\documentclass[twoside,reqno,a4paper,11pt]{memoir}

\setlrmargins{*}{*}{1}
\checkandfixthelayout

\title{\Huge{\textsc{Spatial discretizations of}}\\ \Huge{\textsc{generic dynamical systems}}}
\author{\huge{\textsc{Pierre-Antoine Guihéneuf}}}
\date{\Large{\today}}
\usepackage[latin1]{inputenc}
\usepackage[english,french]{babel}
\usepackage[T1]{fontenc}
\usepackage{amsfonts}
\usepackage{amsmath}
\usepackage{amssymb}
\usepackage{stmaryrd}
\usepackage{amsthm}
\usepackage{enumerate}
\usepackage{pgf,tikz}
\usetikzlibrary{decorations.pathreplacing, shapes.multipart, arrows, matrix, shapes}
\usepackage{color}
\usepackage{float}
\usepackage{makeidx}\makeindex
\usepackage[english,french]{minitoc}
\maxsecnumdepth{subsection}
\usepackage[mathcalasscript,frenchstyle]{kpfonts} 
\usepackage{calc}
\usepackage[columns=3]{idxlayout}
\usepackage{tocloft}
\usepackage[pdftex,colorlinks=true,linkcolor=blue,citecolor=blue,urlcolor=blue]{hyperref}

\DeclareRobustCommand{\SkipTocEntry}[5]{}
\newtheorem{lemme}{Lemma}[chapter]
\newtheorem{theoreme}[lemme]{Theorem}
\newtheorem{prop}[lemme]{Proposition}
\newtheorem{coro}[lemme]{Corollary}
\newtheorem{conj}[lemme]{Conjecture}
\newtheorem{add}[lemme]{Addendum}
\newtheorem{theorem}{Theorem}
\newtheorem{propo}[theorem]{Proposition}
\newtheorem*{conjec}{Conjecture}
\newtheorem*{ques}{Question}
 
\newtheorem{theo}{Theorem}
\newtheorem{proposition}[theo]{Proposition}

\theoremstyle{definition}
\newtheorem{definition}[lemme]{Definition}
\newtheorem{notation}[lemme]{Notation}

\theoremstyle{remark}
\newtheorem{rem}[lemme]{Remark}
\newtheorem{app}[lemme]{Application}
\newtheorem{ex}[lemme]{Example}

\newcounter{counconst}
\newenvironment{constat}{\bigskip\noindent\stepcounter{counconst}\sffamily\thecounconst)}{\bigskip}

\newcommand{\C}{\mathbf{C}}

\newcommand{\N}{\mathbf{N}}
\newcommand{\Pb}{\mathcal{P}}
\newcommand{\R}{\mathbf{R}}
\newcommand{\Sn}{\mathfrak{S}}
\newcommand{\T}{\mathbf{T}}

\newcommand{\V}{\mathcal{V}}
\newcommand{\Q}{\mathbf{Q}}
\newcommand{\Z}{\mathbf{Z}}
\newcommand{\varep}{\varepsilon}
\newcommand{\Hom}{\operatorname{Homeo}}

\newcommand{\Leb}{\operatorname{Leb}}
\newcommand{\Diff}{\operatorname{Diff}}
\newcommand{\ud}{\,\operatorname{d}}
\newcommand{\Prb}{\mathcal{P}}
\newcommand{\tT}{\mathcal{T}}

\newcommand{\card}{\operatorname{Card}}
\newcommand{\dist}{\operatorname{dist}}
\newcommand{\Obs}{\operatorname{Obs}}
\newcommand{\diam}{\operatorname{diam}}
\newcommand{\conv}{\operatorname{conv}}
\newcommand{\1}{\mathbf 1}

\newcommand{\Id}{\operatorname{Id}}
\newcommand{\im}{\operatorname{im}}
\newcommand{\Vect}{\operatorname{span}}
\newcommand{\Sp}{\mathbf{S}}
\newcommand{\End}{\operatorname{End}}
\newcommand{\ind}{{\boldsymbol{i}}}
\newcommand{\Disc}{\operatorname{Disc}}
\newcommand{\len}{\operatorname{length}}
\newcommand{\fat}{\operatorname{father}}
\newcommand{\Ll}{\mathcal{L}}
\newcommand{\Lip}{\operatorname{Lip}}
\newcommand{\E}{\mathbf{E}}

\addto{\captionsfrench}{%
  
}

\makeatletter

\newskip\@bigflushglue \@bigflushglue = -100pt plus 1fil

\makechapterstyle{VZ21}{

\renewcommand\chaptitlefont{\huge\bfseries\centering}
\renewcommand\printchaptertitle[1]{%
\setlength\tabcolsep{7pt}
\settowidth\@tempdimc{\chaptitlefont ##1}%
\setlength\@tempdimc{\textwidth-\@tempdimc-2\tabcolsep}%
\chaptitlefont
\ifdim\@tempdimc > 0pt\relax
\begin{tabular}{c}
\toprule ##1\\ \bottomrule
\end{tabular}
\else
\begin{tabular}{%
>{\chaptitlefont\arraybackslash}p{\textwidth-2\tabcolsep}}
\toprule ##1\\ \bottomrule
\end{tabular}
\fi
\vspace{15pt}
}}

\renewcommand\parttitlefont{\Huge\bfseries\centering}

\renewcommand\printparttitle[1]{%
\setlength\tabcolsep{7pt}
\settowidth\@tempdimc{\parttitlefont #1}%
\setlength\@tempdimc{\textwidth-\@tempdimc-2\tabcolsep}%
\parttitlefont
\ifdim\@tempdimc > 0pt\relax
\begin{tabular}{c}
\toprule #1\\ \bottomrule
\end{tabular}
\else
\begin{tabular}{%
>{\parttitlefont\arraybackslash}p{\textwidth-2\tabcolsep}}
\toprule #1\\ \bottomrule
\end{tabular}
\fi
\vspace{15pt}
}

\makeatother

\begin{document}
\frontmatter
\sloppy

\doparttoc

\clearpage

\thispagestyle{empty}
\mainmatter
\maketitle

\newpage
\chapterstyle{VZ21}

\nouppercaseheads

\let\oldcleartorecto\cleartorecto 
\let\cleartorecto\newpage         

\chapter*{Résumé}

Dans quelle mesure peut-on lire les propriétés dynamiques (quand le temps tend vers l'infini) d'un système sur des simulations numériques ? Pour tenter de répondre à cette question, on étudie dans cette thèse un modèle rendant compte de ce qui se passe lorsqu'on calcule numériquement les orbites d'un système à temps discret $f$ (par exemple un homéomorphisme). L'ordinateur travaillant à précision numérique finie, il va remplacer $f$ par une \emph{discrétisation spatiale} de $f$, notée $f_N$ (où l'ordre de la discrétisation $N$ rend compte de la précision numérique). On s'intéresse en particulier au comportement dynamique des applications finies $f_N$ pour un système $f$ \emph{générique} et pour l'ordre $N$ tendant vers l'infini, où générique sera à prendre dans le sens de Baire (principalement parmi des ensembles d'homéomorphismes ou de $C^1$-difféomorphismes).

La première partie de cette thèse est consacrée à l'étude de la dynamique des discrétisations $f_N$ lorsque $f$ est un homéomorphisme conservatif/dissipatif générique d'une variété compacte. L'étude montre qu'il est illusoire de vouloir retrouver la dynamique du système de départ $f$ à partir de celle d'une seule discrétisation $f_N$ : la dynamique de $f_N$ dépend fortement de l'ordre $N$. Pour détecter certaines dynamiques de $f$ il faut considérer l'ensemble des discrétisations $f_N$, lorsque $N$ parcourt $\N$.

La seconde partie traite du cas linéaire, qui joue un rôle important dans l'étude du cas des $C^1$-difféomorphismes génériques, abordée dans la troisième partie de cette thèse. Sous ces hypothèses, on obtient des résultats similaires à ceux établis dans la première partie, bien que plus faibles et de preuves plus difficiles.

\selectlanguage{english}
\chapter*{Abstract}

How is it possible to read the dynamical properties (ie when the time goes to infinity) of a system on numerical simulations? To try to answer this question, we study in this thesis a model reflecting what happens when the orbits of a discrete time system $f$ (for example an homeomorphism) are computed numerically . The computer working in finite numerical precision, it will replace $f$ by a \emph{spacial discretization} of $f$, denoted by $f_N$ (where the order $N$ of discretization stands for the numerical accuracy). In particular, we will be interested in the dynamical behaviour of the finite maps $f_N$ for a \emph{generic} system $f$ and $N$ going to infinity, where generic will be taken in the sense of Baire (mainly among sets of homeomorphisms or $C^1$-diffeomorphisms).

The first part of this manuscript is devoted to the study of the dynamics of the discretizations $f_N$, when $f$ is a generic conservative/dissipative homeomorphism of a compact manifold. We show that it would be mistaken to try to recover the dynamics of $ f $ from that of a single discretization $f_N$ : its dynamics strongly depends on the order $N$. To detect some dynamical features of $f$, we have to consider all the discretizations $f_N$ when $N$ goes through $\N$.

The second part deals with the linear case, which plays an important role in the study of $C^1$-generic diffeomorphisms, discussed in the third part of this manuscript. Under these assumptions, we obtain results similar to those established in the first part, though weaker and harder to prove.

\newpage

\itshape
This manuscript is an improved version of the thesis \cite{Gui15c} of the author. Some papers have been extracted from it:
\begin{itemize}
\item Chapters~\ref{ChapDissip} and \ref{ChapCons} resume largely the content of \cite{Guih-discr}. Compared to this article, some statements have been improved, others have been added, some misprints have been corrected, additional simulations have been inserted\dots
\item Chapter~\ref{ChapRot} it is almost identical to the article \cite{Guih-rot}.
\item The complete proof of Theorem~\ref{DnZeroIntro} can be found in the article \cite{Gui15a}, which also includes a proof of Theorem~\ref{ConjIntro} concerning the linear case.
\item Still concerning the $C^1$ generic case, the article \cite{Gui15c}aims to prove Theorem~\ref{TheoMesIntroC1}.
\item For its part, the content of Chapter is largely resumed in \cite{Gui15d} and in the permanent preprint \cite{GM} in collaboration with Y. Meyer.
\item A small article \cite{Gui15b} deals with Theorem~\ref{AnswerConjIsom}.
\item Finally, a paper concerning Theorem~\ref{MinkAlm} is in preparation, in collaboration with \'E. Joly \cite{AcEmilien}.
\end{itemize}
This text is aimed to evove with the progression of research in the subject.
\bigskip

The author would like to thank again François Béguin, \'Etienne Ghys, Enrique Pujals, Valérie Berthé, Jérôme Buzzi, Sylvain Crovisier, \'Emilien Joly and Yves Meyer for their input to this manuscript.
\upshape

\let\cleartorecto\oldcleartorecto

\newpage

\selectlanguage{english}
\mtcselectlanguage{english}

\tableofcontents*

\newpage

\chapter{Introduction}

\section{Presentation of the problem}

The goal of this manuscript is to study the dynamics of spatial discretizations of discrete-time dynamical systems. This problem is motivated by numerical considerations: take a dynamics $(X,f)$, where $X$ is a ``good'' space of configurations (the reader can think about a bounded domain of $\R^n$ or the torus $\T^n$), and $f : X\to X$ is a continuous map, given by an explicit formula. We want to understand the dynamics of $f$, i.e. the asymptotical behaviour of the iterates $f^k$ of $f$, by the help of simulations. To do that, the simplest idea consists in taking a point $x\in X$, and asking the computer to compute the images $f^t(x)$, for a time $t$ ``large enough''.

Let us analyse what happens: the computer works with a finite numerical accuracy, for example 10 decimal places\footnote{This is not exactly what happens in general, as the most common number format is floating point. In this format, the number of digits is fixed (instead of the number of decimal places). In our work, we will neglect this aspect.}. So it does not make the computations in the phase space $X$ but rather in a discrete space $E_{10}$, which is the set of points of $X$ whose coordinates have at most $10$ decimal places. The computer also replaces the map $f$ by a \emph{discretization} $f_{10}: E_{10}\to E_{10}$ of $f$, which is an approximation of the dynamics $f$ on the discrete space $E_{10}$. This finite space is a quite good approximation of the continuous space $X$. However, it may happen that the small roundoff errors made at each iteration of $f_{10}$ add up, so that after a while the orbits of a point $x$ under the map $f$ and it discretization $f_{10}$ become very different. Thus, if nothing more is known about the dynamics $f$, one can not deduce the asymptotic behaviour of the orbit of $x$ under $f$ from its orbit under the discretization $f_{10}$, even if a numerical accuracy of 10 numerical places can seem quite good.

Nevertheless, this discrete orbit $\big(f_{10}^t(x)\big)_{t\in\N}$ is close to a true orbit of $f$ on every segment of reasonable length. So one can hope that the \emph{collective} behaviour of the computed orbits -- where by collective we mean that we consider a lot of starting points -- reflects the \emph{global} dynamics of $f$.

Our goal here is to analyse this naive algorithm, which describes what happens when calculi are performed without any precaution. In other words, we want to understand the dynamical phenomena appearing when the phase space of a dynamical system is discretized; we would like to rely the asymptotical behaviour of $f$ (i.e. when the times goes to infinity) with that of the discrete map $f_{10}$. Thus, we will not try to find the best algorithm that allows to detect some dynamical features of some dynamics $f$; our aim will not be to conduct numerical studies aiming to decide if the real dynamics of $f$ is observed in practice either.

Note that this problematic is very different from that of classical numerical analysis. Somehow, given a time $t_0$ and a precision $\delta$, numerical analysts determine a method (in particular a grid) that allows to approach the solution during a time $t_0$ and with a precision $\delta$. Here, we fix a fine enough grid and compare asymptotical behaviours of the real system and the discretized system on this grid.
\bigskip

In this manuscript, we will consider discrete-time dynamical systems $f: X\to X$, where $X$ is a compact manifold (with or without boundary) endowed with a metric $d$, and $f$ is a countinuous map from $X$ into itself. The operation of discretization will be modelled in the following way. Consider a sequence of grids $(E_N)_{N\in\N}$ which are finite subsets of $X$, whose mesh tends to 0\footnote{I.e. for every $\varep>0$ and every $N$ large enough, the grid $E_N$ is $\varep$-dense in $X$.}. We define a projection $P_N : X\to E_N$ which maps every point of $X$ on the closest point\footnote{If there are many closest points, we choose one of them once and for all.} of $E_N$. The \emph{discretization} of a point $x\in X$ is then defined as the point $P_N(x)$, and the \emph{discretization]} of a map $f$ as the map
\[\begin{array}{rcl}
f_N : E_N & \longrightarrow & E_N\\
    x   & \longmapsto     & P_N\big( f(x)\big).
\end{array}\]
To compute $f_N(x)$, we compute $f(x)$ and project it on $E_N$ \emph{via} the projection $P_N$. Remark that this model of the roundoff error is imperfect: in practice, the computer uses also finite numerical precision when doing intermediate calculus to compute $f(x)$; we will neglect this aspect.

We can now specify a little more the question we are interested in.

\begin{ques}
What dynamical properties of $f$ can be read on the (global) dynamics of the discretizations $(f_N)_{N\in \N}$?
\end{ques}
\bigskip

\begin{figure}[!ht]
\makebox[\textwidth]{\parbox{\textwidth}{%
\begin{minipage}[c]{.32\linewidth}
	\includegraphics[width=\linewidth]{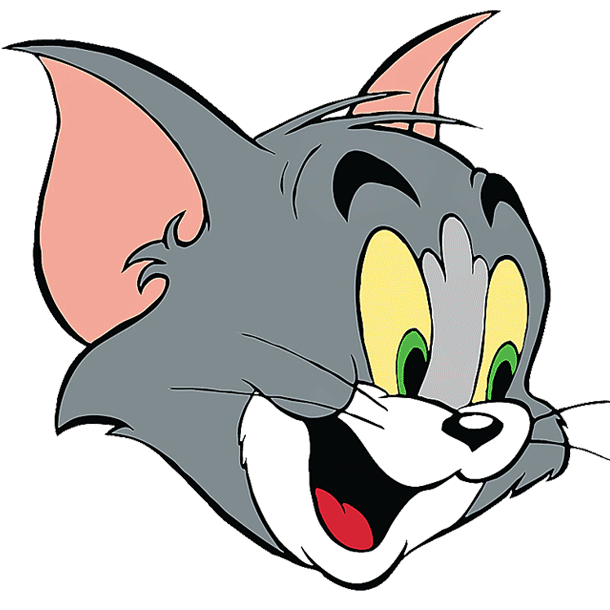}
\end{minipage}\hfill
\begin{minipage}[c]{.32\linewidth}
	\includegraphics[width=\linewidth]{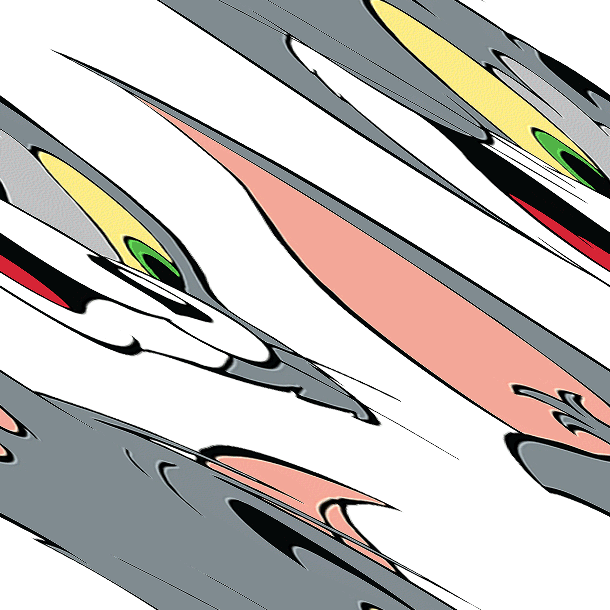}
\end{minipage}\hfill
\begin{minipage}[c]{.32\linewidth}
	\includegraphics[width=\linewidth]{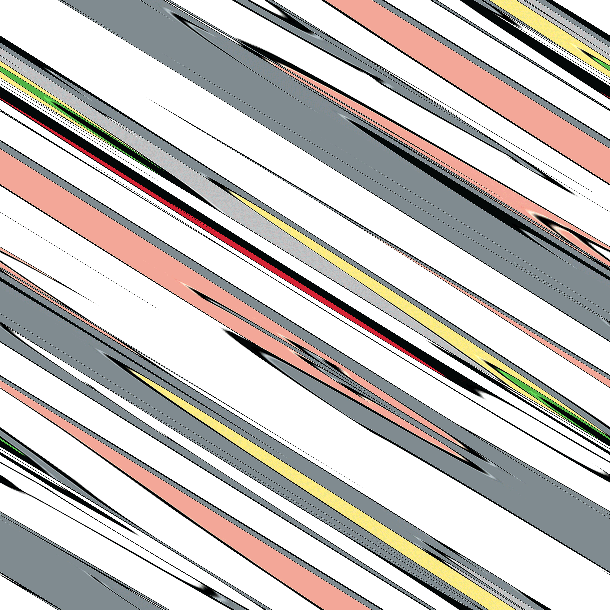}
\end{minipage}
\vspace{.02\linewidth}

\begin{minipage}[c]{.32\linewidth}
	\includegraphics[width=\linewidth]{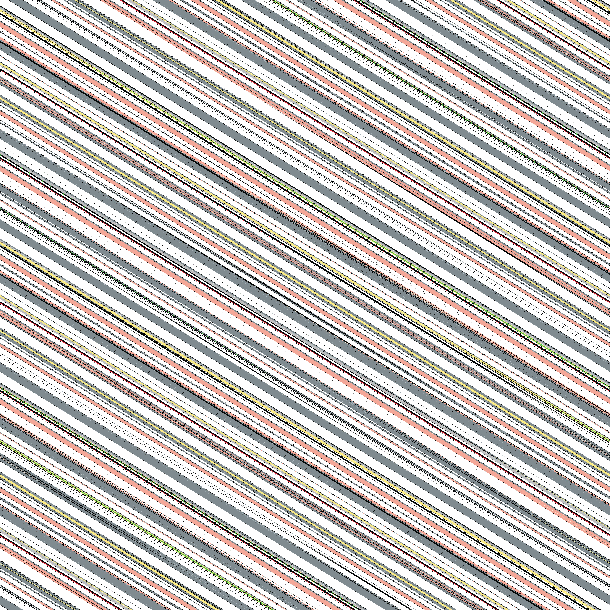}
\end{minipage}\hfill
\begin{minipage}[c]{.32\linewidth}
	\includegraphics[width=\linewidth]{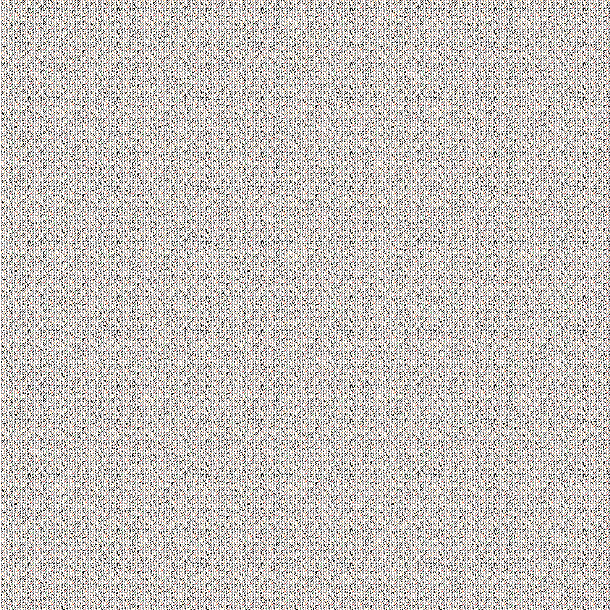}
\end{minipage}\hfill
\begin{minipage}[c]{.32\linewidth}
	\includegraphics[width=\linewidth]{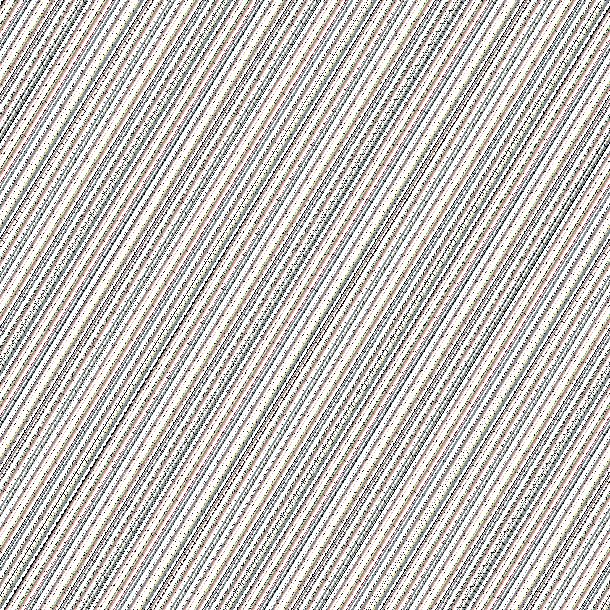}
\end{minipage}
\vspace{.02\linewidth}

\begin{minipage}[c]{.32\linewidth}
	\includegraphics[width=\linewidth]{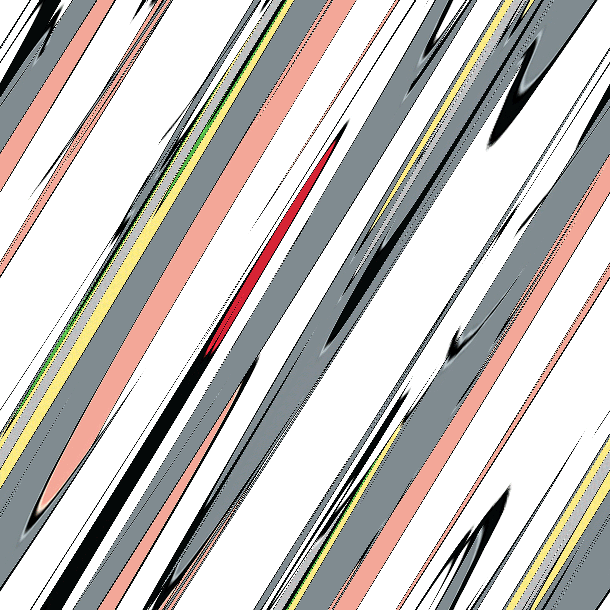}
\end{minipage}\hfill
\begin{minipage}[c]{.32\linewidth}
	\includegraphics[width=\linewidth]{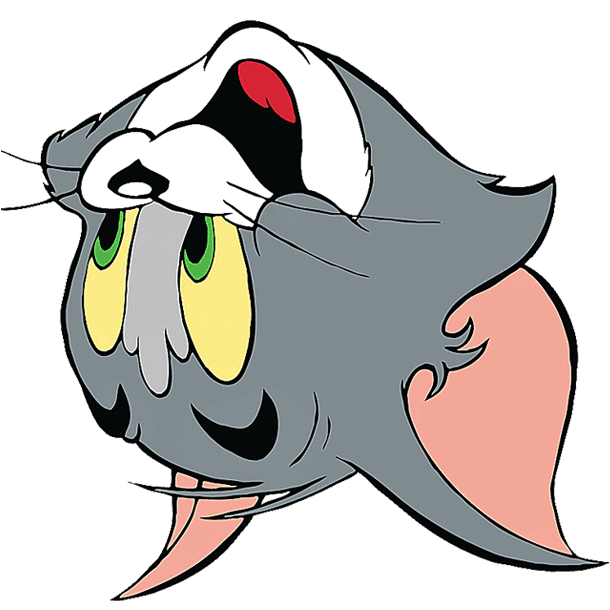}
\end{minipage}\hfill
\begin{minipage}[c]{.32\linewidth}
	\includegraphics[width=\linewidth]{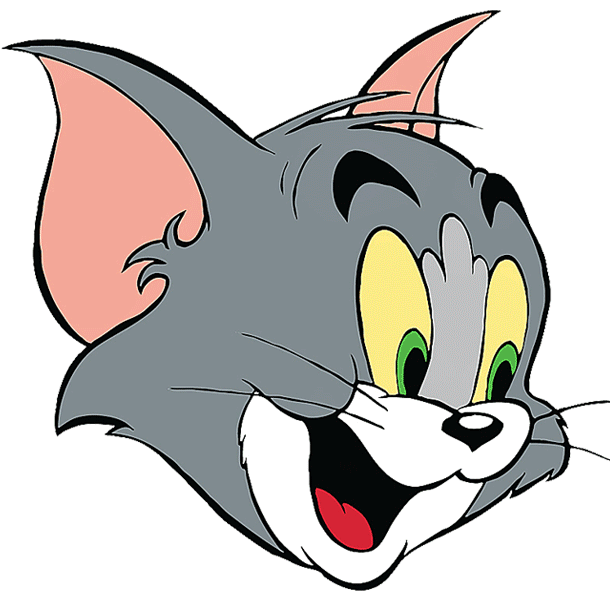}
\end{minipage}
}}
\caption[Discretization of a linear automorphism of the torus]{Images of a cat picture by the discretization of the cat map on the grid of size $610\times 610$. From left to right and top to bottom, times $0$, $2$, $4$, $8$, $14$, $20$, $26$, $30$ et $60$.}\label{Miaou}
\end{figure}

To try to answer this question, the first idea is to study what happens for simple and well-known dynamics. This is done in \'E.~Ghys in \cite{Ghys-vari}, where the author studies the discretizations of the cat map, a linear Anosov automorphism of the torus defined by\label{PageMiaou}
\[\begin{array}{rcl}
A : \R^2/\Z^2 & \longrightarrow & \R^2/\Z^2\\
    (x,y)     & \longmapsto     & (y,x+y),
\end{array}\]
on the grids
\begin{equation}\label{GrillesCano}
E_N = \left\{ \left(\frac{i}{N},\frac{j}{N}\right)\in \R^2/\Z^2 \middle\vert\ 1\le i,j \le N\right\}.
\end{equation}

Figure~\ref{Miaou} shows the iterates of a numerical picture under $A$: such a picture is made of pixels -- here $610\times 610$ --, thus corresponds to a map $\phi : E_{610} \to \llbracket 1,256\rrbracket ^3$. So one can consider its pushforwards by the map $A_{610}$ (recall that $A_N$ denotes the discretization of $A$ on the grid $E_N$). Two very specific phenomena are observed.
\begin{itemize}
\item First, when the picture is iterated, one still obtain pictures. In other words, the colour of each pixel is well defined. In general, there is no reason for the discretizations of a homeomorphism\footnote{Or more generally of a continuous map.} to be bijections. Here, the particular form of $A$ -- it is a linear automorphism with integer coefficients -- forces the discretizations $A_N$ to be permutations of the grids $E_N$. In fact, for the map $A$ there is no need to project: every point $x\in E_N$ is mapped by $A$ to a point of $E_N$.
\item Secondly, the initial map comes back after very few iterations (60 in the example of Figure~\ref{Miaou}). However, the global order of a typical permutation of a set of $N^2$ elements is equivalent to $N^{2\ln N}$ (see for example \cite{Boll-rand}); Applied to $N=610$, this asymptotic gives an order about $5. 10^{35}$: if the map $A_{610}$ were a typical permutation, our poor cat would come back after $10^{35}$ iterations (which is manifestly bigger than $60$)! In fact, the following result has been proved by F.~Dyson and H.~Falk in \cite{MR1176587} : \emph{for every $N\in\N$, the order of the permutation $A_N$ is smaller than $3N$. Moreover, there exists a sequence of integers $(N_k)_{k\in\N}$, tending exponentially to infinity, such that the order of the permutation $A_{N_k}$ is $2k$.}
\end{itemize}

This behaviour of the discretizations contrasts sharply with the actual dynamics of the linear automorphism $A$: this map is a paradigm of exponentially mixing systems. If the dynamics of discretizations reflected that of $A$, the resulting pictures would quickly become uniformly grey (what is observed in the centre of Figure~\ref{Miaou}, at the $14^{\text{th}}$ iteration) and stay forever (unlike what is observed on the last images of Figure~\ref{Miaou}).

It follows from these observations that the example of the linear automorphism $A$ is very particular: the grids $E_N$ somehow resonate with the dynamics $A$, because the grids $E_N$ are regularly spaced, and because $A$ has special arithmetic properties.
\bigskip

To avoid these ``exceptional'' phenomena, \'E.~Ghys proposes (still in \cite{Ghys-vari}) to study the discretizations of \emph{generic} conservative homeomorphisms of the torus; this is what we will do in this manuscript.

\label{DiscussGene}The word ``generic'' has a very precise mathematical sense. Most of reasonable functional spaces are complete (it will always be the case here). In particular, one can apply Baire theorem: \emph{ every countable intersection of open dense sets, is itself dense.} We say that a property (P) is \emph{generic} for the class of functions considered if it is satisfied on at least a countable intersection of open dense sets of such functions. Note that this concept of genericity has two nice properties:
\begin{itemize}
\item a generic property is satisfied on a dense set;
\item if (P) and (Q) are generic, then the property ``satisfying both (P) and (Q)'' is itself generic; this remains true even for a countable family of properties;
\end{itemize}

Thus, we will study the dynamics of discretizations of generic homeomorphisms and $C^1$-diffeomorphisms, both conservative (i.e. under the assumption of preserving a good measure on $X$, fixed a once for all) and dissipative (that is to say without measure preservation assumption).

By abuse of language, we will often use the expression ``generic homeomorphism'', and list its properties. We should have in mind that such a generic homeomorphism does not exist: a bit like in quantum mechanics, where a measurement inevitably perturbs the system studied; try to choose a generic element of a space, thereby it will cease to be. This abuse of language is thus very bad, and when meeting it, the reader should mentally substitute the sentence ``if $f$ is a generic element of this space, then it has the property (P)'' with the more correct one ``the property (P) is generic in this space''.
\bigskip

Let us precise what is meant by ``the dynamics of discretizations''. Every discretization $f_N : E_N\to E_N$ is a map from a finite set into itself. In particular, each of its orbits is pre-periodic, that is, for every point $x\in E_N$, the orbit of $x$ under $f_N$ is eventually periodic: there exist two integers $C$ and $T$ such that for every $k\ge C$, we have $f_N^{k+T}(x) = f_N^k(x)$. Thus, we can partition $E_N$ into two subsets: the \emph{recurrent set} $\Omega(f_N)$ of $f_N$, which is the union of periodic orbits of $f_N$, and its complement called the \emph{wandering set} of $f_N$. Note that the recurrent set $\Omega(f_N)$ is stable under $f_N$, and that the restriction of $f_N$ to this set is a bijection. Moreover, each $x\in E_N$ ``falls'' in $\Omega(f_N)$, i.e. there exists $t\in\N$ such that $f_N^t(x)\in\Omega(f_N)$; and $\Omega(f_N)$ is the smallest subset of $E_N$ having this property.

Thus, from the combinatorial viewpoint, the dynamics of $f_N$ is characterized by a small number of quantities: the number of periodic orbits of $f_N$ and the repartition of their lengths, the sizes of their basins of attraction\footnote{I.e. the number of points whose positive orbit falls into this periodic orbit.}, the stabilization time of $f_N$ (i.e. the smallest $t\in\N$ such that $f_N^t(E_N) = \Omega(f_N)$), etc. We will focus in more detail on the \emph{degree of recurrence}\label{DegRecurIntro} $D(f_N)$ of $f_N$, which is the quotient between the cardinality of the recurrent set $\Omega(f_N)$ and the cardinality of the grid $E_N$.

This combinatorial dynamical properties of $f_N$ does not depend on the geometry of the grids $E_N$. However, we have supposed that grids are good approximations of the space $X$ from the metrical viewpoint. In other words, for every $\varep>0$ and every large enough $N$, the grid $E_N$ is $\varep$-dense in the space $X$. The space $X$ can also be endowed with a good measure $\lambda$, in this case we can also suppose that the uniform measures on the grids $E_N$ tend to $\lambda$. This hypotheses allow to study geometric or ergodic dynamical properties of the discretizations. For instance, one can wonder if the periodic orbits of $f_N$ tend to that of $f$. More generally, one can wonder if the invariant compact subsets of $f_N$ tend (for Hausdorff topology) to the invariant compact subsets of $f$. The same kind of questions can be asked for invariant measures endowed with weak-* topology.
\bigskip

Quite surprisingly, this problem has been only little studied. To my knowledge, apart from works analysing numerically well-chosen examples, or explaining phenomena by heuristic arguments (see Section~\ref{SecHistoGene}), there is very few theoretical works about this problem

Indeed, only P.P.~Flockermann and T. Miernowski have really conducted a systematic study for a large class of systems\footnote{See Chapter~\ref{SecHisto} for a complete historic.}. In his thesis under the supervision of O.E.~lanford, P.P.~Flockermann has considered discretizations of expanding maps of the circle. During his thesis under the supervision of \'E.~Ghys, T. Miernowski has studied the case of circle homeomorphisms (see also \cite{MR2279269}). Basically, he shows that the dynamical properties of a typical (generic or prevalent) circle homeomorphism/diffeomorphism can be read on the dynamics of discretizations. The proofs of all these results use crucially the fact that the phase space is one dimensional (in particular, for homeomorphisms, they depend a lot on the rotation number); thus could not be generalized to other classes of systems.

During his thesis \cite{Mier-dyna}, T.~Miernowski has also studied a bit the case of generic conservative homeomorphisms of the torus. For this purpose, the torus $\T^2 = \R^2/\Z^2$ is endowed with Lebesgue measure and canonical discretization grids (defined by Equation~\eqref{GrillesCano}); the set of homeomorphisms of the torus that preserve Lebesgue measure is denoted by $\Hom(\T^2,\Leb)$. T.~Miernowski shows the following result.

\begin{theo}[Miernowski]\label{TheoMierno}
For a generic conservative homeomorphism $f\in \Hom(\T^2,\Leb)$, there exists a subsequence $f_{N_k}$ of discretizations of $f$, whose elements are permutations of $E_{N_k}$.
\end{theo}

This theorem is somehow the starting point of this thesis: we would like to show other results of this kind, to understand better the dynamics of discretizations of a generic homeomorphism or diffeomorphism.

\section{Some results of this thesis}

This manuscript contains many statements, concerning combinatorial, topological or ergodic properties of discretizations of differentiable or just continuous systems, conservative or dissipative\dots\ In this section, I isolate those which seem the more representative and interesting. For a more systematic presentation of the results, see the introductions of the different parts of this manuscript; also, in Chapter~\ref{Blablablabla}, we will discuss and compare the results obtained in these different contexts.

In this section, for simplicity, we consider the case where $X = \T^2 = \R^2/\Z^2$, endowed if necessary with Lebesgue measure $\Leb$, and with the uniform grids
\begin{equation}\tag{\ref{GrillesCano}}
E_N = \left\{ \left(\frac{i}{N},\frac{j}{N}\right)\ \middle\vert\ 1\le i,j \le N\right\}.
\end{equation}
The results presented in this introduction remain true in more general contexts, explained in the concerning parts of this manuscript.

\subsection{Degree of recurrence}

Let us begin by two theorems concerning the combinatorial dynamics of discretizations, more precisely the degree of recurrence. Recall that the degree of recurrence $D(f_N)$ is defined as the ratio between the cardinality of the recurrent set of $f_N$ (which is the union of the periodic orbits of $f_N$, see page~\pageref{DegRecurIntro}) and that of the grid $E_N$. Since the positive orbit of any point of $E_N$ under the discretization $f_N$ falls in this recurrent set, the degree of recurrence is also equal to  $\card\big(f_N^t(E_N)\big) / \card(E_N)$ for every $t$ large enough.

The first result concerns the behaviour of $D(f_N)$ for a generic conservative\footnote{Recall that ``conservative'' means ``that preserve Lebesgue measure''.} homeomorphism: in Chapter~\ref{ChapCons}, we prove the following result (Corollary~\ref{ConjEt} page~\pageref{ConjEt}).

\begin{theo}\label{IntroAccuDN}
For a generic conservative homeomorphism $f\in\Hom(\T^2,\Leb)$, the sequence $(D(f_N))_{N\ge 0}$ accumulates on the whole segment $[0,1]$.
\end{theo}

The accumulation on $1$ is a trivial corollary of Theorem~\ref{TheoMierno} of T.~Miernowski. Considered independently, Theorem~\ref{TheoMierno} can be seen as very positive: it expresses that an infinite number of discretizations behaves the same way as the homeomorphism, namely is a measure-preserving bijection. In fact, this phenomenon is only a particular case of a very erratic behaviour of the sequence $\big(D(f_N)\big)_{N\ge 0}$: Theorem~\ref{IntroAccuDN} asserts that it accumulates on the biggest set on which it can a priori accumulate.
\bigskip

We then study the case of generic conservative $C^1$-diffeomorphisms. We show the following theorem in Chapter~\ref{ChapDeg} (Theorem~\ref{limiteEgalZero} page~\pageref{limiteEgalZero}).

\begin{theo}\label{DnZeroIntro}
For a generic conservative $C^1$-diffeomorphism $f\in\Diff^1(\T^2,\Leb)$, we have
\[\lim_{N\to +\infty} D(f_N) = 0.\]
\end{theo}

This theorem has to be compared with Theorem~\ref{IntroAccuDN}: it expresses that the global behaviour of discretizations evolves less irregularly when $N$ grows than in the case of homeomorphisms. Nevertheless, the fact that the degree of recurrence tends to 0 means that there is an arbitrarily large loss of information when $N$ is large enough.

To prove the theorem, we link the macroscopic and mesoscopic behaviours of discretizations of a generic diffeomorphism. To begin with, we define the \emph{rate of injectivity} of~$f$
\[\tau_N^t(f) = \frac{\card(f_N^t(E_N))}{\card(E_N)}.\]
By definition, their limit as $t$ tends to infinity is the degree of recurrence $D(f_N)$. Therefore we are led to study these rates of injectivity; in particular we prove a result that connects them to similar quantities concerning the differentials of $f$ (Theorem~\ref{conv} page~\pageref{conv}, see also Theorem~\ref{convBisMieux} page~\pageref{convBisMieux}).

\begin{theo}\label{LocGlobIntro}
For every $r\ge 1$, and for a generic conservative $C^r$-diffeomorphism $f\in\Diff^r(\T^2,\Leb)$, for every $t\in{\N^*}$, we have
\[\lim_{N\to +\infty} \tau_N^t(f) = \int_{\T^2} \tau(D f_{f^{t-1}(x)},\cdots,D f_x) \ud x,\]
where the rate of injectivity of a sequence of matrices is defined similarly to that of a diffeomorphism (see Definition~\ref{DefTaux} page~\pageref{DefTaux}).
\end{theo}

This result allows to reduce the study of the degree of recurrence of a generic diffeomorphism to that of the rate ofinjectivity of a  generic sequence of matrices; this is the subject of Part~\ref{PartII} of this manuscript.

\subsection{Physical measures}

We now turn to the ergodic dynamics of discretizations. To begin with we define the concept of physical measure. We set $\mu_{x,T}^f$ the uniform probability measure on the segment of orbit $\big(x,f(x),\cdots,f^{T-1}(x)\big)$. This allows to set the \emph{basin of attraction} of a probability measure $\mu$ as the set of points $x\in\T^2$ such that the measures $\mu_{x,T}^f$ tend to $\mu$ when $T$ tends to $+\infty$ (in the sense of the weak-* topology). The measure $\mu$ is then called \emph{physical} for $f$ if its basin of attraction has strictly positive Lebesgue measure. A physical measure is a measure that should be detected during physical experiments because it is ``seen'' by a set of points $x$ with positive Lebesgue measure. The question is whether it is possible to identify such measures on numerical experiments.

Note that in this spirit, there are many results on the \emph{stochastic stability} of physical measures, to take into account the fact that in a physical experiment there is always noise (see for example \cite{VianaStoch}). In many cases, it turns out that despite the noise, one can recover the physical measures of the starting dynamics. These studies always assume that the noise is random, and especially independent at each iteration, which is far from being the case of digital truncation. The question is whether such results persists for discretizations. As $f_N$ are finite maps, the measures $\mu_{x,T}^{f_N}$ converge to the 	uniform measure on the periodic orbit in which falls the positive orbit of $x_N$ under $f_N$. We denote $\mu_{x}^{f_N}$ this measure. The goal is to characterize the behaviour of these measures for ``most'' of the points $x\in E_N$. In this regard, we prove the following result (Theorem~\ref{EnsMesInvSimpl} page~\pageref{EnsMesInvSimpl}).

\begin{theo}\label{TheoMesIntro}
For a generic conservative homeomorphism $f\in \Hom(\T^2,\Leb)$, for every $f$-invariant probability measure $\mu$, there exists a subsequence of discretizations $(f_{N_k})_k$ such that $f_{N_k}$ has a unique $f_N$-invariant measure $\mu_k$, tending to $\mu$. In other words, for every $f$-invariant probability measure $\mu$, there exists a sequence of integers $(N_k)_{k\ge 0}$ such that, for every $x\in \T^2$, we have
\[\mu_{x}^{f_{N_k}}\underset{k\to+\infty}{\longrightarrow}  \mu.\]
\end{theo}

Note that the theorem of Oxtoby-Ulam (one of the first results of genericity, see \cite{Oxto-meas}) implies that a generic conservative homeomorphism $f\in \Hom(\T^2,\Leb)$ do has a unique physical measure: Lebesgue measure. Theorem~\ref{TheoMesIntro} expresses that somehow, any $f$-invariant measure is ``physical for an infinity of discretizations'' (it also responds to the question raised by M.~Blank in \cite[p. 114]{MR1440853}, in the case of generic conservative homeomorphisms). This result can be considered as positive or negative depending on the point of view: one can find all the $f$-invariant measures, but for now nothing allows us to detect the physical measures of $f$ on its discretizations. Recall that the very definition of physical measures is supposed to detect the measures that are observable in practice.

\begin{figure}[t]
\begin{minipage}[c]{.31\linewidth}
	\includegraphics[height=4.8cm,trim = 1.5cm .95cm 2.8cm .5cm,clip]{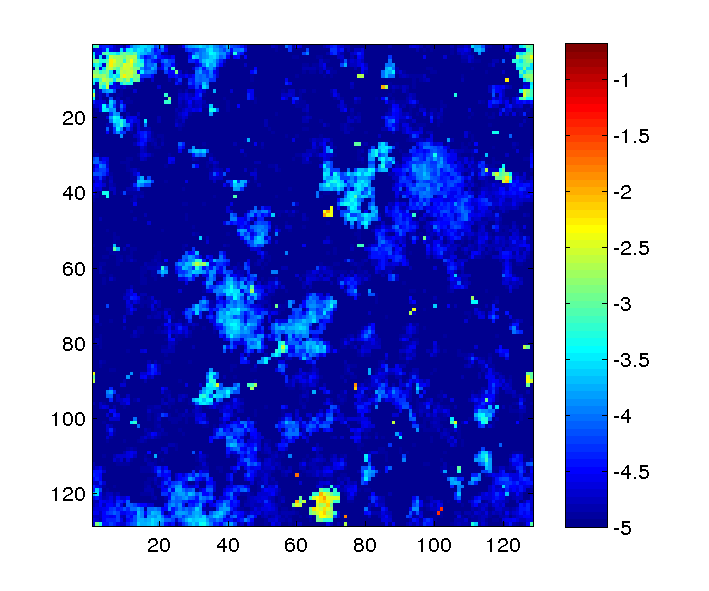}
\end{minipage}\hfill
\begin{minipage}[c]{.31\linewidth}
	\includegraphics[height=4.8cm,trim = 1.5cm .95cm 2.8cm .5cm,clip]{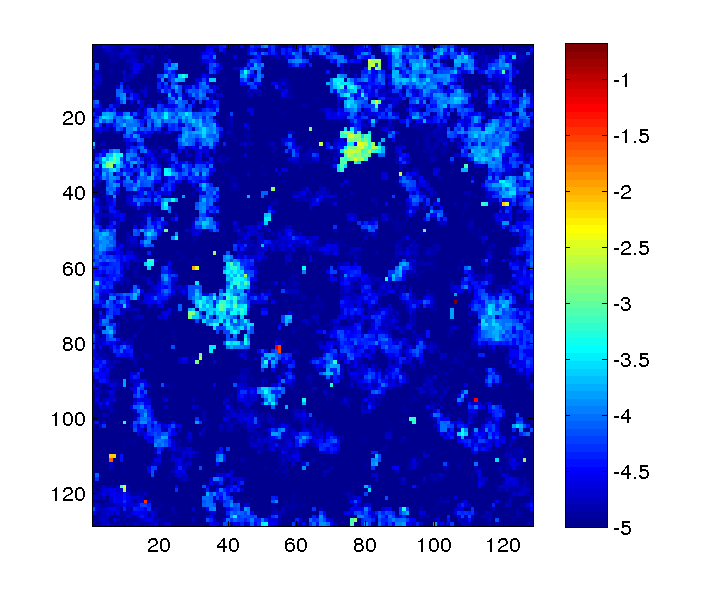}
\end{minipage}\hfill
\begin{minipage}[c]{.37\linewidth}
	\includegraphics[height=4.8cm,trim = 1.5cm .95cm 1cm .5cm,clip]{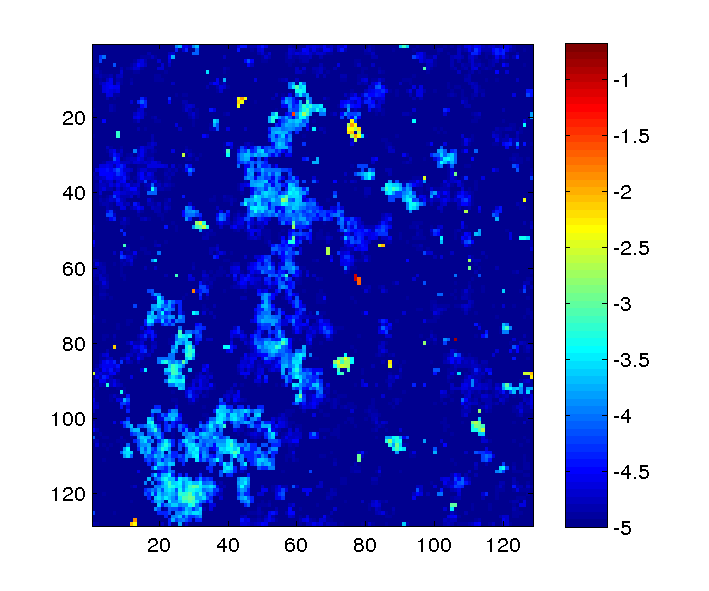}
\end{minipage}
\caption[Density of $\mu_{\T^2}^{f_N}$, for $f$ a $C^0$-perturbation of identity]{Density of the canonical invariant measure $\mu_{\T^2}^{f_N}$ of $f_N$, obtained as a Cesàro limit of the pushforwards by $f_N$ of the uniform measure on $E_N$. The density is represented in logarithmic scale: an orange pixel bears a measure of about $10^{-2}$. The homeomorphism $f$ is a small $C^0$-perturbation of identity (see page~\pageref{PageDefSimulCons}). From left to right, $N = 20\,000,\,20\,001,\,20\,002$.}\label{FigMesIntro}
\end{figure}

This phenomenon even appears on numerical simulations of conservative homeomorphisms: Figure~\ref{FigMesIntro} represents a numerical simulation of measures $\mu_{\T^2}^{f_N}$, where $\mu_{\T^2}^{f_N}$ is the Cesàro limit of the pushforwards by $f_N$ of the uniform measure on $E_N$. We can observe that on this example the measures $\mu_{\T^2}^{f_N}$ do not converge towards Lebesgue measure at all. In addition, they have nothing to do the ones with the others, even if the discretization orders are very close.
\bigskip

In Chapter~\ref{chapPhys} we obtain a similar result for $C^1$-diffeomorphisms. Unfortunately, this result is not strong enough to explain the behaviours of the measures $\mu_{x}^{f_N}$ for all the points $x\in\T^2$; it only express what happens for a generic set of points (Theorem~\ref{TheoMesPhysDiff} page~\pageref{TheoMesPhysDiff}).

\begin{figure}[t]
\begin{minipage}[c]{.31\linewidth}
	\includegraphics[height=4.8cm,trim = 1.5cm .95cm 2.8cm .5cm,clip]{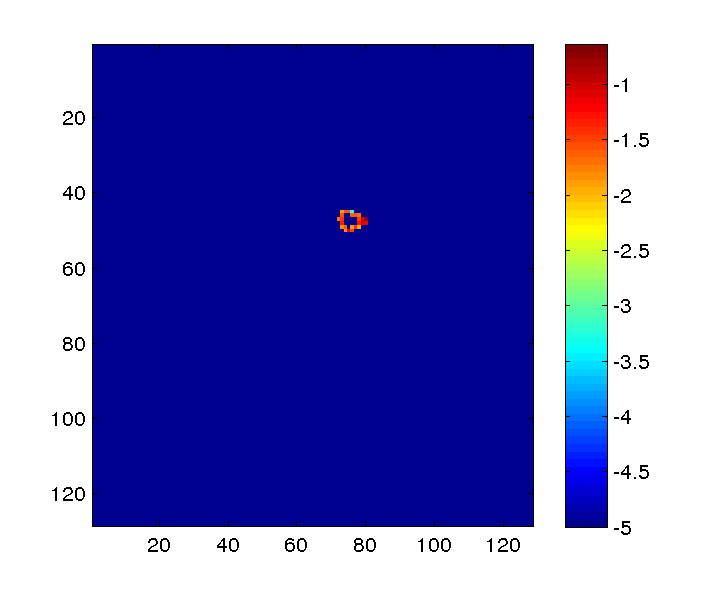}
\end{minipage}\hfill
\begin{minipage}[c]{.31\linewidth}
	\includegraphics[height=4.8cm,trim = 1.5cm .95cm 2.8cm .5cm,clip]{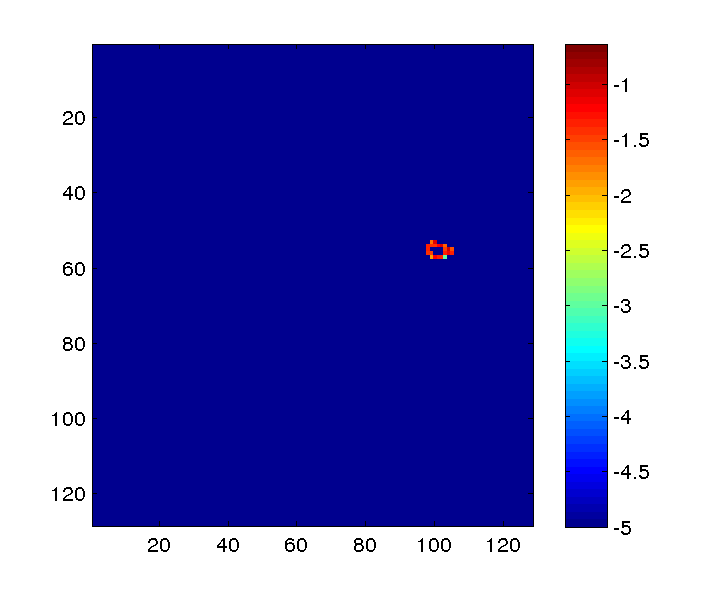}
\end{minipage}\hfill
\begin{minipage}[c]{.37\linewidth}
	\includegraphics[height=4.8cm,trim = 1.5cm .95cm 1cm .5cm,clip]{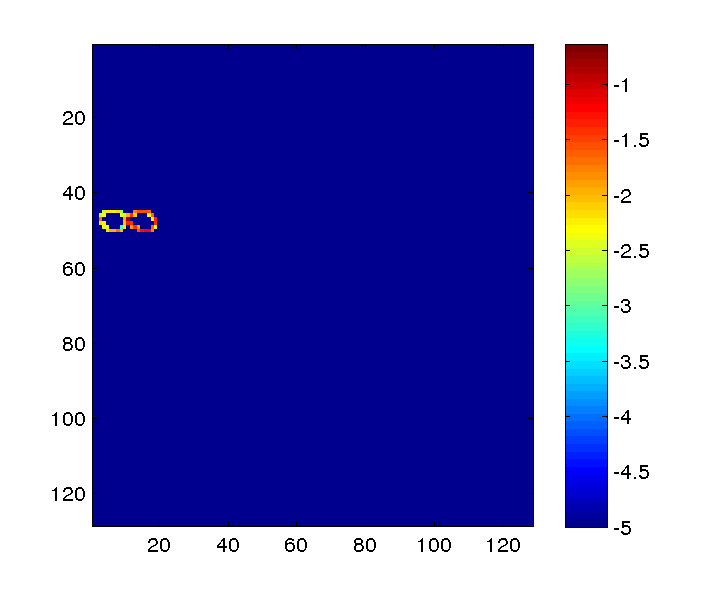}
\end{minipage}
\caption[Density of $\mu_{x}^{f_N}$, for $f$ a $C^1$-perturbation of identity]{Density of the invariant measure $\mu_{x}^{f_N}$ of $f_N$, obtained as a Cesàro limit of the pushforwards by $f_N$ of the measure $\delta_x$, with $x=(1/2,1/2)$. The density is represented in logarithmic scale: an orange pixel bears a measure of about $10^{-2}$. The homeomorphism $f$ is a small $C^1$-perturbation of identity (see section~\ref{NumSimPhys}). From left to right, $N = 2^{20}+1,2^{20}+2,2^{20}+3$.}\label{FigMesIntro2}
\end{figure}

\begin{theo}\label{TheoMesIntroC1}
For a generic conservative $C^1$-diffeomorpihsm $f\in\Diff^1(\T^2,\Leb)$, for a generic point $x\in \T^2$, for every $f$-invariant measure $\mu$ there exists a subsequence $(f_{N_k})_k$ of discretizations such that
\[\mu_{x}^{f_{N_k}} \underset{k\to+\infty}{\longrightarrow} \mu.\]
\end{theo}

Compared to Theorem~\ref{TheoMesIntro}, notice that the sequence $({N_k})_k$ now depends on the point $x$. Still, if one chooses a conservative diffeomorphism $f$ and a point $x\in\T^2$, and we compute numerically the uniform measure on the orbit $x,f(x),\cdots,f^T(x)$ for a large $T$, then Theorem~\ref{TheoMesIntroC1} expresses it may be that one does not find at all a physical measure of $f$, but rather any invariant measure of $f$. This phenomenon can be observed in practice, as shown by numerical simulations of the measures $\mu^{f_N}_x$ (see Figure~\ref{FigMesIntro2}).

The proof of this result is rather long and technical. It uses crucially a linear statement (shown at the end of Part~\ref{PartII}) quite close to the one used to prove Theorem~\ref{DnZeroIntro}. It also uses two classical closing lemmas: the connecting lemma for pseudo-orbits of C.~Bonatti and S.~Crovisier \cite{MR2090361}, and an improvement of an ergodic closing lemma of  F.~Abdenur, C.~Bonatti and S.~Crovisier \cite{MR2811152}.

\subsection{Detecting rare dynamics (from the measure viewpoint)}

In Chapter~\ref{ChapRot}, we study the computation of the rotation set of a generic conservative torus homeomorphism. We take advantage of a phenomenon illustrated by Theorem~\ref{TheoMesIntro}: the discretizations of a generic conservative homeomorphism allow to find all the invariant measures of the homeomorphism, not just its physical ones. Thus we show that, paradoxically, numerical errors are useful (or even necessary) to compute the rotation set of a generic conservative homeomorphism.

The rotation set of a torus homeomorphism is a compact and convex subset of $\R^2$ defined modulo $\Z^2$ (see its definition page~\ref{PageDefRho}). This is a generalization to two dimensions of the concept of rotation number of a circle homeomorphism: this set describes how fast the orbits wind around the torus. For a generic conservative homeomorphism/$C^1$-diffeomorphism, it has nonempty interior (see Propositions~\ref{RotGeneCons} and \ref{RotGeneConsC1} pages~\pageref{RotGeneCons} and \pageref{RotGeneConsC1}).

We will define the \emph{observable rotation set}: a vector $v$ belongs to the observable rotation set if, for every $\varep>0$, there exists a positive Lebesgue measure set of points $x$ whose orbit has a rotation vector $\varep$-close to $v$ (see Definition~\ref{DefRotObs} page~\pageref{DefRotObs}). This definition is meant to represent the rotation set that would be obtained by making exact computations, but by being allowed to make only a finite number of such computations. We will show the following theorem (Theorem~\ref{CoroRotDiscrCons} page~\pageref{CoroRotDiscrCons}).

\begin{theo}\label{EnsRotIntro}
For a generic conservative homeomorphism $f\in\Hom(\T^2,\Leb)$,
\begin{itemize}
\item the observable rotation set is reduced to a point;
\item the superior limit of the rotation sets of discretizations coincides with the rotation set of $f$.
\end{itemize}
\end{theo}

The second assertion of this theorem remains true for a generic conservative $C^1$-diffeomorphism, considering the convex hulls of the rotation sets of discretizations (see Theorem~\ref{DiscrC1} page~\pageref{DiscrC1}). The first assertion remains also true as long as Lebesgue measure is ergodic for $f$; it is conjectured that this is the case for a generic conservative $C^1$-diffeomorphism (and there are open sets of diffeomorphisms in which a generic element is ergodic, see page~\pageref{AvilaErgo}).

\begin{figure}
\makebox[\textwidth]{\parbox{\textwidth}{%
\begin{center}
\begin{minipage}[c]{.35\linewidth}
	\includegraphics[width=\linewidth,trim = .5cm .3cm .6cm .1cm,clip]{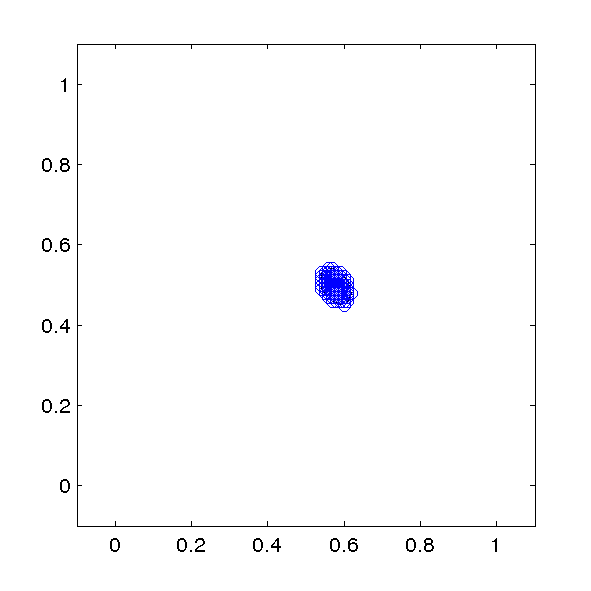}
\end{minipage}\hspace{1.8cm}
\begin{minipage}[c]{.35\linewidth}
	\includegraphics[width=\linewidth,trim = .5cm .3cm .6cm .1cm,clip]{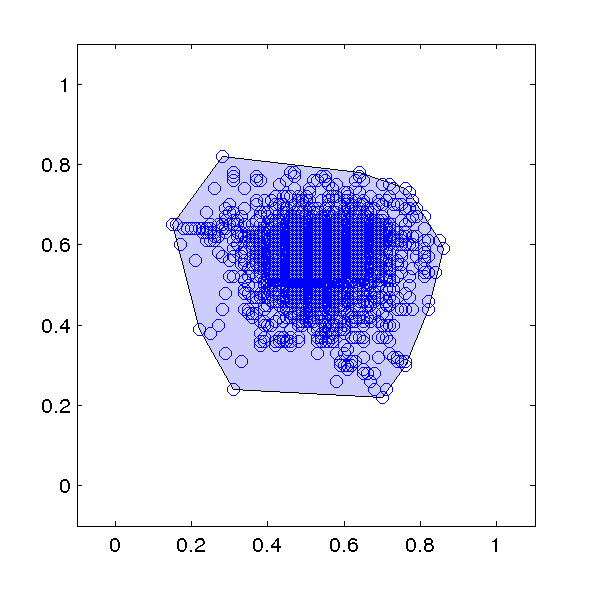}
\end{minipage}
\end{center}
}}
\caption[Comparison of methods of computation of the rotation set]{Comparison of methods of computation of the rotation set on the homeomorphism $f$ (see its definition in section~\ref{SecNumRot}). On the left, the set obtained using the natural algorithm consisting in calculating rotation vectors of very long orbits. On the right, the set obtained by accumulating the rotation sets of discretizations on different grids. Calculation times are similar (about 4 hours).}\label{FigRotIntro}
\end{figure}

This theorem gives an algorithm to approximate the rotation set, consisting in calculating the rotation sets of discretizations. In Figure~\ref{FigRotIntro}, we see that the ``natural'' algorithm only detect a very small part of the whole rotation set, unlike the algorithm using discretizations (as predicted by Theorem~\ref{EnsRotIntro}). In Chapter~\ref{ChapRot}, we will see on other examples that this algorithm is quite efficient. This is an unexpected application of the concept of discretization: truncation, which seems a priori to be an annoying phenomenon, can actually be exploited to detect certain dynamics. Paradoxically, to retrieve these dynamics, it can be useful to make voluntary coarse roundoff errors.
\bigskip

This phenomenon of stabilization of each dynamical invariant of $f$ by an infinite number of discretizations can be explained as follows. Each dynamical invariant of $f$ we will study can be approximated by a periodic orbit $\omega$ (or a finite union of such orbits). Generically, the points of $\omega$ are of Liouville type with respect to the sequence of grids, that is to say, there exists a subsequence of orders $N$ for which each point of $\omega$ is very close to a grid point\footnote{The ``closeness'' being adjusted to the grid mesh and to other quantities like the modulus of continuity of $f$.}, close enough so that for some orders $N$, if $x\in \omega$, then $f_N(x_N) = (f(x))_N$. In this case, there is a discrete orbit close to the actual orbit $\omega$. For a practical system (which therefore does not necessarily verify the hypotheses of the theorems, which assume that the systems are generic), it is hoped that a similar phenomenon occurs: if the orbit $\omega$ is not too long, then it is likely that considering a reasonable number of discretizations on different grids, we recover the orbit $\omega$ on at least one of these grids.

\section{Reading guide}

This manuscript is divided into four parts. Part~\ref{PartOne} concerns discretizations of generic homeomorphisms. Part~\ref{PartII} 
focuses on the discretizations of linear maps. In Part~\ref{PartTri}, we study the discretizations of generic $C^1$-diffeomorpihsms. Finally, Part~\ref{partconcl} analyses and compares the results obtained in the four first parts..

In Part~\ref{PartOne}, Chapter~\ref{ChapNota} consists in notations and preliminaries, which will be used throughout Part~\ref{PartOne}. The reader will also find an index of notations at the end of this memoir. Chapters~\ref{ChapDissip} and \ref{ChapCons} are widely independents, while Chapter~\ref{ChapRot} uses notions developed in Chapters~\ref{ChapDissip} and \ref{ChapCons}.

Part~\ref{PartII} is completely independent from Parts~\ref{PartOne} and \ref{PartTri}. More precisely, Chapter~\ref{ChapAlm} is somehow an introduction to Chapter~\ref{Souris}, where lies most of the original results of this part. Chapter~\ref{Stat0}, is very short, and uses the formalism set up in the two previous chapters.

Part~\ref{PartTri} uses the results of Part~\ref{PartII}. In more details, Chapter~\ref{ChapPerturbLem} is independent from the rest of the manuscript, but we will often compare the results with those obtained in Part~\ref{PartOne}.  Chapter~\ref{ChapDeg} uses crucially the analysis of the linear case made in Part~\ref{PartII}. Finally, in Chapter~\ref{chapPhys}, we combine techniques developped in Part~\ref{PartII} with statements of Chapter~\ref{ChapPerturbLem}.

Finally, Part~\ref{partconcl} is a kind of conclusion. In its first chapter (Chapter~\ref{Blablablabla}), we compare the obtained results from various viewpoints. In the second one (Chapter~\ref{ouaf}), I give examples of open questions that seemed relevant to me.

\chapter{Historic of the topic}\label{SecHisto}

\section{Small survey of classical results in generic dynamics}

Before looking at the dynamics of discretizations of generic maps, it is good to know the actual dynamics of these applications (that is to say, without discretization of the space). In addition, techniques developed for understanding the dynamics of generic homeomorphisms (or diffeomorphisms) will be very useful to study their discretizations.
\bigskip

One of the very firsts results in generic dynamics is Oxtoby-Ulam theorem (published in 1941, see \cite{Oxto-meas}), which states that a generic conservative homeomorphism of a compact manifold of dimension at least two is ergodic. Since then, many works have helped to understand the dynamics of a generic conservative homeomorphism. Note that this dynamics is neither completely trivial nor extremely chaotic: for example, a generic conservative homeomorphism is topologically mixing, but not strongly mixing (in the measurable sense). Note also that the techniques used to prove the results are varied and beautiful: in particular, there is a transfer theorem of ergodic generic properties from the space of automorphisms to that conservative homeomorphisms (due to S.~Alpern see \cite{MR550406}), or a 0-1 law for generic dynamical ergodic properties (due to E.~Glasner and J.~King, see \cite{Glas-zero}). Interested readers can consult the surveys \cite{Alpe-typi} and \cite{MR2931648}. Let us give a sketch of the proof of genericity of transitivity (probably one of the simplest proofs): considering a conservative homeomorphism $f$,
\begin{itemize}
\item we begin by breaking the dynamics $f$ in discretizing the phase space and applying Lax's theorem (see \cite{MR0272983} and theorem~\ref{Lax} page~\pageref{Lax}), which asserts that a conservative homeomorphism is arbitrarily close to a cyclic permutation of a grid of the manifold having some good properties;
\item we then rebuild a conservative homeomorphism from this cyclic permutation, by a $C^0$ closing lemma (Proposition~\ref{extension} page~\pageref{extension}).
\end{itemize}
This proof technique can be used to obtain results on discretizations of generic conservative homeomorphisms; this is what we do in Chapter~\ref{ChapCons}.
\bigskip

The case of generic \emph{dissipative} homeomorphisms has been studied lately. In fact, the systematic study lead by E. Akin, M. Hurley and J. Kennedy in the survey \cite{MR1980335} is a bit disappointing, as all the results go in the same direction : the dynamics of a generic dissipative homeomorphism ``contains'' simultaneously all the possible wild behaviours.

This survey only concerns \emph{topological} dynamics. Quite recently, F. Abdenur and M. Andersson studied \emph{ergodic} properties of generic dissipative homeomorphisms: they consider the behaviour of Birkhoff sums for almost every starting point (see \cite{MR3027586}). For this purpose, they establish a technical result called \emph{shredding lemma} that allows them to show that a generic homeomorpihsm is \emph{weird} (see Definition~\ref{strange} page~\pageref{strange}). We will make use of this technical lemma to establish dynamics of discretizations of a generic dissipative homeomorphism.
\bigskip

Properties of generic $C^1$-diffeomorphisms is a very active subject -- for both conservative and dissipative cases -- and so has became very dense. In the conservative case, the lack of $C^1$ analogue to the property of rebuilding of homeomorphisms forces to establish closing and connecting lemmas. Among others, the closing lemma of C.~Pugh \cite{MR0226669} implies that a generic conservative $C^1$-diffeomorphism has a dense set of periodic points. Also, the connecting lemma for pseudo-orbits of C.~Bonatti and S.~Crovisier \cite{MR2090361} implies that a generic conservative $C^1$-diffeomorphism is transitive. For its part, the question of genericity of ergodicity is still open in this context, despite the recent advance of A.~Avila, S.~Crovisier and A.~Wilkinson \cite{ArturSylvain}: either a generic conservative $C^1$-diffeomorphism has all its Lyapunov exponents null, or is non-uniformly Anosov and ergodic. There is a lot of other results about this subject, the reader could consult the survey \cite{MR2288283} of S.~Crovisier or the introduction of \cite{ArturSylvain}.

For generic dissipative $C^1$ diffeomorphisms, the researches are  by conjectures of J.~Palis which describe what could be the topological and ergodic dynamics of a generic element. Again, the subject is very rich and the reader can consult the Bourbaki seminar \cite{MR2074058} or the ICM report \cite{MR1957538} of C. Bonatti, the habilitation memoir \cite{MR3136194} or the recent survey \cite{SurvolSylvain} of S.~Crovisier.
\bigskip

Compared to the situation of $C^1$-diffeomorphisms, the study of the dynamical properties of generic $C^r$-diffeomorphisms, with $r>1$, is blocked by the absence of a closing lemma (about the closing lemma in $C^1$ topology, see the book \cite{MR1662930} of M.-C. Arnaud). To my knowledge, the only case in which we have results other anecdotal (besides the trivial case of dimension 1) is that of conservative diffeomorphisms of surfaces (see the article \cite{MR1971199} of J. Franks and P. Le Calvez). Note that in larger regularity ($r\ge 4$), KAM theorem implies that on any compact surface, there are open sets of conservative $C^r$-diffeomorphisms (in $C^r$ topology) on which the elements are non transitive (see for example Section 4 of \cite{MR1648127}).

\section[Survey about spatial discretizations of dynamical systems]{Global survey about spatial discretizations of dynamical systems}\label{SecHistoGene}

The first numerical simulations of dynamical systems appear in the 1960s, with for example the famous works of E.~Lorenz and M.~Hénon. Among others, M.~Hénon conducts an extensive digital study of what is now called ``conservative Hénon map'' (see \cite{MR0253513}):
\[f_\alpha(x,y) = \big( x\cos\alpha - (y-x^2)\sin\alpha\ ,\ x\sin\alpha + (y-x^2)\cos\alpha\big).\]
M. Hénon notes that the phenomena appearing on simulations (in particular, ``elliptical islands in a chaotic sea'') are consistent with what we already knew about these applications from a theoretical study. However, he does not actually seem to care about the possible bad effects produced by the digital truncation.

A few years later, P.~Lax remarks in \cite{MR0272983} that the discretization destroys the bijective behaviour of conservative Hénon maps. To overcome this problem, P. Lax shows that any conservative homeomorphism of the torus is arbitrarily well approximated by permutations of ``natural'' grids on the torus. In some sense, even if the discretization of a homeomorphism is not bijective, we know that there is at least one finite map close to this homeomorphism that is. This theorem is to my knowledge the first theoretical result concerning the dynamics of discretizations; its nice proof is essentially combinatorial and based on Hall's marriages lemma (see page~\pageref{mariage}). An improvement of this result is proved by S. Alpern in \cite{Alpe-appr} (this result is in fact due to J. Oxtoby and S. Ulam, see \cite{Oxto-meas}): in the statement of Lax theorem (Theorem~\ref{Lax} page~\pageref{Lax}), the term ``permutation'' can be replaced by ``cyclic permutation''. This result indicates that the dynamics of any conservative homeomorphism is close to a transitive finite map. In the 1990s, in \cite{MR1307740, MR1453713}, Lax theorem is generalized to the case of applications preserving Lebesgue measure (without assumption of bijectivity or continuity); the authors also give examples of actual obtaining of the permutation in dimension 1.
\bigskip

From the late $1970$s, physicists and mathematicians got interested in discretizations of particular dynamical systems with well known behaviour. Some articles perform numerical simulations, and assess the consistency between the results of these simulations and the actual dynamics of the system. Others study this consistency from a more theoretical point of view. For example, in 1978, G. Benettin \emph{et al.} got interested in obtaining numerically the physical measures of some Anosov diffeomorphisms \cite{MR0478237, MR534103}, including the Arnold cat map and some of its perturbations.

\label{refG}In 1983, in \cite{MR700317}, J.-M. Gambaudo and C. Tresser show on well chosen examples that the sinks of a homeomorphism can be undetectable in practice, simply because the size of their basins can be extremely small, even for homeomorphisms with quite simple definition.

Still in the 80s, one sees theoretical work on the dynamics of discretizations of classical dynamical systems appearing. For example, in \cite{MR938160}, the authors are interested in the shadowing property for some parameters of the Hénon map. They study the same property for logistic maps in \cite{MR907194} and for the tent and quadratic maps in \cite{MR929137}. Finally, in \cite{MR835874}, they study the characteristics of the orbits of discretizations of the tent map, namely their number, their lengths and their distribution.

Also, in \cite{MR1169615}, the author studies numerically the number of limit cycles and the length of the longest of these cycles for discretizations of some logistic maps.

In \cite{MR1353178,MR1392078,MR1400185,MR1445902,MR1481914}, the authors note that some quantities related to discretizations of dynamical systems of dimension 1 -- like the proportion of discretization points in the basin of attraction of the fixed point 0, the distribution of the length of the cycles, the stabilization time, etc. -- are similar to the same quantities for random applications with an attractive centre, especially when these quantities are averaged over several successive discretizations. Their study is based mainly on assumptions of convergence of the dynamics of discretizations to that of the initial application, verified experimentally on well-chosen examples.

Finally, in one of its latest papers \cite{MR1678095}, O.E. Lanford studies simulations of discretizations of systems of dimension 1. Among others, he focuses on periods and basins sizes of the periodic orbits of discretizations. He remarks that to retrieve the dynamics of the initial map, the best would be to adjust the time up to which the iterates of discretizations are computed to the mesh of these discretizations. This viewpoint, although very interesting, is very different from that adopted in this manuscript.

These results relate to very specific examples of dynamical systems; they would not be easily generalized to wider classes of dynamics.
\bigskip

The idea of O.E. Lanford consisting in adjusting the length of orbit segments had already been developed in the late $1980$s by A.~Boyarsky. In \cite{MR862028}, he explains heuristically why one usually finds absolutely continuous measures on simulations. His arguments are based on the tracking of long segments of orbits; the only obstacle for the obtaining of a rigorous proof is the lack of uniformity in Birkhoff 's ergodic theorem \footnote{In fact, the intuition of A.~Boyarsky works for a uniquely ergodic homeomorphism (Proposition \ref{UniqErgoMierno} of T.~Miernowski), but is false in more general situations (Theorem~\ref{TheoMesIntro} for generic conservative homeomorphisms, and Theorem~\ref{TheoMesIntroC1} for conservative $C^1$-generic diffeomorphisms.}.

In \cite{MR959419}, A. Boyarsky and P. Góra establish the following result, which also relates to the obtaining of absolutely continuous measures from discretizations.

\begin{theo}
If $f$ has a unique absolutely continuous invariant measure $\mu$, and if $\alpha>0$ is such that there exists a subsequence $f_N$ admitting a segment of orbit of length bigger than $\alpha q_n$ (where $q_n$ is the cardinality of $E_N$), if one set $\nu_N$ the uniform measure on this segment of orbit, then $\nu_N\rightharpoonup \mu$.
\end{theo}

The fact of having an absolutely continuous invariant measure is known for rather large classes of systems ($C^{1+\alpha}$ expanding maps of the circle, $C^{1+\alpha}$ Anosov diffeomorphsims, \dots). on the other hand, the existence of a segment of orbit of length proportional to that of the grid seems to be rarely verified (for example it is not true for a generic expanding map the circle, simply because the degree of recurrence is zero). Despite this, apart from Lax theorem, this is one of the first theoretical results about discretizations of dynamical systems.

Meanwhile, research is carried out by Soviet scientists, including M. Blank in \cite{MR765293,MR875433,MR1037009,MR1299502}, who focuses on the following question (among others): how to retrieve the periodic orbits of the dynamical system (both attractive and repulsive) from the discretizations? The author presents the phenomenon of period doubling: on discretizations, it may appear an orbit whose length is a nontrivial multiple of that of the actual orbit it should represent. Note that M. Blank is the author of a monograph \cite{MR1440853} whose fifth chapter is fully devoted to the problem of spatial discretization; in particular, there is a discussion on the behaviour of invariant measures and a study of the case of torus rotations.
\bigskip

\label{DiaDia}In the 90s, a group of researchers including P. Diamond, P. Kloeden, V. Kozyakin, J. Mustard and A. Pokrovskii published a series of articles about detection on discretizations of some dynamical properties. For example, in 1993, in \cite{MR1201881}, the authors note that any attractor is detected by the discretizations. The study of discretizations leads the authors to define some variants of stable dynamical properties adapted to the study of discretizations. For example, in \cite{MR1265228}, they define a notion of chain transitivity adapted to finite maps, the property of ``weak chain recurrence''. They then state a lot of results explaining how to manipulate this concept; in particular they check that it persists while passing to discretizations. In the same vein, in \cite{MR1331572, MR1354569}, the authors define what they call semi-hyperbolic maps (which include hyperbolic maps); they establish that semi-hyperbolic maps have the shadowing property and others properties such as structural stability of semi-hyperbolicity or semi-continuity of entropy. Finally, the authors define what is the \emph {minimal centre of attraction}, a set containing the non wandering set \cite{MR1387977}; by its very definition, this set is stabilized when passing to discretizations.

This work deserves credit for formally showing that some dynamical properties can be detected on spatial discretizations. However, one should note that they only concern properties that are by definition robust (existence of attractive orbits, weak chain recurrence, etc.), and therefore are naturally inherited by discretizations.
\bigskip

One of the simplest examples of chaotic dynamical systems are torus linear Anosov automorphisms. In 1992, F. Dyson and H. Falk in \cite{MR1176587} make a fairly complete theoretical study of the global period of discretizations of the Arnold's cat map $\begin{pmatrix} 1&1\\1&2\end{pmatrix}$. In particular, the dynamics of discretizations is completely different from that of the linear Anosov automorphism (see the beginning of this introduction, page~\pageref{PageMiaou}). This work was taken up and generalized to all linear Anosov automorphisms of $\T^2$ by \'E. Ghys in \cite{Ghys-vari}. In this paper, \'E. Ghys notes that the global period of discretizations is very low compared to that of a random map on a set with the same number of elements as the grid. He attributes this phenomenon to the strong arithmetic properties of linear automorphisms; this may suggest that it is in some sense exceptional. This leads him to propose to study the behaviour of discretizations of generic homeomorphisms.

This issue was almost completely solved by T. Miernowski in the case of generic homeomorphisms of the circle (see \cite{MR2279269} and \cite{Mier-dyna}). We will survey these interesting results in Section \ref{SurvolMierno}. In \cite{MR2279269}, T. Miernowski also establishes a theorem about convergence of the dynamics of discretizations to that of the initial map, from an ergodic point of view, under the rather restrictive assumption that the considered homeomorphism is uniquely ergodic:

\begin{proposition}(Miernowski)\label{UniqErgoMierno}
Let $f : M \to M$ be a homeomorphism having a unique invariant probability measure $\mu$. For every $N\in\N$, let $\gamma_N$ be a periodic cycle of the $N$-th discretization $f_N$ of $f$. Let $\nu_N$ be the uniform probability measure on the cycle $\gamma_N$. Then the measures $\nu_N$ converge weakly towards $\mu$, independently from the choice of each cycle $\gamma_N$.
\end{proposition}

This convergence of the discrete dynamics to the continuous one is again highlighted in the case of generic conservative torus homeomorphisms by Theorem 2.2.2 of \cite{Mier-dyna} (that we already met earlier this introduction, see Theorem~\ref{TheoMierno}): for a  generic conservative homeomorphism, there is a subsequence of discretizations that are permutations of the grids. The author nevertheless conjecture that this phenomenon is no longer typical in the case of generic hyperbolic $C^r$-diffeomorphisms of the torus, $r\ge 1$. According to him, the behaviour of the discretizations of such maps should approximate that of a random map of a set having the same number of elements as the discretization grid.
\bigskip

The question of the behaviour of discretizations of a generic system was also addressed by P.P.~Flockermann in his thesis \cite{Flocker} under the direction of O.E.~Lanford. P.P.~Flockermann is interested in $C^r$ expanding maps of the circle for $r\ge 2$. He basically shows that in finite time, the discretizations of such applications behave the same way as stochastic perturbations (see Section~\ref{Floque}). These results, although interesting and nontrivial, concern only the finite time behaviour of discretizations and say nothing about their dynamics.
\bigskip

There are many other points of view about modelling numerical simulations. For example, given a family of maps and a dynamical property, one can wonder if a machine can decide in finite time if an element of this family has this property. For example, in \cite{MR2078664}, A. Arbieto and C. Matheus show that the properties ``having positive Lyapunov exponents'' and ``having an SRB measure'' are undecidable among quadratic and Hénon families, and that the property ``having positive topological entropy'' is decidable among the quadratic family.

A study on the modelling of a discretization of a homomorphism by multivalued maps was made in \cite{MR2776399} (see also \cite{MR1403460} and \cite{MR1387977}). This article contains a  very interesting discussion on the transition from continuous to discrete: \emph{a priori}, one can recover the original system properties from a finite number of discretizations only when these properties are robust. This leads the authors to define a ``finite resolution to property'': in particular, these are dynamical properties that can be checked in finite time by computational methods.
\bigskip

For different viewpoints on the history of discretizations, see the short survey of J.~Buzzi \cite{SurvJerome}, or the more complete survey \cite{MR2863582} of S.~Galatolo, M.~Hoyrup and C.~Rojas.
\bigskip

In a completely different context, the discretizations of linear maps have been studied for applications to image processing. The goal is for example to answer the question: how to rotate a digital image without losing quality? In particular, we would like to estimate the loss of information caused by non bijectivity of discretizations of rotations, or try to avoid aliasing phenomena (the appearance of resonances between a regular pattern on the image and the grid made of the pixels). Existing work focuses primarily on the study of discretisations of linear maps in short time; in particular they examine the local behaviour of the iterations of $\Z^2$ by these discretizations. In addition, the linear maps studied are always supposed with rational coefficients. We can get an idea of the topic by consulting the theses \cite{nouvel:tel-00444088} and \cite{thibault:tel-00596947} (see also page~\pageref{BibliLin} of this manuscript for a more complete bibliography).

\section[Works of T.~Miernowski]{Works of T.~Miernowski about discretizations of circle homeomorphisms}\label{SurvolMierno}

In this sub-section, we give the main results concerning the dynamics of uniform discretizations of circle homeomorphisms, obtained by T.~Miernowski in \cite{MR2279269}. By ``uniform'', we mean that these discretizations are taken with respect to the grids $E_N = \{i/N\mid 0\le i \le N-1\}$ (where the circle $\Sp^1$ is identified with $[0,1[$).

The preservation of a cyclic order on the circle implies that if $f$ is an increasing homeomorphism, then all the periodic orbits of a discretization $f_N$ have the same length, denoted by $q_N$. Then, there exists an integer $p_N$ such that if $\widetilde f$ is a lift of $f$ to $\R$, and $\widetilde f_N$ is the discretization of $\widetilde f$ on the uniform grid $E_N$, then for every point $\widetilde x \in E_N$ projecting in a periodic orbit $x$ of $f_N$, we have $\widetilde f_N^{q_N} (\widetilde x) - \widetilde x = p_N$. Then, it can be easily shown that $p_N/q_N$ converges to the rotation number $\rho(f)$ of~$f$.
\bigskip

The simplest example of circle homeomorphism is that of a rotation (there are even many cases when a homeomorphism is conjugated to an irrational rotation, see for example \cite{MR538680}). T.~Miernowski notes that the discretization of a rotation is itself a rotation; in particular it is bijective. If the number of periodic cycles of the rotation of angle $\alpha$ is denoted by $r_N$, the following trivial equality holds $N = r_N q_N$: the dynamics of $f_N$ is completely determined by the number $r_N$. In the case where $\alpha\in\R\setminus\Q$, the behaviour of $r_N$ can be easily deduced from a known result counting rational points in a region of the plane (see \cite{MR0057917}, see also \cite{MR1400185}), we get:
\[\lim_{M\to +\infty} \frac{\card\{N\le M \mid r_N = k\}}{M} = \frac{6}{\pi^2 k^2}.\]
Proposition~3.2 of \cite{MR2279269} extends this result to the rational case. With these results, we can deduce that the length $q_N$ of periodic orbits converges in the sense of Cesàro

\begin{proposition}
If $\alpha\in\R\setminus\Q$, then
\[\frac{1}{N}\sum_{i=1}^{N} q_i \sim N\frac{3}{\pi^2}\zeta(3).\]

If $\alpha = \frac{p}{q}\in\Q$, with coprime $p$ and $q$ there are two cases (with $\varphi$ the Euler indicator function):
\[\frac{1}{N}\sum_{i=1}^{N} q_i \sim \left\{\begin{array}{l}
\displaystyle N\frac{q-1}{q^2}\sum_{c=1}^{\frac{q-1}{2}} \frac{1}{c^2}\sum_{r=1}^{\frac{q-1}{2c}}\frac{\varphi(r)}{r} \quad\text{if }q\text{ odd},\\
\displaystyle N\frac{q-1}{q^2}\sum_{c=1}^{\frac{q}{2}} \frac{1}{c^2}\sum_{r=1}^{\frac{q}{2c}}\frac{\varphi(r)}{r} + \frac{N}{2}\frac{(q-1)^2(q+1)}{q} \varphi(q/2) \quad\text{if }q\text{ even.}
\end{array}\right.\]
\end{proposition}

Therefore, in some sense, there is a convergence of the dynamics of discretizations for Cesàro averages. This is completely opposite to what occurs in the case of generic conservative homeomorphisms, where even if we try to average over successive orders of discretizations, combinatorial dynamics of discretizations continues to be highly erratic (see Theorem~\ref{MoyCesaroIntro} page~\pageref{MoyCesaroIntro}).
\bigskip

T.~ Miernowski then looks at the case of generic circle homeomorphisms. The results describing the behaviour of the length of periodic orbits of discretizations is the following.

\begin{proposition}[Miernowski]\label{mier1}
An increasing generic circle homeomorphism $f$ has at least one stable periodic point. In particular, its rotation number is rational, denoted by $\rho(f) = p/q$. Then, there exists $N_0\in\N$ such that for all $N\ge N_0$, we have $p_N = p$ and $q_N = q$.
\end{proposition}

This easy result is somewhat disappointing: the dynamical behaviour of discretizations is the same for all large enough $N$, and actually stems from the fact that the starting homeomorphism has attractors (it is Morse-Smale). That is why T.~Miernowski try to avoid this type of dynamics: he considers the more difficult case of generic homeomorphisms among those with rational rotation number but no attractive orbit, or those having an irrational rotation number.

In the first case, T.~Miernowski shows the following theorem.

\begin{theo}[Miernowski]\label{mier2}
Among semi-stables homeomorphisms\footnote{I.e., the lift of $f^q$ to $\R$ is either always bigger than (or equal to), or always smaller than (or equal to) the translation of vector $p$.} with rotation number $p/q$, there exists an open and dense set for which there exists a subsequence $(N_k)_{k\in\N}$ such that $q_{N_k}(f) = q$.
\end{theo}

The proof of this result uses crucially the natural order on the circle.

In the case of irrational rotation number, T.~Miernowski proposes two approaches. The first concerns the \emph{generic} diffeomorphisms.

\begin{theo}[Miernowski]\label{mier3}
Let $\Psi : \N \to \R_+$ be an increasing map tending to $+\infty$ in $+\infty$. Let $r\in\N\cup\{\infty,\omega\}$ (in particular, one can have $r=0$). Then, for an generic irrational $C^r$-diffeomorphism $f$,
\[\underset{N\to+\infty}{\underline\lim} \frac{q_N(f)}{\Psi(N)} = 0.\]
\end{theo}

In other words, the sequence $q_N(f)$ tends to $+\infty$, but has a subsequence which does arbitrarily slowly. The proof of this result is pretty nice because it is mainly based on the study conducted in the semi-stable case and the approximation of any irrational diffeomorphism by rational semi-stable ones.

In the second approach, T.~Miernowski is interested in the case of \emph{prevalent} diffeomorphisms. The result is very different from what happens for generic diffeomorphisms: the sequence $q_N(f)$ tends to $+\infty$ almost as fast as $\sqrt N$.

\begin{theo}[Miernowski]\label{mier4}
Let $3\le r \le \omega$, and $f$ be a prevalent irrational $C^r$-diffeomorphism. Then, for every $\varep>0$, there exists $K>0$ such that $q_N(f) \ge KN^{1/(2+\varep)}$.
\end{theo}

The proof of this theorem is mostly based on a classic theorem about prevalent circle dynamics, i.e. the differentiable conjugation to rotation theorem due to M.~Herman (see \cite{MR0458480}).

\section[Works of P.P.~Flockermann and O.E.~Lanford]{Works of P.P.~Flockermann and O.E.~Lanford about discretizations of circle expanding maps}\label{Floque}

In his thesis \cite{Flocker}, P.P.~Flockermann studies the behaviour of discretizations\footnote{Still for the grids $E_N = [0.1 [\cap \Z / N$.} of expanding circle map of class (at least) $C^2$. Besides fairly extensive numerical studies, P.P.~Flockermann shows two main results; they express that for a fixed time $k$, and under assumptions of linear independence of derivatives, the discretization operation of a typical circle expanding map behaves essentially as a random perturbation.
\bigskip

First of all, one can make statistics about numerical error done at each iteration: for $x\in E_N$, we denote $\varep_k(x,N) = N\big(f_N^k(x) - f(f_N^{k-1}(x))\big)$ the normalized error made at time $k$ (this error is between $-1/2$ and $1/2$). The theorem of P.P.~Flockermann is the following (Theorem~9 of \cite{Flocker}).

\begin{theo}[Flockermann]
Let $f$ be a degree 2 expanding map of $\Sp^1$ such that
\begin{itemize}
\item either $f$ is generic among $C^r$ expanding maps, $r\ge 2$;
\item or $f$ is analytic, but different from $x\mapsto 2x$.
\end{itemize}
Let $\vartheta : \N\to\R$ be a map such that $1/N \ll \vartheta(N) \ll 1$.
Then there exists a countable set $D\subset \Sp^1$ such that for all $x_0\in\Sp^1\setminus D$ and all $k_0\in\N$, the distributions of the $k_0$-tuples $\big(\varep_k(x,N)\big)_{k\le k_0}\in[-1/2,1/2]^{k_0}$, for $x\in E_N \cap [x_0-\vartheta(N),x_0+\vartheta(N)]$, converge to the uniform distribution on $[-1/2,1/2]^{k_0}$.
\end{theo}

In other words, at an intermediate (or mesoscopic) space scale $\vartheta(N)$, around a generic point $x_0$, the roundoff errors made at each iteration of the discretization are time-independent and follow the uniform law on $[-1/2,1/2]$. The proof of this result is divided into two independent parts. The first is a result leading to the conclusions of the theorem under the hypothesis of linear independence of the successive derivatives. The second shows that these assumptions are verified by generic $C^r$-expanding maps or analytic expanding maps different from $x\mapsto 2x$.

Using the formalism of model sets, we will give a generalization of this theorem in arbitrary dimension (Proposition~\ref{EquidistribErr} page \pageref{EquidistribErr}).
\bigskip

Secondly, P.P.~Flockermann studies the distribution of $\card(f_N^{-k}(x))$ for fixed $k$ and a lot of points $x$ in a small region of $\Sp^1$. For simplicity, we will enunciate the theorem in the case where time $k$ is 1; the reader will easily infer a general statement (this statement is still quite tricky to formulate, see Theorem 12 and Corollary 3 of \cite{Flocker}).

\begin{theo}[Flockermann]
Let $r\ge 2$, and $f$ be a generic $C^r$-expanding circle map of degree 2. Let $\vartheta : \N\to\R$ be a map such that $1/N \ll \vartheta(N) \ll 1$. Then there exists a countable set $D\subset \Sp^1$ such that for every $x_0\in\Sp^1\setminus D$, if we denote $\{y_1,y_2\} = f^{-1}(x_0)$, when $N$ goes to $+\infty$,
\begin{itemize}
\item the proportion of points $x\in E_N \cap [x_0-\vartheta(N),x_0+\vartheta(N)]$ which do not have preimage under $f_N$ tends to
\[\left(1-\frac{1}{f'(y_1)}\right)\left(1-\frac{1}{f'(y_2)}\right);\]
\item the proportion of points $x\in E_N \cap [x_0-\vartheta(N),x_0+\vartheta(N)]$ which have exactly one preimage under $f_N$ tends to
\[\left(1-\frac{1}{f'(y_1)}\right)\frac{1}{f'(y_2)} + \frac{1}{f'(y_1)}\left(1-\frac{1}{f'(y_2)}\right);\]
\item the proportion of points $x\in E_N \cap [x_0-\vartheta(N),x_0+\vartheta(N)]$ which have two preimages under $f_N$ tends to
\[\frac{1}{f'(y_1)}\frac{1}{f'(y_2)}.\]
\end{itemize}
\end{theo}

As the previous one, this theorem states that in short time, the operation of discretization behaves as a random perturbation: given a dilatation $h_\lambda$ of $\R$ of ratio $\lambda$, then the proportion of points in the image of the discretization of $h_\lambda$ is $1/\lambda$. But the graph of $f$, in the neighbourhood of $f^{-1}(x_0)$, has exactly two branches; the first one is well approximated by the dilatation $h_{f'(y_1)}$ and the second by the dilatation $h_{f'(y_2)} $. If the family $(1,f'(y_1),f'(y_2))$ is $\Q$-free, then the events ``have an inverse image under $f_N$ in the first branch'' and ``have an inverse image by $f_N$ in the second branch'' become independent, which provides the theorem. This implies the following dynamic result.

\begin{theo}[Flockermann, Lanford, unpublished]\label{Ceasar}
Let $r\ge 2$, and $f$ be a generic $C^r$-expanding circle map of degree 2. Then the degree of recurrence of $f$ is 0.
\end{theo}

We will prove this result in Section~\ref{SecRateExpand} (Corollary~\ref{Tau0Dilat}). This section also contains a statement linking the local and global behaviours of discretizations (Theorem~\ref{TauxExpand}).

\part{Discretizations of generic homeomorphisms}\label{PartOne}

\parttoc

\chapter*[Introduction]{Introduction}

In this part, we consider the dynamics of discretizations of generic homeomorphisms, \emph{i.e.} we tackle the following question:

\begin{ques}
Which dynamical properties of a generic homeomorphism $f$ can be read on the dynamics of its discretizations $(f_N)_{N\geq 0}$?
\end{ques}

We will establish properties for both \emph{dissipative}, \emph{i.e.} arbitrary homeomorphisms of $X$, and \emph{conservative} homeomorphisms, \emph{i.e.} homeomorphisms of $X$ that preserve a given good probability measure. In this part, our results concern generic homeomorphisms of a compact manifold (with boundary) of dimension $n\geq 2 $. 

\bigskip

We will prove many results, concerning various aspects of the dynamics of the discretizations, adapting some classical tools of study of the generic dynamics of homeomorphisms. Morally, our results express that in the dissipative generic case, the dynamics of the discretizations tends to the ``physical'' dynamics of the initial homeomorphism\footnote{That is, the dynamics that occur for almost every point with respect to Lebesgue measure.} whereas in the conservative generic setting, the dynamics of the discretizations accumulates on all the possible dynamics of the initial homeomorphism, and moreover the physical dynamics cannot be detected on discretizations. In the rest of this introduction, we try to organize our results according to some ``lessons'':

\begin{constat}
The dynamics of discretizations of a generic dissipative homeomorphism tends to the ``physical dynamics'' of the initial homeomorphism.
\end{constat}

We first study properties of discretizations of generic \emph{dissipative} homeomorphisms\footnote{Without assumption of preservation of a given measure.}. The ergodic behaviour of such a generic homeomorphism is deduced from the \emph{shredding lemma} of F. Abdenur and M. Andersson \cite{MR3027586} (Lemma \ref{déchet}), which implies that a generic homeomorphism has a ``attractor dynamics'' (see Corollaries \ref{convattra} and \ref{MesDissip}):

\begin{theorem}
For a generic homeomorphism $f$, the closure $\overline{A_0}$ of the set of Lyapunov-stable periodic points is a Cantor set of dimension 0 which attracts almost every point of $X$. Moreover, the measure $\mu^f_X$ (see Definition \ref{defmes} page \pageref{defmes}) is well defined, atomless and is supported by $\overline{A_0}$.
\end{theorem}

This behaviour easily transmits to discretizations, for example every attractor of the homeomorphism can be seen on all the fine enough discretizations (Proposition \ref{Hausmodif}).

\begin{theorem}
For a generic homeomorphism $f$, the recurrent set of its discretization $\Omega(f_N)$ tends to $\overline{A_0}$ in the following weak sense: for all $\varep>0$, there exists $N_0\in \N$ such that for all $N\ge N_0$, there exists a subset $\widetilde E_N$ of $E_N$, stabilized by $f_N$, such that, noting $\widetilde \Omega(f_N)$ the corresponding recurrent set, we have $\frac{\card(\widetilde E_N)}{\card(E_N)}>1-\varep$ and $d_H(A_0,\widetilde \Omega(f_N))<\varep$.
\end{theorem}

Moreover, the Cesàro limit of the pushforwards of the uniform measures on the grid by the discretizations tend to the Cesàro limit of the pushforwards of $\lambda$ by the homeomorphism. Indeed, we will prove the following (Theorem \ref{convmesdissip}).

\begin{theorem}
For a generic homeomorphism $f$, the measures $\frac 1m \sum_{i=0}^{m-1}{f}_*^i \lambda$ converge to a measure that we denote by $\mu^f_X$ (see Definition \ref{defmes}).

Concerning the discretizations, for a generic homeomorphism $f$, the measures $\mu^{f_N}_X$ tend to the measure $\mu^f_X$ when $N$ goes to infinity, where $\mu^{f_N}_X$ is the measure on the periodic orbits of $f_N$ such that the global measure of each periodic orbit is proportional to the size of its basin of attraction (see Definition \ref{defmes}).
\end{theorem}

Moreover, there is \emph{shadowing} of the dynamics of $f$ by that of its discretizations $f_N$ (Corollary \ref{stabilis}).

\begin{theorem}
For a generic homeomorphism $f$, for all $\varep>0$ and all $\delta>0$, there is a full measure dense open subset $O$ of $X$ such that for all $x\in O$, all $\delta>0$ and all $N$ large enough, the orbit of $x_N$ under $f_N$ $\delta$-shadows\footnote{That is, for all $k\in\N$, $d(f_N^k(x_N),f^k(x))<\delta$.} the orbit of $x$ under $f$.
\end{theorem}

Thus, it is possible to detect on discretizations the ``physical'' dynamics of a generic dissipative homeomorphism, that is the dynamics that can be seen by almost every point of $X$. This dynamics is mainly characterized by the position of the attractors and of the corresponding basins of attraction.

\begin{constat}
The dynamics of a single discretization of a generic conservative homeomorphism cannot be inferred from the dynamics of the initial homeomorphism.
\end{constat}

We then turn to the study of the conservative case. The starting point of our study is a question from \'E. Ghys (see \cite[Section 6]{Ghys-vari}): for a generic conservative homeomorphism of the torus, what is the asymptotical behaviour of the sequence of degrees of recurrence of $f_N$? A partial answer to this question was obtained by T. Miernowski in the second chapter of his thesis (see Corollary \ref{typlax}).

\begin{theorem}[Miernowski]
For a generic conservative homeomorphism\footnote{That is, there is a $G_\delta $ dense subset of the set of conservative homeomorphisms of the torus on which the conclusion of the theorem holds.} $f$, there are infinitely many integers $N$ such that the discretization $f_N$ is a cyclic permutation.
\end{theorem}

To prove this theorem, T. Miernowski combines a genericity argument with a quite classical technique in generic dynamics: Lax's theorem \cite{MR0272983} (Theorem \ref{Lax}), which states that any conservative homeomorphism can be approximated by cyclic permutations of the discretization grids. In fact this proof can be generalized to obtain many results about discretizations. We will establish some variants of Lax theorem; each of them, combined with a genericity argument, leads to a result for discretizations of generic homeomorphisms. For instance, we will prove the following (Corollary \ref{corovar2}).

\begin{theorem}
For a generic conservative homeomorphism $f$, there exists infinitely many integers $N$ such that the cardinality of $\Omega(f_N)$ is equal to the smallest period of periodic points of $f$.
\end{theorem}

Note that the combination of these two theorems answer \'E. Ghys's question: for a generic homeomorphism $f$, the sequence of the degrees of recurrence of $f_N$ accumulates on both $0$ and $1$; we can even show that it accumulates on the whole segment $[0,1]$ (Corollary~\ref{ConjEt}).

Another variation of Lax's theorem leads to a theorem that throws light on the behaviour of the discretizations on their recurrent set (Corollary \ref{corovar3}).

\begin{theorem}
For a generic conservative homeomorphism $f$ and for all $M\in\N$, there are infinitely many integers $N$ such that $f_N$ is a permutation of $E_N$ having at least $M$ periodic orbits.
\end{theorem}

To summarize, generically, infinitely many discretizations are cyclic permutations, but also infinitely many discretizations are highly non-injective or else permutations with many cycles. In particular, it implies that for all $x\in X$, there exists infinitely many integers $N$ such that the orbit of $x_N$ under $f_N$ does not shadow the orbit of $x$ under $f$: in this sense, generically, the dynamics of discretizations does not reflect that of the homeomorphism. Note that this behaviour is in the opposite of the dissipative case, here the individual behaviour of discretizations does not indicate anything about the actual dynamics of the homeomorphism.

\begin{constat}
A dynamical property of a generic conservative homeomorphism cannot be deduced from the frequency it appears on discretizations either.
\end{constat}

The previous theorems express that the dynamics of a single discretization does not reflect the actual dynamics of the homeomorphism. However, we might reasonably expect that the properties of the homeomorphism are transmitted to many discretizations. More precisely, we may hope that given a property $(P)$ about discretizations, if there are many integers $N$ such that the discretization $f_N$ satisfies $(P)$, then the homeomorphism satisfies a similar property. It is not so, for instance, we will prove the following result (Theorem~\ref{propdemin}).

\begin{theorem}
For a generic conservative homeomorphism $f$, when $M$ goes to infinity, the proportion of integers $N$ between $1$ and $M$ such that $f_N$ is a cyclic permutation accumulates on both $0$ and $1$.
\end{theorem}

In fact, for all the properties considered in the previous paragraph, the frequency with which they appear on discretizations of orders smaller than $M$ accumulates on both $0$ and $1$ when $M$ goes to infinity. Remark that these result imply that, even by looking at the frequency at which some properties occur, the discretizations of a generic conservative homeomorphism do not behave like typical random maps, as for a random map of a set with $q$ elements, the average number of periodic orbits is asymptotically $\log q$ (see for example \cite[XIV.5]{Boll-rand}).

\begin{constat}
Many dynamical properties of a generic conservative homeomorphism can be detected by looking at the dynamics of \emph{all} the discretizations.
\end{constat}

We have observed that we cannot detect the dynamics of a generic homeomorphism when looking at the dynamics of its discretizations, or even at the frequency with which some dynamics appears on discretizations. Nevertheless, the dynamical properties of a generic conservative homeomorphism can be deduced from the analogous dynamical properties of \emph{all} the discretizations. More precisely, we have a shadowing property of the dynamical properties of the homeomorphism: for each dynamical property of the homeomorphism, its discrete analogue can be seen on an infinite number of discretizations. It is worthwhile to note the intriguing fact that this shadowing property occurs for all the dynamical properties of a generic conservative homeomorphism, independently of the measure $\lambda$, while for a generic dissipative homeomorphism the dynamics of the discretizations converges to the ``physical'' dynamics of the homeomorphism (that is, the dynamics depending of $\lambda$).

This idea of convergence of the dynamics when looking at arbitrary large precisions can be related to the work of P. Diamond \emph{et al} (see page \pageref{DiaDia}). For instance, we will prove that the periodic orbits of a homeomorphism can be detected by looking at the periodic orbits of its discretizations (Theorem~\ref{corovar2}).

\begin{theorem}
For a generic homeomorphism $f$, for every $\varep>0$ and every periodic orbit of $f$, this periodic orbit is $\varep$-shadowed by an infinite number of periodic orbits of the same period of the discretizations
\end{theorem}

We will also prove a theorem in the same vein for invariant measures (respectively invariant compact sets), which expresses that the set of invariant measures (respectively compact sets) of the homeomorphism can be deduced from the sets of invariant measures (respectively invariant sets) of its discretizations. More precisely, we will prove the following result (Theorems \ref{EnsMesInv} and \ref{mesinv}, see Theorems \ref{CompactInv} and \ref{EquivCompact} for the compact versions).

\begin{theorem}
For a generic conservative homeomorphism $f$ and for every convex closed set (for Hausdorff topology) $\mathcal M$ of $f$-invariant Borel probability measures there exists an increasing sequence of integers $N_k$ such that the set of $f_{N_k}$-invariant probability measures tends to $\mathcal M$ (for Hausdorff topology).

Moreover, if $\mathcal M$ is reduced to a single measure, then $f_{N_k}$ can be supposed to bear a unique invariant measure.
\end{theorem}

In the third chapter of this part, we will present an application of the notion of discretization to the practical problem of computing numerically the rotation set of a torus homeomorphism. In particular, we will prove a theorem which expresses the shadowing property of the rotation set of a generic conservative torus homeomorphism by the rotation sets of its discretizations (Corollary \ref{CoroRotDiscrCons}).

\begin{theorem}
If $f$ is a generic conservative homeomorphism, then there exists a subsequence $f_{N_i}$ of discretizations such that $\rho_{N_i}(F)$ tends to $\rho(F)$ for Hausdorff topology (where $F$ is a lift of $f$ to $\R^2$).
\end{theorem}

This will give us a convenient way to compute the rotation set in practice, as we will prove that if we compute the rotation set corresponding to a starting point $x\in \T^2$ without roundoff error, then for almost every $x\in\T^2$, the obtained rotation set is reduced to a single vector, that is the mean rotation vector of $f$ (Proposition \ref{ExGeneCons}).

\begin{theorem}\label{I2}
Let $f$ be a generic conservative homeomorphism of the torus $\T^2$ and $F$ a lift of $f$ to $\R^2$. Then for almost every $\tilde x\in\R^2$, the corresponding rotation set
\[\overline \rho(x) = \bigcap_{M\in\N}\overline{\bigcup_{m\ge M} \left\{\frac{F^m(\tilde x)-\tilde x}{m}\right\}}\]
is reduced to the mean rotation vector with respect to the measure $\lambda$.
\end{theorem}

This shows that if we try to compute the rotation set by calculating segments of orbits \emph{without} making any roundoff error, we will only find the mean rotation vector of the homeomorphism. In Chapter \ref{ChapRot}, we will introduce the notion of \emph{observable rotation set}, which expresses which rotation vectors can be found by looking at almost every periodic orbit. We will compute this set for some examples, in particular for both generic conservative (that is, Theorem \ref{I2}) and dissipative (Proposition \ref{GeneDissip}) homeomorphisms.

\begin{constat}
The ``physical dynamics'' of a generic conservative homeomorphism plays no particular role for discretizations.
\end{constat}

The heuristic idea underlying the concept of physical measure is that these measures are the invariant ones which can be detected ``experimentally'' (since many initial conditions lead to these measures). Indeed, some experimental results on specific examples of dynamical systems show that they are actually the measures that are detected in practice for these examples (see for example \cite{MR0478237,MR534103} or \cite{MR862028, MR959419}). Moreover, if the dynamical system is uniquely ergodic, then the invariant measure appears naturally on discretizations (see \cite[Proposition 8.1]{MR2279269} and Proposition \ref{Miern}).

According to this heuristic and these results, we could expect from physical measures to be the only invariant measures that can be detected on discretizations of generic conservative homeomorphisms. This is not the case: for a generic conservative homeomorphism, there exists a unique physical measure, namely $\lambda$ (that follows directly from the celebrated Oxtoby-Ulam's theorem \cite{Oxto-meas}). According to the previous theorem, invariant measures of the discretizations accumulate on all the invariant measures of the homeomorphism and not only on Lebesgue measure.

However, we could still hope to distinguish the physical measure from other invariant measures. For this purpose, we define the canonical physical measure $\mu^{f_N}_X$ associated to a discretization $f_N$: it is the limit in the sense of Cesàro of the images of the uniform measure on $E_N$ by the iterates of $f_N$: if $\lambda_N$ is the uniform measure on $E_N$, then
\[\mu^{f_N}_X = \lim_{m\to\infty}\frac{1}{m}\sum_{i=0}^{m-1}(f_N^i)_*\lambda_N.\]
This measure is supported by the recurrent set $\Omega(f_N)$; it is uniform on every periodic orbit and the total weight of a periodic orbit is proportional to the size of its basin of attraction. The following theorem expresses that these measures accumulate on the whole set of $f$-invariant measures: physical measures cannot be distinguished from other invariant measures on discretizations, at least for generic homeomorphisms (see Theorems \ref{mesinv} and \ref{MesPhys2}).

\begin{theorem}\label{Mattt}
For a generic conservative homeomorphism $f$, the set of limit points of the sequence $(\mu^{f_N}_X)_{N\in\N}$ is the set of all $f$-invariant measures. Also, for every $f$-invariant measure $\mu$, there exists a subsequence $f_{N_k}$ of discretizations such that for every $x\in X$, the sequence of measures\footnote{Recall that by Definition \ref{defmes}, $\mu^{f_N}_x$ is the Cesàro limit of the pushforwards of the Dirac measure $\delta_{x_N}$ by the discretization $f_N$.} $\mu^{f_{N_k}}_x$ tends to $\mu$.
\end{theorem}

The same phenomenon appears for compact invariant subsets: the the recurrent subsets of the discretizations $f_N$ accumulate on the whole set of invariant compact subsets for $f$ (Proposition \ref{EquivCompact}).

\begin{constat}
On the numerical experiments we performed, the dynamics of a dissipative homeomorphism can be detected on discretizations, and a lot of different dynamical behaviours can be observed on discretizations of a conservative homeomorphism.
\end{constat}

We will compare our theoretical results with the reality of numerical simulations. Indeed, it is not clear that the behaviour predicted by our results can be observed on computable discretizations of a homeomorphism defined by a simple formula. On the one hand, all our results are valid ``for \emph{generic} homeomorphisms''; nothing indicates that these results apply to actual examples of homeomorphisms defined by simple formulas. On the other hand, results such as ``there are infinitely many integers $N$ such that the discretization of order $N$\dots'' provide no control over the integers $N$ involved; they may be so large that the associated discretizations are not computable in practice.

We first carried out simulations of dissipative homeomorphisms. The results of discretizations of a small perturbation of identity (in $C^0$ topology) may seem disappointing at first sight: the trapping regions of the initial homeomorphisms cannot be detected, and there is little difference with the conservative case. This behaviour is similar to that highlighted by J.-M. Gambaudo and C. Tresser in \cite{MR700317} (see page \ref{refG}). That is why it seemed to us useful to test a homeomorphism which is $C^0$ close to the identity, but whose basins are large enough. In this case the simulations point to a behaviour that is very similar to that described by theoretical results, namely that the dynamics converges to the dynamics of the initial homeomorphism. In fact, we have actually observed behaviours as described by theorems only for examples of homeomorphisms with a very small number of attractors.

For conservative homeomorphisms, our numerical simulations produce mixed results. From a quantitative viewpoint, the behaviour predicted by our theoretical result cannot be observed on our numerical simulations. For example, we do not observe any discretization whose degree of recurrence is equal to $1$ (i.e. which is a permutation). This is nothing but surprising: the events pointed out by the theorems are \emph{a priori} very rare. For instance, there is a very little proportion of bijective maps among maps from a given finite set into itself. From a more qualitative viewpoint, the behaviour of the simulations is quite in accordance with the predictions of the theoretical results. For example, for a given conservative homeomorphism, the degree of recurrence of a discretization depends a lot on the size of the grid used for the discretization. Similarly, the canonical invariant measure associated with a discretization of a homeomorphism $f$ does depend a lot on the size of the grid used for the 
discretization.

\begin{constat}
Discretizations can actually be very useful and efficient to compute some dynamical invariants like the rotation set of a torus homeomorphism.
\end{constat}

We have also performed numerical simulations of rotation sets. To obtain numerically an approximation of the observable rotation set, we have calculated rotation vectors of long segments of orbits for a lot of starting points, these points being chosen randomly fore some simulations and being all the points of a grid on the torus for other simulations. For the numerical approximation of the asymptotic discretized rotation set we chosed a fine enough grid on the torus and have calculated the rotation vectors of periodic orbits of the discretization of the homeomorphism on this grid.

We made these simulations on an example where the rotation set is known to be the square $[0,1]^2$. It makes us sure of the shape of the rotation set we should obtain numerically, however it limits a bit the ``genericity'' of the examples we can produce.

In the dissipative case we made attractive the periodic points which realize the vertex of the rotation set $[0,1]^2$. It is obvious that these rotation vectors, which are realized by attractive periodic points with basin of attraction of reasonable size, will be detected by the simulations of both observable and asymptotic discretized rotation sets; that is we observe in practice: we can recover quickly the rotation set in both cases.

In the conservative setting we observe the surprising behaviour predicted by the theory: when we compute the rotation vectors of long segments of orbits we obtain mainly rotation vectors which are quite close to the mean rotation vector, in particular we do not recover the initial rotation set. More precisely, when we perform simulations with less than one hour of calculation we only obtain rotation vectors close to the mean rotation vector, and when we let three hours to the computer we only recover one vertex of the rotation set $[0,1]^2$. On the other hand, the rotation set is detected very quickly by the convex hulls of discretized rotation sets (less than one second of calculation). Moreover, when we calculate the union of the discretized rotation sets over several grids to obtain a simulation of the asymptotic discretized rotation set, we obtain a set which is quite close to $[0,1]^2$ for Hausdorff distance. As for theoretical results, this suggests the following lesson:

\emph{When we compute segments of orbits with very good precision it is very difficult to recover the rotation set. However, when we decrease the number of digits used in computations we can obtain quickly a very good approximation of the rotation set. In fact, we have to adapt the precision of the calculation to the number of orbits we can obtain numerically.}

This phenomenon can be explained by the fact that each grid of the torus is stabilized by the corresponding discretization of the homeomorphism. Thus, there exists an infinite number of grids such that every periodic point of the homeomorphism is shadowed by some periodic orbits of the discretizations on these grids.
\bigskip



\chapter{Grids, discretizations, measures}\label{ChapNota}


\section{The manifold $X$ and the measure $\lambda$}\label{DefMani}

The results stated in the introduction for the torus $\T^n$ and the Lebesgue measure $\Leb$ actually extend to any smooth connected manifold $X$\index{$X$} with dimension $n\ge 2$, compact and possibly with boundary, endowed with a Riemannian metric~$d$. \emph{We fix once and for all such a manifold $X$ endowed with the metric $d$.} In the general case, Lebesgue measure on $\T^n$ can be replaced by a \emph{good measure} $\lambda$ on $X$:

\begin{definition}\label{bonne mesure}
A Borel probability measure $\lambda$\index{$\lambda$} on $X$ is called a \emph{good measure}, or an \emph{Oxtoby-Ulam} measure, if it is non-atomic, it has total support (it is positive on every non-empty open sets) and it is zero on the boundary of $X$.
\end{definition}

These restrictions are supported by Oxtoby-Ulam theorem (Theorem \ref{Brown-mesure}). \emph{We fix once and for all a good measure $\lambda$ on $X$.}

\begin{notation}
We denote by $\Hom(X)$\index{$\Hom(X)$} the set of homeomorphisms of $X$, endowed by the metric $d$ defined by:
\[d(f,g) = \sup_{x\in X} d\big(f(x),g(x)\big).\]
We denote by $\Hom(X,\lambda)$\index{$\Hom(X,\lambda)$} the subset of $\Hom(X)$ made of the homeomorphisms that preserve the measure $\lambda$ (\emph{i.e.} for every Borel set $A$, $\lambda(f^{-1}(A)) = \lambda(A)$), endowed with the same metric $d$. Elements of $\Hom(X)$ will be called \emph{dissipative} homeomorphisms, and elements of $\Hom(X,\lambda)$ \emph{conservative} homeomorphisms. 
\end{notation}

\section{Generic properties}

All the functional spaces we will use are Baire spaces, \emph{i.e.} the intersection of every countable collection of dense open sets is dense. In particular, the topological spaces $\Hom(X)$ and $\Hom(X,\lambda)$ are Baire spaces. Indeed, we easily check that the map
\[\delta(f,g) = d(f,g) + d(f^{-1},g^{-1})\]
defines a distance on $\Hom(X)$ and $\Hom(X,\lambda)$ which defines the same topology as $d$, and is complete in $\Hom(X)$ and $\Hom(X,\lambda)$ (see for example \cite[Appendice A.1]{MR2931648}).

We will sometimes abusively use the phrase ``for a generic element of $\mathcal E$, we have the property $(P)$''. By that we will mean that ``the property $(P)$ is generic in $\mathcal E$'', \emph{i.e.} ``there exists $G_\delta$ dense subset $G$ of $\mathcal E$, such that every $f\in G$ satisfy the property $(P)$''.

\section{$C^0$ extension of finite maps}

The basic tool we will use to perturb a homeomorphism in the $C^0$ topology is the finite map extension proposition. It asserts that if we have given a (conservative) homeomorphism and a pseudo-orbit of this homeomorphism, we can perturb this homeomorphism such that this pseudo-orbit becomes a real orbit.

\begin{prop}[Finite map extension]\label{extension}
Let $E\subset X\setminus \partial X$ be a finite set and $\sigma : E \to X\setminus \partial X$ be an injective map. Then there exists an homeomorphism $g\in \Hom(X,\lambda)$ such that $g_{|E} = \sigma$. Moreover, if $d(\sigma,\operatorname{Id}_E)<\delta$, then the homeomorphism $g$ can be chosen such that $d(g,\operatorname{Id}_X)<\delta$.
\end{prop}

The last condition ensures that if the finite map $\sigma$ is close to a given homeomorphism $g$ -- in other words $\sigma$ is a good approximation of $g$ --, then the extension of this finite map can be chosen close to the map $g$. 

The idea of the proof is quite simple: it suffices to compose correctly homeomorphisms whose support's size is smaller than $\varep$ and which are central symmetries within this support. We denote $E=\{x_1,\cdots,x_q\}$ and $\sigma(x_i) = y_i$. For each $i$, we construct a sequence $(z_j)_{1 \le j \le k}$ such that $z_1 = x_i$, $y_i = z_k$ and $d(z_j, z_{j+1}) <\varep/10$. Composing $k-1$ homeomorphisms as above, such that each one sends $z_j$ on $z_{j+1}$, we obtain a conservative homeomorphism which sends $x_i$ on $y_i$. Implementing these remarks is then essentially technical and boring. A detailed proof can be found in \cite{Oxto-meas} or in \cite[Section 2.2]{MR2931648}. This proposition is valid in the neighbourhood of the identity; to have a result concerning the neighbourhood of any homeomorphism $f$ it suffices to compose by the inverse of $f$. The following corollary also includes the case where the map $\sigma$ is not injective (as a non injective map is arbitrarily close to an injective one).

\begin{coro}\label{CoroExtension}
Let $E\subset X\setminus \partial X$ be a finite set, $\sigma : E \to X\setminus \partial X$, $f\in\Hom(X)$ (respectively $f\in\Hom(X,\lambda)$) and $\delta>0$. If $d(f,\sigma)<\varep$, then there exists an homeomorphism $g\in \Hom(X)$ (respectively $\Hom(X,\lambda)$) such that $d(g_{|E}, \sigma)<\delta$ and $d(f,g)<\varep$. Moreover, if $\sigma$ is injective, we can suppose that $g_{|E} = \sigma$.
\end{coro}

\section{Discretization grids, discretizations of a homeomorphism}

We now define a more general notion of discretization grid than in the introduction. Depending on the case, some additional assumptions about these grids will be needed (see also the next section).

\begin{definition}[Discretization grids]\label{grillmiam}
A \emph{sequence of discretization grids} on $X$ is a sequence $(E_N)_{N\in\N}$\index{$E_N$} of discrete subsets of $X\setminus\partial X$, called \emph{grids}, such that the mesh of these grids tends to 0: for all $\varep>0$, there exists $N_0\in\N$ such that for all $N\ge N_0$, the grid $E_N$ is $\varep$-dense\footnote{That is, for every point $x\in X$ there exists a point $y\in E_N$ such that $d(x,y)<\varep$.}. We denote by $q_N$\index{$q_N$} the cardinality of $E_N$.
\end{definition}

\emph{We fix once and for all a sequence $(E_N)_{N\in \N}$ of discretization grids on $X$.}

To each grid is associated a discretization map in the following way.

\begin{definition}[Discretizations]\label{notA}
Let $P_N$\index{$P_N$} be a projection of $X$ on $E_N$ (the projection of $x_0\in X$ on $E_N$ is some $y_0\in E_N$ minimizing the distance $d(x_0,y)$ when $y$ runs through $E_N$). Such a projection is uniquely defined out of the set $E'_N$\index{$E'_N$} consisting of the points $x\in X$ for which there exists at least two points minimizing the distance between $x$ and $E_N$. On $E'_N$ the map $P_N$ is chosen arbitrarily (nevertheless measurably). For $x\in X$ we denote by $x_N$ the \emph{discretization of order $N$} of $x$, defined by $x_N = P_N(x)$. For $f\in\Hom(X)$ we denote by $f_N : E_N\to E_N$ the \emph{discretization of order $N$} of $f$, defined by $f_N = P_N\circ f$. We denote by $\mathcal{D}_N$\index{$\mathcal{D}_N$} be the set of homeomorphisms $f$ such that $f(E_N)\cap E'_N=\emptyset$.

If $\sigma : E_N\to E_N$ and $f\in \Hom(X)$, we denote by $d_N(f,\sigma)$\index{$d_N$} the distance between $f_{|E_N}$ and $\sigma$, considered as maps from $E_N$ into $X$.
\end{definition}

\begin{rem}
We might wonder why the points of the discretization grids are supposed to be \emph{inside} $X$. The reason is simple: a homeomorphism $f$ sends $\partial X$ on $\partial X$. Putting points of some grids on the edge could perturb the dynamics of discretizations $f_N$. In particular, it would introduce at least one orbit with length smaller than $\card(E_N \cap\partial X) $.
\end{rem}

\begin{rem}\label{MerciReferee}
As the exponential map is a local diffeomorphism, the sets $E'_N$ are closed and have empty interior for every $N$ large enough. Subsequently, we will implicitly suppose that the union $\bigcup_{N\in\N} E'_N$ is an $F_\sigma$ with empty interior. It is not a limiting assumption as we will focus only on the behaviour of the discretizations for $N$ going to $+\infty$. It will allow us to restrict the study to the $G_\delta$ dense set $\bigcap_{N\in\N} \mathcal{D}_N$, which is the set of homeomorphisms whose $N$-th discretization is uniquely defined for all $N\in\N$.
\end{rem}

\section{Probability measures on $X$}

In the sequel, we will study some \emph{ergodic} properties of discretizations of $f$. Denote by $\Prb$\index{$\mathcal P$} the set of Borel probability measures on $X$ endowed with the weak-star topology: a sequence $(\nu_m)_{m\in\N}$ of $\Prb$ tends to $\nu\in\Prb$ (denoted by $\nu_m\rightharpoonup\nu$) if for every continuous function $\varphi : X\to\R$,
\[\lim_{m\to\infty}\int_{X}\varphi\,\ud \nu_m = \int_{X}\varphi\,\ud \nu.\]
It is well known that under these conditions, the space $\Prb$ is metrizable and compact, therefore separable (Prohorov, Banach-Alaoglu-Bourbaki theorem); we equip this space with a distance $\dist$\index{$\operatorname{dist}$} which is compatible with this topology.


Let $f\in\Hom(X)$. For $x\in \T^2$, we denote by $p\omega (x)$\index{$p\omega (x)$}\label{pomega} the set of limit points of the sequence
\[\left\{\frac{1}{n}\sum_{k=0}^{n-1}\delta_{f^k(x)}\right\}_{n\in\N^*}.\]
It is a compact subset of the set $\mathcal M^f$\index{$\mathcal M^f$} of $f$-invariant Borel probability measures.
\bigskip

We will pay a particular attention to physical measures 

\begin{definition}\label{sport}
A Borel probability measure $\mu$ is called \emph{physical} (sometimes called SRB, see \cite{MR1933431}) for the map $f$ if its basin of attraction has positive $\lambda$-measure, where the \emph{basin of attraction} of $\mu$ for $f$ is the set
\[\left\{x\in X\ \big\vert\ \frac{1}{M} \sum_{m=0}^{M-1} \delta_{f^m(x)} \underset{M\to+\infty}{\longrightarrow} \mu\right\}\]
of points whose Birkhoff's limit coincides with $\mu$.
\end{definition}

Heuristically, the basin of a measure is the set of points that can see the measure, and physical measures are the ones that can be seen in practice.

To study ergodic properties of homeomorphisms and their discretizations, we define natural invariant probability measures associated with these maps:

\begin{definition}\label{defmes}
For any non-empty open subset $U$ of $X$, we denote by $\lambda_{U}$ the normalized restriction of $\lambda$ on $U$, i.e. $\lambda_U=\frac{1}{\lambda(U)} \lambda_{|U}$. We also denote by $\lambda_{N,U}$\index{$\lambda_{N,U}$} the uniform probability measure on $E_N\cap U$.

For $x\in X$ we denote by $\mu^f_{x,m}$ the uniform measure on the segment of orbit $x,f(x),\cdots,f^{m-1}(x)$:\index{$\mu^f_{x,m}$}
\[\mu^f_{x,m} = \frac 1m \sum_{i=0}^{m-1}{f}_*^i \delta_x.\]
When the limit exists, we denote by $\mu^f_x$ the \emph{Birkhoff limit} of $x$ by:\index{$\mu^f_x$}
\[\mu^f_x = \lim_{m\to\infty} \mu^f_{x,m}.\]
Similarly, for $f_N$ and $x\in X$ (recall that $x_N$ is the projection of $x$ on the grid $E_N$):\index{$\mu^{f_N}_{x,m}$}	
\[\mu^{f_N}_{x,m} = \frac 1m \sum_{i=0}^{m-1}(f_N)_*^i \delta_{x_N};\]
in this case the Birkhoff limit is always well defined: it is the uniform measure on the periodic orbit on which the positive orbit of $x_N$ under $f_N$ falls:\index{$\mu^{f_N}_{x}$}	
\[\mu^{f_N}_{x} = \lim_{m\to\infty}\mu^{f_N}_{x,m},\]

We also define a similar quantity for sets of points. We set\index{$\mu^f_{U,m}$}
\[\mu^f_{U,m} = \frac 1m \sum_{i=0}^{m-1}{f}_*^i \lambda_U,\]
and when it is well defined,\index{$\mu^f_U$}
\[\mu^f_U = \lim_{m\to\infty} \mu^f_{U,m}.\]
Likewise,\index{$\mu^{f_N}_{U,m}$}
\[\mu^{f_N}_{U,m} = \frac 1m \sum_{i=0}^{m-1}(f_N)_*^i\, \lambda_{N,U},\]
and in this case the Birkhoff limit always exists: it is a measure supported by the periodic orbits of $f_N$, such that the measure of each periodic orbit is proportional to the number of points of $E_N\cap U$ whose orbit under $f_N$ eventually belongs to this periodic orbit:\index{$\mu^{f_N}_U$}
\[\mu^{f_N}_{U} = \lim_{m\to\infty} \mu^{f_N}_{U,m}.\]
\end{definition}

We have defined two types of invariant measures: one from a point $x$, the other from the uniform measure $\lambda$. The link between them is done by the following remark which easily follows from the dominated convergence theorem.

\begin{rem}\label{convdom}
When $U$ is an open set almost every point (for $\lambda_U$) in which admits a Birkhoff limit, the measure $\mu^f_U$ is well defined and satisfies, for every continuous map $\varphi : X\to\R$, 
\[\int_X\varphi\, \ud\mu^f_U = \int_U\left(\int_X \varphi\, \ud\mu^f_x\right)\ud\lambda_U.\]
Similarly,
\[\int_X\varphi\, \ud\mu^{f_N}_U = \int_U\left(\int_X \varphi\, \ud\mu^{f_N}_{x}\right)\ud\lambda_{N,U}.\]
\end{rem}

\section{Hypothesis on discretization grids}\label{SecHyp}

Previously, we have given a very general definition of the concept of sequence of discretization grids. In some cases, we will need additional technical assumptions about these sequences of grids; of course all of them will be satisfied by good sequences of uniform discretization grids on the torus (as defined in the introduction).

The first assumption is useful for proving Lax's theorem (Theorem \ref{Lax}), and therefore necessary only in the part concerning conservative homeomorphisms.

\begin{definition}[Well distributed and well ordered grids]\label{Ashe}
We say that a sequence of discretization grids $(E_N)_{N\in\N}$ is \emph{well distributed} if we can associate to each $x\in E_N$ a subset $C_{N,x}$\index{$C_{N,x}$} of $X$, which will be called a \emph{cube of order $N$}, such that:
\begin{itemize}
\item for all $N$ and all $x\in E_N$, $x\in C_{N,x}$;
\item for all $N$, $\{C_{N,x}\}_{x\in E_N}$ is a measurable partition of $X$: $\bigcup_{x\in E_N} C_{N,x}$ has full measure and for two distinct points $x,y$ of $E_N$, the intersection $C_{N,x}\cap C_{N,y}$ is null measure;
\item for a fixed $N$, the cubes $\{C_{N,x}\}_{x\in E_N}$ have all the same measure;
\item the diameter of the cubes of order $N$ tends to $0$: $\max_{x\in E_N}\, \diam(C_{N,x})\underset{N\to+\infty}{\longrightarrow}0$.
\end{itemize}

If $(E_N)_{N\in\N}$ is well distributed and if $\{C_{N,x}\}_{N\in\N,x\in E_N}$ is a family of cubes as above, we will say that $(E_N)_{N\in\N}$ is \emph{well ordered} if, for a fixed $N$, the cubes  $\{C_{N,x}\}_{x\in E_N}$ can be indexed by $\Z/q_N\Z$ (recall that $q_N = \card(E_N)$) such that two consecutive cubes (in $\Z/q_N\Z$) are close to each other: 
\[\max_{i\in \Z/q_N\Z}\, \diam(C_N^i\cup C_N^{i+1})\underset{N\to+\infty}{\longrightarrow} 0\]
(in particular, it is true when the boundaries of two consecutive cubes overlap).
\end{definition}

At first glance, it can seem surprising that there is no link between the cubes and the projections. In fact, the existence of such cubes expresses that the grids ``fit'' the measure $\lambda$.

The following definition describe assumptions on grids that will be useful for obtaining properties on average.

\begin{definition}[self similar grids]\label{Ashe'}
We say that a sequence of discretization grids $(E_N)_{N\in\N}$ is \emph{sometimes self similar} if for all $\varep>0$ and all $N_0\in\N$, there exists two integers $N_1$ and $N_2$ satisfying $N_2\ge N_1\ge N_0$, such that the set $E_{N_2}$ contains disjoint subsets $\widetilde E_{N_2}^1, \cdots, \widetilde E_{N_2}^{\alpha}$\index{$\widetilde E_{N_2}^{i}$} whose union fills a proportion greater than $1-\varep $ of $E_{N_2}$:
\[\frac{\card\big(\widetilde E_{N_2}^1 \cup \cdots \cup \widetilde E_{N_2}^{\alpha}\big)}{\card(E_{N_2})}>1-\varep,\]
and such that for all $i$, the set $\widetilde E_N^i$ is the image of the grid $E_{N_1}$ by a bijection $h_i$ which is $\varep$-close to identity.

We say that a sequence of discretization grids $(E_N)_{N\in\N}$ is \emph{sometimes strongly self similar} if it is sometimes self similar for the parameter $\varep = 0$ (that is, the sets $\widetilde E_{N_2}^i$ form a partition of $E_{N_2}$).

We say that a sequence of discretization grids $(E_N)_{N\in\N}$ is \emph{always self similar} if for all $\varep>0$ and all $N_0\in\N$, there exists $N_1\ge N_0$ such that for all $N_2\ge N_1$, the set $E_{N_2}$ contains disjoint subsets $\widetilde E_{N_2}^1, \cdots, \widetilde E_{N_2}^{\alpha_{N_2}}$ whose union fills a proportion greater than $1-\varep $ of $E_{N_2}$:
\[\frac{\card\big(\widetilde E_{N_2}^1 \cup \cdots \cup \widetilde E_{N_2}^{\alpha_{N_2}}\big)}{\card(E_{N_2})}>1-\varep,\]
and such that for all $i$, $\widetilde E_N^i$ is the image of the grid $E_{N_1}$ by a bijection $h_i$ which is $\varep$-close to identity, and such that for all $i,j$ and all $N_2, N_2'\ge N_1$, either $\widetilde E_{N_2}^i\cap \widetilde E_{N_2'}^j = \emptyset$, or $\widetilde E_{N_2}^i = \widetilde E_{N_2'}^j$.

We say that a sequence of discretization grids $(E_N)_{N\in\N}$ is \emph{always strongly self similar} if it is always self similar for the parameter $\varep = 0$ (that is, the sets $\widetilde E_{N_2}^i$ form a partition of $E_{N_2}$).
\end{definition}

\section{Some examples of discretization grids}\label{exgrilles}

In the previous section we set properties about discretizations -- being well distributed, being well ordered, being sometimes/always (strongly) self similar -- for a later use. In this section, we study these properties for some natural examples of grids.

\subsection*{Uniform discretization grids on the torus}\label{GridUnif}

\begin{figure}
\begin{center}
\begin{minipage}[c]{.32\linewidth}
\begin{center}
\begin{tikzpicture}[scale=.6]
\draw[thick] (.5,.5) rectangle (5.5,5.5);
\foreach \k in {0,...,5}
 {\foreach \l in {0,...,5}
  {\draw[color=red!70!black] node at (\k+.5,\l+.5){\footnotesize$\bullet$};}}
\clip (.5,.5) rectangle (5.5,5.5);
\draw[blue] (0,0) grid (6,6);
\end{tikzpicture}
\end{center}
\end{minipage}
\hfill
\begin{minipage}[c]{.32\linewidth}
\begin{center}
\begin{tikzpicture}[scale=.7714]
\draw[blue] (0,0) grid (4,4);
\draw[thick] (0,0) rectangle (4,4);
\foreach \k in {1,...,4}
 {\foreach \l in {1,...,4}
  {\draw[color=red!70!black] node at (\k*4/5,\l*4/5){\footnotesize$\bullet$};}}
\end{tikzpicture}
\end{center}
\end{minipage}
\hfill
\begin{minipage}[c]{.32\linewidth}
\begin{center}
\begin{tikzpicture}[scale=.7714]
\draw[blue] (0,0) grid (4,4);
\draw[thick] (0,0) rectangle (4,4);
\foreach \k in {1,...,4}
 {\foreach \l in {1,...,4}
  {\draw[color=red!70!black] node at (\k-.5,\l-.5){\footnotesize$\bullet$};}}
\end{tikzpicture}
\end{center}
\end{minipage}
\end{center}
\caption[Examples of grids of discretization]{Uniform discretization grids of order $5$ on the torus $\T^2$ ($E_N$, left) and on the cube $I^2$ ($E_N^0$, middle, and $E_N^1$, right) and their associated cubes}
\label{Grill}
\end{figure}
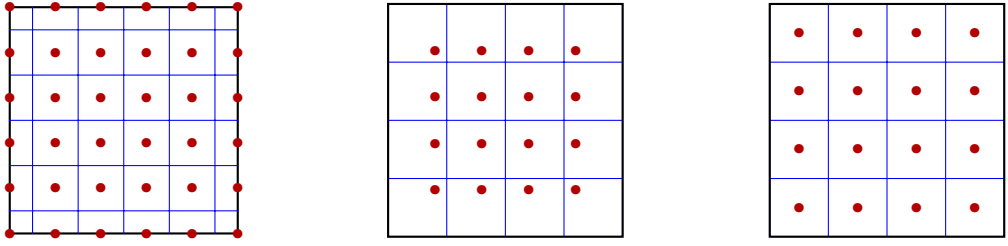

The simplest example of grid of discretization, which will be used for the simulations, is that of the torus $\T^n = \R^n/\Z^n$ of dimension $n\ge 1$ endowed with discretizations called \emph{uniform} discretizations, defined from the fundamental domain $I^n=[0,1]^n$ of $\T^n$: take an increasing sequence of integers $(k_N)_{N\in\N}$ and set (see Figure \ref{Grill}, left)
\[E_N = \left\{\left(\frac{i_1}{k_N},\cdots,\frac{i_n}{k_N}\right)\in \T^n \big|\ \forall j,\, 0\le i_j\le {k_N}-1\right\},\]
\[C_{N,(i_1/N,\cdots,i_n/N)} = \prod_{j=1}^n \left[\frac{i_j-1/2}{k_N}\, ,\, \frac{i_j+1/2}{k_N}\right].\]

We easily verify that this sequence of grids is well distributed, well ordered and always self similar. If we assume that for for every $N_0\in\N$, there exists $N_2\ge N_1\ge N_0$ such that $k_{N_1}$ divides $k_{N_2}$, then the sequence is sometimes strongly self similar. If we further assume that for any $N\in\N$, $k_N$ divides $k_{N+1}$ (which is true when $k_N = p^N$ with $p\ge 2$), then the sequence is always strongly self similar. When $k_N = p^N$ with $p=2$ (respectively, $p=10$) the discretization performs what we can expect from a numerical simulation: doing a binary (respectively, decimal) discretization at order $N$ is the same as truncating each binary (respectively decimal) coordinate of the point $x\in I^n$ to the $N$-th digit, \emph{i.e.} working with a fixed digital precision\label{pagefin}\footnote{In practice, the computer works in floating point format, so that the number of decimal places is not the same when the number is close to $0$ or not.}.

\subsection*{Uniform discretization grids on the cube}

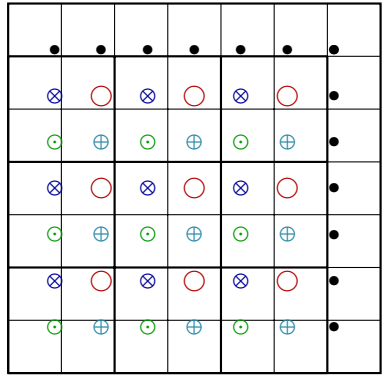
\begin{figure}
\begin{center}
\begin{tikzpicture}[scale=.7]
\draw (0,0) grid (7,7);
\draw[thick] (0,0) rectangle (7,7);
\draw[thick] (0,0) grid[step=2] (6,6);
\foreach \k in {1,...,7}
 {\draw[color=black] node at (\k*7/8,49/8){\footnotesize$\bullet$};
  \draw[color=black] node at (49/8,\k*7/8){\footnotesize$\bullet$};}
\foreach \k in {1,...,3}
 {\foreach \l in {1,...,3}
  {\draw[color=green!60!black] node at (\k*7/4-7/8,\l*7/4-7/8){\footnotesize$\odot$};}}
\foreach \k in {1,...,3}
 {\foreach \l in {1,...,3}
  {\draw[color=blue!60!black] node at (\k*7/4-7/8,\l*7/4){\footnotesize$\otimes$};}}
\foreach \k in {1,...,3}
 {\foreach \l in {1,...,3}
  {\draw[color=cyan!60!black] node at (\k*7/4,\l*7/4-7/8){\footnotesize$\oplus$};}}
\foreach \k in {1,...,3}
 {\foreach \l in {1,...,3}
  {\draw[color=red!70!black] node at (\k*7/4,\l*7/4){\footnotesize$\bigcirc$};}}
\end{tikzpicture}
\caption{Self similarity of grids $E_N^0$}\label{soucubhihi}
\end{center}
\end{figure}

In the case of the cube $I^n = [0,1]^n$, there are two natural types of grids of discretization. In the first one the grids are obtained by taking the same definition as the grids $E_N$ (but avoiding taking points on the boundary of the cube):\index{$E_N^0$}
\[E_N^0 = \left\{\left(\frac{i_1}{k_N},\cdots,\frac{i_n}{k_N}\right)\in I^n \big|\ \forall j,\, 1\le i_j\le {k_N}-1\right\},\]
with an increasing sequence $(k_N)_{N\in\N}$ of integers. To these grids are associated the cubes\footnote{Be careful, these cubes have their vertices on the grid of order $k_N-1$.} (see Figure \ref{Grill}, middle)
\[C_{N,(i_1/N,\cdots,i_n/N)} = \prod_{j=1}^n\left[\frac{i_j-1}{k_N-1}\, ,\, \frac{i_j}{k_N-1}\right].\]
As before, we easily verify that this sequence of grids is well distributed, well ordered and always self similar (see Figure \ref{soucubhihi}). If we assume that for every $N_0\in\N$, there exists $N_2\ge N_1\ge N_0$ such that $k_{N_1}-1$ divides $k_{N_2}-1$, then the sequence is sometimes strongly self similar. If we further assume that for any $N\in\N$, $k_N-1$ divides $k_{N+1}-1$ (that is, for every $N$, there exists an integer $\ell_N>1$ such that $k_{N+1} = k_N + \ell_N k_N(k_N+1)$), and either $k_N$ always divides $k_{N+1}$ or $k_N$ never divides $k_{N+1}$, then the sequence is always strongly self similar.

We can also take discretizations according to the centres of the cubes:\label{pagecubes1}\index{$E_N^1$}
\[E_N^1 = \left\{\left(\frac{i_1+1/2}{k_N},\cdots,\frac{i_n+1/2}{k_N}\right)\in I^n \big|\ \forall j,\, 0\le i_j\le {k_N}-1\right\}.\]
This time, the cubes are much more natural (see Figure \ref{Grill}, right):
\[C_{N,(i_1/N,\cdots,i_n/N)} = \prod_{j=1}^n\left[\frac{i_j}{k_N}\, ,\, \frac{i_j+1}{k_N}\right].\]
Again, we easily verify that this sequence of grids is well distributed, well ordered and always self similar. Moreover, it is sometimes strongly self similar if for every $N_0\in\N$, there exists $N_2\ge N_1\ge N_0$ such that $k_{N_2}/k_{N_1}$ is an odd integer, and always strongly self similar if $k_{N+1}/k_N$ is an odd integer for all $N\in\N$.

\subsection*{Discretization grids on an arbitrary manifold}

In fact, discretizations on arbitrary manifolds can be easily obtained from discretization grids on the cube $I^n$ (for example the grids $E_N^0$) by applying the Oxtoby-Ulam theorem (see \cite{Oxto-meas} or \cite[Corollaire 1.4]{MR2931648}). This theorem asserts that every smooth connected compact manifold of dimension $n$, possibly with boundary, equipped with a good measure $\lambda$, is obtained from the unit cube $I^n$ by gluing together some parts of the boundary of the cube. Moreover, this gluing map can be chosen such that the image of Lebesgue measure on the cube is $\lambda$.

\begin{theoreme}[Oxtoby, Ulam]\label{Brown-mesure}
Under the assumptions that have been made on $X$ and $\lambda$, there exists a map $\phi : I^n\to X$ such that:
\begin{enumerate}
\item $\phi$ is surjective;
\item $\phi_{|\mathring{I^n}}$ is a homeomorphism on its image;
\item $\phi(\partial I^n)$ is a closed subset of $X$ with empty interior and disjoint from $\phi (\mathring{I^n})$;
\item $\lambda(\phi(\partial I^n))=0$;
\item $\phi^*(\lambda) = \Leb$, where $\Leb$ is Lebesgue measure.
\end{enumerate}
\end{theoreme}

This theorem allows to define discretization grids on $X$ from discretizations on the cube; this is outlined by the following informal proposition.

\begin{prop}
If there is a sequence $(E_N)_N$ of grids on $I^n$ whose elements are not on the edge of the cube, then its image by $\phi$ defines a sequence of grids on $X$ which satisfy the same properties as the initial grid on the cube.
\end{prop}

\begin{rem}
The example also includes the case where $X=I^n$, $\lambda=\Leb$ and where the grids are the images of the grids $E_N$ by a unique homeomorphism of $X$ preserving Lebesgue measure.
\end{rem}

\chapter{Discretizations of a generic dissipative homeomorphism}\label{ChapDissip}

In this chapter, we study the dynamical behaviour of the discretizations of generic \emph{dissipative} homeomorphisms, \emph{i.e.} without the assumption of preservation of a given measure.

The dynamics of a generic dissipative homeomorphism is very stable: the \emph{shredding lemma} of F. Abdenur and M. Andersson (Lemma \ref{déchet}) implies that for a generic homeomorphism, there exists a Cantor set of dimension 0 which attracts almost every point of $X$ (see Corollary \ref{convattra}). From that we will deduce easily the dynamics of the discretizations of a generic homeomorphism from that of the homeomorphism itself.

This study of the dynamics of a generic homeomorphism suggests that in the generic dissipative setting, ``the physical dynamics of the discretizations of a homeomorphism $f$ converges to the physical dynamics of $f$'', both from topological and measurable viewpoints. Indeed, we will prove that for a generic homeomorphism and for almost every point $x\in X$, the orbit of $x_N$ by the discretization $f_N$ shadows the orbit of $x$ by $f$ (Corollary \ref{stabilis}). Moreover, we will prove that the recurrent set $\Omega(f_N)$ tends to the closure of the set of Lyapunov stable periodic points of $f$ in a weak sense defined in Proposition~\ref{Hausmodif}; from a combinatorial viewpoint, the cardinality of this recurrent set is as small as possible (Corollary \ref{totsingdiscr}) and the length of the longest periodic orbit tends to $+\infty$ when $N$ goes to $+\infty$ (Corollary \ref{OrbPerDissip}). Finally, we will prove that the measures $\mu^{f_N}_X$ tend to the measure $\mu^f_X$, which is the 
Cesàro limit of the sequence of pushforwards of $\lambda$ by $f$ (Theorem \ref{convmesdissip}).

The end of this chapter will be dedicated to the results of the numerical simulations.

\emph{Throughout this chapter, we fix a compact manifold $X$ of dimension $n\ge 2$. Here, homeomorphisms are not supposed to be conservative; more formally, we recall that we denote by $\Hom(X)$ the set of all \emph{dissipative} homeomorphisms of $X$ (without assumption of conservation of a given measure). The manifold $X$ is equipped with a sequence $(E_N)_{N\in\N}$ of discretization grids and a good measure $\lambda$ as defined in Section~\ref{DefMani} (we will use this measure to look at the physical dynamics, that is the dynamics of almost every point), we suppose that the sequence of grids of discretization satisfy $\sum_{x\in E_N}\delta_x \underset{N\to +\infty}{\longrightarrow} \lambda$. We are interested in properties of discretizations of generic elements of $\Hom(X)$ with respect to the sequence of grids $(E_N)$.}

\section{Dynamics of a generic dissipative homeomorphism}\label{SecShre}

In their article \cite{MR3027586}, F. Abdenur and M. Andersson try to identify the generic ergodic properties of continuous maps and homeomorphisms of compact manifolds. More precisely, they study the behaviour of Birkhoff limits\footnote{We recall that the Birkhoff limit $\mu_x^f$ is the limit in the Cesàro sense of the pushforwards of the Dirac measure $\delta_x$ by the homeomorphism $f$.} $\mu^f_x$ for a generic homeomorphism $f \in \Hom(X)$ and almost every point $x$ for $\lambda$. To do this, they define some interesting behaviours of homeomorphisms related to Birkhoff limits, including one they call \emph{weird}.

\begin{definition}\label{strange}
A homeomorphism $f$ is said \emph{weird} if almost every point $x\in X$ (for $\lambda$) has a Birkhoff's limit $\mu^f_x$, and if $f$ is \emph{totally singular} (\emph{i.e.} there exists a Borel set with total measure whose image by $f$ is null measure) and does not admit any physical measure.
\end{definition}

This definition is supported by their proof, based on the shredding lemma, that a generic homeomorphism is weird. We state an improvement of this lemma, whose main consequence is that a generic homeomorphism has many open attractive sets, all of small measure, and decomposable into a small number of small diameter open sets:

\begin{lemme}[Shredding lemma, F. Abdenur, M. Andersson, \cite{MR3027586}]\label{déchet}
For every homeomorphism $f\in\Hom(X)$, for all $\varep,\delta>0$, there exists integers $\ell$ and $\ell_1,\cdots,\ell_\ell$, bigger than $1/\varep$, and a family of regular pairwise disjoints open sets\footnote{An open set is said \emph{regular} if it is equal to the interior of its closure.} $U_1,\cdots,U_\ell$ such that for all $\varep'>0$, there exists $g\in\Hom(X)$ such that $d(f,g)<\delta$ and:
\begin{enumerate}[(i)]
\item $g(\overline{U_j})\subset U_j$,
\item $\lambda(U_j)<\varep$,
\item $\lambda\left(\bigcup_{j=1}^\ell U_j\right)>1-\varep$,
\item $\lambda(g(U_j))<\varep'\, \lambda(U_j)$,
\item for all $j\le\ell$, there exists open sets $W_{j,1},\cdots,W_{j,\ell_j}$ such that:
  \begin{enumerate}
  \item $\diam(W_{j,i})<\varep'$ for all $i\in\{1,\cdots,\ell_j\}$,
  \item $g(\overline{W_{j,i}})\subset W_{j,i+1}$, for every $i\in\{1,\cdots,\ell_j-1\}$ and $g(\overline{W_{j,\ell_j}})\subset W_{j,1}$,
  \item \[\overline{U_j}\subset\bigcup_{m\ge 0} g^{-m}\big(\bigcup_{i=1}^{\ell_j}W_{j,i}\big).\]
	\item The sets $W_{j,i}$ have disjoints attractive sets, \emph{i.e.} for all $j\neq j'$, we have\footnote{These two sets are decreasing intersections of compact sets, thus compact and non-empty.} 
\[\left(\bigcap_{m\ge 0}\,\bigcup_{m'\ge m} g^{m'}(W_{j,1})\right) \cap \left(\bigcap_{m\ge 0}\,\bigcup_{m'\ge m} g^{m'}(W_{j,i'})\right) = \emptyset.\]
  \item We can further assume that each set $W_{j,i}$ contains a Lyapunov stable periodic point (see Lemma \ref{LyapGene}).
  \end{enumerate}
\end{enumerate}
Moreover, these properties remain true on a neighbourhood of the homeomorphism $g$.
\end{lemme}

We outline the proof of this lemma: by using arguments such as the Oxtoby-Ulam theorem (Theorem \ref{Brown-mesure}) or the concept of uniform grid on the cube $I^n$, it is possible to shorten a little the arguments of \cite{MR3027586}.

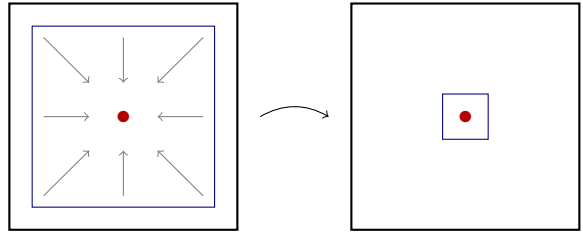
\begin{figure}[t]
 \begin{center}
  \begin{tikzpicture}[scale=1.5]
   \draw[thick] (-3,-1) rectangle (-1,1);
   \draw[red!70!black] (-2,0) node{$\bullet$};
   \draw[blue!50!black] (-2.8,-.8) rectangle (-1.2,.8);
   \draw[->, gray] (-2,.7) -- (-2,.3);
   \draw[->, gray] (-2,-.7) -- (-2,-.3);
   \draw[->, gray] (-2.7,0) -- (-2.3,0);
   \draw[->, gray] (-1.3,0) -- (-1.7,0);
   \draw[->, gray] (-2.7,.7) -- (-2.3,.3);
   \draw[->, gray] (-2.7,-.7) -- (-2.3,-.3);
   \draw[->, gray] (-1.3,-.7) -- (-1.7,-.3);
   \draw[->, gray] (-1.3,.7) -- (-1.7,.3);
   \draw[->] (-.8,0) to[bend left] (-.2,0);

   \draw[thick] (0,-1) rectangle (2,1);
   \draw[red!70!black] (1,0) node{$\bullet$};
   \draw[blue!50!black] (0.8,-.2) rectangle (1.2,.2);
  \end{tikzpicture}
 \end{center}
 \caption{Local perturbation for Lemma \ref{déchet}}\label{PerturbShredding}
\end{figure}

\begin{proof}[Proof of Lemma \ref{déchet}]
Let $f\in\Hom(X)$ and $\delta,\varep>0$. To begin with, we endow $X$ with a collection of ``cubes'': for $N$ large enough, the images of the cubes associated to the the grid $E_N^1$ on the cube $I^n$ (see page \pageref{pagecubes1}) by the map given by Oxtoby-Ulam theorem (Theorem \ref{Brown-mesure}) give us a collection of cubes $(C_i)$ such that their union has full measure, their boundaries have null measure, they all have the same measure and their diameters are smaller than $\delta$. To each cube $C_i$ we associate its centre $e_i$. These properties assert that this collection of cubes behaves like the collection of cubes on $I^n$ endowed with Lebesgue measure. Thus, in the sequel, we will only treat the case of $I^n$ endowed with $\Leb$.

We define a finite map $\sigma : \{e_i\}_i \to \{e_i\}_i$ such that $f(e_i)\in C_{\sigma(e_i)}$. Making a perturbation if necessary, we can suppose that $\sigma$ does not contain any periodic orbit of length smaller than $1/\varep$ (simply because if $N$ is large enough, then for all $e_i$ there exists at least $1/\varep$ other points $e_j$ such that $d(e_i,e_j)<\delta$). We then use Corollary \ref{CoroExtension} to move the points $f(e_i)$ closer to the points $e_j$, in other words we build a homeomorphism $g_1$ which is $\delta$-close to $f$ and such that $\max_i \big(\min_j d(g_1(e_i), e_j)\big) < \min_i \diam(C_i)/10$.

We then set $g_2 = g_1\circ h$, where $h$ is a homeomorphism close to the identity whose restriction to every cube is a huge contraction (see Figure \ref{PerturbShredding}): such a contraction can be easily expressed on the cube $[-1,1]^n$: let $\alpha$ be a ``big'' number (determined by $\min(\varep,\varep')$) and
\[\phi_\alpha(x) = \left\{\begin{array}{ll}
\|x\|_\infty^\alpha x\quad & \text{if}\ \|x\|_\infty\le 1\\
x\quad & \text{otherwise}.
\end{array}\right.\]
We then build easily $h$ by composing such contractions on each cube, each of them being obtained by conjugating by a translation and a homothety.

The ``physical'' dynamics of the map $g_2$ is then close to that of the finite map $\sigma$ which maps $e_i$ on the unique centre $e_j$ such that $g(e_i)\in C_j$. In particular, the periodic sets $W_{j,1},\cdots,W_{j,\ell_j}$ are obtained as neighbourhoods of a given periodic orbit of $\sigma$, and the corresponding basins of attraction $U_j$ are unions of subsets of the cubes\footnote{More precisely, for every cube $C_j$, the sets $U_j$ contain $\{x\in C\mid d(x,\partial C_j)>(1-\varep)\diam(C)\}$.} whose centres have some iterates by $\sigma$ which fall in this periodic orbit. We then easily check that these sets satisfy the conclusions of the lemma, apart from the point \emph{(ii)}.

Thus, we have to prove that the basins of attraction can be supposed to have small measure. To do that, we consider the cubes of order $M=kN$ with $k>1/\varep$, and partition this set of cubes into $k^n$ sets of cubes in the way defined by the fact that these grids are strongly self similar (see page \pageref{pagecubes1}). We then use this self similarity to perturb $f$ such that $\sigma$ stabilizes each one of these subgrids, and apply the same technique as previously.
\end{proof}

\begin{definition}\label{echo}
For a homeomorphism $f\in\Hom(X)$, we denote by $A_0$\index{$A_0$} the set of \emph{Lyapunov stable} periodic points of $f$, \emph{i.e.} the set of periodic points $x$ such that for all $\delta>0$, there exists $\eta>0$ such that if $d(x,y)<\eta$, then $d(f^m(x),f^m(y))<\delta$ for all $m\in\N$.

For two compact sets $K$ and $K'$, by $K'\subset\subset K$ we mean that there exists an open set $O$ such that $K'\subset O\subset K$. In the sequel the set $K$ will be called \emph{strictly periodic} if there exists an integer $i>0$ such that $f^i(K)\subset\subset K$.
\end{definition}

The following lemma ensures that for a generic homeomorphism, every set $W_{j,i}$ of the shredding lemma contains at least one (in fact an infinite number of) Lyapunov stable periodic point.

\begin{lemme}\label{LyapGene}
For a generic homeomorphism $f\in\Hom(X)$, for every strictly periodic topological ball $O$ (\emph{i.e.} there exists $i>0$ such that $f^i(O)\subset\subset O$), there exists a Lyapunov stable periodic point $x\in O$.
\end{lemme}

\begin{proof}[Proof of Lemma \ref{LyapGene}]
We begin by choosing a countable basis of closed sets of $X$: for example we can take $\mathcal K_N$ the set of unions of the closures of the cubes of order $N$. We also denote by $\mathcal B$ the set of all closed topological balls of $X$. We define $\mathcal U_{k,\varep,N}$ as the set of homeomorphisms such that each large enough strictly periodic ball contains a smaller strictly periodic ball with the same period\footnote{For a compact set $K$, $\diam_{int}(K)$ denotes the diameter of the biggest euclidean ball included in $K$.}:
\[\mathcal U_{k,\varep,N} = \left\{ f\in\Hom(X)\ \middle\vert\ \begin{array}{l}
 \forall K\in\mathcal K_N\cap\mathcal B \text{ s.t. } \exists i\le k \text{ s.t. }\\
 f^i(K) \subset\subset K \text{ and } \diam_{int}(K)>\varep,\\
 \exists K'\subset K, K'\in\mathcal B \text{ s.t.\,} \diam(K')<\varep/2\\
 \text{and } f^i(K') \subset\subset K'
\end{array}\right\}.\]
Then for every $k,\varep,N$, it is straightforward that the set $\mathcal U_{k,\varep,N}$ is an open subset of $\Hom(X)$. To show that it is dense it suffices to apply Brouwer's theorem to each $K$ such that $f^i(K) \subset\subset K$ and to make the obtained periodic point attractive.

We now prove that every $f\in \bigcap_{k,\varep,N} \mathcal U_{k,\varep,N}$ satisfies the conclusions of the lemma. First of all, remark that for every topological ball $K$ with non-empty interior which is strictly $i$-periodic, there exits $N\in\N$ and a smaller topological ball $\tilde K\subset K$ which is strictly $i$ periodic such that $\tilde K \in \mathcal K_N$. It implies that if $f$ belongs to the $G_\delta$ dense set $\bigcap_{k,\varep,N} \mathcal U_{k,\varep,N}$, then for every topological ball $K$ with non-empty interior which is strictly $i$-periodic, there exits $N\in\N$ and a topological ball $K'\subset \tilde K\subset K$ which is strictly $i$ periodic and at least twice smaller. Taking the intersection of such balls, we obtain a periodic point with period $i$ which is Lyapunov stable by construction.
\end{proof}

The shredding lemma tells us a lot about the dynamics of a generic homeomorphism, which becomes quite clear: there are many attractors whose basins of attraction are small and attract almost all the manifold $X$. Moreover there is convergence of the attractive sets of the shredding lemma to the closure of the set of Lyapunov stable periodic points of $f$.

\begin{coro}\label{convattra}
Let $f\in\Hom(X)$ verifying the conclusions of the shredding lemma for all $\varep=\varep'>0$, and $W_{j,i,\varep}$ be the corresponding open sets. Such homeomorphisms form a $G_\delta$ dense subset of $\Hom(X)$. Then the sets\index{$A_\varepsilon$}
\[A_\varep = \overline{\bigcup_{j,i} W_{j,i,\varep}}\]
converge for Hausdorff distance when $\varep$ tends to $0$ to a closed set which coincides generically with the set $\overline{A_0}$ (see Definition \ref{echo}).

Moreover, generically, the set $A_0$ is a Cantor set (that is, it is compact, without any isolated point and totally disconnected) whose Hausdorff dimension\footnote{And better, if we are given a countable family $(\lambda_m)_{m\in\N}$ of good measures, generically the Hausdorff dimension of this set with respect to these measures is zero.} is 0.
\end{coro}

\begin{rem}
Thus, for a generic homeomorphism $f\in\Hom(X)$, the $\omega$-limit set of almost every point $x\in X$ is included in $\overline{A_0}$.
\end{rem}

\begin{proof}[Proof of Corollary \ref{convattra}]
Let $f$ verifying the hypothesis of the corollary. We want to show that the sets $A_\varep$ tend to $A_0$ for Hausdorff distance when $\varep$ goes to $0$. This is equivalent to show that for all $\delta>0$, there exists $\varep_0>0$ such that for all $\varep<\varep_0$, $A_0\subset B(A_\varep,\delta)$ and $A_\varep\subset B(A_0,\delta)$ (where $B(A,\delta)$ denotes the set of points of $X$ whose distance to $A$ is smaller than $\delta$). Subsequently we will denote by $U_{j,\varep}$ and by $W_{j, i,\varep}$ the open sets given by the shredding lemma for the parameters $\varep=\varep'$.

Let $\delta> 0$. We start by taking $x\in X$ whose orbit is periodic with period $p$ and Lyapunov stable. Then there exists $\eta>0$ such that if $d(x,y)<\eta$, then $d(f^m(x),f^m(y))<\delta/2$ for all $m\in\N$; we note $O=B(x, \eta)$. As $\lambda(O)>0$ there exists $\varep_0>0$ such that for all $\varep\in]0,\varep_0[$, there exists $j\in\N$ such that the intersection between $O$ and $U_{j,\varep}$ is non-empty. Let $y$ be an element of this intersection. By compactness, there exists a subsequence of $(f^{pm}(y))_{m\in\N}$ which tends to a point $x_0$; moreover $f^m(y)\in \bigcup_{i}W_{j,i}$ eventually and $d(f^{pm}(x),f^{pm}(y))<\delta/2$, thus at the limit $m \to +\infty$, $d(x,x_0)<\delta/2$. We deduce that $x\in  B\big(\overline{\bigcup_{j,i} W_{j,i,\varep_0}},\delta/2\big)$ for all $\varep_0$ small enough. Since $\overline{A_0} $ is compact, it is covered by a finite number of balls of radius $\delta/2$ centred at some points $x_i$ whose orbits attract non-empty open sets. Taking $\varep_0'$ the 
minimum of all the $\
\varep_0$ associated to the $x_i$, the inclusion $A_0\subset B(A_\varep,\delta)$ occurs for all $\varep <\varep_0'$.

Conversely, let $\delta>0$, $\varep<\delta$ and focus on the set $W_{j,i,\varep}$. By Lemma \ref{LyapGene}, we can suppose that there exists $x\in W_{j,i,\varep}$ whose orbit is periodic and Lyapunov stable. Thus $x\in A_0$ and since the diameter of $W_{j,i,\varep}$ is smaller than $\delta$, $W_{j,i,\varep}\subset B(A_0,\delta)$.
\bigskip

We now prove that the set $\overline {A_0}$ has no isolated point. Let $x\in \overline {A_0}$ and $\delta>0$, we want to find another point $y\in \overline {A_0}$ such that $d(x,y)<\delta$. If $x\in \overline{A_0} \setminus A_0$ then it is trivial that $x$ is an accumulation point of $A_0$; thus we suppose that $x\in A_0$. As $x$ is a Lyapunov-stable periodic point, there exists a neighbourhood $O$ of $x$ whose diameter is smaller than $\delta$ and which is periodic (\emph{i.e.} there exists $t>0$ such that $f^t(0)\subset O$). For $\varep$ small enough, the open set $O$ meets at least two sets $U_{j,\varep}$ and $U_{j',\varep}$ of the shredding lemma. Thus, applying the shredding lemma, we deduce that $O$ contains at least two different strictly periodic sets $W_{j,i,\varep}$ and $W_{j',i',\varep}$; by Lemma \ref{LyapGene} they both contain a Lyapunov-stable periodic point, and by construction their distance to $x$ is smaller than $\delta$.

Finally we prove that generically, the set $\overline {A_0}$ has 0 Hausdorff dimension. Let $s>0$. We consider the set of homeomorphisms verifying the conclusions of the shredding lemma for $\varep'$ such that $\sum_j \ell_j {\varep'}^s <1$. This equality implies that for all $\varep>0$, the $s$-Hausdorff measure of the set $\bigcap_{\delta\in]0,\varep[} A_\delta$ is smaller than $1$. As we have
\[\overline{A_0} = \bigcup_{k\ge 0} \bigcap_{\delta\in]0,1/k[} A_\delta,\]
the set $\overline{A_0}$ is a countable union of sets of Hausdorff dimension smaller than $s$ for every $s>0$. Thus, $\overline{A_0}$ has zero Hausdorff dimension; this also proves that $\overline{A_0}$ is perfect.
\end{proof}

\begin{rem}\label{RemPtsPerDissip}
We can also prove that for every $\tau>0$, the set of periodic points with period $\tau$ of a generic homeomorphism is either empty, either a Cantor set with zero Hausdorff dimension.
\end{rem}

From the shredding lemma, we can easily deduce that the sequence of pushforwards of the uniform measure on an open set $U$ is Cauchy. Thus, it converges to a measure $\mu_U^f$, which is supported by $\overline A_0$ and is atomless (because the $\mu_U^f$-measure of each set $\bigcup_i W_{j,i,\varep}$ is between $\varep$ and $2\varep$).

\begin{coro}\label{MesDissip}
For a generic homeomorphism $f\in\Hom(X)$, for every open subset $U$ of $X$, the measure $\mu_U^f$ is well defined, atomless and is supported by $\overline{A_0}$.
\end{coro}

\section{Dynamics of discretizations of a generic homeomorphism}\label{label}

We now establish a discrete counterpart of the shredding lemma. Since having basins of attraction is stable by perturbation, we have a similar statement for discretizations of a generic homeomorphism. In what follows these arguments are developed.

To each point $x_{N} \in E_{N}$, we associate a closed set $\overline{P_N^{-1}(\{x_N\})}$, made of the points in $X$ one of whose projections on $E_{N}$ is $x_{N}$. The closed sets $\overline{P_N^{-1}(\{x_N\})}$ form a basis of the topology of $X$ when $N$ runs through $\N$ and $x_{N}$ runs through $E_{N}$.\label{CNxN} Let $f\in\Hom(X)$, $N\in\N$ and $\vartheta : \N\to\R_+^*$ be a function that tends to $+\infty$ at $+\infty$. Let $\delta,\varep>0$ and $U_j$ be the sets obtained by the shredding lemma for $f$, $\delta$ and $\varep$.

For all $j\le \ell$ we denote by $\widetilde U_j^{N}$ the union over $x_{N}\in E_{N}$ of the closed sets $\overline{P_N^{-1}(\{x_N\})}$ whose intersection with $U_j$ are non-trivial. Then $\widetilde U_j^{N}$ tends to $U_j$ for the metric $d(A,B) = \lambda(A\Delta B)$, and for the Hausdorff metric; the shredding lemma states that these convergences are independent from the choice of $\varep'$. Thus, for all $k$ big enough, properties $(i)$ to $(iii)$ of the shredding lemma remain true for the discretizations $g_{N}$ (for arbitrary $g$ satisfying the properties of the lemma). Taking $\varep'$ small enough and modifying a little $g$ if necessary, there exists sets $W_{j,i}$ and $g\in\Hom(X)$ such that $d(f,g)<\delta$, $W_{j,i}\subset\overline{P_N^{-1}(\{x_N\})}$ and $\card(\bigcup_{j,i} W_{j,i}\cap E_N)\le\vartheta(N)$. The others estimations over the sizes of the sets involved in the lemma are obtained similarly. Finally, we obtain the following lemma:

\begin{lemme}[Discrete shredding lemma]\label{déchetdiscr}
For all $f\in \Hom(X)$, for all $\varep,\delta>0$ and all function $\vartheta : \N\to\R_+^*$ that tends to $+\infty$ at $+\infty$, there exists integers $N_0$, $\ell$ and $\ell_1,\cdots,\ell_\ell$ bigger than $1/\varep$, and a homeomorphism $g\in\Hom(X)$ such that $d(f,g)<\delta$ and for all $N\ge N_0$, there exists a family of subsets $U_1^{N},\cdots,U_\ell^{N}$ of $E_{N}$ such that:
\begin{enumerate}[(i)]
\item $g_{N}(  U_j^{N})\subset   U_j^{N}$,
\item $\card(U_j^{N})<\varep q_{N}$,
\item $\card\left(\bigcup_{j=1}^\ell   U_j^{N}\right)>(1-\varep)q_{N}$,
\item $\card(g_{N}(U_j^{N}))<\card(U_j^{N})/(\ell \vartheta(N_0))$,
\item for all $j$, there exists subsets $W_{j,1}^{N},\cdots,W_{j,\ell_j}^{N}$ of $E_N$ such that
  \begin{enumerate}
  \item $\card(\bigcup_{j,i} W_{j,i}^N)\le q_N / \vartheta(N_0)$,
  \item $g_N(W_{j,i}^N) \subset W_{j,i+1}^N$, for every $i\in\{1,\cdots,\ell_j-1\}$ and $g_N(W_{j,\ell_j}^N) \subset W_{j,1}^N$,
  \item \[U_j^{N}\subset\bigcup_{m\ge 0} g_{N}^{-m}\big(\bigcup_{i=1}^{\ell_j}  W_{j,i}^{N}\big),\]
  \end{enumerate}
\item for all $j$ and all $i$, there exists sets $U_j$ and $W_{j,i}$ satisfying properties $(i)$ to $(v)$ of the shredding lemma such that $U_j\subset P_N^{-1}(U_j^{N})$ and $W_{j,i}\subset P_N^{-1}(W_{j,i}^{N})$.
\item if $N\ge N_0$, then for all $j$ and all $i$, we have
\[\sum_j d_H(\overline{U_j},\overline{P_N^{-1}(U_j^{N})})<\varep\quad\text{and}\quad\sum_{i,j} d_H(\overline{W_{j,i}},\overline{P_N^{-1}(W_{j,i}^{N})})<\varep\]
(where $d_H$\index{$d_H$} is the Hausdorff metric), and
\[\sum_j \lambda(U_j\Delta P_N^{-1}(U_j^{N}))<\varep\quad\text{and}\quad\sum_{i,j} \lambda(W_{j,i}\Delta P_N^{-1}(W_{j,i}^{N}))<\varep.\]
\end{enumerate}
\end{lemme}

\begin{rem}
Properties \emph{(i)} to \emph{(v)} are discrete counterparts of properties \emph{(i)} to \emph{(v)} of the continuous shredding lemma, but the properties \emph{(vi)} and \emph{(vii)} reflect the convergence of the dynamics of discretizations to that of the original homeomorphism.
\end{rem}

This lemma implies that we can theoretically deduce the behaviour of a generic homeomorphism from the dynamics of its discretizations. The next section details this remark.
\bigskip

To begin with, we deduce from the shredding lemma that the dynamics of discretizations $f_N$ tends to that of the homeomorphism $f$. More precisely almost all orbits of the homeomorphism are \emph{$\delta $-shadowed} by the orbits of the corresponding discretizations.

\begin{definition}
Let $f$ and $g$ be two maps from a metric space $X$ into itself, $x,y\in X$ and $\delta>0$. We say that the orbit of $x$ by $f$ \emph{$\delta$-shadows} the orbit of $y$ by $g$ if for all $m\in\N$, $d(f^m(x),g^m(y))<\delta$.
\end{definition}

\begin{coro}\label{stabilis}
For a generic homeomorphism $f\in\Hom(X)$, for all $\varep>0$ and all $\delta>0$, there exists an open set $A$ such that $\lambda(A)>1-\varep$ and $N_0\in\N$, such that for all $N\ge N_0$ and all $x\in A$, the orbit of $x_N = P_N(x)$ by $f_N$ $\delta$-shadows the orbit of $x$ by $f$.

Therefore, for a generic homeomorphism $f$, there exists a full measure dense open set $O$ such that for all $x\in O$, all $\delta>0$ and all $N$ large enough, the orbit of $x_N$ by $f_N$ $\delta$-shadows that of $x$ by $f$.
\end{coro}

\begin{proof}[Proof of Corollary \ref{stabilis}]
This easily follows from the discrete shredding lemma, and especially from the fact that the sets $W_{j,i}^N$ tend to the sets $W_{j,i}$ for Hausdorff metric, in particular $O = \bigcup_{j,\varep} U_{j,\varep}$.
\end{proof}

This statement is a bit different from the genericity of shadowing (see \cite{MR1711347}): here the starting point is not a pseudo-orbit but a point $x\in X$; Corollary \ref{stabilis} expresses that we can ``see'' the dynamics of $f$ on that of $f_N$, with arbitrarily high precision, provided that $N$ is large enough. Among other things, this allows us to observe the basins of attraction of the neighbourhoods of the Lyapunov stable periodic points of $f$ on discretizations. Better yet, to each family of attractors $(W_{j,i})_i$ of the basin $U_j$ of the homeomorphism corresponds a unique family of sets $(W_{j,i}^N)_i$ that are permuted cyclically by $f_N$ and attract a neighbourhood of $U_j$. Thus, attractors are shadowed by cyclic orbits of $f_N$ and we can detect the ``period'' of the attractor (\emph{i.e.} the integer $\ell_j$) on discretizations: when $N=N_0$, the sets $W_{j,i}^N$ each contain only one point; a phenomenon of period multiplication might appear for $N$ bigger \cite{MR1037009,MR875433,
MR765293,MR1299502,MR1980335}. This behaviour is the opposite of what happens in the conservative case, where discretized orbits and true orbits are very different for most points.

Again, in order to show that the dynamics of discretizations converges to that of the initial homeomorphism, we establish the convergence of attractive sets of $f_N$ to that of $f$. Recall that $A_0$ is the closure of the set of Lyapunov stable periodic points of $f$ (see Definition \ref{echo}).

\begin{prop}\label{Hausmodif}
For a generic homeomorphism $f\in\Hom(X)$, the recurrent sets $\Omega(f_N)$ tend weakly to $A_0$ in the following sense: for all $\varep>0$, there exists $N_0\in \N$ such that for all $N\ge N_0$, there exists a subset $\widetilde E_N$ of $E_N$, stabilized by $f_N$, such that, noting $\widetilde \Omega(f_N)$ the corresponding recurrent set, we have $\frac{\card(\widetilde E_N)}{\card(E_N)}>1-\varep$ and $d_H(A_0,\widetilde \Omega(f_N))<\varep$.
\end{prop}

\begin{proof}[Proof of Proposition \ref{Hausmodif}]
Let $\varep>0$. For all $N\in\N$, let $\widetilde E_N$ be the union of the sets $U_j^{N}$ of Lemma \ref{déchetdiscr} for the parameter $\varep$. This lemma ensures that $\widetilde E_N$ is stable by $f_N$ and fills a proportion greater than $1-\varep$ of $E_N$. We also denote by $\widetilde \Omega(f_N)$ the associated recurrent set:
\[\widetilde \Omega(f_N) = \bigcup_{x\in\bigcup_{j=1}^\ell \widetilde U_j^{N}} \omega_{f_N}(x).\]
Property \emph{(viii)} of Lemma \ref{déchetdiscr} ensures that
\[\underset{N\to +\infty}{\overline\lim} d_H(A_\varep,\widetilde \Omega(f_N))<\varep.\]
To conclude, it suffices to apply Corollary \ref{convattra} which asserts that $A_\varep\to A_0$ for Hausdorff distance.
\end{proof}

We now set a final consequence of Lemma \ref{déchetdiscr}, which reflects that the ratio between the cardinality of the image of discretizations and which of the grid is smaller and smaller:

\begin{coro}\label{totsingdiscr}
Let $\vartheta:\R_+ \to \R_+^* $ be a function that tends to $+\infty$ at $+\infty$. Then for a generic homeomorphism $f\in\Hom(X)$,
\[\varliminf_{N\to+\infty} \vartheta(N)\frac{\card(f_{N}(E_{N}))}{\card(E_{N})} = 0\, ;\]
more precisely, for every $M\in\N$, there exists$N_0\ge M$ such that for every $N\ge N_0$, we have
\[\frac{\card(f_{N}(E_{N}))}{\card(E_N)} \le \frac{1}{\vartheta(N_0)}.\]
In particular, the degree of recurrence satisfies $\lim_{N\to +\infty} D(f_N) = 0$.
\end{coro}

Remark that as $\Omega(f_N)\subset f_N(E_N)$, the same estimation holds for the recurrent set. This corollary can be seen as a discrete analogue of the fact that a generic homeomorphism is totally singular, \emph{i.e.} that there exists a Borel set of full measure whose image under $f$ is zero measure. Again, it reflects the regularity of the behaviour of the discretizations of a dissipative homeomorphism: generically, the behaviour of \emph{all} (sufficiently fine) discretizations is the same as the physical behaviour of the initial homeomorphism. This is very different from the conservative case, where sometimes $f_N(E_N) =  E_N$ and sometimes $\card(f_N(E_N)) \le \vartheta(N)$ where $\vartheta:\R_+ \to \R_+^* $ is a given map that tends to $+\infty$ at $+\infty$.

The discrete shredding lemma also allows us to have an estimation about the combinatorial behaviour of ${f_N}_{\Omega(f_N)}$: $f_N$ has a lot of periodic orbits, and along these periodic orbits a lot have long lengths.

\begin{coro}\label{OrbPerDissip}
For a generic homeomorphism $f\in\Hom(X)$ and for every $M\in\N$, there exists $N_0\in\N$ such that for every $N\ge N_0$, there exists a subset $\widetilde E_N$ of $E_N$, stabilized by $f_N$, such that we have $\frac{\card(\widetilde E_N)}{\card(E_N)}>1-\varep$ (as in Proposition \ref{Hausmodif}), such that all the periodic orbits of ${f_N}_{|\widetilde E_N}$ have length bigger than $M$, and such that the number of such orbits is bigger than $M$.
\end{coro}

It remains to study the behaviour of measures $\mu^{f_N}_U$ (see Definition \ref{defmes}). To do that, we have to suppose that the sequence of grids well behaves with respect to the measure $\lambda$. Again, the results are very different from the conservative case: for any open set $U$, the measures $\mu^{f_N}_U$ tend to a single measure, say $\mu^f_U$.

\begin{theoreme}\label{convmesdissip}
For a generic homeomorphism $f\in\Hom(X)$ and an open subset $U$ of $X$, the measure $\mu^f_U$ is well defined\footnote{In other words, a generic homeomorphism is weird, see Definition \ref{strange}, see also \cite{MR3027586}.} and is supported by the set $\overline{A_0}$. Moreover the measures $\mu^{f_N}_U$ tend weakly to $\mu^f_U$.
\end{theoreme}

\begin{proof}[Sketch of proof of Theorem \ref{convmesdissip}]
The proof of this theorem is based on the shredding lemma: the set of homeomorphisms which satisfy the conclusions of the lemma is a $G_\delta$ dense, so it suffices to prove that such homeomorphisms $f$ satisfy the conclusion of the proposition. Let $U$ be an open subset of $X$ and $\varphi : X\to\R$ be a continuous function.

We want to show that on the one hand the integral $\int_X \varphi\, \ud\mu^f_U$ is well defined, \emph{i.e.} that the Birkhoff limits for the function $\varphi$
\[\lim_{m\to+\infty} \frac{1}{m}\sum_{i=0}^{m-1} \varphi(f^i(x))\]
are well defined for almost every $x\in U$; and on the other hand we have the convergence
\[\int_X  \varphi \,\ud\mu^{f_N}_U\underset{N\to+\infty}{\longrightarrow}\int_X \varphi\, \ud\mu^f_U,\]

For the first step, the idea of the proof is that most of the points (for $\lambda$) eventually belong to a set $W_{j,i}$. Since the iterates of the sets $W_{j,i}$ have small diameter, by uniform continuity, the function $\varphi$ is almost constant on the sets $f^m(W_{j,i})$. Thus the measure $\mu^f_x$ is well defined and almost constant on the set of points whose iterates eventually belong to $W_{j,i}$. And by the same construction, since the dynamics of $f_N$ converge to that of $f$, and in particular that the sets $U_j^N$ and $\{w_{j,i}^N\}$ converge to the sets $U_j$ and $W_{j,i}$, the measures $\mu_{U}^{f_N}$ tend to the measures $\mu_U^f$.
\end{proof}

\section{Numerical simulations}\label{partietrois}

We now present some numerical simulations of dissipative homeomorphisms. Again, our aim is to compare the theoretical results with the reality of numerical simulations: for simple homeomorphisms and reasonable orders of discretization, do we have convergence of the dynamics of the discretizations to that of the homeomorphism, as suggested by the above theorems?

We simulate homeomorphisms of the form
\[f(x,y) = (R\circ Q\circ P)(x,y),\]
where $P$ and $Q$ are two homeomorphisms of the torus that modify only one coordinate:
\[P(x,y) = \big(x,y+p(x)\big)\quad\text{and}\quad Q(x,y) = \big(x+q(y),y\big),\]
so that the homeomorphism $Q\circ P$ preserves Lebesgue measure.
We discretize these examples according to the uniform grids on the torus
\[E_N = \left\{\left(\frac{i_1}{N},\cdots,\frac{i_n}{N}\right)\in \T^n \big|\ \forall j,\, 0\le i_j\le {N}-1\right\}.\]

We have tested two homeomorphisms:
\begin{itemize}\label{defbis}
\item To begin with we studied $f_1 = R_1\circ Q\circ P$, with
\[p(x) = \frac{1}{209}\cos(2\pi\times 187x)+\frac{1}{271}\sin(2\pi\times 253 x)-\frac{1}{703}\cos(2\pi\times 775 x),\]
\[q(y) = \frac{1}{287}\cos(2\pi\times 241y)+\frac{1}{203}\sin(2\pi\times 197 y)-\frac{1}{841}\sin(2\pi\times 811 y)\]
and
\begin{align*}
\big(R_1(x,y)\big)_x = x & - 0.00227\sin(2\pi\times 95(x+\alpha))\\
                         & + 0.000224\cos(2\pi\times 197(y+\alpha))\\
												 & - 0.00111\sin(2\pi\times 343(x+\alpha))
\end{align*}
\begin{align*}
\big(R_1(x,y)\big)_y = y & - 0.00376\sin(2\pi\times 107(y+\beta))\\
                         & - 0.000231\cos(2\pi\times 211(x+\beta))\\
												 & + 0.00107\cos(2\pi\times 331(y+\beta)),
\end{align*}
with $\alpha = 0.00137$ and $\beta = 0.00159$.
This dissipative homeomorphism is a small $C^0$ perturbation of the identity, whose derivative has many oscillations whose amplitudes are close to $1$. That creates many fixed points which are attractors, sources or saddles.
\item It has also seemed to us useful to simulate a homeomorphism close to the identity in $C^0$ topology, but with a small number of attractors. Indeed, as explained heuristically by J.-M. Gambaudo and C. Tresser in \cite{MR700317}, a homeomorphism like $f_1$ can have a large number of attractors whose basins of attraction are small. It turns out that the dissipative behaviour of $f_1$ cannot be detected for reasonable orders discretization. We therefore defined another homeomorphism close to the identity in $C^0$ topology, but with much less attractors, say $f_2 =R_2 \circ Q \circ P $, with $P$ and $Q$ identical to those used for $f_1$, and $R_2$ defined by
\[R_2\begin{pmatrix} x \\ y \end{pmatrix} = \begin{pmatrix} x - 0.00227\sin(2\pi\times 14(x+\alpha)) + 0.000324\cos(2\pi\times 33(y+\alpha))\\
y - 0.00376\cos(2\pi\times15(y+\beta)) - 0.000231\sin(2\pi\times 41(x+\beta))\end{pmatrix},\]
with $\alpha = 0.00137$ and $\beta = 0.00159$.
\end{itemize}

Remark that we have chosen to define the homeomorphisms we compute with lacunary trigonometric series, to ``mimic'' the action of Baire theorem.

\subsection{Algorithm used for the calculus of invariant measures}\label{Blabla128}

The algorithm we used to conduct simulations is quite fast (in fact, it is linear in the number of points of the grid). It detects all the periodic orbits of the discretizations $f_N$ in the following way. It takes a first point $x_1\in E_N$ and iterates it until the orbit meets a point that has ever been visited by the orbit. The points belonging to the orbit of $x_1$ are labelled as falling into the periodic orbit number 1. The algorithm also notes the number of points that have been attracted by this orbit, and the coordinates of the points of the periodic orbit. It then takes another point $x_2\in E_N$ which does not belong to the orbit of $x_1$. There are two cases: either an iterate of $x_2$ is equal to an iterate of $x_1$, and in this case it updates the number of points which fall into the periodic orbit number 1; or an iterate of $x_2$ meets another iterate of $x_2$, and in this case it creates a periodic orbit number 2, which attracts all the orbit of $x_2$. This procedure is iterated until there is no more points of $E_N$ that have not been visited. Remark that this algorithm computes the image of a point at most twice.

This algorithm allows to compute quantities like the cardinality of the recurrent set $\Omega(f_N)$, the number of periodic orbits of $f_N$, the maximal size of a periodic orbit of $f_N$, etc. It also allows to represent the invariant measure $\mu^{f_N}_X$ of $f_N$. Recall that this measure is defined as the limit in the Cesàro sense of the push forward of the uniform measures on $E_N$ by the discretizations $f_N$. It is supported by the union $\Omega(f_N)$ of the periodic orbits of $f_N$; the measure of each of these periodic orbits is proportional to the size of its basin of attraction.

We present images of sizes $128\times 128$ pixels representing in logarithmic scale the density of the measures $\mu^{f_N}_X$: each pixel is coloured according to the measure carried by the set of points of $E_N$ it covers. Blue corresponds to a pixel with very small measure and red to a pixel with very high measure. Scales on the right of each image corresponds to the measure of one pixel on the $\log 10$ scale: if green corresponds to $-3$, then a green pixel will have measure $10^{-3}$ for $\mu^{f_N}_X$. For information, when Lebesgue measure is represented, all the pixels have a value about $-4.2$.\label{pagealgo}
 
We also compute the distance between the measure $\mu^{f_N}_X$ and Lebesgue measure. The distance we have chosen is given by the formula
\[d(\mu,\nu) = \sum_{k=0}^\infty \frac{1}{2^k} \sum_{i,j=0}^{2^k-1} \big| \mu(C_{i,j,k}) - \nu(C_{i,j,k})\big|\in[0,2],\]
where
\[C_{i,j,k} = \left[\frac{i}{2^k},\frac{i+1}{2^k}\right] \times \left[\frac{j}{2^k},\frac{j+1}{2^k}\right].\]
This distance spans the weak-* topology, which makes compact the set of probability measures on $\T^2$. In practice, we have computed an approximation of this quantity by summing only on the $k\in\llbracket 0,7 \rrbracket$.

From a practical point of view, we have restricted ourselves to grids of sizes smaller than $2^{15}\times 2^{15}$: the initial data become quickly very large, and the algorithm creates temporary variables that are of size of the order of five times the size of the initial data. For example, for a grid $2^{15}\times 2^{15}$, the algorithm needs between $25$ and $30$ Go of RAM, and takes about two days of calculus for a single order of discretization.

\subsection{Combinatorial behaviour}\label{simulgrafdissip}

\begin{figure}[ht]
\begin{center}
\makebox[0.8\textwidth]{\parbox{0.8\textwidth}{%
\begin{minipage}[c]{.49\linewidth}
	\includegraphics[width=\linewidth,trim = .5cm .3cm .6cm .1cm,clip]{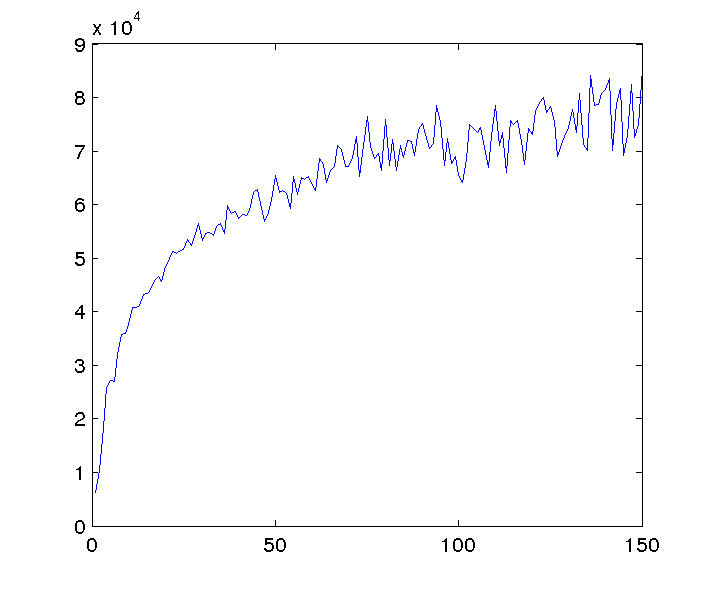}
\end{minipage}\hfill
\begin{minipage}[c]{.49\linewidth}
	\includegraphics[width=\linewidth,trim = .5cm .3cm .6cm .1cm,clip]{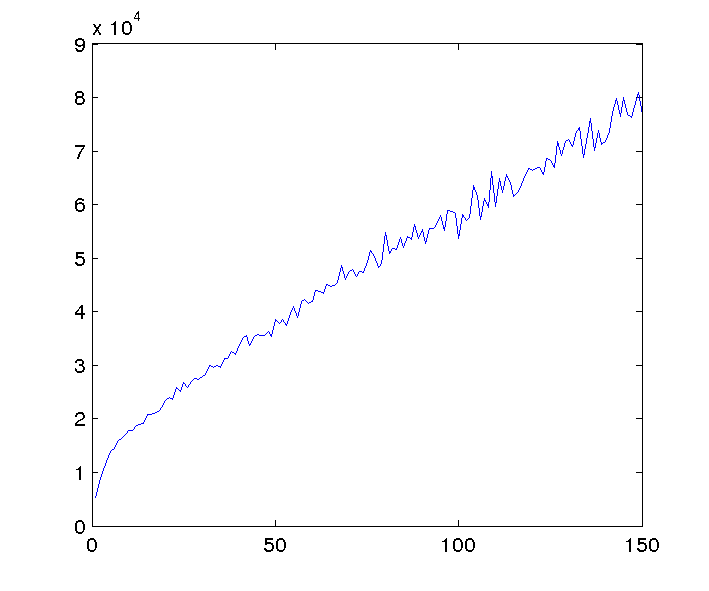}
\end{minipage}

\begin{minipage}[c]{.49\linewidth}
	\includegraphics[width=\linewidth,trim = .5cm .3cm .6cm .1cm,clip]{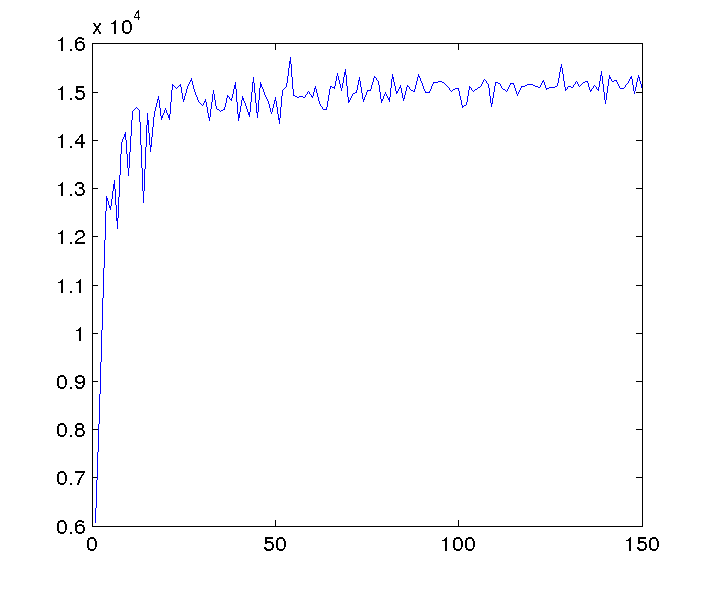}
\end{minipage}\hfill
\begin{minipage}[c]{.49\linewidth}
	\includegraphics[width=\linewidth,trim = .5cm .3cm .6cm .1cm,clip]{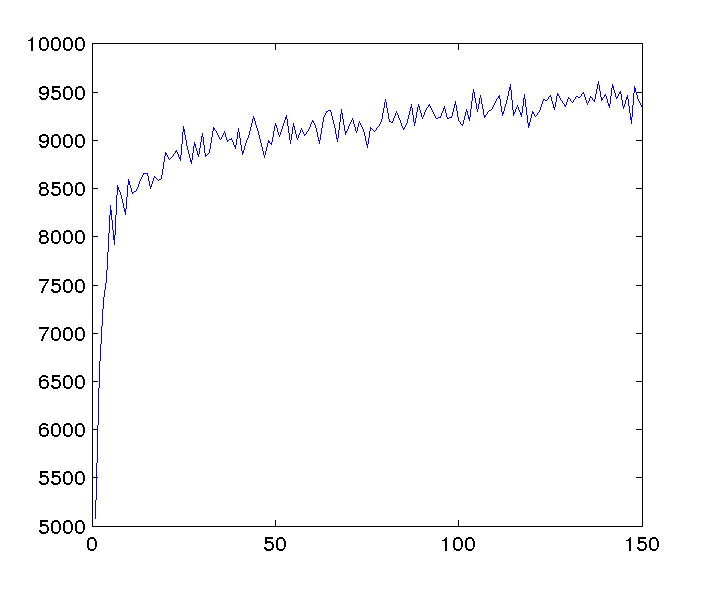}
\end{minipage}

\begin{minipage}[c]{.49\linewidth}
	\includegraphics[width=\linewidth,trim = .5cm .3cm .6cm .1cm,clip]{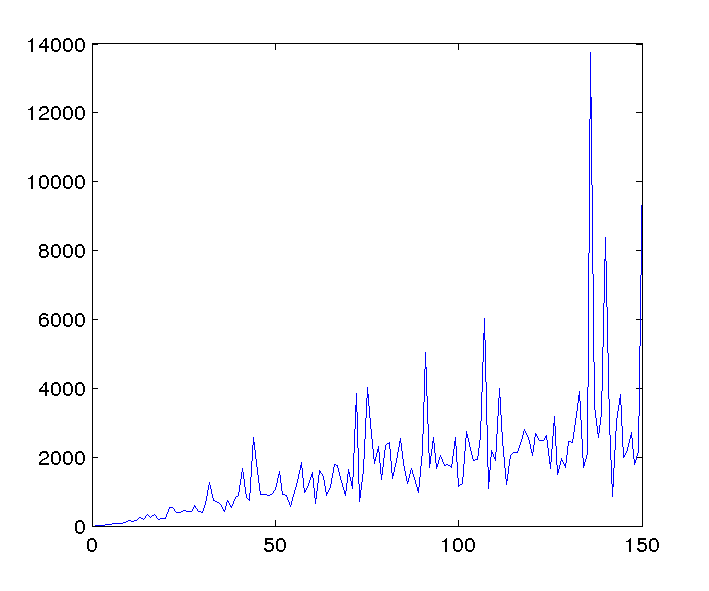}
\end{minipage}\hfill
\begin{minipage}[c]{.49\linewidth}
	\includegraphics[width=\linewidth,trim = .5cm .3cm .6cm .1cm,clip]{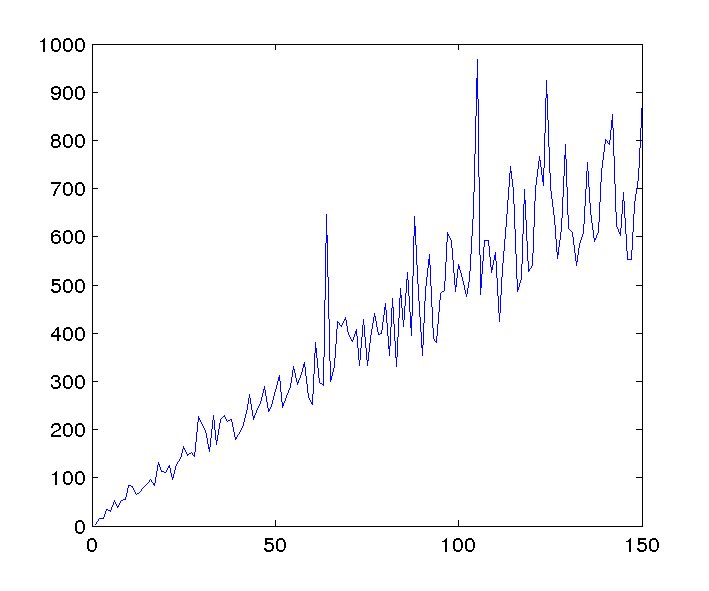}
\end{minipage}
}}
\end{center}
\caption[Simulations of the combinatorial behaviour of 2 examples of dissipative diffeomorphisms]{Size of the recurrent set $\Omega((f_i)_N)$ (top), number of periodic orbits of $(f_i)_N$ (middle) and length of the largest periodic orbit of $(f_i)_N$ (bottom) depending on $N$, for $f_1$ (left) and $f_2$ (right), on the grids $E_N$ with $N=128k$, $k=1,\cdots,150$.}\label{GrafDissip}
\end{figure}

We simulated some quantities related to the combinatorial behaviour of discretizations of homeomorphisms, namely:
\begin{itemize}
\item the cardinality of the recurrent set $\Omega(f_N)$,
\item the number of periodic orbits of $f_N$,
\item the maximal size of a periodic orbit of $f_N$.
\end{itemize}
We calculated these quantities for discretizations of orders $128k$ for $k$ from $1$ to $150$, and represented it graphically (see Figure \ref{GrafDissip}). For information, if $N =128\times150$, then $q_N \simeq 3.6. 10^8$.

Theoretically, the degree of recurrence, that is, the ratio between the cardinality of $\Omega(f_N)$ and $q_N$, should tend to $0$ (see Corollary~\ref{totsingdiscr}); this is what we observe on simulations. This is not really surprising: we will even see it for discretizations of conservative homeomorphisms (see Figure~\ref{GrafCons}). In this context, it is interesting to compare the behaviour of $\Omega(f_N)$ in the conservative and the dissipative case. The result is a little disappointing: the graphic for $f_1$, the dissipative homeomorphism, is very similar to that of $f_3$, the corresponding conservative homeomorphism, while in theory they should be very different. This is quite different for $f_2$, where the cardinality of $\Omega((f_2)_N)$ is more or less linear in $N$. We have no explanation to the linear shape of this function; if the maps $f_N$ were typical random maps, their degree of recurrence would be linear in $N$, with a value would close to $2.4.10^4$ for $N=150\times 128$ (here the value is about three times bigger).
\bigskip

The theoretical results assert that the number of periodic orbits of $f_N$ should tend to $+\infty$ (as a generic dissipative homeomorphism has an infinite number of attractors). We can hope that this quantity reflects the fact that the dynamics converges to that of the initial homeomorphism: among others, we can test if it is of the same order as the number of attractors of the homeomorphism. In practice, this number of periodic orbits of $f_N$ first increases rather quickly, to stabilize to around a value of $1.5.10^4$ for $f_1$ and $9.10^3$ for $f_2$. We could be tempted to interpret this phenomenon by the fact that after a while, the discretization has detected all the attractors of $f$ and thus, the number of periodic of the discretizations reflects the number of attractors of $f$. This idea may be reasonable for $f_2$ (as we will see in observing the invariant measures of $(f_2)_N$ in Figure~\ref{MesC1IdDissipSer}), but if we compare these graphs in the dissipative case with that of the conservative case (Figure~\ref{GrafCons}), we see that they are as alike as two peas in a pod. Thus, this is not clear at all that this behaviour is due to the dissipative character of the homeomorphism or not.
\bigskip

Since the dynamics of discretizations is assumed to converge to that of the initial homeomorphism, we could expect that the length of the longest periodic orbit of discretizations $(f_i)_N$ is almost always a multiple of that of an attractive periodic orbit of $f_i$. The graphic of this length for $f_1$ looks like the conservative case (see Figure~\ref{GrafCons}), so we can say that the dissipative behaviour of this homeomorphism is not detected in practical by this quantity. For $f_2$, the length of the longest orbit is much smaller than for $f_1$ (up to a factor 10), and seem to increase linearly in $N$. This may be imputed to the fact that $f_2$ is ``almost conservative'' around its attractive periodic points; thus it has a conservative behaviour, but at a smaller scale than $f_1$.

\subsection{Behaviour of invariant measures}

We have computed the invariant measures $\mu^{(f_i)_N}_{\T^2}$ of dissipative homeomorphisms $f_1$ and $f_2$ as defined on page \pageref{defbis}. Our aim is to test whether Theorem \ref{convmesdissip} applies in practice or if there are technical constraints such that this behaviour cannot be observed on these examples. For a presentation of the representations of the measures, see page~\pageref{pagealgo}.

\begin{figure}[ht]
\begin{center}
\makebox[0.8\textwidth]{\parbox{0.8\textwidth}{%
\begin{minipage}[c]{.49\linewidth}
	\includegraphics[width=\linewidth,trim = .5cm .3cm .6cm .1cm,clip]{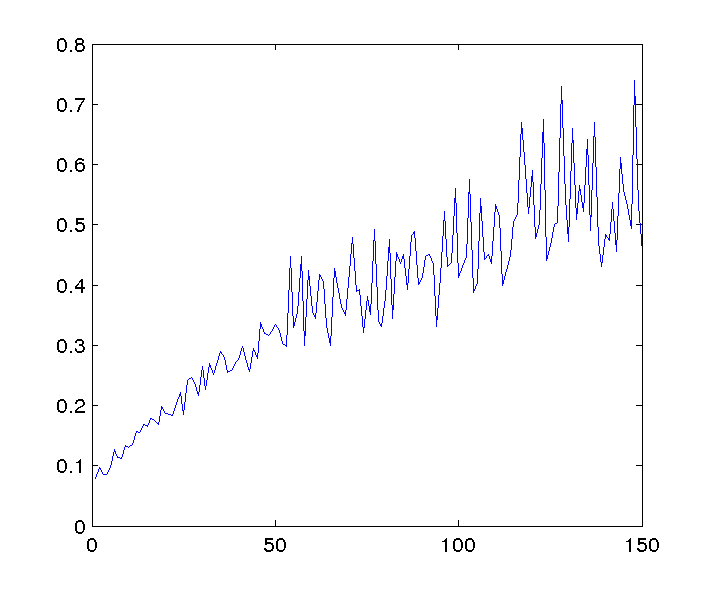}
\end{minipage}\hfill
\begin{minipage}[c]{.49\linewidth}
	\includegraphics[width=\linewidth,trim = .5cm .3cm .6cm .1cm,clip]{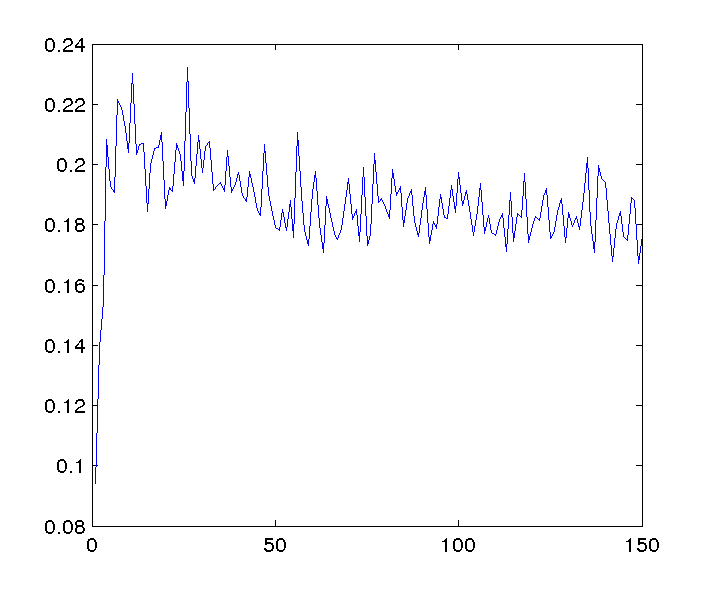}
\end{minipage}
}}
\end{center}
\caption[Distance between $\mu^{f_N}_{\T^2}$ and $\Leb$ for 2 examples of dissipative diffeomorphisms]{Distance between Lebesgue measure and the measure $\mu^{f_N}_{\T^2}$ depending on $N$ for $f_1$ (left) and $f_2$ (right), on the grids $E_N$ with $N=128k$, $k=1,\cdots,150$.}\label{GrafMaxmeasDissip}
\end{figure}

On Figure \ref{GrafMaxmeasDissip}, we have represented the distance between the measure $\mu^{(f_i)_N}_{\T^2}$ and Lebesgue measure. Theorem \ref{convmesdissip} says that for a generic dissipative homeomorphism, this quantity converges to the distance between $\mu_{\T^2}^f$ and Lebesgue measure. Clearly, this is not what happens in practice for $f_1$: the distance between both measures globally increases when $N$ increases. Locally, the behaviour of the map $N\mapsto \dist(\mu^{f_N}_{\T^2},\Leb)$ is quite erratic: there is no sign of convergence to any measure. When we compare this with the conservative case (Figure~\ref{GrafDistLebCons}), we see that both behaviours are very similar. In other words, we do not see the dissipative nature of $f_1$ on simulations. The behaviour of the distance between the measure $\mu^{(f_2)_N}_{\T^2}$ and Lebesgue measure is much more interesting. First of all, we observe that this distance is smaller than for $f_1$. Moreover, it is globally slightly decreasing in $N$, and it seem to converge to a value close to $0.2$. So it suggests that the measures $\mu^{(f_2)_N}_{\T^2}$ converge to a given measure whose distance to $\Leb$ is close to $0.2$, as predicted by the theory. In the view of these graphics, we can say that the behaviour of $f_2$ is much more close to that of a generic dissipative homeomorphism that that of $f_1$. This can be interpreted in the view of the article of J.-M. Gambaudo and C. Tresser \cite{MR700317}: for $f_1$, the attractors are much too small to be observed in practice for reasonable orders of discretization.

\begin{figure}[ht]
\begin{minipage}[c]{.31\linewidth}
	\includegraphics[height=4.8cm,trim = 1.5cm .95cm 2.8cm .5cm,clip]{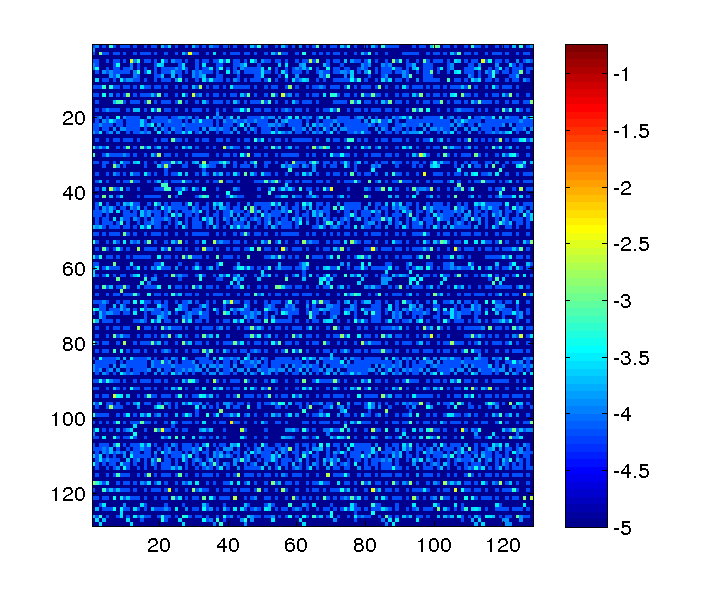}
\end{minipage}\hfill
\begin{minipage}[c]{.31\linewidth}
	\includegraphics[height=4.8cm,trim = 1.5cm .95cm 2.8cm .5cm,clip]{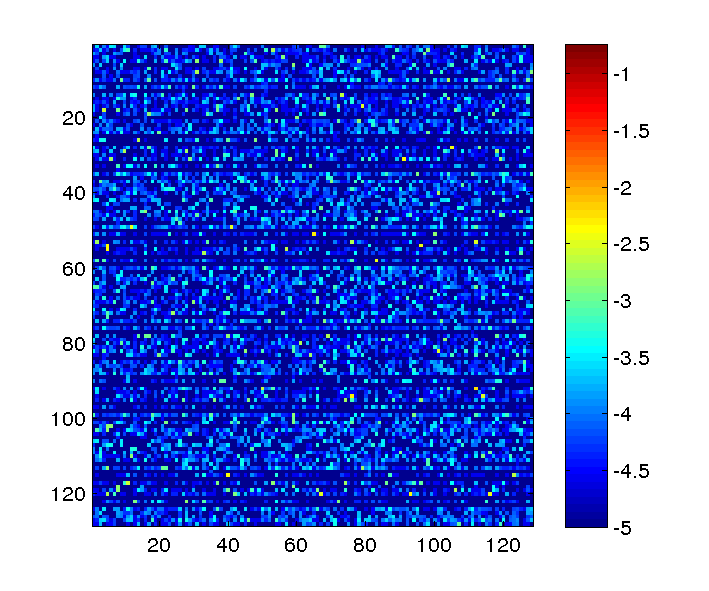}
\end{minipage}\hfill
\begin{minipage}[c]{.37\linewidth}
	\includegraphics[height=4.8cm,trim = 1.5cm .95cm 1cm .5cm,clip]{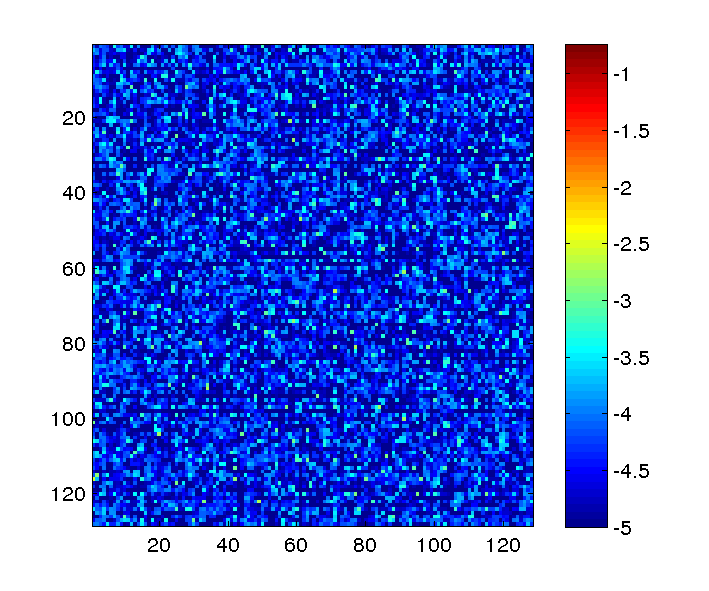}
\end{minipage}

\begin{minipage}[c]{.31\linewidth}
	\includegraphics[height=4.8cm,trim = 1.5cm .95cm 2.8cm .5cm,clip]{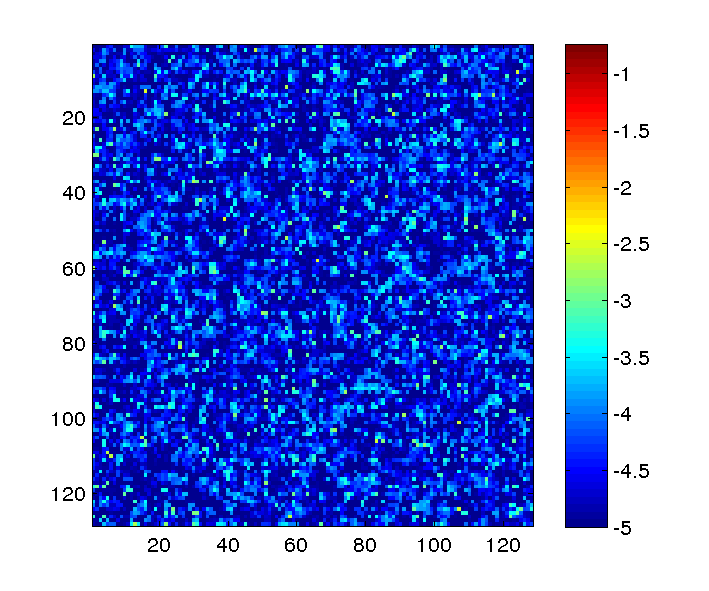}
\end{minipage}\hfill
\begin{minipage}[c]{.31\linewidth}
	\includegraphics[height=4.8cm,trim = 1.5cm .95cm 2.8cm .5cm,clip]{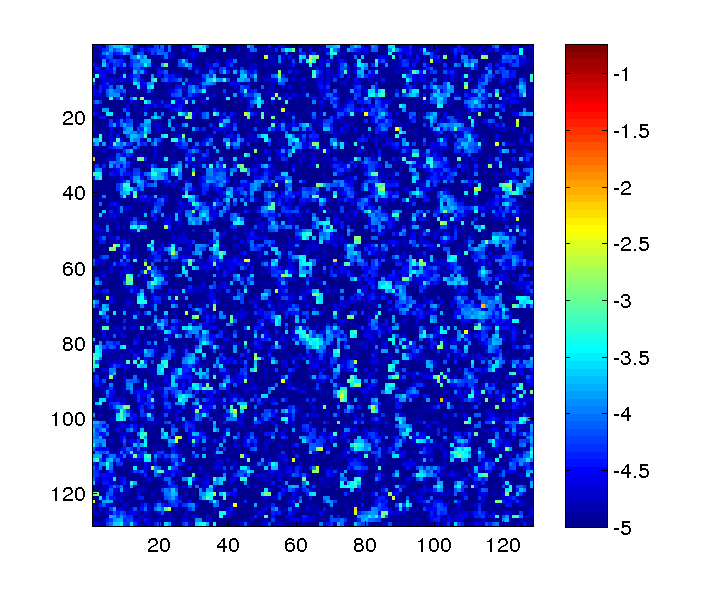}
\end{minipage}\hfill
\begin{minipage}[c]{.37\linewidth}
	\includegraphics[height=4.8cm,trim = 1.5cm .95cm 1cm .5cm,clip]{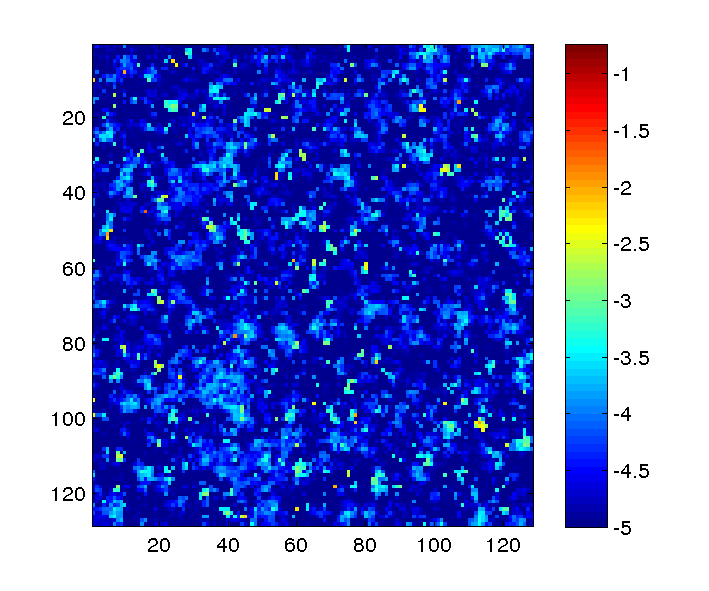}
\end{minipage}

\begin{minipage}[c]{.31\linewidth}
	\includegraphics[height=4.8cm,trim = 1.5cm .95cm 2.8cm .5cm,clip]{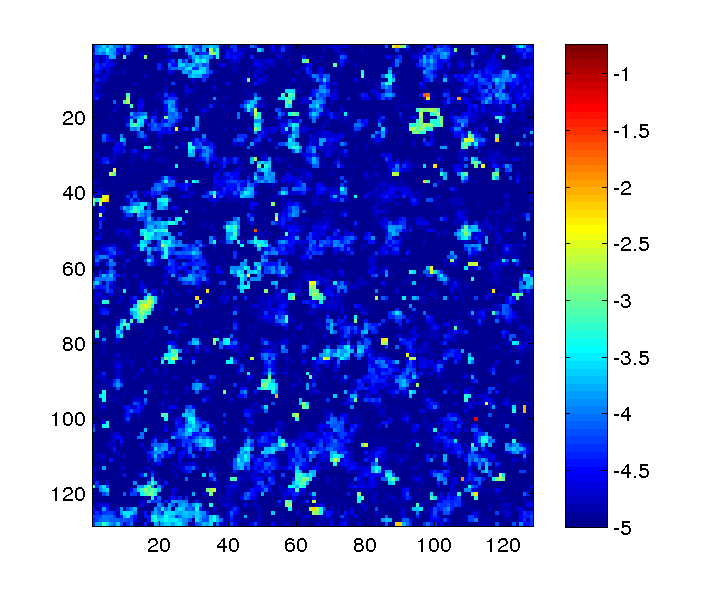}
\end{minipage}\hfill
\begin{minipage}[c]{.31\linewidth}
	\includegraphics[height=4.8cm,trim = 1.5cm .95cm 2.8cm .5cm,clip]{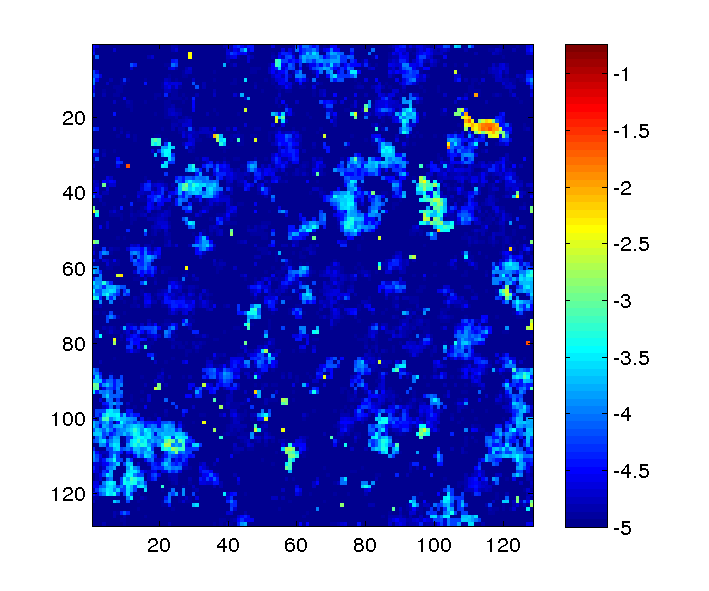}
\end{minipage}\hfill
\begin{minipage}[c]{.37\linewidth}
	\includegraphics[height=4.8cm,trim = 1.5cm .95cm 1cm .5cm,clip]{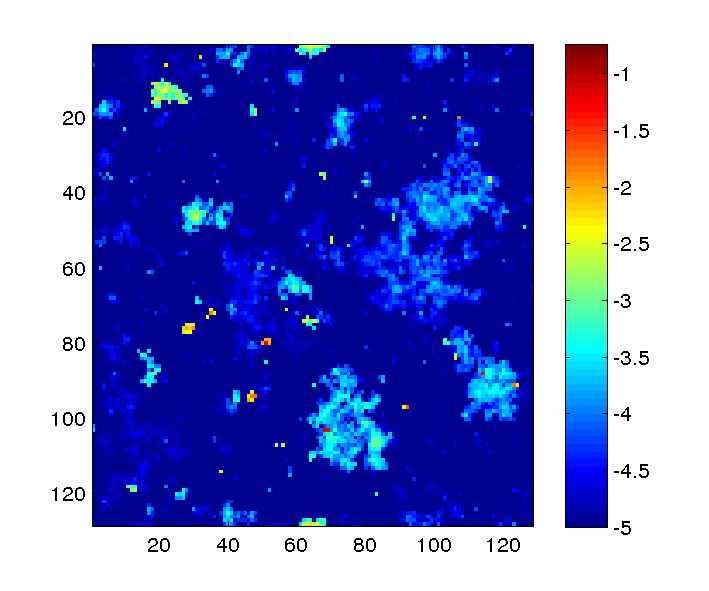}
\end{minipage}
\caption[Simulations of $\mu^{(f_1)_N}_{\T^2}$ on the grids $E_N$, with $N=2^k$, $k= 7,\cdots,15$]{Simulations of $\mu^{(f_1)_N}_{\T^2}$ on grids $E_N$, with $N=2^k$, $k= 7,\cdots,15$ (from left to right and top to bottom).}\label{MesC0IdDissipSer}
\end{figure}

The behaviour of invariant measures for $f_1$ (see Figure~\ref{MesC0IdDissipSer}), which is a small $C^0$ dissipative perturbation of identity, is relatively similar to that of invariant measures for $f_3$ \emph{i.e.} the corresponding conservative case (see Figure~\ref{MesC0IdCons2p}): when the order discretization is large enough, there is a strong variation of the measure $\mu^{(f_1)_N}_{\T^2}$. Moreover, this measure has a significant absolutely continuous component with respect to Lebesgue measure. This is very different from what is expected from the theoretical results (in particular Theorem~\ref{convmesdissip}), which say that for a generic dissipative homeomorphism, the measures $\mu^{f_N}_{\T^2}$ must converge to the measure $\mu_{\T^2}^f$. Thus, we can say that it is impossible to detect the dissipative character of $f_1$ on these simulations. Again, as noted by J.-M. Gambaudo and C. Tresser in \cite{MR700317}, the size of the attractors of $f_1$ can be very small compared to the numbers involved in the definition of $f_1$. So even in orders discretization such as $2^{15}$, the dissipative nature of the homeomorphism is undetectable on discretizations. This is why the discretizations of $f_1$ are very similar to those of its conservative counterpart $f_3$.

\begin{figure}[ht]
\begin{minipage}[c]{.31\linewidth}
	\includegraphics[height=4.8cm,trim = 1.5cm .95cm 2.8cm .5cm,clip]{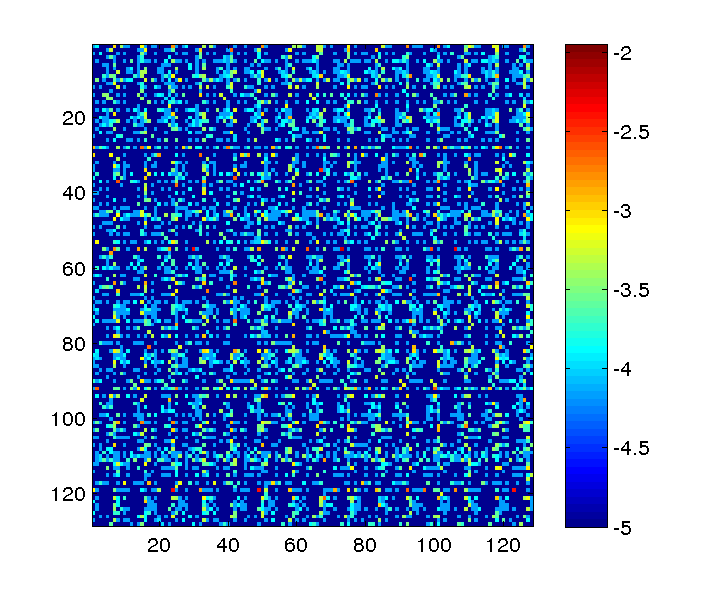}
\end{minipage}\hfill
\begin{minipage}[c]{.31\linewidth}
	\includegraphics[height=4.8cm,trim = 1.5cm .95cm 2.8cm .5cm,clip]{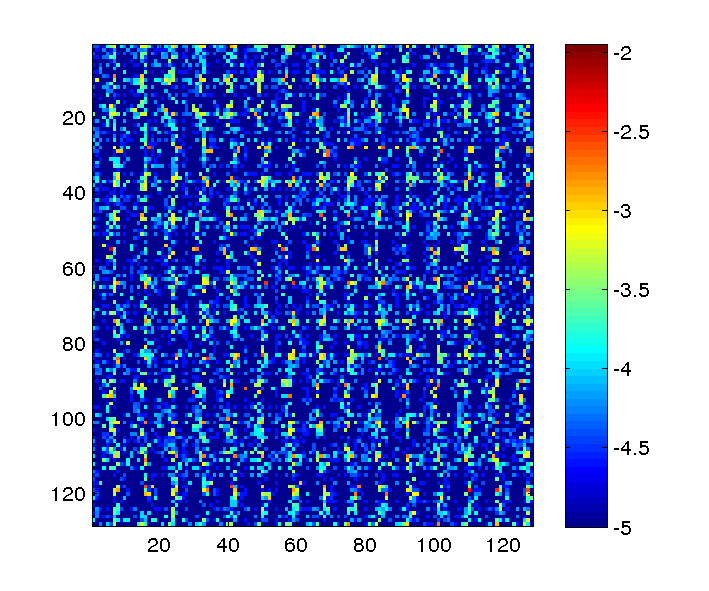}
\end{minipage}\hfill
\begin{minipage}[c]{.37\linewidth}
	\includegraphics[height=4.8cm,trim = 1.5cm .95cm 1cm .5cm,clip]{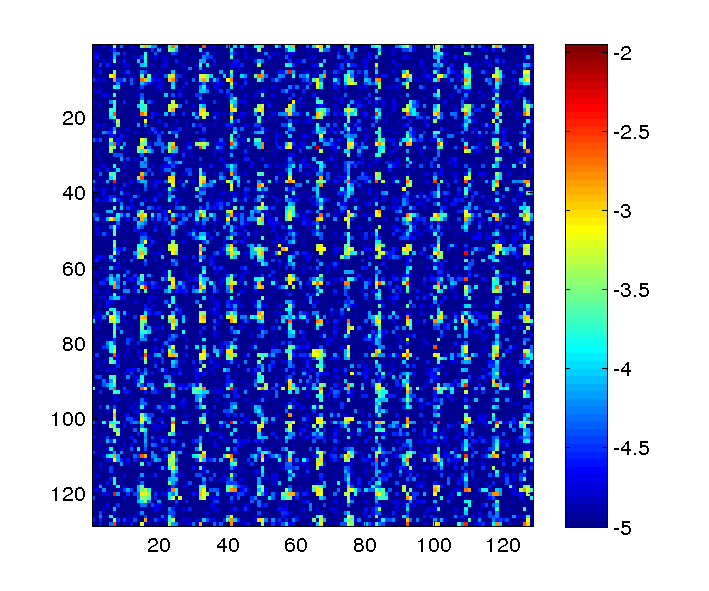}
\end{minipage}

\begin{minipage}[c]{.31\linewidth}
	\includegraphics[height=4.8cm,trim = 1.5cm .95cm 2.8cm .5cm,clip]{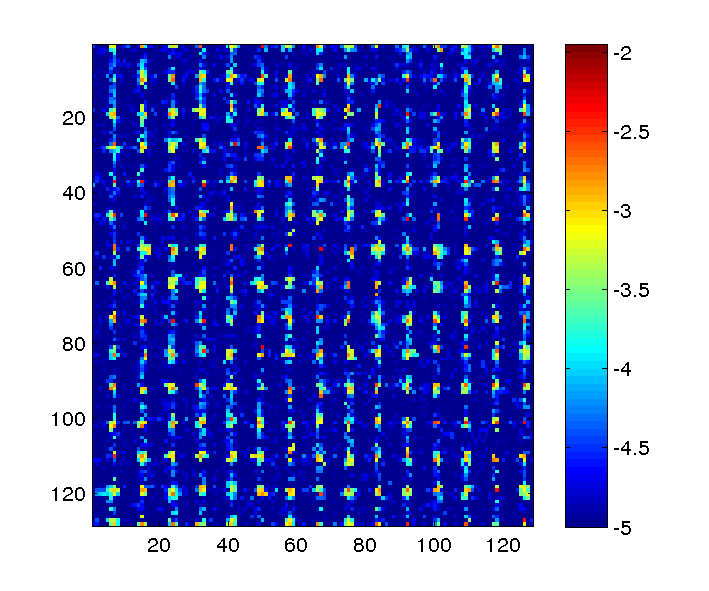}
\end{minipage}\hfill
\begin{minipage}[c]{.31\linewidth}
	\includegraphics[height=4.8cm,trim = 1.5cm .95cm 2.8cm .5cm,clip]{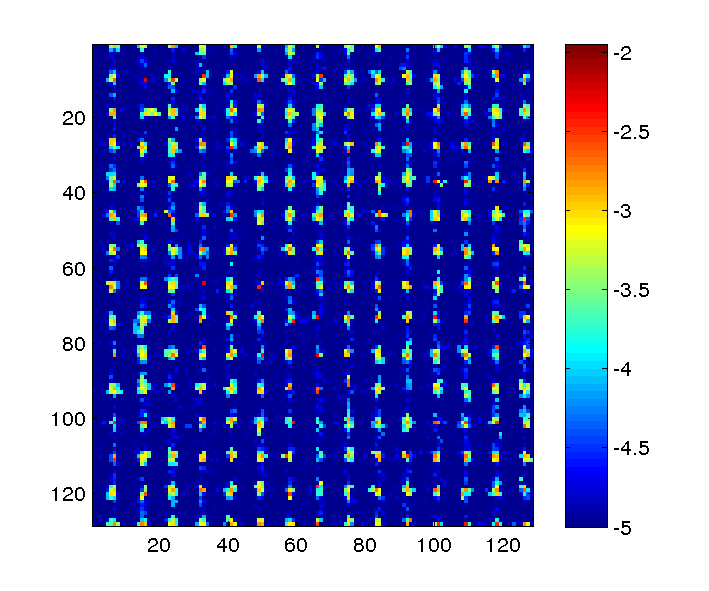}
\end{minipage}\hfill
\begin{minipage}[c]{.37\linewidth}
	\includegraphics[height=4.8cm,trim = 1.5cm .95cm 1cm .5cm,clip]{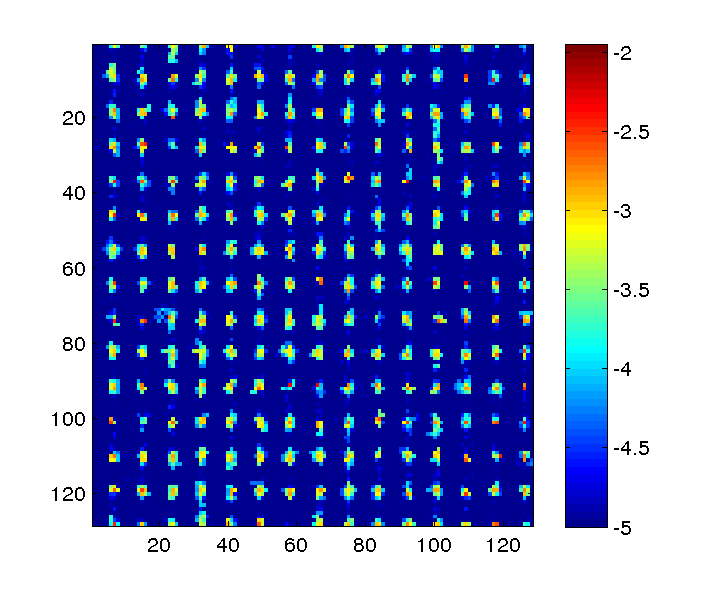}
\end{minipage}

\begin{minipage}[c]{.31\linewidth}
	\includegraphics[height=4.8cm,trim = 1.5cm .95cm 2.8cm .5cm,clip]{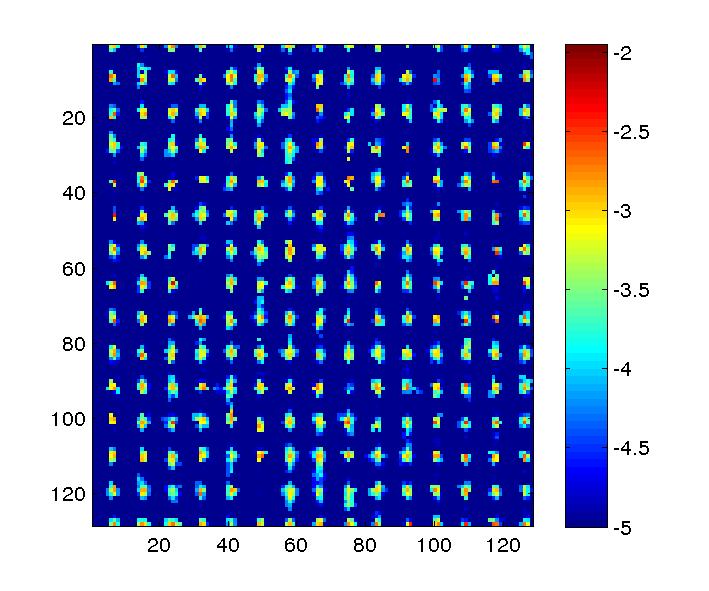}
\end{minipage}\hfill
\begin{minipage}[c]{.31\linewidth}
	\includegraphics[height=4.8cm,trim = 1.5cm .95cm 2.8cm .5cm,clip]{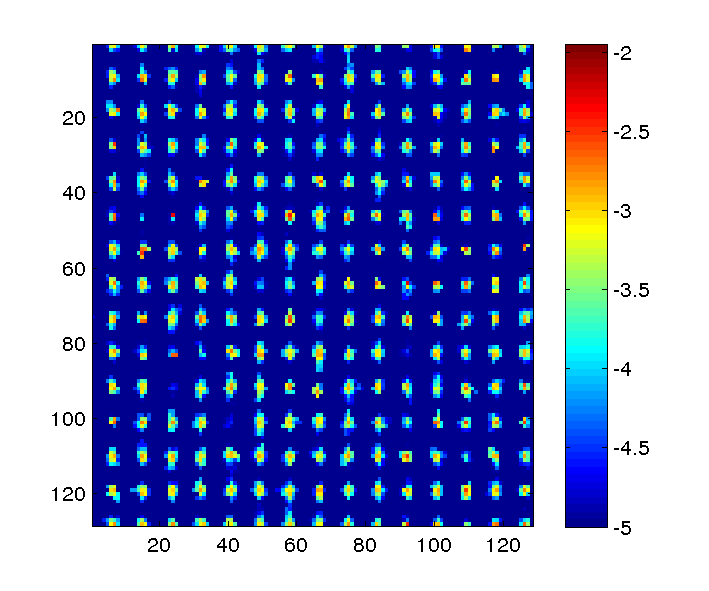}
\end{minipage}\hfill
\begin{minipage}[c]{.37\linewidth}
	\includegraphics[height=4.8cm,trim = 1.5cm .95cm 1cm .5cm,clip]{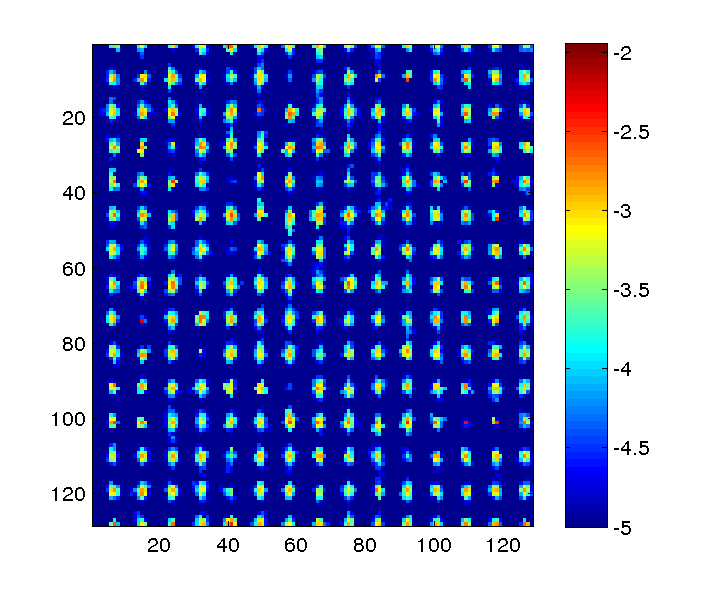}
\end{minipage}
\caption[Simulations of $\mu^{(f_2)_N}_{\T^2}$ on the grids $E_N$, with $N=20\,000,\cdots,20\,008$]{Simulations of $\mu^{(f_2)_N}_{\T^2}$ on grids $E_N$, with $N=2^k$, $k= 7,\cdots,15$  (from left to right and top to bottom).}\label{MesC1IdDissipSer}
\end{figure}

\begin{figure}[ht]
\begin{minipage}[c]{.31\linewidth}
	\includegraphics[height=4.8cm,trim = 1.5cm .95cm 2.8cm .5cm,clip]{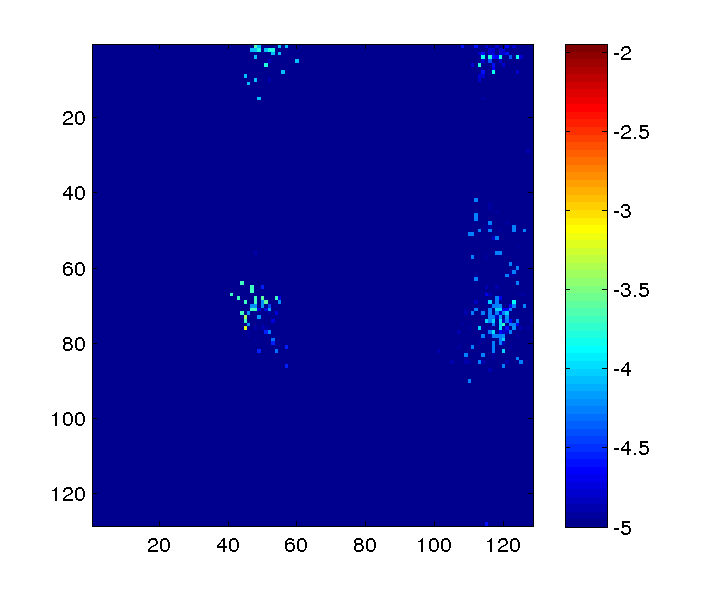}
\end{minipage}\hfill
\begin{minipage}[c]{.31\linewidth}
	\includegraphics[height=4.8cm,trim = 1.5cm .95cm 2.8cm .5cm,clip]{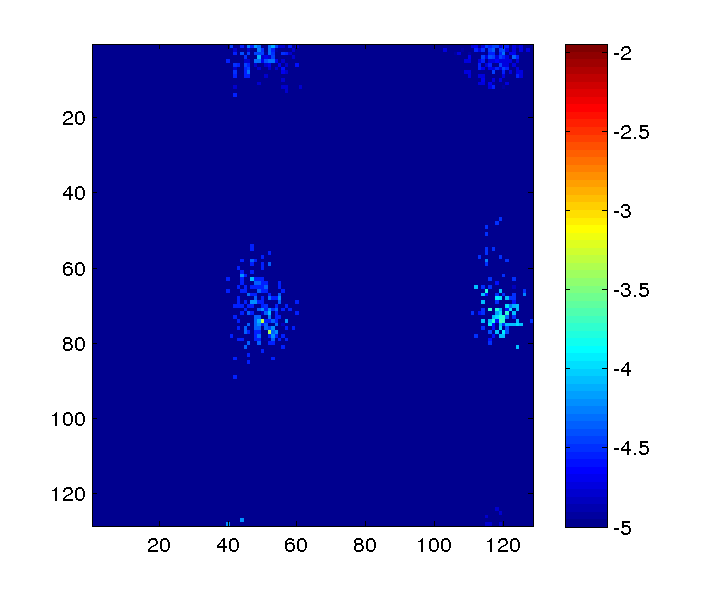}
\end{minipage}\hfill
\begin{minipage}[c]{.37\linewidth}
	\includegraphics[height=4.8cm,trim = 1.5cm .95cm 1cm .5cm,clip]{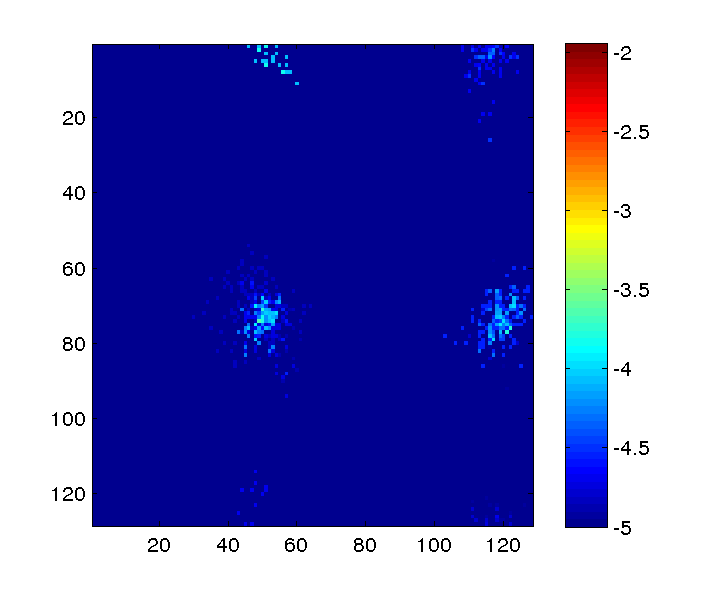}
\end{minipage}
\caption[Zoom on $\mu^{(f_2)_N}_{\T^2}$ on the grids $E_N$, with $N=2^k$, $k= 13,14,15$]{Zoom on the density of the measures $\mu^{(f_2)_N}_{\T^2}$ on grids $E_N$, with $N=2^k$, $k= 13,14,15$ (from left to right); the zoom is made on the top right of the representation of Figure~\ref{MesC1IdDissipSer}, on a square of size $1/8 \times 1/8$.}\label{MesC1IdDissipSerZoom}
\end{figure}

Recall that what happens for the dissipative homeomorphism $f_1$ is rather close to what happens for the conservative homeomorphism $f_3$. For their part, simulations of invariant measures for $f_2$ on grids of size $2^k \times 2^k$ (see Figure~\ref{MesC1IdDissipSer}) highlight that we expect from a generic dissipative homeomorphism: the measures $\mu^{(f_2)_N}_{\T^2}$ seem to tend to a single measure (that is also observed on a series of simulations), which is carried by the attractors of $f_2$. The fact that $f_2$ has few attractors allows the discretizations of reasonable orders (typically $2^{11}$) to find the actual attractors of the initial homeomorphism, contrary to what we observed for $f_1$. We also present a zoom of the density of the measures $\mu^{(f_2)_N}_{\T^2}$ (Figure~\ref{MesC1IdDissipSerZoom}). On these simulations, we can see that the attracting regions that we observe on the simulations of Figure~\ref{MesC1IdDissipSer} are in fact crumbled, it particular they are not connected. This is what is predicted by the theory: the set of Lyapunov periodic points of a generic dissipative homeomorphism is a Cantor set. On these zoomed simulations, we also observe that the dynamics of the discretizations has not completely converged at the order $N = 2^{15}$: locally, the measures $\mu^{(f_2)_N}_{\T^2}$ are quite different for different orders $N$; locally, the homeomorphism still behaves like in the conservative one.

\chapter{Discretizations of a generic conservative homeomorphism}\label{ChapCons}

We now study the conservative case, \emph{i.e.} we study the dynamical behaviour of discretizations of homeomorphisms of $X$ which are generic among those which preserve a given good measure $\lambda$. In the previous chapter, we proved that the dynamics of the discretizations of a generic \emph{dissipative} homeomorphism $f$ converges to the dynamics of $f$. This was due to the fact that such a generic homeomorphism possesses attractors whose basins cover almost all the phase space $X$. Now, the hypothesis of preservation of a measure $\lambda$ precludes this kind of behaviours, as conservative maps do not have attractors and thus are chain transitive.

We will begin by introducing the concept of \emph{dense type of discrete approximation} (Definition \ref{DensTyp}), which expresses that a given set of  finite maps on the grids can approach any conservative homeomorphism. Then, using the proposition of finite maps extension (Proposition \ref{extension}), we prove that any dense type of approximation occurs infinitely many times in the discretizations of a generic homeomorphism (Theorem \ref{génécycl}). Those dense types of approximations will be obtained by using Lax's theorem (Theorem \ref{Lax}). This strategy of proof will lead us to many results concerning the dynamical properties of discretizations of generic conservative homeomorphisms. In particular, for a generic conservative homeomorphism $f$, under the assumptions that the grids are well distributed and well ordered:
\begin{itemize}
\item there are infinitely many discretizations $f_N$ such that $f_N$ is a cyclic permutation (Corollary \ref{typlax});
\item for every $\varep>0$, there are infinitely many discretizations $f_N$ such that $f_N$ is $\varep$-topologically weakly mixing (Corollary \ref{méldiscr}, see Definition \ref{epmélfaibl});
\item for every period $p$ of a periodic orbit of $f$, there are infinitely many discretizations $f_N$ such that $f_N$ has a unique injective orbit whose corresponding periodic orbit has length $p$ (Corollary \ref{corovar2});\newpage
\item for every map $\vartheta : \N\to\R_+^*$ which tends to $+\infty$ at $+\infty$, there are infinitely many discretizations $f_N$ such that $\card(f_{N}(E_{N})) \le \vartheta(N)$ (Corollary \ref{crush});
\item if we further assume that the grids are self similar, for every map $\vartheta : \N\to\R$ such that $\vartheta(N) = o(q_N)$, there are infinitely many discretizations $f_N$ such that $f_{N}$ has at least $\vartheta(N)$ cycles which are pairwise conjugated (Corollary \ref{corovar3}).
\end{itemize}

Moreover, we prove that we can not hope to observe the actual dynamics of a generic homeomorphism by looking at the frequency with which some properties appear on discretizations. For example, if the sequence of grids is strongly self similar, we will prove that for a generic conservative homeomorphism, for every $\varep>0$, there are infinitely many $M\in\N$ such that the proportion of discretizations $f_N$, with $1\le N \le M$, such that $f_N$ is a cyclic permutation, is greater than $1-\varep$ (Theorem \ref{propdemin}). We also prove that under the assumption that the sequence of grids is self similar, the same property holds when we replace ``$f_N$ is a cyclic permutation'' by ``$\card(\Omega(f_N)) = o(q_N)$'' (Corollary \ref{CoroAver1}), ``$f_{N}$ is $\varep$-topologically weakly mixing'' (Corollary \ref{CoroAver2}), ``$\frac{\card(f_{N}(E_{N}))}{\card(E_{N})}<\varep$'' (Corollary \ref{petitepermieux}) or ``$f_N$ has at least $M$ periodic orbits'' (Corollary \ref{CoroAver3}).

However, we will prove that it is possible to recover the dynamics of the initial homeomorphism by looking at the dynamics of \emph{all} its discretizations. For example, it is possible to detect the periodic orbits of $f$ and their periods (Theorem \ref{corovar2}), to recover the set of invariant measures of $f$ (Theorem \ref{EnsMesInv}), or to recover the set of compact subsets of $X$ that are invariant under $f$ (Theorem~\ref{CompactInv}).

From this point of view of the shadowing property of the dynamics of the homeomorphism by that of all its discretizations, the ``physical dynamics'', that is the dynamical properties that depend on the good measure $\lambda$, plays non important role among all the dynamical properties of $f$. For instance, we will prove that if the sequence of grids is self similar, then for a generic homeomorphism $f$, the sequence of $f_N$-physical\footnote{Recall that $\mu^{f_N}_X$ is the limit in the sense of Cesàro of the pushforwards of the uniform measure on $E_N$ by the iterates of $f_N$} measures $(\mu^{f_N}_X)_{N\in\N}$ accumulates on the the whole set of $f$-invariant measures (Theorem \ref{mesinv}). And with the same techniques of proof we will show a similar result for invariant compact sets (Proposition \ref{EquivCompact}).

The end of this chapter is dedicated to the results of the numerical simulations.
\bigskip

\emph{Recall that we have fixed once and for all a manifold $X$, a good measure $\lambda$ on $X$ and a sequence of discretization grids $(E_N)_{N\in\N}$ (see Definition \ref{grillmiam}). We further assume that this sequence of grids is well distributed and well ordered (see Chapter \ref{ChapNota}). In this chapter, we will focus on discretizations of elements of $\Hom (X,\lambda)$, so homeomorphisms will always be supposed conservative.}

\section{Dense types of approximation}\label{Grosse}

To begin with, we define the notion of \emph{type of approximation}.

\begin{definition}\label{DensTyp}
A \emph{type of approximation} $\tT = (\tT_N)_{N\in\N}$ is a sequence of subsets of the set $\mathcal F (E_N,E_N)$ of applications from $E_N$ into itself.

Let $\mathcal{U}$ be an open subset of $\Hom (X,\lambda)$ (therefore $\mathcal U$ is a set of homeomorphisms). A type of approximation $\tT$ is said to be \emph{dense} in $\mathcal{U}$ if for all $f\in \mathcal{U}$, all $\varepsilon>0$ and all $N_0\in\N$, there exists $N\ge N_0$ and $\sigma_N\in \tT_N$ such that $d_N(f,\sigma_N)<\varepsilon$ (recall that $d_N$ is the distance between $f_{|E_N}$ and $\sigma_N$ considered as maps from $E_N$ into $X$).
\end{definition}

The goal of this paragraph is to obtain the following result: every dense type of approximation appears on infinitely many discretizations of a generic homeomorphism. This will lead us to a systematic study of dense types of approximation.

\begin{theoreme}\label{génécycl}
Let $\tT$ be a type of approximation which is dense in $\Hom(X,\lambda)$. Then for a generic homeomorphism $f\in\Hom(X,\lambda)$ and for all $N_0\in\N$, there exists $N\ge N_0$ such that $f_{N}\in \tT_{N}$.
\end{theoreme}

Thus, we have decomposed the obtaining of generic properties of discretizations into two steps.
\begin{enumerate}
\item The first step consists in proving that a given type of approximation is dense. This part will be done in Section \ref{partie 1.3}, the systematic use of Lax's theorem (see Section \ref{SecLax}) will reduce the proof of density of types of approximation to combinatorial problems.
\item The second step is always the same, it simply consists in applying Theorem~\ref{génécycl}.
\end{enumerate}

\begin{rem}
However, some properties can be proved in different ways: for instance Corollary \ref{crush} can be proved by inserting horseshoes in the dynamics of a given homeomorphism, using the local modification theorem (Theorem \ref{extension-sphères}, for a presentation of the technique in another context see \cite[Section 3.3]{MR2931648}); also a variation of Corollary \ref{corovar1} can be shown in perturbing any given homeomorphism such that it has a periodic orbit whose distance to the grid $E_N$ is smaller than the modulus of continuity of $f$, so that the actual periodic orbit and the discrete orbit fit together.
\end{rem}

To prove Theorem \ref{génécycl}, we start from a dense type of approximation -- in other words a sequence of discrete applications -- and we want to get properties of homeomorphisms. The tool that allows us to restore a homeomorphism from a finite map $\sigma_N : E_N\to E_N$ is the finite map extension proposition (Proposition \ref{extension}).

\begin{lemme}\label{lemmetrans}
Let $\mathcal U$ be an open subset\footnote{In most cases we will take $\mathcal U = \Hom(X,\lambda)$.} of $\Hom(X,\lambda)$ and $\tT=(\tT_N)_{N\in\N}$ be a type of approximation which is dense in $\mathcal U$. Then for all $N_0\in\N$, there exists an open and dense subset of $\mathcal U$ on which every homeomorphism $f$ satisfies: there exists $N\ge N_0$ such that the discretization $f_N$ belongs to $\tT_N$.
\end{lemme}

\begin{proof}[Proof of Lemma \ref{lemmetrans}] Let $f\in\mathcal U$, $N_0\in\N$ and $\varepsilon>0$. Since $\tT$ is a type of approximation which is dense in $\mathcal U$, there exists $N\ge N_0$ and $\sigma_{N}\in \tT_{N}$ such that $d_{N}(f,\sigma_{N})< \varepsilon$. Let $x_1,\cdots,x_{q_{N}}$ be the elements of $E_N$. Then for all $\ell$, $d(f(x_\ell),\sigma_{N}(x_\ell))<\varepsilon$. Now, the problem is that $\sigma_N$ is not injective. To solve this problem, we take advantage of the fact that $P_N$ is highly non injective to modify $\sigma_N$ into a injection $\sigma'_N : E_N\to X$ whose discretization on $E_N$ equals to $\sigma_N$: set $\sigma'_N(x_1) = \sigma_N(x_1)$, $\sigma'_N$ is defined by induction by choosing $\sigma'_N(x_\ell)$ such that $\sigma'_N(x_\ell)$ is different from $\sigma'_N(x_i)$ for $i\le\ell$, such that $P_N(\sigma'_N(x_\ell)) = P_N(\sigma_N(x_\ell))$ and that $d(f(x_\ell),\sigma'_N(x_\ell))<\varep$.

Since $\sigma'_N$ is a bijection, Proposition \ref{extension} can be applied to $f(x_\ell)$ and $\sigma'_{N}(x_\ell)$; this gives a measure-preserving homeomorphism $\varphi$ such that $\varphi(f(x_\ell)) = \sigma'_{N}(x_\ell)$ for all $\ell$ and such that $d(\varphi,\mathrm{Id})< \varepsilon$. Set $f' = \varphi\circ f$, we have $d(f,f')\leq \varepsilon$ and $f'(x_\ell) = \sigma'_{N}(x_\ell)$ for all $\ell$, and therefore $f_N = \sigma_N$.
\end{proof}

Now, proof of Theorem \ref{génécycl} is straightforward.

\begin{proof}[Proof of Theorem \ref{génécycl}]
Let $(x_{{N},\ell})_{1\le \ell\le q_{N}}$ be the elements of $E_{N}$ and consider the set (where $\rho_{N}$ is the minimal distance between two distinct points of $E_{N}$)
\[\bigcap_{N_0\in\N}\,\bigcup_{\substack{N\ge N_0\\\sigma_{N}\in \tT_{N}}}\,\left\{f\in\Hom(X,\lambda) \mid\forall \ell,\, d_{N}\big(f(x_{{N},\ell}),\sigma_{N}(x_{{N},\ell})\big)<\frac{\rho_{N}}{2}\right\}.\]
This set is clearly a $G_\delta$ set and its density follows directly from Lemma \ref{lemmetrans}. Moreover we can easily see that its elements satisfy the conclusions of the theorem.
\end{proof}

\section{Lax's theorem}\label{SecLax}

Now that we have shown Theorem \ref{génécycl}, it remains to obtain dense types of approximation. Systematic use will be made of a theorem due to Oxtoby and Ulam \cite{Oxto-meas} but classically referred as Lax's theorem \cite{MR0272983}, and more precisely its improvement as stated by S. Alpern in \cite{Alpe-newp}, which allows to approach any conservative homeomorphism by a cyclic permutation of the elements of a discretization grid.

\begin{theoreme}[Lax, Alpern]\label{Lax}
Let $f\in\mathrm{Homeo}(X,\lambda)$ and $\varepsilon>0$. Then there exists $N_0\in\N$ such that for all $N\ge N_0$, there exists a cyclic permutation $\sigma_N$ of $E_N$ such that $d_N(f,\sigma_N)<\varepsilon$.
\end{theoreme}

It is worthwhile to note that this theorem was initially stated in the view of numerical approximation of conservative homeomorphisms: P. Lax noticed that in some pioneering works showing numerical simulations of chaotic dynamical systems (e.g. \cite{MR0253513}), the discretizations of conservative maps were only approximately measure-preserving. His theorem states that among all possible approximations of a conservative homeomorphism, some of them are bijective, which is a discrete counterpart of measure preservation. This theorem was later improved by S. Alpern in \cite{Alpe-newp, Alpe-appr} to establish generic properties of conservative homeomorphisms (see also \cite{MR1307740,MR1453713} for a generalisation and some simulations in dimension 1). This theorem is now one of the keystones of the theory of generic properties of homeomorphisms (see \cite{MR2931648}), together with other theorems of approximation by permutations\footnote{See also \cite{MR0219697,Kato-metr} for approximations in weak topology 
or \cite[page 65]{MR0097489} for 
approximations of conservative automorphisms.}.

As a compact metric set can be seen as a ``finite set up to $\varep$'', Lax's theorem allows us to see a homeomorphism as a ``cyclic permutation up to $\varep$''. Genericity of transitivity in $\Hom(X,\lambda)$ follows easily from Lax's theorem together with the finite map extension proposition (see \cite{Alpe-typi} or part 2.4 of \cite{MR2931648}). In our case, applying Theorem \ref{génécycl}, we deduce that infinitely many discretizations of a generic homeomorphism are cyclic permutations. The purpose of Section \ref{partie 1.3} is to find variations of Lax's theorem (which are at the same time corollaries of it) concerning other properties of discretizations.

We give briefly the beautiful proof of Lax's theorem, which is essentially combinatorial and based on marriage lemma and on a lemma of approximation of permutations by cyclic permutations due to S. Alpern. Readers wishing to find proofs of these lemmas may consult \cite[Section 2.1]{MR2931648}. It is here that we use the fact that discretizations grids $E_N$ are well distributed and well ordered (see Section \ref{SecHyp}).

\begin{lemme}[Marriage lemma]\label{mariage}
Let $E$ and $F$ be two finite sets and $\approx$ a relation between elements of $E$ and $F$. Suppose that the number of any subset $E'\subset E$ is smaller than the number of elements in $F$ that are associated with an element of $E'$, i.e.:
\[\forall E'\subset E,\quad \card (E')\leq \card\left\{f\in F\mid \exists e \in E' : e\approx f \right\},\]
then there exists an injection $\Phi : E\to F$ such that for all $e\in E$, $e\approx \Phi(e)$.
\end{lemme}

\begin{lemme}[Cyclic approximations in $\Sn_q$, \cite{Alpe-newp}]\label{Pioure}
Let $q\in\N^*$ and $\sigma \in \mathfrak{S}_q$ ($\Sn_q$ is seen as the permutation group of $\Z/q\Z$). Then there exists $\tau \in \mathfrak{S}_q$ such that $|\tau(k)-k|\leq 2$ for all $k$ (where $|.|$ is the distance in $\Z/q\Z$) and such that the permutation $\tau\sigma$ is cyclic.
\end{lemme}

\begin{proof}[Proof of Theorem \ref{Lax}]
Let $f\in\mathrm{Homeo}(X,\lambda)$ and $\varep>0$. Consider $N_0\in\N$ such that for all $N\ge N_0$, the diameter of the cubes of order $N$ (given by the hypothesis ``$E_N$ is well distributed'') and their images by $f$ is smaller than $\varepsilon$. We define a relation $\approx$ between cubes of order $N-1$: $C\approx C'$ if and only if $f(C)\cap C'\neq\emptyset$. Since $f$ preserves $\lambda$, the image of the union of $\ell$ cubes intersects at least $\ell$ cubes, so the marriage lemma (Lemma \ref{mariage}) applies: there exists an injection $\Phi_N$  from the set of cubes of order $N$ into itself (then a bijection) such that for all cube $C$, $f(C)\cap\Phi_N(C)\neq\emptyset$. Let $\sigma_N$ be the application that maps the centre of any cube $C$ to the centre of the cube $\Phi_N(C)$, we obtain:
\[d_N(f,\sigma_N) \leq \sup_{C}\big(\diam(C)\big)+\sup_{C}\big(\diam(f(C))\big)\le 2\varepsilon.\]

It remains to show that $\sigma_N$ can be chosen as a cyclic permutation. Increasing $N$ if necessary, and using the hypothesis that the grids are well ordered, we number the cubes such that the diameter of the union of two consecutive cubes is smaller than $\varep$. Then we use Lemma \ref{Pioure} to obtain a cyclic permutation $\sigma_N'$ whose distance to $\sigma_N$ is smaller than $\varep$. Thus we have found $N_0\in\N$ and for all $N\ge N_0$ a cyclic permutation $\sigma_N'$ of $E_N$ whose distance to $f$ is smaller than $3\varepsilon$.
\end{proof}

\begin{rem}
In this proof, the marriage lemma can be replaced by Birkhoff's theorem about bistochastic matrices: we define the transition matrix $M$ for the partition of $X$ by the cubes by $m_{i,j} = \lambda\big(C_i\cap f(C_j)\big)/\lambda(C_i)$; the hypothesis of preservation of $\lambda$ implies that $M$ is bistochastic\footnote{That is, its entries are nonnegative and the sum of the elements of each of its columns and rows is equal to 1.}. Then, Birkhoff's theorem ensures that bistochastic matrices form a convex set whose extremal points are permutation matrices; thus $M$ can be written as a convex combination of permutation matrices and each one of the corresponding permutation provides an approximating permutation for $f$.
\end{rem}

\section{Individual behaviour of discretizations}\label{partie 1.3}

Now we have set the technical Theorems \ref{génécycl} and \ref{Lax}, we can establish results concerning the behaviour of discretizations of a generic conservative homeomorphism. Here we study \emph{individual} behaviour of discretizations, \emph{i.e.} properties about only one order of discretization. As has already been said, applying Theorem \ref{génécycl}, it suffices to find dense types of approximation to obtain properties about discretizations. In practice, these dense types of approximations are obtained from variations of Lax's theorem (Theorem \ref{Lax}).

Recall that the sequence $(E_N)_{n\in\N}$ of discretization grids is well distributed and well ordered (see Definition \ref{Ashe}), we denote by $f_N$ the discretization of a homeomorphism $f$ and by $\Omega(f_N) $ the recurrent set of $f_N$ (\emph{i.e.} the union of periodic orbits of $f_N$).

We will show that for a lot of simple dynamical properties $(P)$ about finite maps and for a generic conservative homeomorphism $f$, infinitely many discretizations $f_N$ satisfy $(P)$ as well as infinitely many discretizations satisfy its contrary. For instance, for a generic homeomorphism $f$, the recurrent set $\Omega(f_N)$ is sometimes as large as possible, \emph{i.e.} $\Omega(f_N) = E_N$ (Corollary \ref{typlax}), sometimes very small (Corollary \ref{corovar1}) and even better sometimes the number of elements of the image of $E_N$ is small (Corollary \ref{crush}). In the same way stabilization time\footnote{That is, the smallest integer $k$ such that $f_N^k(E_N) = \Omega(f_N)$.} is sometimes zero (Corollary \ref{typlax}, for example), sometimes around $\card (E_N)$ (Corollary \ref{corovar2}). Finally, concerning the dynamics of ${f_N}_{|\Omega(f_N)}$, sometimes it is a cyclic permutation (Corollary \ref{typlax}) or a bicyclic permutation (Corollary \ref{méldiscr}), sometimes it has many orbits (Corollary \ref{corovar3})\dots
\bigskip

Firstly, we deduce directly from Lax's theorem that cyclic permutations of the sets $E_N$ form a dense type of approximation in $\Hom (X,\lambda)$. Combining this with Theorem \ref{génécycl}, we obtain directly:

\begin{coro}[Miernowski, \cite{Mier-dyna}]\label{typlax}
For a generic homeomorphism $f\in\Hom(X,\lambda)$, for every $N_0\in\N$, there exists $N\ge N_0$ such that $f_{N}$ is a cyclic permutation\footnote{In fact, T. Miernowski proves ``permutation'' but his arguments, combined with S. Alpern's improvement of the Lax's theorem, show ``cyclic permutation''.}.
\end{coro}

This theorem states that for every generic conservative homeomorphism, there exists a subsequence of discretizations which are ``transitive''. Recall
that generic homeomorphism are transitive (see \cite{Oxto-Note}). So, in a certain sense, transitivity can be detected on discretizations. Remark that this result implies that the discretizations of a generic conservative homeomorphism do not behave like typical random maps, as for a random map of a set with $q$ elements the average number of periodic orbits is asymptotically $\log q$ (see for example \cite[XIV.5]{Boll-rand}).
\bigskip

There exists a variation of Lax's Theorem for bicyclic permutations, which are permutations having exactly two orbits whose lengths are relatively prime (see \cite[lemme 2.9 page 28]{MR2931648}); this variation leads to a proof of genericity in $\Hom(X,\lambda)$ of topological weak mixing (see for example \cite{Alpe-typi} or  \cite[Part 2.4]{MR2931648}). In the following we define a discrete analogue to topological weak mixing and state that this property occurs infinitely often on discretizations of a generic homeomorphism.

\begin{definition}
A homeomorphism $f$ is said to be \emph{topologically weakly mixing} if for all non-empty open sets $(U_i)_{i\le M}$ and $(U'_i)_{i\le M}$, there exists $m\in\N$ such that $f^m(U_i)\cap U'_i$ is non-empty for all $i\le M$.
\end{definition}

The proof of the genericity of topological weak mixing starts by an approximation of every conservative homeomorphism by another having $\varep$-dense periodic orbits whose lengths are relatively prime. The end of the proof primarily involves the use of Baire's theorem and Bézout's identity. In the discrete case, the notion of weak mixing is replaced by the following.

\begin{definition}\label{epmélfaibl}
Let $\varep>0$. A finite map $\sigma_N$ is said \emph{$\varep$-topologically weakly mixing} if for all $M\in\N$ and all balls $(B_i)_{i\le M}$ and $(B'_i)_{i\le M}$ with diameter $\varep$, there exists $m\in\N$ such that for all $i$
\[\sigma_N^m(B_i\cap E_{N})\cap (B'_i\cap E_{N})\neq\emptyset.\]
\end{definition}

The first step of the proof is replaced by the following variation of Lax's theorem:

\begin{prop}[First variation of Lax's theorem]\label{melfaibl}
Let $f\in\Hom(X,\lambda)$ be a homeomorphism such that all the iterates of $f$ are topologically transitive. Then for all $\varep>0$ and all $M\in\N^*$, there exists $N_0\in\N$ such that for all $N\ge N_0$, there exists $\sigma_{N} : E_{N}\to E_{N}$ which has $M$ $\varep$-dense periodic orbits whose lengths are pairwise relatively prime, and such that $d_{N}(f,\sigma_{N})<\varep$.
\end{prop}

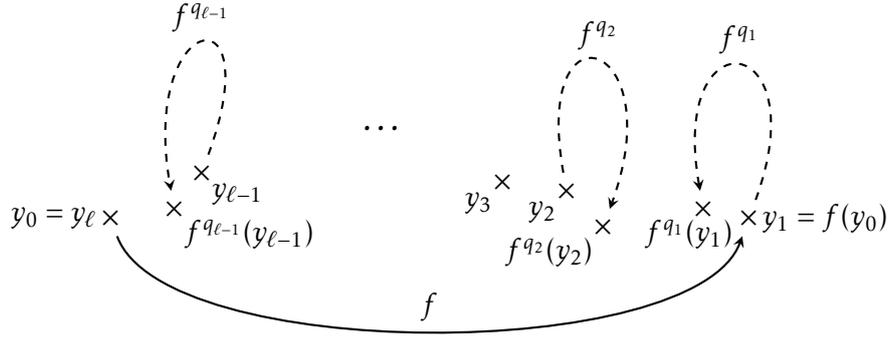
\begin{figure}
\begin{center}
\begin{tikzpicture}[scale=1.2]
\node (A) at (0,0) {\large$\times$};\node at (A) [left]{$y_0=y_\ell\,$};
\node (B) at (7,0) {\large$\times$};\node at (B) [right]{$\,y_1 = f(y_0)$};
\node (C) at (5,.3) {\large$\times$};\node at (C) [below left]{$y_2$};
\node (D) at (4.3,.4) {\large$\times$};\node at (D) [below left]{$y_3$};
\node (F) at (1,.5) {\large$\times$};\node at (F) [below right]{$y_{\ell-1}$};

\node (b) at (6.5,.1){\large$\times$};\node at (b)[below]{$f^{q_1}(y_1)\quad$};
\node (c) at (5.4,-.1){\large$\times$};\node at (c)[below left]{$f^{q_2}(y_2)$};
\node (f) at (.7,.1){\large$\times$};\node at (f)[below right]{$f^{q_{\ell-1}}(y_{\ell-1})$};
\node at (3,1){\Large$\dots$};

\draw[->,>=stealth,thick] (A) to [out=-70,in=-110,distance=1.5cm] node[above]{$f$}(B);
\draw[->,>=stealth,thick,dashed] (B) to [out=70,in=100,distance=2cm] node[above]{$f^{q_1}$}(b);
\draw[->,>=stealth,thick,dashed] (C) to [out=100,in=70,distance=2cm] node[above]{$f^{q_2}$}(c);
\draw[->,>=stealth,thick,dashed] (F) to [out=70,in=100,distance=2cm] node[above]{$f^{q_{\ell-1}}$}(f);

\end{tikzpicture}
\caption{Construction of the sequence $(y_m)_{1\le m\le \ell}$}\label{constrigrec}
\end{center}
\end{figure}

\begin{proof}[Proof of Proposition \ref{melfaibl}]
We prove the proposition in the case where $M=2$, the other cases being easily obtained by an induction. Let $\varep>0$ and $f$ be a homeomorphism whose all iterates are topologically transitive. Then there exists $x_0\in X$ and $p\in\N^*$ such that $\{x_0,\cdots,f^{p-1}(x_0)\}$ is $\varep$-dense and $d(x_0,f^p(x_0))<\varep/2$. Since transitive points of $f^p$ form a dense $G_\delta$ subset of $X$, while the orbit of $x_0$ form a $F_\sigma$ set with empty interior, the set of points whose orbit under $f^p$ is dense and disjoint from that of $x_0$ is dense. So we can pick such a transitive point $y_0$. Set $y_1 = f(y_0)$. Then there exists a multiple $q_1$ of $p$ such that the orbit $\{y_1,\cdots,f^{q_1-1}(y_1)\}$ is $\varep$-dense and $d(y_1,f^{q_1}(y_1))<\varep/2$. Again, by density, we can choose a transitive point $y_2$ whose orbit is disjoint from that of $x_0$ and $y_1$, with $d(y_1,y_2)<\varep/2$ and $d(y_0,y_1)-d(y_0,y_2)>\varep/4$. Then there exists a multiple $q_2$ of $p$ such that $d(y_2,f^{q_2}(y_
2))<\varep/2$. And so on, we construct a sequence $(y_m)_{1\le m\le \ell}$ such that (see figure \ref{constrigrec}):
\begin{enumerate}[(i)]
\item for all $m$, there exists $q_m>0$ such that $p| q_m$ and $d(y_m, f^{q_m}(y_m))<\varep/2$,
\item the orbits $\{x_0,\cdots,f^{p-1}(x_0)\}$ and $\{y_m,\cdots,f^{q_m-1}(y_m)\}$ ($m$ going from $0$ to $\ell-1$) are pairwise disjoints,
\item for all $m$, $d(y_m,y_{m+1})<\varep/2$ and $d(y_0,y_m)-d(y_0,y_{m+1})>\varep/4$,
\item $y_\ell = y_0$.
\end{enumerate}

Let $\sigma_{N}$ be a finite map given by Lax's theorem. For all $N$ large enough, $\sigma_{N}$ satisfies the same properties (i) to (iii) than $f$. Changing $\sigma_{N}$ at the points $\sigma_{N}^{q_m-1}((y_m)_{N})$ and $\sigma_{N}^{p-1}((x_0)_{N})$, we obtain a finite map $\sigma'_{N}$ such that ${\sigma'_{N}}^{q_N}((y_m)_{N}) = (y_{m+1})_{N}$ and $\sigma_{N}^{p}((x_0)_{N}) = (x_0)_{N}$. Thus the orbit of $(x_0)_N$ under $\sigma'_{N}$ is $2\varep$-dense and has period $p$ and the orbit of $(y_0)_{N}$ under $\sigma'_{N}$ is $2\varep$-dense, disjoint from which of $(x_0)_{N}$ and has period $1+q_1+\cdots+q_{\ell-1}$ relatively prime to $p$.
\end{proof}

\begin{coro}\label{méldiscr}
For a generic homeomorphism $f\in \Hom(X,\lambda)$, for all $\varep>0$ and all $N_0\in\N$, there exists $N\ge N_0$ such that $f_{N}$ is $\varep$-topologically weakly mixing.
\end{coro}

\begin{proof}[Proof of Corollary \ref{méldiscr}]
Again, we prove the corollary in the case where $M=2$, other cases being easily obtained by induction. Let $\varep>0$ and $N_0\in\N$. All iterates of a generic homeomorphism $f$ are topologically transitive: it is an easy consequence of the genericity of transitivity (see for example Corollary \ref{typlax} or \cite[Theorem 2.11]{MR2931648}); we pick such a homeomorphism. Combining Theorem \ref{génécycl} and Proposition \ref{melfaibl}, we obtain $N\ge N_0$ such that $f_{N}$ has two $\varep/3$-dense periodic orbits whose lengths $p$ et $q$ are coprime. We now have to prove that $f_N$ is $\varep$-topologically weakly mixing. Let $B_1$, $B_2$, $B'_1$ and $B'_2$ be balls with diameter $\varep$. Since each one of these orbits is $\varep/3$-dense, there exists $x_N\in X$ which is in the intersection of the orbit whose length is $p$ and $B_1$, and $y_{N}\in X$ which is in the intersection of the orbit whose length is $q$ and $B_2$. Similarly, there exists two integers $a$ and $b$ such that $f_{N}^a(x_{N})\in B'_1$ 
and $f_{N}^b(y_{N})\in B'_2$.

Recall that we want to find a power of $f_{N}$ which sends both $x_{N}$ in $B'_1$ and $y_{N}$ in $B'_2$. It suffices to pick $m\in\N$ such that $m = a+\alpha p = b+\beta q$. Bézout's identity states that there exists two integers $\alpha$ and $\beta$ such $\alpha p-\beta q = b-a$. Set $m = a+\alpha p$, adding a multiple of $pq$ if necessary, we can suppose that $m$ is positive. Thus $f_{N}^m(x_{N})\in B'_1$ and $f_{N}^m(y_{N})\in B'_2$.
\end{proof}
\bigskip

For now, the two approximation types we studied concern analogues of properties that are generic among $\Hom(X,\lambda)$. We now show that some discrete analogues of properties that are not generic among conservative homeomorphisms also occur infinitely often in the discretizations of generic homeomorphisms. The second variation of Lax's theorem concerns the approximation of applications whose recurrent set is small.

\begin{prop}[Second variation of Lax's theorem]\label{var1}
Let $f\in\Hom(X,\lambda)$. Then for all $\varepsilon, \varep'>0$, there exists $N_0\in\N$ such that for all $N\ge N_0$, there exists a map $\sigma_{N} : E_{N}\to E_{N}$ such that $d_{N}(f,\sigma_{N})<\varepsilon$ and
\[\frac{\card(\Omega(\sigma_N))}{\card(E_{N})} = \frac{\card(\Omega(\sigma_N))}{q_{N}}<\varepsilon',\]
and such that $E_N$ is made of a unique (pre-periodic) orbit of $\sigma_N$.
\end{prop}

\begin{proof}[Proof of Proposition \ref{var1}]
Let $f\in\Hom(X,\lambda)$, $\varepsilon>0$ and a recurrent point $x$ of $f$. There exists $\tau\in\N^*$ such that $d(x,f^\tau(x))<\frac{\varepsilon}{8}$; this inequality remains true for fine enough discretizations: there exists $N_1\in\N$ such that if $N\ge N_1$, then
\[d(x,x_{N})<\frac{\varep}{8},\quad d\big(f^\tau(x),f^\tau(x_{N})\big)<\frac{\varepsilon}{8}\quad \text{and}\quad \frac{\tau}{q_{N}}<\varepsilon'.\]
Using the modulus of continuity of $f^\tau$ and Lax's theorem (Theorem \ref{Lax}), we obtain an integer $N_0\ge N_1$ such that for all $N\ge N_0$, there exists a cyclic permutation $\sigma_{N}$ of $E_{N}$ such that $d_{N}(f,\sigma_{N})<\frac{\varepsilon}{2}$ and $d_{N}(f^\tau,\sigma_{N}^\tau)<\frac{\varepsilon}{8}$. Then
\begin{align*}
d(x_{N},\sigma_N^\tau(x_{N})) \le & d(x_{N},x) + d\big(x,f^\tau(x)\big) + d\big(f^\tau(x),f^\tau(x_{N})\big)\\
                                  & + d\big(f^\tau(x_{N}),\sigma_N^\tau(x_{N})\big)\\
                               <  & \frac{\varepsilon}{2}.
\end{align*}

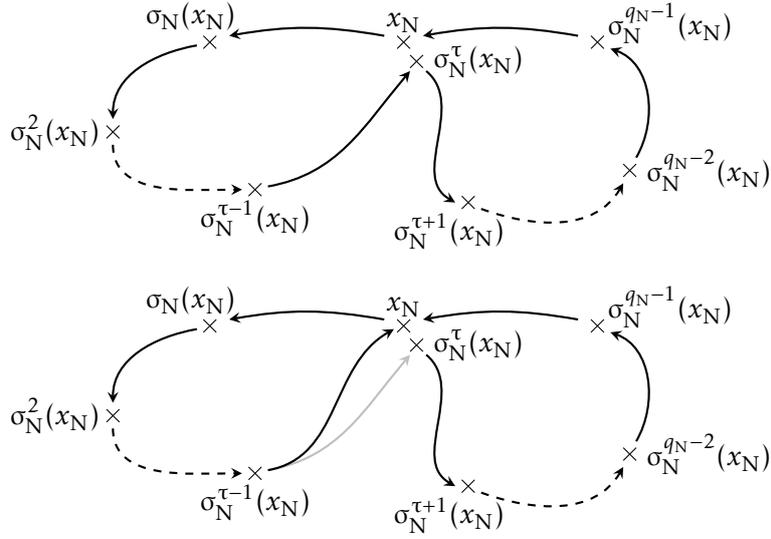
\begin{figure}
\begin{center}
\begin{tikzpicture}[scale=.85]
\draw (-.1,-.1) -- (.1,.1); \draw (-.1,.1) -- (.1,-.1); \draw (0,.3) node{$x_N$};
\draw (-3.1,-.1) -- (-2.9,.1);\draw (-3.1,.1) -- (-2.9,-.1);\draw (-3.3,.4) node{$\sigma_N(x_N)$};
\draw (-4.6,-1.5) -- (-4.4,-1.3);\draw (-4.6,-1.3) -- (-4.4,-1.5);\draw (-4.55,-1.4) node[left]{$\sigma_N^2(x_N)$};
\draw (-2.4,-2.4) -- (-2.2,-2.2);\draw (-2.4,-2.2) -- (-2.2,-2.4);\draw (-2.3,-2.75) node{$\sigma_N^{\tau-1}(x_N)$};
\draw (.1,-.4) -- (.3,-.2); \draw (.1,-.2) -- (.3,-.4);\draw (.3,-.3) node[right]{$\sigma_N^{\tau}(x_N)$};
\draw (.9,-2.4) -- (1.1,-2.6); \draw (.9,-2.6) -- (1.1,-2.4);\draw (.7,-2.95) node{$\sigma_N^{\tau+1}(x_N)$};
\draw (3.4,-1.9) -- (3.6,-2.1); \draw (3.4,-2.1) -- (3.6,-1.9);\draw (3.6,-2) node[right]{$\sigma_N^{q_N-2}(x_N)$};
\draw (3.1,-.1) -- (2.9,.1);\draw (3.1,.1) -- (2.9,-.1);\draw (3,.3) node[right]{$\sigma_N^{q_N-1}(x_N)$};
\draw[->,>=stealth,thick] (-.3,.1) to [out=170,in=10] (-2.7,.1);
\draw[->,>=stealth,thick] (-3.2,-.05) to [out=190,in=90] (-4.5,-1.2);
\draw[->,>=stealth,thick,dashed] (-4.5,-1.6) to [out=-90,in=180] (-2.5,-2.3);
\draw[->,>=stealth,thick] (-2.1,-2.25) to [out=10,in=-130] (.1,-.45);
\draw[->,>=stealth,thick] (.35,-.45) to [out=-40,in=150] (.8,-2.5);
\draw[->,>=stealth,thick,dashed] (1.2,-2.6) to [out=-20,in=-120] (3.4,-2.2);
\draw[->,>=stealth,thick] (3.6,-1.8) to [out=60,in=-20] (3.2,-.1);
\draw[->,>=stealth,thick] (2.7,.1) to [out=170,in=10] (.3,.1);
\draw (0,-3.6);
\end{tikzpicture}

\begin{tikzpicture}[scale=.85]
\draw (-.1,-.1) -- (.1,.1); \draw (-.1,.1) -- (.1,-.1); \draw (0,.3) node{$x_N$};
\draw (-3.1,-.1) -- (-2.9,.1);\draw (-3.1,.1) -- (-2.9,-.1);\draw (-3.3,.4) node{$\sigma_N(x_N)$};
\draw (-4.6,-1.5) -- (-4.4,-1.3);\draw (-4.6,-1.3) -- (-4.4,-1.5);\draw (-4.55,-1.4) node[left]{$\sigma_N^2(x_N)$};
\draw (-2.4,-2.4) -- (-2.2,-2.2);\draw (-2.4,-2.2) -- (-2.2,-2.4);\draw (-2.3,-2.75) node{$\sigma_N^{\tau-1}(x_N)$};
\draw (.1,-.4) -- (.3,-.2); \draw (.1,-.2) -- (.3,-.4);\draw (.3,-.3) node[right]{$\sigma_N^{\tau}(x_N)$};
\draw (.9,-2.4) -- (1.1,-2.6); \draw (.9,-2.6) -- (1.1,-2.4);\draw (.7,-2.95) node{$\sigma_N^{\tau+1}(x_N)$};
\draw (3.4,-1.9) -- (3.6,-2.1); \draw (3.4,-2.1) -- (3.6,-1.9);\draw (3.6,-2) node[right]{$\sigma_N^{q_N-2}(x_N)$};
\draw (3.1,-.1) -- (2.9,.1);\draw (3.1,.1) -- (2.9,-.1);\draw (3,.3) node[right]{$\sigma_N^{q_N-1}(x_N)$};
\draw[->,>=stealth,thick] (-.3,.1) to [out=170,in=10] (-2.7,.1);
\draw[->,>=stealth,thick] (-3.2,-.05) to [out=190,in=90] (-4.5,-1.2);
\draw[->,>=stealth,thick,dashed] (-4.5,-1.6) to [out=-90,in=180] (-2.5,-2.3);
\draw[->,>=stealth,thick,color=gray!50] (-2.1,-2.25) to [out=10,in=-130] (.1,-.45);
\draw[->,>=stealth,thick] (-2.1,-2.25) to [out=10,in=-150] (-.15,-.05);
\draw[->,>=stealth,thick] (.35,-.45) to [out=-40,in=150] (.8,-2.5);
\draw[->,>=stealth,thick,dashed] (1.2,-2.6) to [out=-20,in=-120] (3.4,-2.2);
\draw[->,>=stealth,thick] (3.6,-1.8) to [out=60,in=-20] (3.2,-.1);
\draw[->,>=stealth,thick] (2.7,.1) to [out=170,in=10] (.3,.1);
\end{tikzpicture}
\caption[Modification of a cyclic permutation]{Modification of a cyclic permutation in the proof of Proposition \ref{var1}}\label{trajectoire}
\end{center}
\end{figure}

We compose $\sigma_{N}$ by the (non bijective) application mapping $\sigma_{N}^{\tau-1}(x_{N})$ on $x_{N}$ and being identity anywhere else (see Figure \ref{trajectoire}), in other words we consider the application
\[\sigma'_{N}(x) = \left\{\begin{array}{ll}
x_{N} \quad & \text{if}\ x=\sigma_{N}^{\tau-1}(x_{N})\\
\sigma_{N}(x) \quad        & \text{otherwise.}
\end{array}\right.\]
The map $\sigma'_{N}$ has a unique injective orbit whose associated periodic orbit $\Omega(\sigma'_N)$ has length $\tau$ (it is $(x_{N},\sigma_N(x_{N}),\cdots,\sigma_{N}^{\tau-1}(x_{N}))$). Since $d(f,\sigma'_N)<\varep$, the map $\sigma'_{N}$ verifies the conclusions of the proposition.
\end{proof}

A direct application of Theorem \ref{génécycl} leads to the following corollary.

\begin{coro}\label{corovar1}
For a generic homeomorphism $f\in\Hom(X,\lambda)$,
\[\underline\lim_{N\to +\infty}\frac{\card(\Omega(f_N))}{\card(E_{N})} =0.\]
Specifically for all $\varepsilon>0$ and all $N_0\in\N$, there exists $N\ge N_0$ such that $\card(\Omega(f_N)) / \card(E_{N}) <\varepsilon$ and such that $E_N$ is made of a unique (pre-periodic) orbit of $f_N$.
\end{coro}

We now improve this corollary in stating that for a generic homeomorphism, there exists $C>0$ such that we have $\card(\Omega(f_N))\le C$ for an infinite number of orders $N$; in particular these discretizations are highly non transitive, which is the opposite behaviour to the dynamics of a generic homeomorphism.

The same kind of idea than in the proof of Proposition \ref{var1} leads to the following variation of Lax's theorem.

\begin{prop}\label{propvar2}
Let $f\in\Hom(X,\lambda)$ having at least one periodic point $x$ of period $p$. Then for all $\varep>0$, there exists $N_0\in\N$ such that for all $N\ge N_0$, there exists an application $\sigma_{N} : E_{N}\to E_{N}$ with $d_{N}(f,\sigma_{N})<\varepsilon$, such that $E_N$ is made of a unique (pre-periodic) orbit of $\sigma_N$ and such that the unique periodic orbit of $f_N$ is of length $p$ and $\varep$-shadows the $f$-orbit of $x$.
\end{prop}

\begin{proof}[Proof of proposition \ref{propvar2}]
Simply replace the recurrent point by a periodic point of period $p$ in the proof of Proposition \ref{var1}.
\end{proof}

We will use this proposition in Section~\ref{Sec8} to obtain Theorem \ref{corovar2}. In particular, it will imply the following statement.

\begin{coro}\label{corovarN}
For a generic homeomorphism $f\in\Hom(X,\lambda)$, for every period $p$ of a periodic point of $f$ and for infinitely many integers $N$, $f_{N}$ has a unique periodic orbit, whose length is $p$, and such that $E_N$ is covered by a single (pre-periodic) orbit of $f_N$.
\end{coro}

\begin{rem}
In particular, if $p_0$ denotes the minimal period of the periodic points of $f$, then there are infinitely many discretizations such that the cardinality of their recurrent set is equal to $p_0$. Note that however, the shortest period of periodic points of generic homeomorphisms has no global upper bound in $\Hom(X,\lambda)$: for example, for all $p\in\N$, there is an open set of homeomorphisms of the torus without periodic point of period less than $p$ (\emph{e.g.} the neighbourhood of an irrational rotation) and this property remains true for discretizations.
\end{rem}

Corollary \ref{corovarN} states that for a generic homeomorphism, $\card(\Omega(f_N)) / \card(E_{N})$ is as small as possible for an infinite number of orders $N$. The next result states that this loss of injectivity can even occur from the first iteration of $f_N$.

\begin{prop}[Third variation of Lax's theorem]\label{crunch}
Let $f\in\Hom(X,\lambda)$ and $\vartheta : \N\to\R_+^*$ a map which tends to $+\infty$ at $+\infty$. Then for all $\varep>0$, there exists $N_0\in\N$ such that for all $N\ge N_0$, there exists a map $\sigma_{N} : E_{N}\to E_{N}$ such that $\card(\sigma_{N}(E_N))<\vartheta({N})$ and $d_{N}(f,\sigma_{N})<\varep$.
\end{prop}

\begin{proof}[Proof of Proposition \ref{crunch}]
Let $f\in\Hom(X,\lambda)$, $\vartheta : \N\to\R_+^*$ a map which tends to $+\infty$ at $+\infty$ and $\varep>0$. By Lax's theorem (Theorem \ref{Lax}) there exists $N_1\in\N$ such that for all $N\ge N_1$, there exists a cyclic permutation $\sigma_{N} : E_{N}\to E_{N}$ whose distance to $f$ is smaller than $\varep/2$. For $N\ge N_1$, set $\sigma'_{N} = P_{N_1}\circ \sigma_{N}$. Increasing $N_1$ if necessary we have $d(f,\sigma'_{N})<\varep$, regardless of $N$. Moreover $\card(\sigma'_{N}(E_{N}))\le q_{N_1}$; if we choose $N_0$ large enough such that for all $N\ge N_0$ we have $q_{N_1}<\vartheta({N})$, then $\card(\sigma'_{N}(E_{N}))\le \vartheta({N})$. We have shown that the map $\sigma'_{N}$ satisfies the conclusions of proposition for all $N\ge N_0$.
\end{proof}

\begin{coro}\label{crush}
Let $\vartheta : \N\to\R_+^*$ a map which tends to $+\infty$ at $+\infty$. Then for a generic homeomorphism $f\in\Hom(X,\lambda)$, 
\[\underset{N\to+\infty}{\underline\lim}\ \frac{\card(f_{N}(E_{N}))}{\vartheta({N})}=0.\]
In particular, generically, $\underline{\lim}\frac{\card(f_{N}(E_{N}))}{\card(E_{N})}=0$.
\end{coro}

\begin{proof}[Proof of Corollary \ref{crush}]
Remark that if we replace $\vartheta({N})$ by $\sqrt{\vartheta({N})}$, it suffices to prove that for a generic homeomorphism, $\underline\lim\frac{\card(f_{N}(E_{N}))}{\vartheta({N})}\le 1$. This is easily obtained by combining Theorem \ref{génécycl} and Proposition \ref{crunch}.
\end{proof}

So far all variations of Lax's theorem have constructed finite maps with a small number of cycles. With the additional assumption that the sequence of grids is sometimes self similar, we show a final variation of Lax's theorem, approaching every homeomorphism by a finite map with a large number of orbits.

\begin{prop}[Fourth variation of Lax's theorem]\label{killing}
Assume that the sequence of grids $(E_N)_{N\in\N}$ is sometimes self similar. Let $f\in\Hom(X,\lambda)$ and $\vartheta : \N\to\R$ such that $\vartheta(N) = o(q_N)$. Then for all $\varepsilon>0$, there exists $N_1\in\N$ such that for all $N\ge N_1$, there exists a permutation $\sigma_{N}$ of $E_{N}$ such that $d_{N}(f,\sigma_{N})<\varepsilon$ and that the number of cycles of $\sigma_{N}$ is greater than $\vartheta({N})$. Moreover, $\vartheta({N})$ of these cycles of $\sigma_{N}$ are conjugated to a cyclic permutation of $E_{N_0}$ by bijections whose distance to identity is smaller than $\varep$.
\end{prop}

\begin{proof}[Proof of Proposition \ref{killing}]
Let $\varep>0$, for all $N_0\in\N$ large enough, Lax's theorem gives us a cyclic permutation $\sigma'_{N_0}$ of $E_{N_0}$ whose distance to $f$ is smaller than $\varep$. Since the grids are sometimes self similar, there exists $N_1\in\N$ such that for all $N\ge N_1$, the set $E_N$ contains $q_N/q_{N_0}\ge\vartheta(N)$ disjoint subsets $\widetilde E_N^j$, each one conjugated to a grid $E_{N_0}$ by a bijection $h_j$ whose distance to identity is smaller than $\varep$. On each $\widetilde E_N^j$, we define $\sigma_N$ as the conjugation of $\sigma'_{N_0}$ by $h_j$; outside these sets we just pick $\sigma_N$ such that $d_N(f,\sigma_N)<\varep$. Since the distance between $h_j$ and identity is smaller than $\varep$, we have $d_N(f,\sigma_N)<2\varep$. Moreover, $\sigma_N$ has at least $\vartheta(N)$ cycles which are conjugated to a cyclic permutation of $E_{N_0}$ by the maps $h_j$; this completes the proof.
\end{proof}

The application of Theorem \ref{génécycl} leads to the following corollary, which ensures that an infinite number of discretizations $f_N$ have a lot of periodic orbits.

\begin{coro}\label{corovar3}
We still assume that the sequence of grids $(E_N)_{N\in\N}$ is sometimes self similar. Let $\vartheta : \N\to\R$ such that $\vartheta(N) = o(q_N)$. Then for a generic homeomorphism $f\in\Hom(X,\lambda)$ and for infinitely many integers $N$, the discretization $f_{N}$ of $f$ has at least $\vartheta(N)$ periodic orbits (which all have the same period).
\end{coro}

If we compose the map obtained by Proposition \ref{killing} by an appropriate map of the finite set $\{h_i(x_0)\}_i$ into itself ($x_0$ being fixed), we can prove the following result.

\begin{prop}
Let $f\in\Hom(X,\lambda)$, $\varepsilon>0$ and $]a,b[\subset [0,1]$. Then there exists $N_1\in\N$ such that for all $N\ge N_1$, there exists a permutation $\sigma_{N}$ of $E_{N}$ such that $d_{N}(f,\sigma_{N})<\varepsilon$ and that the degree of recurrence of $\sigma_N$ (that is, the ratio between the cardinality of the recurrent set of $\sigma_N$ and the cardinality of $E_N$) belongs to $]a,b[$.
\end{prop}

This proposition implies trivially the following corollary.

\begin{coro}\label{ConjEt}
For a generic conservative homeomorphism $f\in\Hom(X,\lambda)$, the sequence $D(f_N)$ of the degrees of recurrence of the discretizations accumulates on the whole segment $[0,1]$.
\end{coro}

\section{Average behaviour of discretizations}\label{bofbof}

We now want to study the average behaviour of discretizations of a generic homeomorphism. For example, we could imagine that even if for a generic homeomorphism $f$, the event ``$f_N$ is a cyclic permutation'' appears for infinitely many orders $N$, it is statistically quite rare.

More precisely, we study the frequency of occurrence of properties related to the discretizations of generic homeomorphisms: given a property $(P)$ concerning discretizations, we are interested in the behaviour of the proportion between $1$ and $M$ of discretizations satisfying the property $(P)$, when $ M $ goes to infinity. For this study, we assume that the sequence of discretization grids is always self similar (which is true for example for the torus equipped with discretizations upon uniform grids of orders powers of an integer, see Section \ref{exgrilles}). This prevents us from having to deal with tricky arithmetic problems about overlay of grids.

\begin{definition}
Let $f\in\Hom(X,\lambda)$. We say that a property $(P)$ concerning discretizations is \emph{satisfied in average} if for all $N_0\in\N$ and all $\varep>0$, there exists $N\ge N_0$ such that the proportion of integers $M\in \{ 0,\cdots,N\}$ such that $f_M$ satisfies $(P)$ is greater than $1-\varep$, \emph{i.e.}
\[\underset{N\to+\infty}{\overline \lim}\,\frac{1}{N+1}\card\big\{M\in\{ 0,\cdots,N\}\mid f_M\text{ satisfies }(P)\big\}=1.\]
\end{definition}

We will show that most of the dynamical properties studied in the previous section are actually satisfied on average for generic homeomorphisms. To start with, we set out a technical lemma:

\begin{lemme}\label{recopie}
Let $\tT$ be a dense type of approximation in $\Hom(X,\lambda)$. Then for a generic homeomorphism $f\in\Hom(X,\lambda)$, for all $\varep>0$ and all $\alpha>0$, the property $(P)$: ``$E_N$ contains at least $\alpha$ disjoints subsets which fills a proportion greater than $1-\varep$ of $E_N$, each of them stabilized by $f_N$ and such that the restriction of $f_N$ to each one is conjugated to a map of $\tT$ by a bijection whose distance to identity is smaller than $\varep$'' is satisfied in average.
\end{lemme}

In practice, this lemma provides many properties satisfied on average, for instance:
\begin{itemize}
\item quantitative properties on discretizations, such as possessing at least $M$ periodic orbits;
\item properties of existence of sub-dynamics on discretizations, such as possessing at least one dense periodic orbit.
\end{itemize}

\begin{proof}[Proof of Lemma \ref{recopie}]
Let us consider the set
\[\mathcal{C} = \bigcap_{\substack{\varepsilon>0\\ N_0\in\N}}\,\bigcup_{N\ge N_0}\,\left\{
\begin{array}{c}
f\in\Hom(X,\lambda)\ \big|\\ \frac{1}{N+1}\card\big\{M\in\{ 0,\cdots,N\}\mid f_M\text{ satisfies }(P)\big\}>1-\varep
\end{array}\right\}.\]
We want to show that $\mathcal{C}$ contains a dense $G_\delta$ of $\Hom(X,\lambda)$. The set $\mathcal{C}$ is a $G_\delta$ of the generic set  $\bigcap_{N\in\N} \mathcal{D}_{N}$, it suffices to prove that it is dense in $\bigcap_N \mathcal{D}_{N}$. Let $f\in \Hom(X,\lambda)$, $N_0\in\N$, $\delta>0$ and $\varepsilon>0$. To prove it we want to find a homeomorphism $g$ whose distance to $f$ is smaller than $\delta$ and an integer $N\ge N_0$ such that
\[\frac{1}{N+1}\card\big\{M\in\{ 0,\cdots,N\}\mid g_M\text{ satisfies }(P)\big\}>1-\varep.\]
It is simply obtained in combining the density of the type of approximation $\tT$ and the fact that the sequence of grids is always self similar.
\end{proof}

This lemma allows us to obtain some properties about the average behaviour of discretizations. For instance here is an improvement of Corollary~\ref{corovar2}.

\begin{coro}\label{CoroAver1}
For a generic homeomorphism $f\in\Hom(X,\lambda)$, the property ``$f_{N}$ has a $\varep$-dense periodic orbit and the cardinality of $\Omega(f_N)$ satisfies  $\card(\Omega(f_N)) = o(q_N)$'' is satisfied in average.
\end{coro}

Or an improvement of Corollary \ref{méldiscr}.

\begin{coro}\label{CoroAver2}
For a generic homeomorphism $f\in\Hom(X,\lambda)$ and for all $\varep>0$, the property ``$f_{N}$ is $\varep$-topologically weakly mixing'' (see Definition \ref{epmélfaibl}) is satisfied in average.
\end{coro}

Or even an improvement of Corollary \ref{crush}.

\begin{coro}\label{petitepermieux}
For a generic homeomorphism $f\in\Hom(X,\lambda)$ and for all $\varep>0$, the property ``$\frac{\card(f_{N}(E_{N}))}{\card(E_{N})}<\varep$'' is satisfied in average.
\end{coro}

And an improvement of Corollary \ref{corovar3}.

\begin{coro}\label{CoroAver3}
For a generic homeomorphism $f\in\Hom(X,\lambda)$ and for all $M\in\N$, property ``$f_N$ has at least $M$ periodic orbits'' is satisfied in average.
\end{coro}

However, note that the most simple property about discretizations, \emph{i.e.} being a cyclic permutation, cannot be proved by using Lemma \ref{recopie}. To do this, we need a slightly more precise result, that requires the hypothesis that the grids are always strongly self similar.

\begin{theoreme}\label{propdemin}
For a generic homeomorphism $f\in\Hom (X,\lambda)$, the property ``$f_{N}$ is a cyclic permutation'' is satisfied in average.
\end{theoreme}

Note that this implies that most of the discretizations of $f$ does not behave like a random map of a set of cardinality $q_N$, as a random map of a set with $q$ elements has in average $\log q$ periodic orbits (see for example \cite[XIV.5]{Boll-rand}).

\begin{lemme}\label{demin}
Let $f\in\Hom(X,\lambda)$ and $\varep>0$. The sequence of grids is supposed to be always strongly self similar. Then there exists $N_0\in\N$ such that for all $N\ge N_0$, there exists a permutation $\sigma_N$ of $E_N$, such that $d_N(f,\sigma_N)<\varep$, ${\sigma_N}_{|E_{N_0}}$ is a cyclic permutation of $E_{N_0}$ and for all $M\in\{ N_0,\cdots,N-1\}$, $\sigma_N$ is a cyclic permutation of $E_{M+1}\setminus E_{M}$.
\end{lemme}

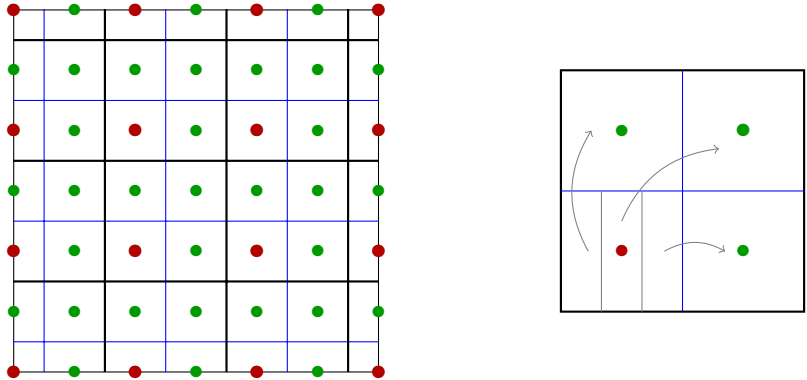
\begin{figure}[ht]
\begin{center}
\begin{tikzpicture}[scale=.8]
\draw (.5,.5) rectangle (6.5,6.5);
\foreach \k in {0,...,6}
 {\foreach \l in {0,...,6}
  {\draw[color=green!60!black] node at (\k+.5,\l+.5){$\bullet$};}}
\foreach \k in {0,...,3}
 {\foreach \l in {0,...,3}
  {\draw[color=red!70!black] node at (2*\k+.5,2*\l+.5){\large$\bullet$};}}

\draw[thick] (9.5,1.5) rectangle (13.5,5.5);
\draw[blue] (11.5,1.5) -- (11.5,5.5);
\draw[blue] (9.5,3.5) -- (13.5,3.5);
\draw[color=green!60!black] node at (12.5,4.5){\large$\bullet$};
\draw[color=green!60!black] node at (12.5,2.5){$\bullet$};
\draw[color=green!60!black] node at (10.5,4.5){$\bullet$};
\draw[color=red!70!black] node at (10.5,2.5){$\bullet$};
\draw[gray] (10.8333,1.5) -- (10.8333,3.5);
\draw[gray] (10.1666,1.5) -- (10.1666,3.5);

\draw[->,gray] (9.95,2.5) to[bend left] (10,4.5);
\draw[->,gray] (10.5,3) to[bend left] (12.1,4.2);
\draw[->,gray] (11.2,2.5) to[bend left] (12.2,2.5);

\clip (.5,.5) rectangle (6.5,6.5);
\draw[blue] (0,0) grid (7,7);
\draw[step=2, thick] (0,0) grid (7,7);
\end{tikzpicture}
\caption[Construction of Lemma \ref{demin}]{Construction\!\protect\footnotemark\ of Lemma \ref{demin} for two grids of orders $3$ and $6$ and zoom around a point of a grid of order $3$}
\label{soucub}
\end{center}
\end{figure}
\footnotetext{Do not forget that we identify some points of the boundary!}

\begin{proof}[Proof of Lemma \ref{demin}]
Applying Lax's theorem (Theorem \ref{Lax}), we can find an integer $N_0$ and a cyclic permutation $\sigma_{N_0}$ of $E_{N_0}$ such that $d_{N_0}(f,\sigma_{N_0})<\varep$. As the sequence of grids is always strongly self similar, we can suppose that $N_0$ is big enough to verify the conclusions of Definition~\ref{Ashe'}.

Let us observe what happens for the order $N_0+1$. We will define an application $\sigma'_{N_0+1}$ on $E_{N_0+1}$, which will be close to the homeomorphism $f$. On $E_{N_0}$, we define $\sigma'_{N_0+1}$ as being equal to $\sigma_{N_0}$. The idea is to repeat the proof of Lax's theorem for the elements of $E_{N_0+1}\setminus E_{N_0}$. To do that, we have to find a partition of $X$ into sets with the same measure, such that every element $A_x$ of this partition is one to one associated to a point $x$ of $E_{N_0+1}\setminus E_{N_0}$, and has ``small'' diameter. So it suffices to split equitably the mass of the cubes  $E_{N_0}$ to the other cubes.

For this purpose, we cut each cube of $E_{N_0+1}$ corresponding to a point $x_{N_0}$ of the grid $E_N$ into $\alpha_N-1$ subsets of the same measure $1/(q_{N_0+1}(\alpha_N-1))$ (see Figure \ref{soucub}). Each of these subsets is assigned to the cube corresponding to the point $h_i(x_{N_0})$. For each $x\in E_{N_0+1}\setminus E_{N_0}$, we define $A_x$ as the union of the cube $C_x$ with the subset that is associated to it. Then all the sets $A_x$ have the same measure and have a diameter smaller than twice the diameter of a cube of $E_{N_0+1}$ plus $\varep$. We can now apply Lax's theorem to the homeomorphism $f$ and the partition $\{A_x\}_{x\in E_N}$ of $X$, this partition being numbered the same way as the partition $E_{N_0+1}$ (some numbers are skipped); this gives a cyclic permutation of $E_{N_0+1}\setminus E_{N_0}$ which is close to $f$. This cyclic permutation defines the permutation $\sigma'_{N_0+1}$ where it was not yet. This $\sigma'_{N_0+1}$ is close to $f$ and permutes cyclically the elements of $E_
{N_0}$ and these of $E_{N_0+1}\setminus E_{N_0}$.

We finish the proof by an induction, iterating $N-N_0$ times.
\end{proof}

\begin{coro}\label{bidouille}
Let $f\in\Hom(X,\lambda)$ and $\varep>0$ (we recall that the sequence of grids is supposed to be always strongly self similar). Then there exists $N_0\in\N$ such that for all $N\ge N_0$, there exists a homeomorphism $g$ whose distance to $f$ is smaller than $\varep$ such that for all $M\in\{ N,\cdots,N_0\}$, the discretization $g_{M}$ is a cyclic permutation of $E_{M}$. Moreover, this property can be supposed to be verified on a neighbourhood of $g$.
\end{coro}

\begin{proof}[Proof of Corollary \ref{bidouille}]
Let $f\in\Hom(X,\lambda)$ and $\varep>0$. Lemma \ref{demin} gives us an integer $N_0\in\N$ such that for all $N\ge N_0$, there exists a permutation $\sigma_{N}$ of $E_{N}$, whose distance to $f$ is smaller than $\varep/2$, such that $\sigma_{N}$ is a cyclic permutation of $E_{N_0}$ and for all $M\in\{ N_0+1,\cdots,N\}$, $\sigma_{M}$ is a cyclic permutation of $E_{M}\setminus E_{M-1}$. Set $N\ge N_0$, the idea is to modify slightly the orbits of $\sigma_{N}$ such that $\sigma_{N}$ becomes a cyclic permutation on all the sets $E_{M}$.

We begin by choosing two points $x_{N_0}\in E_{N_0}$ and $x'_{N_0+1}\in E_{N_0+1}$ such that $x'_{N_0+1}$ belongs to the cube corresponding to the point $x_{N_0}$. Then, we modify $\sigma_{N}$ in interchanging the points $x_{N_0}$ et $x'_{N_0+1}$, in other words we set
\[\sigma_{N}^{N_0}(x) = \left\{\begin{array}{ll}
\sigma_{N}(x'_{N_0+1}) \quad & \text{if}\ x=x_{N_0}\\
\sigma_{N}(x_{N_0}) \quad      & \text{if}\ x = x'_{N_0+1}\\
\sigma_{N}(x) \quad        & \text{otherwise.}
\end{array}\right.\]
Thus, $\sigma_{N}^{N_0}$ is a cyclic permutation of $E_{N_0+1}$ and the discretization of order $N_0+1$ of $\sigma_{N}^{N_0}$ is also cyclic. We build the same way some maps $\sigma_{N}^{N_0+1},\cdots,\sigma_{N}^{N}$ in interchanging the images of two adjacent points lying in the grids $E_M$ et $E_ {M+1}$. Then, for all $M\in\{ N_0,\cdots,N\}$, the discretization of $\sigma_{N}^{M-1}$ on $E_{M}$ is a cyclic permutation. Moreover, the distance between the map $\sigma_{N}^{N}$ and $\sigma_N$ is smaller than $\varep/2$. Combined with Proposition \ref{extension}, this allows us to interpolate $\sigma_{N}^{N}$ to a homeomorphism $g$ whose distance to $f$ is smaller than $\varep$ such that for all $M\in\{ N_0,\cdots,N\}$, the discretization of order $M$ of $g$ is a cyclic permutation of $E_M$. A careful reader would notice that it may happen that in this construction, the discretizations of $g$ are not uniquely defined; depending on the choices made during the definition of $P_N$ on $E'_N$, these discretizations 
may not be cyclic permutations. If we want to avoid this problem, it suffices to modify slightly the map $\sigma_{N}^{M}$; moreover this ensures that the conclusions of the corollary are verified on a whole neighbourhood of $g$.
\end{proof}

\begin{proof}[Proof of Theorem \ref{propdemin}]
Let the set
\[\mathcal{C} = \bigcap_{\substack{\varepsilon>0\\ N_0\in\N}}\,\bigcup_{N\ge N_0}\,\left\{\begin{array}{c}
\rule[-0.3cm]{0cm}{.8cm} f\in\Hom(\T^n,\Leb)\mid\\
\frac{1}{N+1}\card\big\{M\in\{ 0,\cdots, N\}\mid f_{M}\text{ cyclic permutation}\big\}>1-\varep
\end{array}\right\},\]
we want to show that $\mathcal{C}$ contains a $G_\delta$ dense subset of $\Hom(X,\lambda)$. Let $f\in \Hom(X,\lambda)$, $N_0\in\N$, $\delta>0$ and $\varepsilon>0$. This is equivalent to find a homeomorphism $g$ whose distance to $f$ is smaller than $\delta$, and an integer $N\ge N_0$, such that for all homeomorphism $g'$ close enough to $g$, we have
\[\frac{1}{N+1}\card\big\{M\in\{ 0,\cdots,N\}\mid g'_{M}\in P\big\}>1-\varep.\]
Such a homeomorphism $g$ is simply given by Corollary \ref{bidouille} for $N\ge N_0$ such that $\varep k_N>k_{N_0}$.
\end{proof}

\begin{rem}
However, the property of approximation by bicyclic permutations in average cannot be proven directly with this technique.
\end{rem}

\begin{rem}
A simple calculation shows that everything that has been done in this section also applies to the behaviour of discretizations in average of Cesàro average, in average of average of Cesàro average etc., \emph{i.e.} when studying quantities
\[\frac{1}{N_2+1}\sum_{N_1=0}^{N_2}\frac{1}{N_1+1}\card\big\{M\in\{ 0,\cdots,N_1\}\mid f_M\text{ satisfies }(P)\big\},\]
\[\frac{1}{N_3+1}\sum_{N_2=0}^{N_3}\frac{1}{N_2+1}\sum_{N_1=0}^{N_2}\frac{1}{N_1+1}\card\big\{M\in\{ 0,\cdots,N_1\}\mid f_M\text{ satisfies }(P)\big\}\dots\]
\end{rem}

\section{Behaviour of all the discretizations}\label{Sec8}

In the previous sections, we showed that the dynamical behaviour of discretizations depends drastically on the order of discretization: even when looking at the frequency a property appears on discretizations, a lot of different behaviours can occur. In contrast, the dynamics of a generic homeomorphism is well known (see for example \cite{MR2931648}) and ``independent from the homeomorphism''. We even have a 0-1 law on $\Hom(X,\lambda)$ (see \cite{Glas-zero} or the final chapter of \cite{MR2931648}) which states that either a given ergodic property is generic among conservative homeomorphisms, or its contrary is generic. Thus, the dynamics of a generic homeomorphism and that of its discretizations seem quite of independent. In this chapter, we will see that the variations of the dynamics of the discretizations are as large as possible. Indeed, in general, from the convergence $d_N(f,f_N)\to 0$, it can be deduced by a compactness argument that the accumulation set of the dynamical invariants (for example, 
invariant measures) of $f_N$ is included in that of $f$. We will see that in the generic case, the inverse inclusion is also true: every dynamical invariant of $f$ is an accumulation point of the corresponding dynamical invariants of $f_N$. Thus, in a certain sense, it is possible to deduce some dynamical features of a generic homeomorphism from the corresponding dynamical features of \emph{all} its discretizations. In other words, there is a shadowing property of the dynamical features of the homeomorphism by that of its discretizations: for each dynamical feature of the homeomorphism, its discrete analogues appear on an infinite number of discretizations.

For example, we will show that every periodic orbit of $f$ is shadowed by a periodic orbit of an infinite number of $f_N$ (Theorem \ref{corovar2}). We will also show that every $f$-invariant measure is the limit of a subsequence of $f_N$-invariant measures (Theorem \ref{EnsMesInvSimpl}), and every $f$-invariant chain transitive compact set is the limit of a subsequence of $f_N$-invariant sets (Theorem \ref{CompactInvSimpl}). Moreover, in Chapter \ref{ChapRot}, we will give an application of this convergence of all the dynamics of the discretizations to that of the homeomorphism. This will give an algorithm to obtain numerically the rotation set of a generic conservative homeomorphism: the upper limit (for the inclusion) of the rotation sets of the discretizations of a generic conservative homeomorphism of the torus is equal to the rotation set of the homeomorphism itself (see Corollary \ref{CoroRotDiscrCons}).

\subsection{Periodic orbits}

Firstly, we deduce from Proposition \ref{propvar2} that every periodic orbit of a generic homeomorphism is shadowed by a periodic orbit with the same period of an infinite number of discretizations. We will see later that this property of shadowing is true for a larger class of compact sets (namely the chain transitive invariant compact sets, see Theorem \ref{CompactInvSimpl}), but the following theorem is also concerned with the period of the periodic orbits of the discretizations. Moreover, the case of periodic orbits seems natural enough to be handled separately.

\begin{theoreme}\label{corovar2}
Let $f\in\Hom(X,\lambda)$ be a generic homeomorphism. Then, for every periodic point $x$ of period $p$ for $f$, for every $\delta>0$ and for every $N_0\in\N$, there exists $N\ge N_0$ such that $f_N$ has a unique periodic orbit, whose length is $p$, which $\delta$-shadows the $f$-orbit of $x$, and such that $E_N$ is covered by a single (pre-periodic) orbit of $f_N$ (in particular, this implies that this periodic orbit attracts the whole set $E_N$).
\end{theoreme}

In particular, this theorem implies that it is possible to detect the periods of the periodic orbits of a generic homeomorphism from that of its discretizations.

\begin{coro}\label{détecpér}
Let $f$ be a generic homeomorphism of $\Hom(X,\lambda)$ and $p$ be an integer. Then $f$ has a periodic orbit of period $p$ if and only if there exists infinitely many integers $N$ such that $f_N$ has a periodic orbit with period~$p$.
\end{coro}

\begin{proof}[Proof of Corollary \ref{détecpér}]
The first implication is Corollary \ref{corovar2}, and the other follows easily from a compactness argument.
\end{proof}

F. Daalderop and R. Fokkink have proved in \cite{Daal-chao} (see also \cite[Part 3.2]{MR2931648}, the proof is similar to that of Theorem \ref{corovar2}) that for a generic element of $f\in\Hom(X,\lambda)$, the set of $f$-periodic points is dense in $X$.

To prove Theorem \ref{corovar2}, the concept of \emph{persistent point} will allow us to define open properties about periodic points.

\begin{definition}\label{rétine}
Let $f\in\Hom(X)$. A periodic point $x$ of $f$ with period $p$ is said \emph{persistent} if for all neighbourhood $U$ of $x$, there exists a neighbourhood $\V$ of $f$ in $\Hom(X)$ such that every $\widetilde{f}\in \V$ has a periodic point $\widetilde{x}\in U$ with period $p$.
\end{definition}

\begin{ex}\label{linéaire}
The endomorphism $h = \operatorname{Diag}(\lambda_1,\cdots,\lambda_n)$ of $\R^n$, with $\prod \lambda_i = 1$ and $\lambda_i\neq 1$ for all $i$, is measure-preserving and has a persistent fixed point at the origin (see for example \cite[page 319]{MR1326374}). Let $s$ be a reflection of $\R^n$, the application $h\circ s$ is also measure-preserving and has a persistent fixed point at the origin.
\end{ex}

To create persistent periodic points, we use the theorem of local modification of conservative homeomorphisms, which allows to replace locally a homeomorphism by another. Although it seems ``obvious'' and has an elementary proof in dimension two, the proof in higher dimensions uses the (difficult) \emph{annulus theorem}. For more details we may refer to \cite{Daal-chao} or \cite[Part 3.1]{MR2931648}.

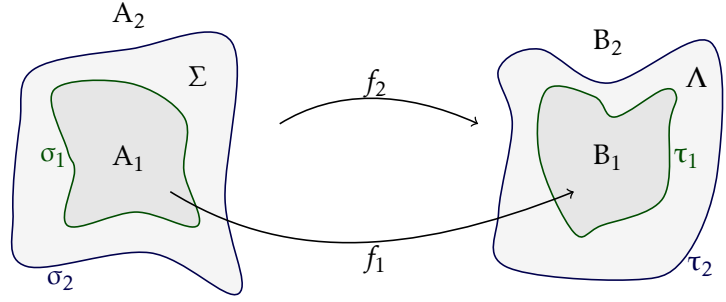
\begin{figure}
\begin{center}
\begin{tikzpicture}[scale=.9]
\draw[color=blue!30!black, fill=gray!8!white, semithick] plot[tension=0.6, smooth cycle] coordinates{(-1.4,1.3) (0,1.6) (1.5,1.7) (1.4,-.6) (1.6,-2) (0.3,-1.4) (-1.5,-1.6) (-1.7,-.8) (-1.6,.5)};
\draw[color=green!30!black, fill=gray!20!white, semithick] plot[tension=0.6, smooth cycle] coordinates{(-1.1,1) (0,1.1) (.8,.6) (.8,-.2) (1,-1) (0.1,-.8) (-.9,-1) (-.8,-.2) (-.9,.1)};
\draw[color=blue!30!black, fill=gray!8!white, semithick] plot[tension=0.6, smooth cycle] coordinates{(5.4,1.3) (6,1.6) (7,1.1) (8.5,1.7) (8.6,-.2) (8.1,-1.5) (7.3,-1.8) (5.5,-1.6) (5.4,-.8) (5.6,.5)};
\draw[color=green!30!black, fill=gray!20!white, semithick] plot[tension=0.6, smooth cycle] coordinates{(6.1,1) (6.8,.9) (7.1,.6) (7.8,1) (8,.9) (7.9,0.5) (7.8,-.6) (6.9,-1) (6.5,-1.1) (6,0)};
\draw[->, semithick] (.6,-.5) .. controls (2.2,-1.5) and (4,-1.5) .. (6.5,-.5);
\draw[->, semithick] (2.2,.5) .. controls (3,1) and (4,1) .. (5.1,.5);
\node at (0,0) {$A_1$};
\node[color=green!30!black] at (-1.1,0) {$\sigma_1$};
\node[color=blue!30!black] at (-1,-1.8) {$\sigma_2$};
\node at (1,1.2) {$\Sigma$};
\node at (0,2.1) {$A_2$};
\node at (3.6,-1.5) {$f_1$};
\node at (3.6,1.1) {$f_2$};
\node at (7,0) {$B_1$};
\node[color=green!30!black] at (8.15,0) {$\tau_1$};
\node[color=blue!30!black] at (8.35,-1.6) {$\tau_2$};
\node at (7,1.7) {$B_2$};
\node at (8.3,1.2) {$\Lambda$};
\end{tikzpicture}\caption{Technique of local modification}\label{dessin-extension}
\end{center}
\end{figure}

\begin{theoreme}[Local modification]\label{extension-sphères}
Let $\sigma_1$, $\sigma_2$, $\tau_1$ and $\tau_2$ be four bicollared embeddings\footnote{An embedding $i$ of a manifold $M$ in $\R^n$ is said \emph{bicollared} if there exists an embedding  $j :M\times [-1,1] \to \R^n$ such that $j_{M\times \{0\}} = i$.} of $\mathbf{S}^{n-1}$ in $\R^n$, such that $\sigma_1$ is in the bounded connected component\footnote{By the Jordan-Brouwer theorem the complement of a set which is homeomorphic to $\mathbf{S}^{n-1}$ has exactly two connected components: one bounded and one unbounded.} of $\sigma_2$ and $\tau_1$ is in the bounded connected component of $\tau_2$. Let $A_1$ be the bounded connected component of $\R^n\setminus\sigma_1$ and $B_1$ the bounded connected component of $\R^n\setminus\tau_1$; $\Sigma$ be the connected component of $\R^n\setminus(\sigma_1\cup\sigma_2)$ having $\sigma_1 \cup \sigma_2$ as boundary and $\Lambda$ the connected component of $\R^n\setminus(\tau_1\cup\tau_2)$ having $\tau_1 \cup \tau_2$ as boundary; $A_2$ be the unbounded connected component 
of $\R^n\setminus\sigma_2$ and $B_2$ the unbounded connected component of $\R^n\setminus\tau_2$ (see Figure \ref{dessin-extension}).

Suppose that $\mathrm{Leb}(A_1) = \mathrm{Leb}(B_1)$ and $\mathrm{Leb}(\Sigma) = \mathrm{Leb}(\Lambda)$. Let $f_i : A_i\to B_i$ be two measure-preserving homeomorphisms such that either each one preserves the orientation, or each one reverses it. Then there exists a measure-preserving homeomorphism $f : \R^n \to \R^n$ whose restriction to $A_1$ equals $f_1$ and restriction to $A_2$ equals $f_2$.
\end{theoreme}

Since this theorem is local, it can be applied to an open space $O$ instead of $\R^n$ or, even better, together with the Oxtoby-Ulam theorem (Theorem \ref{Brown-mesure}), to any domain of chart of a manifold $X$ instead of $\R^n$ and measure $\lambda$ instead of Lebesgue measure.

By perturbing a homeomorphism, we can make every periodic point persistent, as stated by the following lemma.

\begin{lemme}\label{lemper}
Let $f\in \Hom(X,\lambda)$, $\varep>0$, $\delta>0$ and $P\in\N$. Then there exists $g\in\Hom(X,\lambda)$ such that $d(f,g)<\delta$ and such that for every $f$-periodic point $x$ of period smaller than $P$, there exists a persistent periodic point of $g$ with the same period which $\varep$-shadows the $f$-orbit of $x$. 
\end{lemme}

\begin{proof}[Proof of Lemma \ref{lemper}]
Indeed, by compactness of the set of compact subsets of $X$ endowed with Hausdorff topology, it suffices to prove it only for a finite number of $f$-periodic orbits. If for each of these periodic orbits, we perform the perturbation of $f$ outside of a neighbourhood of the other periodic orbits, we reduce the proof to the case of a single periodic orbit. More formally, we choose a neighbourhood $D$ of $x$ such that the sets $D$, $f(D),\cdots, f^{p-1}(D)$ are pairwise disjoints (and disjoint from a neighbourhood of the others periodic orbits) and we locally replace $f$ by the map $h\circ f^{-(p-1)}$ (where $h$ is one of the two maps of Example \ref{linéaire}, depending of whether $f^{-(p-1)}$ 
preserves orientation or not) in the neighbourhood of $f^{-1}(x)$, using the theorem of local modification (Theorem \ref{extension-sphères}), such that $f$ do not change outside the union of the sets $f^i(D)$. 
\end{proof}

\begin{proof}[Proof of Theorem \ref{corovar2}]
We choose a basis $\{C_k\}$ of the topology of $X$. Now, for $p,N,k\in\N$, we define the set $\mathcal{U}_{p,N,k}$ as the set of $f\in\Hom(X,\lambda)$ such that if $f$ possesses a persistent periodic point of period $p$ which intersects $C_k$, then $f_N$ has a unique injective orbit whose periodic orbit intersects $C_k$ and has length $p$. Note that the fact that we consider \emph{persistent} periodic points ensures that these sets are open.

Then, the combination of Proposition \ref{propvar2}, Lemma \ref{lemmetrans} and Lemma \ref{lemper} implies that the set
\[\bigcap_{p\in\N^*}\bigcap_{M\in\N}\bigcup_{N\ge M}\mathcal{U}_{p,N}\]
is a dense $G_\delta$ of $\Hom(X,\lambda)$, consisting of the homeomorphisms that verify the conclusions of the theorem.
%
\end{proof}

\subsection{Invariant measures}

We now try to obtain information about invariant measures of a generic homeomorphism from invariant measures of its discretizations. More precisely, given all the invariant measures of discretizations of a generic homeomorphism, what can be deduced about invariant measures of the initial homeomorphism? A first step in this study was taken in 2006 by T. Miernowski in \cite[Part 8]{MR2279269}:

\begin{prop}\label{Miern}
Let $f\in\Hom(X)$. For all $N\in\N$ let $\gamma_N\subset E_N$ be a periodic cycle of $f_N$ and $\nu_N$ the uniform probability measure on $\gamma_N$. Then, any limit point of the sequence of measures $(\nu_N)_N$ is $f$-invariant. 

In particular, if $f$ is uniquely ergodic, whose unique invariant probability measure is denoted by $\mu$, then the sequence of measures $(\nu_N)_N$ tends to $\mu$.
\end{prop}

The proof of this proposition essentially consists in an appropriate use of the compactness of the set of probability measures on $X$. We now state the same kind of result for \emph{generic} homeomorphisms.

\begin{theoreme}\label{EnsMesInvSimpl}
Let $f\in\Hom(X,\lambda)$ be a generic homeomorphism and suppose that the sequence of grids $(E_N)_{N\in\N}$ is sometimes strongly self similar. Then for every $f$-invariant measure $\mu$, there exists a subsequence $\nu_{N_k}$ of $f_{N_k}$-invariant periodic measures such that $(\nu_{N_k})_k$ tends to $\mu$.
\end{theoreme}

This result is a particular case of the more general Theorem \ref{EnsMesInv}, for which we will need the following notation.

\begin{notation}
For $f\in\Hom(X)$, we denote by $\mathcal M^f$\index{$\mathcal M^f$} the set of Borel probability measures on $X$ which are invariant under $f$, and by $\mathcal M^f_N$\index{$\mathcal M^f_N$} the set of Borel probability measures on $E_N$ which are invariant under $f_N$
\end{notation}


\begin{theoreme}\label{EnsMesInv}
Let $f\in\Hom(X,\lambda)$ be a generic homeomorphism and suppose that the sequence of grids $(E_N)_{N\in\N}$ is sometimes strongly self similar. Let $\mathcal M_N$\index{$\mathcal M_N$} be the set of probability measures on $E_N$ that are invariant under $f_N$. Then for every closed and convex subset $\mathcal M$ of $f$-invariant measures, there exists a subsequence of $\{\mathcal M_N\}$ which tends to $\mathcal M$ for Hausdorff topology.
\end{theoreme}

Proposition \ref{Miern} asserts that the upper limit of the sets of measures that are invariant under $f_N$ is included in the set of $f$-invariant measures. As the set of $f_N$-invariant measures is always convex, Theorem~\ref{EnsMesInv} asserts that generically, the sets of $f_N$-invariant measures accumulate on all the possible sets of measures.

Before giving a detailed proof of Theorem \ref{EnsMesInv}, let us give its main arguments. An \emph{ad hoc} application of Baire's theorem, together with compactness arguments, reduces the proof to that of the following variation of Lax's theorem, which asserts that  for every finite collection $\mathcal M$ of measures that are invariant under $f$, there exists a homeomorphism $g$ which is close to $f$ and a big order of discretization $N$ such that the set of periodic measures that are invariant under $g_N$ is close to $\mathcal M$ for Hausdorff distance.

\begin{lemme}[Ergodic variation of Lax's theorem]\label{Laxergod}
Suppose that the sequence of grids $(E_N)_{N\in\N}$ is sometimes strongly self similar. Then, for all $f\in\Hom(X,\lambda)$, for all collection of $f$-invariant measures $\nu_1,\cdots,\nu_\ell$, for all $\varep>0$ and all $k_0,N_0\in\N$, there exists $g\in\Hom(X,\lambda)$ and $N\ge N_0$ such that $d(f,g)<\varep$, and that $g_N$ possesses exactly $\ell$ periodic orbits, and the the invariant periodic measures $\{\nu_{1,N}^g,\cdots, \nu_{\ell,N}^g\}$ supported by these periodic orbits satisfy $\dist(\nu_i,\nu_{i,N}^g)<\frac{1}{k_0}$ for all $i$. Moreover, we can suppose that these properties are still true on a whole neighbourhood of $g$.
\end{lemme}

For the proof of Theorem \ref{EnsMesInv}, we will need a slightly weaker result; however this lemma will also be useful in the next part concerning the physical dynamics.

Before giving the formal proof of Lemma \ref{Laxergod}, we sketch its main arguments. Suppose first that $\ell=1$ and that $\nu_1$ is ergodic. For this purpose we apply Birkhoff's theorem to $f$, $\nu_1$ and a recurrent point $x$: for all $M$ large enough,
\[\frac{1}{M}\sum_{k=0}^{M-1}\delta_{f^k(x)} \simeq \nu_1.\]
Since $x$ is recurrent we can choose an integer $M$ large enough such that $x\simeq f^M(x)$. First we approximate $f$ by a cyclic permutation $\sigma_N$ given by Lax's theorem, then we slightly modify $\sigma_N$ into a map $\sigma'_N$ by choosing $\sigma_N'(\sigma_N^M(x_N)) = x_N$, as in Proposition \ref{var1}. The measure $\nu^{\sigma_N'}_N$ is the uniform measure on the orbit $x_N,\cdots, \sigma_N^{M-1}(x_N)$, so it is close to $\nu_1$; we then just have to apply the proposition of extension of finite maps to finish the proof in the case where $\ell=1$ and $\nu_1$ ergodic. The proof of the lemma when $\nu_1$ is not ergodic but only invariant is obtained by approximating the invariant measure by a finite convex combination of ergodic measures. We use the hypothesis of self similarity of the grids to apply the previous process on each subgrid; we approximate each ergodic measure of the convex combination on a number of subgrids proportional to the coefficient of this measure in the convex combination. We 
then make a final perturbation to merge the periodic orbits of the different subgrids and apply the proposition of extension of finite maps; this proves the case $\ell=1$. Finally, the general case (when there are several measures to approximate) is obtained is using (again) the fact that the grids are sometimes strongly self similar.

\begin{proof}[Proof of Theorem \ref{EnsMesInv}]
Recall that we denote by $\Prb$ the set of all Borel probability measures on $X$ and that $\Prb$ is equipped with a distance $\dist$ compatible with the weak-* topology. The set of compact sets of Borel probability measures on $X$ is equipped with the Hausdorff distance $\dist_H$ which makes this set compact; we consider a sequence $\{\mathcal M_\ell\}_\ell$ of compact sets of Borel probability measures which is dense for this topology. Thereafter homeomorphisms will be taken in the generic set
\[\bigcap_{N\in\N} \mathcal{D}_{N},\]
made of homeomorphisms whose $N$-th discretization is uniquely defined for all $N$. Consider
\[\mathcal{A} = \bigcap_{(\ell,N_0,k_0)\in\N^3} \mathcal{O}_{\ell,N_0,k_0},\]
where $\mathcal{O}_{\ell,N_0,k_0}$ is the set of homeomorphisms $f\in \bigcap_{N\in\N} \mathcal{D}_{N}$ such that if there exists a closed convex set $\mathcal M^f$ of $f$-invariant measures such that $\dist_H(\mathcal M_\ell,\mathcal M^f)\le\frac{1}{k_0}$, then there exists $N\ge N_0$ such that the set $\mathcal M^f_N$ of $f_N$-invariant measures satisfies $\dist_H(\mathcal M_\ell,\mathcal M^f_N)<\frac{2}{k_0}$.
\bigskip

We easily check that if $f\in \mathcal{A}$, then it satisfies the conclusions of the theorem.

The fact that the sets $\mathcal{O}_{\ell,N_0,k_0}$ are open follows from the upper semi-continuity of $f\mapsto \mathcal M^f$ and the fact that in the set $\bigcap_{N\in\N} \mathcal{D}_{N}$, the map $g\mapsto \mathcal M_N^g$ is locally constant.

It remains to show that the sets $\mathcal{O}_{\ell,N_0,k_0}$ are dense, it follows from Lemma~\ref{Laxergod} in the following way. Let $f\in \bigcap_{N\in\N} \mathcal{D}_{N}$ and $\ell$, $N_0$ and $k_0$. If $\dist_H(\mathcal M_\ell,\mathcal M^f)>\frac{1}{k_0}$, there is nothing to prove. So we suppose that $\dist_H(\mathcal M_\ell,\mathcal M^f)\le\frac{1}{k_0}$. By compactness, there exists a finite collection of $f$-invariant measures $\nu_1,\cdots,\nu_{\ell}$ such that $\mathcal M^f\subset \bigcup_i B\big(\nu_i,\frac{1}{k_0}\big)$ (where $B\big(\nu_i,\frac{1}{k_0}\big)$ denotes the set of measures whose distance to $\nu_i$ is smaller than $1/k_0$). So it suffices to find $g$ close to $f$ such that the set of $g_N$-invariant measures is included in $\bigcup_i B\big(\nu_i,\frac{1}{k_0}\big)$, but this fact is implied by Lemma \ref{Laxergod}.
\end{proof}

\begin{proof}[Proof of Lemma \ref{Laxergod}]
To begin with, we prove the lemma for only one measure $\nu=\nu_1$ which is ergodic. We want to show that there exists a homeomorphism $g$ whose distance to $f$ is smaller than $\varep$, and an integer $N\ge N_0$ such that $\dist(\nu,\mu^{g_N}_X)<\frac{2}{k_0}$.

Since $\nu$ is ergodic, for all continuous map $\varphi$, by Birkhoff's theorem,
\begin{equation}\label{eqergo}
\frac{1}{M}\sum_{m=0}^{M-1} \delta_{f^m(x)}\underset{M\to+\infty}{\longrightarrow} \nu
\end{equation}
for $\nu$-almost every $x$. Let $x\in X$ be a recurrent point for $f$ satisfying equation (\ref{eqergo}) (such points form $\nu$-full measure set). There exists $M_0\in\N$ such that for all $M\ge M_0$ we have
\begin{equation}\label{eqergo2}
\dist\left(\frac{1}{M}\sum_{m=0}^{M-1}\delta_{f^m(x)},\, \nu \right)\le \frac{1}{2k_0}.
\end{equation}

Since $x$ is recurrent, there exists $\tau\ge M_0$ such that $d(x,f^\tau(x))\le\varep/4$. Let $\sigma_N$ be a map from $E_N$ into itself given by Lax's theorem: it is a cyclic permutation and its distance to $f$ is smaller than $\varep/2$. For all $N$ large enough the orbit $(\sigma_{N}^m(x_{N}))_{0\le m\le \tau}$ shadows the orbit $(f^m(x))_{0\le m\le \tau}$, thus $d(x_{N},\sigma_N^\tau(x_{N}))<\varep/2$. Then, we ``close'' the orbit of  $x_{N}$ between the points $x_{N}$ and $\sigma_{N}^\tau(x_{N})$, \emph{i.e.} we set (as in Proposition \ref{var1}, see also Figure \ref{trajectoire})
\[\sigma'_{N}(y_{N}) = \left\{\begin{array}{ll}
x_{N}\quad & \text{if}\ \,y_{N} = \sigma_{N}^{\tau-1}(x_N)\\
\sigma_{N}(y_N)\quad & \text{otherwise.} 
\end{array}\right.\]
Then $d_{N}(\sigma'_{N},f)<\varep$ and $(x_N,\cdots,\sigma_{N}'^{\tau-1}(x_{N}))$ is a periodic orbit for $\sigma'_N$ whose basin of attraction is the whole set $E_N$. If $\varep$ is small enough, then
\begin{equation}\label{machintruc}
\dist\left(\frac{1}{\tau}\sum_{m=0}^{\tau-1}\delta_{f^m(x)},\, \frac{1}{\tau}\sum_{m=0}^{\tau-1}\delta_{{\sigma'}_N^m(x_N)} \right)\le \frac{1}{2k_0}.
\end{equation}

Since the periodic orbit $(\sigma_{N}'^m(x))_{0\le m\le \tau}$ attracts $E_{N}$, for all $M'$ large enough and $y_{N}\in E_{N}$ we have
\begin{equation}\label{eqergo4}
\frac{1}{M'}\sum_{m=0}^{M'-1}\delta_{\sigma_{N}'^m(y_{N})}\underset{M'\to\infty}{\longrightarrow}  \frac{1}{\tau}\sum_{m=0}^{\tau-1}\delta_{ \sigma_{N}'^m(x_{N})}.
\end{equation}
With Proposition \ref{extension}, the same way as in Lemma \ref{lemmetrans}, we construct a homeomorphism $g'$ from the map $\sigma'_N$ such that the discretization $g'_N$ and $\sigma'_N$ fit together, and such that $g'$ and the restriction of $g'_N$ to the orbit of $x_N$ fit together. Since $d_{N}(f,\sigma_{N}')<\varep$, we can furthermore assume that $d(f,g')<\varep$. Thus, Equation (\ref{eqergo4}) implies that
\begin{equation}\label{defmug}
\mu^{g'_N}_X = \frac{1}{\tau}\sum_{m=0}^{\tau-1}\delta_{\sigma'^m_N(x_N)}.
\end{equation}

We now have all the necessary estimations to compute the distance between $\mu^{g'}_N$ and $\nu$: applying Equations \eqref{eqergo2}, \eqref{machintruc}, \eqref{defmug} and triangle inequality, we obtain:
\[\dist(\mu^{g'_N}_X,\nu)<\frac{1}{k_0}.\]
\bigskip

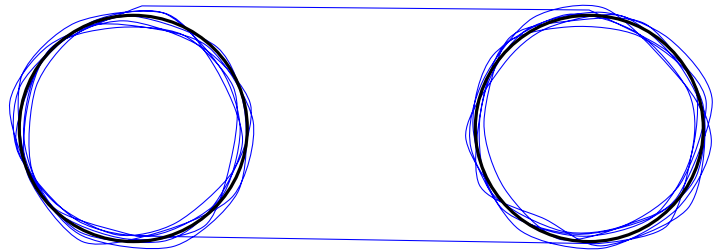
\begin{figure}[b]
\begin{center}
\begin{tikzpicture}[scale=1.5]
\draw [blue,domain=90:7*360-90, samples=61, smooth]  plot ({cos(\x)+.1*rand},{sin(\x)+.1*rand}) -- plot ({4+cos(\x-180)+.1*rand},{sin(\x-180)+.1*rand}) -- cycle;
\draw[very thick] (0,0) circle (1);
\draw[very thick] (4,0) circle (1);
\end{tikzpicture}
\caption[Construction of Lemma \ref{Laxergod}]{Construction of Lemma \ref{Laxergod}: the measure carried by the blue cycle is close to a convex combination of the measures carried by the black cycles}\label{FigrAlea}
\end{center}
\end{figure}

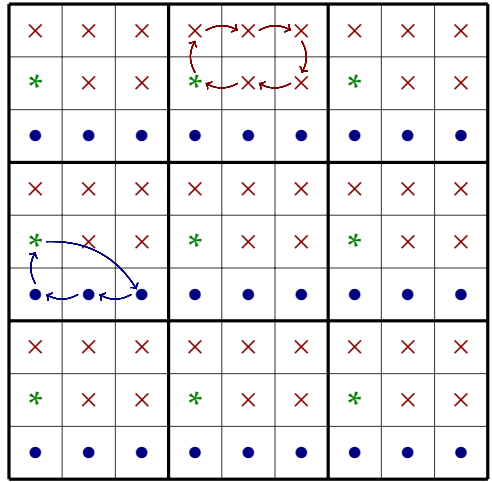
\begin{figure}[t]
\begin{center}
\begin{tikzpicture}[scale=.7]
\draw[gray!50!black] (0,0) grid (9,9);
\draw[step=3cm, very thick] (0,0) grid (9,9);
\foreach\i in {0,...,2}{
\foreach\j in {0,...,2}{
\draw[blue!50!black] (3*\i+.5,3*\j+.5) node {$\bullet$};
\draw[blue!50!black] (3*\i+1.5,3*\j+.5) node {$\bullet$};
\draw[blue!50!black] (3*\i+2.5,3*\j+.5) node {$\bullet$};
\draw[green!50!black] (3*\i+.5,3*\j+1.5) node {\Large$\ast$};
\draw[red!50!black] (3*\i+1.5,3*\j+1.5) node {$\times$};
\draw[red!50!black] (3*\i+2.5,3*\j+1.5) node {$\times$};
\draw[red!50!black] (3*\i+.5,3*\j+2.5) node {$\times$};
\draw[red!50!black] (3*\i+1.5,3*\j+2.5) node {$\times$};
\draw[red!50!black] (3*\i+2.5,3*\j+2.5) node {$\times$};
\draw[->, red!50!black] (3.7,8.5) to[bend left](4.3,8.5);
\draw[->, red!50!black] (4.7,8.5) to[bend left](5.3,8.5);
\draw[->, red!50!black] (5.5,8.3) to[bend left](5.5,7.7);
\draw[->, red!50!black] (5.3,7.5) to[bend left](4.7,7.5);
\draw[->, red!50!black] (4.3,7.5) to[bend left](3.7,7.5);
\draw[->, red!50!black] (3.5,7.7) to[bend left](3.5,8.3);
\draw[->, blue!50!black] (.5,3.7) to[bend left](.5,4.3);
\draw[->, blue!50!black] (.7,4.5) to[bend left](2.4,3.6);
\draw[->, blue!50!black] (2.3,3.5) to[bend left](1.7,3.5);
\draw[->, blue!50!black] (1.3,3.5) to[bend left](.7,3.5);
}}
\end{tikzpicture}
\caption[Construction of Lemma \ref{Laxergod}]{Construction of Lemma \ref{Laxergod}: before the perturbation, the red (crosses) subgrids carry a measure close to $\nu_1^e$, the blue (points) subgrids carry a measure close to $\nu_1^e$ and the green (stars) subgrid carries a cyclic permutation; we compose $\sigma'_N$ by two cyclic permutations (red and blue arrows) to create a long periodic orbit which attracts the whole grid $E_N$ and carries a measure close to $\nu'$}\label{ReunCycl}
\end{center}
\end{figure}

In the general case the measure $\nu$ is only invariant (and not ergodic). So in the second step of the proof we treat the case where there is only one measure that is not ergodic. It reduces to the ergodic case by the fact that the set of invariant measures is a compact convex set whose extremal points are exactly the ergodic measures: by the Krein-Milman theorem, for all $M\ge 1$ there exists an $f$-invariant measure $\nu'$ which is a finite convex combination of ergodic measures: 
\[\nu' = \sum_{j=1}^{r} \lambda_j\nu_j^e,\]
and whose distance to $\nu$ is smaller than $\frac{1}{k_0}$. To simplify the proof, we treat the case where $r=2$, the other cases being treated the same way. So $\nu'=\lambda_1\nu_1^e + \lambda_2\nu_2^e$. We use the hypothesis of self similarity of the grids to build a permutation of a grid $E_N$ which is close to $f$ and has a unique periodic orbit whose associated measure is close to $\nu'$ in the way described by Figure \ref{FigrAlea}. More precisely,  the self similarity of the grids imply that for all $\varep>0$, there exists $N_0\in\N$ such that for all $N\ge N_0$, the set $E_N$ is partitioned by disjoint subsets $\widetilde E_N^1, \cdots, \widetilde E_N^{P}$ such that for all $i$, $\widetilde E_N^i$ is the image of the grid $E_{N_0}$ by a bijection $h_i$ which is $\varep$-close to identity; in particular we have $P = q_N/q_{N_0}$. Then, Lax's theorem (Theorem \ref{Lax}) imply that there exists a permutation $\sigma_N$ whose distance to $f$ is smaller than $\varep$ and which is a cyclic permutation on 
each set $\widetilde E_N^i$. We denote by $P$ the number of such cycles and approximate the $\lambda_i$ by multiples of $1/P$:
\[\lambda'_1 = \frac{\lfloor\lambda_1 P\rfloor}{P}\]
and
\[\lambda'_2 = 1-\lambda'_1-\frac{1}{P}.\]
Then we have $|\lambda_1-\lambda'_1|\le \frac{1}{P}$ and $|\lambda_2-\lambda'_2|\le \frac{3}{P}$. Increasing $N_0$ if necessary, we make two different modifications of the map $\sigma_N$ on the sets $\widetilde E_N^i$ to obtain a map $\sigma'_N$:
\begin{enumerate}
\item on the sets $\widetilde E_N^i$ with $1\le i\le P\lambda'_1$, we do the previous construction concerning ergodic measures: on each of these sets, ${\sigma'_N}_{|\widetilde E_N^i}$ carries an unique invariant probability measure ${\nu_i^e}'$ satisfying $\dist({\nu_i^e} , {\nu_1^e})\le\frac{1}{10k_0}$;
\item on the sets $\widetilde E_N^i$ with $2+P\lambda'_1\le i\le P$, we do the previous construction concerning ergodic measures: on each of these sets, ${\sigma'_N}_{|\widetilde E_N^i}$ carries an unique invariant probability measure ${\nu_i^e}'$ satisfying $\dist({\nu_i^e} , {\nu_2^e})\le\frac{1}{10k_0}$.
\end{enumerate}

It remains to merge the periodic orbits on the sets $\widetilde E_N^i$ to create a periodic orbit which carries a measure which is close to $\nu'$ (see Figure \ref{ReunCycl}). We first consider a point $x_1\in E_N$ such that $x_1$ belongs to the periodic orbit of ${\sigma'_N}_{|\widetilde E_N^1}$, and compose $\sigma'_N$ by a cyclic permutation of the set $\{h_i(h_1^{-1}(x_1))\}_{1\le i\le P(\lambda'_1+1)}$. Then, we consider a point $x_2\in E_N$ such that $x_2$ belongs to the periodic orbit of ${\sigma'_N}_{|\widetilde E_N^P}$, and compose $\sigma'_N$ by a cyclic permutation of the set $\{h_i(h_1^{-1}(x_2))\}_{P(\lambda'_1+1)\le i\le P}$.

Thus, by construction, the obtained map $\sigma''_N$ carries a single periodic orbit (thus this orbit attracts the whole grid $E_N$); we denote by $\mu^{\sigma''_N}$. Moreover, if we consider a point $x_N$ belonging to this orbit, its orbit distributes like the measure $\nu_1$ during a time proportional to $\lambda'_1$, like the measure $\nu_2$ during a time proportional to $\lambda'_2$, whereas it covers the rest of the orbit during a (small) time proportional to $1/P$. Thus, increasing $N$ (therefore $P$) if necessary, it means that the distance between the measures $\mu^{\sigma''_N}$ and $\nu'$ is smaller than $1/k_0$. To finish the proof of the lemma when there is only one measure to approximate, it suffices to interpolate the map $\sigma''_N$ by a homeomorphism close to $f$ in applying the proposition of extension of finite maps (Proposition \ref{extension}).
\bigskip

To treat the case where there are several measures to approximate, we use once more the fact that the grids are sometimes strongly self similar. On each of the subgrids, (and considering sub-subgrids if necessary), we apply the process of the previous part of the proof, and apply the proposition of extension of finite maps to find a homeomorphism which satisfies the conclusions of the lemma.
\end{proof}

\subsection{Invariant compact sets}

There is a similarity between invariant measures and invariant compact sets, thus the previous theorem is also true for invariant compact sets. As in the previous section, we begin by a particular case of the main theorem of this section (Theorem \ref{CompactInv}).

\begin{theoreme}\label{CompactInvSimpl}
Let $f\in\Hom(X,\lambda)$ be a generic homeomorphism and suppose that the sequence of grids $(E_N)_{N\in\N}$ is sometimes strongly self similar. Then for every $f$-invariant chain transitive compact set $K\subset X$, there exists a subsequence $K_{N_k}$ of $f_{N_k}$-periodic orbits which tends to $K$ for Hausdorff topology.
\end{theoreme}

In the ergodic case, the set of invariant measures of a given map had to be closed and convex; concerning compact sets, there are also conditions for a collection of compact sets to be the set of compact invariant sets of a map. This motivates the following definition.

\begin{definition}
Let $f\in\Hom(X)$. We say that a set $\mathcal K$ of compact subsets of $X$ that are invariant under $f$ is \emph{admissible for $f$} if there exists $\mathcal K_0\subset \mathcal K$ such that:
\begin{enumerate}[(i)]
\item every $K_0\in \mathcal K_0$ is chain transitive;
\item every $K\in\mathcal K$ is a finite union of elements of $\mathcal K_0$.
\end{enumerate}
\end{definition}

\begin{theoreme}\label{CompactInv}
Let $f\in\Hom(X,\lambda)$ be a generic homeomorphism and suppose that the sequence of grids $(E_N)_{N\in\N}$ is sometimes strongly self similar. Let $\mathcal K_N$ be the set of compact subsets of $E_N$ that are invariant under $f_N$. Then for every collection $\mathcal K$ of compact subsets of $X$ that is admissible for $f$, there exists a subsequence of $\mathcal K_N$ which tends to $\overline{\mathcal K}$ for Hausdorff topology\footnote{That is, Hausdorff topology on the sets of compacts subsets of $X$ endowed with Hausdorff distance.}.
\end{theoreme}

Like in the ergodic case, this theorem expresses that the sets of invariant compact sets for $f_N$ accumulate on all the admissible sets of invariant sets for~$f$. 

The proof of Theorem \ref{CompactInv} is very similar to that of Theorem \ref{EnsMesInv}. An appropriate application of Baire's theorem reduces the proof to that of the following variation of Lax's theorem.

\begin{lemme}[Compact variation of Lax's theorem]\label{CompacLax}
Suppose that the sequence of grids $(E_N)_{N\in\N}$ is sometimes strongly self similar. For all $f\in\Hom(X,\lambda)$, for all collection $\mathcal K$ of compact subsets of $X$ that is admissible for $f$, for all $\varep>0$ and $k_0,N_0\in\N$, there exists $g\in\Hom(X,\lambda)$ and $N\ge N_0$ such that $d(f,g)<\varep$, and that the set of $g_N$-invariant sets\footnote{Which consists of unions of periodic orbits of $g_N$.} on $E_N$ is $1/k_0$-close to $\overline{\mathcal K}$ for Hausdorff distance. Moreover, we can suppose that these properties are satisfied on a whole neighbourhood of $g$.
\end{lemme}

\begin{proof}[Sketch of proof of Lemma \ref{CompacLax}]
The proof is very similar to that of Lemma~\ref{Laxergod}. We take a collection $\mathcal K$ of compact subsets of $X$ which is admissible for $f$. By compactness of $\overline{\mathcal K}$, it is close for Hausdorff topology to a finite set $\{K_1\}\cup\cdots\cup\{K_p\}$. Then, every compact set $K_i$ is close to a finite union $\tilde K_i^1\cup\cdots\cup \tilde K_i^{\ell_i}$ of chain transitive compact sets that are invariant under $f$. As each set $\tilde K_i^j$ is chain transitive, there exists a $f$-pseudo orbit $\tilde \omega_i^j$ which is close to $\tilde K_i^j$ for Hausdorff distance. We then apply the same proof strategy as for Lemma~\ref{Laxergod} to the pseudo orbits $\tilde \omega_i^j$ to obtain directly the conclusions of Lemma~\ref{CompacLax}.
\end{proof}

\section{Physical dynamics}\label{grobra}

In the previous section, we proved that generically, the upper limit of the sets of invariant measures of discretizations is the set of invariant measures of the initial homeomorphism. This expresses that the sets of invariant measures of discretizations accumulate on ``the biggest possible set''. However, we might expect that physical measures -- that is, Borel measures which attracts a lot of points with respect to $\lambda$ (see Definition \ref{sport} page \pageref{sport}) -- play a specific part: their definition expresses that they are the measures that can be observed in practice on numerical experiments, because they governs the ergodic behaviour of $\lambda$ almost every point\footnote{Recall that $\lambda$ is a good measure, in particular it is positive on every non-empty open set.}; in the case of a generic conservative homeomorphism, Oxtoby-Ulam theorem \cite{Oxto-meas} implies that $\lambda$ is the unique physical measure. Moreover, results of stochastic stability are known to be true in 
various contexts (for example, expanding maps \cite{MR884892},\cite{MR874047}, \cite{MR685377}, uniformly hyperbolic attractors \cite{MR874047}, \cite{MR857204}\dots). These theorems suggest that the physical measures can always be observed in practice, even if the system is noisy. Given this background, we are tempted to think that the natural invariant measures $\mu^{f_N}_X$ of $f_N$, which can be seen as the physical measures of $f_N$, converge to the physical measures of $f$.

Recall that $\mu^{f_N}_X$ is the limit in the sense of Cesàro of the pushforwards by iterates of $f_N$ of uniform measure $\lambda_N$\index{$\lambda_N$} on $E_N$ (see Definition \ref{defmes} page \pageref{defmes}):
\[\mu^{f_N}_X = \lim_{m\to\infty}\frac 1m \sum_{i=0}^{m-1}(f_N)_*^i \lambda_N.\]
The measure $\mu^{f_N}_X$ is supported by the recurrent set of $f_N$; it is uniform on every periodic orbit and the total weight of a periodic orbit is proportional to the size of its basin of attraction.

The expectation that the measure $\lambda$ plays a specific part is supported by the following variation of Proposition~\ref{Miern}, obtained by replacing $\nu_N$ by $\mu^{f_N}_X$: \emph{if $f$ is uniquely ergodic, then the measures $\mu^{f_N}_X$ converge weakly to the only measure $\mu$ that is invariant under $f$}.

Unfortunately, we show that this is not at all the case: the sequence of measures $(\mu^{f_N}_X)$ accumulates of the whole set of $f$-invariant measures. More precisely we have the following theorem:

\begin{theoreme}\label{mesinv}
If the sequence of grids $(E_N)_{N\in\N}$ is sometimes strongly self similar, then for a generic homeomorphism $f\in\Hom(X,\lambda)$, the set of limit points of the sequence $(\mu^{f_N}_X)_{N\in\N}$ is exactly the set of $f$-invariant measures.
\end{theoreme}

This theorem can be seen as a discretized version of the following conjecture.

\begin{conj}[F. Abdenur, M. Andersson, \cite{MR3027586}]
A homeomorphism $f$ which is generic in the set of homeomorphisms of $X$ (without measure preserving hypothesis) that do not have any open trapping set is \emph{wicked}, i.e. it is not uniquely ergodic and the measures
\[\frac{1}{m}\sum_{k=0}^{m-1}f_*^k(\Leb)\]
accumulate on the whole set of invariant measures under~$f$.
\end{conj}

The behaviour described in this conjecture is the opposite of that consisting of possessing a physical measure.

\begin{proof}[Sketch of proof of theorem \ref{mesinv}]
The proof is similar to which of theorem \ref{EnsMesInv}: the set $\mathcal A$ is replaced by
\[\mathcal{A}' = \bigcap_{(\ell,N_0,k_0)\in\N^3} \mathcal{O}_{\ell,N_0,k_0},\]
where
\[\mathcal{O}_{\ell,N_0,k_0} = \left\{f\in \bigcap_{N\in\N} \mathcal{D}_{N}\  \middle\vert\
\begin{array}{l}
\Big(\exists\nu f\text{-inv.}: \dist(\nu,\tilde\nu_\ell)\le\frac{1}{k_0}\Big)\implies\\
\Big(\exists N\ge N_0 : \dist(\tilde\nu_\ell,\mu^{f_N}_X)<\frac{2}{k_0}\Big)
\end{array}\right\},\]
and $\{\nu_\ell\}$ is a countable set of Borel probability measures that is dense in the whole set of Borel probability measures. A direct application of  Lemma \ref{Laxergod} leads to the fact that every set $\mathcal{O}_{\ell,N_0,k_0}$ is open and dense.
\end{proof}

We have also another corollary of Lemma \ref{Laxergod} about physical dynamics, which can be seen as a combination of Theorem \ref{mesinv} and of a discrete Birkhoff's theorem.

\begin{theoreme}\label{MesPhys2}
If the sequence of grids $(E_N)_{N\in\N}$ is sometimes strongly self similar, then for a generic homeomorphism $f\in\Hom(X,\lambda)$ and for every $f$-invariant measure $\mu$, there exists a subsequence $f_{N_k}$ of discretizations such that for every $x\in X$, the sequence of measures\footnote{Recall that by Definition \ref{defmes}, $\mu^{f_N}_x$ is the Cesàro limit of the pushforwards of the Dirac measure $\delta_{x_N}$ by the discretization $f_N$.} $\mu^{f_{N_k}}_x$ tends to $\mu$.
\end{theoreme}

In particular, this theorem asserts that it is impossible to detect the physical measure of a generic conservative homeomorphism by looking at the physical measures of its discretizations: 

\begin{rem}
This seems to contradict the empirical observations made by A. Boyarsky in 1986 (see \cite{MR862028} or \cite{MR959419}): when a homeomorphism $f$ has only one ergodic measure $\mu$ which is absolutely continuous with respect to Lebesgue measure, ``most of'' the measures $\mu^{f_N}_{x}$ tend to measure $\mu$. However, the author does not specify in what sense he means ``most of the points'', or if his remark is based on a tacit assumption of regularity for $f$.
\end{rem}

Note that as in the previous section, we have a compact counterpart of Theorem~\ref{mesinv}:

\begin{theoreme}\label{EquivCompact}
If the sequence of grids $(E_N)_{N\in\N}$ is sometimes strongly self similar, then for a generic homeomorphism $f\in\Hom(X,\lambda)$, the set of limit points of the sequence $(\Omega(f_{N}))_{N\in\N}$ is exactly the set of limit points of finite unions of $f$-invariant chain transitive compact sets.
\end{theoreme}

\section{Addendum: generic conjugates and generic grids}\label{AddendParti1}

In this addendum we want to adapt the proofs of this chapter to the case of discretizations of conjugated generic homeomorphisms. More precisely, we fix a conservative homeomorphism $f$ and we look at the dynamics of the discretizations of $hfh^{-1}$, where $h$ is a generic conservative homeomorphism. This is equivalent to fix a dynamics $f$ and to wonder what is the dynamics of the discretizations of generic realizations of $f$. Another point of view is to see these discretizations as discretizations of a fixed homeomorphism on a generic sequence of grids, where by definition a generic sequence of grids is the image of a good sequence of grids by a generic homeomorphism.

The following property is a finite maps extension result for conjugations of a given homeomorphism. It implies that under some weak hypothesis, the results are the same for generic conjugations of a fixed homeomorphism and for generic conservative homeomorphisms.

\begin{prop}\label{seule}
Let $f$ be a conservative homeomorphism whose set of fixed points has empty interior, and $\sigma : E_N\to E_N$ be such that $d_N(f,\sigma)<\varep$. Then there exists a conservative homeomorphism $h$ such that $d(h,\operatorname{Id})<\varep$ and such that the homeomorphism $g = hfh^{-1}$satisfies $g_N = \sigma$ (and these properties remain true on a whole neighbourhood of $h$).
\end{prop}

\begin{proof}[Proof of Proposition \ref{seule}]
Let $\delta = \min\{1/(2N),\varep-d(h,\operatorname{Id})\}$. Let $x_1,\cdots,x_{q_N}$ be the points of $E_N$. It may happen that some of these points are fixed points of $f$, thus we define an alternative sequence of points $(x'_1,y'_1,\cdots,x'_{q_N},y'_{q_N})$ by induction, which is close to the sequence $(x_1,f(x_1),\cdots,x_{q_N},f(x_{q_N}))$ and contains no fixed points. Suppose that we have constructed the sequence until $x'_{\ell-1}$ and, $y'_{\ell-1}$. Then there exists $x'_\ell\in X$ $\delta$-close to $x_\ell$, which is not a fixed point of $f$ and which is different from all the points $x'_1,f(x'_1),f^{-1}(x'_1),\cdots,x'_{\ell-1},f(x'_{\ell-1}),f^{-1}(x'_{\ell-1})$. Similarly, there exists a point $y'_\ell$ $\delta$-close to $\sigma(x_\ell)$ and different from all the points $y_1,\cdots,y_{\ell-1}$ and $x_1,\cdots,x_{q_N}$. So on, we define a set $(x'_1,y'_1,\cdots,x'_{q_N},y'_{q_N})$.

We can now build the homeomorphism $h$ with the proposition of extension of finite maps (Proposition \ref{extension}): we choose a conservative homeomorphism $h$ such that for all $\ell$, we have $h(x'_\ell) = x_\ell$ and $h(f(x'_\ell)) = y'_\ell$. As for all $\ell$, the point $x'_\ell$ is $\varep$-close to $x_\ell$, the point $\sigma(x_\ell)$ is $\varep-d(h,\operatorname{Id})$-close to $y'_\ell$ and the point $f(x_\ell)$ is $d(h,\operatorname{Id})$-close to $\sigma(x_\ell)$, we can choose $h$ such that $d(h,\operatorname{Id})<\varep$. Moreover, if we set $g = hfh^{-1}$, then by construction we have $g(x_\ell) = h(f(x'_\ell)) = y'_\ell$. As $y'_\ell$ is $1/(2N)$-close to $\sigma(x_\ell)$, we obtain $g_N(x_\ell) = \sigma(x_\ell)$ for all $\ell$.
\end{proof}

As in this chapter, the only perturbation result we used was the proposition of finite maps extension (Proposition \ref{extension}), this proposition implies the following corollary.

\begin{coro}
All the results stated in this chapter concerning the dynamical behaviour of the discretizations of generic conservative homeomorphisms remain true for the dynamical behaviour of the discretizations of conservative homeomorphisms which are generic among those conjugated to a given conservative homeomorphism whose set of fixed points has empty interior.
\end{coro}

\section{Numerical simulations}\label{partietroisb}

We now present the results of numerical simulations of conservative homeomorphisms. We are interested in whether it is possible to observe the behaviours described by our theoretical results on actual simulations. It is not clear that the orders of discretization involved in these results can be reached in practice, or if simple examples of homeomorphisms behave the same way as generic homeomorphisms. For instance, can we observe a lot of discretizations that are cyclic permutations of the grid, or that have a small degree of recurrence? Given a homeomorphism $f$, is it possible to see the invariant measures of the discretizations accumulating on the set of invariant measures of $f$? Is it possible to recover the periodic points of the homeomorphism? 

We simulate homeomorphisms $f(x,y) = (Q\circ P)(x,y)$, where both $P$ and $Q$ are homeomorphisms of the torus that modify only one coordinate:
\[P(x,y) = \big(x,y+p(x)\big)\quad\text{and}\quad Q(x,y) = \big(x+q(y),y\big),\]
so that $P$ and $Q$ both preserve Lebesgue measure.
We discretize these examples according to the uniform grids on the torus
\[E_N = \left\{\left(\frac{i}{N},\frac{j}{N}\right)\in \T^2 \big|\  0\le i,j\le {N}-1\right\}.\]
We perform simulations of two conservative homeomorphisms which are small perturbations of identity or of the standard Anosov automorphism $A : (x,y)\mapsto (x+y,x+2y)$.\label{PageDefSimulCons}
\begin{itemize}
\item To begin with we study $f_3 = Q\circ P$, with\label{defex}
\[p(x) = \frac{1}{209}\cos(2\pi\times 187x)+\frac{1}{271}\sin(2\pi\times 253x)-\frac{1}{703}\cos(2\pi\times775x),\]
\[q(y) = \frac{1}{287}\cos(2\pi\times 241y)+\frac{1}{203}\sin(2\pi\times 197y)-\frac{1}{841}\sin(2\pi\times811y).\]
This conservative homeomorphism is a small $C^0$ perturbation of identity. Experience shows that even dynamical systems with fairly simple definitions behave quite chaotically (see for example \cite{MR700317}). We can expect that a homeomorphism such as $f_3$ has a complex dynamical behaviour and even more, behaves essentially as a generic homeomorphism, at least for orders of discretization that can be reached in practice. Note that we choose coefficients that have (virtually) no common divisors to avoid arithmetical phenomena such as periodicity or resonance. Also, we have chosen to simulate a homeomorphism close to the identity to avoid phenomenons like ``every orbit of the discretization $f_N$ is $\varep$-dense for all $N$ large enough'': in this case it would be difficult to see something interesting in the images of the invariant measures.
\item We also simulate $f_4$ the composition of $f_3$ with the standard Anosov automorphism
\[A = \begin{pmatrix} 2 & 1 \\ 1 & 1 \end{pmatrix},\]
say $f_4 = Q\circ P\circ A$. This makes it a small $C^0$-perturbation of $A$. Thus $f_4$ is semi-conjugated to $A$ but not conjugated: to each periodic orbit of $A$ corresponds many periodic orbits of $f_4$. As for $f_3$, we define $f_4$ in the hope that the behaviour of its discretizations is fairly close to that of discretizations of a generic homeomorphism.
\end{itemize}

We have chosen to define the homeomorphisms we compute with lacunary trigonometric series, to ``mimic'' the action of Baire theorem.

Note that the homeomorphisms $f_3$ and $f_4$ have at least one fixed point (for $f_3$, simply make simultaneously $p(x)$ and $q(y)$ vanish; for $f_4$, note that $0$ is a persistent fixed point of $A$). Theoretical results indicate that for a generic homeomorphism $f$ which has a fixed point, infinitely many discretizations has a unique fixed point; moreover a subsequence of $(\mu^{f_N}_{\T^2})_{N\in\N}$ tends to an invariant measure under $f$ supported by a fixed point; we can test if these properties are true on simulations or not.

We refer to the page~\pageref{Blabla128} for a presentation of the algorithm we used for the simulations.

\subsection{Combinatorial behaviour}\label{simulgrafcons}

As a first step, we are interested in some quantities related to the combinatorial behaviour of discretizations of homeomorphisms. These quantities are:
\begin{itemize}
\item the cardinality of the recurrent set $\Omega(f_N)$;
\item the number of periodic orbits of the map $f_N$;
\item the maximal size of a periodic orbit of $f_N$.
\end{itemize}
Recall that according to the theoretical results we obtained, for a generic homeomorphism, the degree of recurrence $\frac{\card(\Omega(f_N))}{q_N}$ should accumulate on the whole segment $[0,1]$ (Corollary~\ref{ConjEt}), the number of orbits of $f_N$ should be $1$ for some $N$ and should be bigger than (for example) $\sqrt{q_N}$ for other orders $N$ (Corollaries~\ref{typlax} and \ref{corovar3}) and the maximal size of a periodic orbit of $f_N$ should be sometimes $1$, and sometimes equal to $q_N$ (Corollaries~\ref{typlax} and \ref{corovar3}).

\begin{figure}[h!]
\begin{center}
\makebox[0.8\textwidth]{\parbox{0.8\textwidth}{%
\begin{minipage}[c]{.49\linewidth}
	\includegraphics[width=\linewidth,trim = .5cm .3cm .6cm .1cm,clip]{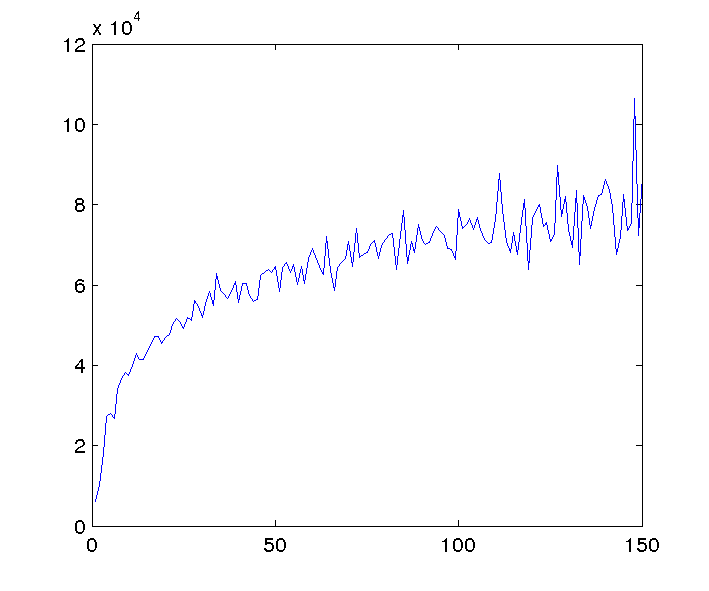}
\end{minipage}\hfill
\begin{minipage}[c]{.49\linewidth}
	\includegraphics[width=\linewidth,trim = .5cm .3cm .6cm .1cm,clip]{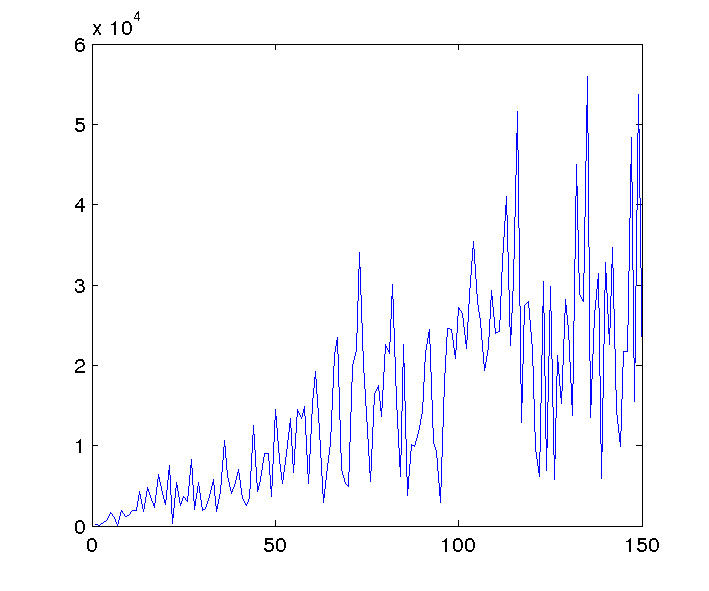}
\end{minipage}

\begin{minipage}[c]{.49\linewidth}
	\includegraphics[width=\linewidth,trim = .5cm .3cm .6cm .1cm,clip]{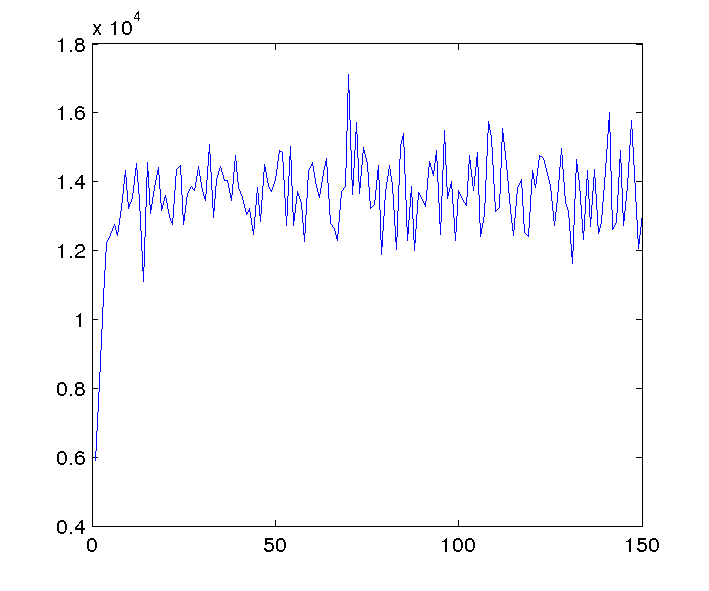}
\end{minipage}\hfill
\begin{minipage}[c]{.49\linewidth}
	\includegraphics[width=\linewidth,trim = .5cm .3cm .6cm .1cm,clip]{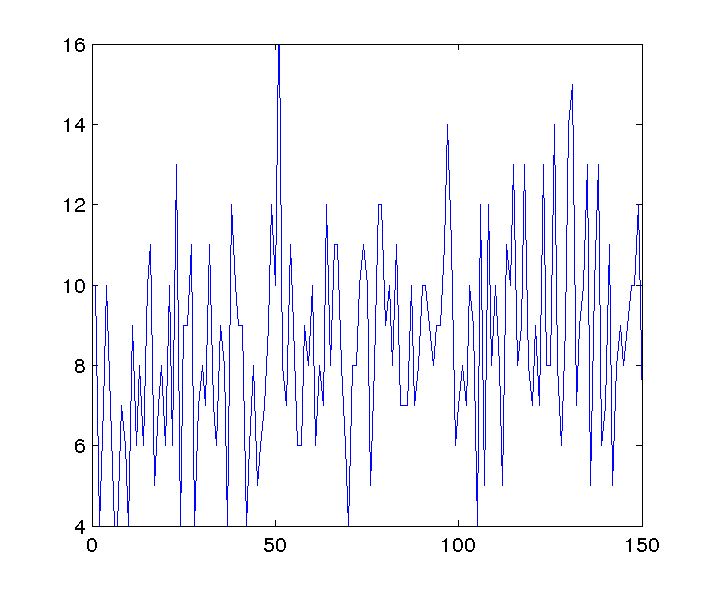}
\end{minipage}

\begin{minipage}[c]{.49\linewidth}
	\includegraphics[width=\linewidth,trim = .5cm .3cm .6cm .1cm,clip]{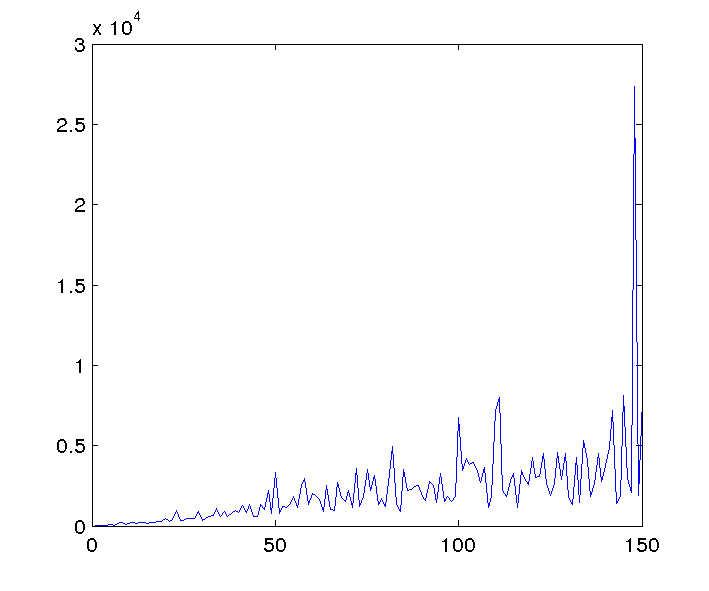}
\end{minipage}\hfill
\begin{minipage}[c]{.49\linewidth}
	\includegraphics[width=\linewidth,trim = .5cm .3cm .6cm .1cm,clip]{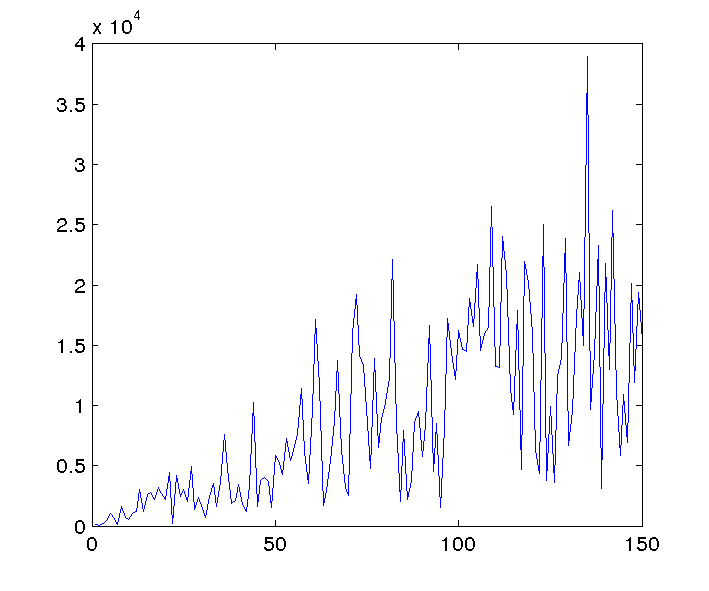}
\end{minipage}
}}\end{center}
\caption[Simulations of the combinatorial behaviour of 2 examples of conservative homeomorphisms]{Size of the recurrent set $\Omega((f_i)_N)$ (top), number of periodic orbits of $(f_i)_N$ (middle) and length of the largest periodic orbit of $(f_i)_N$ (bottom) depending on $N$, for $f_3$ (left) and $f_4$ (right), on the grids $E_N$ with $N=128k$, $k=1,\cdots,150$.}\label{GrafCons}
\end{figure}

We have calculated these quantities for discretizations of orders $128 k$, for $k$ from $1$ to $150$ and have represented them graphically (see Figure \ref{GrafCons}). For information, if $N =128\times 150$, then $q_N\simeq 3.7\times10^8$.

We begin with the cardinality of the recurrent set $\Omega(f_N)$ (Figure~\ref{GrafCons}, top). Contrary to what the theoretical results provide for a generic homeomorphism, for all simulations, the degree of recurrence $\frac{\card(\Omega(f_N))}{q_N}$ tends clearly to $0$ as $N$ increases. More specifically, the cardinality of $\Omega(f_N)$ evolves much more regularly for $f_3$ than for $f_4$: for $f_3$ the value of this cardinality seems to be the sum of a smooth increasing function and a random noise, but for $f_4$ this value seems to be the product of a smooth increasing function with a random noise. We have no explanation for the parabolic shape of the curve for $f_3$: it reflects the fact that the cardinality of $\Omega((f_3)_N) $ evolves in the same way as $\sqrt N = \sqrt[4]{q_N}$ (whereas for a random map of a finite set with $q$ elements into itself, it evolves as $\sqrt{q}$). Finally, it is interesting to note that for $f_4$, the size of the recurrent set is distributed more or less around the size of the recurrent set of a random map of a set with $q_N$ elements into itself, which depends (asymptotically) linearly of $N$ (because it is of the order of $\sqrt{q_N}$) and is worth about $24\,000$ for $N=128 \times 150$ (see \cite[XIV.5]{Boll-rand} or the Theorem~2.3.1 of \cite{Mier-dyna}).

According to the results of the previous sections, for a generic conservative homeomorphism $f$, the number of periodic orbits of $f_N$ should be sometimes equal to $1$ and sometimes bigger than (for example) $\sqrt{q_N}$. It is clearly not the case for these simulations. In addition, its behaviour is clearly not the same for $f_3$ and for $f_4$ (Figure \ref{GrafCons}, middle): for $f_3$ the number of orbits reaches quickly a value around $1.5\times 10^4$ to stagnate thereafter, while for $f_4$ it oscillates between $1$ and $16$, regardless the order of discretization. Note that contrary to what we have observed for the cardinality of the recurrent set, this seem to contradict the fact that the discretizations of $f_4$ behave like random maps: if this were the case, the behaviour of this quantity should be proportional to $\log N$, with a value close to $25$ for $N = 128\times 150$. These graphics can be compared with those representing the size of the recurrent set $\Omega((f_i)_N)$: the number of periodic orbits and the size of the recurrent set are of the same order of magnitude for $f_3$ (up to a factor $5$), which means that the average period of a periodic orbit is small (which is not surprising, since $f_3$ is a small perturbation of identity). They differ by a factor roughly equal to $10^3$ for $f_4$, which means that this time the average period of a periodic orbit is very large. This can be explained partly by the fact that the standard Anosov automorphism tends to mix what happens in the neighbourhood of identity. The fact remains that these simulations (such as the size of the recurrent set $\Omega ((f_i)_N)$) suggest that the behaviour in the neighbourhood of identity and of the standard Anosov automorphism are quite different, at least for reasonable orders of discretization.

Regarding the maximum size of a periodic orbit of $f_N$ (Figure~\ref{GrafCons}, bottom), again its behaviour does not correspond to that of a generic homeomorphism: it should oscillate between $1$ and $q_N$, and this is not the case. However, it varies widely depending on $N$, especially when $N$ is large. The qualitative behaviours are quite similar for the three simulations, but there are some quantitative differences: the maximum of the maximal length of a periodic orbit is twice greater for $f_4$ than for $f_3$. For these simulations, there is no significant difference between the behaviours of the discretizations of $f_3$: if we remove the very big value that is attained by the maximal length of the periodic orbit of $f_3$ (for $N$ close to $145$), the graphics are very similar.

\subsection{Behaviour of invariant measures}

\begin{figure}[h!]
\begin{center}
\makebox[0.8\textwidth]{\parbox{0.8\textwidth}{%
\begin{minipage}[c]{.49\linewidth}
	\includegraphics[width=\linewidth,trim = .5cm .3cm .6cm .1cm,clip]{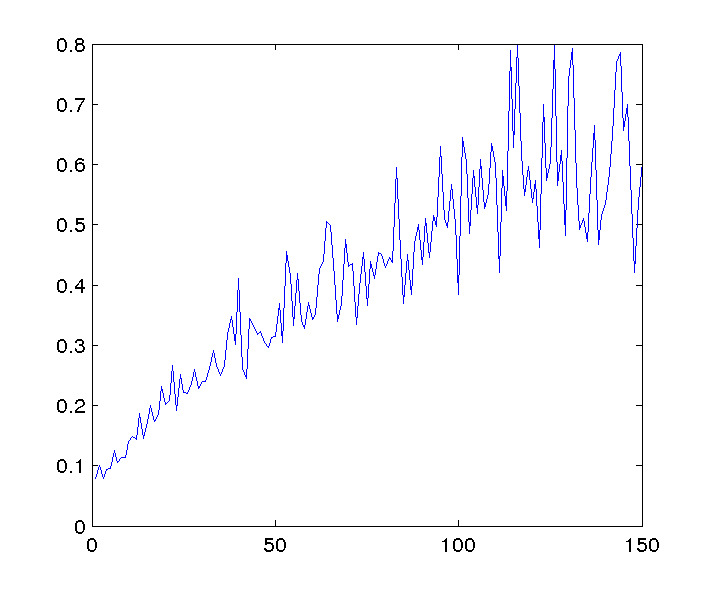}
\end{minipage}\hfill
\begin{minipage}[c]{.49\linewidth}
	\includegraphics[width=\linewidth,trim = .5cm .3cm .6cm .1cm,clip]{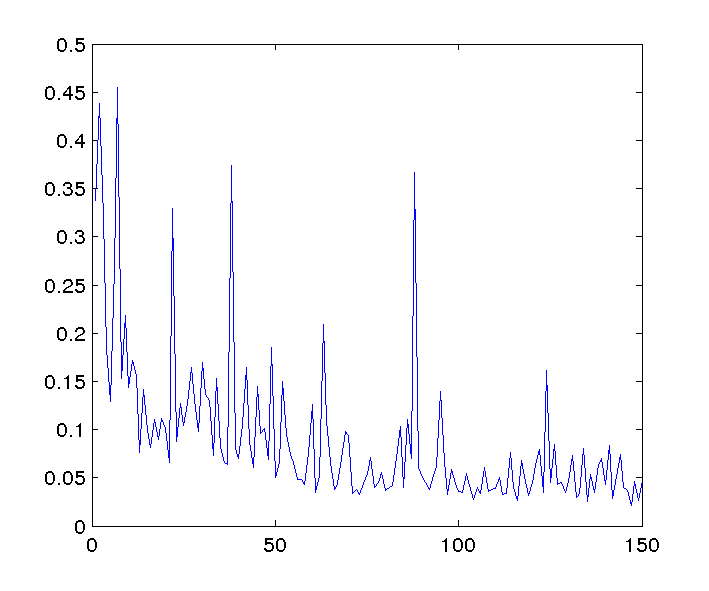}
\end{minipage}
}}\end{center}
\caption[Simulation of the distance between $\Leb$ and $\mu^{(f_i)_N}_{\T^2}$ for 3 examples of conservative homeomorphisms]{Distance between Lebesgue measure and the measure $\mu^{(f_i)_N}_{\T^2}$ depending on $N$ for $f_3$ (left) and $f_4$ (right), on the grids $E_N$ with $N=128k$, $k=1,\cdots,150$.}\label{GrafDistLebCons}
\end{figure}

We have simulated the measure $\mu^{(f_i)_N}_{\T^2}$ for the two examples of conservative homeomorphisms $f_3$ and $f_4$ as defined page~\pageref{defex}. Our purpose is to test whether phenomena as described by Theorem~\ref{mesinv} can be observed in practice or not. We obviously cannot expect to see the sequence $(\mu^{(f_i)_N}_{\T^2})_{N\in\N}$ accumulating on \emph{all} the invariant probability measures of $f$, since these measures generally form an infinite-dimensional convex set\footnote{To see that, simply observe that the set of periodic points is a union of Cantor sets; the uniform measures on these periodic orbits form an uncountable free family.}, but we can still test if it seems to converge or not. In particular, we can try to detect whether Lebesgue measure (which is the unique physical measure by Oxtoby-Ulam theorem \cite{Oxto-meas}) play a particular role for these invariant measures or not.

We present images of sizes $128\times 128$ pixels representing in logarithmic scale the density of the measures $\mu^{f_N}_{\T^2}$: each pixel is coloured according to the measure carried by the set of points of $E_N$ it covers. Blue corresponds to a pixel with very small measure and red to a pixel with very high measure. Scales on the right of each image corresponds to the measure of one pixel on the $\log 10$ scale: if green corresponds to $-3$, then a green pixel will have measure $10^{-3}$ for $\mu^{f_N}_{\T^2}$. For information, when Lebesgue measure is represented, all the pixels have a value about $-4.2$.
\bigskip

Firstly, we have computed the distance between the measure $\mu^{f_N}_{\T^2}$ and Lebesgue measure. The distance we have chosen is given by the formula
\[d(\mu,\nu) = \sum_{k=0}^\infty \frac{1}{2^k} \sum_{i,j=0}^{2^k-1} \big| \mu(C_{i,j,k}) - \nu(C_{i,j,k})\big|\in[0,2],\]
where
\[C_{i,j,k} = \left[\frac{i}{2^k},\frac{i+1}{2^k}\right] \times \left[\frac{j}{2^k},\frac{j+1}{2^k}\right].\]
In practice, we have computed an approximation of this quantity by summing only on the $k\in\llbracket 0,7 \rrbracket$.
Theoretically, the distance between the measure $\mu^{f_N}_{\T^2}$ and $\Leb$ (Figure~\ref{GrafDistLebCons}) should oscillate asymptotically between $0$ and the supremum over the set of $f$-invariant measures $\mu$ of the distance between Lebesgue measure and $\mu$. For $f_3$ (left of Figure~\ref{GrafDistLebCons}), wee see that the distance between $\mu^{f_N}_{\T^2}$ and $\Leb$ seem to increase with $N$: the measure $\mu^{f_N}_X$ is close to Lebesgue measure for small values of $N$, and then is more and more far away from Lebesgue measure. It can be explained by the fact that for a small $N$, the discretization $f_N$ does not see the irregularities of $f$ and is more or less close to the discretization of identity. The behaviour for $f_4$ is very different: overall, the distance between $\mu^{f_N}_{\T^2}$ and $\Leb$ is much smaller and decreasing, but we observe peaks of this distance: there are a few values of $N$ for which $\mu^{f_N}_{\T^2}$ is far away from $\Leb$.

\begin{figure}[ht]
\begin{minipage}[c]{.31\linewidth}
	\includegraphics[height=4.8cm,trim = 1.5cm .95cm 2.8cm .5cm,clip]{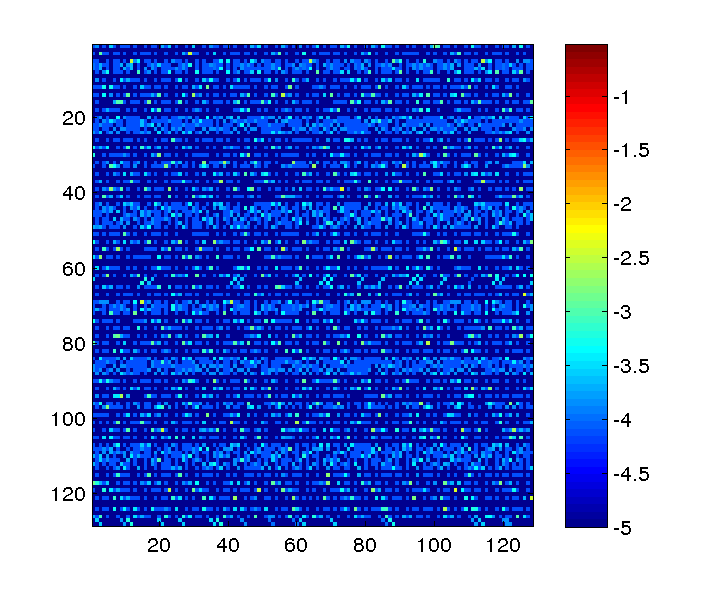}
\end{minipage}\hfill
\begin{minipage}[c]{.31\linewidth}
	\includegraphics[height=4.8cm,trim = 1.5cm .95cm 2.8cm .5cm,clip]{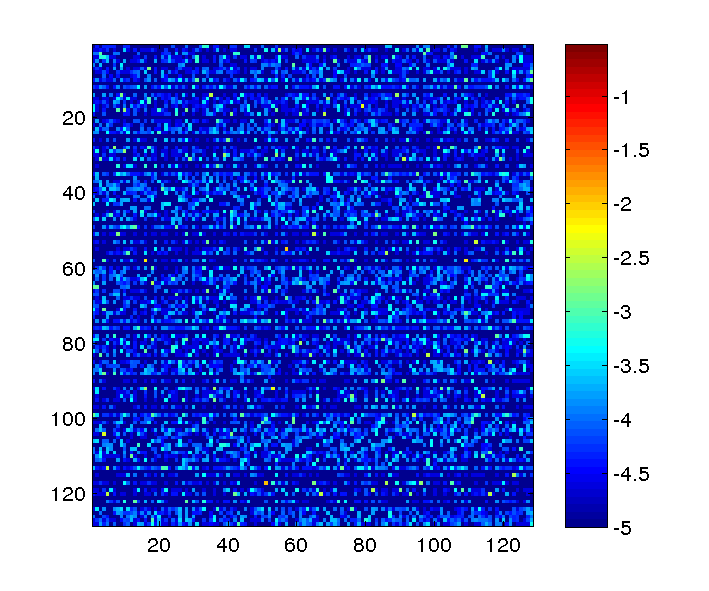}
\end{minipage}\hfill
\begin{minipage}[c]{.37\linewidth}
	\includegraphics[height=4.8cm,trim = 1.5cm .95cm 1cm .5cm,clip]{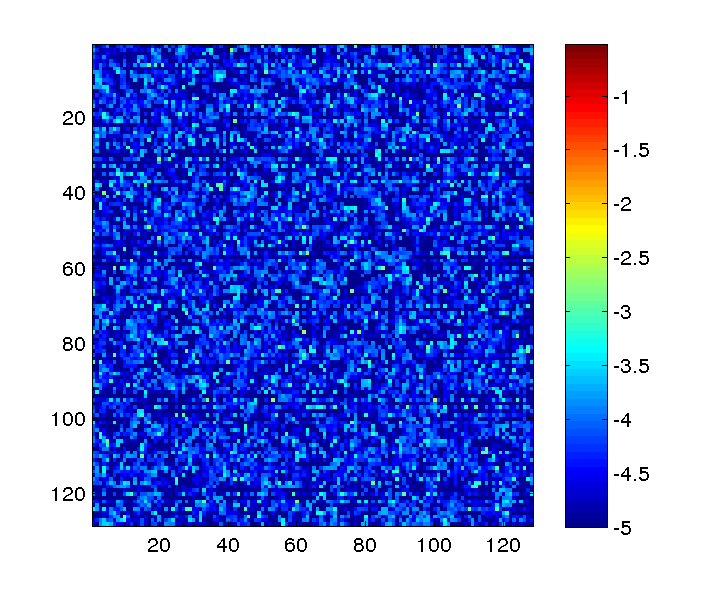}
\end{minipage}

\begin{minipage}[c]{.31\linewidth}
	\includegraphics[height=4.8cm,trim = 1.5cm .95cm 2.8cm .5cm,clip]{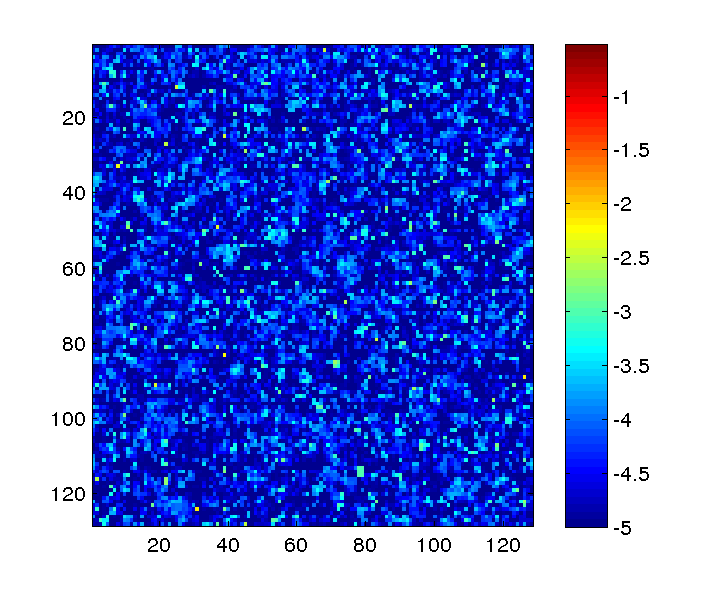}
\end{minipage}\hfill
\begin{minipage}[c]{.31\linewidth}
	\includegraphics[height=4.8cm,trim = 1.5cm .95cm 2.8cm .5cm,clip]{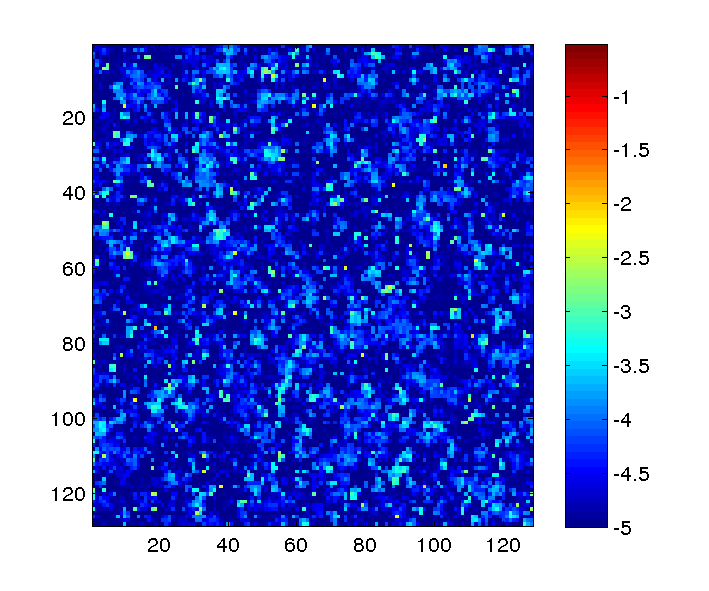}
\end{minipage}\hfill
\begin{minipage}[c]{.37\linewidth}
	\includegraphics[height=4.8cm,trim = 1.5cm .95cm 1cm .5cm,clip]{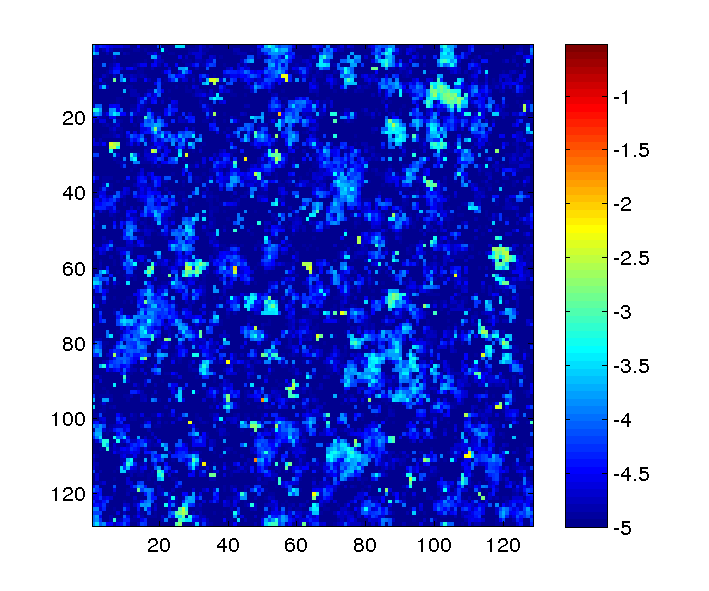}
\end{minipage}

\begin{minipage}[c]{.31\linewidth}
	\includegraphics[height=4.8cm,trim = 1.5cm .95cm 2.8cm .5cm,clip]{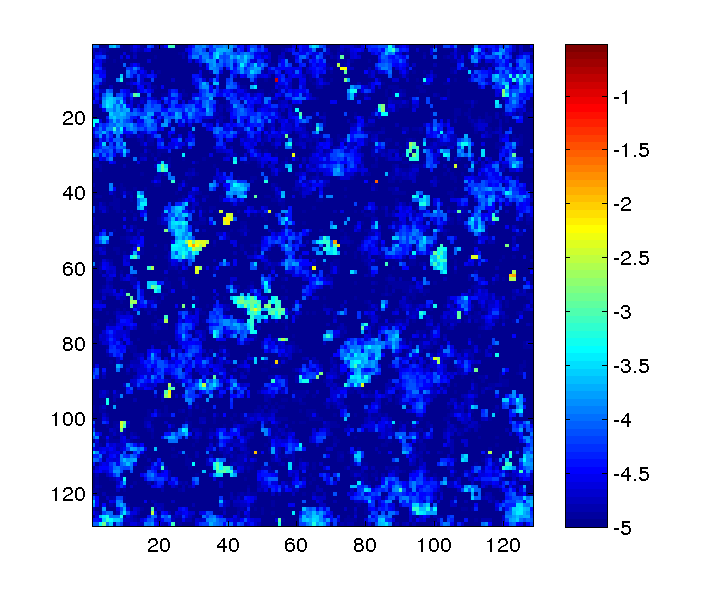}
\end{minipage}\hfill
\begin{minipage}[c]{.31\linewidth}
	\includegraphics[height=4.8cm,trim = 1.5cm .95cm 2.8cm .5cm,clip]{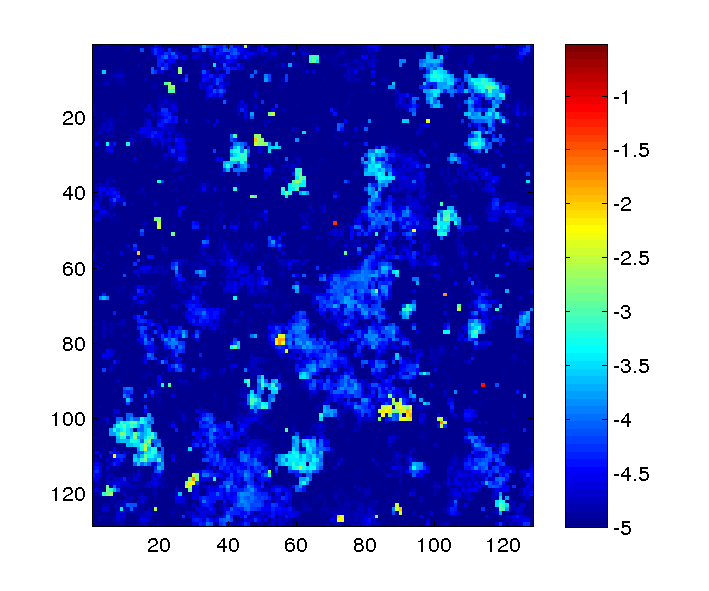}
\end{minipage}\hfill
\begin{minipage}[c]{.37\linewidth}
	\includegraphics[height=4.8cm,trim = 1.5cm .95cm 1cm .5cm,clip]{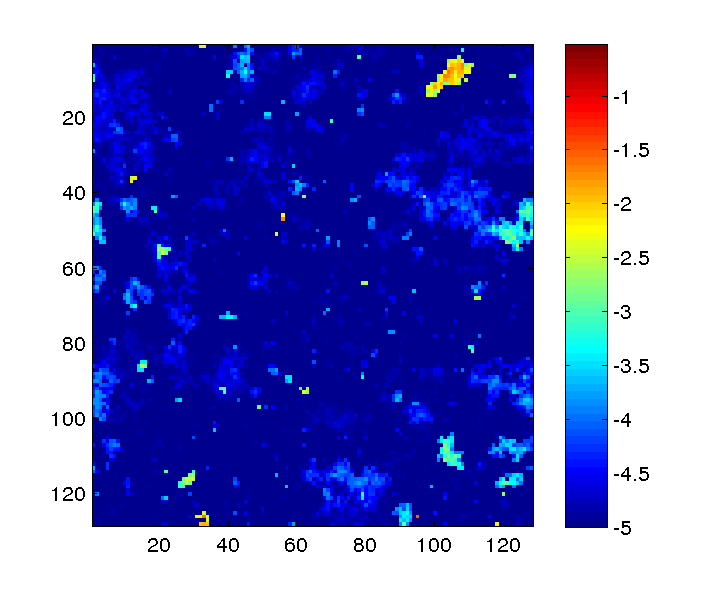}
\end{minipage}
\caption[Simulations of $\mu^{(f_3)_N}_{\T^2}$ on the grids $E_N$, with $N=2^k$, $k= 7,\cdots,15$]{Simulations of invariant measures $\mu^{(f_3)_N}_{\T^2}$ on the grids $E_N$, with $N=2^k$, $k= 7,\cdots,15$ (from left to right and top to bottom).}\label{MesC0IdCons2p}
\end{figure}

\begin{figure}[ht]
\begin{minipage}[c]{.31\linewidth}
	\includegraphics[height=4.8cm,trim = 1.5cm .95cm 2.8cm .5cm,clip]{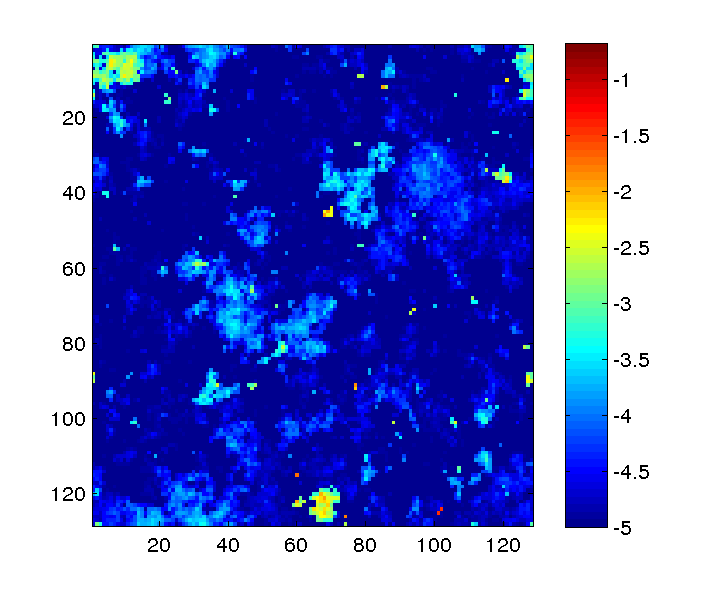}
\end{minipage}\hfill
\begin{minipage}[c]{.31\linewidth}
	\includegraphics[height=4.8cm,trim = 1.5cm .95cm 2.8cm .5cm,clip]{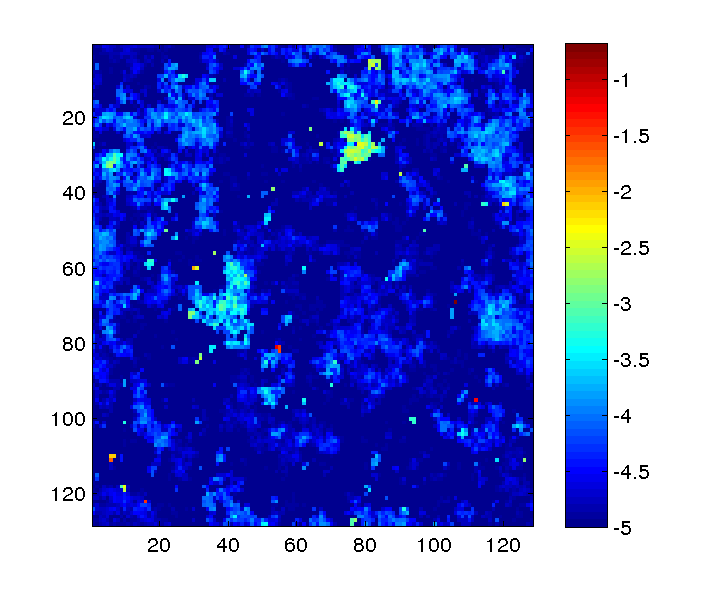}
\end{minipage}\hfill
\begin{minipage}[c]{.37\linewidth}
	\includegraphics[height=4.8cm,trim = 1.5cm .95cm 1cm .5cm,clip]{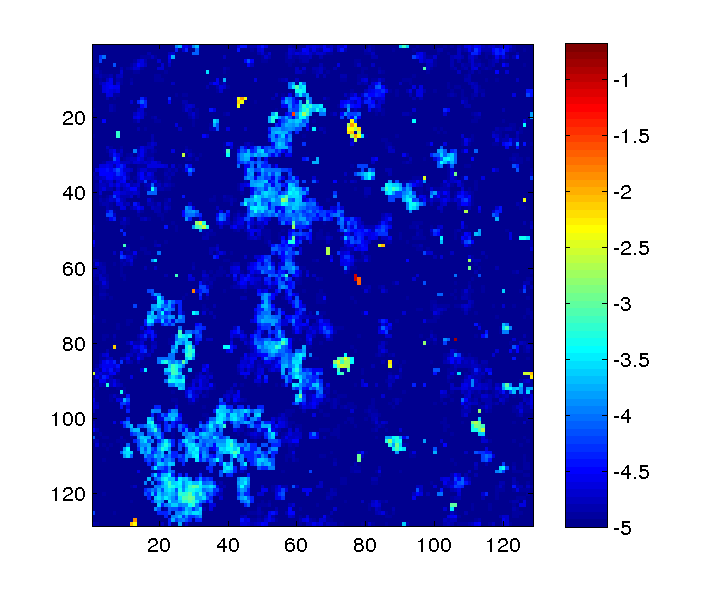}
\end{minipage}

\begin{minipage}[c]{.31\linewidth}
	\includegraphics[height=4.8cm,trim = 1.5cm .95cm 2.8cm .5cm,clip]{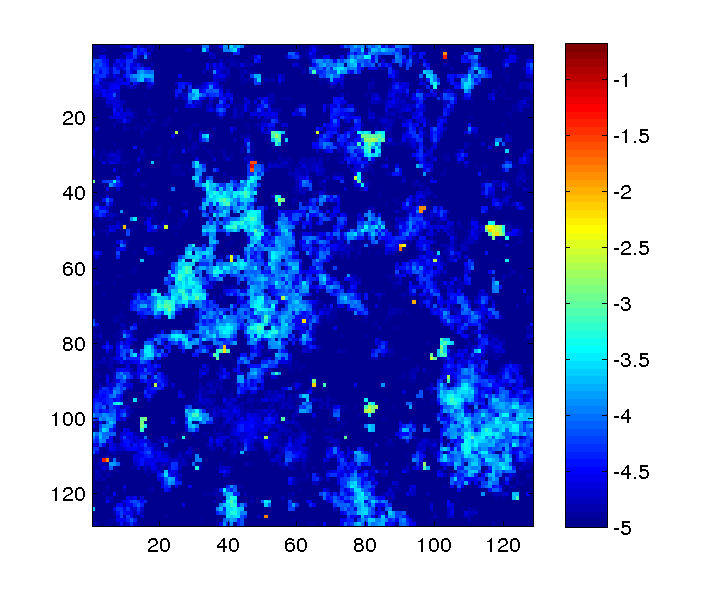}
\end{minipage}\hfill
\begin{minipage}[c]{.31\linewidth}
	\includegraphics[height=4.8cm,trim = 1.5cm .95cm 2.8cm .5cm,clip]{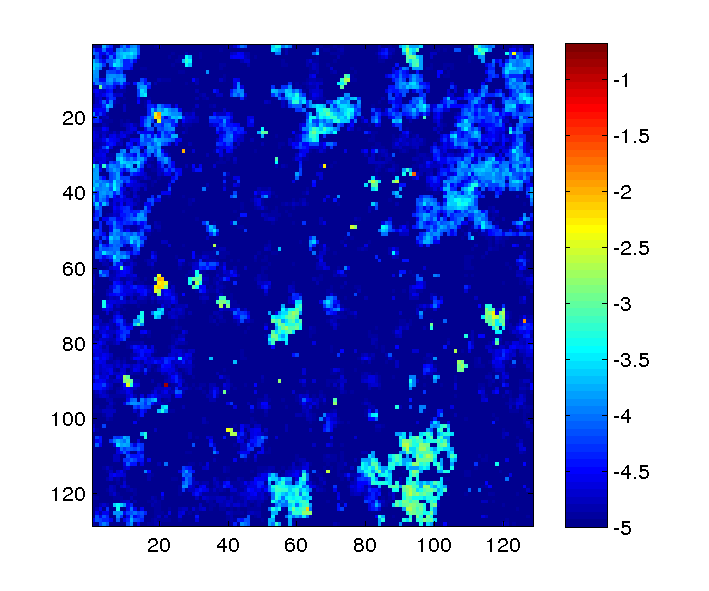}
\end{minipage}\hfill
\begin{minipage}[c]{.37\linewidth}
	\includegraphics[height=4.8cm,trim = 1.5cm .95cm 1cm .5cm,clip]{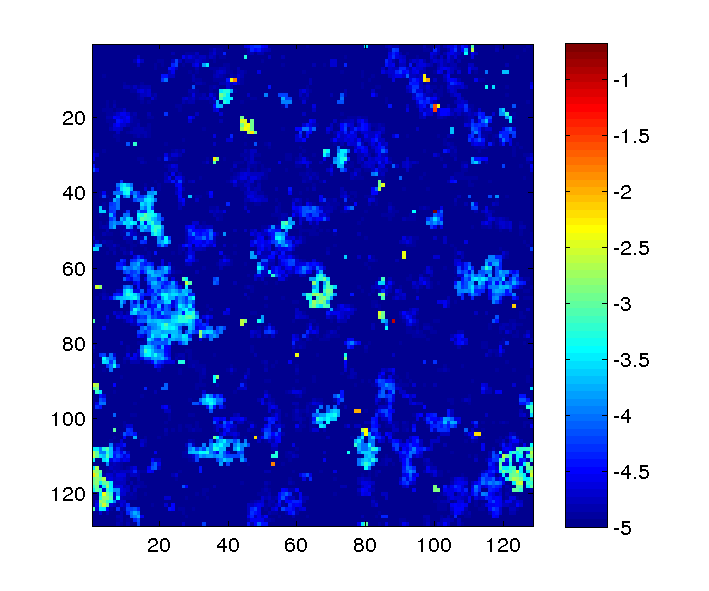}
\end{minipage}

\begin{minipage}[c]{.31\linewidth}
	\includegraphics[height=4.8cm,trim = 1.5cm .95cm 2.8cm .5cm,clip]{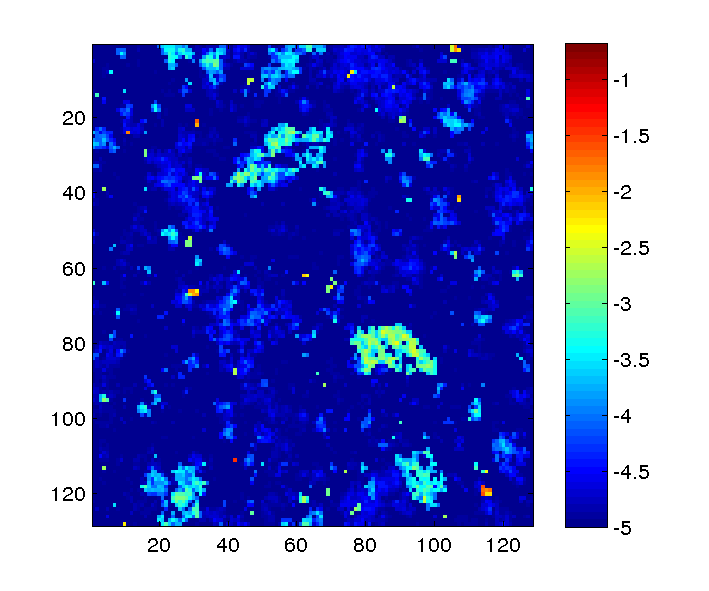}
\end{minipage}\hfill
\begin{minipage}[c]{.31\linewidth}
	\includegraphics[height=4.8cm,trim = 1.5cm .95cm 2.8cm .5cm,clip]{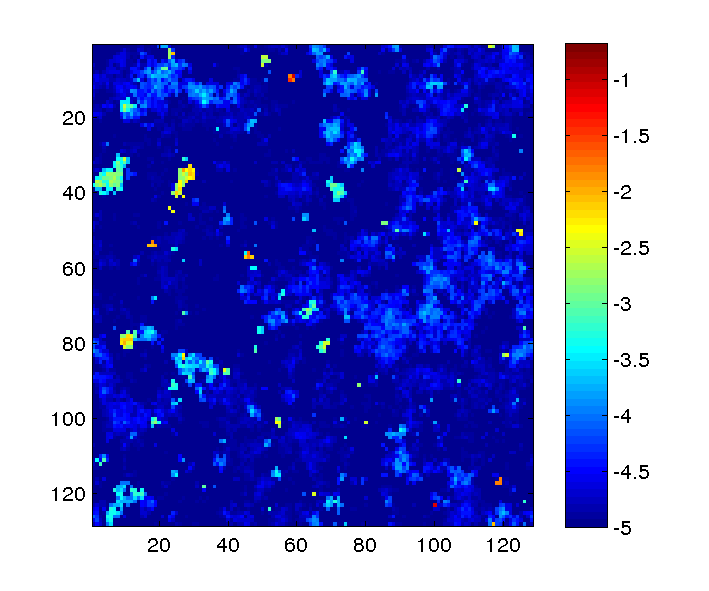}
\end{minipage}\hfill
\begin{minipage}[c]{.37\linewidth}
	\includegraphics[height=4.8cm,trim = 1.5cm .95cm 1cm .5cm,clip]{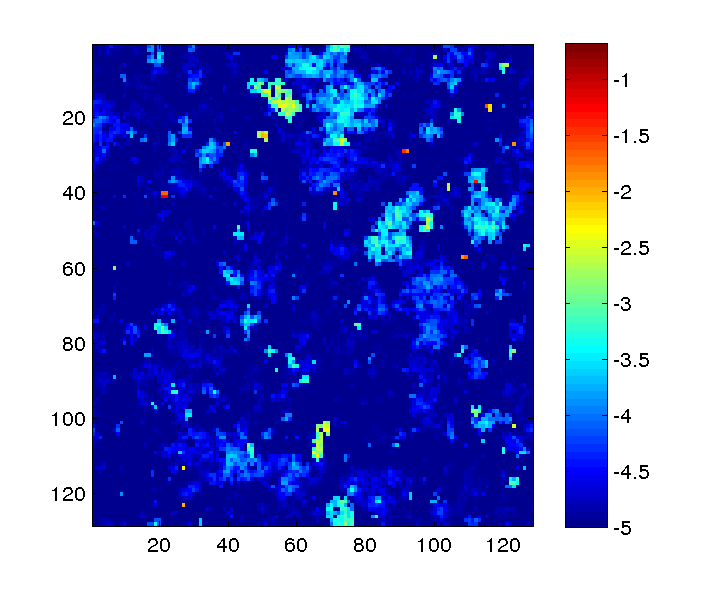}
\end{minipage}
\caption[Simulations of $\mu^{(f_3)_N}_{\T^2}$ on the grids $E_N$, with $N=20\,000,\cdots,20\,008$]{Simulations of invariant measures $\mu^{(f_3)_N}_{\T^2}$ on the grids $E_N$, with $N=20\,000,\cdots,20\,008$ (from left to right and top to bottom).}\label{MesC0IdConsSer}
\end{figure}

The results of simulations of invariant measures of discretizations of $f_3$ (which is a $C^0$ conservative perturbation of identity) are fairly positive: they agree with theoretical results about discretizations of generic conservative homeomorphisms, in particular with Theorem \ref{mesinv}. Indeed, when we compute discretizations of the homeomorphism $f_3$ on the grids of size $2^k\times 2^k$ (Figure~\ref{MesC0IdCons2p}), we first observe that $\mu^{f_N}_{\T^2}$ is fairly well distributed, say for $2^k=128,256,512$. When the order of the discretization increases, we can observe places where the measure accumulates; moreover these places changes a lot when the order of discretization varies (see Figure~\ref{MesC0IdConsSer}): for the orders $N$ we tested, the measures $\mu^{f_N}_{\T^2}$ and $\mu^{f_{N+1}}_{\T^2}$ are always completely different, and do not have anything in common with Lebesgue measure. This agrees fairly with what happens in the $C^0$ generic case, where we have proved that the measure $\mu^{f_N}_X$ depends very much on the order of discretization, rather than on the homeomorphism itself. There is also an other phenomenon: when the size of the grid is large enough (around $10^{12}\times 10^{12}$), some areas uniformly charged by the measure $\mu^{f_N}_{\T^2}$ appear; their sizes seem to be inversely proportional to the common mass of the pixels of the area; we do not know why this phenomenon appears\dots

\begin{figure}[ht]
\begin{minipage}[c]{.31\linewidth}
	\includegraphics[height=4.8cm,trim = 1.5cm .95cm 2.8cm .5cm,clip]{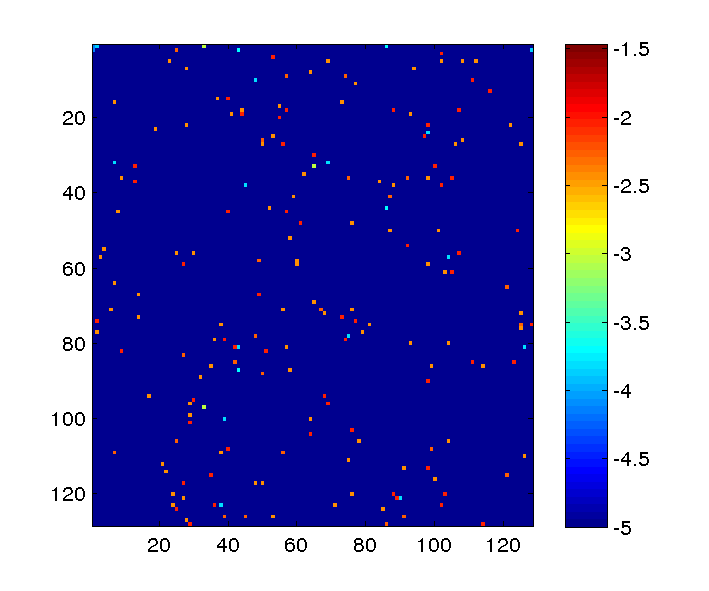}
\end{minipage}\hfill
\begin{minipage}[c]{.31\linewidth}
	\includegraphics[height=4.8cm,trim = 1.5cm .95cm 2.8cm .5cm,clip]{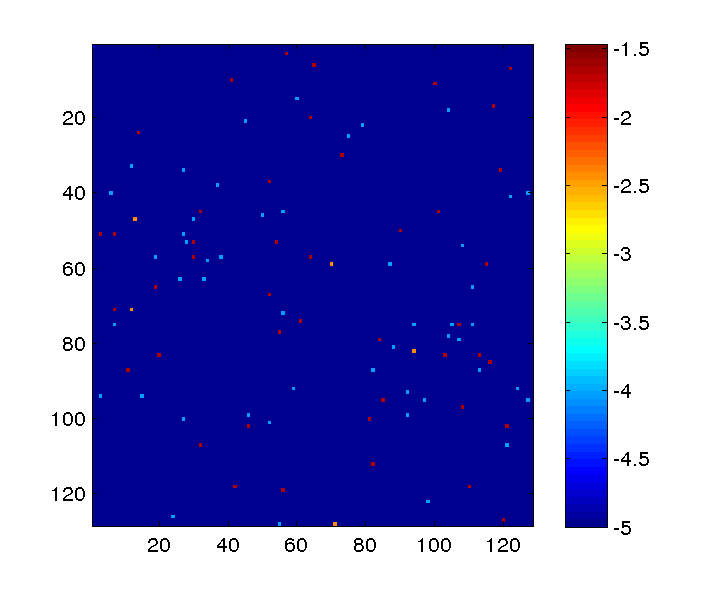}
\end{minipage}\hfill
\begin{minipage}[c]{.37\linewidth}
	\includegraphics[height=4.8cm,trim = 1.5cm .95cm 1cm .5cm,clip]{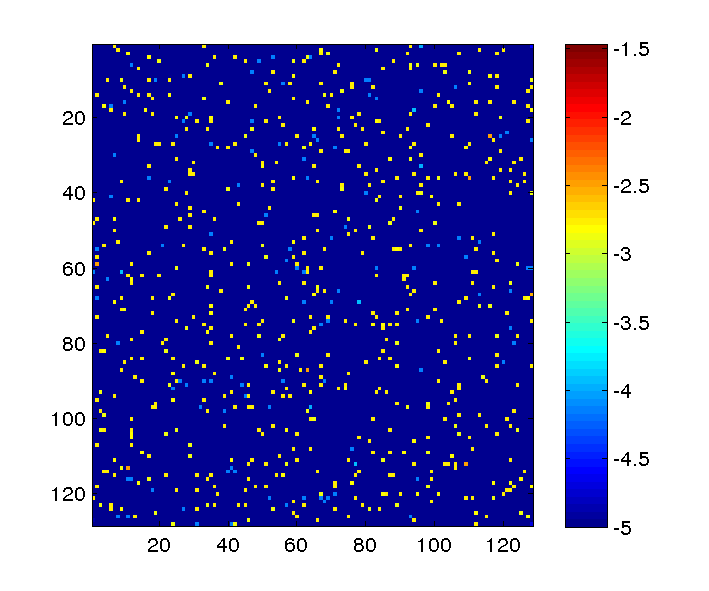}
\end{minipage}

\begin{minipage}[c]{.31\linewidth}
	\includegraphics[height=4.8cm,trim = 1.5cm .95cm 2.8cm .5cm,clip]{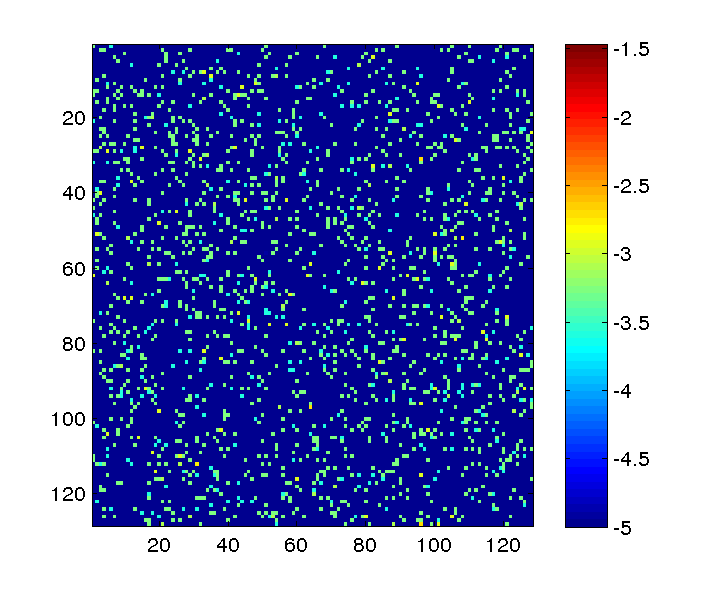}
\end{minipage}\hfill
\begin{minipage}[c]{.31\linewidth}
	\includegraphics[height=4.8cm,trim = 1.5cm .95cm 2.8cm .5cm,clip]{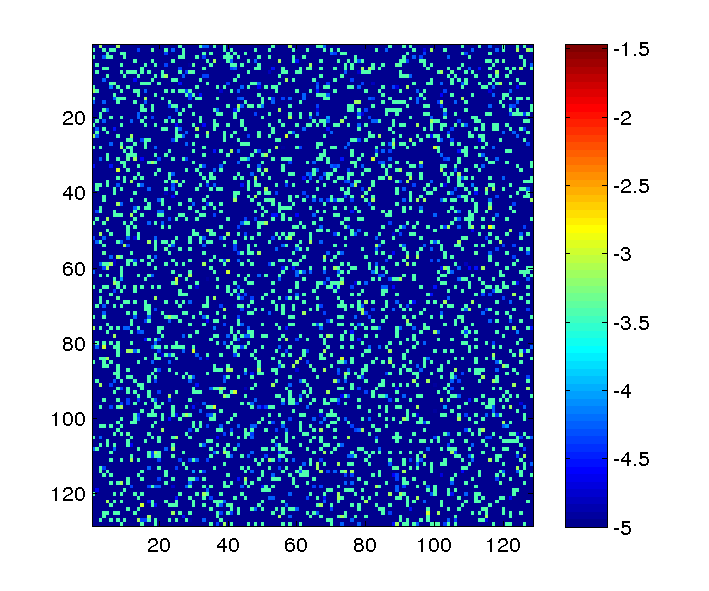}
\end{minipage}\hfill
\begin{minipage}[c]{.37\linewidth}
	\includegraphics[height=4.8cm,trim = 1.5cm .95cm 1cm .5cm,clip]{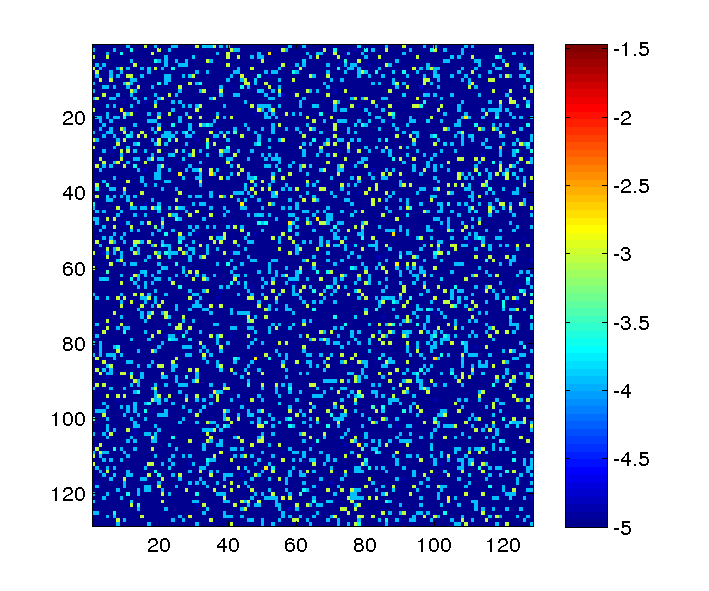}
\end{minipage}

\begin{minipage}[c]{.31\linewidth}
	\includegraphics[height=4.8cm,trim = 1.5cm .95cm 2.8cm .5cm,clip]{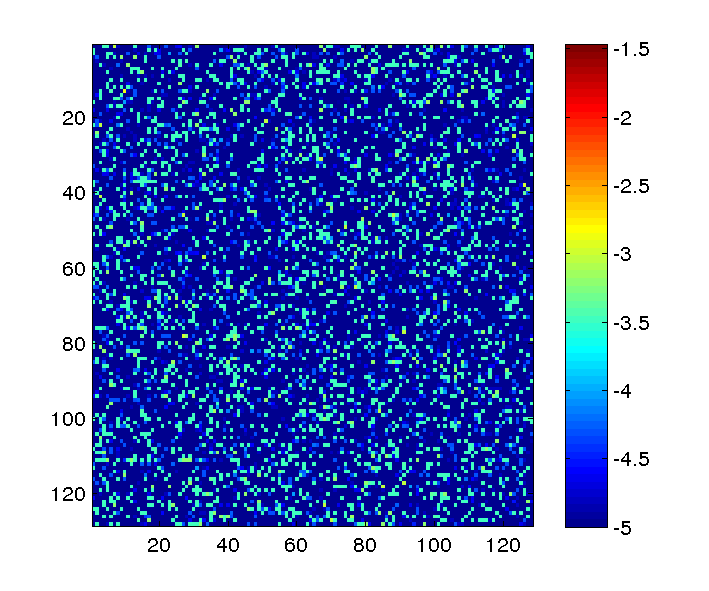}
\end{minipage}\hfill
\begin{minipage}[c]{.31\linewidth}
	\includegraphics[height=4.8cm,trim = 1.5cm .95cm 2.8cm .5cm,clip]{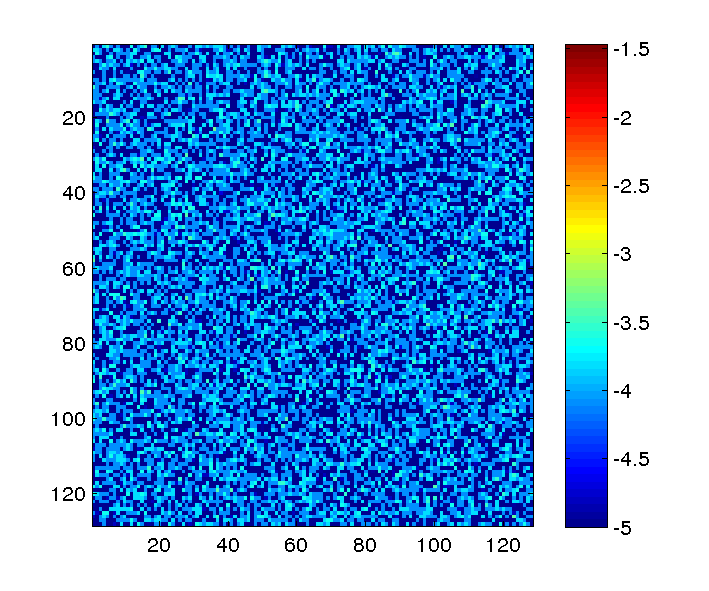}
\end{minipage}\hfill
\begin{minipage}[c]{.37\linewidth}
	\includegraphics[height=4.8cm,trim = 1.5cm .95cm 1cm .5cm,clip]{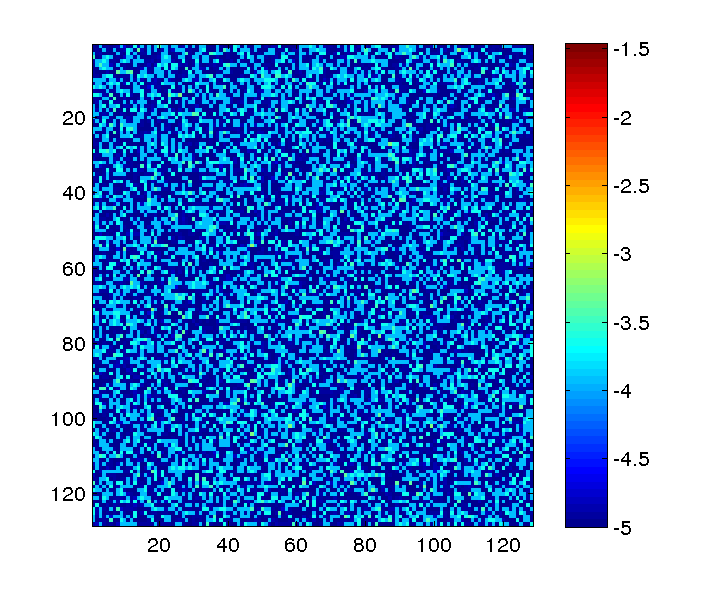}
\end{minipage}
\caption[Simulations of $\mu^{(f_4)_N}_{\T^2}$ on the grids $E_N$, with $N=2^k$, $k= 7,\cdots,15$]{Simulations of invariant measures $\mu^{(f_4)_N}_{\T^2}$ on the grids $E_N$, with $N=2^k$, $k= 7,\cdots,15$ (from left to right and top to bottom).}\label{MesC0AnoCons2p}
\end{figure}

\begin{figure}[ht]
\begin{minipage}[c]{.31\linewidth}
	\includegraphics[height=4.8cm,trim = 1.5cm .95cm 2.8cm .5cm,clip]{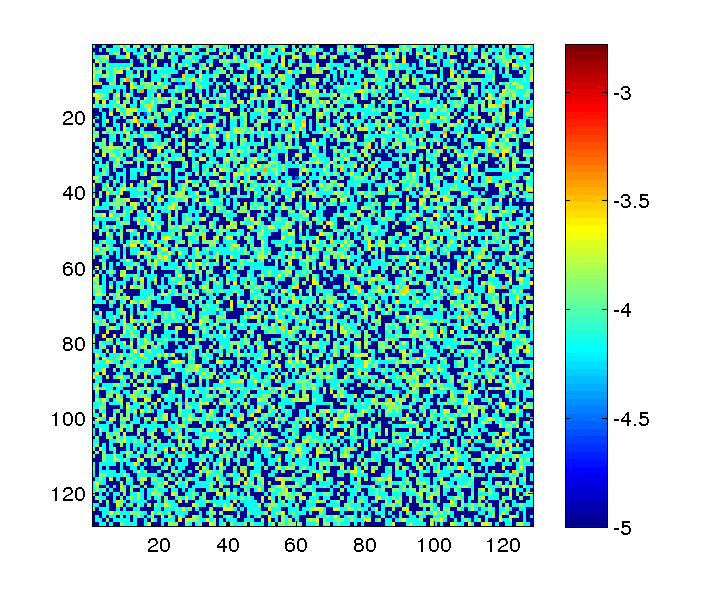}
\end{minipage}\hfill
\begin{minipage}[c]{.31\linewidth}
	\includegraphics[height=4.8cm,trim = 1.5cm .95cm 2.8cm .5cm,clip]{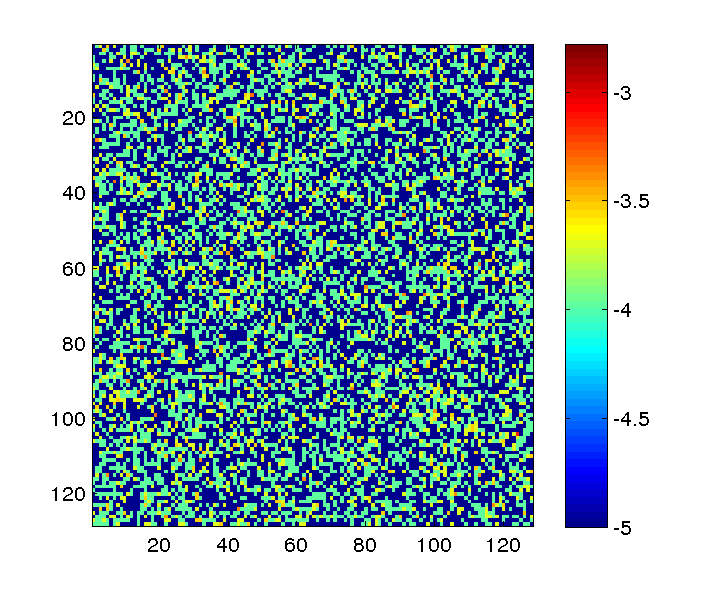}
\end{minipage}\hfill
\begin{minipage}[c]{.37\linewidth}
	\includegraphics[height=4.8cm,trim = 1.5cm .95cm 1cm .5cm,clip]{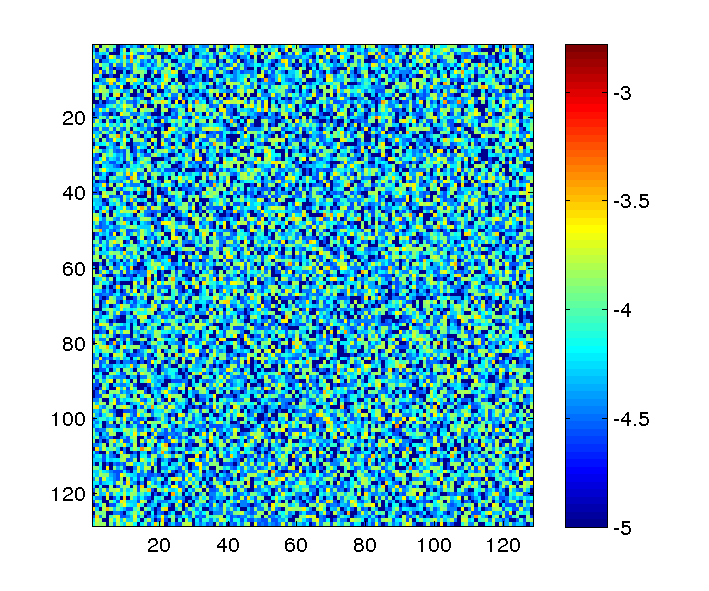}
\end{minipage}

\begin{minipage}[c]{.31\linewidth}
	\includegraphics[height=4.8cm,trim = 1.5cm .95cm 2.8cm .5cm,clip]{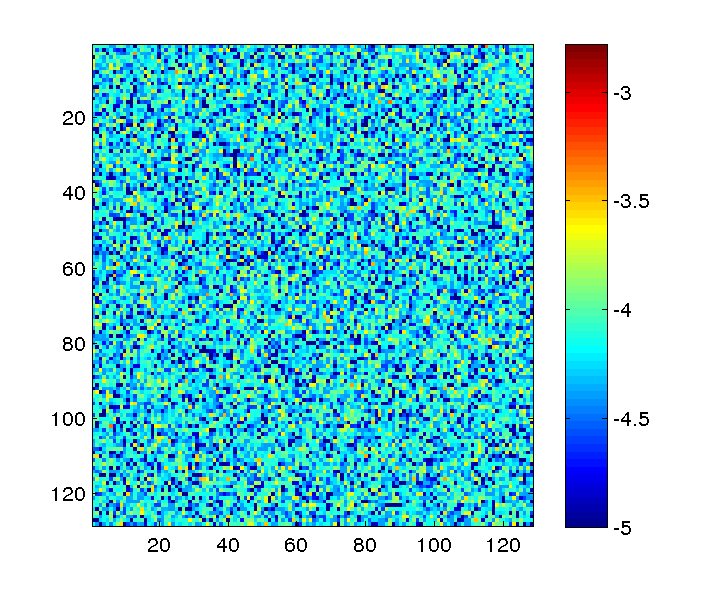}
\end{minipage}\hfill
\begin{minipage}[c]{.31\linewidth}
	\includegraphics[height=4.8cm,trim = 1.5cm .95cm 2.8cm .5cm,clip]{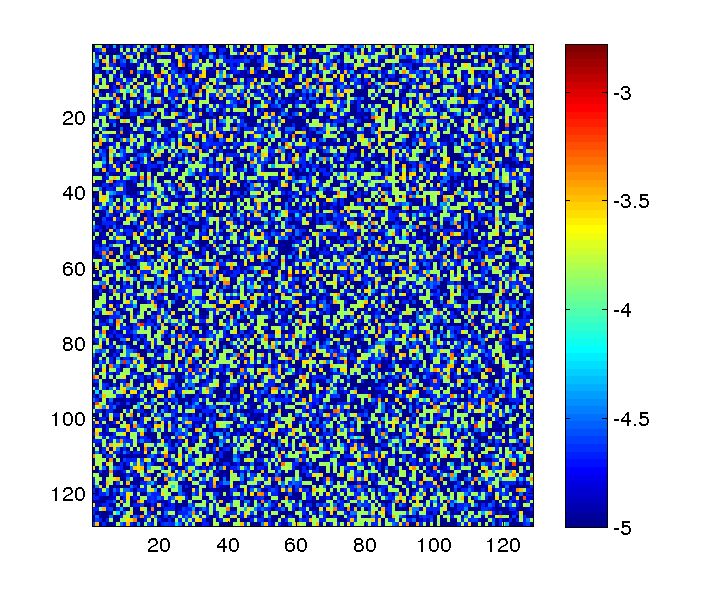}
\end{minipage}\hfill
\begin{minipage}[c]{.37\linewidth}
	\includegraphics[height=4.8cm,trim = 1.5cm .95cm 1cm .5cm,clip]{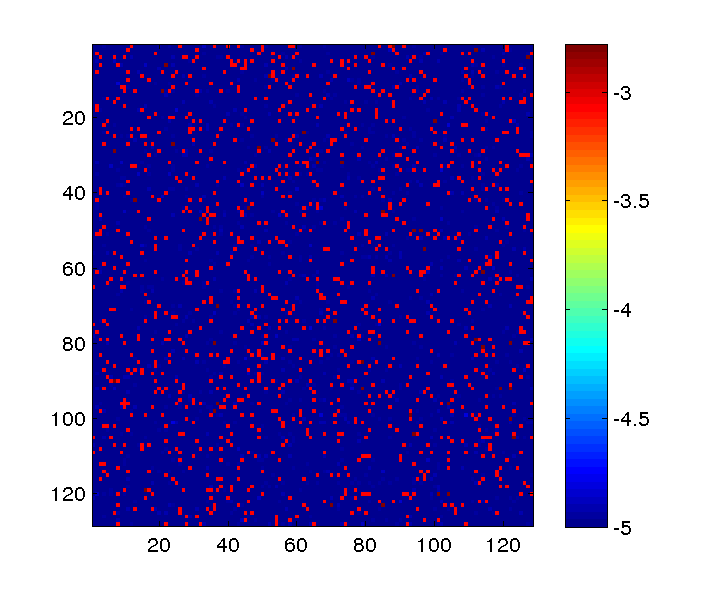}
\end{minipage}

\begin{minipage}[c]{.31\linewidth}
	\includegraphics[height=4.8cm,trim = 1.5cm .95cm 2.8cm .5cm,clip]{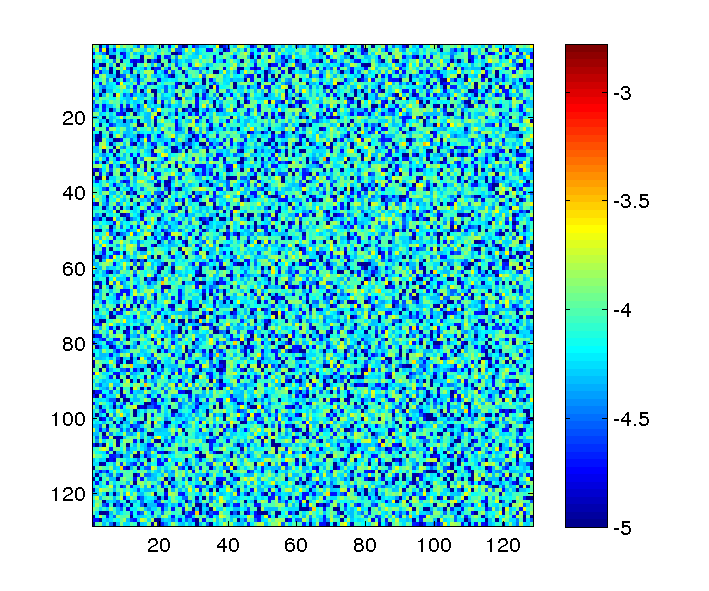}
\end{minipage}\hfill
\begin{minipage}[c]{.31\linewidth}
	\includegraphics[height=4.8cm,trim = 1.5cm .95cm 2.8cm .5cm,clip]{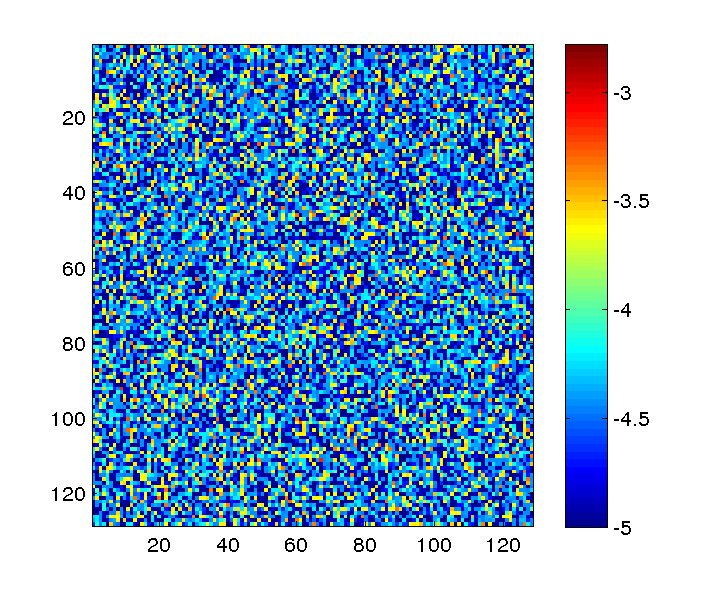}
\end{minipage}\hfill
\begin{minipage}[c]{.37\linewidth}
	\includegraphics[height=4.8cm,trim = 1.5cm .95cm 1cm .5cm,clip]{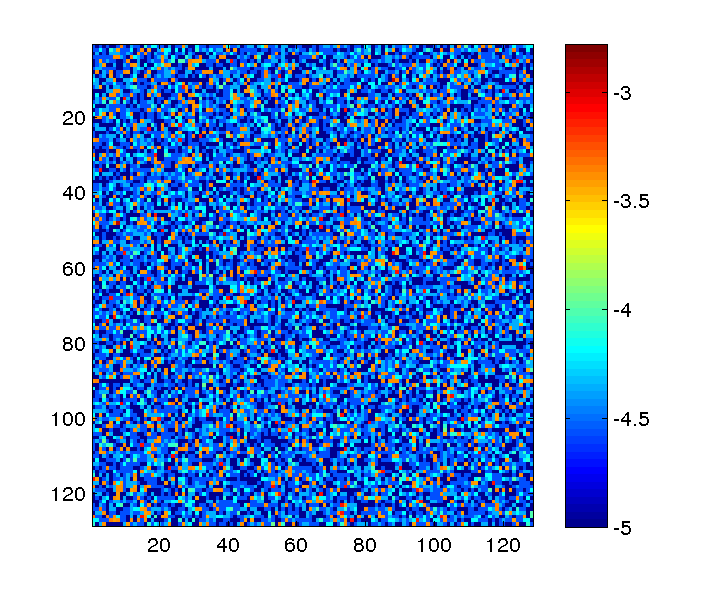}
\end{minipage}
\caption[Simulations of $\mu^{(f_4)_N}_{\T^2}$ on the grids $E_N$, with $N=22\,395,\cdots,22\,403$]{Simulations of invariant measures $\mu^{(f_4)_N}_{\T^2}$ on the grids $E_N$, with $N=22\,395,\cdots,22\,403$ (from left to right and top to bottom).}\label{MesC0AnoConsSer}
\end{figure}

For the discretizations of $f_4$, which is a $C^0$ conservative perturbation of the linear Anosov automorphism $A$, the simulations on grids of size $2^k\times 2^k$ might suggest that the measures $\mu^{(f_4)_N}_{\T^2}$ tend to Lebesgue measure (Figure~\ref{MesC0AnoCons2p}). In fact, making a large number of simulations, we realize that there are also strong variations of the behaviour of the measures (Figure \ref{MesC0AnoConsSer}): the measure is often well distributed in the torus, and sometimes quite singular with respect to Lebesgue measure (as it can be seen in Figure \ref{GrafDistLebCons}). For example, when we discretize on the grid of size $22\,400\times 22\,400$ (middle right of Figure \ref{MesC0AnoConsSer}), we observe an orbit of length $369$ which mass $84\%$ of the total measure.
In fact the behaviour of discretizations looks the same that in the neighbourhood of identity, modulo the fact that the linear Anosov automorphism $A$ tends to spread the attractive periodic orbits of the discretizations on the entire torus: for many values of $N$, composing by $A$ spreads the behaviour of the measure $\mu^{(f_4)_N}_{\T^2}$, but sometimes (in fact, seldom) a fixed point of $(f_3)_N$ which attracts a large part of $E_N$ is located around one of the few periodic points of small period for $A$. This then creates a periodic orbit for $(f_4)_N$ with a big measure for $\mu^{(f_4)_N}_{\T^2}$.

\clearpage

\subsection{Periodic points}

\begin{figure}[ht]
\begin{minipage}[c]{.33\linewidth}
	\includegraphics[width=\linewidth,trim = .5cm .3cm .5cm .3cm,clip]{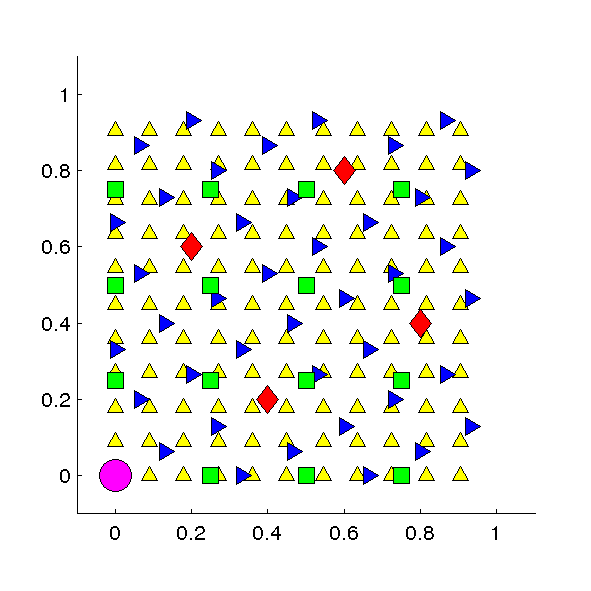}
\end{minipage}\hfill
\begin{minipage}[c]{.33\linewidth}
	\includegraphics[width=\linewidth,trim = .5cm .3cm .5cm .3cm,clip]{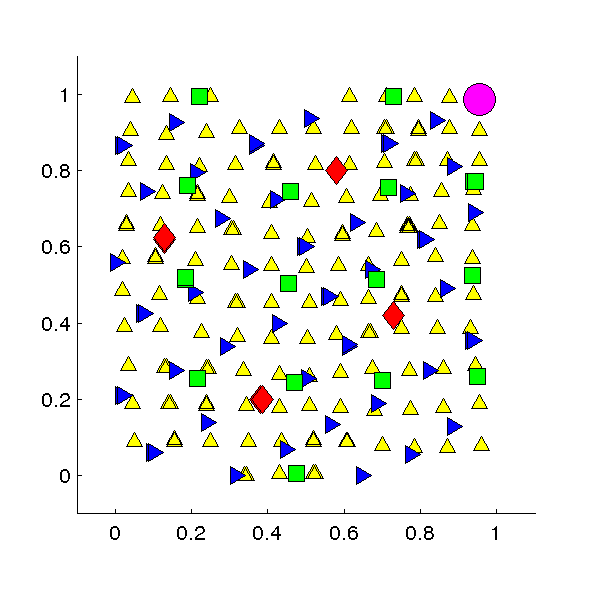}
\end{minipage}\hfill
\begin{minipage}[c]{.33\linewidth}
	\includegraphics[width=\linewidth,trim = .5cm .3cm .5cm .3cm,clip]{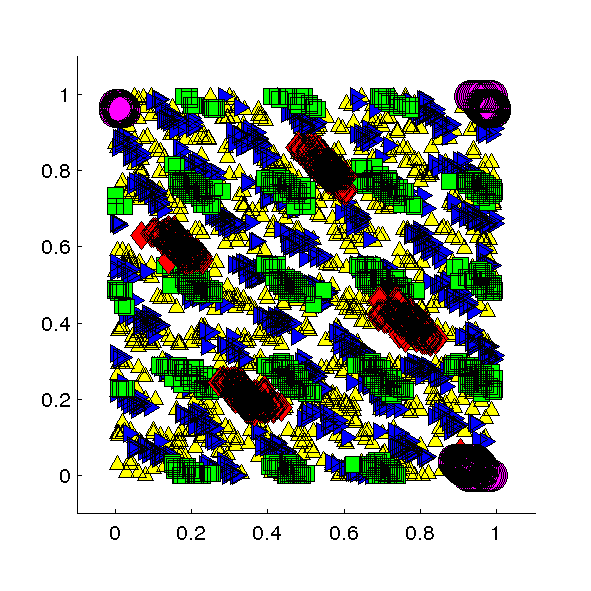}
\end{minipage}
\caption[Periodic points of discretizations of perturbations of $A$]{Unions of the sets of periodic points of period smaller than 5 of all the discretizations on grids $E_N$ with $1\,000\le N\le 2\,000$ for $A$ (left), a small $C^1$ perturbation of $A$ (middle) and a small $C^0$ perturbation of $A$ (right); period 1: purple circle, period 2: red diamond, period 3: green square, period 4: blue triangle pointing right, period 5: yellow triangle pointing up.}\label{FigPtsPer}
\end{figure}

We have also conducted simulations of the periodic orbits of period smaller than 5 of the discretizations for three conservative homeomorphisms:
\begin{itemize}
\item the linear map
\[A = \begin{pmatrix}
      2 & 1\\
      1 & 1
      \end{pmatrix};\]
\item a small $C^1$ perturbation of $A$, say $g_1 = Q_1\circ P_1\circ A$, with
\[P_1(x,y) = \left(x-0.01975\ , \ y+\frac{1}{41}\cos(2\pi\times x)+\frac{1}{351}\sin(2\pi\times 7 x)+0.02478\right),\]
\[Q_1(x,y) = \left(x+\frac{1}{47}\cos(2\pi\times y)+\frac{1}{311}\sin(2\pi\times 5 y)+0.01237\ , \ y+0.00975\right);\]
\item and a small $C^0$ perturbation of $A$, say $g_2 = Q_2\circ P_2\circ A$, with
\[P_2(x,y) = \left(x-0.01975\ , \ y+\frac{1}{31}\cos(2\pi\times 17x)+\frac{1}{41}\sin(2\pi\times 233 x)+0.02478\right),\]
\[Q_2(x,y) = \left(x+\frac{1}{27}\cos(2\pi\times 23y)+\frac{1}{37}\sin(2\pi\times 217 y)+0.01237\ , \ y+0.00975\right).\]
\end{itemize}

We represent the accumulation of the periodic points of period smaller than 5 of the discretizations $f_N$, for $1\,000 \le N \le 2\,000$.

It can be easily seen that the periodic points of $A$ are rational points; it is also easy to compute by hand the coordinates of the periodic points of period smaller than 5. As $A$ is Anosov and as Anosov maps are structurally stable (see \cite{MR1326374}, Theorem 2.6.3 for the linear maps of the torus and Corollary 18.2.2 for the general case), a small $C^1$ perturbation of $A$ is ($C^0$) conjugated to $A$: if $d_{C^1}(A,g)$ is small enough, then there exists a homeomorphism $h\in\Hom(\T^2)$ which maps bijectively each set of periodic points of a given period of $A$ to the set of periodic points with the same period of $g$. If the perturbation of $A$ is only $C^0$, then if $d_{C^0}(A,g)$ is small enough but $d_{C^1}(A,g)$ is big, then it can happen that the homeomorphism $g$ is only semi-conjugated to $A$. But for a generic homeomorphism, the set of periodic points of a given period is either empty, either a Cantor set with zero Hausdorff dimension. Thus, to any periodic point of the linear map $A$ will correspond a Cantor set of periodic points of $g$ (in other words, a generic small $C^0$ perturbation blows up each periodic point).

The theoretical results about discretizations assert that any periodic point of a generic homeomorphism is shadowed by a periodic orbit with the same period of an infinite number of discretizations (see Theorem~\ref{corovar2}); the same result holds for a generic $C^1$-diffeomorphism (see Lemma~\ref{PropShadowC1}). Thus, on the simulations, we should recover all the periodic points of the initial map for both a small $C^1$ and $C^0$ generic perturbation.

On the left of Figure~\ref{FigPtsPer}, we see that we recover all the periodic points of period $\le 5$ of the linear Anosov map $A$. It is logical: these periodic points have rational coordinates with small denominators, thus they are located on a lot of grids of order $N\in \llbracket 1\,000,2\,000\rrbracket$. 

In the case of a small $C^1$ perturbation $g_1$ of $A$ (middle of Figure~\ref{FigPtsPer}), we also recover all the periodic points of period $\le 5$ of the map: we can see that we detect one (and almost always only one) periodic point of $g_1$ in the neighbourhood of each periodic point af $A$. In a certain sense, the simulation detects the conjugation between $A$ and $g_1$.

For the small $C^0$ perturbation $g_2$ of $A$ (right of Figure~\ref{FigPtsPer}), we observe the phenomenon of blow up of the periodic points of $A$: to each periodic point of $A$ corresponds a lot of periodic points of the discretizations of $g_2$. Thus, this simulation show two phenomenon. Firstly, the fact that the set of periodic points of a generic homeomorphism is a Cantor set is illustrated by some simple examples of homeomorphisms. Secondly, the behaviour predicted by Theorem~\ref{corovar2} can be observed on practice: we recover (all) the periodic points of the homeomorphism with small period by looking at all the periodic points of the discretizations on grids of reasonable orders.

\newpage


\chapter[How roundoff errors help to compute the rotation set]{How roundoff errors help to compute the rotation set of torus homeomorphisms}\label{ChapRot}

%
%
%
%
%

\section*{Introduction}

We now present a practical application of the concept of discretization to the numerical computation of rotation set of torus homeomorphisms.

In this chapter, to be able to talk about the rotation set of the homeomorphisms, we will only consider homeomorphisms of the torus $\T^2$ which are homotopic to the identity. In other words, the space phase will always be $X = \T^2$. We fix once for all a good measure $\lambda$ on $\T^2$ (see Definition~\ref{bonne mesure}), and a sequence of grids $E_N$ on $\T^2$, which is sometimes strongly self-similar (see Definition~\ref{Ashe'}). As an example of such grids, when $\lambda=\Leb$, the reader can think about the canonical grids
\[E_N^0 = \left\{\left(\frac{i}{N},\frac{j}{N}\right)\in \T^2 \big|\ 1\le i,j\le N\right\}\]
(see Section~\ref{exgrilles} for other examples of grids).

We begin with the presentation of the concept of rotation set of a torus homeomorphism. 

The concept of rotation number for circle homeomorphisms was introduced by H. Poincaré in 1885. In \cite{Poincare}, he states the theorem of classification of orientation-preserving circle homeomorphisms: if a homeomorphism has a rational rotation number, then it posses a periodic point and all its periodic points have the same period; moreover the $\omega$-limit set of every point is a periodic orbit (the dynamics is asymptotically periodic). On the contrary, if a homeomorphism has an irrational rotation number $\alpha$, then it is semi-conjugated to the rigid rotation of angle~$\alpha$ (the dynamics contains that of the irrational rotation). Ever since, the rotation number has been the fundamental tool in the study of the dynamics of circle homeomorphisms (see for example \cite{MR538680}).

About a century after was introduced a generalisation to dimension 2 of this rotation number, the rotation set for homeomorphisms of the torus which are homotopic to the identity. Due to the loss of natural cyclic order on the phase space, there is no longer a single speed of rotation for orbits; informally the rotation set is defined as the set of all possible rotation speeds of all possible orbits. Like in dimension 1, this topological invariant gives precious informations about the dynamics of the homeomorphism; for example, depending of the shape of this set, we can ensure the existence of periodic points of a given period (\cite{MR967632}, \cite{MR958891}). Moreover, the size of this convex set gives lower bounds on the topological entropy of the homeomorphism (\cite{MR1101087} and \cite{MR1213082} for an explicit estimation)\dots\ See the course \cite{Begu-ens} of F.~Béguin for a quite complete survey of the results about the rotation sets.

The aim of this chapter is to tackle the question of numerical approximation of the rotation set: given a homeomorphism of the torus homotopic to the identity, is it possible to compute numerically its rotation set? In particular, is it possible to detect its dimension? Is it possible to approximate it in Hausdorff topology? And what algorithm shall we use to compute it?
\bigskip

First of all, we build a theoretical model of what happens when we try to calculate the rotation set of a homeomorphism with a computer. To do that, we first take into account the fact that the computer can calculate only a finite number of orbits; in particular it will detect only phenomenon that occur on $\lambda$-positive measure sets. This leads us to the notion of \emph{observable rotation set}: a rotation vector is called observable if it is the rotation vector of an observable measure in the sense given by E. Catsigeras and H. Enrich in \cite{MR2852870}; more precisely, a measure $\mu$ is observable if for every $\varep>0$, the set of points which have a Birkhoff limit whose distance to $\mu$ is smaller than $\varep$ has $\lambda$-positive measure (see Definition \ref{DefObs}).

However, this notion of observable measure does not take into account the fact that the computer uses finite precision numbers and can calculate only finite length orbits; this observation leads to the definition leads to the definition of the \emph{asymptotic discretized rotation set} in the following way. We fix a sequence of finite grids on the torus with precision going to 0; the discretized rotation set on one of these grids is the collection of rotation vectors of periodic orbits of the discretization of the homeomorphism on this grid (see section~\ref{sectionNotRot}); the asymptotic discretized rotation set is then the upper limit of these discretized rotation sets on the grids.

We focus mainly on the generic behaviour of both observable and asymptotic discretized rotation sets. We recall that a result of A. Passeggi states that for a generic dissipative homeomorphism of the torus the rotation set is a polygon with rational vertices, possibly degenerated\footnote{Namely it can be a segment or a singleton. However there are open sets of homeomorphisms where the rotation set has non-empty interior.} \cite{rata}. In this chapter we will prove the following result about generic dissipative homeomorphisms.

\begin{theorem}
For a generic dissipative homeomorphism,
\begin{enumerate}
\item the observable rotation set is the closure of the set of rotation vectors corresponding to Lyapunov stable periodic points (Lemma \ref{ExGeneDissip});
\item the convex hull of the observable rotation set, the convex hull of the asymptotic discretized rotation set and the rotation set are equal;
\item if the rotation set has non-empty interior, there is no need to take convex hulls, \emph{i.e.} both observable and asymptotic discretized rotation sets coincide with the rotation set (Propositions \ref{GeneDissip} and \ref{RotDiscrDissip}).
\end{enumerate}
\end{theorem}

Thus, it is possible to obtain the rotation set of a generic dissipative homeomorphism from the observable or the asymptotic discretized rotation set. In other words, from the theoretical point of view, it is possible to recover numerically the rotation set of a generic homeomorphism. The generic conservative setting is quite different.

\begin{theorem}
For a generic conservative homeomorphism,
\begin{enumerate}
\item the rotation set has non-empty interior (Proposition \ref{RotGeneCons});
\item the observable rotation set consists in a single point: the mean rotation vector (Proposition \ref{ExGeneCons}). On the other hand, the asymptotic discretized rotation set coincides with the rotation set (Corollary \ref{CoroRotDiscrCons}).
\end{enumerate}
\end{theorem}

These results suggest the quite surprising moral that to recover the rotation set of a conservative homeomorphisms, it is better to do coarse roundoff errors at each iteration. More precisely, if we compute a finite number of orbits with arbitrarily good precision and long length, we will find only the mean rotation vector of the homeomorphism; but if we make roundoff errors while computing, we will be able to retrieve the whole rotation set.
\bigskip

We have performed numerical simulations to see whether these behaviours can be observed in practice or not. To obtain numerically an approximation of the observable rotation set, we have calculated rotation vectors of long segments of orbits for a lot of starting points with high precision (these points being chosen randomly). For the numerical approximation of the asymptotic discretized rotation set we have chosen a fine enough grid on the torus and have calculated the rotation vectors of periodic orbits of the discretization of the homeomorphism on this grid.

We have chosen to make these simulations on some examples where the rotation set is known to be the square $[0,1]^2$. It makes us sure of the shape of the rotation set we should obtain numerically, however it limits a bit the ``genericity'' of the examples we can produce. We also produced simulations for a homeomorphism for which we do not know the shape of the rotation set.

In the dissipative case, we made attractive the periodic points which realize the vertex of the rotation set $[0,1]^2$. It is obvious that these rotation vectors, which are realized by attractive periodic points with basin of attraction of reasonable size, will be detected by the simulations of both observable and asymptotic discretized rotation sets; that is we observe in practice: we can recover quickly the rotation set in both cases (Figures \ref{RotDissipAlea} and \ref{RotDissipDiscr}).

\begin{figure}[t]
\begin{minipage}[c]{.33\linewidth}
	\includegraphics[width=\linewidth,trim = .6cm .6cm .6cm .6cm,clip]{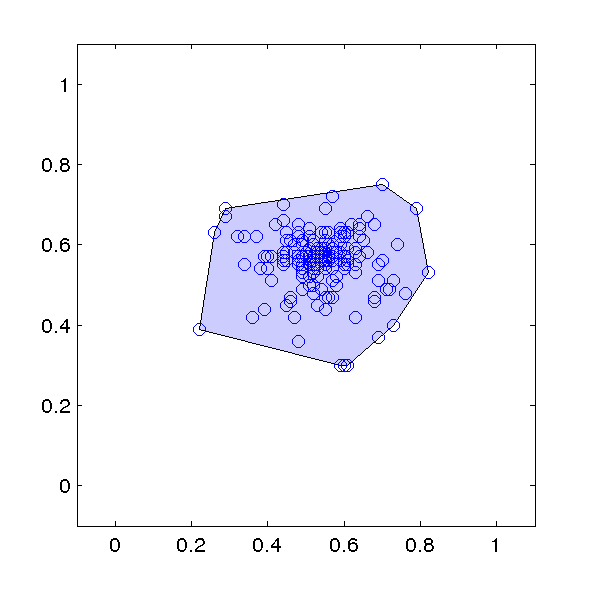}
\end{minipage}\hfill
\begin{minipage}[c]{.33\linewidth}
	\includegraphics[width=\linewidth,trim = .6cm .6cm .6cm .6cm,clip]{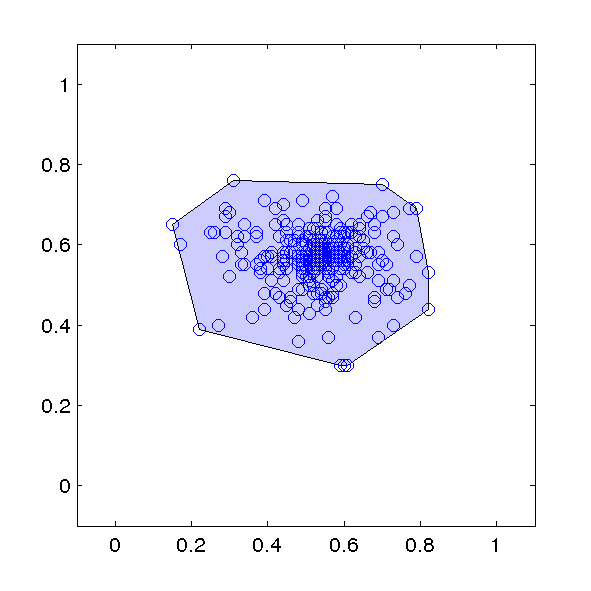}
\end{minipage}\hfill
\begin{minipage}[c]{.33\linewidth}
	\includegraphics[width=\linewidth,trim = .6cm .6cm .6cm .6cm,clip]{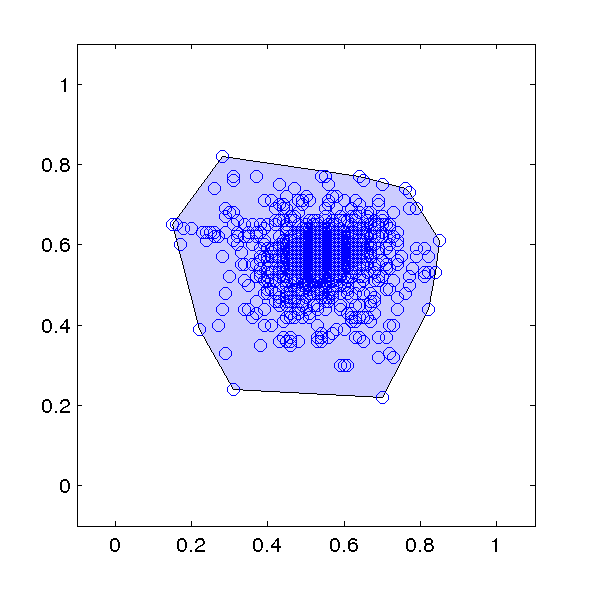}
\end{minipage}

\begin{minipage}[c]{.33\linewidth}
	\includegraphics[width=\linewidth,trim = .6cm .6cm .6cm .6cm,clip]{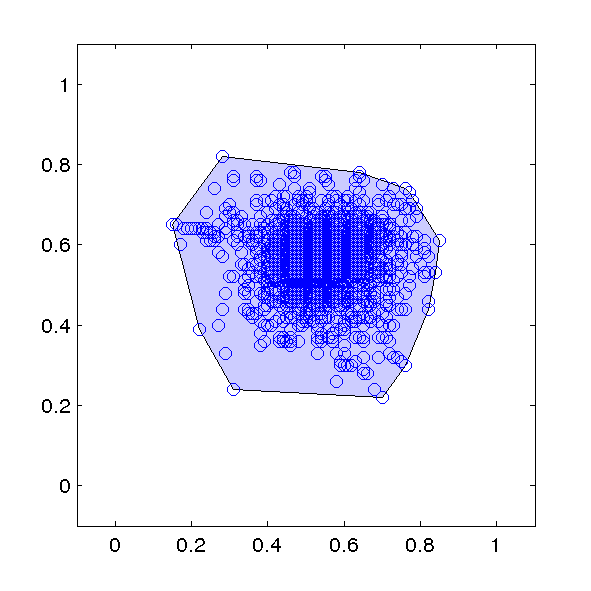}
\end{minipage}\hfill
\begin{minipage}[c]{.33\linewidth}
	\includegraphics[width=\linewidth,trim = .6cm .6cm .6cm .6cm,clip]{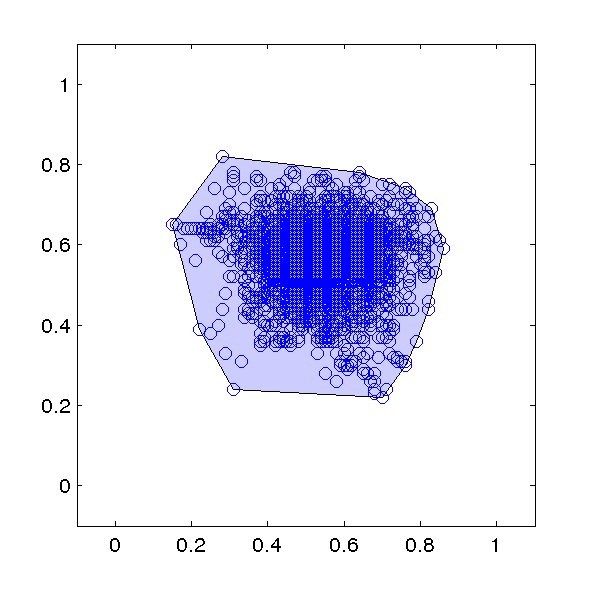}
\end{minipage}\hfill
\begin{minipage}[c]{.33\linewidth}
	\includegraphics[width=\linewidth,trim = .6cm .6cm .6cm .6cm,clip]{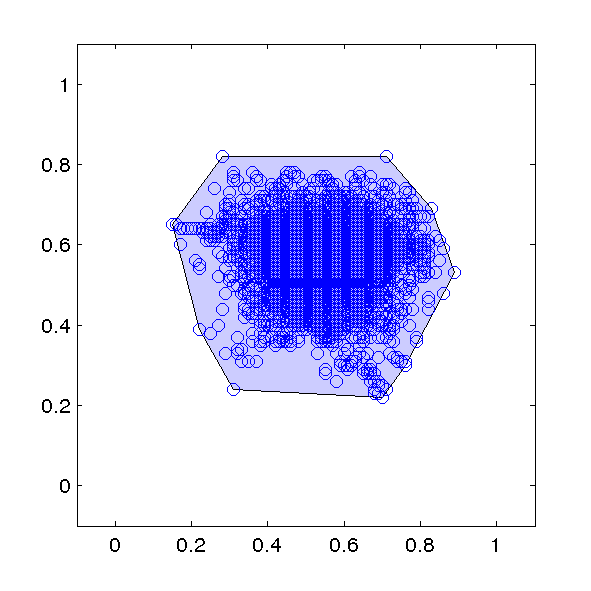}
\end{minipage}
\caption[Simulations of the asymptotic discretized rotation set of $g_3$]{Asymptotic discretized rotation set of $g_3$ as the union of the discretized rotation sets on grids $E_N$ with $100\le N\le M$ with $M = 150$ (top left), $M=200$ (top middle), $M = 500$ (top right), $M=1\,000$ (bottom left), $M=2\,000$ (bottom middle) and $M=3\,950$ (bottom right)}\label{RotLoinSerie}
\end{figure}

In the conservative setting, we observe the surprising behaviour predicted by the theory: when we compute the rotation vectors of long segments of orbits we obtain mainly rotation vectors which are quite close to the mean rotation vector, in particular we do not recover the initial rotation set. More precisely, when we perform simulations with three hours of calculation we only obtain rotation vectors close to the mean rotation vector (Figure \ref{RotConsC0Alea}). On the other hand, when we calculate the union of the discretized rotation sets over several grids to obtain a simulation of the asymptotic discretized rotation set, the rotation set is detected very quickly by the convex hulls of discretized rotation sets (less than one minute of calculation) and when we compute more orders of discretizations, we obtain a set which is quite close to $[0,1]^2$ for Hausdorff distance (Figure \ref{RotConsC0Serie}). Moreover, when we compute the observable rotation set of a homeomorphism whose rotation set is unknown, 
we obtain a single rotation vector (Figure \ref{RotLoinAlea}); but when we simulate the asymptotic discretized rotation set, then we obtain a sequence of ``thick'' sets whose convex hulls seem to converge (Figure \ref{RotLoinSerie}). As for theoretical results, this suggests the following lesson:

\emph{When we compute segments of orbits with very good precision it is very difficult to recover the rotation set. However, when we decrease the number of digits used in computations we can obtain quickly a very good approximation of the rotation set.}

This phenomenon can be explained by the fact that each grid of the torus is stabilized by the corresponding discretization of the homeomorphism. Thus, there exists an infinite number of grids such that every periodic point of the homeomorphism is shadowed by some periodic orbits of the discretizations on these grids.

\section{Notations and preliminaries}\label{sectionNotRot}

\subsection{Rotation sets}

The definition of the rotation set is made to mimic the rotation number for homeomorphisms of the circle. At first sight the natural generalisation to dimension 2 of this notion is the point rotation set, defined as follows. For every homeomorphism $f$ of the torus $\T^2$ homotopic to the identity\footnote{In this chapter every homeomorphism will be supposed homotopic to the identity.} we take a lift $F : \R^2\to\R^2$ of $f$ to the universal cover $\R^2$ of $\T^2$. The difference with the one dimensional case is that as we lose the existence of a total order on our space, the sequence $\frac{F^n(\tilde x)-\tilde x}{n}$ no longer need to converge. Thus, we have to consider all the possible limits of such sequences, called \emph{rotation vectors}; the set of rotation vectors associated to $\tilde x\in\R^2$ will be denoted by $\overline\rho(\tilde x)$: \index{$\overline\rho(\tilde x)$}
\[\overline\rho(\tilde x) = \bigcap_{N_0\in\N}\overline{\bigcup_{n\ge N_0}\left\{\frac{F^n(\tilde x)-x}{n}\right\}}.\]

Then, the \emph{point rotation set} is defined as $\rho_{pts}(F) = \bigcup_{\tilde x\in\R^2} \overline\rho(\tilde x)$\index{$\rho_{pts}(F)$}.\label{PageDefRho} Unfortunately this definition is not very convenient and it turns out that when we interchange the limits in the previous definition, we obtain the \emph{rotation set}\index{$\rho(F)$}
\[\rho(F) = \bigcap_{M\in\N}\overline{\bigcup_{m\ge M} \left\{\frac{F^m(\tilde x)-\tilde x}{m}\mid \tilde x\in\R^2\right\}}\]
which has much better properties and is easier to manipulate. In particular, it is compact and convex (see \cite{MR1053617}), and it is the convex hull of $\rho_{pts}(F)$. Moreover, it coincides with the \emph{measure rotation set}:  if we denote by $D(F)$ the \emph{displacement function}, defined on $\T^2$ by $D(F)(x) = F(\tilde x)-\tilde x$, where $\tilde x$ is a lift of $x$ to $\R^2$ (we easily check that this quantity does not depend of the lift), then (recall that $\mathcal M^f$ is the set of $f$-invariant probability measures)
\[\rho(F) = \left\{\int_{\T^2} D(F)(x)\ \ud\mu\mid \mu\in\mathcal M^f\right\}.\]

Finally, for a homeomorphism $f$ preserving $\lambda$, we denote by $\rho_{mean}(F)$\index{$\rho_{mean}(F)$} the \emph{mean rotation vector} of $F$:
\[\rho_{mean}(F) = \int_{`\T^2} D(F)(x)\ \ud\lambda(x).\]

The geometry of the rotation set of a generic dissipative homeomorphism is given by a recent result published by A. Passeggi:

\begin{theoreme}[Passeggi, \cite{rata}]\label{ThRata}
On an open and dense set of homeomorphisms $f\in\Hom(\T^2)$, the rotation set is locally constant around $f$ and is equal to a rational polygon. 
\end{theoreme}

We end this paragraph by giving a proof that if $f$ is a generic conservative homeomorphism of the torus, then $\rho(F)$ has non-empty interior.

\begin{prop}\label{RotGeneCons}
If $f$ is generic\footnote{In fact on a open dense subset of $\Hom(\T^2,\lambda)$.} among $\Hom(\T^2,\lambda)$, then $\rho(F)$ has non-empty interior.
\end{prop}

\begin{rem}
We do not know the shape of the boundary of the rotation set of a generic conservative homeomorphism. In particular we do not know if it is a polygon or not.
\end{rem}

\begin{proof}[Proof of Proposition \ref{RotGeneCons}]
We use an argument due to S. Crovisier. If $\rho(F)$ consists in a single point, we use classical perturbation techniques for conservative homeomorphisms to create a persistent periodic point $x$ for $f$. Then, by composing by a small rotation of the torus, we can move a little the mean rotation vector; in particular as the rotation set still contains the rotation vector of the periodic point $x$, it is not reduced to a single point. Now if the rotation set is a segment, by a $C^0$ ergodic closing lemma we can create a persistent periodic point whose rotation vector is close to the mean rotation vector in the following way. A small perturbation allows us to suppose that the homeomorphism we obtained, still denoted by $f$, is ergodic (it is the Oxtoby-Ulam theorem, see \cite{Oxto-meas}). We then choose a recurrent point $y\in\T^2$ which verifies the conclusion of Birkhoff's theorem: for $N$ large enough, the measure $\frac{1}{N}\sum_{k=0}^{N-1} \delta_{f^k(y)}$ is close to the measure $\lambda$. As this point 
is recurrent, by making a little perturbation, we can make it periodic (like in the proof of Lemma \ref{Laxergod}) and even persistent (see Definition \ref{rétine} and Theorem \ref{extension-sphères} page \pageref{extension-sphères}); by construction $\overline\rho (y)$ is close to the mean rotation vector. We now have two persistent periodic points, say $x$ and $y$, whose rotation vectors are different. It then suffices to compose by an appropriate rotation such that the mean rotation vector goes outside of the line generated by these two rotation vectors, and to repeat the construction to find a persistent periodic point whose rotation vector is close to this new mean rotation vector. Thus, we obtain a homeomorphism $g$ which is arbitrarily close to $f$ and possesses three periodic points $x$, $y$ and $z$ whose rotation vectors are non-aligned; therefore the rotation set of this homeomorphism has nonempty interior. Moreover, as the periodic points $x$, $y$ and $z$ are persistent, this property remains true 
on a neighbourhood of $g$, which concludes the proof for $\Hom(\T^2,\Leb)$.
\end{proof}

\subsection{Observable measures}

From the ergodic point of view, we could be tempted to define the observable rotation set to be the set of rotation vectors associated to physical measures (see Definition \ref{sport}), which are defined to express which measures can be observed in practice. However, such measures do not need to exist for every dynamical system, in this case the associated observable rotation set would be empty. To solve this problem of non existence of physical measures, E. Catsigeras and H. Enrich have defined in \cite{MR2852870} the weaker notion of \emph{observable measure}:

\begin{definition}\label{DefObs}
A probability measure $\mu$ is \emph{observable} for $f$ if, for every $\varep>0$, the set\index{$A_\varep(\mu)$}
\begin{equation}\label{defAep}
A_\varep(\mu) = \{x\in \T^2\mid \exists \nu\in p\omega(x) : \dist(\nu,\mu)<\varep \}
\end{equation}
has $\lambda$-positive measure (recall that $p\omega(x)$ is the set of Birkhoff limits of $x$, see page~\pageref{pomega}). The set of observable measures is denoted by $\Obs(f)$\index{$\Obs(f)$}.
\end{definition}

The interesting property of these measures is that, unlike physical measures, they always exist. More precisely, the set $\Obs(f)$ is a non-empty compact subset of the set of invariant measures of $f$ containing the set of physical measures (see \cite{MR2852870}).

\begin{rem}\label{remConj}
The behaviour of observable measures is compatible with topological conjugacy in the following sense: if $\mu$ is observable for $f$ and $h$ is a homeomorphism which preserves null sets, then $h^*\mu$ is observable for $hfh^{-1}$.
\end{rem}

\begin{ex}\label{ExObsMes}
\begin{enumerate}
\item If $f = \operatorname{Id}$, then $\Obs(f) = \{\delta_x\mid x\in X\}$, but $f$ has no physical measure.
\item If a dynamical system possesses a collection of physical measures whose basins of attraction cover almost all the phase space $X$ (for example if it is ergodic with respect to a smooth measure), then the set of physical measures coincides with the set of observable measures.
\end{enumerate}
\end{ex}

\begin{prop}\label{PointDissip}
If $f$ is generic among $\Hom(\T^2)$, then
\[\Obs(f) = \operatorname{Cl}\{\delta_\omega\mid \omega\text{ is a Lyapunov stable periodic orbit}\},\]
where $\operatorname{Cl}$ denotes the closure.
\end{prop}

Thus, a generic homeomorphism $f$ has a lot of observable measures\footnote{The set of Lyapunov stable periodic orbits is a Cantor set.}, but no physical measure (it is a direct consequence of the shredding lemma \ref{déchet}, see \cite{MR3027586}).

\begin{proof}[Proof of Proposition \ref{PointDissip}]
The first inclusion is easy: it suffices to remark that every stable measure supported by a Lyapunov stable periodic orbit is observable.

For the other inclusion , let $f$ be a generic dissipative homeomorphism, $\mu\in \Obs(f)$ and $\varep>0$. By hypothesis $\lambda(A_\varep(\mu))>0$ (see Equation \eqref{defAep}), then $\varep' = \frac12 \min(\varep, \lambda(A_\varep(\mu)))>0$. As $f$ is generic, it satisfies the conclusions of the shredding lemma (see Section \ref{SecShre}) applied to $f$ and $\varep'$, in particular there exists a Borel set $B\subset A_\varep(\mu)$ and an open set $O\subset \T^2$ such that:
\begin{itemize}
\item $\lambda (B) > 0$;
\item $O$ is strictly periodic: $\exists i>0 : f^i(O)\subset\subset O$;
\item $\diam (O)<\varep'$,
\item every orbit of every point of $B$ belongs to $O$ eventually.
\end{itemize}
By Lemma \ref{LyapGene} page \pageref{LyapGene}, $O$ contains a Lyapunov stable periodic point whose orbit is denoted by $\omega$; thus for every $x\in B$ and every $\nu\in p\omega(x)$, we have $\dist(\nu,\delta_\omega)<\varep'$. But by hypothesis $\dist(\nu,\mu)<\varep$, then $\dist(\mu, \delta_\omega)<2\varep$, with $\omega$ a Lyapunov stable periodic orbit.
\end{proof}

\begin{lemme}\label{PointCons}
If $f$ is generic among $\Hom(\T^2,\lambda)$, then $\Obs(f) = \{\lambda\}$ coincide with the set of physical measures.
\end{lemme}

\begin{proof}[Proof of Lemma \ref{PointCons}]
A classical theorem of J. Oxtoby and S. Ulam \cite{Oxto-meas} states that a generic conservative homeomorphism $f\in \Hom(\T^2,\lambda)$ is ergodic with respect to the measure $\lambda$. But Remark 1.8 of \cite{MR2852870} states that if the measure $\lambda$ is ergodic, then $\Obs(f) = \{\lambda\}$.
\end{proof}

\section{Observable rotation sets}

\subsection{Definitions}

As said before, from the notion of observable measure, it is easy to define a notion of observable ergodic rotation set. Another definition, more topologic, seemed reasonable to us for observable rotation sets:

\begin{definition}\label{DefRotObs}
\[\rho^{obs}(F) = \Big\{ v\in\R^2 \mid \forall \varep>0,\, \lambda\big\{x\mid \exists u\in \overline\rho (x) : d(u,v)<\varep\big\}>0\Big\}.\]\index{$\rho^{obs}(F)$}
\[\rho^{obs}_{mes}(F) = \left\{\int_{\T^2} D(F)(x) \ud \mu(x) \mid \mu\in \Obs(f)\right\}.\]\index{$\rho^{obs}_{mes}(F)$}
\end{definition}

These two sets are non-empty compact subsets of the classical rotation set, and the first one is even a subset of $\rho_{pts}(F)$. The next lemma states that these two definitions coincide:

\begin{lemme}\label{equal}
$\rho^{obs}_{mes}(F) = \rho^{obs}(F)$.
\end{lemme}

\begin{proof}[Proof of Lemma \ref{equal}]
We first prove that $\rho^{obs}_{mes}(F) \subset \rho^{obs}(F)$. Let $v\in \rho^{obs}_{mes}(F)$ and $\varep>0$. Then there exists $\mu\in \Obs(f)$ such that $v = \int_{\T^2} D(F) \ud \mu$, in particular $\lambda(A_{\varep/2}(\mu))>0$. But if $x\in A_{\varep/2}(\mu)$, then there exists a strictly increasing sequence of integers $(n_i(x))_i$ such that for every $i\ge 0$,
\[\dist \left(\frac{1}{n_i(x)}\sum_{k=0}^{n_i(x)-1} \delta_{f^k(x)},\mu\right)<\varep.\]
Thus,
\[\left|\frac{1}{n_i(x)}\sum_{k=0}^{n_i(x)-1} D(F)(f^k(x)) - \int_{\T^2} D(F) \ud \mu\right|<\varep,\]
in other words the inequality
\[\left|\frac{F^{n_i(x)}(x)-x}{n_i(x)} - v\right|<\varep\]
holds for every $i$ and on a $\lambda$-positive measure set of points $x$.
\bigskip

For the other inclusion, let $v\in \rho^{obs}(F)$ and set
\[\tilde A_\varep(v) = \{x\in \T^2 \mid \exists u\in\overline\rho(x) : d(u,v)<\varep\}.\]
By hypothesis, $\lambda(\tilde A_\varep(v))>0$ for every $\varep>0$. To each $x\in \tilde A_\varep(v)$ we associate the set $p\omega_\varep^v(x)$ of limit points of the sequence of measures 
\[\frac{1}{n_i(x)}\sum_{k=0}^{n_i(x)-1} \delta_{f^k(x)},\]
where  $(n_i(x))_i$ is a strictly increasing sequence such that
\[\left| \frac{F^{n_i(x)}(x)-x}{n_i(x)}-v\right| <\varep.\]
By compactness of $\Pb$, the set $p\omega_\varep^v(x)$ is non-empty and compact. In the sequel we will use the following easy remark: if $0<\varep<\varep'$ and $x\in \tilde A_\varep$, then $p\omega_\varep^v(x)\subset p\omega_{\varep'}^v(x)$.

By contradiction, suppose that for every $\mu\in\Pb$, there exists $\varep_\mu>0$ such that
\[\lambda\big\{ x\in \tilde A_{\varep_\mu}(v) \mid \exists \nu\in p\omega_{\varep_\mu}^v(x) : \dist(\nu,\mu)<\varep_\mu \big\} =0.\]
By compactness, $\Pb$ is covered by a finite number of balls $B(\mu_j,\varep_{\mu_j})$. Taking $\varep = \min \varep_{\mu_j}$, for every $j$ we have
\[\lambda\big\{ x\in \tilde A_\varep(v) \mid \exists\nu\in p\omega_\varep^v(x) : \dist(\nu,\mu)<\varep_{\mu_j}\big\} =0,\]
thus, as balls $B(\mu_j,\varep_{\mu_j})$ cover $\Pb$,
\[\lambda\big\{ x\in \tilde A_\varep(v) \mid p\omega_\varep^v(x)\cap \Pb\neq\emptyset \big\} =0,\]
which is a contradiction.

Therefore, there exists $\mu_0\in\Pb$ such that for every $\varep>0$,
\[\lambda\big\{ x\in \tilde A_{\varep}(v) \mid \exists\nu\in p\omega_{\varep}^v(x) : \dist(\nu,\mu_0)<\varep \big\} > 0,\]
in particular $\mu_0\in\Obs(f)$. Furthermore, for $\varep>0$, there exists $x\in \tilde A_{\varep}(v)$ and $\mu_x\in p\omega_\varep^v(x)$ such that $\dist(\mu_x,\mu_0)<\varep$. As $\mu_x\in p\omega_\varep^v(x)$, there exists a sequence $(n_i(x))_i$ such that
\[\dist\left(\mu_x\,,\, \frac{1}{n_i(x)}\sum_{k=0}^{n_i(x)-1} \delta_{f^k(x)}\right)<\varep \qquad \text{and} \qquad \left|\frac{F^{n_i(x)}(x)-x}{n_i(x)} - v\right| < \varep.\]
Thus,
\[\dist \left(\mu_0\,,\, \frac{1}{n_i(x)}\sum_{k=0}^{n_i(x)-1} \delta_{f^k(x)}\right)<2\varep.\]
Integrating this estimation according to the function $D(F)$, we obtain:
\[\left|\int_{\T^2} D(F)\ud\mu_0 - \frac{F^{n_i(x)}(x)-x}{n_i(x)}\right| < 2\varep,\]
so
\[\left|\int_{\T^2} D(F)\ud\mu_0 - v\right| < 3\varep,\]
for every $\varep>0$, in other words,
\[v = \int_{\T^2} D(F)\ud\mu_0.\]
\end{proof}


\subsection{Properties of the observable rotation set}

We begin by giving two lemmas which state the dynamical behaviour of the observable rotation sets.

\begin{lemme}\label{iter}
For every $q\in \N$, $\rho^{obs}(F^q) = q\rho^{obs}(F)$.
\end{lemme}

\begin{proof}[Proof of Lemma \ref{iter}]
It suffices to remark that $\overline\rho_{F^q}(x) = q\overline\rho_{F}(x)$ (one inclusion is trivial and the other is easily obtained by Euclidean division).
\end{proof}

\begin{rem}
In general $\rho^{obs}(F^{-1}) \neq -\rho^{obs}(F)$: see for instance the point \ref{F-1}. of Example \ref{exrot}.
\end{rem}

\begin{lemme}\label{Dynamic}
If $H$ is a homeomorphism of $\R^2$ commuting with integral translations and preserving null sets, then $\rho^{obs}(H\circ F\circ H^{-1}) = \rho^{obs}(F)$.
\end{lemme}

\begin{proof}[Proof of Lemma \ref{Dynamic}]
It follows easily from the fact that the notion of observable measure is stable by conjugacy (see Remark \ref{remConj}).
\end{proof}

We now give a few simple examples of calculation of observable rotation sets.

\begin{ex}\label{exrot}
\begin{enumerate}
\item If $f = \operatorname{Id}$, then $\rho^{obs}(F) = \{(0,0)\}$.
\item If
\[F(x,y) = \left( x+\cos(2\pi y)\, ,\, y\right),\]
then $\rho_{pts}(F) = \rho^{obs}(F) = [-1,1]\times \{0\}$.
\item\label{F-1} If
\[F(x,y) = \left( x+\cos(2\pi y)\, ,\, y+\frac{1}{100}\sin(2\pi y)\right),\]
then $\rho_{pts}(F) = \{(0,-1),(0,1)\}$, but $\rho^{obs}(F) = \{(0,-1)\}$ and $\rho^{obs}(F^{-1}) = \{(0,1)\}$.
\item Let
\[P\begin{pmatrix} x\\ y \end{pmatrix} = \begin{pmatrix} x + \frac12 \cos(2\pi y)+1\\ y \end{pmatrix} \quad \text{and} \quad Q\begin{pmatrix} x\\ y \end{pmatrix} = \begin{pmatrix} x\\ y + \frac12 \cos(2\pi x)+1\end{pmatrix}.\]
Then the rotation set of the (conservative) homeomorphism $F= P\circ Q$ is equal to $[0,1]^2$. Moreover, we can perturb $F$ into a (conservative) homeomorphism $\tilde F$ such that $\tilde F$ is the identity on the neighbourhoods of the points whose coordinates belong to $1/2 \Z$ (applying for example Theorem \ref{extension-sphères}). Then, the vertices of the square $[0,1]^2$ belong to the observable rotation set of $\tilde F$.
\item Let $P$ be a convex polygon with rational vertices. In \cite{MR1176627}, J. Kwapisz has constructed an axiom A diffeomorphism $f$ of $\T^2$ whose rotation set is the polygon $P$. It is possible to modify slightly Kwapisz's construction so that all the sinks of $f$ are fixed points, and so that the union of the basins of these sinks have $\lambda$-full measure. Hence, the observable rotation set of $f_P$ is reduced to $\{(0,0)\}$.
\end{enumerate}
\end{ex}
\bigskip

We now give the results about the link between the rotation set and the observable rotation set in the generic setting. We begin by the dissipative case.

\begin{prop}\label{GeneDissip}
If $f$ is generic among $\Hom(\T^2)$, then $\rho(F) = \conv(\rho^{obs}(F))$. If moreover $f$ is generic with a non-empty interior rotation set, then $\rho(F) = \rho^{obs}(F)$.
\end{prop}

To prove this proposition, we will use the following lemma, which is a direct consequence of Proposition \ref{PointDissip}.

\begin{lemme}\label{ExGeneDissip}
If $f$ is generic among $\Hom(\T^2)$, then
\[\rho^{obs}(F) = \operatorname{Cl}\{\overline\rho(\tilde x)\mid x\text{ is a Lyapunov stable periodic point}\}.\]
\end{lemme}

We will also need a theorem of realization of rotation vectors by periodic points.

\begin{theoreme}[J. Franks, \cite{MR958891}]\label{Franks}
For every $f\in\Hom(\T^2)$, every rational point of the interior of $\rho(F)$ is realized as the rotation vector of a periodic point of the homeomorphism~$f$.
\end{theoreme}

\begin{proof}[Proof of Proposition \ref{GeneDissip}]
Theorem \ref{ThRata} states that for an open dense set of homeomorphisms, the rotation set is a rational polygon. Then, a theorem of realization of J. Franks \cite[Theorem 3.5]{MR967632} implies that every vertex of this polygon is realized as the rotation vector of a periodic point of the homeomorphism, which can be made attractive by a little perturbation of the homeomorphism. Then generically we can find a Lyapunov stable periodic point which shadows the previous periodic point (by Lemma \ref{LyapGene}), in particular it has the same rotation vector. Thus every vertex of $\rho(F)$ belongs to $\rho^{obs}(F)$ and $\rho(F) = \conv(\rho^{obs}(F))$.

For $\varep>0$, we can find a finite $\varep$-dense subset $R_\varep$ of $\rho(F)$ made of rational points. Thus, Theorem \ref{Franks} associates to each of these rational vectors a periodic point of the homeomorphism which realizes this rotation vector; we can even make these periodic points of the homeomorphism attractive. Thus, for every $\varep>0$, the set $O_\varep$ made of the homeomorphisms such that every vector of $R_\varep$ is realized by a strictly periodic open subset of $\T^2$ is open and dense in the set of homeomorphisms with non-empty interior rotation set. Applying Lemma \ref{LyapGene} we find a $G_\delta$ dense subset of $O_\varep$ on which every strictly periodic open subset of $\T^2$ contains a Lyapunov stable periodic point; on this set the Hausdorff distance between $\rho(F) = \rho^{obs}(F)$ is smaller than $\varep$. The conclusion of the proposition then easily follows from Baire theorem.
\end{proof}

\begin{rem}
It is not true that $\rho(F) = \rho^{obs}(F)$ holds for a generic homeomorphism: see for instance the point \ref{F-1} of Example \ref{exrot}, where on a neighbourhood of $f$ the set $\rho^{obs}$ is contained in a neighbourhood of the points $(0,-1)$ and $(0,1)$.
\end{rem}

For the conservative case, we recall the result of Proposition \ref{RotGeneCons}: the rotation set of a generic conservative homeomorphism has non-empty interior. The following result states that in this case the observable rotation set is much smaller, more precisely it consists in a single vector, namely the mean rotation vector.

\begin{prop}\label{ExGeneCons}
If $f$ is generic among $\Hom(\T^2,\lambda)$, then $\rho^{obs}(F) = \{\rho_\lambda(F)\}$, where $\rho_\lambda(F)$ is the mean rotation vector with respect to the measure $\lambda$.
\end{prop}

Thus, for almost every $x\in\T^2$ (with respect to the measure $\lambda$), the set $\overline\rho(x)$ is reduced to a single point which is the mean rotation vector.

\begin{proof}[Proof of Proposition \ref{ExGeneCons}]
It is easily implied by the fact that the measure $\lambda$ is the only observable measure (Lemma \ref{PointCons}, which easily follows from Oxtoby-Ulam theorem).
\end{proof}

\section{Discretized rotation sets}\label{DiscRot}

We now take into account the fact that the computer has a finite digital precision. It will be the occasion to apply the techniques of proof  explained in Chapter \ref{ChapCons}.

The \emph{discretized rotation set} is defined as follows. Consider a lift $F : \R^2 \to\R^2$ of $f$ and a lift $\tilde E_N$ of the grid $E_N$ to $\R^2$. Then\index{$\rho(F_N)$}
\[\rho(F_N) = \bigcap_{M\in\N}\overline{\bigcup_{m\ge M} \left\{\frac{F_N^m(\tilde x)-\tilde x}{m}\mid \tilde x\in\R^2\right\}}.\]
Remark that this set coincides with the set of rotation vectors of the periodic orbits of $f_N$. Then the \emph{asymptotic discretized rotation set} is the upper limit of the sets $\rho(F_N)$:\index{$\rho^{discr}(F)$}
\[\rho^{discr}(F) = \bigcap_{M\in\N}\overline{\bigcup_{N\ge M} \rho(F_N)}.\]

The first result is that for every homeomorphism $f$, the discretized rotation set $\rho(F_N)$ is almost included in the rotation set $\rho(F)$ when $N$ is large enough. This property follows easily with a compactness argument from the convergence of the sequence $f_N$ to the homeomorphism $f$ (for example for the Hausdorff distance on the graphs of these maps).

\begin{prop}\label{InclPts}
For every homeomorphism $f$ and every $\varep>0$, it exists $N_0\in\N$ such that for every $N\ge N_0$, we have $\rho(F_N) \subset B(\rho(F),\varep)$, where $B(\rho(F),\varep)$ denotes the set of points whose distance to $\rho(F)$ is smaller than $\varep$. In particular $\rho^{discr}(F)\subset\rho(F)$.
\end{prop}

\begin{proof}[Proof of Proposition \ref{InclPts}]
By definition of the rotation set, for $\varep>0$ there exists $m\in\N$ such that
\[\left\{\frac{F^m(\tilde x)-\tilde x}{m}\mid \tilde x\in\R^2\right\} \subset B(\rho(F),\varep).\]
Then there exists $N_0\in\N$ such that for every $N\ge N_0$, 
\[\left|\frac{F^m(\tilde x)-\tilde x}{m} - \frac{F_N^m(\tilde x_N)-\tilde x_N}{m}\right| \le \varep.\]
This allows us to handle the case of long periodic orbits of the discretizations: by euclidean division, each periodic orbit of $f_N$ of length bigger than $m/\varep$ will be in the $\varep$ neighbourhood of the convex hull of the set
\[\frac{F_N^m(\tilde x_N)-\tilde x_N}{m},\]
so in the $3\varep$-neighbourhood of the rotation set $\rho(f)$.

For short orbits we argue by contradiction: suppose that there exists $\varep>0$ such that for every $N_0\in\N$ there exists $N\ge N_0$ and $x_N\in E_N$ which is periodic under $f_N$ with period smaller than $m/\varep$ and whose associated rotation vector is not in $B(\rho(F),\varep)$. Then up to take subsequences these periodic points $x_N$ have the same period and converge to a periodic point $x\in\T^2$ whose associated rotation vector (for $F$) is not in $B(\rho(F),\varep)$, which is impossible.
\end{proof}

The other inclusion depends on the properties of the map $f$. We begin by the dissipative case.

\begin{prop}\label{RotDiscrDissip}
If $f$ is generic among $\Hom(\T^2)$, then $\conv\big(\rho(F_N)\big)$ tends to $\rho(F)$ for the Hausdorff topology. In particular $\conv\big(\rho^{discr}(F)\big) = \rho(F)$. Moreover, if $\rho(F)$ has nonempty interior, then there is no need to take convex hulls.
\end{prop}

\begin{proof}[Proof of Proposition \ref{RotDiscrDissip}]
The fact that the upper limit of $\rho(F_N)$ is included in $\rho(F)$ follows directly from Lemma \ref{InclPts}.

It remains to prove that the lower limit of $\conv\big(\rho(F_N)\big)$ contains $\rho(F)$. To do that, we prove that $\rho(F_N)$ converges to $\rho^{obs}(F)$. First of all the rotation set is the closure of the convex hull of the rotation vectors of Lyapunov stable periodic points (Propositions \ref{PointDissip} and \ref{GeneDissip}). To each one of these points we can associate a periodic closed set $K$ with non-empty interior and with period $\tau$ which has the same rotation vector. Then there exists an open set $O\subset K$ such that for $N$ large enough and $x\in K$ we also have $f_N^\tau(x_N) \in O\subset K$. Thus there exists $i\in\N^*$ such that $f_N^{\tau i}(x_N) = f_N^{2\tau i}(x_N)$ and $f_N^{\tau i}(x_N)$ has the same rotation vector as $K$, thus the same rotation vector as the initial Lyapunov stable periodic point.
\end{proof}

For the conservative case, with the same techniques as in Chapter \ref{ChapCons}, we can prove the following result:

\begin{lemme}\label{EnsRotDiscrCons}
If $f$ is generic among $\Hom(\T^2,\lambda)$, then for every finite collection of rotation vectors $\{v_1,\cdots,v_n\}$, each one realized by a periodic orbit of $f$, there exists a subsequence $f_{N_i}$ of discretizations such that for every~$i$, $\rho_{N_i}(f) = \{v_1,\cdots,v_n\}$.
\end{lemme}

\begin{proof}[Proof of Lemma \ref{EnsRotDiscrCons}]
We denote by $\mathcal D_q$ the set of subsets of $\Q^2$ made of elements whose coordinates are of the type $p'/q'$, with $0<q'<q$ and $-q^2<p'<q^2$. Consider the set
\[\bigcap_{q,N_0}\bigcap_{D\in\mathcal D_q}\bigcup_{N\ge N_0} \left\{\begin{array}{r}
f\in\Hom(\T^2,\lambda) \mid(\forall v\in D,\,v\text{ is realised by a}\\
\text{persistent periodic point of $f$}) \implies \rho(F_N) = D 
\end{array}\right\}.
\]
To prove the lemma it suffices to prove that this set contains a $G_\delta$ dense. This is obtained with the same kind of proof as for Proposition~\ref{var1}. 

Let $f\in\Hom(\T^2,\lambda)$, $\varep>0$,  $q, N_0\in\N$ and $D\in \mathcal D_q$. We suppose that for all $v\in D$, $v$ is realizable by a persistent periodic orbit $\omega_i$ of $f$. For all of these orbits $\omega_1,\cdots,\omega_\ell$, we denote by $p_i$ the length of the orbit $\omega_i$ and choose a point $x_i$ belonging to $\omega_i$. We then apply Lax's theorem (Theorem \ref{Lax}): if $N$ is large enough, then there exists a cyclic permutation $\sigma_N$ of $E_N$ such that $d_N(f,\sigma_N)<\varep$. If $N$ is large enough, then the families 
\[\big\{(x_1)_N,\cdots,\sigma_N^{p_1-1}((x_1)_N)\big\},\, \cdots,\, \big\{(x_\ell)_N,\cdots,\sigma_N^{p_\ell-1}((x_\ell)_N)\big\}\]
are disjoint and satisfy $d\big((x_i)_N, \sigma_N^{p_i-1}((x_i)_N))<\varep$ for all $i$. We then use the same technique as in the proof of Proposition~\ref{var1} (see also Figure \ref{trajectoire}) to close each orbit $\{(x_i)_N,\cdots,\sigma_N^{p_i-1}((x_i)_N)\}$. The discrete map $\sigma'_N$ we obtain has then exactly $\ell$ periodic orbits, and each of them has the same rotation vector as the corresponding real periodic orbit of $f$. As in the proof of Lemma \ref{lemmetrans}, we then use the proposition of finite maps extension (Proposition \ref{extension}) to build a homeomorphism $g$ which is $\varep$-close to $f$ and whose discretization $g_N$ satisfies $\rho(G_N) = D$, and moreover we can suppose that this occurs on a whole neighbourhood of $g$.
\end{proof}

The combination of the realisation theorem of J. Franks \cite[Theorem 3.2]{MR958891} and the fact that for a generic conservative homeomorphism the rotation set has non-empty interior (Proposition \ref{RotGeneCons}) leads to the following corollary.

\begin{coro}\label{CoroRotDiscrCons}
If $f$ is generic among $\Hom(\T^2,\lambda)$, then for every compact subset $K$ of the rotation set of $F$ there exists a subsequence $f_{N_i}$ of discretizations such that $\rho_{N_i}(F)$ tends to $K$ for the Hausdorff topology. In particular $\rho^{discr}(F) = \rho(F)$.
\end{coro}

\section{Numerical simulations}\label{SecNumRot}

We have conducted numerical simulations of the rotation sets associated to both dissipative and conservative homeomorphisms. For the first examples we treated, we have made the deliberate choice to choose homeomorphisms whose rotation set is known to be the square $[0,1]^2$. Of course these homeomorphisms are not the best candidates for ``generic'' homeomorphisms, but at least we are sure of what is the shape of the rotation set we want to obtain.

As an example of dissipative homeomorphism we have taken $f_1 = R_1\circ Q_1\circ P_1$, and for the conservative homeomorphism we have chosen the very similar expression $g_1 = Q_1\circ P_1$, where
\begin{align*}
P_1(x,y) = & \Big(x\ , \ y+\frac12\big(\cos(2\pi(x+\alpha))+1\big)\\
           & + 0.0234\sin^2(4\pi(x+\alpha))\big(\sin(6\pi(x+\alpha))+0.3754\cos(26\pi(x+\alpha))\big)\Big),
\end{align*}
\begin{align*}
Q_1(x,y) = & \Big(x + \frac12\big(\cos(2\pi (y+\beta))+1\big)\\
			     & + 0.0213\sin^2(4\pi (y+\beta))\big(\sin(6\pi(y+\beta))+0.4243\cos(22\pi(y+\beta))\big)\  , \ y\Big),
\end{align*}
\begin{align*}
R_1(x,y) = & \big(x-0.0127\sin(8\pi (x+\alpha)) + 0.000324\sin(33\pi (x+\alpha))\ ,\\
           & y - 0.0176\sin(12\pi (y+\beta)) + 0.000231\sin(41\pi y)\big),
\end{align*}
with $\alpha = 0.00137$ and $\beta = 0.00159$.

The homeomorphisms $P_1$ and $Q_1$ are close to the homeomorphisms
\[\tilde P(x,y) = \Big(x\ , \ y+\frac12\big(\cos(2\pi(x+\alpha))+1\big)\Big)\]
and
\[\tilde Q(x,y) = \Big(x + \frac12\big(\cos(2\pi (y+\beta))+1\big)\ ,\ y\Big);\]
it can easily be seen that the rotation set of the homeomorphism $\tilde Q\circ \tilde P$ is the square $[0,1]^2$, whose vertices are realized by the points $(0,0)$, $(0,1/2)$, $(1/2,0)$ and $(1/2,1/2)$. The perturbations $P_1$ and $Q_1$ of $\tilde P$ and $\tilde Q$ are small enough (in $C^2$ topology) to ensure that the rotation set remains the square $[0,1]^2$; these perturbations are made in order to make $f_1$ ``more generic'' (in particular, the periodic orbits whose rotation vectors realize the vertices of the square do not belong to the grids). The key property of the homeomorphism $R_1$ is that is has the fixed points of $f_1$ which realize the vertices of $[0,1]^2$ as fixed attractive points; this creates fixed attractive points which realize the vertices of the rotation set. 

We have chosen $R_1$ to be very close to the identity in $C^1$-topology to ensure that the basins of the sinks and sources are large enough. Indeed, J.-M. Gambaudo and C. Tresser have shown in \cite{MR700317} that, even for dissipative diffeomorphisms defined by very simple formulas, sinks and sources are often undetectable in practice because the size of the their basins are too small.

We have conducted other series of simulations for two other examples of conservative homeomorphisms. The first one has an expression which is very similar to that of $g_1$, but the cosines are replaced by a piecewise affine map with the same following properties: $s$ is 1-periodic, $s(0)=1$, $s(1/2)=0$ and $s$ is affine between $0$ and $1/2$ and between $1/2$ and $1$. More precisely, we set $g_2 = Q_2\circ P_2$, with
\begin{align*}
 P_2(x,y) = \Big(x\ ,\ y & + 2s(x+\alpha) + 0.0234s(2(x+\alpha))\\
                         & + 0.0167s(10(x+\alpha)\Big);
\end{align*}
\begin{align*}
 Q_2(x,y) = \Big(x & + 2s(y+\beta) + 0.0213s(2(y+\beta))\\ 
                   & + 0.0101s(6(y+\beta))\ ,\ y \Big).
\end{align*}

The properties of $s$ imply that the rotation set of $g_2$ is also the square $[0,1]^2$; the difference with $g_1$ is that the vertices of this rotation set are no longer realized by elliptic periodic points, which makes them harder to detect.

For the last conservative homeomorphism we tested, we made ``random'' choices of the coefficients; we do not know \emph{a priori} what is its rotation set. More precisely, we took $g_3 = Q_3\circ P_3$, with 
\begin{align*}
 P_3(x,y) = \big(x\ ,\ y & + 0.3\sin(2\pi(x+0.34137))\\
                         & + 0.2\sin(3\pi(x+0.21346)) + 0.578675)\big);
\end{align*}
\begin{align*}
  Q_3(x,y) = \big(x & + 0.25\sin(2\pi(y+0.9734))\\
                    & + 0.35\sin(3\pi(y-0.20159))+0.551256\ ,\ y \big).
\end{align*}
We will test on simulations whether the computed rotation sets seem to converge or not. If so, it could be a good indication that the rotation set we obtained is close to the actual rotation set.
\bigskip

We have made two kinds of simulations of the rotation set.
\begin{itemize}
\item In the first one we have computed the rotation vectors of segments of orbits of length $1\,000$ with good precision (52 binary digits); in other words for $N$ random starting points $x\in \T^2$, we have computed $\frac{F^{1000}(x)-x}{1000}$. We have made these tests for $N= 100$, which takes around 1s of calculation, $N=10\,000$, which takes about 2min of calculation, and $N=1\,000\,000$, which takes about 4h of calculation. This is maybe the most simple process that can be used to find numerically the rotation set. It should lead to a good approximation of the observable rotation set; in particular, Proposition \ref{GeneDissip} suggests that, for the dissipative homeomorphism $f_1$, we should obtain a set which is close (for Hausdorff distance) to the square $[0,1]^2$, and if not at least a set whose convex hull is this square. On the other hand, for the conservative homeomorphisms $g_i$, Proposition \ref{ExGeneCons} suggests that we should only obtain the mean rotation vector, which is close to $(1/2,
1/2)$.
\item In the second kind of simulations we have computed the rotation vectors of the periodic orbits of the discretization $(f_i)_{N}$ on a grid $N\times N$; these simulations calculate the discretized rotation sets. For each homeomorphism we have represented these sets for $N=499$, $N=500$ and $N=501$, each calculation taking about 2s of calculation. We have also computed the union of the discretized rotation sets for $100\le N\le M$, which represents the asymptotic discretized rotation set. We represent these sets for $M = 100$ ($\simeq$ 0.5s of calculation), $M=150$ ($\simeq$ 15s of calculation), $M = 200$ ($\simeq$ 45s of calculation), $M=500$ ($\simeq$ 13min of calculation), $M=1\,000$ ($\simeq$ 1h 45min of calculation), $M=2\,000$ ($\simeq$ 14h of calculation) and for $g_3$, $M=3950$ ($\simeq$ 100h of calculation). The theory tells us that in both conservative and dissipative cases, for some $N$, the discretized rotation set should be close (for Hausdorff distance) to the square $[0,1]^2$; a weaker 
property would be that its convex hull should be close to this square. Moreover this should also be true for the asymptotic discretized rotation sets.
\end{itemize}

We shall notice that these two methods are formally the same: making simulations on a grid $N\times N$ is equivalent to calculate with $-\log_2 N$ binary digits (for example about 10 for $N=1\,000$). The only difference is that for the second method we use deliberately a very bad numerical precision, which allows us to detect the actual dynamics of the discretizations.

Moreover, in practice, for a given calculation time, the calculation of the rotation set by discretization (\emph{i.e.} by the second method) allows to compute much more orbits than the other method. More precisely, the algorithm we have used to compute the asymptotic discretized rotation set visits each point of the grid $N\times N$ once. Thus, for $N^2$ starting points we only have to compute $N^2$ images of the discretization of the homeomorphism on the grid; the number of rotation vectors we obtain is simply the number of periodic orbits of the discretization. So in a certain sense this second algorithm is much faster than the naive algorithm consisting in computing long segments of orbits. All the simulations have been performed on a computer equipped with a processor Intel Core I5 2.40GHz.
\bigskip

We shall notice that the calculated rotation sets we obtained for $f_1$, $g_1$ and $g_2$ are always contained inside of the square $[0,1]^2$ (see the figures below). Indeed, the very definition of the rotation set ensures that if $T$ is large enough, then the rotation vector of every segment of orbit of length $T$ belongs to a neighbourhood of the rotation set; therefore the computed observable rotation set should be included in a small neighbourhood of the rotation set provided we have chosen a large enough length of orbit. Concerning the discretized rotation set, Proposition \ref{InclPts} ensures that if the order of discretization is large enough, then the discretized rotation set is included in a neighbourhood of the actual rotation set. However, there is no global estimation of these integers, as can be seen in the following example. Set $f : \T^2 \to \T^2$ defined by
\[f (x,y) = \big(x+\alpha\, ,\, y+\sin(2\pi x)\big),\]
with $\alpha$ an irrational number close to (say) $1/(2T)$. As $\alpha\in\R\setminus\Q$, we easily obtain that the rotation set of $f$ is equal to $(\alpha,0)$; but if we compute the rotation numbers of all the segment of orbit of length $T$, they will form a set which is close to the segment $\big[(\alpha,-2/\pi),(\alpha,2/\pi)\big]$.
\bigskip

\begin{figure}[ht]
\begin{minipage}[c]{.33\linewidth}
	\includegraphics[width=\linewidth,trim = .6cm .6cm .6cm .6cm,clip]{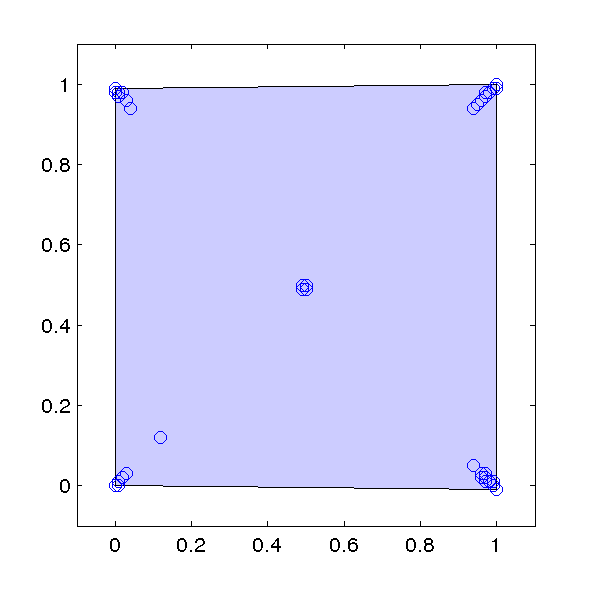}
\end{minipage}\hfill
\begin{minipage}[c]{.33\linewidth}
	\includegraphics[width=\linewidth,trim = .6cm .6cm .6cm .6cm,clip]{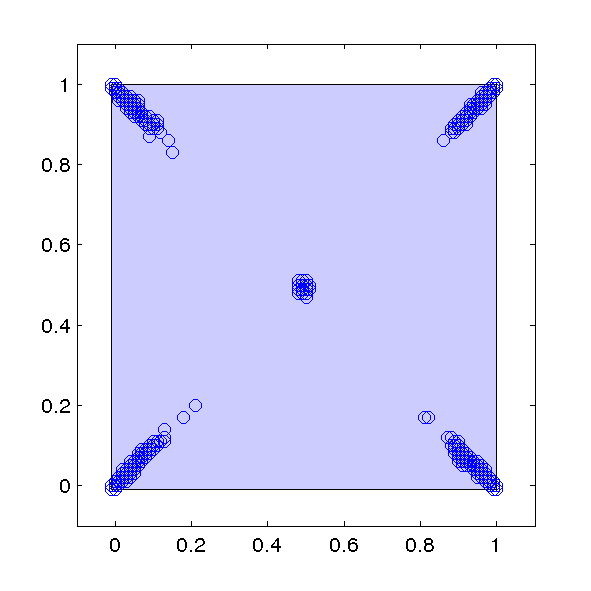}
\end{minipage}\hfill
\begin{minipage}[c]{.33\linewidth}
	\includegraphics[width=\linewidth,trim = .6cm .6cm .6cm .6cm,clip]{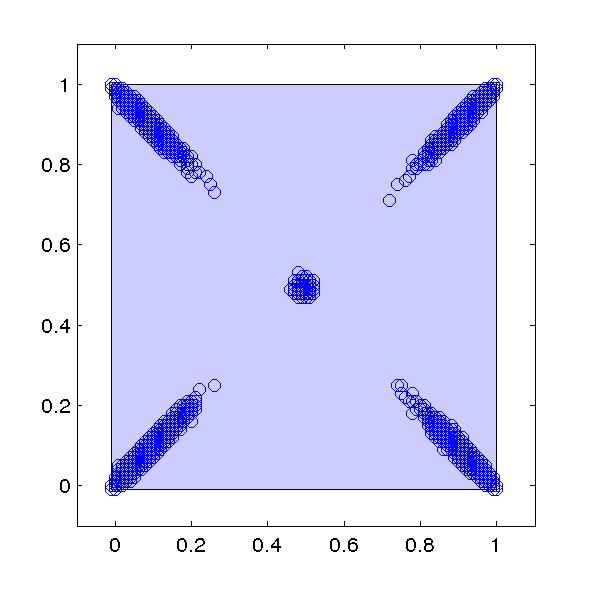}
\end{minipage}
\caption[Simulations of the observable rotation set of $f_1$]{Observable rotation set of $f_1$, $k$ orbits of length $1\,000$ with random starting points with $k=100$ (left), $10\,000$ (middle) and $1\,000\,000$ (right)}\label{RotDissipAlea}
\end{figure}

In the dissipative case, a lot of the obtained rotation vectors are close to one of the vertices of the real rotation set $[0,1]^2$ of $f_1$, the others being located around $(1/2,1/2)$ (see Figure \ref{RotDissipAlea}). This is what is predicted by the theory, in particular by Lemma \ref{ExGeneDissip}: we detect rotation vectors realized by Lyapunov stable periodic points. The fact that the rotation vectors are not located exactly on the vertices of $[0,1]^2$ can be explained by the slow convergence of the orbits to the attractive points: it may take a while until the orbit become close to one of the Lyapunov stable periodic points. We will see that this behaviour is very different from the one in the conservative case, even if the homeomorphism $f_1$ is very close to $g_1$ (approximately $10^{-2}$ close).

\begin{figure}[ht]
\begin{minipage}[c]{.33\linewidth}
	\includegraphics[width=\linewidth,trim = .6cm .6cm .6cm .6cm,clip]{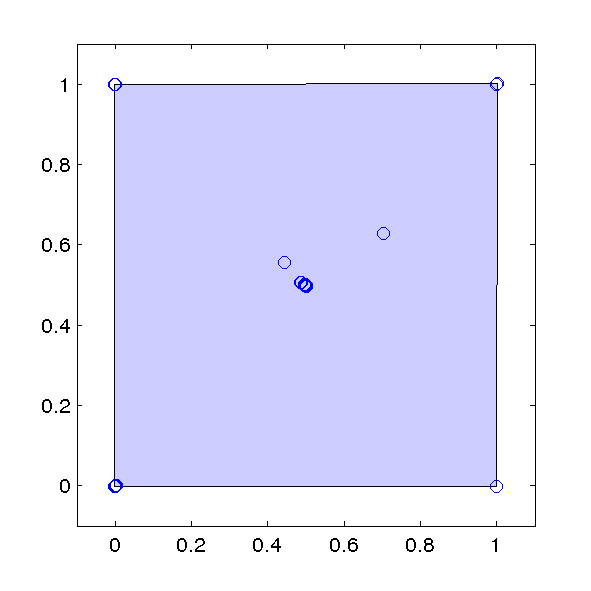}
\end{minipage}\hfill
\begin{minipage}[c]{.33\linewidth}
	\includegraphics[width=\linewidth,trim = .6cm .6cm .6cm .6cm,clip]{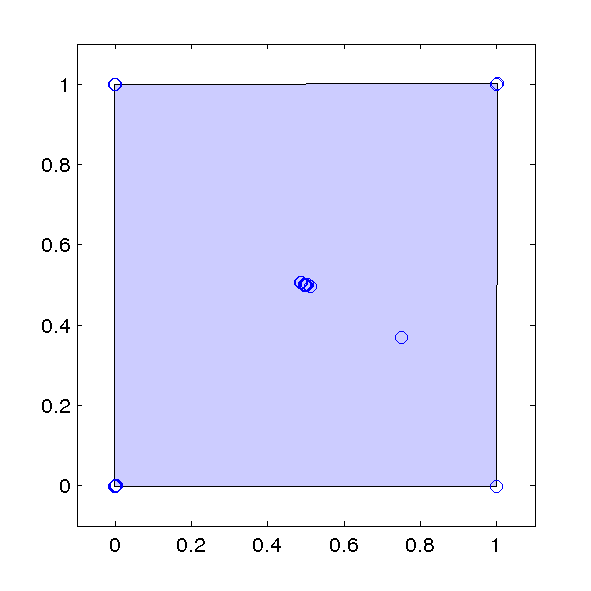}
\end{minipage}\hfill
\begin{minipage}[c]{.33\linewidth}
	\includegraphics[width=\linewidth,trim = .6cm .6cm .6cm .6cm,clip]{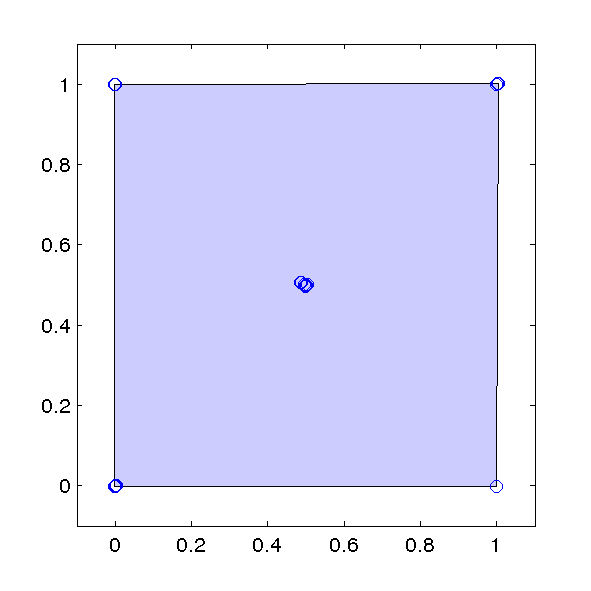}
\end{minipage}

\caption[Simulations of the discretized rotation set of $f_1$]{Discretized rotation set of $f_1$ on grids $E_N$, with $N=499$ (left), $N=500$ (middle) and $N=501$ (right)}\label{RotDissipDiscr}
\end{figure}

\begin{figure}[t]
\begin{minipage}[c]{.33\linewidth}
	\includegraphics[width=\linewidth,trim = .6cm .6cm .6cm .6cm,clip]{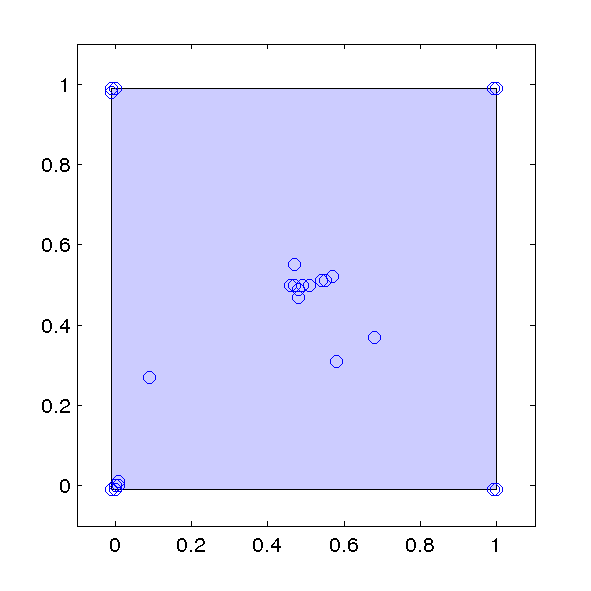}
\end{minipage}\hfill
\begin{minipage}[c]{.33\linewidth}
	\includegraphics[width=\linewidth,trim = .6cm .6cm .6cm .6cm,clip]{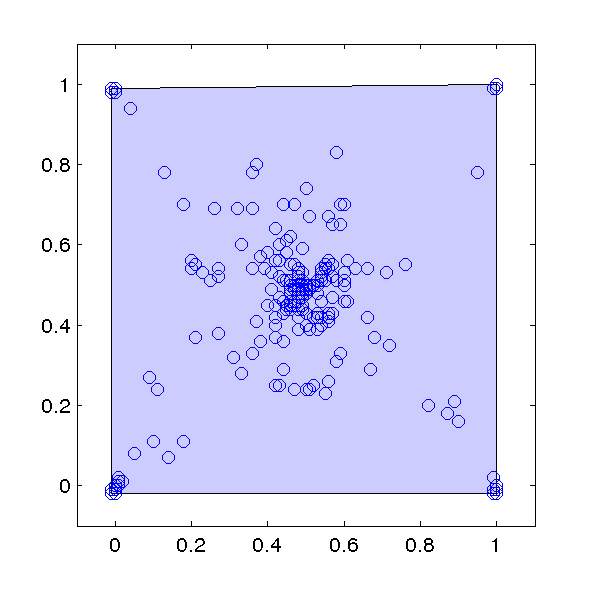}
\end{minipage}\hfill
\begin{minipage}[c]{.33\linewidth}
	\includegraphics[width=\linewidth,trim = .6cm .6cm .6cm .6cm,clip]{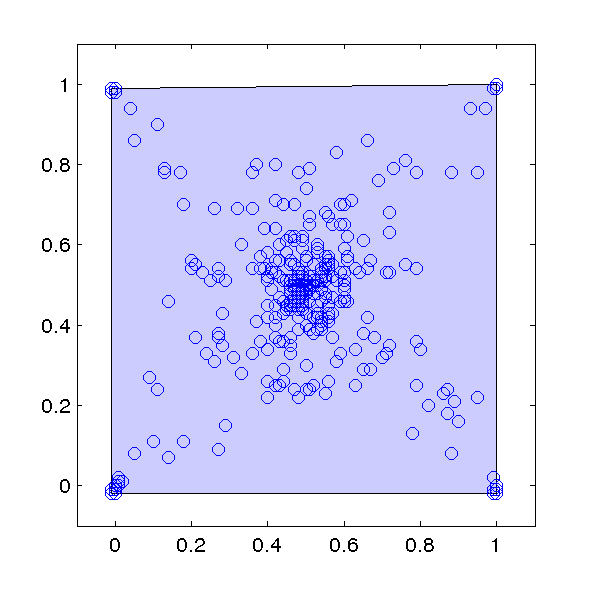}
\end{minipage}

\begin{minipage}[c]{.33\linewidth}
	\includegraphics[width=\linewidth,trim = .6cm .6cm .6cm .6cm,clip]{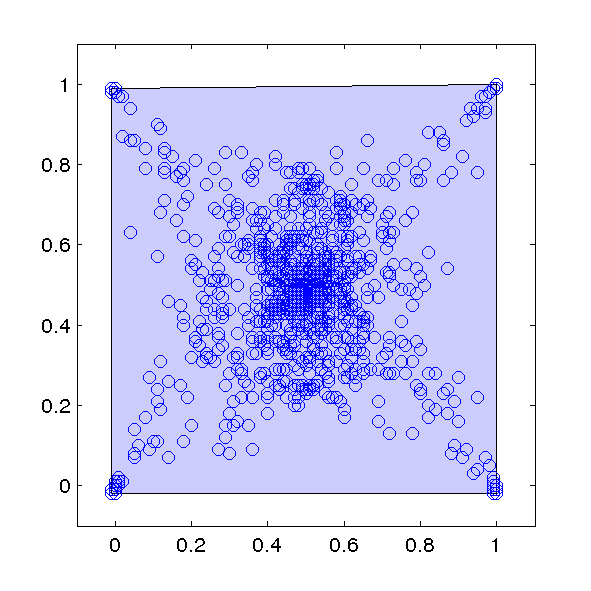}
\end{minipage}\hfill
\begin{minipage}[c]{.33\linewidth}
	\includegraphics[width=\linewidth,trim = .6cm .6cm .6cm .6cm,clip]{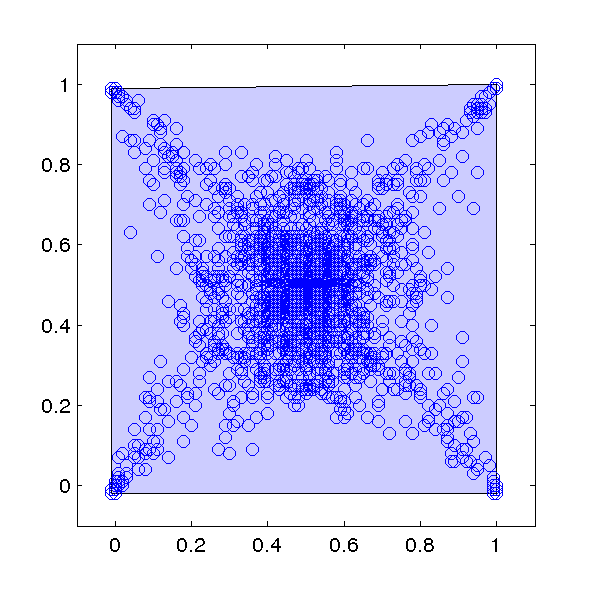}
\end{minipage}\hfill
\begin{minipage}[c]{.33\linewidth}
	\includegraphics[width=\linewidth,trim = .6cm .6cm .6cm .6cm,clip]{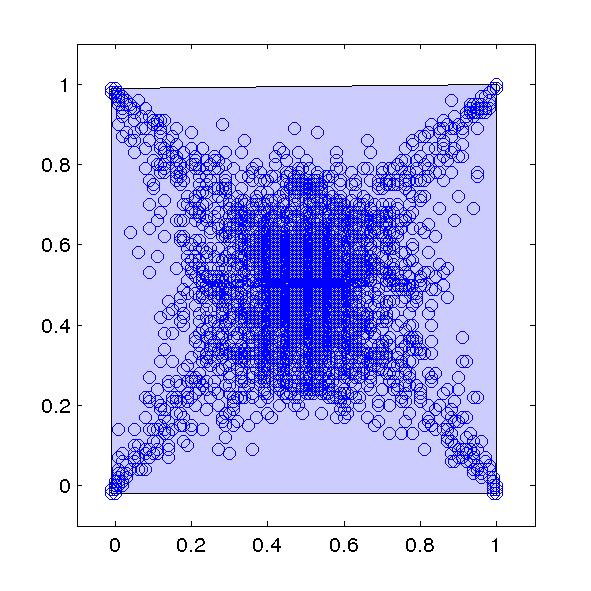}
\end{minipage}

\caption[Simulations of the asymptotic discretized rotation set of $f_1$]{Asymptotic discretized rotation set of $f_1$ as the union of the discretized rotation sets on grids $E_N$ with $100\le N\le M$ with $M = 100$ (top left), $M=150$ (top middle), $M = 200$ (top right), $M=500$ (bottom left), $M=1\,000$ (bottom middle) and $M=2\,000$ (bottom right)}\label{RotDissipSerie}
\end{figure}

For the discretized rotation set for $f_1$, the vertices of $[0,1]^2$ are also detected, and we only have a few points in the interior of the square (see Figure \ref{RotDissipDiscr}). However, when we compute the asymptotic discretized rotation set (see Figure \ref{RotDissipSerie}), we observe that the computed rotation vectors fill a great proportion of the square $[0,1]^2$, as predicted by the theory.
\bigskip

\begin{figure}[ht]
\begin{minipage}[c]{.33\linewidth}
	\includegraphics[width=\linewidth,trim = .6cm .6cm .6cm .6cm,clip]{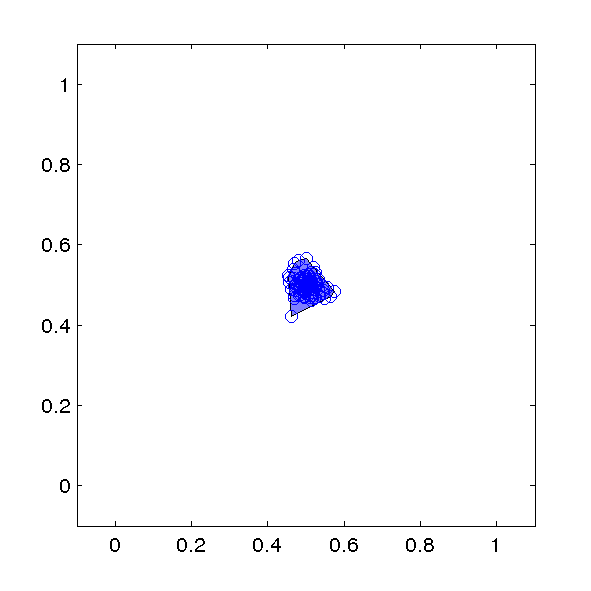}
\end{minipage}\hfill
\begin{minipage}[c]{.33\linewidth}
	\includegraphics[width=\linewidth,trim = .6cm .6cm .6cm .6cm,clip]{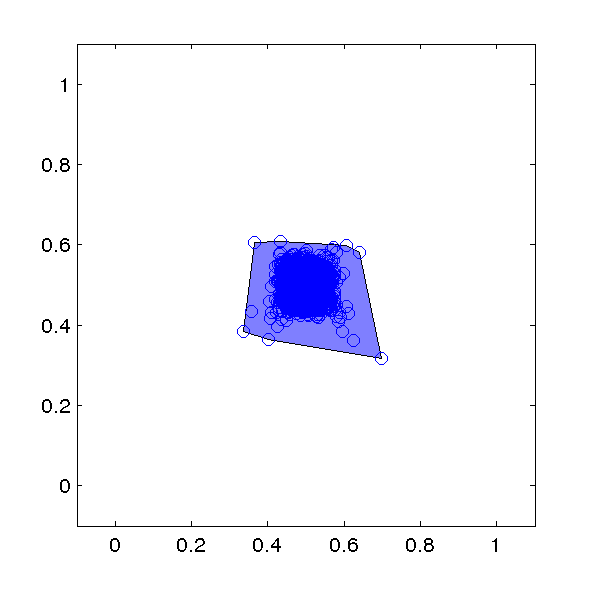}
\end{minipage}\hfill
\begin{minipage}[c]{.33\linewidth}
	\includegraphics[width=\linewidth,trim = .6cm .6cm .6cm .6cm,clip]{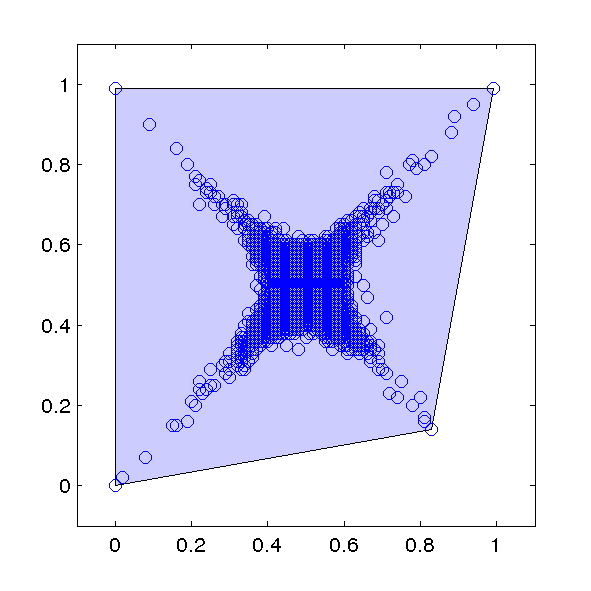}
\end{minipage}
\caption[Simulations of the observable rotation set of $g_1$]{Observable rotation set of $g_1$, $k$ orbits of length $1\,000$ with random starting points with $k=100$ (left), $10\,000$ (middle) and $1\,000\,000$ (right)}\label{RotConsC1Alea}
\end{figure}

In the conservative case, the rotation vectors of the observable rotation set are mainly quite close to the mean rotation vector of $g_1$, as predicted by Proposition \ref{ExGeneCons}. In particular in Figure \ref{RotConsC1Alea}, left, all the 100 rotation vectors of the computed observable rotation set are in the neighbourhood of $(1/2,1/2)$. Thus, the behaviour of these vectors is governed by Birkhoff's ergodic theorem with respect to the ergodic measure $\Leb$; \emph{a priori} this behaviour is quite chaotic and converges slowly: a typical orbit will visit every measurable subset with a frequency proportional to the measure of this set, so the rotation vectors will take time to converge. When the number of computed orbits increases (Figure \ref{RotConsC1Alea}, middle and rignt), we observe that a few rotation vectors are not close to the mean rotation vector; for $10^6$ different orbits we even detect three of the vertices of the actual rotation set. Anyway, even after 4 hours of calculation, we are 
unable to recover completely the initial rotation set of the homeomorphism.

\begin{figure}[ht]
\begin{minipage}[c]{.33\linewidth}
	\includegraphics[width=\linewidth,trim = .6cm .6cm .6cm .6cm,clip]{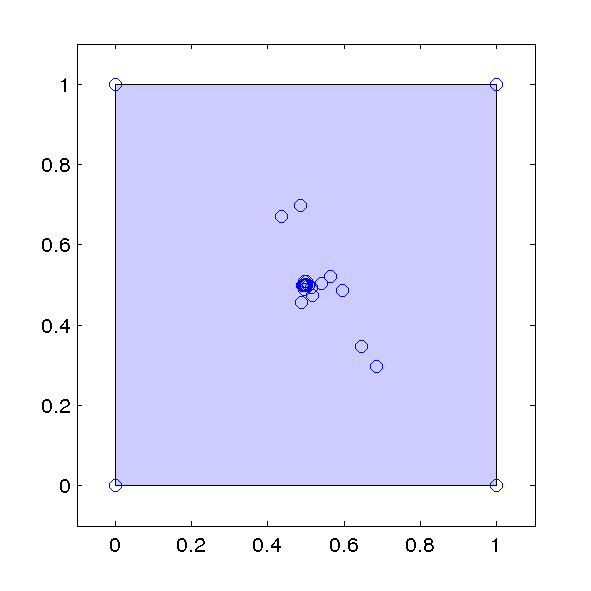}
\end{minipage}\hfill
\begin{minipage}[c]{.33\linewidth}
	\includegraphics[width=\linewidth,trim = .6cm .6cm .6cm .6cm,clip]{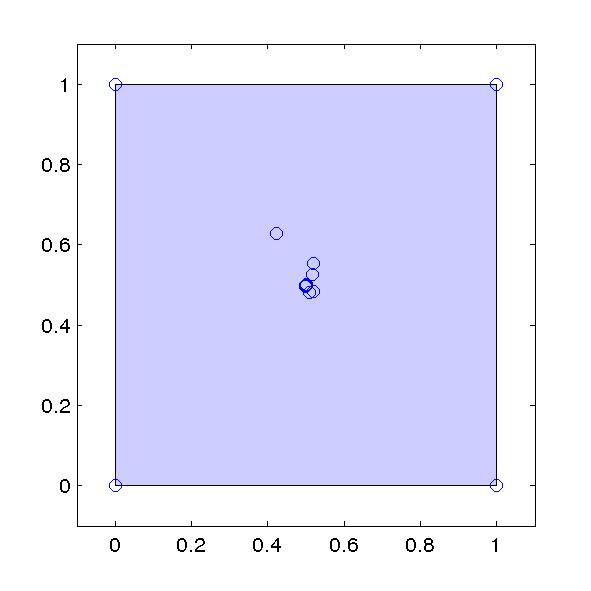}
\end{minipage}\hfill
\begin{minipage}[c]{.33\linewidth}
	\includegraphics[width=\linewidth,trim = .6cm .6cm .6cm .6cm,clip]{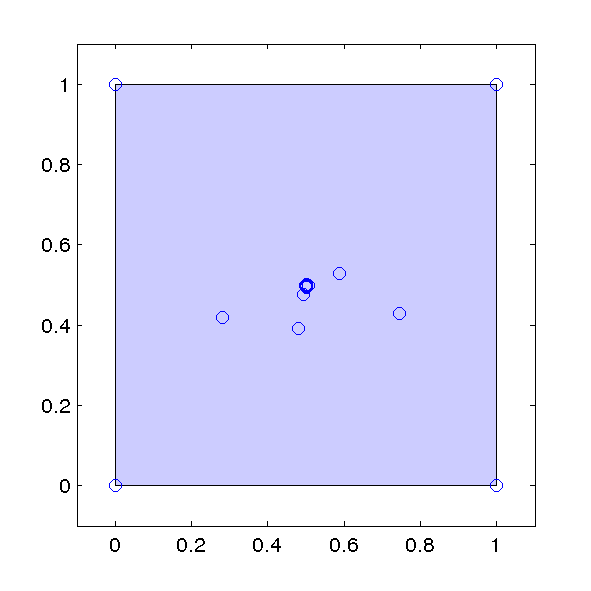}
\end{minipage}

\caption[Simulations of the discretized rotation set of $g_1$]{Discretized rotation set of $g_1$ on grids $E_N$, with $N=499$ (left), $N=500$ (middle) and $N=501$ (right)}\label{RotConsC1Discr}
\end{figure}

\begin{figure}[t]
\begin{minipage}[c]{.33\linewidth}
	\includegraphics[width=\linewidth,trim = .6cm .6cm .6cm .6cm,clip]{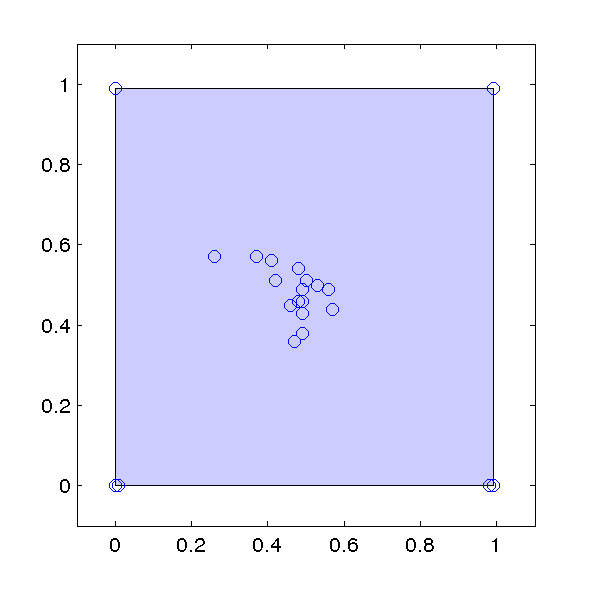}
\end{minipage}\hfill
\begin{minipage}[c]{.33\linewidth}
	\includegraphics[width=\linewidth,trim = .6cm .6cm .6cm .6cm,clip]{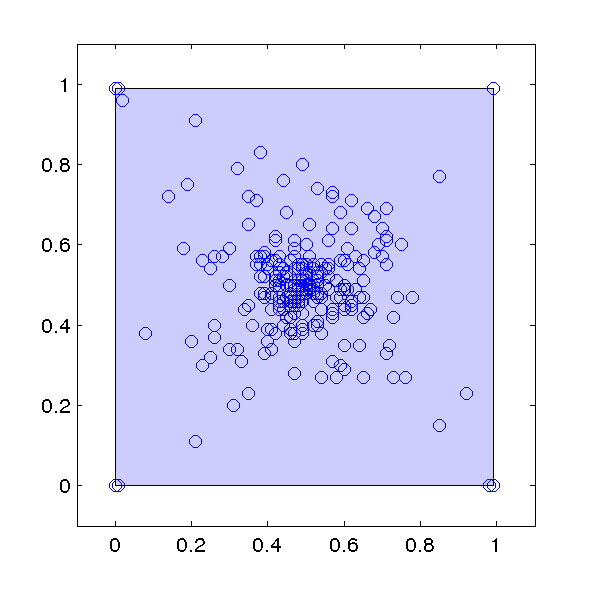}
\end{minipage}\hfill
\begin{minipage}[c]{.33\linewidth}
	\includegraphics[width=\linewidth,trim = .6cm .6cm .6cm .6cm,clip]{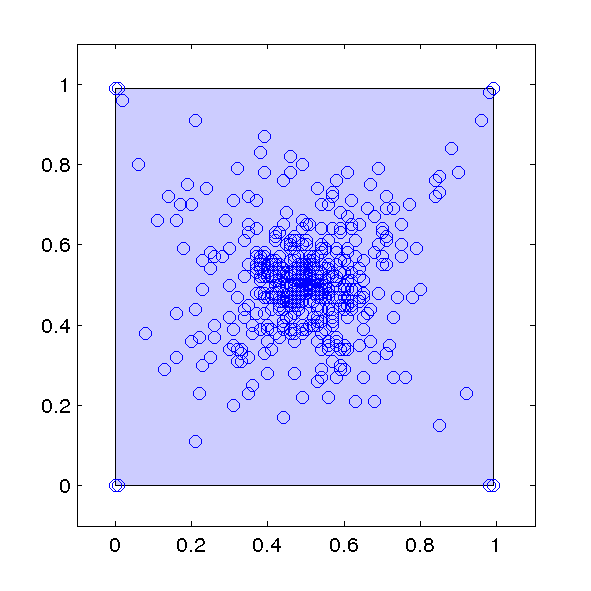}
\end{minipage}

\begin{minipage}[c]{.33\linewidth}
	\includegraphics[width=\linewidth,trim = .6cm .6cm .6cm .6cm,clip]{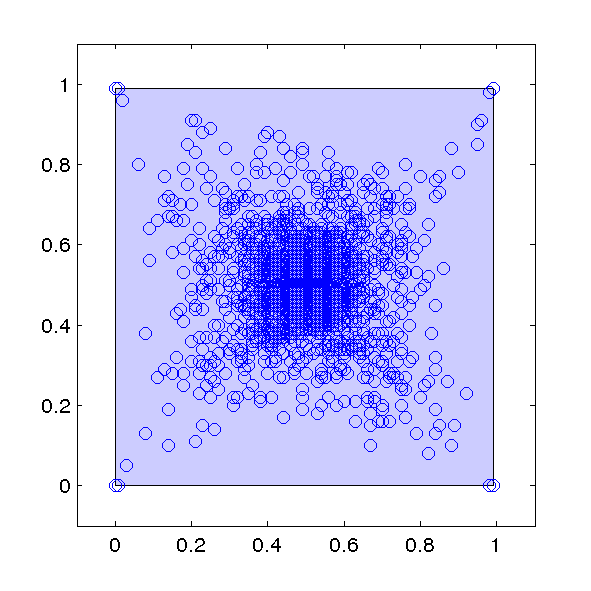}
\end{minipage}\hfill
\begin{minipage}[c]{.33\linewidth}
	\includegraphics[width=\linewidth,trim = .6cm .6cm .6cm .6cm,clip]{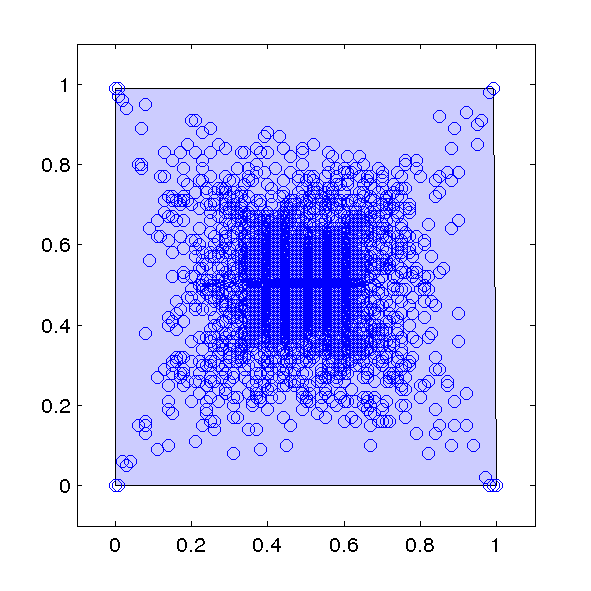}
\end{minipage}\hfill
\begin{minipage}[c]{.33\linewidth}
	\includegraphics[width=\linewidth,trim = .6cm .6cm .6cm .6cm,clip]{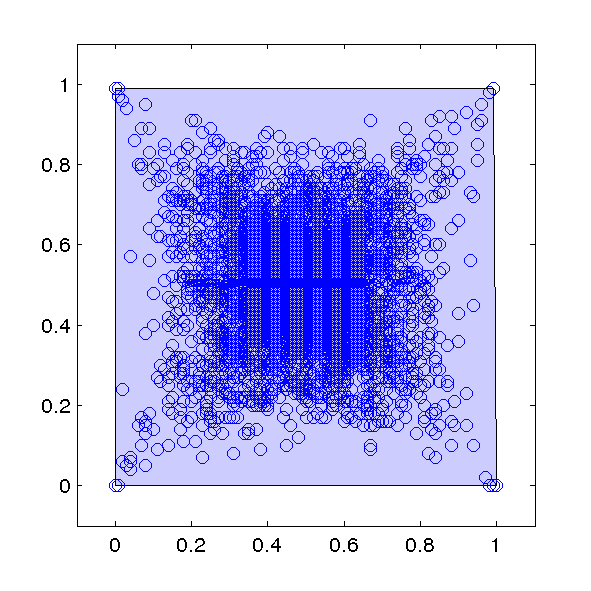}
\end{minipage}
\caption[Simulations of the asymptotic discretized rotation set of $g_1$]{Asymptotic discretized rotation set of $g_1$ as the union of the discretized rotation sets on grids $E_N$ with $100\le N\le M$ with $M = 100$ (top left), $M=150$ (top middle), $M = 200$ (top right), $M=500$ (bottom left), $M=1\,000$ (bottom middle) and $M=2\,000$ (bottom right)}\label{RotConsC1Serie}
\end{figure}

On the other hand, the convex hull of the discretized rotation set gives quickly a very good approximation of the rotation set. For example on a grid $500\times 500$ (Figure \ref{RotConsC1Discr}), with 2s of calculation (and even on a grid $100\times 100$ and $0.2s$ of calculation), we obtain a rotation set whose convex hull is already very close to $[0,1]^2$. However, for a single size of grid, we do not obtain exactly the conclusions of Corollary \ref{CoroRotDiscrCons} which states that for some integers $N$ the discretized rotation set should be close to the rotation set for Hausdorff distance; here for each $N$ we only have a few points in the interior of $[0,1]^2$. That is why we represented the union of the discretized rotation sets on grids $N\times N$ with $100\le N\le 1\,000$ (Figure \ref{RotConsC1Serie}). In this case we recover almost all the rotation set of $g_1$, except from the points which are close to one edge of the square but far from its vertices. The fact that we can obtain very easily 
the vertices of the rotation set can be due to the fact that in our example $f_1$ these vertices are realized by elliptic periodic points of the homeomorphism (in fact the derivative on this points is the identity). That is why we also conducted simulations of the homeomorphism $g_2$ which rotation set is also the square $[0,1]^2$ whose vertices are realized by non-elliptic periodic points.
\bigskip

\begin{figure}[ht]
\begin{minipage}[c]{.33\linewidth}
	\includegraphics[width=\linewidth,trim = .6cm .6cm .6cm .6cm,clip]{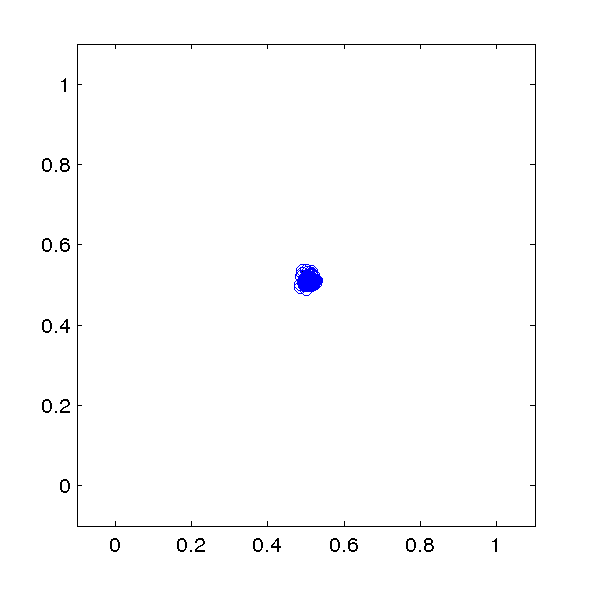}
\end{minipage}\hfill
\begin{minipage}[c]{.33\linewidth}
	\includegraphics[width=\linewidth,trim = .6cm .6cm .6cm .6cm,clip]{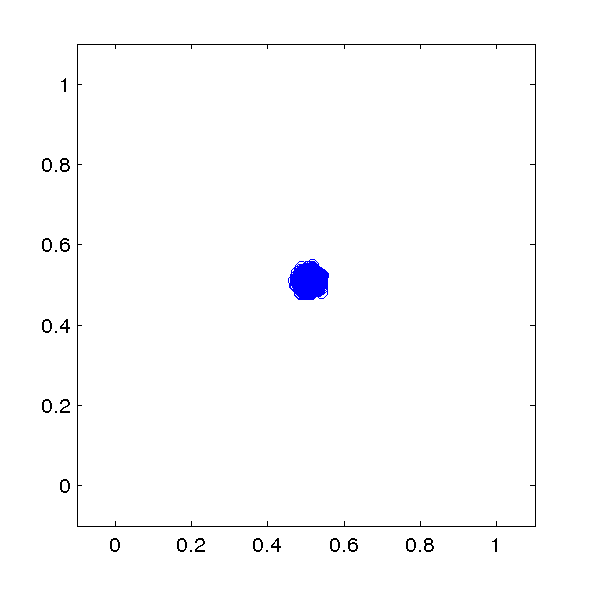}
\end{minipage}\hfill
\begin{minipage}[c]{.33\linewidth}
	\includegraphics[width=\linewidth,trim = .6cm .6cm .6cm .6cm,clip]{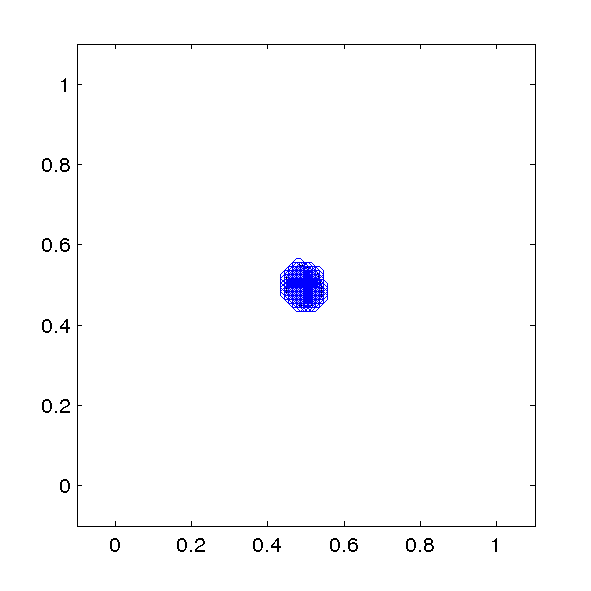}
\end{minipage}
\caption[Simulations of the observable rotation set of $g_2$]{Observable rotation set of $g_2$, $k$ orbits of length $1\,000$ with random starting points with $k=100$ (left), $10\,000$ (middle) and $1\,000\,000$ (right)}\label{RotConsC0Alea}
\end{figure}

In fact, when we compute the observable rotation set for $g_2$ (Figure \ref{RotConsC0Alea}), we only find rotation vectors which are close to the mean rotation vector, even after 4h of calculation. As the periodic points which realize the vertices of the rotation set are no longer elliptic, they are much more unstable and thus they are not detected by these simulations.

\begin{figure}[ht]
\begin{minipage}[c]{.33\linewidth}
	\includegraphics[width=\linewidth,trim = .6cm .6cm .6cm .6cm,clip]{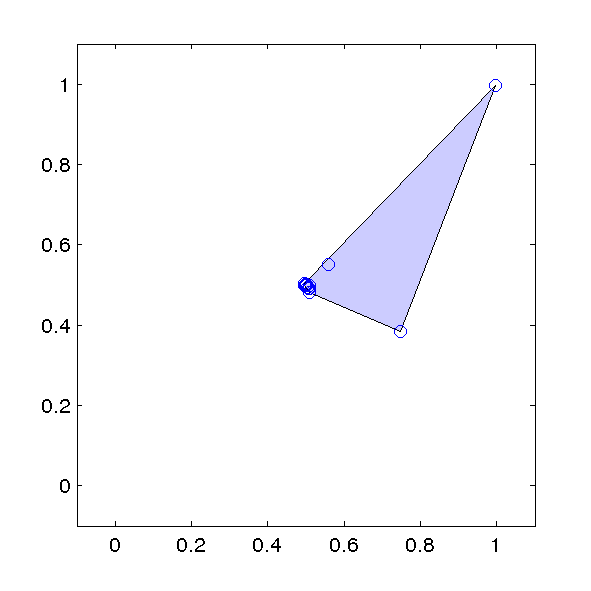}
\end{minipage}\hfill
\begin{minipage}[c]{.33\linewidth}
	\includegraphics[width=\linewidth,trim = .6cm .6cm .6cm .6cm,clip]{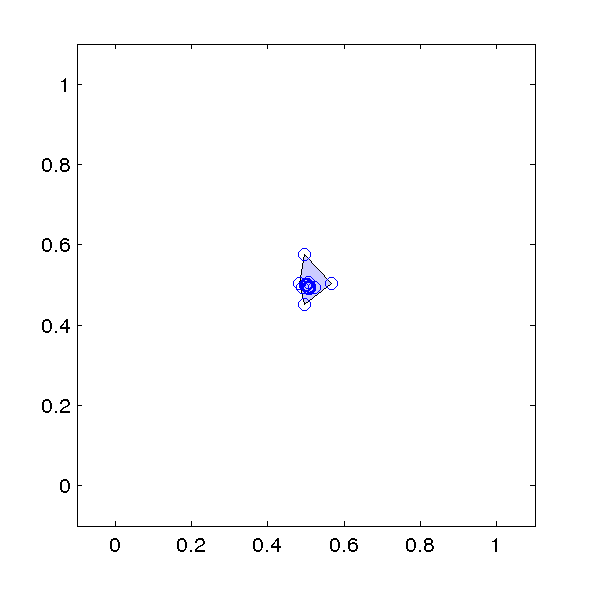}
\end{minipage}\hfill
\begin{minipage}[c]{.33\linewidth}
	\includegraphics[width=\linewidth,trim = .6cm .6cm .6cm .6cm,clip]{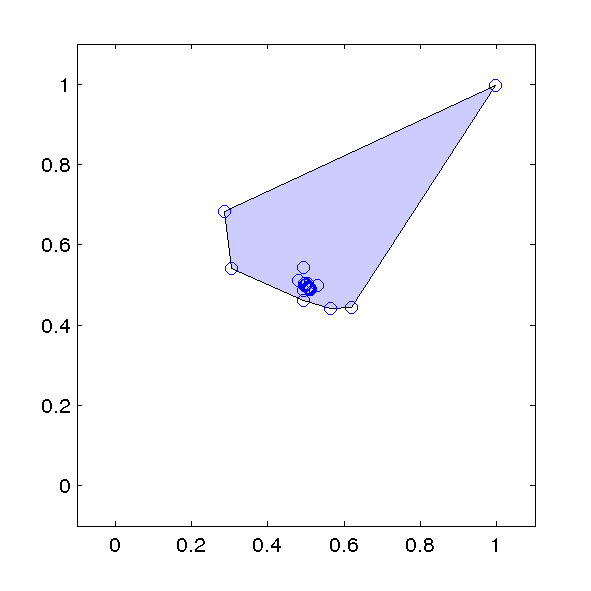}
\end{minipage}
\caption[Simulations of the discretized rotation set of $g_2$]{Discretized rotation set of $g_2$ on grids $E_N$, with $N=499$ (left), $N=500$ (middle) and $N=501$ (right)}\label{RotConsC0Discr}
\end{figure}

\begin{figure}[t]
\begin{minipage}[c]{.33\linewidth}
	\includegraphics[width=\linewidth,trim = .6cm .6cm .6cm .6cm,clip]{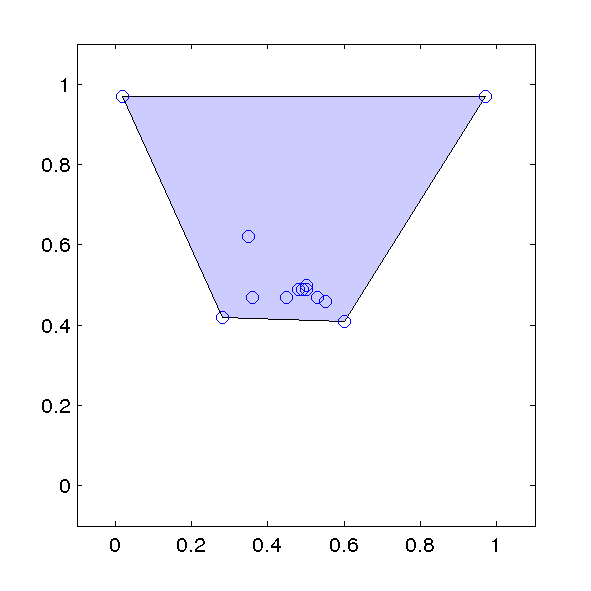}
\end{minipage}\hfill
\begin{minipage}[c]{.33\linewidth}
	\includegraphics[width=\linewidth,trim = .6cm .6cm .6cm .6cm,clip]{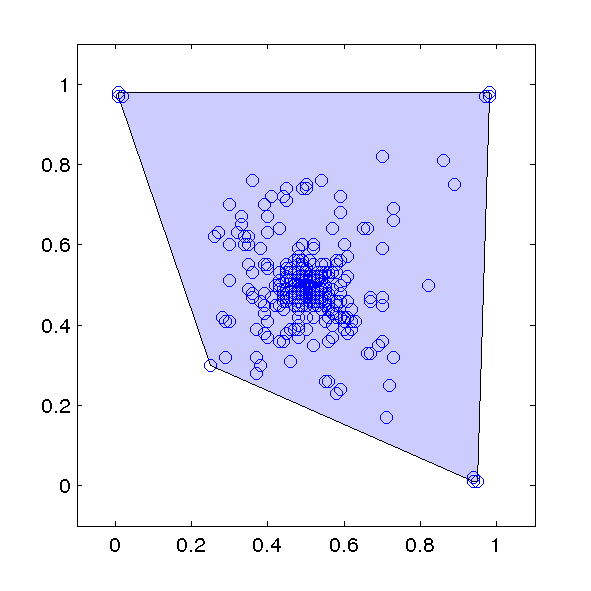}
\end{minipage}\hfill
\begin{minipage}[c]{.33\linewidth}
	\includegraphics[width=\linewidth,trim = .6cm .6cm .6cm .6cm,clip]{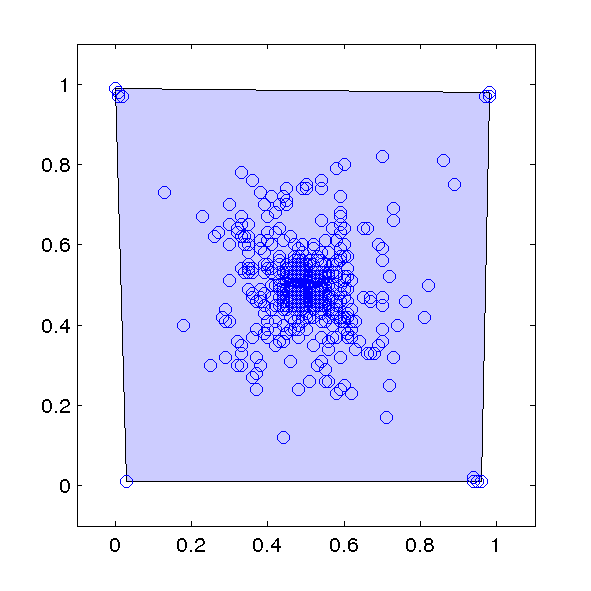}
\end{minipage}

\begin{minipage}[c]{.33\linewidth}
	\includegraphics[width=\linewidth,trim = .6cm .6cm .6cm .6cm,clip]{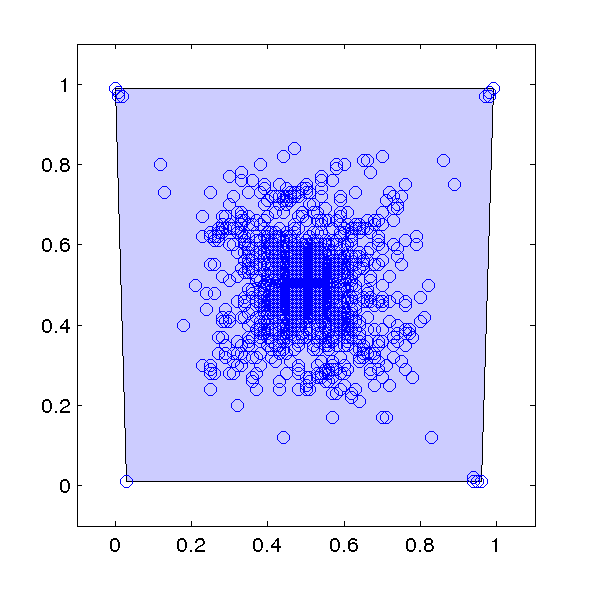}
\end{minipage}\hfill
\begin{minipage}[c]{.33\linewidth}
	\includegraphics[width=\linewidth,trim = .6cm .6cm .6cm .6cm,clip]{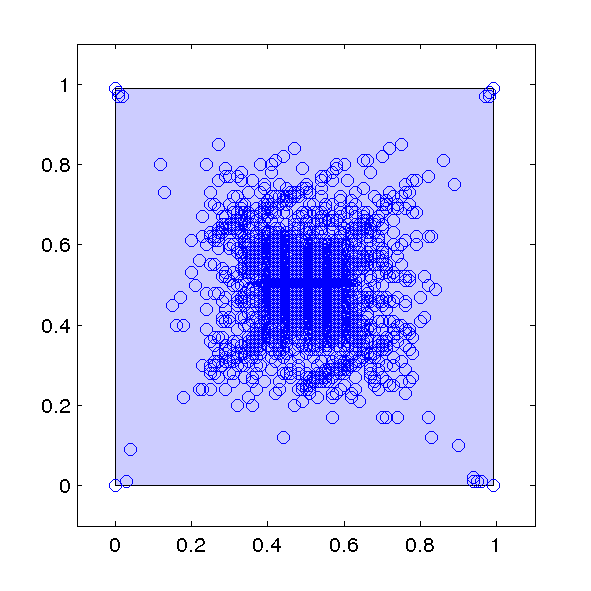}
\end{minipage}\hfill
\begin{minipage}[c]{.33\linewidth}
	\includegraphics[width=\linewidth,trim = .6cm .6cm .6cm .6cm,clip]{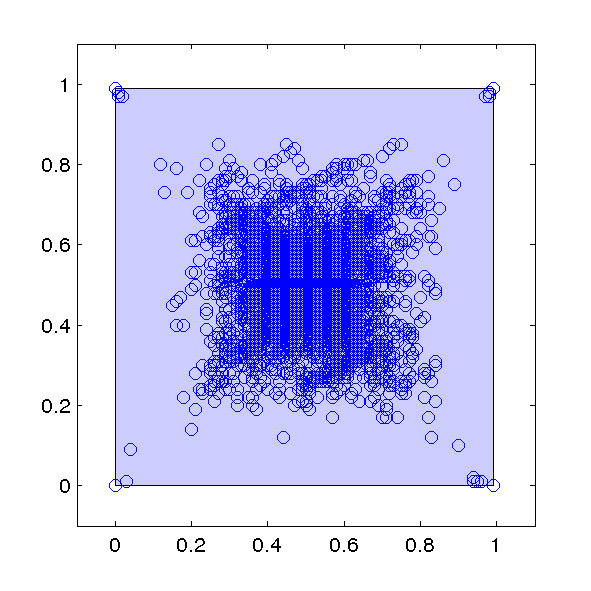}
\end{minipage}
\caption[Simulations of the asymptotic discretized rotation set of $g_2$]{Asymptotic discretized rotation set of $g_2$ as the union of the discretized rotation sets on grids $E_N$ with $100\le N\le M$ with $M = 100$ (top left), $M=150$ (top middle), $M = 200$ (top right), $M=500$ (bottom left), $M=1\,000$ (bottom middle) and $M=2\,000$ (bottom right)}\label{RotConsC0Serie}
\end{figure}

The sets detected by the discretized rotation sets of order $499$, $500$ and $501$ (Figure \ref{RotConsC0Discr}) are quite bigger than those detected by the simulations of the observable rotation set, even if the time of calculation is much smaller. However, we do not recover the whole rotation set of the homeomorphism (we conducted simulations for higher orders around $N=1\,000$ and $N=2\,000$ and the behaviour is similar). By contrast, the simulations of the asymptotic discretized rotation set (Figure \ref{RotConsC0Serie}) allows us to see the actual rotation set of the homeomorphism: when we represent all the rotation vectors of the discretizations of order $100\le N\le M$ with $M=200$ (which takes about 45s of calculation) we obtain a set which is very close to the square $[0,1]^2$; for $M=100$ ($\simeq$ 1h 45min of calculation) we recover almost exactly the initial rotation set.
\bigskip

\begin{figure}[ht]
\begin{minipage}[c]{.33\linewidth}
	\includegraphics[width=\linewidth,trim = .6cm .6cm .6cm .6cm,clip]{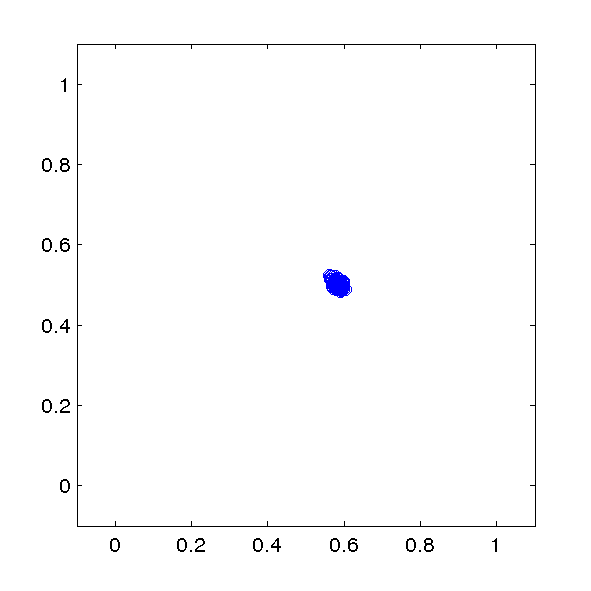}
\end{minipage}\hfill
\begin{minipage}[c]{.33\linewidth}
	\includegraphics[width=\linewidth,trim = .6cm .6cm .6cm .6cm,clip]{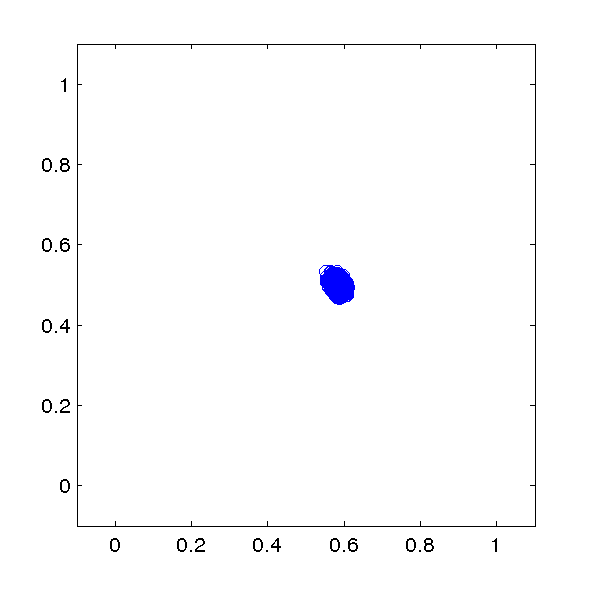}
\end{minipage}\hfill
\begin{minipage}[c]{.33\linewidth}
	\includegraphics[width=\linewidth,trim = .6cm .6cm .6cm .6cm,clip]{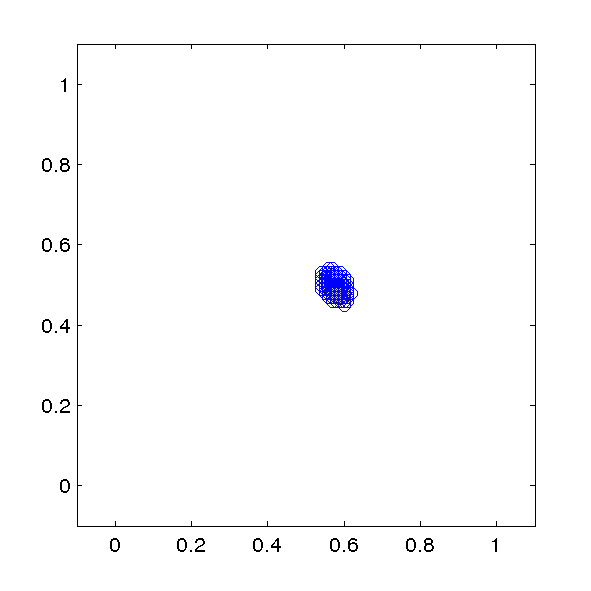}
\end{minipage}
\caption[Simulations of the observable rotation set of $g_3$]{Observable rotation set of $g_3$, $k$ orbits of length $1\,000$ with random starting points with $k=100$ (left), $10\,000$ (middle) and $1\,000\,000$ (right)}\label{RotLoinAlea}
\end{figure}

Finally, the behaviour of the observable rotation set of $g_3$ is very similar to that of $g_2$ (see Figure \ref{RotLoinAlea}): even when we compute $1\,000\,000$ different orbits with random starting points, we only obtain rotation vectors which are close to $(0.55,0.5)$, which should be a good approximation of the mean rotation vector.

\begin{figure}[ht]
\begin{minipage}[c]{.33\linewidth}
	\includegraphics[width=\linewidth,trim = .6cm .6cm .6cm .6cm,clip]{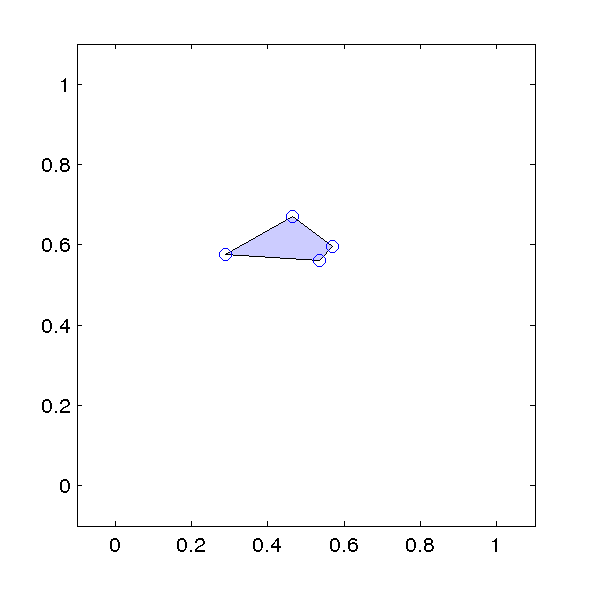}
\end{minipage}\hfill
\begin{minipage}[c]{.33\linewidth}
	\includegraphics[width=\linewidth,trim = .6cm .6cm .6cm .6cm,clip]{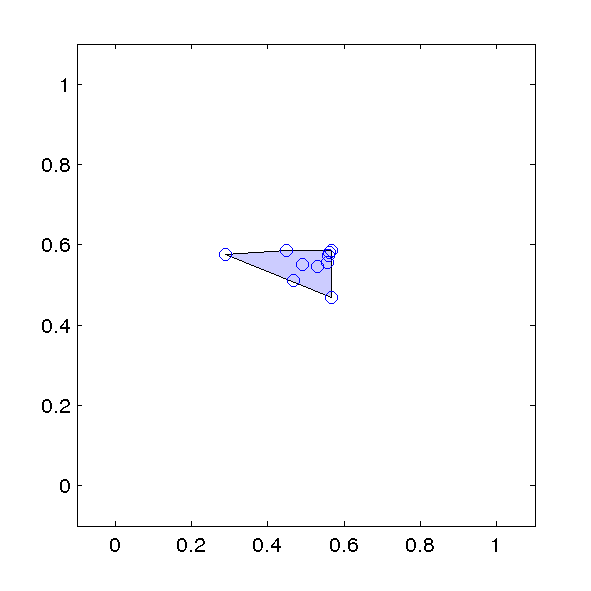}
\end{minipage}\hfill
\begin{minipage}[c]{.33\linewidth}
	\includegraphics[width=\linewidth,trim = .6cm .6cm .6cm .6cm,clip]{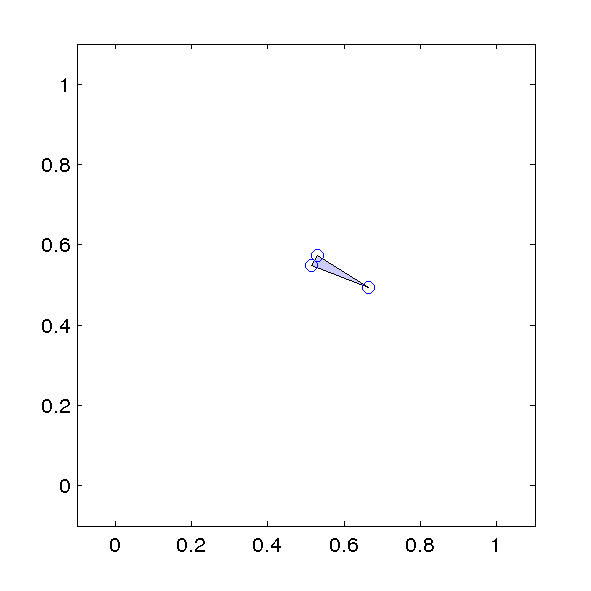}
\end{minipage}
\caption[Simulations of the discretized rotation set of $g_3$]{Discretized rotation set of $g_3$ on grids $E_N$, with $N=499$ (left), $N=500$ (middle) and $N=501$ (right)}\label{RotLoinDiscr}
\end{figure}

Like for $g_2$, the simulations of the discretized rotation sets for the grids $E_N$ with $N\in\{499,500,501\}$ (Figure \ref{RotLoinDiscr}) are not very convincing: the sets do not seem to converge to anything. We have to compute the asymptotic discretized rotation sets (Figure \ref{RotLoinDiscr}) to see something that looks like a convergence  for the Hausdorff topology of the computed rotation sets. However, this convergence in practical is just an indication that the set we compute is close to the actual rotation set of $g_3$. To our knowledge, it is impossible to ensure that for a given order of discretization, the asymptotic rotation set computed to this order is close to the rotation set of $g_3$.


\part{Discretizations of linear maps}\label{PartII}

\parttoc

\chapter*[Introduction]{Introduction}

In the second part of this manuscript, we will consider the dynamical behaviour of the discretizations of \emph{linear maps}. We were led to this question by the study of the dynamics of the discretizations of generic \emph{$C^1$-diffeomorphisms}; for exmaple, we will see in the next part that the study of the degree of recurrence of a generic conservative $C^1$-diffeomorphism requires a good understanding of the dynamics of the discretizations of linear maps. However, it seemed to us that the study of the linear case could form a whole part. Indeed, it appeared that this subject is itself quite interesting and does not use the tools developed in the other parts of this manuscript. Moreover, the techniques involved in this second part of the thesis are very different from that used in the others: we will use tools like almost periodic patterns, model sets or lattice tilings of the Euclidean space by unit cubes.
\bigskip

The definition of the discretization of a linear map of $\R^n$ ($n\ge 2$) is made to mimic that of the discretization of a map of a compact space. The idea underlying this definition is to take the grid $\frac{1}{N}\Z^n \cap [-1,1]^n$ for a large $N$, and to consider the images of this grid by the discretization of linear maps on the grid $\frac{1}{N}\Z^n$. As the action of a linear map is invariant under homothety, it is equivalent to look at the grids $\Z^n \cap [-N,N]^n$ and to consider the images of these grids by the discretization of linear maps on the grid $\Z^n$. Thus, we define a  projection $\pi$ of $\R^n$ onto $\Z^n$, which maps any point of $\R^n$ onto (one of) the closest point of the lattice $\Z^n$ (see Definition~\ref{DefDiscrLin}); given $A\in GL_n(\R)$, the \emph{discretization} of $A$ is the map
\[\widehat A = \pi\circ A : \Z^n\to\Z^n.\]
Given a sequence $(A_k)_{k\ge 1}$ of matrices of $GL_n(\R)$, we want to study the dynamics of the sequence $(\widehat{A_k})_{k\ge 1}$, and in particular, the \emph{density} of the set $\Gamma_k = (\widehat{A_k}\circ\cdots\circ\widehat{A_1})(\Z^n)$. This density is defined in the following way, where $B_R$ denotes the infinite ball $B_\infty(0,R)$ :
\begin{equation}\label{premEq}
D^+(\Gamma_k) =  \underset{R\to +\infty}{\overline\lim} \frac{\card (\Gamma_k \cap B_R)}{\card (\Z^n \cap B_R)}.
\end{equation}

First of all, we study the structure of the image sets $\Gamma_k$. It appears that there is a kind of ``regularity at infinity'' of the behaviour of $\Gamma_k$. More precisely, this set is an \emph{almost periodic pattern}: for every $\varep>0$ there exists a relatively dense\footnote{A set $\mathcal N$ is \emph{relatively dense} if there exists $R>0$ such that every ball of radius $R$ contains at least one point of $\mathcal N$.} set of $\varep$-translations of $\Gamma_k$, where a vector $v\in\Z^n$ is an $\varep$-translation of $\Gamma_k$ if $D^+\big((\Gamma_k-v)\Delta \Gamma_k\big)\le\varep$ (see Definition~\ref{DefAlmPer}). Roughly speaking, for $R$ large enough, the set $\Gamma_k \cap B_R$ determines the whole set $\Gamma_k$ up to an error of density smaller than $\varep$. We prove that the image of an almost periodic pattern by the discretization of a linear map is still an almost periodic pattern (Theorem~\ref{imgquasi}); thus, given a sequence $(A_k)_{k\ge 1}$ of invertible matrices, the sets $(\widehat{A_k}\circ\cdots\circ\widehat{A_1}) (\Z^n)$ are almost periodic patterns. In particular, these sets possess a uniform density: the superior limit in Equation~\eqref{premEq} is in fact a limit (Corollary~\ref{corolimitexist}). This allows us to define the rate of injectivity: given a sequence $(A_k)_{k\ge 1}$ of linear maps, the \emph{rate of injectivity in time $k$} of this sequence is the quantity (see Definition~\ref{DefTaux})
\[\tau^k(A_1,\cdots,A_k) = \lim_{R\to +\infty} \frac{\card \big((\widehat{A_k}\circ\cdots\circ\widehat{A_1}) (B_R\cap\Z^n)\big)}{\card (B_R\cap\Z^n)}\in]0,1],\]
where $B_R = B_\infty(0,R)$. As these quantities are decreasing in $k$, we can also define the \emph{asymptotic rate of injectivity}
\[\tau^\infty\big((A_k)_{k\ge 1}\big) = \lim_{k\to +\infty}\tau^k(A_1,\cdots,A_k)\in[0,1].\]
These rates on injectivity can be seen as the quantity of information we lose when we apply a sequence of discretizations of matrices.

The link between the rate of injectivity and the density of the images sets is made by the following formula (Proposition~\ref{TauDens}):
\[\tau^k(A_1,\cdots,A_k) = \det(A_1)\cdots \det(A_k) D^+\big((\widehat{A_k}\circ\cdots\circ\widehat{A_1})(\Z^n)\big).\]
\bigskip

The goal of the second chapter of this part is to study the behaviour of the asymptotic rate of injectivity of a generic sequence of matrices of $SL_n(\R)$ (in fact, the same holds for matrices with determinant $\pm 1$). It is given by the main theorem of this chapter (Theorem~\ref{ConjPrincip}).

\begin{theorem}\label{ConjIntro}
For a generic sequence $(A_k)_{k\ge 1}$ of matrices of $SL_n(\R)$, we have $\tau^\infty\big((A_k)_k\big) = 0$.
\end{theorem}

If this result can seem quite natural, its proof is far from being trivial. First of all, the sets $\Gamma_k = (\widehat{A_k}\circ\cdots\widehat{A_1})(\Z^n)$ are ``more and more complex'' when $k$ increases: \emph{a priori}, the radius $R_0$ for which $\Gamma_k\cap B_{R_0}$ determines almost all $\Gamma_k$ is more and more large when $k$ increases; it is very difficult to have an idea of the local geometry of these sets, thus to decide which is the best possible perturbation of the matrices (especially since we have to get practical estimates on the loss of injectivity). Moreover, once the set $\Gamma_k$ has a density smaller than $1/2$, it may be impossible to make the rate decrease in one step of time. For example, if we set
\[A_1 = \begin{pmatrix} 4 & \\ & 1/4 \end{pmatrix}\quad \text{and}\quad A_2 = \begin{pmatrix} 1/2 & \\ & 2 \end{pmatrix},\]
then $\big(\widehat{A_2}\circ \widehat{A_1}\big)(\Z^2) = (2\Z)^2$, and for every $B_3$ close to the identity, we have $\tau^3(A_1,A_2,B_3) = \tau^2(A_1,A_2) = 1/4$. In such cases, it seems difficult to have a long term strategy to make the rate decrease\dots

The idea of the proof of Theorem~\ref{ConjIntro} is to take advantage of the fact that for a generic sequence of matrices, the coefficients of all the matrices are rationally independent (in fact, we will only need a weaker assumption about he independence). For example, for a generic matrix $A\in SL_n(\R)$, the set $A(\Z^n)$ is uniformly distributed modulo $\Z^n$. We then remark that the local pattern of the image set $\widehat A(\Z^n)$ around $\widehat A(x)$ is only determined by $A$ and the the remainder of $Ax$ modulo $\Z^n$: the global behaviour of $\widehat A(Z^n)$ is coded by the quotient $\R^n/\Z^n$. This somehow reduces the study to a local problem.

As a first application of this remark, we state that the rate of injectivity in time 1 can be seen as the area of an intersection of cubes (Proposition~\ref{FormTau1}). This observation, combined with considerations about the frequency of differences $\rho_{\Gamma_k}(v) = D\big(({\Gamma_k}-v)\cap{\Gamma_k}\big)$, allows to prove a weak version of Theorem~\ref{ConjIntro}: the asymptotic rate of injectivity of a generic sequence of matrices of $SL_n(\R)$ is smaller than $1/2$ (Theorem~\ref{PerLin1}). It also leads to a proof of the following result.

\begin{theorem}\label{DD}
Let $(P_k)_{k\ge 1}$ be a generic sequence of matrices of $O_n(\R)$. Then $\tau^\infty((P_k)_k) = 0$.
\end{theorem}

\begin{figure}[t]
\begin{center}
\begin{minipage}[c]{.35\linewidth}
	\includegraphics[width=\linewidth]{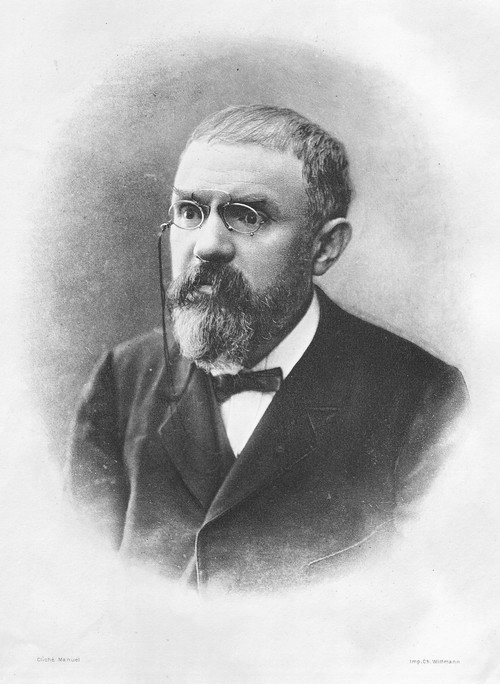}
\end{minipage}\hspace{50pt}
\begin{minipage}[c]{.35\linewidth}
	\includegraphics[width=\linewidth]{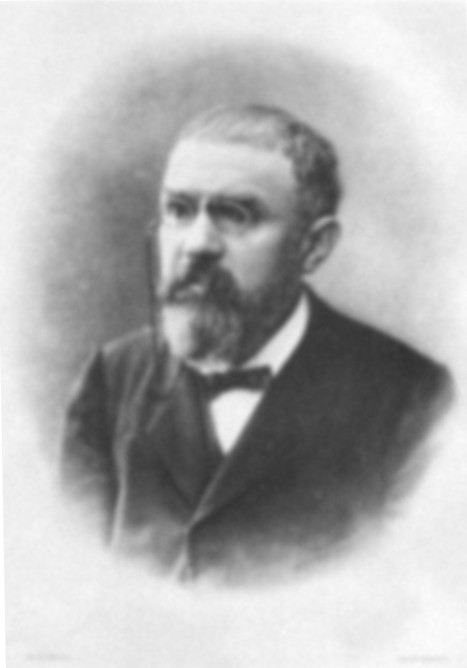}
\end{minipage}
	\caption[40 successive rotations of an image]{Original image (left) and 40 successive random rotations of this image, obtained with a consumer software.}\label{PoincareRot}
\end{center}	
\end{figure}

In particular, for a generic sequence of rotations of the plane, we have a total loss of information when we apply successively the discretizations of these rotations (see Figure~\ref{PoincareRot}).

At the end of the second chapter of this part, we make a full use of the equidistribution property to see the rate of a sequence $A_1,\cdots,A_k$ of matrices of $SL_n(\R)$ in terms of areas of intersections of cubes in $\R^{nk}$ (Proposition~\ref{CalculTauxModel}). The proof of this formula is based on the notion of model set\footnote{Sometimes called ``cut-and-project'' set.} (Definition~\ref{DefModel}), which is a particular class of almost periodic patterns. Using this formula, we replace the iteration by a passage in high dimension. These considerations allow to prove Theorem~\ref{ConjIntro}, without having to make ``clever'' perturbations of the sequence of matrices (that is, the perturbations made a each iteration are chosen independently from that made in the past or in the future).

Finally, the last chapter of this part is devoted to the study of the statistics of the roundoff errors induced by the discretization process. The main result is the following (Proposition~\ref{EquidistribErr}).

\begin{propo}
For a generic sequence $(A_k)_{k\ge 1}$ of matrices of $GL_n(\R)$, or $SL_n(\R)$, or $O_n(\R)$, for every fixed integer $k$, the finite sequence of errors $\varep_x= \big(\varep_x^1,\cdots,\varep_x^k\big)$ is equidistributed in $(\R^n/\Z^n)^k$ when $x$ ranges over $\Z^n$.
\end{propo}

This chapter presents a work in progress; the goal of these considerations is to study in more depth a conjecture of O.E. Lanford (Conjecture~\ref{Lalan}) concerning the physical measures of expanding maps of the circle.
\bigskip

\begin{figure}[t]
\begin{center}
\begin{minipage}[c]{.4\linewidth}
	\includegraphics[width=\linewidth]{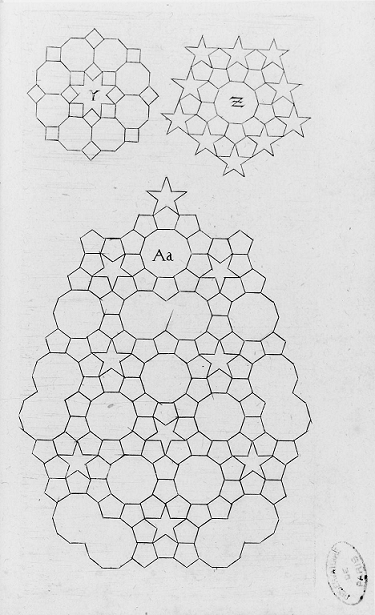}
	\caption[Some drawings of J. Kepler]{Some drawings of J.~Kepler (1619).}\label{Kepler}
\end{minipage}\hfill
\begin{minipage}[c]{.55\linewidth}
	\includegraphics[width=\linewidth]{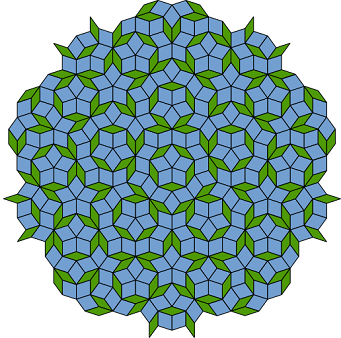}
	\caption[Penrose tiling]{A Penrose tiling.}\label{Penrose}
\end{minipage}
\end{center}	
\end{figure}

\begin{figure}[t]
\begin{center}
\begin{minipage}[c]{.33\linewidth}
	\includegraphics[width=\linewidth]{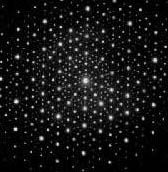}
	\caption[Figure of diffraction of a quasicrystal]{The figure of diffraction behind the discovery of quasicrystals; this figure possesses a symmetry of order 5, which is impossible for regular crystals.}\label{Diffraction}
\end{minipage}\hfill
\begin{minipage}[c]{.6\linewidth}
	\includegraphics[width=\linewidth]{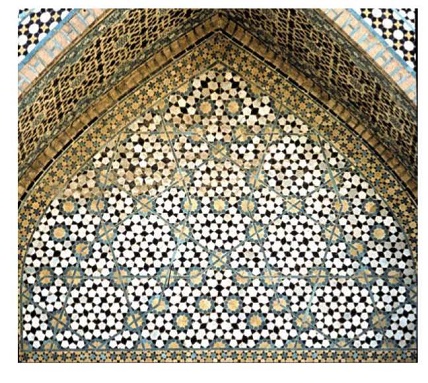}
	\caption[Islamic art]{Mosaic in a madrassa of Boukhara (Ouzbekistan, 15th century). In 2007, P. Lu discovered that it has the same structure as a Penrose tiling.}\label{Islam}
\end{minipage}
\end{center}	
\end{figure}
\bigskip

Throughout this part, we will use several notions of almost periodicity: Delone sets, almost periodic patterns, weakly almost periodic sets (which is, to our knowledge, an original concept), model sets, Bohr and Besicovitch almost periodic functions. We now give a quick overview of the history of these different concepts. The history of almost periodic sets is made of many independent rediscoveries, made in different areas of mathematics and physics. The first traces of studies of some notions of almost periodicity date back to 1893 with the master's thesis of P. Bohl ``\emph{Über die Darstellung von Funktionen einer Variabeln durch trigonometrische Reihen mit mehreren einer Variabeln proportionalen Argumenten}'', where the author introduces quasi-periodic functions. This class of functions were rediscovered about ten years later by E. Esclangon, in the view of the study of celestial mechanic. Almost periodic functions have been studied in more depth by H. Bohr in 1924. In the paper \cite{MR1555192}, he proves in particular the fundamental theorem stating that any uniform almost periodic function (now called Bohr almost periodic function) is the limit of generalized trigonometric series (and reciprocally). Two years later, A. Besicovitch defined a weaker notion of almost periodicity \cite{MR1575297}; the results of these papers were summarized later in the books \cite{MR0020163} and \cite{MR0068029}. In the early seventies, Y.~Meyer \cite{MR0485769} introduced \emph{model sets} in the context of harmonic analysis (more precisely, he wanted to study Pisot and Salem numbers). These sets were rediscovered independently in 1981 by N. G. de Bruijn (see \cite{MR609465}) to study Penrose aperiodic tilings. These aperiodic tilings we introduced in 1974 by the physicist R. Penrose \cite{Penrose} (see Figure~\ref{Penrose}), inspired by the work of J. Kepler (see Figure~\ref{Kepler}). In 1984 was made the fundamental discovery of quasicrystals \cite{1984PhRvL..53.1951S} (see Figure~\ref{Diffraction}), that is, solids that are ordered but not periodic. Quickly, the link was made with the previous works of H. Bohr and A. Besicovitch. Since then, the various notions of almost periodicity play an important role in the study of quasicrystals and aperiodic tilings, but also in many other parts of mathematics (see for example the survey \cite{Moody25}). Various mathematical formulations of quasicrystals have been proposed, such as Meyer sets (see \cite{MR1420415}), and the harmonic study of these sets has been investigated (see \cite{MR2876415}). Quite recently, P.~Lu, an American physicist, discovered that some very old Uzbek mosaics have the same structure than the Penrose tiling (see Figure~\ref{Islam}). It is amusing to note that in our case, the study of almost periodic sets arises from a branch of mathematics which is still different: the discretizations of generic diffeomorphisms of the torus.
\bigskip

\label{BibliLin}The particular problem of the discretization of linear maps has been only little studied. To our knowledge, what has been made in this direction has been initiated by image processing. The goal of these studies is to try to answer to the question: what is the best way to define the action of a linear map on the lattice $\Z^2$? In particular, how can we compute the image of a numerical image by a linear map? More precisely, we want to avoid phenomenons like loss of information (due to the fact that discretizations of linear maps are not injective) or aliasing (the apparition of undesirable periodic patterns in the image, due for example to a resonance between a periodic pattern in the image and the discretized map). For example, in Figure~\ref{PoincareRot}, we have applied 40 successive random rotations to a $500\times 684$ pixels picture, using a consumer software. These discretized rotations induce a very strong blur in the resulting image, thus a big loss of information. To our knowledge, the existing studies are mostly interested in the linear maps with \emph{rational coefficients} (see for example \cite{Jacob25}, \cite{MR1382839} or \cite{MR1832794}), and also in the specific case of \emph{rotations} (see for example \cite{A1996_1075}, \cite{nouvel:tel-00444088}, \cite{thibault:tel-00596947}, \cite{MR1782038}). These works mainly focus on the local behaviour of the images of $\Z^2$ by discretizations of linear maps: given a radius $R$, what pattern can follow the intersection of this set with any ball of radius $R$? What is the number of such patterns, what are their frequencies? Are they complex (in a sense to define) or not? Are these maps bijections? In particular, the thesis \cite{nouvel:tel-00444088} of B.~Nouvel gives a characterization of the angles for which the discrete rotation is a bijection (such angles are countable and accumulate only on $0$).

Here, our point of view is quite different: we want to determine the dynamical behaviour of discretizations of \emph{generic} linear maps; in particular generic matrices are totally irrational (by this, we mean that the image of the lattice $\Z^n$ by the matrix is equidistributed modulo $\Z^n$). We will see that the behaviour of discretizations of generic linear maps is in a certain sense smoother than that of rational linear maps: for example the rate of injectivity is continuous when restricted to totally irrational matrices (Proposition~\ref{ThMeanRate}), while it is not on some rational matrices (Proposition~\ref{rateEx}). However, we study the rate of injectivity of a sequence of matrices, which measures the loss of information we have when we apply several discretized linear maps to an image; Theorem~\ref{ConjIntro} expresses that in the general case, when we apply a lot of linear maps, then we lose most of the information contained in the image, and Theorem~\ref{DD} states that this phenomenon also appears for a generic sequence of rotations.

\chapter{Almost periodic sets}\label{ChapAlm}

In this chapter, we introduce the basic notions that we will use during the study of discretizations of linear maps of $\R^n$, $n\ge 1$. First of all, we introduce the notion of almost periodic pattern: roughly speaking, a set $\Gamma\subset \Z^n$ is an almost periodic pattern if for $R$ large enough, the set $\Gamma \cap [-R,R]^n$ determines the whole set $\Gamma$ up to an error of density smaller than $\varep$ (see Definition~\ref{DefAlmPer}). It can be easily seen that an almost periodic pattern $\Gamma$ possesses a uniform density, that is, that the limit 
\[D(\Gamma) = \lim_{R\to+\infty}\frac{\card(\Gamma\cap [-R,R]^n)}{\card(\Z^n\cap [-R,R]^n)}\]
is well defined (Corollary \ref{corolimitexist}). 

Recall that the \emph{discretization} of a linear map $A\in GL_n(\R)$ is the map $\widehat A = \pi\circ A : \Z^n\to\Z^n$, where $\pi$ is a projection from $\R^n$ to the nearest element of $\Z^n$ (see Definition~\ref{DefDiscrLin}). The definition of an almost periodic pattern is supported by Theorem~\ref{imgquasi}, which states that the image of an almost periodic pattern by the discretization of a linear map is still an almost periodic pattern. In particular, the successive images of the integer lattice $\Z^n$ by the discretizations of linear maps are almost periodic patterns, and possess a uniform density.

We then study the frequency of any difference $v$ in the almost periodic pattern $\Gamma$, defined by $\rho_\Gamma(v) = D\big((\Gamma-v)\cap\Gamma\big)$ (Definition~\ref{DefDiff}); in particular we prove a Minkowski-like theorem for these differences (Theorem \ref{MinkAlm}, obtained in collaboration with \'E.~Joly): if $S$ is a centrally symmetric convex body, then the sum of the frequency of the differences $v\in S$ is bigger than a quantity which depends linearly on the measure of $S$.

Finally, we introduce the notion of model set (Definition \ref{DefModel}): the model set modelled on a lattice $\Lambda$ of $\R^{m+n}$ and on a ``regular'' set $W\subset\R^m$ is the projection on the $n$ last coordinates of the points of $\Lambda$ whose projection on the $m$ first coordinates belongs to $W$. We state that the images sets $(\widehat{A_k}\circ\cdots\circ\widehat{A_1}) (\Z^n)$ are model sets, and prove that these sets are almost periodic patterns (Theorem~\ref{ModelAlmost}, obtained in collaboration with Y.~Meyer).
\bigskip

We fix once for all an integer $n\ge 1$. We will denote by $\llbracket a, b \rrbracket$\index{$\llbracket\cdot\rrbracket$} the integer segment $[a,b]\cap\Z$. In this part, every ball will be taken with respect to the infinite norm; in particular, for $x = (x_1,\cdots,x_n)$, we will have\index{$B(x,R)$}
\[B(x,R) = B_\infty(x,R) = \big\{y=(y_1,\cdots,y_n)\in\R^n\mid \forall i\in \llbracket 1, n\rrbracket, |x_i-y_i|<R\big\}.\]
We will also denote $B_R = B(0,R)$\index{$B_R$}. Finally, we will denote by $\lfloor x \rfloor$\index{$\lfloor \cdot \rfloor$} the biggest integer that is smaller than $x$ and $\lceil x \rceil$\index{$\lceil \cdot \rceil$} the smallest integer that is bigger than $x$. For a set $B\subset \R^n$, we will denote $[B] = B\cap \Z^n$.\index{$B$@$[B]$}

\section[Almost periodic patterns]{Almost periodic patterns: definitions and first properties}

In this section, we define the notion of almost periodic pattern and prove that these sets possess a uniform density.

\begin{definition}
Let $\Gamma$ be a subset of $\R^n$.
\begin{itemize}
\item We say that $\Gamma$ is \emph{relatively dense} if there exists $R_\Gamma>0$ such that each ball with radius at least $R_\Gamma$ contains at least one point of $\Gamma$.
\item We say that $\Gamma$ is a \emph{uniformly discrete} if there exists $r_\Gamma>0$ such that each ball with radius at most $r_\Gamma$ contains at most one point of $\Gamma$.
\end{itemize}
The set $\Gamma$ is called a \emph{Delone} set if it is both relatively dense and uniformly discrete.
\end{definition}

\begin{definition}
For a discrete set $\Gamma\subset \R^n$ and $R\ge 1$, we define the uniform $R$-density:\index{$D_R^+$}
\[D_R^+(\Gamma) = \sup_{x\in\R^n} \frac{\card\big(B(x,R)\cap \Gamma\big)}{\card\big(B(x,R)\cap\Z^n\big)},\]
and the uniform upper density:\index{$D^+$}
\[D^+(\Gamma) = \underset{R\to +\infty}{\overline\lim} D_R^+(\Gamma).\]
\end{definition}

Remark that if $\Gamma\subset \R^n$ is Delone for the parameters $r_\Gamma$ and $R_\Gamma$, then its upper density satisfies:
\[\frac{1}{(2R_\Gamma+1)^n} \le D^+(\Gamma) \le \frac{1}{(2r_\Gamma+1)^n}.\]

We can now define the notion of almost periodic pattern that we will use throughout this chapter. Roughly speaking, an almost periodic pattern $\Gamma$ is a set for which there exists a relatively dense set of translations of $\Gamma$, where a vector $v$ is a translation of $\Gamma$ if $\Gamma-v$ is equal to $\Gamma$ up to a set of upper density smaller than $\varep$. More precisely, we state the following definition.

\begin{definition}\label{DefAlmPer}\index{$\mathcal N_\varep$}
A Delone set $\Gamma$ is an \emph{almost periodic pattern} if for every $\varep>0$, there exists $R_\varep>0$ and a relatively dense set $\mathcal N_\varep$, called the \emph{set of $\varep$-translations} of $\Gamma$, such that
\begin{equation}\label{EqAlmPer}
\forall R\ge R_\varep,\  \forall v\in\mathcal N_\varep,\  D_R^+\big( (\Gamma+v)\Delta \Gamma \big) <\varep.
\end{equation}
\end{definition}

Remark that if $\Gamma$ is a subset of $\Z^n$ with positive upper density and which satisfies the condition of this definition, then it is a Delone set. In the sequel, we will only use this definition for subsets of the lattice $\Z^n$. Remark that by \cite[Theorem 3]{MR2869161}, an almost periodic pattern which is also a Meyer set\footnote{A set $\Gamma$ is called a \emph{Meyer set} if $\Gamma-\Gamma$ is a Delone set. It is equivalent to ask that there exists a finite set $F$ such that $\Gamma-\Gamma \subset \Gamma + F$ (see \cite{MR1400744}).} is ``almost included'' in a finite union of lattices.

Of course, every lattice, or every finite union of translates of a given lattice, is an almost periodic pattern. We will see in next section a large class of examples of almost periodic patterns: images of $\Z^n$ by discretizations of linear maps.

This definition is stronger than the following one, that we initially used for this study.

\begin{definition}\label{wap}
We say that a Delone set $\Gamma$ is \emph{weakly almost periodic} if for every $\varep>0$, there exists $R>0$ such that for every $x,y\in\R^n$, there exists $v\in\R^n$ such that
\begin{equation}\label{EqWeakAlmPer}
\frac{\card\Big( \Big(B(x,R)\cap\Gamma\Big) \Delta \Big(\big(B(y,R)\cap\Gamma\big)-v\Big) \Big)}{\card(B_R\cap\Z^n)} \le \varep.
\end{equation}
Remark that \emph{a priori}, the vector $v$ is different from $y-x$.
\end{definition}

We had defined this concept because it seemed to us that it was the weakest to imply the existence of a uniform density. Unfortunately, this notion is not very convenient to manipulate and we have not succeeded to prove that it is stable under the action of discretization of linear maps. Of course, we have the following result (see also the addendum \cite{GM} of \cite{MR2876415} for more details on the subject).

\begin{prop}\label{WeakAlmPer}
Every almost periodic pattern is weakly almost periodic.
\end{prop}

We do not know if the converse is true or not (that is, if there exists weakly almost periodic sets that are not almost periodic patterns).

\begin{proof}[Proof of Proposition \ref{WeakAlmPer}]
We prove that an almost periodic pattern satisfies Equation \eqref{EqWeakAlmPer} for $x=0$, the general case being obtained by applying this result twice.

Let $\Gamma$ be an almost periodic pattern and $\varep>0$. Then by definition, there exists $R_\varep>0$ and a relatively dense set $\mathcal N_\varep$ (for a parameter $R_{\mathcal N_\varep}>0$) such that
\begin{equation}\label{eq2.1}
\forall R\ge R_\varep,\  \forall v\in\mathcal N_\varep,\  D_{R}^+\big( (\Gamma+v)\Delta \Gamma \big) <\varep.
\end{equation}
Moreover, as $\Gamma$ is Delone, there exists $r_\Gamma>0$ such that each ball with radius smaller than $r_\Gamma$ contains at most one point of $\Gamma$.

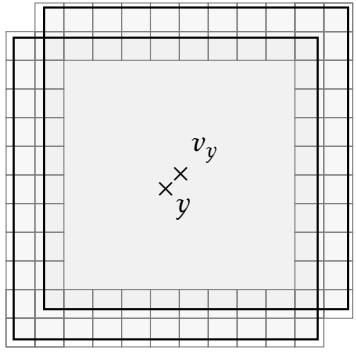
\begin{figure}[t]
\begin{center}

\begin{tikzpicture}[scale=1]
\draw[fill=gray!20!white, opacity=0.3] (-1.8,-1.8) rectangle (2.2,2.2);
\draw[fill=gray!20!white, opacity=0.3] (-1.4,-1.4) rectangle (2.6,2.6);
\node (A) at (.2,.2) {$\times$};
\draw (A) node[below right]{$y$};
\node (B) at (.4,.4) {$\times$};
\draw (B) node[above right]{$v_y$};
\foreach\i in {-3,...,7}{
\foreach\j in {6,...,7}{
\draw[color=gray] (\i*.38,\j*.38) rectangle (\i*.38-.38,\j*.38-.38);
\draw[color=gray] (\i*.38-0.38,\j*.38-3.8) rectangle (\i*.38-.76,\j*.38-3.8-.38);
}}
\foreach\i in {-4,...,-3}{
\foreach\j in {-4,...,6}{
\draw[color=gray] (\i*.38,\j*.38) rectangle (\i*.38-.38,\j*.38-.38);
\draw[color=gray] (\i*.38+3.8-0.38,\j*.38+0.38) rectangle (\i*.38+3.8,\j*.38);
}}
\draw[thick] (-1.8,-1.8) rectangle (2.2,2.2);
\draw[thick] (-1.4,-1.4) rectangle (2.6,2.6);
\end{tikzpicture}
\caption[Covering the set $B(y,R)\Delta B(v_y,R)$ by cubes of radius $r_\Gamma$]{Covering the set $B(y,R)\Delta B(v_y,R)$ by cubes of radius $r_\Gamma$.}\label{DiffCub}
\end{center}
\end{figure}

As $\mathcal N_\varep$ is relatively dense, for every $y\in \R^n$, there exists $v_y\in -\mathcal N_\varep$ such that $d_\infty(y, v_y)<R_{\mathcal N_\varep}$. This $v_y$ is the vector $v$ we look for to have the property of Definition~\ref{wap}. Indeed, by triangle inequality, for every $R\ge R_\varep$, we have
\begin{align}\label{EqTriangleIn}
\card\Big(\big(B(0,R)\cap\Gamma\big) & \Delta \big((B(y,R)\cap\Gamma)-v_y\big)\Big)\nonumber\\
         \le & \card\Big(\big(B(0,R)\cap\Gamma\big) \Delta\big((B(v_y,R)\cap\Gamma) - v_y\big)\Big)\\
				     & + \card\Big(\big(B(v_y,R)\cap\Gamma\big) \Delta\big(B(y,R)\cap\Gamma\big)\Big).\nonumber
\end{align}
By Equation~\eqref{eq2.1}, the first term of the right side of the inequality is smaller than $\varep\card\big(B(0,R)\cap\Z^n\big)$. It remains to bound the second one.

For every $y\in \R^n$, as $d_\infty(y, v_y)<R_{\mathcal N_\varep}$, the set $B(y,R)\Delta B(v_y,R)$ is covered by
\[\frac{2n(R+r_\Gamma)^{n-1}(R_{\mathcal N_\varep} + r_\Gamma)}{r_\Gamma^n}\]
disjoint cubes of radius $r_\Gamma$ (see Figure \ref{DiffCub}). Thus, as each one of these cubes contains at most one point of $\Gamma$, this implies that
\[\card\Big(\big(B(y,R)\Delta B(v_y,R)\big)\cap \Gamma\Big) \le 2n\frac{(R+r_\Gamma)^{n-1}(R_{\mathcal N_\varep} + r_\Gamma)}{r_\Gamma^n}.\]
Increasing $R_\varep$ if necessary, for every $R\ge R_\varep$, we have 
\[2n\frac{(R+r_\Gamma)^{n-1}(R_{\mathcal N_\varep} + r_\Gamma)}{r_\Gamma^n} \le \varep \card\big(B(0,R) \cap \Z^n\big),\]
so,
\[\card\Big(\big(B(y,R)\Delta B(v_y,R)\big)\cap \Gamma\Big) \le \varep \card\big(B(0,R) \cap \Z^n\big).\]
This bounds the second term of Equation~\eqref{EqTriangleIn}. We finally get
\[\card\Big(\big(B(0,R)\cap\Gamma\big) \Delta\big(B(y,R)\cap\Gamma-v_y\big)\Big) \le 2 \varep \card\Big(B(0,R) \cap \Z^n\Big),\]
which proves the proposition.
\end{proof}

Weakly almost periodic sets -- and in particular almost-periodic patterns -- have a regular enough behaviour at the infinity to possess a density.

\begin{prop}\label{limitexist}
Let $\Gamma$ be a weakly almost periodic set. Then the uniform upper density of $\Gamma$ is a limit, more precisely:\index{$D$}
\[D^+(\Gamma) = \underset{R\to +\infty}{\lim} D_R^+(\Gamma) = D(\Gamma).\]
Moreover, we have uniformity of the convergence of the density with respect to the base point, that is: for every $x\in\R^n$, we have
\[D(\Gamma) = \lim_{R\to +\infty}\frac{\card\big(B(x,R)\cap \Gamma\big)}{\card\big(B(x,R)\cap \Z^n\big)}.\]
In this case, we say that $D(\Gamma)$ is the \emph{uniform density} of $\Gamma$.
\end{prop}

Combined with Proposition \ref{WeakAlmPer}, this proposition directly implies the following corollary.

\begin{coro}\label{corolimitexist}
Let $\Gamma$ be an almost periodic pattern. Then the uniform upper density of $\Gamma$ is a limit; in other words $\Gamma$ possesses a uniform density.
\end{coro}

\begin{rem}\label{Jordan}
The same proof also shows that the same property holds if instead of considering the density $D^+$, we take a \emph{Jordan-measurable}\footnote{We say taht a set $J$ is Jordan-measurable if for every $\varep>0$, there exists $\eta>0$ such that there exists two disjoint unions $\mathcal C$ and $\mathcal C'$ of cubes of radius $\eta$, such that $\mathcal C\subset J\subset\mathcal C'$, and that $\Leb(\mathcal C'\setminus\mathcal C)<\varep$.} set $J$ and consider the density $D_J^+(\Gamma)$ of a set $\Gamma\subset \Z^n$ defined by\index{$D_A^+$}
\[D_J^+(\Gamma) = \underset{R\to +\infty}{\overline\lim} \sup_{x\in\R^n} \frac{\card\big(J_R\cap \Gamma\big)}{\card\big(J_R\cap\Z^n\big)},\]
where $J_R$ denotes the set of points $x\in\R^n$ such that $x/R\in J$.
\end{rem}

\begin{proof}[Proof of Proposition \ref{limitexist}]
Let $\Gamma$ be a weakly almost periodic set and $\varep>0$. Then, by definition, there exists $R>0$ such that for all $x,y\in\R^n$, there exists $v\in \R^n$ such that Equation \eqref{EqWeakAlmPer} holds. We take a ``big'' $M\in\R$, $x\in\R^n$ and $R'\ge MR$. We use the tiling of $\R^n$ by the collection of squares $\{B(Ru,R)\}_{u\in (2\Z)^n}$ and the Equation \eqref{EqWeakAlmPer} (applied to the radius $R'$ and the points $0$ and $Ru$) to find the number of points of $\Gamma$ that belong to $B(x,R')$: as $B(x,R')$ contains at least $\lfloor M\rfloor^n$ disjoint cubes $B(Ru,R)$ and is covered by at most $\lceil M\rceil^n$ such cubes, we get (recall that $B_R = B(0,R)$)
\begin{flalign*}
\frac{\lfloor M\rfloor^n \big(\card (B_R\cap\Gamma)-2\varep\card(B_R\cap\Z^n)\big)}{\lceil M\rceil^n\card(B_R\cap\Z^n)} & & &
\end{flalign*}
\[\le \frac{\card\big(B(x,R')\cap \Gamma\big)}{\card\big(B(x,R')\cap \Z^n\big)} \le\]
\begin{flalign*}
 & & & \frac{\lceil M\rceil^n \big(\card (B_R\cap\Gamma)+2\varep\card (B_R \cap \Z^n)\big)}{\lfloor M\rfloor^n\card (B_R\cap \Z^n )},
\end{flalign*}
thus
\begin{flalign*}
\frac{\lfloor M\rfloor^n}{\lceil M\rceil^n}\left(\frac{\card (B_R\cap\Gamma)}{\card (B_R\cap \Z^n)}-2\varep\right) & & &
\end{flalign*}
\[\le \frac{\card\big(B(x,R')\cap \Gamma\big)}{\card\big(B(x,R')\cap \Z^n\big)} \le\]
\begin{flalign*}
 & & & \qquad \qquad \qquad \qquad \frac{\lceil M\rceil^n}{\lfloor M\rfloor^n}\left(\frac{\card (B_R\cap\Gamma)}{\card (B_R\cap \Z^n)}+2\varep\right).
\end{flalign*}
For $M$ large enough, this ensures that for every $R'\ge MR$ and every $x\in \R^n$, the density
\[\frac{\card\big(B(x,R')\cap \Gamma\big)}{\card\big(B(x,R')\cap \Z^n\big)}\quad  \text{is close to}  \quad \frac{\card (B_R\cap\Gamma)}{\card (B_R\cap \Z^n)};\]
this finishes the proof of the proposition.
\end{proof}

We end this section by an easy lemma which asserts that for $\varep$ small enough, the set of translations $\mathcal N_\varep$ is ``stable under additions with a small number of terms''.

\begin{lemme}\label{arithProg}
Let $\Gamma$ be an almost periodic pattern, $\varep>0$ and $\ell\in\N$. Then if we set $\varep'=\varep/\ell$ and denote by $\mathcal N_{\varep'}$ the set of translations of $\Gamma$ and $R_{\varep'}>0$ the corresponding radius for the parameter $\varep'$, then for every $k\in\llbracket 1,\ell \rrbracket$ and every $v_1,\cdots,v_\ell\in\mathcal N_{\varep'}$, we have
\[\forall R\ge R_{\varep'},\  D_R^+\Big( \big(\Gamma+\sum_{i=1}^\ell v_i\big)\Delta \Gamma \Big) <\varep.\]
\end{lemme}

\begin{proof}[Proof of Lemma \ref{arithProg}]
Let $\Gamma$ be an almost periodic pattern, $\varep>0$, $\ell\in\N$, $R_0>0$ and $\varep'=\varep/\ell$. Then there exists $R_{\varep'}>0$ such that
\[\forall R\ge R_{\varep'},\  \forall v\in\mathcal N_{\varep'},\  D_R^+\big( (\Gamma+v)\Delta \Gamma \big) <\varep'.\]
We then take $1\le k\le\ell$, $v_1,\cdots,v_k\in\mathcal N_{\varep'}$ and compute
\begin{align*}
D_R^+\Big( \big(\Gamma+\sum_{i=1}^k v_i\big)\Delta \Gamma \Big) & \le \sum_{m=1}^k D_R^+\Big( \big(\Gamma+\sum_{i=1}^m v_i\big)\Delta \big(\Gamma+\sum_{i=1}^{m-1} v_i\big) \Big)\\
             & \le \sum_{m=1}^k D_R^+\Big( \big((\Gamma+v_m)\Delta \Gamma\big) + \sum_{i=1}^{m-1} v_i \Big).
\end{align*}
By the invariance under translation of $D_R^+$, we deduce that
\begin{align*}
D_R^+\Big( \big(\Gamma+\sum_{i=1}^k v_i\big)\Delta \Gamma \Big) & \le \sum_{m=1}^k D_R^+ \big((\Gamma+v_m)\Delta \Gamma\big)\\
						 & \le k\varep'.
\end{align*}
As $k\le \ell$, this ends the proof.
\end{proof}

\begin{rem}\label{arithProg2}
In particular, this lemma implies that each set $\mathcal N_\varep$ contains arbitrarily large arithmetical progressions. More precisely, for every almost periodic pattern $\Gamma$, $\varep>0$ and $\ell\in\N$, if we set $\varep'=\varep/\ell$, then for every $k\in \llbracket 1,\ell\rrbracket$ and every $v\in \mathcal N_{\varep'}$, we have
\[\forall R\ge R_{\varep'},\  D_R^+\big( (\Gamma+kv )\Delta \Gamma \big) <\varep.\]

It also implies that the set $\mathcal N_\varep$ contains arbitrarily large patches of lattices of $\R^n$: for every almost periodic pattern $\Gamma$, $\varep>0$ and $\ell\in\N$, there exists $\varep'>0$ such that for every $k_i \in \llbracket -\ell,\ell\rrbracket$ and every $v_1,\cdots,v_n\in\mathcal N_{\varep'}$, we have
\[\forall R\ge R_{\varep'},\  D_{R}^+\Big( \big(\Gamma+\sum_{i=1}^n k_iv_i \big)\Delta \Gamma \Big) <\varep.\]
\end{rem}

\section{Almost periodic patterns and linear maps}

In this section, we prove that the notion of almost periodic pattern is invariant under discretizations of linear maps: the image of an almost periodic pattern by the discretization of a linear map is still an almost periodic pattern. First of all, we define precisely what we mean by discretization of a linear map.

\begin{definition}\label{DefDiscrLin}
The map $P : \R\to\Z$\index{$P$} is defined as a projection from $\R$ onto $\Z$. More precisely, for $x\in\R$, $P(x)$ is the unique\footnote{Remark that the choice of where the inequality is strict and where it is not is arbitrary.} integer $k\in\Z$ such that $k-1/2 < x \le k + 1/2$. This projection induces the map\index{$\pi$}
\[\begin{array}{rrcl}
\pi : & \R^n & \longmapsto & \Z^n\\
 & (x_i)_{1\le i\le n} & \longmapsto & \big(P(x_i)\big)_{1\le i\le n}
\end{array}\]
which is an Euclidean projection on the lattice $\Z^n$. Let $A\in M_n(\R)$. We denote by $\widehat A$ the \emph{discretization}\index{$\widehat A$} of the linear map $A$, defined by 
\[\begin{array}{rrcl}
\widehat A : & \Z^n & \longrightarrow & \Z^n\\
 & x & \longmapsto & \pi(Ax).
\end{array}\]
\end{definition}

The main result of this section is the following theorem.

\begin{theoreme}\label{imgquasi}
Let $\Gamma\subset\Z^n$ be an almost periodic pattern and $A\in GL_n(\R)$. Then the set $\widehat A(\Gamma)$ is an almost periodic pattern.
\end{theoreme}

\begin{figure}[t]
\begin{minipage}[c]{.33\linewidth}
	\includegraphics[width=\linewidth, trim = 1.5cm .5cm 1.5cm .5cm,clip]{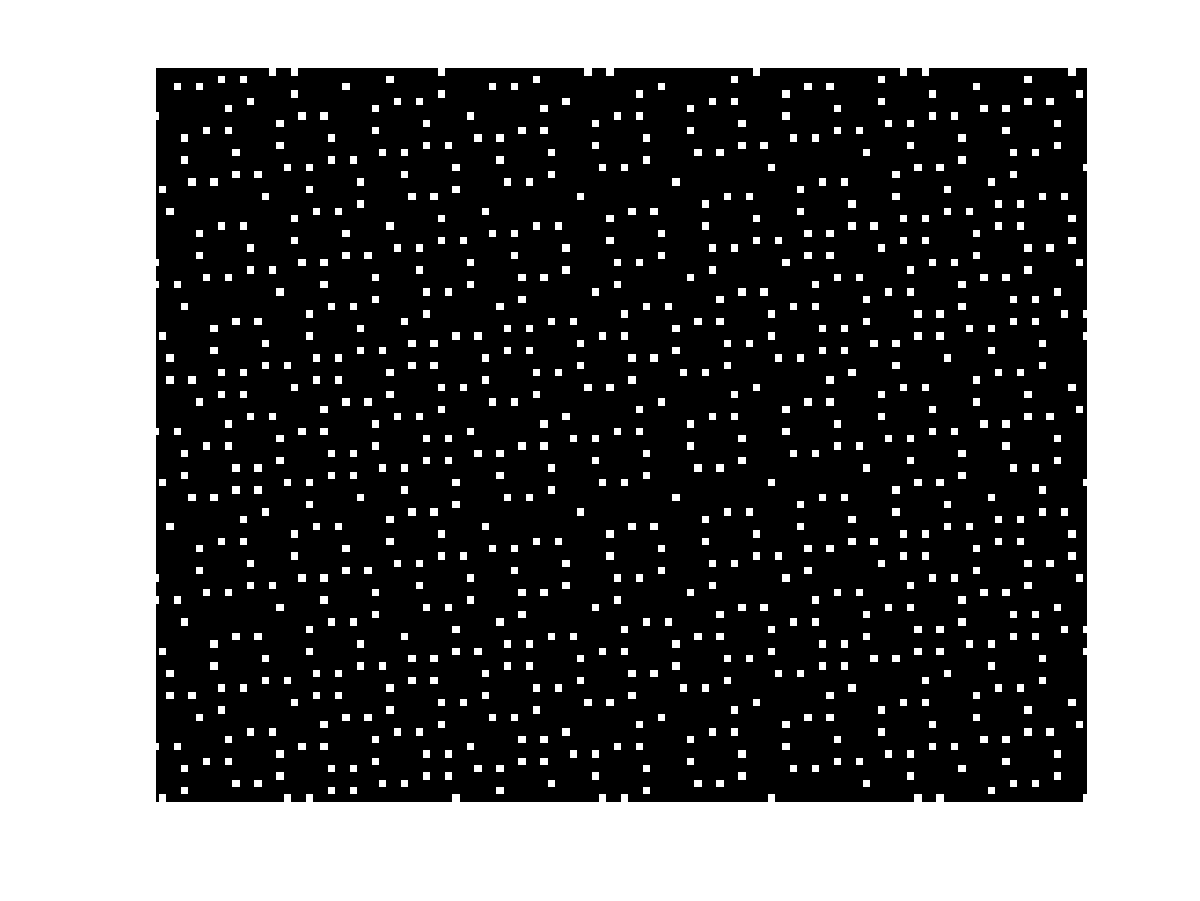}
\end{minipage}\hfill
\begin{minipage}[c]{.33\linewidth}
	\includegraphics[width=\linewidth, trim = 1.5cm .5cm 1.5cm .5cm,clip]{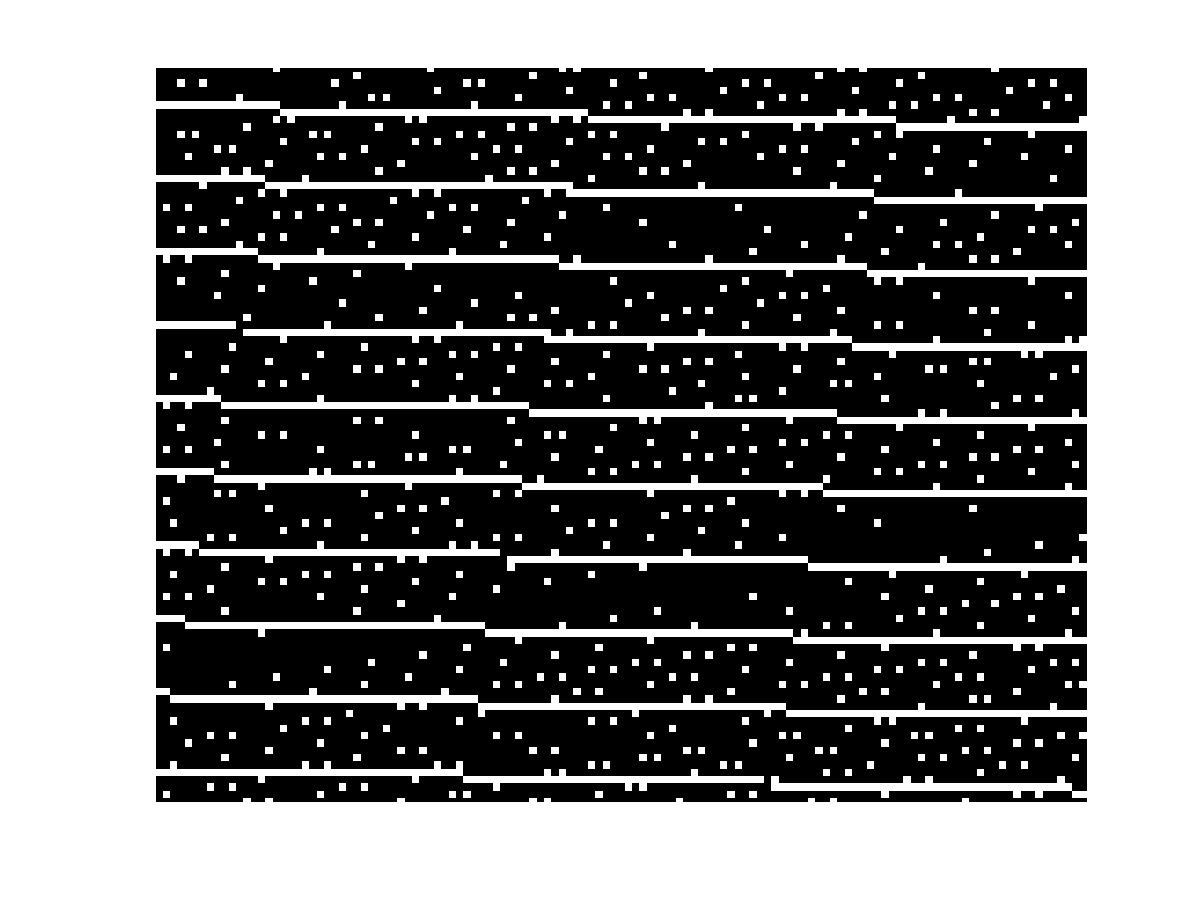}
\end{minipage}\hfill
\begin{minipage}[c]{.33\linewidth}
	\includegraphics[width=\linewidth, trim = 1.5cm .5cm 1.5cm .5cm,clip]{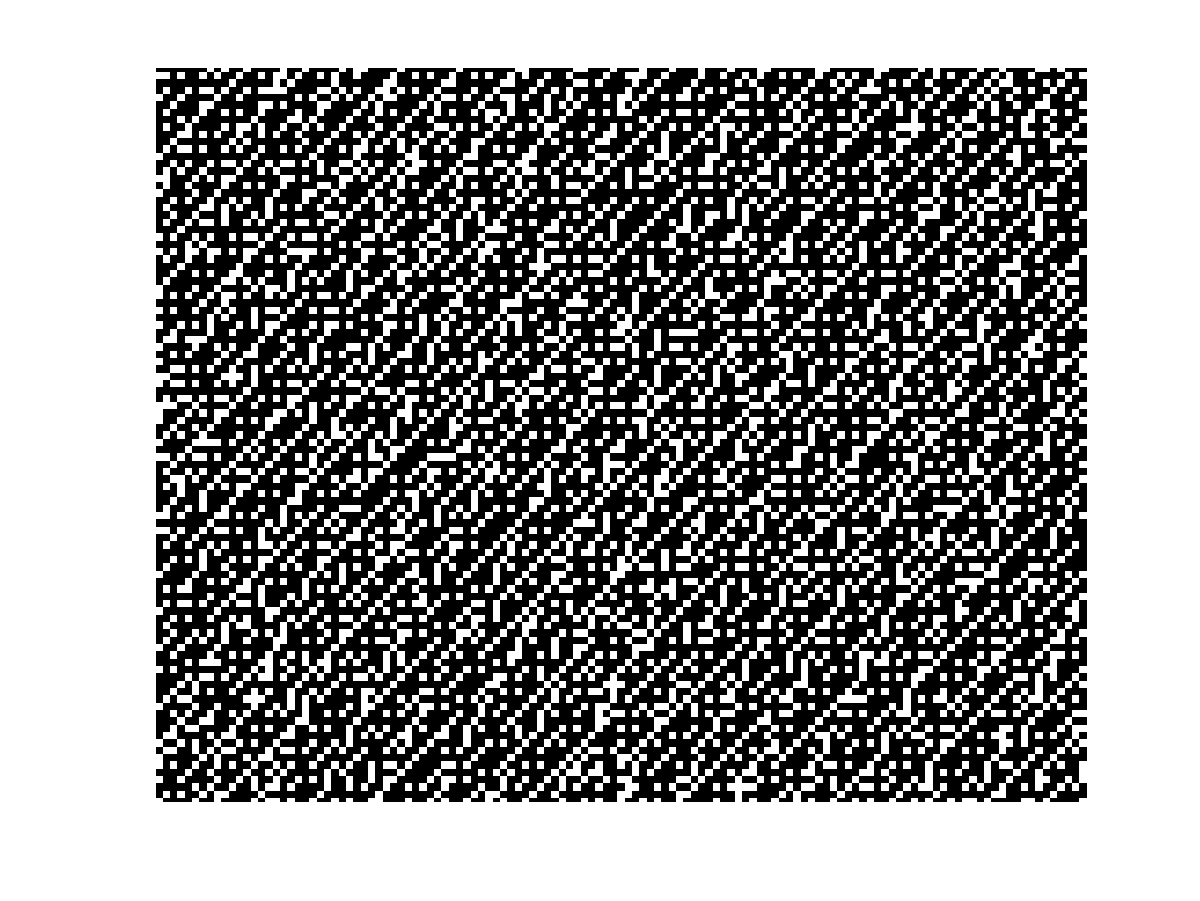}
\end{minipage}

\begin{minipage}[c]{.33\linewidth}
	\includegraphics[width=\linewidth, trim = 1.5cm .5cm 1.5cm .5cm,clip]{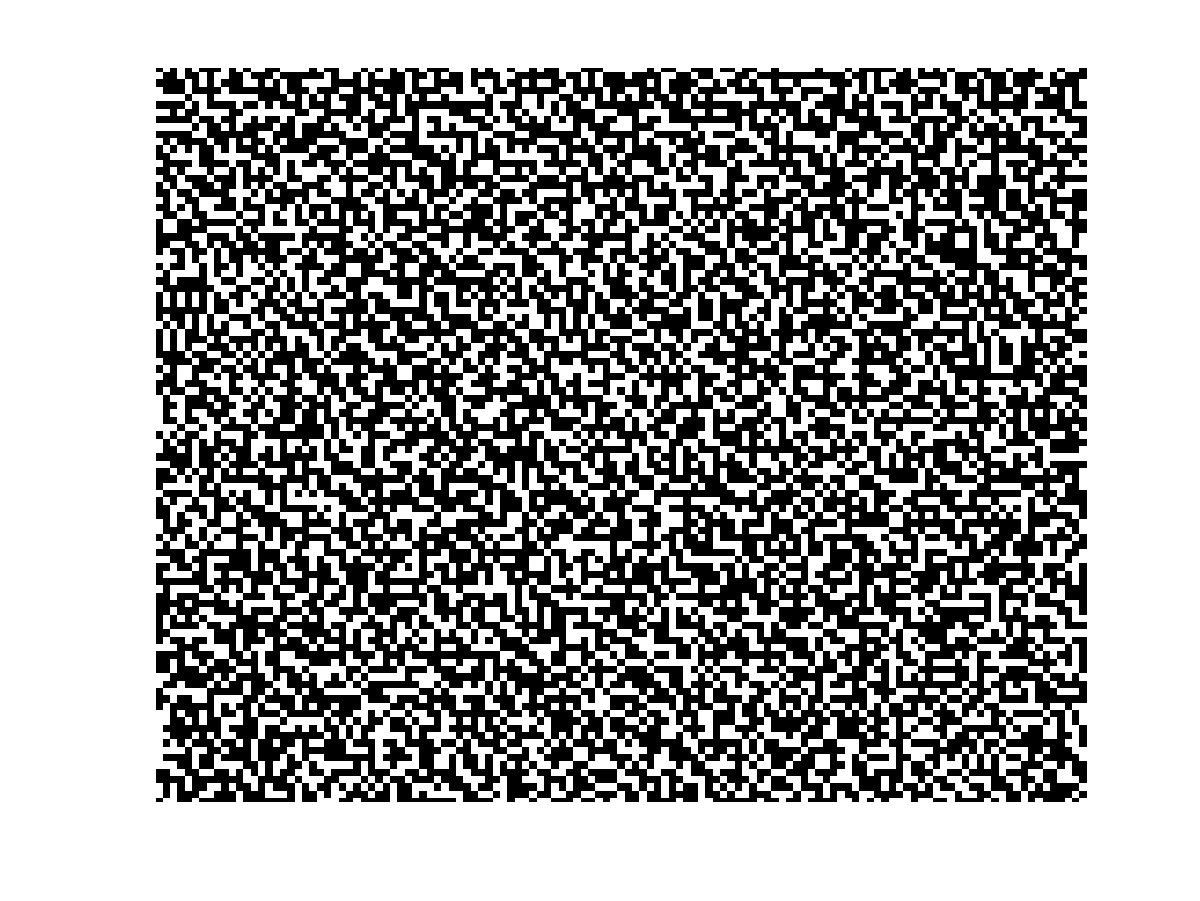}
\end{minipage}\hfill
\begin{minipage}[c]{.33\linewidth}
	\includegraphics[width=\linewidth, trim = 1.5cm .5cm 1.5cm .5cm,clip]{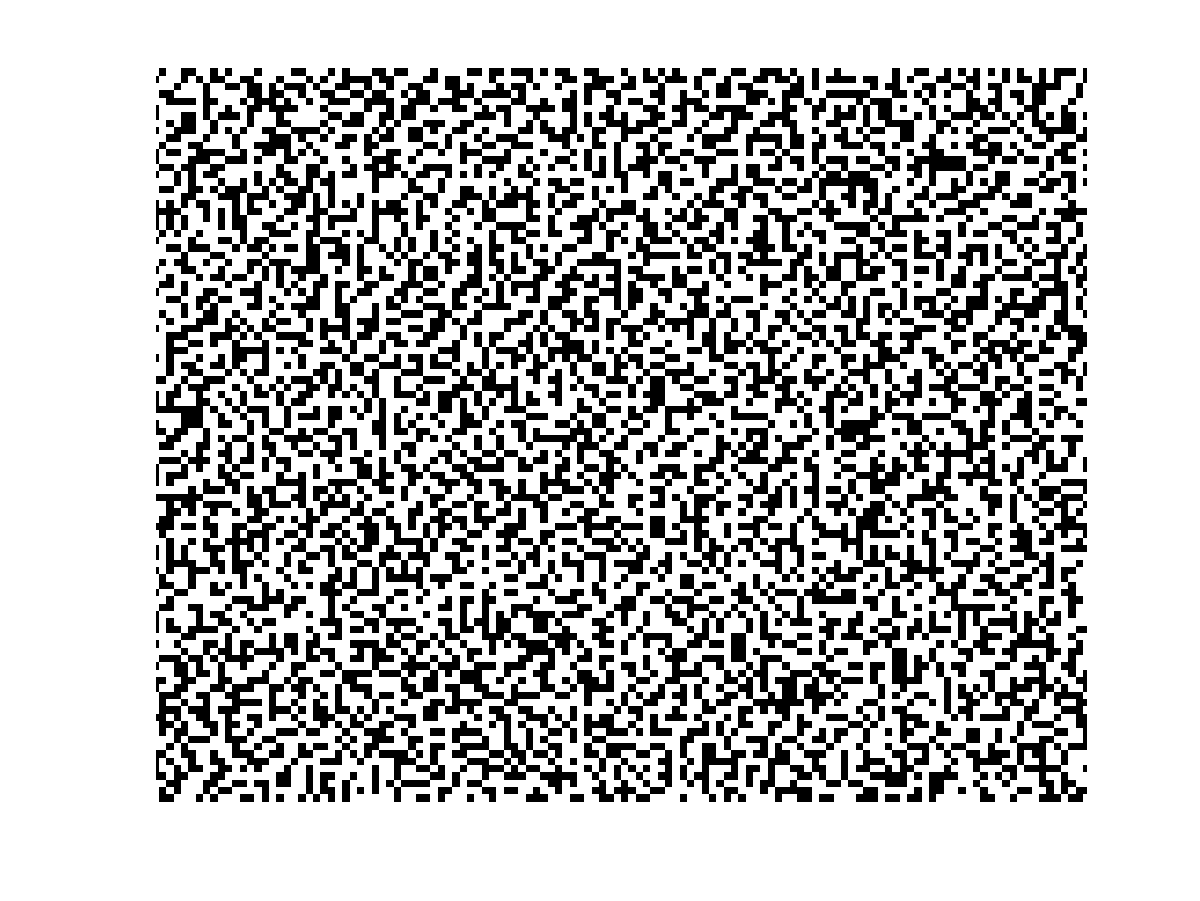}
\end{minipage}\hfill
\begin{minipage}[c]{.33\linewidth}
	\includegraphics[width=\linewidth, trim = 1.5cm .5cm 1.5cm .5cm,clip]{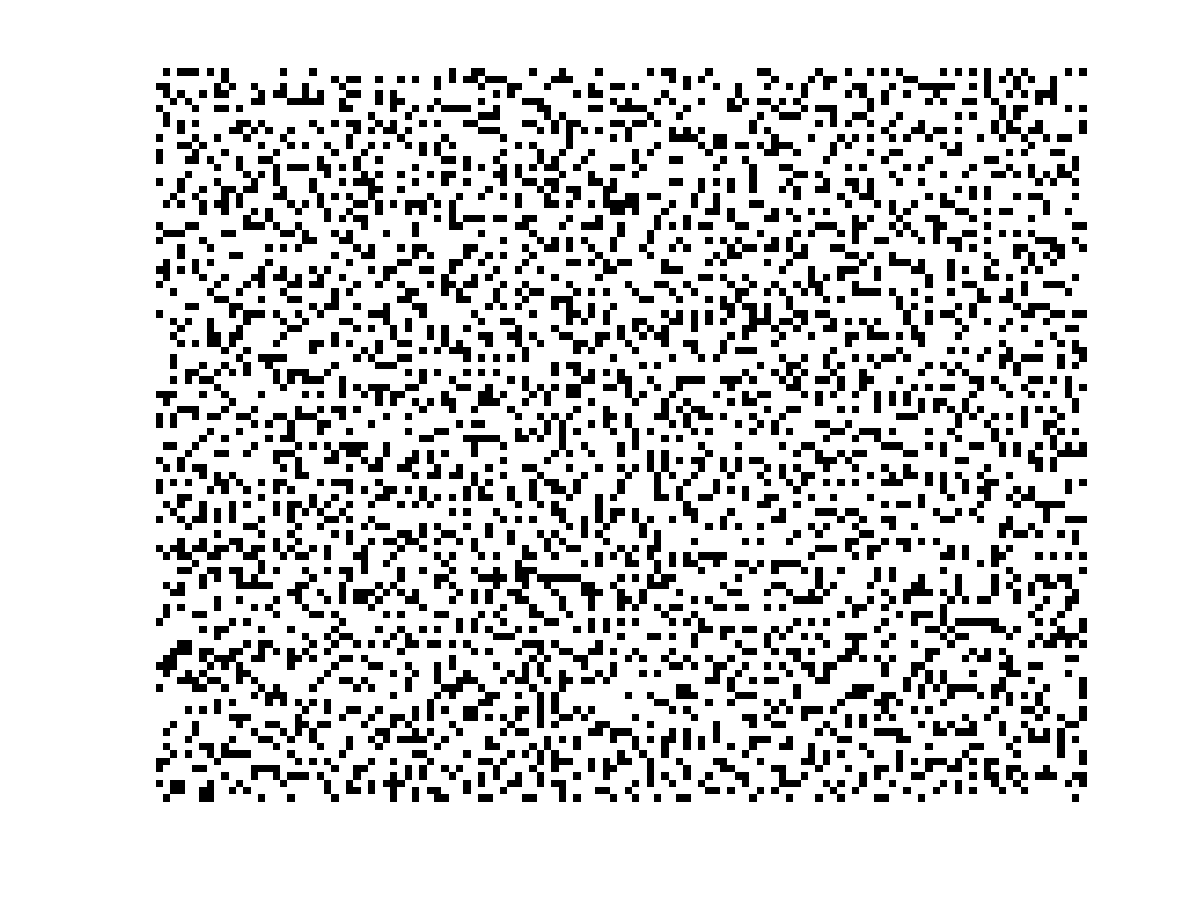}
\end{minipage}
\caption[Successive images of $\Z^2$ by discretizations of random matrices]{Successive images of $\Z^2$ by discretizations of random matrices in $SL_2(\R)$, a point is black if it belongs to $(\widehat{A_k}\circ\cdots\circ\widehat{A_1})(\Z^2)$. The $A_i$ are chosen randomly, using the singular value decomposition: they are chosen among the matrices of the form $R_\theta D_t R_{\theta'}$, with $R_\theta$ the rotation of angle $\theta$ and $D_t$ the diagonal matrix $\operatorname{Diag}(e^t,e^{-t})$, the $\theta$, $\theta'$ being chosen uniformly in $[0,2\pi]$ and $t$ uniformly in $[-1/2,1/2]$. From left to right and top to bottom, $k=1,\, 2,\, 3,\, 5,\, 10,\, 20$.}\label{ImagesSuitesMat}
\end{figure}

In particular, for every lattice $\Lambda$ of $\R^n$, the set $\pi(\Lambda)$ is an almost periodic pattern. More generally, given a sequence $(A_k)_{k\ge 1}$ of invertible matrices of $\R^n$, the successive images $(\widehat{A_k}\circ\cdots\circ \widehat{A_1})(\Z^n)$ are almost periodic patterns. See Figure~\ref{ImagesSuitesMat} for an example of the successive images of $\Z^2$ by a random sequence of bounded matrices of $SL_2(\R)$.

\begin{notation}\label{intelligent}
For $A\in GL_n(\R)$, we denote $A=(a_{i,j})_{i,j}$. We denote by $I_\Q(A)$\index{$I_\Q(A)$} the set of indices $i$ such that $a_{i,j}\in\Q$ for every $j\in\llbracket 1,n\rrbracket$
\end{notation}

The proof of Theorem~\ref{imgquasi} relies on the following remark:

\begin{rem}\label{pmtriv}
If $a\in\Q$, then there exists $q\in\N^*$ such that $\{ax\mid x\in\Z\}\subset \frac{1}{q}\Z$. On the contrary, if $a\in\R\setminus\Q$, then the set $\{ax\mid x\in\Z\}$ is equidistributed in $\R/\Z$.
\end{rem}

Thus, in the rational case, the proof will lie in an argument of periodicity. On the contrary, in the irrational case, the image $A(\Z^n)$ is equidistributed modulo $\Z^n$: on every large enough domain, the density does not move a lot when we perturb the image set $A(\Z^n)$ by small translations. This reasoning is formalized by Lemmas~\ref{tiroir} and \ref{équi}. 

More precisely, for $R$ large enough, we would like to find vectors $w$ such that $D^+_R\big((\pi(A\Gamma) +w)\Delta \pi(A\Gamma)\big)$ is small. We know that there exists vectors $v$ such that $D^+_R\big((\Gamma+v)\Delta\Gamma\big)$ is small; this implies that $D^+_R\big((A\Gamma+Av)\Delta A\Gamma\big)$ is small, thus that $D^+_R\big(\pi(A\Gamma+Av)\Delta \pi(A\Gamma)\big)$ is small. The problem is that in general, we do not have $\pi(A\Gamma+Av) = \pi(A\Gamma)+\pi(Av)$. However, this is true if we have $Av\in\Z^n$. Lemma~\ref{tiroir} shows that in fact, it is possible to suppose that $Av$ ``almost'' belongs to $\Z^n$, and Lemma~\ref{équi} asserts that this property is sufficient to conclude.

The first lemma is a consequence of the pigeonhole principle.

\begin{lemme}\label{tiroir}
Let $\Gamma\subset \Z^n$ be an almost periodic pattern, $\varep>0$, $\delta>0$ and $A\in GL_n(\R)$. Then we can suppose that the elements of $A(\mathcal N_\varep)$ are $\delta$-close to $\Z^n$. More precisely, there exists $R_{\varep,\delta}>0$ and a relatively dense set $\widetilde{\mathcal N}_{\varep,\delta}$\index{$\widetilde{\mathcal N}_{\varep,\delta}$} such that 
\[\forall R\ge R_{\varep,\delta},\  \forall v\in\widetilde{\mathcal N}_{\varep,\delta},\  D_R^+\big( (\Gamma+v)\Delta \Gamma \big) <\varep,\]
and that for every $v\in\widetilde{\mathcal N}_{\varep,\delta}$, we have $d_\infty(Av,\Z^n)<\delta$. Moreover, we can suppose that for every $i\in I_\Q(A)$ and every $v\in\widetilde{\mathcal N}_{\varep,\delta}$, we have $(Av)_i\in \Z$.
\end{lemme}

The second lemma states that in the irrational case, we have continuity of the density under perturbations by translations.

\begin{lemme}\label{équi}
Let $\varep>0$ and $A\in GL_n(\R)$. Then there exists $\delta>0$ and $R_0>0$ such that for all $w\in B_\infty(0,\delta)$ (such that for every $i\in I_\Q(A)$, $w_i=0$), and for all $R\ge R_0$, we have
\[D_R^+\big(\pi(A\Z^n) \Delta \pi(A\Z^n+w) \big) \le \varep.\]
\end{lemme}

\begin{rem}
In Section \ref{SecCont} of Chapter \ref{Souris}, we will present an example which shows that the assumption ``for every $i\in I_\Q(A)$, $v_i=0$'' is necessary to obtain the conclusion of the lemma.
\end{rem}

\begin{rem}\label{RemContTrans}
When $I_\Q(A) = \emptyset$, and in particular when $A$ is totally irrational (see Definition~\ref{DefMeanRate}), the map $v\mapsto \tau(A+v)$ is continuous in 0; the same proof as that of this lemma implies that this function is globally continuous.
\end{rem}

We begin by the proofs of both lemmas, and prove Theorem~\ref{imgquasi} thereafter.

\begin{proof}[Proof of Lemma \ref{tiroir}]
Let us begin by giving the main ideas of the proof of this lemma. For $R_0$ large enough, the set of remainders modulo $\Z^n$ of vectors $Av$, where $v$ is a $\varep$-translation of $\Gamma$ belonging to $B_{R_0}$, is close to the set of remainders modulo $\Z^n$ of vectors $Av$, where $v$ is any $\varep$-translation of $\Gamma$. Moreover (by the pigeonhole principle), there exists an integer $k_0$ such that for each $\varep$-translation $v\in B_{R_0}$, there exists $k\le k_0$ such that $A(k v)$ is close to $\Z^n$. Thus, for every $\varep$-translation $v$ of $\Gamma$, there exists a $(k_0-1)\varep$-translation $v' = (k-1)v$, belonging to $B_{k_0 R_0}$, such that $A(v+v')$ is close to $\Z^n$. The vector $v+v'$ is then a $k_0\varep$-translation of $\Gamma$ (by additivity of the translations) whose distance to $v$ is smaller than $k_0 R_0$.
\bigskip

We now formalize these remarks. Let $\Gamma$ be an almost periodic pattern, $\varep>0$ and $A\in GL_n(\R)$. First of all, we apply the pigeonhole principle. We partition the torus $\R^n/\Z^n$ into squares whose sides are smaller than $\delta$; we can suppose that there are at most  $\lceil 1/\delta\rceil^n$ such squares. For $v\in \R^n$, we consider the family of vectors $\{A(kv)\}_{0\le k\le \lceil 1/\delta\rceil^n}$ modulo $\Z^n$. By the pigeonhole principle, at least two of these vectors, say $A(k_1v)$ and $A(k_2v)$, with $k_1<k_2$, lie in the same small square of $\R^n/\Z^n$. Thus, if we set $k_v = k_2-k_1$ and $\ell = \lceil 1/\delta\rceil^n$, we have
\begin{equation}\label{eqdistZ}
1\le k_v\le \ell \quad \text{and} \quad d_\infty\big(A(k_vv),\Z^n\big)\le\delta.
\end{equation}
To obtain the conclusion in the rational case, we suppose in addition that $v\in q\Z^n$, where $q\in\N^*$ is such that for every $i\in I_\Q(A)$ and every $1\le j\le n$, we have $q\, a_{i,j}\in\Z$ (which is possible by Remark~\ref{arithProg2}).

We set $\varep'=\varep/\ell$. By the definition of an almost periodic pattern, there exists $R_{\varep'}>0$ and a relatively dense set ${\mathcal N}_{\varep'}$ such that Equation \eqref{EqAlmPer} holds for the parameter $\varep'$:
\begin{equation}\label{EqAlmPer3}\tag{\ref{EqAlmPer}'}
\forall R\ge R_{\varep'},\  \forall v\in\mathcal N_{\varep'},\  D_R^+\big( (\Gamma+v)\Delta \Gamma \big) <\varep',
\end{equation}

We now set
\[P = \big\{Av\operatorname{mod} \Z^n \mid v\in {\mathcal N}_{\varep'}\big\} \quad \text{and} \quad P_R = \big\{Av\operatorname{mod} \Z^n \mid v\in \mathcal N_{\varep'}\cap B_R\big\}.\]
We have $\bigcup_{R>0} P_R = P$, so there exists $R_0>R_{\varep'}$ such that $d_H(P,P_{R_0})<\delta$ (where $d_H$\index{$d_H$} denotes Hausdorff distance). Thus, for every $v\in\mathcal N_{\varep'}$, there exists $v'\in \mathcal N_{\varep'}\cap B_{R_0}$ such that
\begin{equation}\label{eq666}
d_\infty(Av-Av',\Z^n)<\delta.
\end{equation}

We then remark that for every $v'\in {\mathcal N}_{\varep'}\cap B_{R_0}$, if we set $v'' = (k_{v'}-1)v'$, then by Equation \eqref{eqdistZ}, we have
\[d_\infty(Av' + Av'',\Z^n) = d_\infty\big(A(k_{v'}v'),\Z^n\big)\le\delta.\]
Combining this with Equation~\eqref{eq666}, we get
\[d_\infty(Av + Av'',\Z^n)\le 2\delta,\]
with $v''\in B_{\ell R_0}$. 

On the other hand, $k_{v'}\le \ell$ and Equation \eqref{EqAlmPer3} holds, so Lemma \ref{arithProg} (more precisely, the first point of Remark~\ref{arithProg2}) implies that $v''\in \mathcal N_\varep$, that is
\[\forall R\ge R_{\varep'},\  D_{R}^+\big( (\Gamma+ v'')\Delta \Gamma \big) <\varep.\]

In other words, for every $v\in\mathcal N_{\varep'}$, there exists $v''\in \mathcal N_\varep \cap B_{\ell R_0}$ (with $\ell$ and $R_0$ independent from $v$) such that $d_\infty\big(A(v+v''),\Z^n\big)<2\delta$. The set $\widetilde{\mathcal N}_{2\varep,2\delta}$ we look for is then the set of such sums $v+v''$.
\end{proof}

\begin{proof}[Proof of Lemma \ref{équi}]
Under the hypothesis of the lemma, for every $i\notin I_\Q(A)$, the sets
\[\left\{\sum_{j=1}^n a_{i,j} x_j\mid (x_j)\in\Z^n\right\},\]
are equidistributed modulo $\Z$. Thus, for all $\varep>0$, there exists $R_0>0$ such that for every $R\ge R_0$,
\[D_R^+\big\{v\in\Z^n \,\big|\, \exists i\notin I_\Q(A) : d\big((Av)_i,\Z+\frac12\big)\le \varep\big\} \le 2(n+1)\varep.\]
As a consequence, for all $w\in\R^n$ such that $\|w\|_\infty\le\varep/(2(n+1))$ and that $w_i=0$ for every $i\in I_\Q(A)$, we have
\[D_R^+\big(\pi(A\Z^n) \Delta \pi(A(\Z^n+w))\big)\le\varep.\]
Then, the lemma follows from the fact that there exists $\delta>0$ such that $\|A(w)\|_\infty\le \varep/(2(n+1))$ as soon as $\|w\|\le\delta$.
\end{proof}

\begin{proof}[Proof of Theorem \ref{imgquasi}]
Let $\varep>0$. Lemma \ref{équi} gives us a corresponding $\delta>0$, that we use to apply Lemma \ref{tiroir} and get a set of translations $\widetilde{\mathcal N}_{\varep,\delta}$. Then, for every $v\in \widetilde{\mathcal N}_{\varep,\delta}$, we write $\pi(Av) = Av + \big(\pi(Av)-Av\big) = Av + w$. The conclusions of Lemma~\ref{tiroir} imply that $\|w\|_\infty <\delta$, and that $w_i=0$ for every $i\in I_\Q(A)$.

We now explain why $\hat Av = \pi(Av)$ is a $\varep$-translation for the set $\widehat A(\Gamma)$. Indeed, for every $R\ge \max(R_{\varep,\delta},MR_0)$, where $M$ is the maximum of the greatest modulus of the eigenvalues of $A$ and of the greatest modulus of the eigenvalues of $A^{-1}$, we have
\begin{align*}
D^+_R \Big(\pi(A\Gamma) \Delta \big(\pi(A \Gamma)+\widehat Av\big)\Big) \le &\ D^+_R \Big(\pi(A\Gamma) \Delta \big(\pi(A \Gamma)+w\big)\Big)\\
                  & + D^+_R \Big(\big(\pi(A\Gamma) + w\big) \Delta \big(\pi(A \Gamma)+\widehat Av\big)\Big)
\end{align*}
(where $w=\pi(Av)-Av$). By Lemma \ref{équi}, the first term is smaller than $\varep$. For its part, the second term is smaller than
\[D^+_R\big((A\Gamma + Av) \Delta A \Gamma\big) \le M^2 D^+_{RM}\big((\Gamma + v) \Delta \Gamma\big),\]
which is smaller than $\varep$ because $v\in\mathcal N_\varep$.
\end{proof}

Theorem \ref{imgquasi} motivates the following definition.

\begin{definition}\label{DefTaux}
Let $A\in GL_n(\R)$. The \emph{rate of injectivity} of $A$ is the quantity\footnote{For the definition of the discretization $\widehat A$, see Definition~\ref{DefDiscrLin}. By definition, $[B_r] = B_\infty(0,R)\cap \Z^n$.}\index{$\tau$}
\[\tau(A) = \lim_{R\to +\infty} \frac{\card (\widehat A [B_R])}{\card [B_R]}\in]0,1].\]
More generally, for $A_1,\cdots,A_k \in GL_n(\R)$, we set\index{$\tau^k$}
\[\tau^k(A_1,\cdots,A_k) = \lim_{R\to +\infty} \frac{\card \big((\widehat{A_k}\circ\cdots\circ\widehat{A_1}) [B_R]\big)}{\card [B_R]}\in]0,1],\]
and for an infinite sequence $(A_k)_{k\ge 1}$ of invertible matrices, as the previous quantity is decreasing, we can define the \emph{asymptotic rate of injectivity}\index{$\tau^\infty$}
\[\tau^\infty\big((A_k)_{k\ge 1}\big) = \lim_{k\to +\infty}\tau^k(A_1,\cdots,A_k)\in[0,1].\]
\end{definition}

An easy calculation shows that the rate of injectivity can be deduced from the uniform density of the image $\widehat A(\Z^n)$ (in particular, it uses the fact that for every $A\in GL_n(\R)$ and every $R>0$, the set $A(B_R)$ is Jordan-measurable, see Remark~\ref{Jordan}).

\begin{prop}\label{TauDens}
For every matrix $A\in GL_n(\R)$, we have
\begin{equation}\label{EqTauDens}
\tau(A) = |\det(A)| D\big(\widehat A(\Z^n)\big).
\end{equation}
More generally, for every $A_1,\cdots,A_k \in GL_n(\R)$, we have
\[\tau^k(A_1,\cdots,A_k) = |\det(A_1)|\cdots |\det(A_k)| D\big((\widehat{A_k}\circ\cdots\circ\widehat{A_1})(\Z^n)\big).\]
\end{prop}

The same convergence holds for every affine map, thus we will also use the notion of rate of injectivity (and the notation $\tau$) in this more general context.

The quantity $\tau(A)$ does not change when we multiply $A$ on his right by an element of $SL_n(\Z)$. Similarly, $\tau^k(A_1,\cdots,A_k)$ in invariant under the multiplication of $A_1$ on his right by an element of $SL_2(\Z)$. It is difficult to detect more simple invariants for the rate $\tau$ (or even $\tau^k$), for example it has no obvious good behaviour with respect to the geodesic flow (or horocyclic flow) on the modular surface $SL_2(\R)/SL_2(\Z)$ (see Figure~\ref{FlowMean})

\begin{figure}[t]
\begin{center}
\includegraphics[width=.6\linewidth]{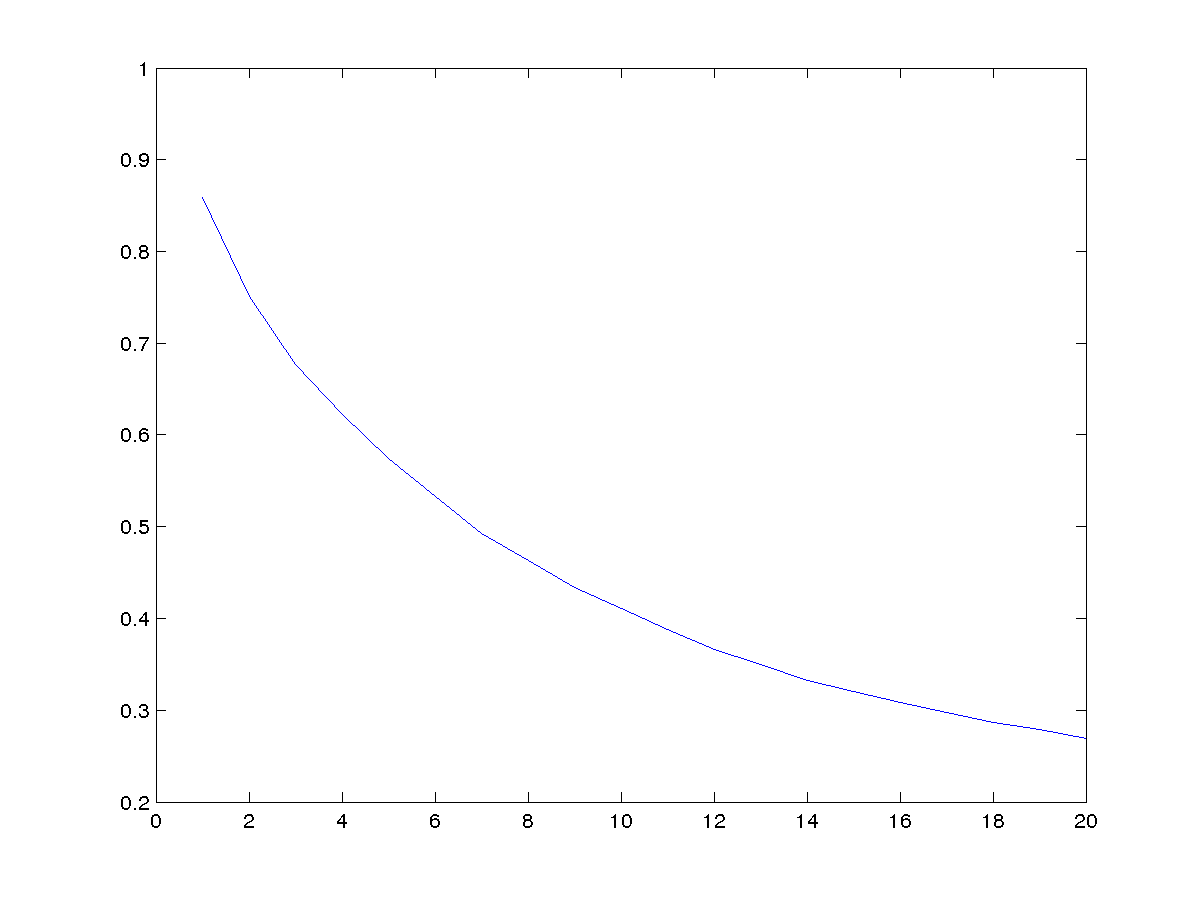}
\caption[Expectation of the rate of injectivity of a random sequences of matrices]{Expectation of the rate of injectivity of random sequences of matrices: the graphic represents the mean of the rate of injectivity $\tau^k(A_1,\cdots,A_k)$ depending on $k$, $1\le k\le 20$, for 50 random draws of matrices $A_i$. The $A_i$ are chosen randomly using the singular value decomposition: they are chosen among the matrices of the form $R_\theta D_t R_{\theta'}$, with $R_\theta$ the rotation of angle $\theta$ and $D_t$ the diagonal matrix $\operatorname{Diag}(e^t,e^{-t})$, and $\theta$, $\theta'$ chosen uniformly in $[0,2\pi]$ and $t$ chosen uniformly in $[-1/2,1/2]$. Note that the behaviour is experimentally not exponential. Also note that contrary to what happens in the case of isometries (see Figure~\ref{TauxSuiteRotations}), some hyperbolic-like phenomenons can occur in a random sequences of matrices of $SL_2(\R)$; so that in practical, it becomes difficult to plot a graph up to a long time (20 iterations for $SL_2(\R)$ versus 200 iterations for $O_2(\R)$) for reasons of memory constraints.}\label{TauxSuiteMat}
\end{center}
\end{figure}

The goal of Chapter \ref{Souris} will be to study in detail this quantity $\tau$, and in particular to prove that the behaviour of the rate suggested by Figure~\ref{TauxSuiteMat} actually holds, namely that the asymptotic rate of injectivity of a generic sequence of $SL_n(\R)$ is zero (Theorem~\ref{ConjPrincip}).

\bigskip

We end this section by a technical lemma that we will use a lot in next chapter. It expresses that given an almost periodic pattern $\Gamma$, a generic matrix $A\in GL_n(\R)$ is non resonant with respect to $\Gamma$.

\begin{lemme}\label{passifacil}
Let $\Gamma\subset \Z^n$ be an almost periodic pattern with positive uniform density. Then the set of $A\in GL_n(\R)$ (respectively $SL_n(\R)$, $O_n(\R)$) such that $A(\Gamma)$ is equidistributed modulo $\Z^n$ is generic. More precisely, for every $\varep>0$, there exists an open and dense set of $A\in GL_n(\R)$ (respectively $SL_n(\R)$, $O_n(\R)$) such that there exists $R_0>0$ such that for every $R>R_0$, the projection on $\R^n/\Z^n$ of the uniform measure on $A(\Gamma\cap B_R)$ is $\varep$-close to Lebesgue measure on $\R^n/\Z^n$.
\end{lemme}

\begin{rem}\label{Rempassifacil}
The proof also allows to suppose that the radius $R_0$ is uniform in a whole neighbourhood of every matrix $A$ of this open set of matrices.
\end{rem}

\begin{proof}[Proof of Lemma \ref{passifacil}]
During this proof, we consider a distance $\dist$ on $\Prb(\R^n/\Z^n)$ which is invariant under translations. We also suppose that this distance satisfies the following convexity inequality: if $\mu, \nu_1,\cdots,\nu_d\in\Prb(\R^n/\Z^n)$, then
\[\dist\left(\mu,\frac{1}{d}\sum_{i=1}^d \nu_i\right) \le \frac{1}{d}\sum_{i=1}^d \dist(\mu_,\nu_i).\]
For the simplicity of the notations, when $\mu$ and $\nu$ have not total mass 1, we will denote by $\dist(\mu,\nu)$ the distance between the normalizations of $\mu$ and $\nu$.

We consider the set $\mathcal U_\varep$ of matrices $A\in GL_n(\R)$ satisfying: there exists $R_0>0$ such that for all $R\ge R_0$,
\[\dist \left(\Leb_{\R^n/\Z^n}, \sum_{x\in B_R\cap \Gamma} \bar\delta_{Ax} \right) <\varep,\]
where $\bar\delta_x$\index{$\bar\delta$} is the Dirac measure of the projection of $x$ on $\R^n/\Z^n$. We show that for every $\varep>0$, $\mathcal U_\varep$ contains an open dense set. Then, the set $\bigcap_{\varep>0}\mathcal U_\varep$ will be a $G_\delta$ dense set made of matrices $A\in GL_n(\R)$ such that $A(\Gamma)$ is well distributed.

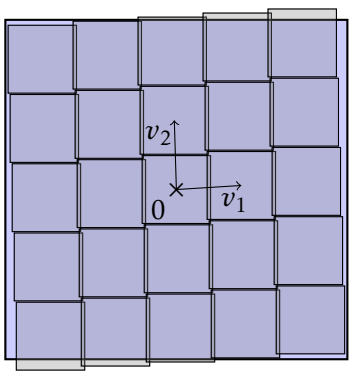
\begin{figure}[t]
\begin{center}
\begin{tikzpicture}[scale=.9]
\draw[fill=blue!20!white, thick] (-2.5,-2.5) rectangle (2.5,2.5);
\foreach\i in {-2,...,2}{
\foreach\j in {-2,...,2}{
\draw[fill=gray,opacity=.3] (0.95*\i-0.03*\j-.5,0.06*\i+1.02*\j-.5) rectangle (0.95*\i-0.03*\j+.5,0.06*\i+1.02*\j+.5) ;
\draw (0.95*\i-0.03*\j-.5,0.06*\i+1.02*\j-.5) rectangle (0.95*\i-0.03*\j+.5,0.06*\i+1.02*\j+.5) ;
}}
\draw (0,0) node {$\times$};
\draw[->] (0,0) to (0.95,0.06);
\draw (.5,.1) node[below right] {$v_1$};
\draw[->] (0,0) to (-0.03,1.02);
\draw (.1,.5) node[above left] {$v_2$};
\draw (0,0) node[below left]{$0$};
\end{tikzpicture}
\caption[``Almost tiling'' of $B_{\ell R_0}$]{``Almost tiling'' of $B_{\ell R_0}$ by cubes $B(\sum_{i=1}^n k_i v_i, R_0)$, with $-\ell\le k_i\le\ell$.}\label{AlmTil}
\end{center}
\end{figure}

Let $\varep>0$, $\delta>0$, $\ell>0$ and $A\in GL_n(\R)$. We apply the second part of Remark~\ref{arithProg2} to obtain a parameter $R_0>0$ and a family $v_1,\cdots,v_n$ of $\varep$-translations of $\Gamma$ such that the family of cubes $\big(B(\sum_{i=1}^n k_i v_i, R_0)\big)_{-\ell\le k_i\le \ell}$ is an ``almost tiling'' of $B_{\ell R_0}$ (in particular, each $v_i$ is close to the vector having $2R_0$ in the $i$-th coordinate and $0$ in the others, see Figure~\ref{AlmTil}):
\begin{enumerate}[(1)]
\item\label{pati1} this collection of cubes fills almost all $B_{\ell R_0}$:
\[\frac{\card\Big(\Gamma \cap \big(\bigcup_{-\ell\le k_i\le \ell}B(\sum_{i=1}^n k_i v_i, R_0) \Delta B_{\ell R_0}\big)\Big)}{\card(\Gamma \cap B_{\ell R_0})} \le \varep;\]
\item\label{pati2} the overlaps of the cubes are not too big: for all collections $(k_i)$ and $(k'_i)$ such that $-\ell \le k_i,k'_i \le \ell$, 
\[\frac{\card\Big(\Gamma \cap \big(B(\sum_{i=1}^n k_i v_i, R_0)\Delta B(\sum_{i=1}^n k'_i v_i, R_0)\big)\Big)}{\card(\Gamma \cap B_{\ell R_0})} \le \varep;\]
\item\label{pati3} the vectors $\sum_{i=1}^n k_i v_i$ are translations for $\Gamma$: for every collection $(k_i)$ such that $-\ell \le k_i \le \ell$, 
\[\frac{\card\Big(\big(\Gamma \Delta (\Gamma-\sum_{i=1}^n k_i v_i)\big)\cap B_{R_0}\Big)}{\card(\Gamma \cap B_{R_0})} \le \varep.\]
\end{enumerate}

Increasing $R_0$ and $\ell$ if necessary, there exists $A'\in GL_n(\R)$ (respectively $SL_n(\R)$, $O_n(\R)$) such that $\|A-A'\|\le\delta$ and that we have
\begin{equation}\label{machinrioa}
\dist\left(\Leb_{\R^n/\Z^n}, \sum_{-\ell\le k_i\le\ell}\bar\delta_{A'(\sum_{i=1}^n k_i v_i)}\right) \le \varep.
\end{equation}
Indeed, if we denote by $\Lambda$ the lattice spanned by the vectors $v_1,\cdots,v_n$, then the set of matrices $A'$ such that $A'\Lambda$ is equidistributed modulo $\Z^n$ is dense in $GL_n(\R)$ (respectively $SL_n(\R)$ and $O_n(\R)$).

Then, we have,
\begin{align*}
\dist\Bigg(\Leb_{\R^n/\Z^n}, & \sum_{\substack{-\ell \le k_i \le \ell\\ x\in \Gamma \cap B(\sum_{i=1}^n k_i v_i, R_0)}} \bar\delta_{A'x}\Bigg) \le\\
  & \dist\Bigg(\Leb_{\R^n/\Z^n}, \sum_{\substack{-\ell \le k_i \le \ell\\ x\in \Gamma \cap B(0, R_0)}} \bar\delta_{A'(\sum_{i=1}^n k_i v_i)+ A'x}\Bigg)\\
  & + \dist\Bigg(\sum_{\substack{-\ell \le k_i \le \ell\\ x\in \Gamma \cap B(0, R_0)}} \bar\delta_{A'(\sum_{i=1}^n k_i v_i)+ A'x}, \sum_{\substack{-\ell \le k_i \le \ell\\ x\in \Gamma \cap B(\sum_{i=1}^n k_i v_i, R_0)}} \bar\delta_{A'x}\Bigg)
\end{align*}
By the property of convexity of $\dist$, the first term is smaller than
\[ \frac{1}{\card\big(\Gamma \cap B(0, R_0)\big)}\sum_{x\in \Gamma \cap B(0, R_0)}\dist\Bigg(\Leb_{\R^n/\Z^n}, \sum_{-\ell \le k_i \le \ell} \bar\delta_{A'(\sum_{i=1}^n k_i v_i)+ A'x}\Bigg);\]
by Equation~\eqref{machinrioa} and the fact that $\dist$ is invariant under translation, this term is smaller than $\varep$. As by hypothesis, the vectors $\sum_{i=1}^n k_i v_i$ are $\varep$-translations of $\Gamma$ (hypothesis \eqref{pati3}), the second term is also smaller than $\varep$. 
Thus, we get 
\[\dist\Bigg(\Leb_{\R^n/\Z^n}, \sum_{\substack{-\ell \le k_i \le \ell\\ x\in \Gamma \cap B(\sum_{i=1}^n k_i v_i, R_0)}} \bar\delta_{A'x}\Bigg) \le 2\varep\]

By the fact that the family of cubes $\big(B(\sum_{i=1}^n k_i v_i, R_0)\big)_{-\ell\le k_i\le \ell}$ is an almost tiling of $B_{\ell R_0}$ (hypotheses \eqref{pati1} and \eqref{pati2}), we get, for every $v\in \R^n$,
\[ \dist\left(\Leb_{\R^n/\Z^n}, \sum_{x\in \Gamma \cap B_{\ell R_0}} \bar\delta_{A'x}\right) < 4\varep.\]
Remark that we can suppose that this remains true on a whole neighbourhood of $A'$. We use the fact that $\Gamma$ is an almost periodic pattern to deduce that $A'$ belongs to the interior of $\mathcal U_\varep$.
\end{proof}

\section{Differences in almost periodic patterns}

From now, we suppose that the almost periodic patterns we consider are subsets of $\Z^n$. 

In the sequel, we will use the concept of difference in an almost periodic pattern. The following lemma, whose proof is straightforward, states that the occurrences of a given difference in an almost periodic pattern form an almost periodic pattern.

\begin{lemme}
Let $v\in\R^n$ and $\Gamma$ be an almost periodic pattern. Then the set
\[\{x\in\Gamma\mid x+v\in\Gamma\} = \Gamma\cap(\Gamma-v)\]
is an almost periodic pattern.
\end{lemme}

Then Proposition~\ref{limitexist} allows to associate to each translation of $\Z^n$ the frequency it appears in the almost periodic pattern.

\begin{definition}\label{DefDiff}
For $v\in\Z^n$, we set\index{$\rho_\Gamma$}
\[\rho_\Gamma(v) = \frac{D\{x\in\Gamma\mid x+v\in\Gamma\}}{D(\Gamma)} = \frac{D\big(\Gamma\cap(\Gamma-v)\big)}{D(\Gamma)} \in [0,1]\]
the \emph{frequency} of the difference $v$ in the almost periodic pattern $\Gamma$.
\end{definition}

The function $\rho_\Gamma$ is itself almost periodic in the sense given by H. Bohr (see \cite{MR1555192}).

\begin{definition}
Let $f : \Z^n\to \R$. Denoting by $T_v$\index{$T_v$} the translation of vector $v$, we say that $f$ is \emph{Bohr almost periodic} (also called \emph{uniformly almost periodic}) if for every $\varep>0$, the set
\[\mathcal N_\varep = \big\{v\in\Z^n\mid \|f - f\circ T_v \|_\infty<\varep \big\},\]
is relatively dense.

If $f: \Z^n\to \R$ is a Bohr almost periodic function, then it possesses a \emph{mean}\index{$\mathcal M$} $\mathcal M(f)$ (see for example the historical paper of H. Bohr \cite[Satz VIII]{MR1555192}), which satisfies: for every $\varep>0$, there exists $R_0>0$ such that for every $R\ge R_0$ and every $x\in\R^n$, we have
\[\left|\mathcal M(f) - \frac{1}{\card[B(x,R)]}\sum_{v\in [B(x,R)]} f(v) \right| <\varep.\]
\end{definition}

The fact that $\rho_\Gamma$ is Bohr almost periodic is straightforward.

\begin{lemme}
If $\Gamma$ is an almost periodic pattern, then the function $\rho_\Gamma$ is Bohr almost periodic.
\end{lemme}

In fact, we can compute precisely the mean of $\rho(\Gamma)$.

\begin{prop}\label{IntRho}
If $\Gamma$ is an almost periodic pattern, then we have 
\[\mathcal M (\rho_\Gamma) = D(\Gamma).\]
\end{prop}

\begin{proof}[Proof of Proposition \ref{IntRho}]
This proof lies primarily in an inversion of limits.

Let $\varep>0$. As $\Gamma$ is an almost periodic pattern, there exists $R_0>0$ such that for every $R\ge R_0$ and every $x\in\R^n$, we have
\begin{equation}\label{eqDens}
\left|D(\Gamma) - \frac{\Gamma \cap [B(x,R)]}{\card[B_R]}\right|\le \varep.
\end{equation}

So, we choose $R\ge R_0$, $x\in\Z^n$ and compute
\begin{align*}
\frac{1}{\card[B_R]} & \sum_{v\in[B(x,R)]} \rho_\Gamma(v) = \frac{1}{\card[B_R]}\sum_{v\in[B(x,R)]} \frac{D\big((\Gamma-v)\cap \Gamma\big)}{D(\Gamma)}\\
       = & \frac{1}{\card[B_R]}\sum_{v\in[B(x,R)]} \lim_{R'\to +\infty}\frac{1}{\card[B_{R'}]}\sum_{y\in[B_{R'}]} \frac{\1_{y\in\Gamma-v} \1_{y\in\Gamma}}{D(\Gamma)}\\
       = & \frac{1}{D(\Gamma)}\lim_{R'\to +\infty}\frac{1}{\card[B_{R'}]} \sum_{y\in[B_{R'}]} \1_{y\in\Gamma}\frac{1}{\card[B_R]}\sum_{v\in[B(x,R)]} \1_{y\in\Gamma-v}\\
       = & \frac{1}{D(\Gamma)}\underbrace{\lim_{R'\to +\infty}\frac{1}{\card[B_{R'}]} \sum_{y\in[B_{R'}]} \1_{y\in\Gamma}}_{\text{first term}}\underbrace{\frac{1}{\card[B_R]}\sum_{v'\in[B(y+x,R)]} \1_{v'\in\Gamma}}_{\text{second term}}.
\end{align*}
By Equation \eqref{eqDens}, the second term is $\varep$-close to $D(\Gamma)$. Considered independently, the first term is equal to $D(\Gamma)$ (still by Equation \eqref{eqDens}). Thus, we have
\[\left|\frac{1}{\card[B(x,R)]} \sum_{v\in[B(x,R)]} \rho_\Gamma(v) - D(\Gamma)\right|\le \varep,\]
that we wanted to prove.
\end{proof}

%
%

We now state a Minkowski-type theorem for the map $\rho_\Gamma$. To begin with, we recall the classical Minkowski theorem (see for example \cite[IX.3]{berger2009géométrie} or the whole book \cite{MR893813}).

\begin{theoreme}[Minkowski]\label{Minkowski}
Let $\Lambda$ be a lattice of $\R^n$, $k\in\N$ and $S\subset\R^n$ be a centrally symmetric convex body. If $\Leb(S/2) > k \operatorname{covol}(\Lambda)$, then $S$ contains at least $2k$ distinct points of $\Lambda\setminus\{0\}$.
\end{theoreme}

In particular, if $\Leb(S/2) > \operatorname{covol}(\Lambda)$, then $S$ contains at least one point of $\Lambda\setminus\{0\}$. This theorem is optimal in the following sense: for every lattice $\Lambda$ and for every $k\in\N$, there exists a centrally symmetric convex body $S$ such that $\Leb(S/2) = k\operatorname{covol}(\Lambda)$ and that $S$ contains less than $2k$ distinct points of $\Lambda\setminus\{0\}$. In Chapter~\ref{Souris}, we will make use of Haj\'os theorem (Theorem~\ref{hajos}), which specify what can happen in the case where $S$ is an infinite ball satisfying $\Leb(S/2) = \operatorname{covol}(\Lambda)$.

\begin{proof}[Proof of Theorem \ref{Minkowski}]
We consider the function
\[\varphi = \sum_{\lambda\in\Lambda} \1_{\lambda + S/2}.\]
The hypothesis about the covolume of $\Lambda$ and the volume of $S/2$ imply that the mean of the periodic function $\varphi$ satisfies $\mathcal M(\varphi)>k$. In particular, there exists $x_0\in \R^n$ such that $\varphi(x_0)\ge k+1$ (note that this argument is similar to the pigeonhole principle). Then, there exists $\lambda_0,\cdots,\lambda_k\in \Lambda$, with the $\lambda_i$ sorted in lexicographical order (for a chosen basis), such that the points $x_0-\lambda_i$ belong to $S/2$. As $S/2$ is centrally symmetric, $\lambda_i-x_0$ belongs to $S/2$ and as $S/2$ is convex, $\big((x_0-\lambda_0) + (\lambda_i-x_0)\big)/2 = (\lambda_i - \lambda_0)/2$ also belongs to $S/2$. Then, $\lambda_i-\lambda_0\in\Lambda\setminus\{0\} \cap S$ for every $i\in \llbracket 1,k\rrbracket$. By hypothesis, these $k$ vectors are all different. To obtain $2k$ different points of $S\cap \Lambda\setminus\{0\}$ (instead of $k$ different points), it suffices to consider the points $\lambda_0-\lambda_i$; this collection is disjoint from the 
collection of $\lambda_i-\lambda_0$ by 
the fact that the $\lambda_i$ are sorted in lexicographical order. This proves the theorem.
\end{proof}

With the definitions we introduced, Minkowski theorem can be seen as a result about the function $\rho_\Gamma$: \emph{Let $\Lambda$ be a lattice of $\R^n$, and $S\subset\R^n$ be a centrally symmetric convex body, then}
\[\sum_{u\in S} \rho_\Lambda(u) \ge 2\lceil D(\Lambda) \Leb(S/2)\rceil-1.\]
We now state a similar statement in the more general case of almost periodic patterns.

\begin{theoreme}[in collaboration with \'E. Joly]\label{MinkAlm}
Let $\Gamma\subset \Z^n$ be an almost periodic pattern, and $S\subset\R^n$ be a centrally symmetric convex body. Then (recall that $[S] = S\cap \Z^n$)
\[\sum_{u\in [S]} \rho_\Gamma(u) \ge D(\Gamma) \card[S/2].\]
\end{theoreme}

Remark that by Theorem~\ref{Minkowski}, we have $\card[S/2] \ge 2\lceil \Leb(S/4)\rceil-1$. Thus,
\[\sum_{u\in [S]} \rho_\Gamma(u) \ge D(\Gamma) \big(2\lceil \Leb(S/4)\rceil-1\big).\]

\begin{rem}
It is possible to prove in an easier way a weak version of Theorem~\ref{MinkAlm} which is sufficient for the use we will make (Theorem~\ref{AnswerConjIsom}). However, it seemed nicer to us to state a more optimal property.
\end{rem}

The case of equality in the theorem is attained even in the non trivial case where $\card[S/2]>1$, as shown by the following example.

\begin{figure}[t]
\begin{center}
\begin{tikzpicture}[scale=.75]
\draw (0,0) node[right]{$0$};
\foreach\i in {-1,...,1}{
\foreach\j in {-2,...,2}{
\draw[color=red!70!black] (3*\i,\j) node {$\bullet$};
}}
\foreach\i in {-3,...,3}{
\foreach\j in {-2,...,2}{
\draw[color=black] (\i,\j) node {$\cdot$};
}}
\draw[color=green!70!black,thick] (-2.2,.1) -- (-2.2,-1.2) -- (-1,-1.2) -- (2.2,-.1) -- (2.2,1.2) -- (1,1.2) -- cycle;
\draw[color=green!40!black] (1,1.1) node[above right] {$S$};
\draw[color=blue!70!black] (-1.1,.05) -- (-1.1,-.6) -- (-.5,-.6) -- (1.1,-.05) -- (1.1,.6) -- (.5,.6) -- cycle;
\draw[color=blue!40!black] (1,.5) node[right] {\small$S/2$};
\draw[color=red!60!black] (3,.3) node[right] {$\Gamma$};
\end{tikzpicture}
\caption{Example \ref{OptMink} of equality case in Theorem \ref{MinkAlm} for $k=3$}\label{FigExMink}
\end{center}
\end{figure}
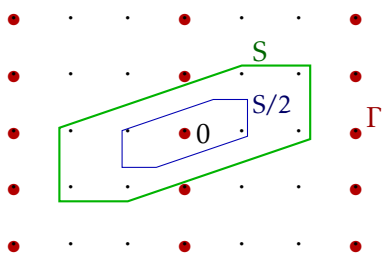

\begin{ex}\label{OptMink}
If $k$ is an odd number, if $\Gamma$ is the lattice $k\Z\times \Z$, and if $S$ is a centrally symmetric convex set such that (see Figure~\ref{FigExMink})
\[S\cap \Z^2 = \{(i,0)\mid i\in\llbracket -(k-1), k-1\rrbracket\} \cup \{\pm(i,1)\mid i\in\llbracket 1, k-1\rrbracket\} ,\]
then $\sum_{u\in S} \rho(u) = 1$, $D(\Gamma) = 1/k$ and $\card(S/2 \cap \Z^n) = k$.
\end{ex}

The strategy of proof of this theorem is similar to that of the classical Minkowski theorem: we consider the set $\Gamma + S/2$ and define a suitable auxiliary function based on this set.

\begin{proof}[Proof of Theorem \ref{MinkAlm}]
We first define 
\[\rho_a^R = \frac{1}{\card[B_R]}\sum_{v\in [B_R]} \1_{v\in[S/2+a]} \sum_{u\in [S]}\frac{\1_{v\in\Gamma} \1_{u+v\in\Gamma}}{D(\Gamma)}.\]
We use this function to apply an argument of double counting. Firstly, we have
\begin{align*}
\sum_{a\in\Z^n} \rho_a^R & = \frac{1}{\card[B_R]}\sum_{v\in [B_R]} \sum_{u\in [S]}\frac{\1_{v\in\Gamma} \1_{u+v\in\Gamma}}{D(\Gamma)}\sum_{a\in\Z^n}\1_{a\in[S/2+v]}\\
    & = \frac{1}{\card[B_R]}\sum_{v\in [B_R]} \sum_{u\in [S]}\frac{\1_{v\in\Gamma} \1_{u+v\in\Gamma}}{D(\Gamma)}\card[S/2+v]\\
    & = \card[S/2]\sum_{u\in [S]} \frac{1}{\card[B_R]}\sum_{v\in [B_R]} \frac{\1_{v\in\Gamma} \1_{u+v\in\Gamma}}{D(\Gamma)}.
\end{align*}
Thus, by the definition of $\rho_\Gamma$, we get 
\begin{equation}\label{res1}
\lim_{R\to+\infty} \sum_{a\in\Z^n} \rho_a^R\, =\, \card[S/2]\sum_{u\in [S]} \rho_\Gamma(u).
\end{equation}
The conclusion of the theorem is an estimate on the right side of this equality; to prove it we compute the left side in another way.

First of all, we remark that as $S$ is a centrally symmetric convex body, $v,w\in S/2$ implies that $w-v\in S$, so (applying this property to $w = u+v$)
\[\1_{v\in S/2+a} \1_{u\in S}  \ge \1_{v\in S/2+a} \1_{u+v\in S/2+a},\]
and thus
\[\1_{v\in[S/2+a]} \1_{u\in [S]} \1_{v\in\Gamma} \1_{u+v\in\Gamma} \ge \1_{v\in\Gamma\cap[S/2+a]} \1_{u+v\in\Gamma\cap[S/2+a]}.\]
We now sum this inequality over $u\in\Z^n$:
\begin{align*}
\sum_{u\in [S]} \1_{v\in[S/2+a]} \1_{v\in\Gamma} \1_{u+v\in\Gamma} & \ge \1_{v\in\Gamma\cap[S/2+a]} \sum_{u\in\Z^n}\1_{u+v\in\Gamma\cap[S/2+a]}\\
     & \ge \1_{v\in\Gamma\cap[S/2+a]} \sum_{u'\in\Z^n}\1_{u'\in\Gamma\cap[S/2+a]}\\
     & \ge \1_{v\in\Gamma\cap[S/2+a]} \card\big(\Gamma\cap[S/2+a]\big);
\end{align*}
so
\[\rho_a^R \ge \frac{1}{D(\Gamma)}\frac{1}{\card[B_R]}\sum_{v\in[B_R]}\1_{v\in\Gamma\cap[S/2+a]} \card(\Gamma\cap[S/2+a]).\]
We denote by $B_R^S$ the $S$-interior of $B_R$, that is
\[B_R^S = \big(B_R^\complement + S\big)^\complement = \{x\in B_R\mid\forall s\in S, x+s\in B_R\}.\]
In particular, $a\in[B_R^S]$ implies that $[S/2+a]\subset [B_R]$, then
\begin{align*}
\sum_{a\in \Z^n} \rho_a^R & \ge \frac{1}{D(\Gamma)} \sum_{a\in \Z^n} \frac{1}{\card[B_R]}\sum_{v\in[B_R]}\1_{v\in\Gamma\cap[S/2+a]} \card\big(\Gamma\cap[S/2+a]\big)\\
    & \ge \frac{1}{D(\Gamma)}\frac{1}{\card[B_R]}\sum_{a\in [B_R^S]}\sum_{v\in[B_R]} \1_{v\in\Gamma\cap[S/2+a]} \card\big(\Gamma\cap[S/2+a]\big)\\
		& \ge \frac{1}{D(\Gamma)}\frac{1}{\card[B_R]}\sum_{a\in [B_R^S]} \card\big(\Gamma\cap[S/2+a]\big)^2.
\end{align*}
We then use the fact that the family $\{B_R\}_{R>0}$ is van Hove when $R$ goes to infinity (see for example \cite[Equation 4]{MR1884143}), that is
\[\lim_{R\to +\infty} \frac{\card[B_R]-\card[B_R^S]}{\card[B_R]} = 0;\]
thus,
\[\lim_{R\to+\infty} \sum_{a\in\Z^n} \rho_a^R \ge \frac{1}{D(\Gamma)}\underset{R\to +\infty}{\overline\lim}\ \frac{1}{\card[B_R]}\sum_{a\in [B_R]} \card\big(\Gamma\cap[S/2+a]\big)^2.\]
Using the convexity of $x\mapsto x^2$, we 
deduce that
\begin{equation}\label{res2}
\lim_{R\to+\infty} \sum_{a\in\Z^n} \rho_a^R \ge \underset{R\to +\infty}{\overline\lim}\ \frac{1}{D(\Gamma)}\left(\frac{1}{\card[B_R]}\sum_{a\in [B_R]} \card\big(\Gamma\cap[S/2+a]\big)\right)^2.
\end{equation}

Now, it remains to compute
\begin{align*}
\frac{1}{\card[B_R]}\sum_{a\in [B_R]} \card\big(\Gamma\cap[S/2+a]\big) & = \frac{1}{\card[B_R]}\sum_{a\in [B_R]}\sum_{s\in [S/2]} \1_{s+a\in\Gamma}\\
      & = \sum_{s\in [S/2]}\frac{1}{\card[B_R]}\sum_{a'\in [B_R]-s} \1_{a'\in\Gamma}.
\end{align*}
But for every $s\in\R^n$, we have
\[\lim_{R\to+\infty} \frac{1}{\card[B_R]}\sum_{a'\in [B_R]-s} \1_{a'\in\Gamma} = D(\Gamma),\]
thus
\[\lim_{R\to+\infty}\frac{1}{\card[B_R]}\sum_{a\in [B_R]} \card\big(\Gamma\cap[S/2+a]\big) = \card[S/2] D(\Gamma).\]
Applied to Equation \eqref{res2}, this gives 
\[\lim_{R\to+\infty} \sum_{a\in\Z^n} \rho_a^R \ge \card[S/2]^2 D(\Gamma).\]
To finish the proof, we combine it with the first estimate in Equation \eqref{res1} and get:
\[\card[S/2]\sum_{u\in S} \rho_\Gamma(u) \ge \card[S/2]^2 D(\Gamma),\]
so
\[\sum_{u\in S} \rho_\Gamma(u) \ge \card[S/2] D(\Gamma),\]
\end{proof}

\nocite{minkowski1910geometrie}

\section{Model sets}\label{SecModel}

We now turn to a more precise notion about almost periodicity: model sets. We begin this section by motivating the introduction of these sets: we give an alternative construction of the sets $\widehat A(\Z^n)$ in terms of model sets.

A point $x\in\Z^n$ belongs to $\widehat A(\Z^n)$ if and only if there exists $y\in\Z^n$ such that $\|x-Ay\|_\infty <\frac 12$. In other words, if we set
\[ M_A = \left(\begin{array}{cc}
A & -\operatorname{Id}\\
  & \operatorname{Id}
\end{array}\right)\in M_{2n}(\R),\]
if we note $p_1$ and $p_2$ the projections of $\R^{2n}$ on respectively the $n$ firsts and the $n$ last coordinates, and if we set $W = ]-\frac 12,\frac 12]^n$, then
\[\widehat A(\Z^n) = \big\{p_2(M_A v)\mid v\in\Z^{2n},\, p_1(M_A v)\in W\big\}.\]
This notion is close to that of model set introduced by Y. Meyer in the early seventies \cite{MR0485769}, but in our case the projection $p_2$ is not injective. Model sets are sometimes called ``cut and project'' sets in the literature.

\begin{figure}[t]
\begin{center}
\begin{tikzpicture}[scale=1.4]
\fill[color=blue!10!white] (-.6,-2) rectangle (.9,2);
\draw[color=blue!80!black] (-.6,-2) -- (-.6,2);
\draw[color=blue!80!black] (.9,-2) -- (.9,2);
\draw[color=blue!80!black, thick] (-.6,0) -- (.9,0);
\draw[color=blue!80!black] (.25,0) node[below] {$W$};
\draw (-3,0) -- (3,0);
\draw (0,-2) -- (0,2);
\clip (-3,-2) rectangle (3,2);

\draw (.866,.364) -- (0,.364);
\draw (-.129,.987) -- (0,.987);
\draw (.737,1.351) -- (0,1.351);
\draw (-.258,1.974) -- (0,1.974);
\draw (.129,-.987) -- (0,-.987);
\draw (.258,-1.974) -- (0,-1.974);

\draw (0,0) node {$\times$};
\draw (0,.364) node {$\times$};
\draw (0,.987) node {$\times$};
\draw (0,1.351) node {$\times$};
\draw (0,1.974) node {$\times$};
\draw (0,-.987) node {$\times$};
\draw (0,-1.974) node {$\times$};

\foreach\i in {-3,...,3}{
\foreach\j in {-3,...,3}{
\draw[color=green!40!black] (.866*\i-.129*\j,.364*\i+.987*\j) node {$\bullet$};
}}
\draw[color=green!40!black] (1.2,-.8) node {$\Lambda$};
\end{tikzpicture}
\caption[Model set]{Construction of a model set.}\label{FigModel}
\end{center}
\end{figure}
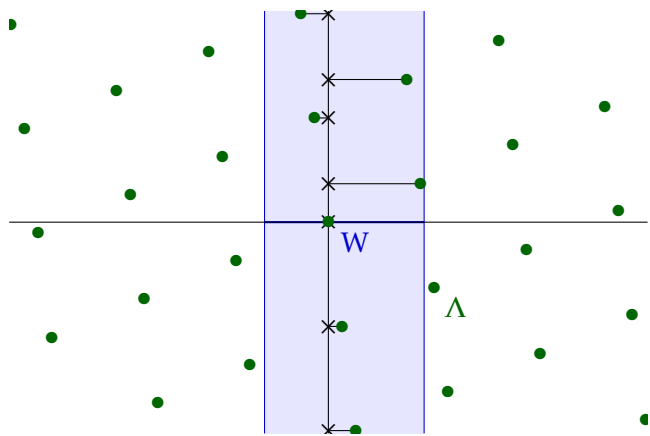

\begin{definition}\label{DefModel}
Let $\Lambda$ be a lattice of $\R^{m+n}$, $p_1$ and $p_2$ the projections of $\R^{m+n}$ on respectively $\R^m\times \{0\}_{\R^n}$ and $\{0\}_{\R^m} \times \R^n$, and $W$ a Riemann integrable subset of $\R^m$. The \emph{model set} modelled on the lattice $\Lambda$ and the \emph{window} $W$ is (see Figure~\ref{FigModel})
\[\Gamma = \big\{ p_2(\lambda)\mid \lambda\in\Lambda,\, p_1(\lambda)\in W \big\}.\]
\end{definition}

Here, we will call \emph{model set} every set of this type, even if the projection $p_2$ is not injective. Indeed, this phenomenon is what is interesting for us, if it did not occur there would be no loss of injectivity when applying discretizations of linear maps. Notice that this definition, which could seem very restrictive for the set $\Gamma$, is in fact quite general: as stated by Y. Meyer in \cite{MR0485769}, every Meyer set\footnote{A set $\Gamma$ is a \emph{Meyer set} if $\Gamma-\Gamma$ is a Delone set. It is equivalent to ask that there exists a finite set $F$ such that $\Gamma-\Gamma \subset \Gamma + F$ (see \cite{MR1400744}).} is a subset of a model set. Conversely, model sets are Meyer sets (see \cite{MR1420415}).

Returning to our problem of images of the lattice $\Z^n$ by discretizations of linear maps, we have the following (trivial) result.

\begin{prop}\label{ImgModel}
Let $A_1,\cdots,A_k \in GL_n(\R)$ be $k$ invertible maps, then the $k$-th image $\Gamma_k = \widehat{A_k}\circ\dots\circ\widehat{A_1} (\Z^n)$ of $\Z^n$ by the discretizations is the model set modelled on the window $W = ]-\frac 12,\frac 12]^{nk}$ and the lattice $M_{A_1,\cdots,A_k}\Z^{n(k+1)}$, where
\[ M_{A_1,\cdots,A_k} = \left(\begin{array}{ccccc}
A_1 & -\Id &        &        & \\
    & A_2  & -\Id   &        & \\
    &      & \ddots & \ddots & \\
    &      &        & A_k    & -\Id\\
    &      &        &        & \Id
\end{array}\right)\in M_{n(k+1)}(\R).\]
\end{prop}

\begin{rem}
This notion has the advantage that it builds the $k$-th image directly: the concept of time disappears, so we will be able ``anticipate''\ the behaviour of successive images. The downside is the increasing of the dimension; thus it will be more difficult to have a geometric intuition\dots
\end{rem}

In the sequel, we will only consider model sets whose window is regular.

\begin{definition}
Let $W$ be a subset of $\R^n$. We say that $W$ is \emph{regular} if for every affine subspace $V\subset\R^n$, we have
\[\Leb_V\Big( B_V\big(\partial (V\cap W),\eta\big)\Big) \underset{\eta\to 0}{\longrightarrow} 0,\]
where $\Leb_V$ denotes the Lebesgue measure on $V$, and $B_V\big(\partial (V\cap W),\eta\big)$ the set of points of $V$ whose distance to $\partial (V\cap W)$ is smaller than $\eta$ (of course, the boundary is also take in restriction to $V$).
\end{definition}

The link with the previous sections is made by the following theorem.

\begin{theoreme}[in collaboration with Y. Meyer]\label{ModelAlmost}
A model set modelled on a regular window is an almost periodic pattern.
\end{theoreme}

In other words, for every $\varep>0$, there exists $R_0>0$ and a relatively dense set $\mathcal N$ such that for every $v\in\mathcal N$ and every $R\ge R_0$, most of the points (\emph{i.e.} a proportion greater than $\varep$) of the model set $\Gamma$ also belong to $v+\Gamma$ (see Definition \ref{DefAlmPer}).

We begin by proving a weak version of this theorem.

\begin{lemme}\label{LemModelAlmost}
A model set modelled on a window with nonempty interior is relatively dense.
\end{lemme}

\begin{proof}[Proof of Lemma \ref{LemModelAlmost}]
We prove this lemma in the specific case where the window is $B_\eta$ (recall that $B_\eta$ is the infinite ball of radius $\eta$ centred at 0). We will use this lemma only in this case (and the general case can be treated the same way).

Let $\Gamma$ be a model set modelled on a lattice $\Lambda$ and a window $B_\eta$. We will use the fact that for any centrally symmetric convex set $S\subset \R^n$, if there exists a basis $e_1,\cdots,e_n$ of $\Lambda$ such that for each $i$, $\lceil n/2 \rceil e_i \in S$, then $S$ contains a fundamental domain of $\R^n/\Lambda$, that is to say, for every $v\in\R^n$, we have $(S+v)\cap \Lambda\neq\emptyset$. This is due to the fact that the parallelepiped spanned by the vectors $e_i$ is included into the simplex spanned by the vectors $\lceil n/2\rceil e_i$.

We set
\[ V = \bigcap_{\eta'>0} \Vect\big(p_1^{-1}(B_{\eta'})\cap \Lambda\big) = \bigcap_{\eta'>0}\Vect\big\{\lambda\in\Lambda\mid d_\infty(\lambda,\ker p_1)\le\eta'\big\},\]
and remark that $\im p_2 = \ker p_1\subset V$, simply because for every vectorial line $D\subset \R^n$ (and in particular for $D\subset \ker p_1$), there exists some points of $\Lambda\setminus\{0\}$ arbitrarily close to $D$. We also take $R>0$ such that
\[V\subset V'\doteq\Vect\big(p_1^{-1}(B_{\eta/\lceil n/2\rceil})\cap \Lambda \cap p_2^{-1}(B_R)\big).\]
We then use the remark made in the beginning of this proof and apply it to the linear space $V'$, the set $S = \big(p_2^{-1}(B_R)\times p_1^{-1}(B_\eta)\big)\cap V'$, and the module $V'\cap \Lambda$. This leads to
\[\forall v\in V,\ \Big(\big(p_1^{-1}(B_\eta) \cap p_2^{-1}(B_R)\big) + v \Big) \cap \Lambda \neq \emptyset,\]
and as $\im p_2 \subset V$, we get
\[\forall v'\in \R^n,\ \big(p_1^{-1}(B_\eta) \cap p_2^{-1}(B_R+v')\big) \cap \Lambda \neq \emptyset;\]
this proves that the model set is relatively dense for the radius $R$.
\end{proof}

\begin{proof}[Proof of Theorem \ref{ModelAlmost}]
Let $\Gamma$ be a model set modelled on a lattice $\Lambda$ and a window $W$.

First of all, we decompose $\Lambda$ into three supplementary modules: $\Lambda = \Lambda_1 \oplus \Lambda_2 \oplus \Lambda_3$, such that (see \cite[Chap. VII, \S 1, 2]{MR1726872}):
\begin{enumerate}
\item $\Lambda_1 = \ker p_1 \cap \Lambda$;
\item $p_1(\Lambda_2)$ is discrete;
\item $p_1(\Lambda_3)$ is dense in the vector space $V$ it spans (and such a vector space is unique), and $V\cap p_1(\Lambda_2) = \{0\}$.
\end{enumerate}
As $\Lambda_1 = \ker p_1 \cap \Lambda = \im p_2 \cap \Lambda$, we have $\Lambda_1 = p_2(\Lambda_1)$. Thus, for every $\lambda_1\in\Lambda_1$ and every $\gamma\in\Gamma$, we have $\lambda_1+\gamma\in\Gamma$. So $\Lambda_1$ is a set of periods for $\Gamma$. Therefore, considering the quotients $\R^n/\Vect\Lambda_1$ and $\Lambda/\Lambda_1$ if necessary, we can suppose that ${p_1}_{|\Lambda}$ is injective (in other words, $\Lambda_1 = \{0\}$).

Under this assumption, the set $p_2(\Lambda_3)$ spans the whole space $\im p_2$. Indeed, as $\ker p_1\cap\Lambda = \{0\}$, we have the decomposition
\begin{equation}\label{eqAgr}
\R^{m+n} = \underbrace{\ker p_1}_{=\im p_2} \oplus \underbrace{\Vect\big(p_1(\Lambda_2)\big) \oplus \Vect\big(p_1(\Lambda_3)\big)}_{=\im p_1}.
\end{equation}
Remark that as $p_1(\Lambda_2)$ is discrete, we have $\dim\Vect\big(p_1(\Lambda_2)\big) = \dim\Lambda_2$; thus, considering the dimensions in the decomposition \eqref{eqAgr}, we get
\begin{equation}\label{EstimDimpt}
\dim\Vect(\Lambda_3) = \dim \Big(\ker p_1 \oplus \Vect\big(p_1(\Lambda_3)\big)\Big).
\end{equation}
The following matrix represents a basis of $\Lambda = \Lambda_2\oplus\Lambda_3$ in a basis adapted to the decomposition \eqref{eqAgr}.
\[\begin{array}{cc}
&
\begin{array}{cc}
\Lambda_2 &  \Lambda_3\\
\overbrace{\phantom{\hspace{1cm}}} & \overbrace{\phantom{\hspace{1cm}}}\end{array}\\
\begin{array}{r}
\ker p_1 = \im p_2           \Big\{\!\!\!\!\!\!\\
\Vect\big(p_1(\Lambda_2)\big)\Big\{\!\!\!\!\!\!\\
\Vect\big(p_1(\Lambda_3)\big)\Big\{\!\!\!\!\!\!
\end{array}
& \renewcommand{\arraystretch}{1.5}
\left(\begin{array}{c|c}
* & *\\ \hline
* & 0\\ \hline
\hspace{.4cm}0 \hspace{.4cm}{} & \hspace{.4cm} *\hspace{.4cm}{}
\end{array}\right)
\renewcommand{\arraystretch}{1}
\end{array}\]

We can see that the projection of the basis of $\Lambda_3$ on $\im p_2 \oplus \Vect\big(p_1(\Lambda_3)\big)$ form a free family; by Equation~\eqref{EstimDimpt}, this is in fact a basis. Thus, $\Vect(\Lambda_3)\supset \ker p_1 = \im p_2$, so $\Vect\big(p_2(\Lambda_3)\big) = \im (p_2)$.
\bigskip

For $\eta>0$, let $\mathcal N(\eta)$ be the model set modelled on $\Lambda$ and $B(0,\eta)$, that is
\[\mathcal N(\eta) = \{p_2(\lambda_3)\mid \lambda_3 \in \Lambda_3,\,\|p_1(\lambda_3)\|_\infty\le \eta\}.\]
Lemma \ref{LemModelAlmost} asserts that $\mathcal N(\eta)$ is relatively dense in the space it spans, and the previous paragraph asserts that this space is $\im p_2$. The next lemma, which obviously implies Theorem \ref{ModelAlmost}, expresses that if $\eta$ is small enough, then $\mathcal N(\eta)$ is the set of translations we look for.

\begin{lemme}\label{Lem2ModelAlmost}
For every $\varep>0$, there exists $\eta>0$ and a regular model set $Q(\eta)$ such that $D^+(Q(\eta))\leq \varep$ and
\[v\in \mathcal N(\eta)\Rightarrow (\Gamma+v)\Delta \Gamma\subset Q(\eta).\]
\end{lemme}

We have now reduced the proof of Theorem~\ref{ModelAlmost} to that of Lemma \ref{Lem2ModelAlmost}.
\end{proof}

\begin{proof}[Proof of Lemma \ref{Lem2ModelAlmost}]
We begin by proving that $(\Gamma+v)\setminus\Gamma\subset Q(\eta)$ when $v\in \mathcal N(\eta)$. As $v\in \mathcal N(\eta)$, there exists $\lambda_0\in \Lambda_3$ such that $p_2(\lambda_0)=v$ and $\|p_1(\lambda_0)\|_\infty\le \eta$.

If $x\in \Gamma+v$, then $x=p_2(\lambda_2+\lambda_3)+p_2(\lambda_0)=p_2(\lambda_2+\lambda_3+\lambda_0)$ where $\lambda_2 \in \Lambda_2$, $\lambda_3 \in \Lambda_3$ and $p_1(\lambda_2+\lambda_3)\in W$. If moreover $x\notin \Gamma$, it implies that $p_1(\lambda_2+\lambda_3+\lambda_0)\notin W$. Thus, $p_1(\lambda_2+\lambda_3+\lambda_0)\in W_\eta$, where (recall that $V = \Vect (p_1(\Lambda_3))$)
\[W_\eta = \big\{k+w\mid k\in\partial W, w\in V\cap B_\eta\big\}.\]
We have proved that $\Gamma+v\setminus\Gamma\subset Q(\eta)$, where
\[Q(\eta) =\big\{p_2(\lambda)\mid \lambda\in \Lambda,\,p_1(\lambda)\in W_\eta\big\}.\]
Let us stress that the model set $Q(\eta)$ does not depend on $v$. We now observe that as $W$ is regular, we have
\[\sum_{\lambda_2\in\Lambda_2} \Leb_{V+p_1(\lambda_2)}\big(W_\eta\cap (V+p_1(\lambda_2))\big)\underset{\eta\to 0}{\longrightarrow}0.\]
As $p_1(\Lambda_3)$ is dense in $V$ (thus, it is equidistributed), the uniform upper density of the model set $Q(\eta)$ defined by the window $W_\eta$ can be made smaller than $\varep$ by taking $\eta$ small enough.

The treatment of $\Gamma\setminus(\Gamma+v)$ is similar; this ends the proof of Lemma~\ref{Lem2ModelAlmost}.
\end{proof}

\chapter{Rate of injectivity of linear maps}\label{Souris}

In this chapter, we focus in more detail on the rate of injectivity of a sequence of linear maps (see Definition~\ref{DefTaux}).

To begin with, study the easier case of the rate of injectivity in time 1. First of all, we study the continuity of the map $\tau = \tau^1$. Unfortunately, this map is not continuous (Proposition~\ref{rateEx}). However, $\tau$ is continuous at every totally irrational matrix\footnote{We say that a matrix $A\in GL_n(\R)$ is \emph{totally irrational} when the set $A(\Z^n)$ is equidistributed modulo $\Z^n$.} (Proposition~\ref{ThMeanRate}). More precisely, we define the \emph{mean rate of injectivity} $\overline\tau(A)$ of a matrix $A$: it is the mean of the rates of injectivity of the affine maps $A+v$, for $v\in\R^n/\Z^n$ (see Definition~\ref{DefMeanRate}). It turns out that this quantity is much more convenient to use: by an argument of equirepartition, it coincides with the rate of injectivity when $A$ is totally irrational, and it is continuous on the whole set $GL_n(\R)$ (Proposition~\ref{ThMeanRate}). We finish this study by estimating the lack of continuity of $\tau$ -- that is, the difference between $\overline\tau$ and $\tau$ -- at the matrices which are not totally irrational: this difference is arbitrarily small out of a topologically small set, that is, out of a locally finite union of hyperplanes (Proposition~\ref{oscill}).

These considerations allow us to focus on the mean rate of injectivity. The equidistribution property we gain by using the mean rate of injectivity allows to have a formula to compute it directly: for $A\in GL_n(\R)$, we have (Proposition~\ref{FormTau1}, see also Proposition~\ref{FormTau2} and Figures~\ref{tourp} and \ref{tourper})
\[\overline\tau(A) = \det(A) D\left(\bigcup_{\lambda\in A\Z^n} B_\infty(\lambda,1/2)\right).\]
As the set $A\Z^n$ is a lattice, this quantity is equal to the area of the intersection between a fundamental domain of $A\Z^n$ and the union of cubes $\bigcup_{\lambda\in A\Z^n} B_\infty(\lambda,1/2)$: this reduces the computation of the mean rate of injectivity to that of the area of a finite union of cubes. This formula proves that the mean rate of injectivity is continuous and piecewise polynomial (Corollary~\ref{superCont}). It also allows to compute the mean rate of injectivity of some practical examples (see Applications~\ref{AppRot} and \ref{AppLattPara}), and gives a quick method to compute numerically the mean rate of injectivity (see Figure~\ref{FigGradMean}).

To characterize the matrices for which the mean rate of injectivity is equal to 1, we combine this formula with the classical Haj\'os theorem (Theorem~\ref{hajos}) and get directly Corollary~\ref{CoroHajos}: a matrix $A\in SL_n(\R)$ has a mean rate of injectivity equal to 1 if and only if there exists $B\in SL_n(\Z)$ such that in a canonical basis of $\R^n$ (that is, permuting coordinates if necessary), the matrix $AB$ is upper triangular with ones on the diagonal. Thus, the set of such matrices forms a locally finite union of manifolds of positive codimension.
\bigskip

We then study the behaviour of the asymptotic rate of injectivity of a generic sequence of matrices of $SL_n(\R)$. It is given by the main theorem of this chapter (Theorem~\ref{ConjIntro}, see also Theorem~\ref{ConjPrincip}): for a generic sequence $(A_k)_{k\ge 1}$ of matrices of $SL_n(\R)$, we have $\tau^\infty\big((A_k)_k\big) = 0$.

The end of the chapter is devoted to the study of this theorem and some variations on it. The first variant is given in the (easy) case where all the matrices are diagonal (Proposition~\ref{CasLin}): for a generic sequence $(A_k)_{k\ge 1}$ of diagonal matrices, we have $\tau^\infty((A_k)_k) = 0$.

In the general case, we give a proof of the following weak version of Theorem~\ref{ConjPrincip} (Theorem~\ref{PerLin1}): if $(A_k)_{k\ge 1}$ is a generic sequence of matrices of $SL_n(\R)$, then $\tau^\infty((A_k)_k)\le 1/2$. This proof involves the study of the differences in an almost periodic pattern: the frequency of the difference $v\in\Z^n$ in the almost periodic pattern $\Gamma$ is the quantity $\rho_\Gamma(v) = D\big((\Gamma-v)\cap\Gamma)$ (see Definition~\ref{DefDiff}). To start the proof, we remark that if $\Gamma\subset\Z^n$ is an almost periodic pattern whose density is bigger than $1/2 + \delta$, then for every $v\in\Z^n$, the frequency of the difference $v$ in $\Gamma$ satisfies $\rho_\Gamma(v) \ge 2\delta$ (Lemma~\ref{majoration}). The theorem results from this remark and from the fact that matrices with a rate of injectivity equal to 1 are very rare (Corollary~\ref{CoroHajos}). The proof also uses crucially the geometric construction used to compute the mean rate of injectivity $\overline\tau$ (Proposition~\ref{FormTau1}).

We then state a variation of Theorem~\ref{ConjIntro} in the case of isometries (Theorem~\ref{AnswerConjIsom}) : if $(P_k)_{k\ge 1}$ is a generic sequence of matrices of $O_n(\R)$, then $\tau^\infty((P_k)_k) = 0$.

\noindent To prove this result, we study in more detail the action of the discretizations of linear maps on the frequencies of the differences of an almost periodic pattern (Proposition~\ref{ActionDiff}). We show that the discretization $\widehat A$ of $A\in GL_n(\R)$ acts smoothly on the frequency of differences: when we want to compute $\widehat A(\Gamma)$, we apply A and then make a projection, whereas when we want to compute $\rho_{\widehat A(\Gamma)}$, we apply $A$ and then make a \emph{weighted} projection. To prove Theorem~\ref{DD}, we combine this computation with the Minkowski-like theorem for the differences in an almost periodic pattern we obtained in the previous chapter (Theorem \ref{MinkAlm}).

Finally, we use the notion of model set\footnote{Sometimes called ``cut-and-project'' set.} (Definition~\ref{DefModel}) to improve Theorem~\ref{PerLin1}. It is defined as follows: given a lattice $\Lambda$ of $\R^{m+n}$ and a ``regular'' set $W\subset\R^m$, we select the points of the lattice whose projection on the $m$ first coordinates belongs to $W$; the projection on the $n$ last coordinates of these points is called the \emph{model set} modelled on $\Lambda$ and $W$. Model sets form a subclass of almost periodic patterns (Theorem~\ref{ModelAlmost}). The introduction of this notion is supported by the fact that for every sequence $(A_k)_{k\ge 1}$ of invertible linear maps, the set $(\widehat{A_k}\circ\cdots\circ\widehat{A_1}) (\Z^n)$ is a model set (Proposition~\ref{ImgModel}). This fact gives us a more precise information about the structure of the successive images of $\Z^n$ by discretizations of linear maps, as the definition of model set is stronger than that of almost periodic pattern. This viewpoint is fruitful to study the question raised in Theorem~\ref{ConjIntro}: first of all, an argument of equidistribution leads to a generalization of the geometric construction to compute the mean rate of injectivity in time 1 to the rate of injectivity in arbitrary time $k$ (Propositions~\ref{CalculTauxModel} and \ref{CalculTauxModel2}). Moreover, these constructions allows us to give an alternative proof of Theorem~\ref{PerLin1}, and to finally prove Theorem~\ref{ConjIntro}.

\section{Study of the continuity of the rate of injectivity}\label{SecCont}

In this section, we study the properties of continuity of the function $\tau$. We first show that there exists some matrices in which $\tau$ is not continuous (Proposition~\ref{rateEx}). However, we show in Proposition~\ref{ThMeanRate} that $\tau$ is continuous at every totally irrational matrix; more precisely, the map $\tau$ coincides with the \emph{mean rate of injectivity} $\overline \tau$ (Definition~\ref{DefMeanRate}) on the set of totally irrational matrices, and the map $\overline \tau$ is continuous on the whole set $GL_n(\R)$. We finish this section in estimating the lack of continuity of $\tau$ at the matrices which are not totally irrational (Proposition~\ref{oscill}).
\bigskip

In the following example, we show that on some rational matrices, the rate of injectivity is not continuous. In particular, it will give us an example where Lemma \ref{équi} is not true when there is no restriction about the value of $v$ on $I_\Q(A)$ (the set of indices corresponding to rational coefficients, see Notation \ref{intelligent}).

\begin{figure}[t]
\begin{center}
\begin{tikzpicture}[scale=.93]
\clip (-1.6,-1.6) rectangle (2.6,2.6);
\draw (-2,-2) grid (3,3);
\draw (.5,.5) node[above]{$0$};
\foreach\i in {-4,...,4}{
\foreach\j in {-2,...,2}{
\draw[color=green!40!black] (.5*\i-\j+.5,.5*\i+\j+.5) node {$\bullet$};
\draw[color=black] (\i+.5,\j+.5) node {$\cdot$};
}}
\end{tikzpicture}\hspace{12pt}
\begin{tikzpicture}[scale=.93]
\clip (-1.6,-1.6) rectangle (2.6,2.6);
\draw (-2,-2) grid (3,3);
\draw (.5,.5) node[above]{$0$};
\foreach\i in {-4,...,4}{
\foreach\j in {-2,...,2}{
\draw[color=green!40!black] (.5*\i-\j+.7,.5*\i+\j+.4) node {$\bullet$};
\draw[color=black] (\i+.5,\j+.5) node {$\cdot$};
}}
\end{tikzpicture}\hspace{12pt}
\begin{tikzpicture}[scale=.93]
\clip (-1.6,-1.6) rectangle (2.6,2.6);
\draw (-2,-2) grid (3,3);
\draw (.5,.5) node[above]{$0$};
\foreach\i in {-4,...,4}{
\foreach\j in {-2,...,2}{
\draw[color=green!40!black] (.5*\i-\j+.7,.5*\i+\j+.6) node {$\bullet$};
\draw[color=black] (\i+.5,\j+.5) node {$\cdot$};
}}\end{tikzpicture}
\caption{Sets $f_0(\Z^n) + v$ for $v=(0,0)$, $v = (0.2\,,\,-0.1)$ and $v = (0.2,\,0.1)$}\label{flemme}
\end{center}
\end{figure}
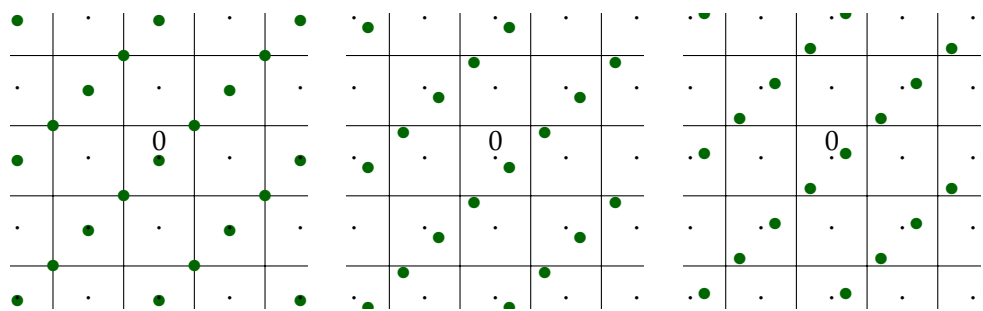

Our counterexample is given by irrational perturbations of the matrix
\[f_0 = \left(\begin{matrix} \frac12 & -1 \\ \frac12 & 1 \end{matrix}\right).\]
It can be easily seen that the rate of injectivity of $f_0$ is 1/2 (see Figure \ref{flemme}). Remark that it depends on the choice made for the projection of $\R^2$ on $\Z^2$: if we make the same choice for both directions (that is what we have chosen in Definition~\ref{DefDiscrLin}) then the rate of injectivity is $1/2$, otherwise this rate of injectivity is $1$. For $\varep$ ``small'', we consider the linear map
\[f_\varep = \left(\begin{matrix} 1+\varep & 0 \\ 0 & \frac{1}{1+\varep} \end{matrix}\right) f_0,\]
which is close to $f_0$ when $\varep$ is small.

\begin{prop}\label{rateEx}
The rate of injectivity of $f_\varep$ tends $\frac 34$ when $\varep\notin\Q$ tends to 0.
\end{prop}

In particular, the rate of injectivity $\tau$ is not continuous in $f_0$.

\begin{lemme}\label{CalcEx}
\begin{itemize}
\item For all $v\in ]0,\frac 12[^2\cup]\frac 12,1[^2$ modulo $\Z^n$, the rate of injectivity of $f_0+v$ is $\frac 12$.
\item For all $v\in ]0,\frac 12[\times ]\frac 12,1[\, \cup\, ]\frac 12,1[\times]0,\frac 12[$ modulo $\Z^n$, the rate of injectivity of $f_0+v$ is $1$.
\end{itemize}
\end{lemme}

\begin{proof}[Proof of Lemma \ref{CalcEx}]
We want to know when it is possible to find two different vectors $x,y\in\Z^n$ such that $\pi(f_0x+v) = \pi(f_0y+v)$. This implies that $d_\infty(f_0x,f_0y)< 1$, a simple calculation (see Figure \ref{flemme}) shows that in this case, $x=y\pm (1,0)$. As $f_0(2\Z\times \Z)\subset \Z^2$, we can suppose that $x=0$ ; we want to know if for $y = \pm (1,0)$, we have $\pi(f_0x + v) = \pi(f_0y + v)$, that is $\pi(v) = \pi(v\pm (1/2,1/2))$. Again, a simple calculation shows that this occurs if and only if $v\in ]0,\frac 12[^2\cup]-\frac 12,0[^2$.
\end{proof}

Thus, for half of the vectors $v$ (for Lebesgue measure), the rate of injectivity of $f_0+v$ is 1 and for the other half of the vectors $v$, the rate of injectivity of $f_0+v$ is $\frac 12$. Proposition~\ref{rateEx} then follows from an argument of equirepartition. To make it rigorous, we define the mean rate of injectivity.

\begin{definition}\label{DefMeanRate}
For $A\in GL_n(\R)$, the quantity\index{$\tau o$ @ $\overline\tau$}
\[\overline{\tau}(A) = \int_{\T^n}\tau(A+v)\ud v\]
is called the \emph{mean rate of injectivity} of $A$.

We say that a matrix $A\in GL_n(\R)$ is \emph{totally irrational} if the image $A(\Z^n)$ is equidistributed\footnote{It is equivalent to require that it is dense instead of equidistributed; it is also equivalent to ask that $\Z^n$ is equidistributed modulo $A(\Z^n)$.} modulo $\Z^n$; in particular, this is true when the coefficients of $A$ form a $\Q$-free family.
\end{definition}

The motivation of this definition is that the mean rate of injectivity is continuous and coincides with the rate of injectivity on totally irrational matrices.

\begin{prop}\label{ThMeanRate}
The mean rate of injectivity $\overline \tau$ is continuous on $GL_n(\R)$. Moreover, if $A$ is totally irrational, then $\overline \tau(A) = \tau(A)$; and even more, $\tau(A) = \tau(A+v)$ for every $v\in\T^n$.
\end{prop}

Thus, the restriction of the rate of injectivity to totally irrational matrices is a continuous function, which extends to $GL_n(\R)$ into a continuous function. In particular, the restriction of this function to any bounded subset of $SL_n(\R)$ is uniformly continuous.

Proposition~\ref{ThMeanRate}, combined with Lemma \ref{CalcEx}, obviously implies Proposition~\ref{rateEx} (as for $\varep\notin\Q$, the map $f_\varep$ is totally irrational).

\begin{proof}[Proof of Proposition \ref{ThMeanRate}]
The continuity of $\overline\tau$ will be obtained later as a direct consequence of Corollary~\ref{superCont}.

Let $A\in GL_n(\R)$ be a totally irrational matrix. We want to prove that $\tau(A) = \overline\tau(A)$. To do that, we show that $\tau(A) = \tau(A+v)$ for every $v\in\T^n$. As $\{Ax\mid x\in \Z^n\}$ is equidistributed modulo $\Z^n$, for every $v\in \T^n$, there exists $x\in\Z^n$ such that $d_\infty(v-Ax,\Z^n)<\varep$. But as the density is a limit  independent from the choice of the centre of the ball (Corollary~\ref{corolimitexist} and Theorem~\ref{imgquasi}), we have
\[\lim_{R\to +\infty}\frac{\card\big(\pi( A(\Z^n))\cap B_R\big)}{\card[B_R]} = \lim_{R\to +\infty}\frac{\card\big(\pi (A(\Z^n)+Ax)\cap B(Ax,R)\big)}{\card[B_R]}.\]
This proves that $\tau(A) = \tau(A+Ax)$. We then use Remark~\ref{RemContTrans}, which states that as $A$ is totally irrational, the function $v\mapsto \tau(A+v)$ is continuous. As the vector $Ax$ is arbitrarily close to the vector $v$, we get $\tau(A) = \tau(A+v)$.
\end{proof}

\begin{rem}\label{ContiTaux}
The same proof shows that if $\Gamma\subset\Z^n$ is an almost periodic pattern, and $A\in GL_n(\R)$ is such that $A(\Gamma)$ is equidistributed modulo $\Z^n$ (which is true for a generic $A$, see Lemma~\ref{passifacil}), then $B\mapsto D\big(\widehat B(\Gamma)\big)$ is continuous in $B=A$.
\end{rem}

We now study the behaviour of $\tau$ at the matrices where it is not continuous. In particular, we prove that the set of matrices on which the rate of injectivity makes ``big jumps'' is small. More precisely, we define the oscillation of a map.

\begin{definition}
For $A\in GL_n(\R)$, the \emph{oscillation} of $\tau$ at $A$ is the quantity\index{$\omega_\tau$}
\[\omega_\tau(A) = \underset{r\to 0}{\overline\lim}\ \sup\big\{|\tau(B_1) - \tau(B_2)|\ \big\vert\ \|A-B_i\|\le r\big\}.\]
\end{definition}

\begin{prop}\label{oscill}
For every $\varep>0$, the set $\big\{ A\in GL_n(\R) \mid \omega_\tau(A)>\varep\big\}$ is locally contained in a finite union of hyperplanes.
\end{prop}

To prove this proposition, we will need a technical lemma which requires the following definition.

\begin{definition}
Let $\Lambda$ be a lattice of $\R^n$. Then \cite[Chap. VII, \S 1, 2]{MR1726872} implies that the $\Z$-module $\Lambda + \Z^n$ can be decomposed into two complementary modules \index{$\Lambda_\text{cont}$}\index{$\Lambda_\text{discr}$}
\[\Lambda + \Z^n = \Lambda_\text{cont} \oplus \Lambda_\text{discr},\]
such that $\Lambda_\text{cont}$ is dense in the vector space it spans and $\Lambda_\text{discr}$ is discrete.
\end{definition}

During the proof of Proposition~\ref{oscill}, we will make use of the following lemma, that we will not prove.

\begin{lemme}\label{LemOscill}
For every $\varep>0$, the set of matrices $A\in GL_n(\R)$ such that $(A\Z^n)_\text{discr}$ has a fundamental domain of diameter smaller than $\varep$ contains the complement of a locally finite union of hyperplanes.
\end{lemme}

Remark that if $A$ is totally irrational, then $(A\Z^n)_\text{discr} = \{0\}$

\begin{proof}[Proof of Proposition \ref{oscill}]
Let $\varep>0$. We consider a matrix $A\in GL_n(\R)$ such that $(A\Z^n)_\text{discr}$ has a fundamental domain $\mathcal D$ of diameter smaller than $\varep$; Lemma~\ref{LemOscill} asserts that such matrices contains the complement of a locally finite union of hyperplanes. We denote $\Lambda = A\Z^n$. We use Proposition~\ref{ThMeanRate}, which states that $\overline\tau$ is continuous and coincides with $\tau$ on totally irrational matrices, to estimate the oscillation of $\tau$ in $A$ : it is smaller than the maximal value of $\big| D(\pi(\Lambda + v)) - D(\pi(\Lambda + v'))\big|$, when the vectors $v$ and $v'$ run through $\T^n$.

First of all, we remark that the map
\[v\mapsto D(\pi(\Lambda + v))\]
is $\Lambda + \Z^n$-periodic. Thus, we only have to estimate $\big| D(\pi(\Lambda + v)) - D(\pi(\Lambda + v'))\big|$ for $v$ and $v'$ in $\mathcal D$. We denote by $(\Z^n)'$ the set of points in $\R^n$ whose discretization is not canonically defined, more precisely,
\[(\Z^n)' = \big\{(x_i)_{1\le i\le n} \in\R^n \mid \exists i : x_i\in \Z+1/2\big\}.\]
Then, for every $v,v'\in\R^n$, we have
\[\big| D(\pi(\Lambda + v)) - D(\pi(\Lambda + v'))\big| \le D\left\{x\in\Lambda\mid \exists w\in\mathcal D : x+w\in (\Z^n)'\right\},\]
and as $\diam(\mathcal D)$ is smaller than $\varep$, we get
\[\big| D(\pi(\Lambda + v)) - D(\pi(\Lambda + v'))\big| \le \frac{\card\big\{x\in\Lambda_\text{discr} \cap [0,1]^n\mid d(x,(\Z^n)')\le \varep\big\}}{\card(\Lambda_\text{discr} \cap [0,1]^n)}.\]
The proposition then follows from the fact that this last quantity is small (uniformly in $\varep$).
\end{proof}

\section{A geometric viewpoint on the mean rate}\label{ptgeom}

In this section, we present two geometric constructions to compute the rate of injectivity of a matrix.

\paragraph{First construction}

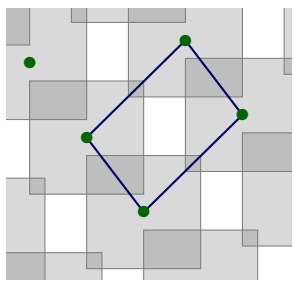
\begin{figure}[t]
\begin{center}
\begin{tikzpicture}[scale=1.5]
\clip (-1.2,-.6) rectangle (1.3,1.8);
\foreach\i in {-1,...,2}{
\foreach\j in {-1,...,3}{
\fill[color=gray,opacity = .3] (.866*\i-.5*\j-.5,.859*\i+.659*\j-.5) rectangle (.866*\i-.5*\j+.5,.859*\i+.659*\j+.5);
\draw[color=gray] (.866*\i-.5*\j-.5,.859*\i+.659*\j-.5) rectangle (.866*\i-.5*\j+.5,.859*\i+.659*\j+.5);
}}
\draw[color=blue!40!black,thick] (0,0) -- (.866,.859) -- (.366,1.518) -- (-.5,.659) -- cycle;
\foreach\i in {-2,...,2}{
\foreach\j in {-2,...,2}{
\draw[color=green!40!black] (.866*\i-.5*\j,.859*\i+.659*\j) node {$\bullet$};
}}
\end{tikzpicture}
\caption[First geometric construction to compute the rate of injectivity]{First geometric construction: the green points are the elements of $\Lambda$, the blue parallelogram is a fundamental domain of $\Lambda$ and the grey squares are centred on the points of $\Lambda$ and have radii $1/2$; their union form the set $U$. A point $x\in \Z^n$ belongs to $\pi(\Lambda)$ if and only if it belongs to at least one grey square.}\label{tourp}
\end{center}
\end{figure}

Let $A\in GL_n(\R)$ and $\Lambda = A(\Z^n)$. The density of $\pi(\Lambda)$ is the proportion of $x\in\Z^n$ belonging to $\pi(\Lambda)$; in other words the proportion of $x\in\Z^n$ such that there exists $\lambda\in \Lambda$ whose distance to $x$ (for $\|\cdot\|_\infty$) is smaller than $1/2$. This property only depends on the value of $x$ modulo $\Lambda$. If we consider the union\index{$U$}
\[U = \bigcup_{\lambda\in\Lambda} B(\lambda,1/2)\]
of balls of radius $1/2$ centred on the points of $\Lambda$ (see Figure~\ref{tourp}), then $x\in\pi(\Lambda)$ if and only if $x\in U\cap\Z^n$. So, if we set $\nu$ the measure of repartition of the $x\in\Z^n$ modulo $\Lambda$, that is\index{$\nu$}
\[\nu = \lim_{R\to+\infty} \frac{1}{\card(B_R\cap \Z^n)} \sum_{x\in B_R\cap \Z^n} \delta_{\operatorname{pr}_{\R^n/\Lambda}(x)},\]
then we obtain the following formula (using also Equation~\eqref{EqTauDens}, which links $\tau(A)$ and $D(\Lambda)$).

\begin{prop}
For every $A\in GL_n(\R)$ (we identify $U$ with its projection of $\R^n/\Lambda$),
\[\tau(A) = |\det(A)| D\big(\pi(\Lambda)\big) = |\det(A)| \nu\big(\operatorname{pr}_{\R^n/\Lambda}(U)\big).\]
\end{prop}

In particular, when the matrix $A$ is totally irrational, the measure $\nu$ is the uniform measure; thus if $\mathcal D$ is a fundamental domain of $\R^n/\Lambda$, then $\tau(A)$ is the area of $\mathcal D \cap U$. The same holds for the mean rate of injectivity and any matrix (not necessarily totally irrational).

\begin{prop}\label{FormTau1}
For every $A\in GL_n(\R)$,
\[\overline\tau(A) = |\det(A)| \Leb\big(\operatorname{pr}_{\R^n/\Lambda}(U)\big) = |\det(A)| \Leb(U\cap \mathcal D).\]
Thus, the mean rate of injectivity can be seen as the area of an intersection of parallelepipeds (see Figure~\ref{tourp}).
\end{prop}

With the same kind of arguments, we easily obtain a formula for $\rho_{\widehat A(\Z^n)}(v)$ (the frequency of the difference $v$ in $\widehat A(\Z^n)$, see Definition~\ref{DefDiff}).

\begin{prop}\label{ActionDiffGeom}
If $A\in GL_n(\R)$ is totally irrational, then for every $v\in\Z^n$,
\[\rho_{\widehat A(\Z^n)}(v) = \Leb\big( B(v,1/2) \cap U\big).\]
\end{prop}

\begin{proof}[Sketch of proof of Proposition~\ref{ActionDiffGeom}]
We want to know which proportion of points $x\in\Gamma = \widehat A(\Z^n)$ are such that $x+v$ also belongs to $\Gamma$. But modulo $\Lambda = A(\Z^n)$, $x$ belongs to $\Gamma$ if and only if $x\in B(0,1/2)$. Similarly, $x+v$ belongs to $\Gamma$ if and only if $x\in B(-v,1/2)$. Thus, by equirepartition, $\rho_{\widehat A(\Z^n)}(v)$ is equal to the area of $B(v,1/2) \cap U$.
\end{proof}

From Proposition~\ref{FormTau1}, we deduce the continuity of $\overline\tau$.

\begin{coro}\label{superCont}
The mean rate of injectivity is a continuous function on $GL_n(\R)$, which is locally polynomial with degree smaller than $n$ in the coefficients of the matrix.
\end{coro}

\begin{figure}[t]
\begin{minipage}[c]{.33\linewidth}
	\includegraphics[width=\linewidth, trim = 1.6cm .5cm 3.4cm 1cm, clip]{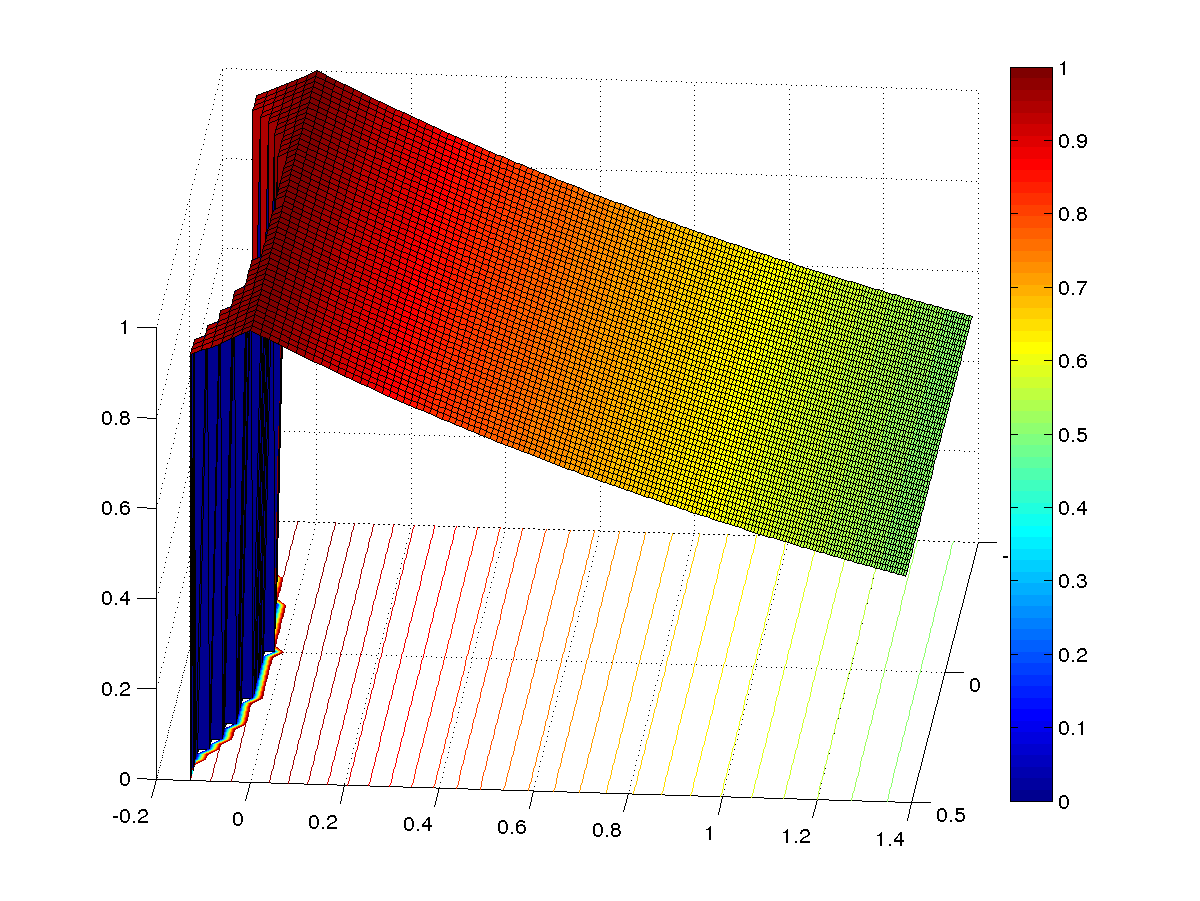}
\end{minipage}\hfill
\begin{minipage}[c]{.33\linewidth}
	\includegraphics[width=\linewidth, trim = 1.6cm .5cm 3.4cm 1cm, clip]{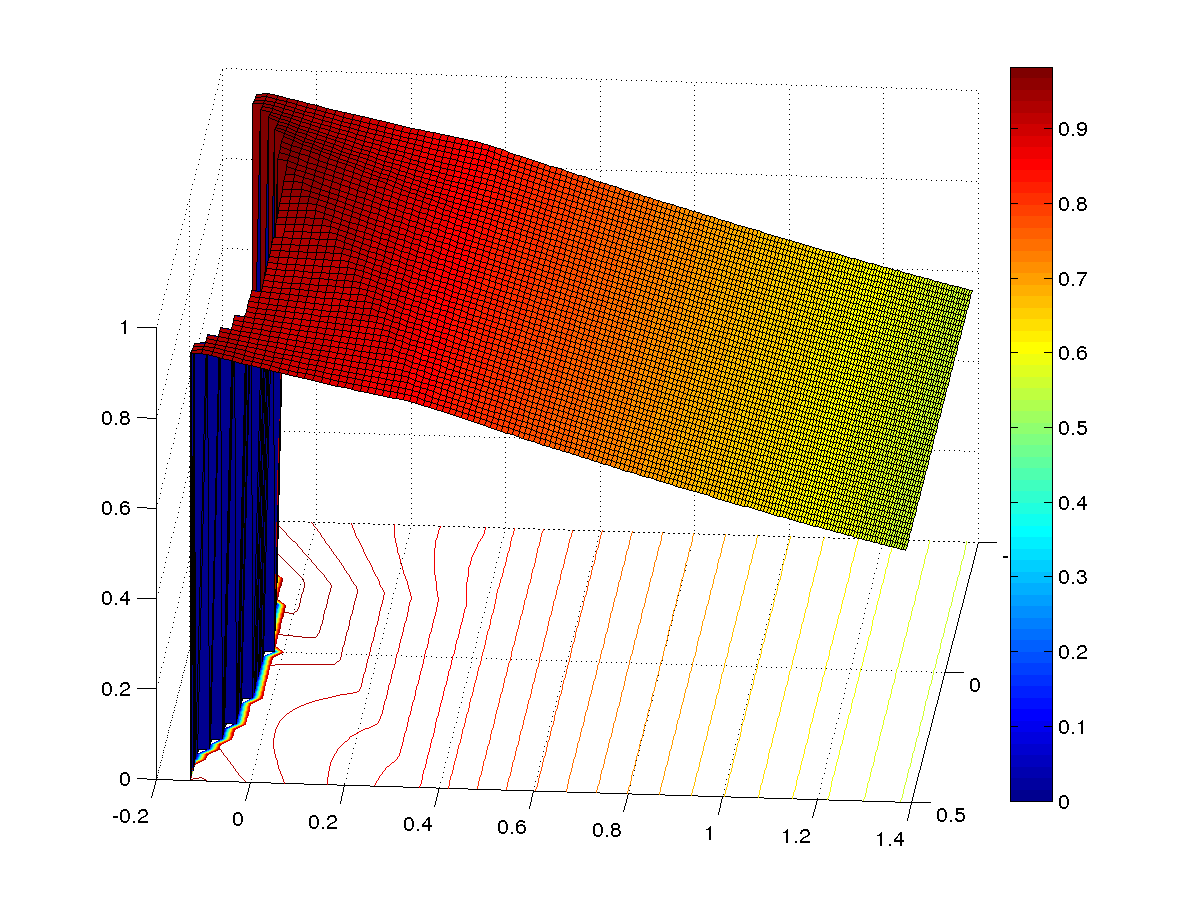}
\end{minipage}\hfill
\begin{minipage}[c]{.33\linewidth}
	\includegraphics[width=\linewidth, trim = 1.6cm .5cm 3.4cm 1cm, clip]{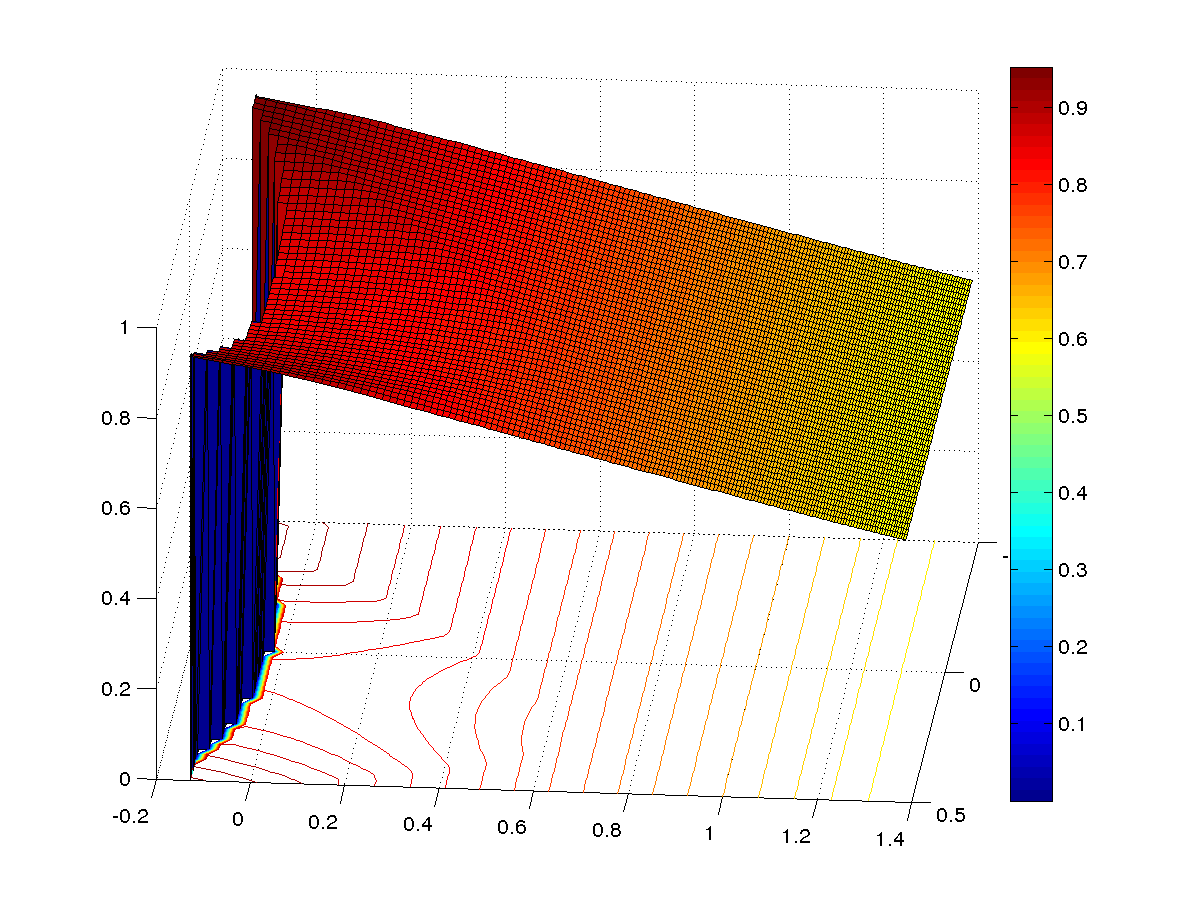}
\end{minipage}

\begin{minipage}[c]{.33\linewidth}
	\includegraphics[width=\linewidth, trim = 1.6cm .5cm 3.4cm 1cm, clip]{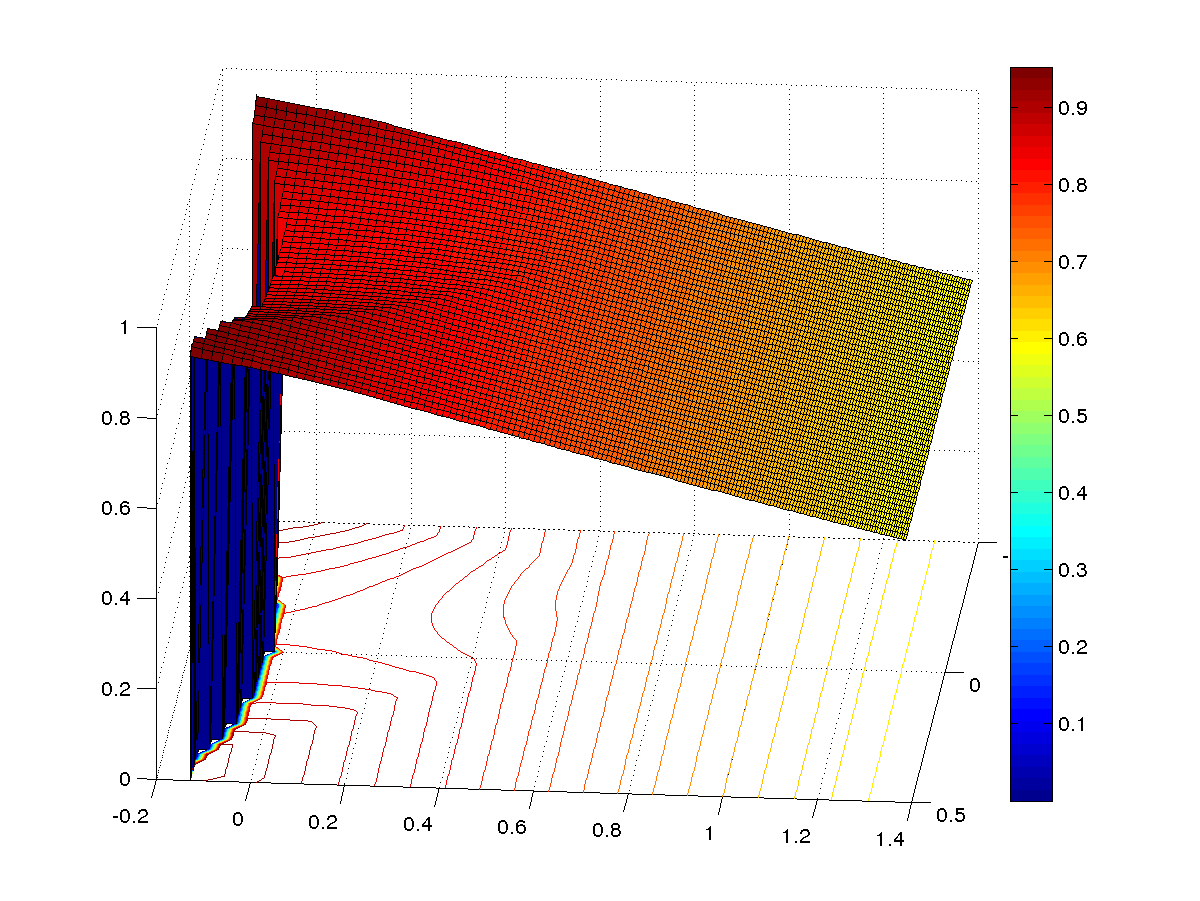}
\end{minipage}\hfill
\begin{minipage}[c]{.33\linewidth}
	\includegraphics[width=\linewidth, trim = 1.6cm .5cm 3.4cm 1cm, clip]{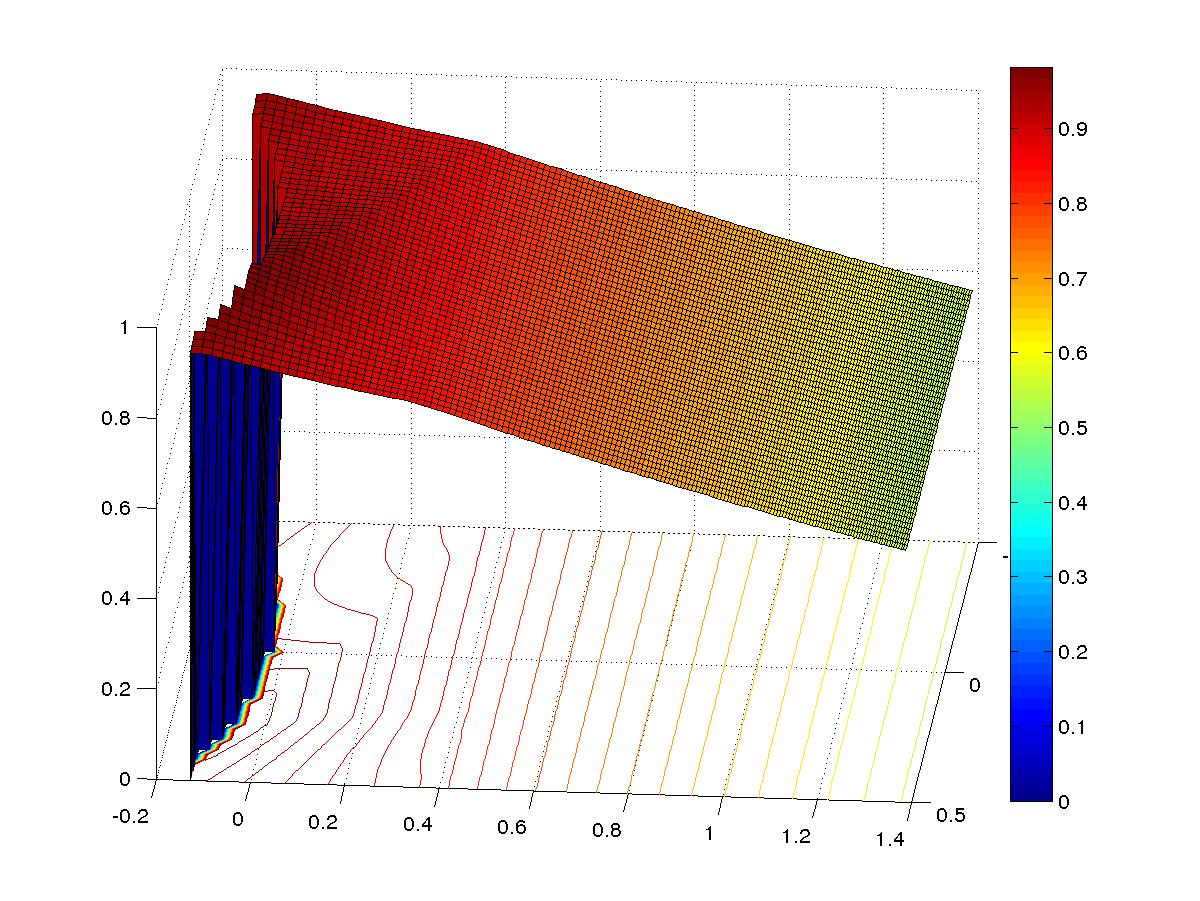}
\end{minipage}\hfill
\begin{minipage}[c]{.33\linewidth}
	\includegraphics[width=\linewidth, trim = 1.6cm .5cm 3.4cm 1cm, clip]{Fichiers/TauxPourAngle25.png}
\end{minipage}
\caption[Mean rate of injectivity on $SL_2(\R)$]{Mean rate of injectivity on $SL_2(\R)$ on the fundamental domain $\mathcal D$ (more precisely, $\mathcal D = \{z\in\C\mid |z|>1,\,\operatorname{Im} z>0,\, -1/2<\operatorname{Re} z<1/2\}$) of the modular surface, for various angles in $T^1\mathcal D$ ($0$, $\pi/10$, $\pi/5$, $3\pi/10$, $2\pi/5$, $\pi/2$).}\label{FigGradMean}
\end{figure}

\begin{figure}
\begin{center}
\includegraphics[width=11cm,clip]{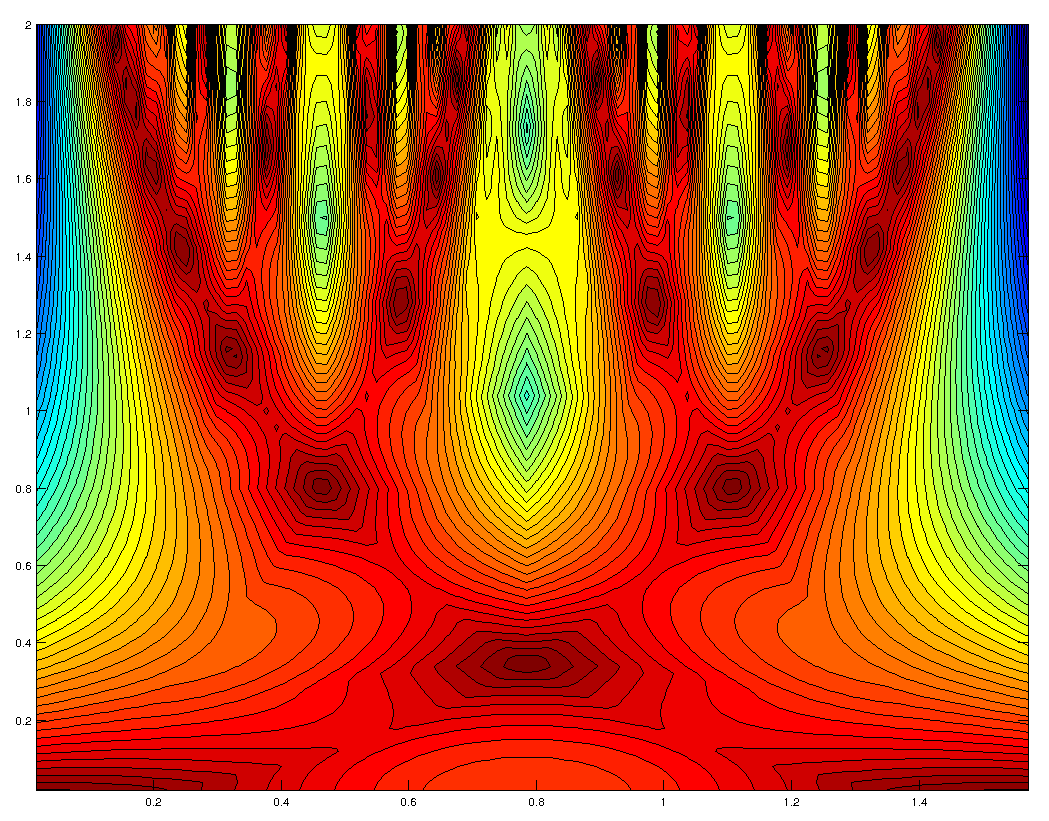}
\caption[Action of the geodesic flow on the mean rate of injectivity]{Action of the geodesic flow on the mean rate of injectivity: the red corresponds to a rate of injectivity close to 1 and the blue to a rate of injectivity close to 0. The figure represents the rate of injectivity on geodesics of the space $SL_2(\R)/SL_2(\Z)$, identified with the unitary tangent space $T^1 S$ of the modular surface $S$ starting from the point $i$, depending on the time (on the vertical axis) and on the starting angle of the geodesic (on the horizontal axis). The mean rate of injectivity of each lattice has been computed in using Proposition~\ref{FormTau2}}\label{FlowMean}
\end{center}
\end{figure}

Moreover, this construction gives a quick algorithm to compute numerically the mean rate of injectivity of some matrices: this algorithm is much more efficient than the naive one consisting in computing the cardinality of the image of a large ball by the discretization, see Figure~\ref{FigGradMean} and Figure~\ref{FlowMean}.

It also allows to compute simply the mean rate of injectivity of some examples of matrices.

\begin{app}\label{AppRot}
For $\theta\in[0,\pi/2]$, the mean rate of injectivity of a rotation of $\R^2$ of angle $\theta$ is (see Figure \ref{RateRot}).
\[\overline\tau(R_\theta) = 1-(\cos(\theta)+\sin(\theta)-1)^2.\]
\end{app}

\begin{figure}
\begin{center}
\begin{minipage}[b]{.47\linewidth}
\begin{center}
\begin{tikzpicture}[scale=1.4]
\clip (-1,-.6) rectangle (1.5,2.2);
\foreach\i in {0,...,1}{
\foreach\j in {0,...,1}{
\fill[color=gray,opacity = .3] (-.5+.866*\i-.5*\j,-.5+.5*\i+.866*\j) rectangle (.5+.866*\i-.5*\j,.5+.5*\i+.866*\j);
\draw[color=gray] (-.5+.866*\i-.5*\j,-.5+.5*\i+.866*\j) rectangle (.5+.866*\i-.5*\j,.5+.5*\i+.866*\j);
}}
\draw (0,0) -- (.6,0);
\draw[thick] (.5,0) arc (0:30:.5);
\draw (.62,.18) node {$\theta$};
\draw[color=red!60!black,thick] (0,.866) -- (0,.5) -- (.366,.5) -- (.366,.866) -- cycle;
\draw[color=blue!40!black,thick] (0,0) -- (.866,.5) -- (.866-.5,.866+.5) -- (-.5,.866) -- cycle;
\foreach\i in {0,...,1}{
\foreach\j in {0,...,1}{
\draw[color=green!40!black] (.866*\i-.5*\j,.5*\i+.866*\j) node {$\bullet$};
}}
\end{tikzpicture}
\caption[Mean rate of injectivity of a rotation]{Computation of the mean rate of injectivity of a rotation of $\R^2$: it is equal to $1$ minus the area of the interior of the red square.}\label{RateRot}
\end{center}
\end{minipage}\hfill
\begin{minipage}[b]{.47\linewidth}
\begin{center}
\begin{tikzpicture}[scale = 1.4]
\draw[color=gray,fill=gray,fill opacity=0.3] (0,0) -- (0,1) -- (1,1) -- (1,0) -- cycle;
\draw[color=gray,fill=gray,fill opacity=0.3] (.7,0) -- (.7,1) -- (1.7,1) -- (1.7,0) -- cycle;
\draw[color=gray,fill=gray,fill opacity=0.3] (.4,1.43) -- (.4,2.43) -- (1.4,2.43) -- (1.4,1.43) -- cycle;
\draw[color=gray,fill=gray,fill opacity=0.3] (1.1,1.43) -- (1.1,2.43) -- (2.1,2.43) -- (2.1,1.43) -- cycle;
\draw[color=blue!40!black,thick] (.5,.5) -- (1.2,.5) -- (1.6,1.93) -- (.9,1.93) -- cycle;
\draw [color=green!40!black] (.5,.5) node{$\bullet$};
\draw [color=green!40!black] (1.2,.5) node{$\bullet$};
\draw [color=green!40!black] (1.6,1.93) node{$\bullet$};
\draw [color=green!40!black] (.9,1.93) node{$\bullet$};
\draw[color=red!60!black,thick] (0.64,1) -- (1.34,1) -- (1.46,1.43) -- (0.76,1.43) -- cycle;
\end{tikzpicture}
\caption[Mean rate of injectivity of a lattice having one vector parallel with the horizontal axis, $\ell<1$]{Computation of the mean rate of injectivity of a lattice having one vector parallel with the horizontal axis with length $\ell<1$.}\label{RateParaLP}
\end{center}
\end{minipage}
\end{center}
\end{figure}

\begin{figure}[t]
\begin{center}
\begin{minipage}[b]{.47\linewidth}
\begin{center}
\begin{tikzpicture}[scale = 1.4]
\draw[color=gray,fill=gray,fill opacity=0.3] (0,0) -- (0,1) -- (1,1) -- (1,0) -- cycle;
\draw[color=gray,fill=gray,fill opacity=0.3] (1.4,0) -- (1.4,1) -- (2.4,1) -- (2.4,0) -- cycle;
\draw[color=gray,fill=gray,fill opacity=0.3] (.6,0.71) -- (.6,1.71) -- (1.6,1.71) -- (1.6,0.71) -- cycle;
\draw[color=gray,fill=gray,fill opacity=0.3] (2,0.71) -- (2,1.71) -- (3,1.71) -- (3,0.71) -- cycle;
\draw[color=blue!40!black,thick] (.5,.5) -- (1.9,.5) -- (2.5,1.21) -- (1.1,1.21) -- cycle;
\draw [color=green!40!black] (.5,.5) node{$\bullet$};
\draw [color=green!40!black] (1.9,.5) node{$\bullet$};
\draw [color=green!40!black] (2.5,1.21) node{$\bullet$};
\draw [color=green!40!black] (1.1,1.21) node{$\bullet$};
\draw[color=red!60!black,thick] (1,.5) -- (1.4,.5) -- (1.4,.71) -- (1,.71) -- cycle;
\draw[color=red!60!black,thick] (1.6,1) -- (2,1) -- (2,1.21) -- (1.6,1.21) -- cycle;
\end{tikzpicture}
\caption[Mean rate of injectivity of a lattice having one vector parallel with the horizontal axis, $\ell>1$ and $x>\ell-1$]{Computation of the mean rate of injectivity of a lattice having one vector parallel with the horizontal axis with length $\ell\in]1,2[$ for a parameter $x>\ell-1$.}\label{RateParaLGXG}
\end{center}
\end{minipage}\hfill
\begin{minipage}[b]{.47\linewidth}
\begin{center}
\begin{tikzpicture}[scale = 1.4]
\draw[color=gray,fill=gray,fill opacity=0.3] (0,0) -- (0,1) -- (1,1) -- (1,0) -- cycle;
\draw[color=gray,fill=gray,fill opacity=0.3] (1.4,0) -- (1.4,1) -- (2.4,1) -- (2.4,0) -- cycle;
\draw[color=gray,fill=gray,fill opacity=0.3] (.2,0.71) -- (.2,1.71) -- (1.2,1.71) -- (1.2,0.71) -- cycle;
\draw[color=gray,fill=gray,fill opacity=0.3] (1.6,0.71) -- (1.6,1.71) -- (2.6,1.71) -- (2.6,0.71) -- cycle;
\draw[color=blue!40!black,thick] (.5,.5) -- (1.9,.5) -- (2.1,1.21) -- (.7,1.21) -- cycle;
\draw [color=green!40!black] (.5,.5) node{$\bullet$};
\draw [color=green!40!black] (1.9,.5) node{$\bullet$};
\draw [color=green!40!black] (2.1,1.21) node{$\bullet$};
\draw [color=green!40!black] (.7,1.21) node{$\bullet$};
\draw[color=red!60!black,thick] (1,.5) -- (1.4,.5) -- (1.4,.71) -- (1,.71) -- cycle;
\draw[color=red!60!black,thick] (1.2,1) -- (1.6,1) -- (1.6,1.21) -- (1.2,1.21) -- cycle;
\draw[color=red!60!black,thick] (1.2,.71) -- (1.4,.71) -- (1.4,1) -- (1.2,1) -- cycle;
\end{tikzpicture}
\caption[Mean rate of injectivity of a lattice having one vector parallel with the horizontal axis, $\ell>1$ and $x<\ell-1$]{Computation of the mean rate of injectivity of a lattice having one vector parallel with the horizontal axis with length $\ell\in]1,2[$ for a parameter $x<\ell-1$.}\label{RateParaLGXP}
\end{center}
\end{minipage}
\end{center}
\end{figure}

\begin{app}\label{AppLattPara}
We consider a lattice of $\R^2$ of covolume 1 with one basis vector parallel with the horizontal axis. We have several cases:
\begin{enumerate}
\item If the length $\ell$ of the vector parallel with the horizontal axis is smaller than 1, then the mean rate of injectivity is equal to $\ell$, independently from the choice of the second basis vector (see Figure \ref{RateParaLP}).
\item If the length $\ell$ is bigger than 1 but smaller than 2, we choose $v = (x,1/\ell)$ another basis vector of the lattice. We can suppose that $-\ell/2\le x\le \ell/2$; by symmetry we only treat the case $x\in[0,\ell/2]$. The mean rate of injectivity is then a piecewise affine map with respect to $x$:
\begin{itemize}
\item if $\ell-1 \le x \le \ell/2$, then the mean rate of injectivity is constant and equal to $1-(\ell-1)(2/\ell-1) = \ell-2+2/\ell$ (see Figure \ref{RateParaLGXG});
\item if $0 \le x \le \ell-1$, then the mean rate of injectivity is equal to $1 - (\ell-1)(2/\ell-1) - (1-1/\ell)(\ell-1-x) = 1/\ell + x - x/\ell$ (see Figure \ref{RateParaLGXP}).
\end{itemize}
\end{enumerate}
We do not treat the other cases where $\ell>2$.
\end{app}

\paragraph{Second construction}

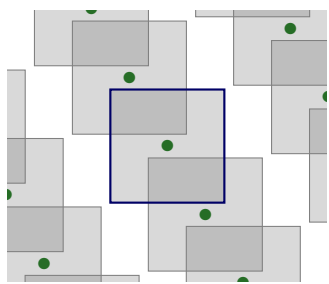
\begin{figure}[ht]
\begin{center}
\begin{tikzpicture}[scale=1.5]
\clip (-.9,-.7) rectangle (1.9,1.7);
\foreach\i in {-2,...,2}{
\foreach\j in {-2,...,2}{
\draw[color=green!40!black] (1.414*\i-.333*\j+.5,.433*\i+.605*\j+.5) node {$\bullet$};
\draw[color=gray, fill=gray,fill opacity=.3] (1.414*\i-.333*\j,.433*\i+.605*\j) rectangle (1.414*\i-.333*\j+1,.433*\i+.605*\j+1);
}}
\draw[color=blue!40!black,thick] (0,0) rectangle (1,1);
\end{tikzpicture}
\caption[Second geometric construction to compute the rate of injectivity]{Second geometric construction: the green points are the elements of $\Lambda$, the blue square is a fundamental domain of $\Z^n$ and the grey squares are centred on the points of $\Lambda$ and have radii $1/2$. The rate of injectivity is equal to the integral on the blue square of the inverse of the number of different squares the point belongs to.}\label{tourper}
\end{center}
\end{figure}

Here, we take another viewpoint to compute the rate of injectivity of $A$. To begin with, we remark that for every map $f$ defined on a finite set $E$, the cardinality of $f(E)$ can be computed by the formula
\[\card\big(f(E)\big)=\sum_{x\in E} \frac{1}{\card\big(f^{-1}(\{f(x)\})\big)};\]
in particular,
\begin{equation}\label{Merciki}
\tau(A)=\lim_{R\to +\infty} \frac{\card\big(\widehat A([B_R])\big)}{\card[B_R]} = \lim_{R\to +\infty} \frac{1}{\card[B_R]} \sum_{x\in [B_R]} \frac{1}{\card\big(\widehat A^{-1}(\{\hat A(x)\}))}.
\end{equation}
Thus, we are reduced to compute $\card\big(\widehat A^{-1}(\{\hat A(x)\})\big)$ for $x\in\Z^n$. This cardinality only depends on the class of $Ax$ modulo $\Z^n$, thus we can suppose that $Ax\in B_{1/2}$ (in other words, $\widehat A x= 0$).

\begin{lemme}\label{montagne}
\[\card\big(\widehat A^{-1}(\{\hat A(x)\})\big) = \psi(Ax),\]
where $\psi$ is the function (different from the function $\1_U$ used previously; recall that $\Lambda = A\Z^n$)\index{$\psi$}
\[\psi = \sum_{\lambda\in\Lambda} \1_{B(\lambda,1/2)}\]
which is the sum of the indicator functions of the balls of radius $1/2$ centred on the points of $\Lambda$ (see Figure~\ref{tourper}).
\end{lemme}

\begin{proof}[Proof of Lemma~\ref{montagne}]
We denote by $\tilde x$ the projection of $Ax$ on the fundamental domain $B_{1/2}$ of $\Z^n$. We want to prove that there exists $k$ different vectors $y\in\Z^n$ such that $\widehat A (x) = \widehat A(y)$ if and only if $\psi(\tilde x) = k$. Indeed, if $\lambda\in \Lambda$ is such that $\tilde x\in B(\lambda,1/2)$, then $\tilde x-\lambda \in B_{1/2}$ and thus $Ax + \lambda$ projects on the same point of $\Z^n$ as $Ax$; the number of such vectors $\lambda$ determines the number of points $y\in\Z^n$ such that $\widehat A x = \widehat A y$.
\end{proof}

If we set $\nu'$\index{$\nu'$} the measure of repartition of the $\lambda\in\Lambda$ modulo $\Z^n$, that is
\[\nu' = \lim_{R\to+\infty} \frac{1}{\card (B_R\cap \Lambda)} \sum_{\lambda\in\Lambda\cap B_R} \delta_{\operatorname{pr}_{\R^n/\Z^n}(\lambda)},\]
then, by combining Equation~\eqref{Merciki} with Lemma~\ref{montagne}, we obtain

\begin{prop}\label{2Fr}
\[\tau(A) = \int_{B_{1/2}} \frac{1}{\psi(\lambda)}\, \ud \nu'(\lambda).\]
\end{prop}

And we have a similar statement for the mean rate of injectivity.

\begin{prop}\label{FormTau2}
\[\overline\tau(A) = \int_{B_{1/2}} \frac{1}{\psi(\lambda)}\, \ud \Leb(\lambda).\]
\end{prop}

These properties can be also directly deduced from the first geometric construction by applying a double counting argument. We state it in more general context because we will need a more precise statement in the sequel (here, we only need the case $m=0$, thus $B_1 = B$); this lemma allows to compute the area of the projection of a set on a quotient by a lattice.

\begin{lemme}\label{DoubleComptage}
Let $\Lambda_1$ be a subgroup of $\R^m$ which is a lattice in the vector space it spans, $\Lambda_2$ be such that $\Lambda_1 \oplus \Lambda_2$ is a lattice of covolume 1 of $\R^m$, and $B$ be a compact subset of $\R^m$. Let $B_1$ be the projection of $B$ on the quotient $\R^m/\Lambda_1$, and $B_2$ the projection of $B$ on the quotient $\R^m/(\Lambda_1\oplus\Lambda_2)$. We denote by
\[a_i = \Leb\big\{x\in B_1 \mid \card\{\lambda_2\in\Lambda_2 \mid x\in B_1+\lambda_2\} = i \big\}.\]
Then,
\[\Leb(B_2) = \sum_{i\ge 1} \frac{a_i}{i}.\]
\end{lemme}

In particular, the area of $B_2$ (the projection on the quotient by the direct sum of both lattices) is smaller than (or equal to) that of $B_1$ (the projection on the quotient by the first lattice). The loss of area is given by the following corollary.

\begin{coro}\label{CoroSansNom}
With the same notations as for Lemma~\ref{DoubleComptage}, if we denote by
\[D_1 = \Leb\big\{x\in B_1 \mid \card\{\lambda_2\in\Lambda_2 \mid x\in B_1+\lambda_2\} \ge 2 \big\},\]
then,
\[\Leb(B_2) \le \Leb(B_1) - \frac{D_1}{2}.\]
\end{coro}

%

\begin{proof}[Proof of Lemma~\ref{DoubleComptage}]
Consider the function
\[\begin{array}{rcl}
\phi : & \R^n/\Lambda_1 & \longrightarrow \N^*\\
       & x              & \longmapsto \left\{
			\begin{array}{ll}
			\big(\card\{\lambda_2\in\Lambda_2 \mid x\in B_1+\lambda_2\}\big)^{-1}\quad & \text{if}\ x\in B_1\\
			0 &  \text{if}\ x\notin B_1.
			\end{array}\right.
\end{array}\]
Obviously, it satisfies $\int_{B_1} \phi = \sum_{i\ge 1} \frac{a_i}{i}$. Moreover, if we set
\[\begin{array}{rcl}
\Phi : & \R^n/(\Lambda_1\oplus\Lambda_2) & \longrightarrow \N^*\\
       & x                               & \longmapsto \sum_{\lambda_2\in\Lambda_2} \phi(x-\lambda_2),
\end{array}\]
then we easily see that on the one hand, $\Phi = \1_{\Lambda_2 + B_1}$, thus $\int \Phi = \Leb(B_2)$, and on the other hand, $\int \Phi = \operatorname{covol}(\Lambda_2)\int_{B_2} \phi = \sum_{i\ge 1} \frac{a_i}{i}$.
\end{proof}

\section{When is the rate close to 1?}

Obviously, the rate of injectivity $\tau(A)$ of any matrix of $A\in SL_n(\Z)$ is equal to 1. In Section \ref{SecCont}, we have found other examples of affine maps with determinant 1 whose rate of injectivity is also 1. Also, in Application \ref{AppLattPara}, taking $\ell=1$, this gives another class of examples where the rate of injectivity is one. In this section, we investigate more in detail this question: what are the matrices with determinant 1 whose mean rate of injectivity is 1? With Proposition~\ref{FormTau1}, we can reformulate this question in terms of intersection of cubes: if $\det(A)=1$, $\overline \tau(A)=1$ if and only if the cubes $B(\lambda,1/2)$, with $\lambda\in A \Z^n$, tile the space $\R^n$. This is a classical problem raised by H.~Minkowski in 1896 (see \cite{minkowski1910geometrie}), and answered by G.~Haj\'os in 1941.

\begin{theoreme}[Haj\'os, \cite{MR0006425}]\label{hajos}
Let $\Lambda$ be a lattice of $\R^n$. Then the collection of squares $\{B(\lambda,1/2)\}_{\lambda\in\Lambda}$ tiles the plane if and only if in a canonical basis of $\R^n$ (that is, permuting coordinates if necessary), $\Lambda$ admits a generating matrix which is upper triangular with ones on the diagonal.
\end{theoreme}

This theorem can be stated geometrically: if the collection of squares tiles the space $\R^n$, then at least two of these squares have a face in common. The conclusion of Theorem~\ref{hajos} can then be inferred by induction, quotienting in the direction which is orthogonal to this face and iterating the argument. The proof of this theorem involves fine results of group theory; we will not prove it here in the general case (see for example the excellent book \cite{MR1311249} for a complete investigation on the subject). Remark that Haj\'os theorem studies a particular case of the case of equality in Minkowski's theorem (Theorem \ref{Minkowski}): it states a necessary and sufficient condition for a lattice with covolume 1 to possess a non-trivial point in the boundary of the unit square centred at 0, but no non-trivial point in its interior. We give an elementary proof of Haj\'os theorem in dimension 2 (for dimension 3, there is also an elementary proof, see for example the book of H.~Minkowski \cite{minkowski1907diophantische}).

\begin{proof}[Proof of Theorem \ref{hajos} in dimension 2]
Let $\Lambda$ be a lattice of $\R^2$ such that the family $\{B(\lambda,1/2)\}_{\lambda\in\Lambda}$ tiles the plane. We consider the point $x_1 = (1/2,0)$. By the hypothesis of tiling, there exists $\lambda_1\in\Lambda\setminus\{0\}$ such that $x_1\in \overline{B(\lambda_1,1/2)}$. Remark that this implies that $\lambda_1$ has the form $\lambda_1 = (1,y)$. We have two cases.
\begin{enumerate}
\item If $\lambda_1 = (1,0)$, then we consider the point $x_2 = (0,1/2)$ and again, there exists $\lambda_2\in\Lambda\setminus\{0\}$ such that $x\in \overline{B(\lambda_2,1/2)}$. The vector $\lambda_2$ has the form $\lambda_2 = (*,1)$, thus the basis $(\lambda_1,\lambda_2)$ of $\Lambda$ has the desired form.
\item If $\lambda_1 = (1,y) \neq (1,0)$, we suppose that $y>0$ (doing a symmetry if necessary). Considering the point $x_2 = (1/2,1/2)$, this implies that $\lambda_2 = (0,1)$ belongs to $\Lambda$ (in other words, $x_2\in B(\lambda_2,1/2)$). The basis $(\lambda_1,\lambda_2)$ of $\Lambda$ has the desired form.
\end{enumerate}
\end{proof}

Combined with the results of the previous section (more precisely, Proposition \ref{FormTau2}), Haj\'os theorem (Theorem~\ref{hajos}) leads to the following corollary.

\begin{coro}\label{CoroHajos}
A matrix $A\in SL_n(\R)$ has a mean rate of injectivity equal to 1 if and only if there exists $B\in SL_n(\Z)$ such that in a canonical basis of $\R^n$ (that is, permuting coordinates if necessary), the matrix $AB$ is upper triangular with ones on the diagonal (equivalently, if and only if there exists a permutation matrix $P$ and a matrix $B\in SL_n(\Z)$ such that $PAB$ is upper triangular with ones on the diagonal). Remark that the set of such matrices is a locally finite union of manifolds of positive codimension.
\end{coro}

\section{Rate of injectivity of a sequence of matrices}

The aim of the end of this chapter is to study the asymptotic rate of injectivity of sequences of matrices with determinant 1. More precisely, we want to estimate the rate of injectivity $\tau^\infty((A_k)_k)$ (that is, the limit of the densities of the sets $(\widehat{A_k}\circ\cdots\circ\widehat{A_1})(\Z^n)$ when $k$ goes to infinity, see Definition~\ref{DefTaux}) when $(A_k)_{k\ge 1}$ is a \emph{generic} sequence of matrices of $SL_n(\R)$ in the following sense.

\begin{definition}\label{DefTopoSL}
We fix once for all a norm $\|\cdot\|$ on $M_n(\R)$. For a bounded sequence $(A_k)_{k\ge 1}$ of matrices of $SL_n(\R)$, we set\index{$\|(A_k)_k\|$}
\[\|(A_k)_k\| = \sup_{k\ge 1} \|A_k\|.\]
In other words, we consider the space $\ell^\infty(SL_n(\R))$ of uniformly bounded sequences of matrices of determinant 1.
\end{definition}

This metric is complete, thus there is a good notion of genericity on the set of bounded sequences of matrices of determinant 1 (see page~\pageref{DiscussGene}). We give the main theorem of this chapter, which states the behaviour of the asymptotic rate of injectivity of a sequence of matrices (stated in the introduction as Theorem~\ref{ConjIntro}).

\begin{theoreme}\label{ConjPrincip}
For a generic\footnote{Generic for the topology $\ell^\infty$ on $(SL_n(\R))^\N$, see Definition~\ref{DefTopoSL}.} sequence of matrices $(A_k)_{k\ge 1}$ of $SL_n(\R)$, we have\footnote{For a definition of $\tau^\infty$, see Definition~\ref{DefTaux}.}
\[\tau^\infty\big( (A_k)_{k\ge 1}\big) = 0.\]
\end{theoreme}

This statement was motivated by some numerical simulations (see Figure~\ref{TauxSuiteMat} page~\pageref{TauxSuiteMat}). In Part~\ref{PartTri}, we will deduce from this theorem a similar statement for generic diffeomorphisms of the torus (Theorem~\ref{limiteEgalZero}).

\begin{rem}
This statement remain true if we replace ``for a generic sequence of matrices $(A_k)_{k\ge 1}$ of $SL_n(\R)$'' by ``for a generic sequence of matrices $(A_k)_{k\ge 1}$ among matrices of determinant $\pm 1$'': the property we will need in the proofs is that the matrices preserve the volume of $\R^n$.
\end{rem}

\begin{rem}\label{ContiTauxInfini}
Remark~\ref{ContiTaux} implies that for a generic sequence $(A_i)_{i\ge 1}$ of matrices of $SL_n(\R)$ and for every $k\in\N^*$, the map $(B_i)_{i\ge 1} \mapsto \tau^k(B_1,\cdots,B_k)$ is continuous in $A_1,\cdots,A_k$. Thus, as the asymptotic rate of injectivity is the infimum of the rates of injectivity in times $k$ for $k\in\N^*$, it is upper semi continuous on every generic sequence of matrices. In particular, for every generic sequence $(A_i)_{i\ge 1}$ of matrices of $SL_n(\R)$, and every $\varep>0$ there exists a neighbourhood of this sequence on which the asymptotic rate of injectivity is smaller than $\varep$.
\end{rem}

We begin the study of Theorem~\ref{ConjPrincip} by handling the easiest case: the result is true when we restrict to diagonal matrices (in the canonical basis). More precisely, the set of diagonal matrices is a Baire space, so we can talk about generic diagonal matrices, and we have the following property.

\begin{prop}\label{CasLin}
Let $(A_k)_{k\ge 1}$ be a generic sequence of diagonal matrices of $SL_n(\R)$. Then $\tau^\infty((A_k)_k)=0$.
\end{prop}

\begin{proof}[Proof of Proposition \ref{CasLin}]
We denote $A_k = \operatorname{Diag}(\lambda_{k,1},\cdots,\lambda_{k,n})$. Since $A_k \in SL_n(\R)$ for every $k$, there exists $i_0\in\llbracket 1,n\rrbracket$ and an infinite number of integers $\ell$ such that $\prod_{k=0}^\ell \lambda_{k,i_0} \le 1$.
We then use the following lemma.

\begin{lemme}\label{LemDiag}
Let $\delta_0>0$, and $(\lambda_1,\cdots,\lambda_\ell)$ be a collection of real numbers such that $\lambda_i>\delta_0$ for every $i$ and $\prod_{k=1}^\ell\lambda_k\le 1$. Then, for every $\delta\in ]0,\delta_0[$, we have 
\[\prod_{k=1}^\ell(\lambda_k-\delta)\le \frac{1}{1+\delta (\ell+1)}.\]
\end{lemme}

\begin{proof}[Proof of Lemma \ref{LemDiag}]
We decompose the difference $\prod_{k=1}^\ell \lambda_k - \prod_{k=1}^\ell(\lambda_k-\delta)$ to have a telescopic sum:
\begin{align*}
\prod_{k=1}^\ell \lambda_k - \prod_{k=1}^\ell(\lambda_k-\delta) = & \sum_{i=1}^{\ell} \left(\prod_{k=1}^i \lambda_k \prod_{k=i+1}^\ell (\lambda_k - \delta) - \prod_{k=1}^{i-1} \lambda_k \prod_{k=i}^\ell (\lambda_k - \delta)\right)\\
= & \sum_{i=1}^{\ell} \left(\prod_{k=1}^{i-1} \lambda_k \prod_{k=i+1}^\ell (\lambda_k - \delta)\right) \big(\lambda_i - (\lambda_i - \delta)\big)\\
= & \sum_{i=1}^{\ell} \left(\prod_{k=1}^{i} \lambda_k \prod_{k=i+1}^\ell (\lambda_k - \delta)\right) \frac{\delta}{\lambda_i}\\
\ge & \delta \left(\prod_{k=1}^\ell (\lambda_k - \delta)\right) \sum_{i=1}^{\ell} \frac{1}{\lambda_i}.
\end{align*}
As a consequence, using the inequality of arithmetic and geometric means,
\[\prod_{k=1}^\ell \lambda_k - \prod_{k=1}^\ell(\lambda_k-\delta) \ge (\ell+1) \delta \left(\prod_{k=1}^\ell (\lambda_k - \delta)\right) \left(\prod_{i=1}^{\ell} \frac{1}{\lambda_i}\right)^{1/\ell},\]
and as $\prod_{i=1}^\ell\lambda_i\le 1$, we get
\[\prod_{k=1}^\ell \lambda_k - \prod_{k=1}^\ell(\lambda_k-\delta) \ge (\ell+1) \delta \left(\prod_{k=1}^\ell (\lambda_k - \delta)\right).\]
This implies that 
\[\prod_{k=1}^\ell (\lambda_k - \delta) \le \frac{\prod_{k=1}^\ell \lambda_k}{(\ell+1) \delta+1}\le \frac{1}{1+\delta(\ell+1)}.\]
\end{proof}

Proposition \ref{CasLin} easily follows from Lemma~\ref{LemDiag} and the following property, which is specific to the case where the matrices are diagonal: the rate of injectivity of $A_1,\cdots,A_k$ is smaller than the smallest eigenvalue of $A_1\circ\cdots \circ A_k$. Indeed, the eigenspace associated to this eigenvalue $\lambda$ is stable by discretization, and on this eigenspace the image of a big ball of radius $R$ is included in a ball of radius $\lambda R + o(R)$.
\end{proof}

\section[Frequency of differences: the rate is smaller than $1/2$]{Frequency of differences and a first result: the asymptotic rate of injectivity is generically smaller than $1/2$}

In this section, we prove a weak version of Theorem~\ref{ConjPrincip}, which states that the rate of injectivity of a generic sequence of matrices is smaller than one half. It is based on the following lemma, which asserts that if a difference appears rarely in $\Gamma$, then the rate of injectivity of $\Gamma$ is small.

\begin{lemme}\label{majoration}
Let $\Gamma\subset \Z^n$ be an almost periodic pattern, and $\rho_0\in]0,1[$. If $D(\Gamma)\ge \frac{1}{2-\rho_0}$, then for every $v\in \Z^n$, we have $\rho_\Gamma(v)\ge\rho_0$.
\end{lemme}

\begin{proof}[Proof of Lemma \ref{majoration}]
We argue by contraposition and suppose that there exists $v_0\in\Z^n$ such that $\rho_\Gamma(v_0)\le\rho_0$. We set $\Gamma_{v_0} = \{x\in\Gamma\mid x+v_0\notin\Gamma\}$, thus we have $D(\Gamma_{v_0})\ge(1-\rho_0)D(\Gamma)$. But $\Gamma \cap (\Gamma_{v_0}+v_0) = \emptyset$, so $D(\Gamma) + D(\Gamma_{v_0}+v_0)\le 1$, which implies that $D(\Gamma)(1+1-\rho_0)\le 1$.
\end{proof}

\begin{theoreme}\label{PerLin1}
Let $(A_k)_{k\ge 1}$ be a generic sequence of matrices of $SL_n(\R)$. Then there exists a parameter $\lambda\in ]0,1[$ such that for every $k\ge 1$, we have $\tau^k(A_1,\cdots,A_k) \le (\lambda^k+1)/2$. In particular, $\tau^\infty((A_k)_k)\le 1/2$.
\end{theoreme}

The proof of this theorem is based on an argument of equirepartition, which translates the problem in terms of areas of intersections of cubes.

\begin{lemme}\label{lambda}
Let $\delta>0$ and $M>0$. Then there exists $V_0>0$ such that for all $A\in SL_n(\R)$ with $\|A\|\le M$, there exists $B\in SL_n(\R)$ totally irrational, with $\|B-A\|\le \delta$, there exist a polygon $P\subset \T^n$ whose volume is bigger than $V_0$, and some $v_0\in\Z^n\setminus\{0\}$, such that if $Bx\in P$ mod $\Z^n$, then $\pi(Bx) = \pi(B(x+v_0))$.
\end{lemme}

\begin{proof}[Proof of Lemma \ref{lambda}]
Haj\'os theorem (Theorem~\ref{hajos}, see also Corollary~\ref{CoroHajos}) proves that for $\varep>0$ small enough, the set of matrices $B\in SL_n(\R)$ such that for every $v\in\Z^n\setminus\{0\}$, $\|Bv\|_\infty > 1-\varep$ is a small neighbourhood of a locally finite union of manifolds. Thus, there exists $\varep = \varep(\delta,M)>0$ such that for all $A\in SL_n(\R)$ with $\|A\|\le M$, there exists $B\in SL_n(\R)$ and $v_0\in\Z^n\setminus\{0\}$ such that $\|B-A\|\le \delta$ and $\|Bv_0\|_\infty < 1-\varep$. As totally irrational matrices are dense among $SL_n(\R)$, we can moreover suppose that $B$ is totally irrational. We now set $P = B(0,1/2) \cap B(-v_0,1/2) \subset \R^n/\Z^n$ (we identify $\R^n/\Z^n$ with its fundamental domain $B(0,1/2)$). The volume of $P$ is bigger than $V_0 = \varep^n$. Reasoning as in the proof of Proposition~\ref{2Fr}, we get that if $Bx\in P$ mod $\Z^n$, then $\pi(Bx) = \pi(B(x+v_0))$.
\end{proof}

\begin{proof}[Proof of Theorem \ref{PerLin1}]
We prove the theorem by induction. Suppose that the theorem is proved for a rank $k\in\N$. We set 
\[\Gamma_k = (\widehat {B_k}\circ\cdots\circ \widehat{B_1}) (\Z^n)\qquad \text{and}\qquad \rho_0 = 2-\frac{1}{D(\Gamma_k)}.\]
By Theorem~\ref{imgquasi}, $\Gamma_k$ is an almost periodic pattern.

Lemma~\ref{lambda} applied to $A=A_{k+1}$, $\delta$ and $M = \|(A_i)_i\| + \delta$ gives us a parameter $V_0>0$ (depending only on $\delta$ and $M$), a polygon $P\subset \T^n$ of volume greater than $V_0$, a matrix $B = B_{k+1}$ and a vector $v_0\in\Z^n\setminus\{0\}$. Lemma~\ref{majoration} implies that $\rho_{\Gamma_k}(v_0) \ge \rho_0$. But Lemma~\ref{passifacil} asserts that for a generic set of matrices $B$, the set $B(\Gamma_k)$ is equidistributed modulo $\Z^n$; perturbing a little $B$ if necessary, we suppose that it is true. These facts imply that
\[D\{x\in\Gamma_k \mid Bx\in P\ \operatorname{mod}\ \Z^n,x+v_0\in\Gamma_k\} \ge V_0\rho_0 D(\Gamma_k),\]
and thus, as $Bx\in P\ \operatorname{mod}\ \Z^n$ implies $\widehat B(x) = \widehat B(x+v_0)$, we get
\[D(\widehat{B_{k+1}}(\Gamma_k))\le D(\Gamma_k) - V_0\rho_0 D(\Gamma_k) = D(\Gamma_k)(1-V_0\rho_0).\]
But $\rho_0 = 2-\frac{1}{D(\Gamma)}$, thus (we can suppose that $V_0 \le 1/2$)
\[D(\Gamma_{k+1}) = D(\widehat{B_{k+1}}(\Gamma_k))\le D(\Gamma_k)\left(1-V_0\left(2-\frac{1}{D(\Gamma_k)}\right)\right) = V_0+D(\Gamma_k)(1-2V_0).\]
Setting $\lambda = 1-2V_0$, we obtain the general term of an arithmetico-geometric sequence, this proves the theorem for the rank $k+1$.
\end{proof}

\section{Diffusion process and the case of isometries}

Unfortunately, problems arise when we try to perturb a sequence of matrices to make its asymptotic rate smaller than $1/2$. First of all, if the density of an almost periodic pattern $\Gamma$ is smaller than $1/2$, then the set of differences may not be the full set $\Z^n$. Even worse, it might happen that this set of differences has big holes, as shown by the following example. We take the almost periodic pattern
\[\Gamma = \bigcup_{i=0}^{39} 100\Z+i.\]
For every $x\in\Gamma$ and $v\in\llbracket 40, 60\rrbracket$, $x+v\notin\Gamma$. In other words, for every $v\in\llbracket 40, 60\rrbracket$, we have $\rho_\Gamma(v) = 0$, whereas $D(\Gamma) = 0.4$. However, the things are not too bad for the frequencies of differences when we are close to 0, as shown by the Minkowski-like theorem for almost periodic patterns (Theorem \ref{MinkAlm}).

\paragraph{Diffusion process} In this paragraph, we study the action of a discretization of a matrix on the set of differences of an almost periodic pattern $\Gamma$; more precisely, we study the link between the functions $\rho_\Gamma$ and $\rho_{\widehat A(\Gamma)}$.

For $u\in\R^n$, we define the function $\varphi_u$, which is a ``weighted projection'' of $u$ on $\Z^n$.

\begin{definition}
Given $u\in\R^n$, the function $\varphi_u = \Z^n\to [0,1]$ is defined by
\[\varphi_u (v) = \left\{\begin{array}{ll}
0 & \ \text{if}\ d_\infty(u,v)\ge 1\\
\prod_{i=1}^n (1-|u_i+v_i|) & \ \text{if}\ d_\infty(u,v)< 1.
\end{array}\right.\]
\end{definition}

\begin{figure}
\begin{minipage}[b]{.4\linewidth}
\begin{center}
\begin{tikzpicture}[scale=.85]
\draw (0,0) rectangle (3,3);
\draw (1.3,0.9) node {$\times$};
\draw (1.3,0.9) node[above left] {$u$};
\draw[color=green!40!black] (0,0) node {$\bullet$};
\draw (0,0) node[below right] {$v$};
\draw[color=green!40!black] (0,3) node {$\bullet$};
\draw (0,3) node[above right] {$v+(0,1)$};
\draw[color=green!40!black] (3,0) node {$\bullet$};
\draw (3,0) node[below right] {$v+(1,0)$};
\draw[color=green!40!black] (3,3) node {$\bullet$};
\draw (3,3) node[above right] {$v+(1,1)$};
\draw[->] (1.35,.75) to[bend right] (2.9,.1);
\draw[->] (1.35,1.05) to[bend left] (2.9,2.9);
\draw[->] (1.25,.75) to[bend left] (.1,.1);
\draw[->] (1.25,1.05) to[bend right] (.1,2.9);
\draw[->, color=red!60!black] (.1,.1) to (1.2,0.8);
\draw[color=blue!80!black] (1.3,0) -- (1.3,3);
\draw[color=blue!80!black] (0,.9) -- (3,.9);
\end{tikzpicture}
\end{center}
\caption[The function $\varphi_u$]{The function $\varphi_u$ in dimension 2: its value on one vertex of the square is equal to the area of the opposite rectangle; in particular, $\varphi_u(v)$ is the area of the rectangle with the vertices $u$ and $v+(1,1)$.}\label{IxMa}
\end{minipage}\hfill
\begin{minipage}[b]{.55\linewidth}
\begin{center}
\begin{tikzpicture}[scale=.8]
\draw[color=blue!80!black] (1.7,0) -- (1.7,3);
\draw[color=blue!80!black] (0,2.1) -- (3,2.1);
\draw[color=blue!80!black] (0,0) rectangle (3,3);
\draw (1.7,2.1) node {$\times$};
\draw (1.7,2.1) node[above left] {$u'$};

\draw[dashed] (0,0) rectangle (6,6);
\draw[dashed] (3,0) -- (3,6);
\draw[dashed] (0,3) -- (6,3);
\draw (1.5,1.5) node[below left] {\small$(0,0)$};
\draw (4.5,1.5) node[below right] {\small$(1,0)$};
\draw (1.5,4.5) node[above left] {\small$(0,1)$};
\draw (4.5,4.5) node[above right] {\small$(1,1)$};
\draw[->, color=red!60!black] (1.8,2.2) to (2.9,2.9);
\draw[color=green!40!black] (1.5,1.5) node {$\bullet$};
\draw[color=green!40!black] (1.5,4.5) node {$\bullet$};
\draw[color=green!40!black] (4.5,1.5) node {$\bullet$};
\draw[color=green!40!black] (4.5,4.5) node {$\bullet$};
\end{tikzpicture}
\end{center}
\caption[Proof of Proposition~\ref{ActionDiff}]{The red vector is equal to that of Figure~\ref{IxMa} for $u=Av$. If $Ax$ belongs to the bottom left rectangle, then $\pi(Ax+Av) = y\in\Z^2$; if $Ax$ belongs to the top left rectangle, then $\pi(Ax+Av) = y+(0,1)$ etc.}\label{IxMa2}
\end{minipage}
\end{figure}

In particular, the function $\varphi_u$ satisfies $\sum_{v\in\Z^n} \varphi_u(v) = 1$, and is supported by the vertices of the integral unit cube\footnote{By definition, an integral cube has its faces parallel to the canonical hyperplanes.} that contains\footnote{More precisely, the support of $\varphi_u$ is the smallest integral unit cube of dimension $n'\le n$ which contains $u$.} $u$. Figure~\ref{IxMa} gives a geometric interpretation of this function $\varphi_u$.

The following property asserts that the discretization $\widehat A$ acts ``smoothly'' on the frequency of differences. In particular, when $D(\Gamma) = D(\widehat A \Gamma)$, the function $\rho_{\widehat A \Gamma}$ is obtained from the function $\rho_\Gamma$ by applying a linear operator $\mathcal A$, acting on each Dirac function $\delta_v$ such that $\mathcal A \delta_u(v) = \varphi_{A(u)}(v)$. Roughly speaking, to compute $\mathcal A \delta_v$, we take $\delta_{Av}$ and apply a diffusion process. In the other case where $D(\widehat A \Gamma) < D(\Gamma)$, we only have inequalities involving the operator $\mathcal A$ to compute the function $\rho_{\widehat A\Gamma}$.

\begin{prop}\label{ActionDiff}
Let $\Gamma\subset \Z^n$ be an almost periodic pattern and $A\in SL_n(\R)$ be a generic matrix.
\begin{enumerate}[(i)]
\item If $D(\widehat A(\Gamma)) = D(\Gamma)$, then for every $u\in\Z^n$,
\[\rho_{\widehat A(\Gamma)}(u) = \sum_{v\in\Z^n} \varphi_{A(v)} (u) \rho_\Gamma(v).\]
\item In the general case, for every $u\in\Z^n$, we have
\[\frac{D(\Gamma)}{D(\widehat A(\Gamma))}\sup_{v\in\Z^n} \varphi_{A(v)} (u) \rho_\Gamma(v) \le \rho_{\widehat A(\Gamma)}(u) \le \frac{D(\Gamma)}{D(\widehat A(\Gamma))}\sum_{v\in\Z^n} \varphi_{A(v)} (u) \rho_\Gamma(v).\]
\end{enumerate}
\end{prop}

\begin{rem}
This proposition expresses that the action of the discretization of a linear map $A$ on the differences is more or less that of a multivalued map.
\end{rem}

\begin{proof}[Proof of Proposition \ref{ActionDiff}]
We begin by proving the first point of the proposition. Suppose that $P\in O_n(\R)$ is generic and that $D(\widehat A(\Gamma)) = D(\Gamma)$. Let $x\in \Gamma\cap (\Gamma-v)$. We consider the projection $y'$ of $y=Px$, and the projection $u'$ of $u=Pv$, on the fundamental domain $]-1/2,1/2]^n$ of $\R^n/\Z^n$. We have
\[P(x+v) = \pi(Px) + \pi(Pv) + y' + u'.\]
Suppose that $y'$ belongs to the parallelepiped whose vertices are $(-1/2,\cdots,-1/2)$ and $u'$ (in bold in Figure~\ref{IxMa2}), then $y'+u'\in [-1/2,1/2[^n$. Thus, $\pi(P(x+v)) = \pi(Px) + \pi(Pv)$. The same kind of results holds for the other parallelepipeds whose vertices are $u'$ and one vertex of $[-1/2,1/2[^n$.

We set $\Gamma = \widehat A(\Z^n)$. The genericity of $P$ ensures that for every $v\in\Z^n$, the set $\Gamma\cap (\Gamma-v)$, which has density $D(\Gamma)\rho_\Gamma(v)$ (by definition of $\rho_\Gamma$), is equidistributed modulo $\Z^n$ (by Lemma~\ref{passifacil}). Thus, the points $x'$ are equidistributed modulo $\Z^n$. In particular, the difference $v$ will spread into the differences which are the support of the function $\varphi_{Pv}$, and each of them will occur with a frequency given by $\varphi_{Pv} (x) \rho_\Gamma(v)$. The hypothesis about the fact that the density of the sets does not decrease imply that the contributions of each difference of $\Gamma$ to the differences of $\widehat A (\Gamma)$ add.

In the general case, the contributions to each difference of $\Gamma$ may overlap. However, applying the argument of the previous case, we can easily prove the second part of the proposition.
\end{proof}

\begin{rem}\label{RemActionDiff}
We also remark that:
\begin{enumerate}[(i)]
\item the density strictly decreases (that is, $D(\widehat A(\Gamma)) < D(\Gamma)$) if and only if there exists $v_0\in\Z^n$ such that $\rho_\Gamma(v_0)>0$ and $\|Av_0\|_\infty <1$;
\item if there exists $v_0\in\Z^n$ such that 
\[\sum_{v\in\Z^n} \varphi_{A(v)} (v_0) \rho_\Gamma(v_0)>1,\]
then the density strictly decreases by at least $\sum_{v\in\Z^n} \varphi_{A(v)} (v_0) \rho_\Gamma(v_0)- 1$;
\item we can compute which differences will go to the difference 0, that is, the differences $u\in\Z^n$ such that there exists $x\in\Gamma\cap(\Gamma-u)$ such that $\widehat A(x) = \widehat A(x+u)$ (that will make the rate of injectivity decrease). The set of such differences $u$ is $A^{-1}(B(0,1))\cap \Z^n$ (recall that $B(0,1)$ is an infinite ball). By Haj\'os theorem (Theorem~\ref{hajos}, see also Corollary~\ref{CoroHajos}), this set contains a non zero vector when $A$ is generic. More generally, the set of differences $u\in \Gamma\cap(\Gamma-u)$ such that $\widehat A(x) + v = \widehat A(x+u)$ is $A^{-1}(B(v,1))\cap \Z^n$. Iterating this process, it is possible to compute which differences will go to 0 in time $t$.
\end{enumerate}
\end{rem}

\paragraph{Rate of injectivity of a generic sequence of isometries} Proposition~\ref{ActionDiff} allows to give an alternative proof of Theorem~\ref{ConjPrincip} for isometries.

\begin{theoreme}\label{AnswerConjIsom}
Let $(P_k)_{k\ge 1}$ be a generic sequence of matrices of $O_n(\R)$. Then $\tau^\infty((P_k)_k) = 0$.
\end{theoreme}

\begin{figure}[ht]
\begin{minipage}[c]{.33\linewidth}
	\includegraphics[width=\linewidth, trim = 1.5cm .5cm 1.5cm .5cm]{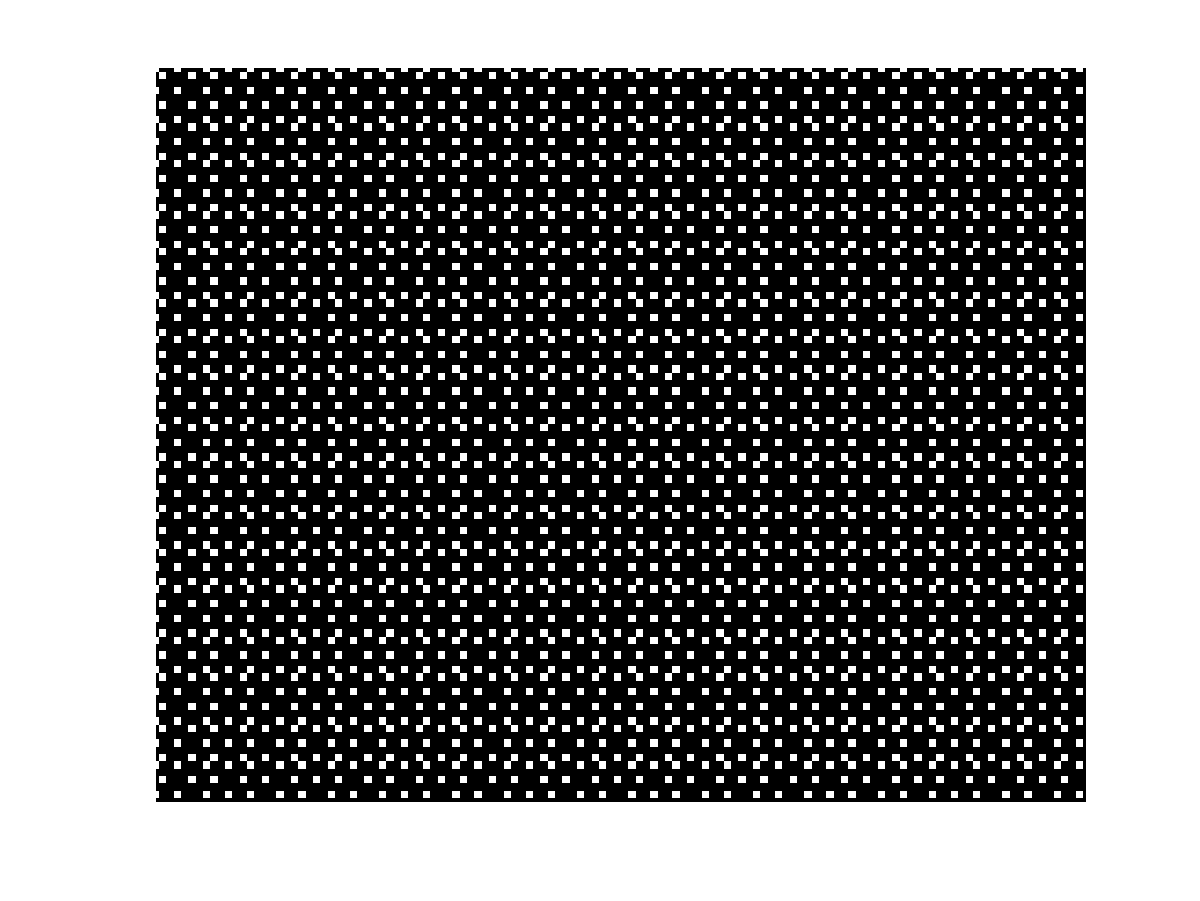}
\end{minipage}\hfill
\begin{minipage}[c]{.33\linewidth}
	\includegraphics[width=\linewidth, trim = 1.5cm .5cm 1.5cm .5cm]{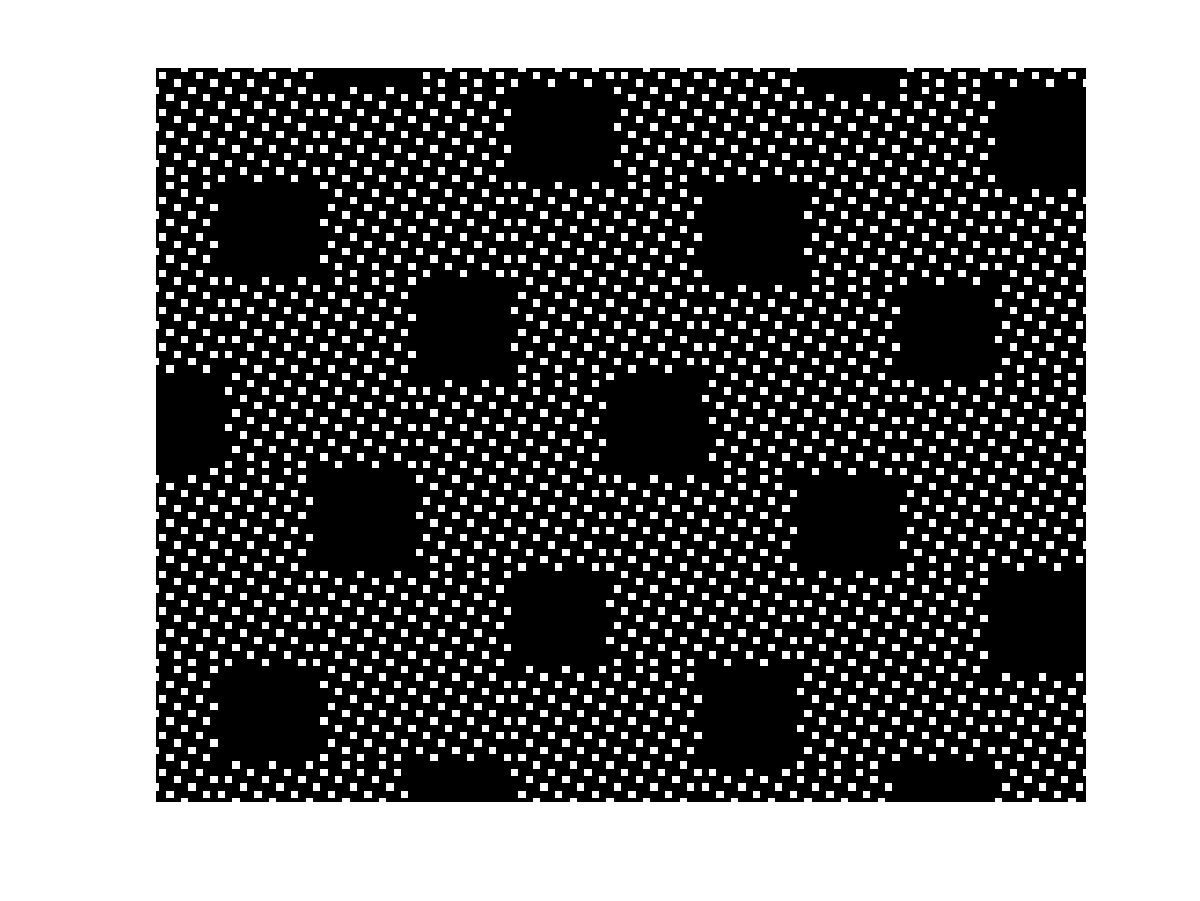}
\end{minipage}\hfill
\begin{minipage}[c]{.33\linewidth}
	\includegraphics[width=\linewidth, trim = 1.5cm .5cm 1.5cm .5cm]{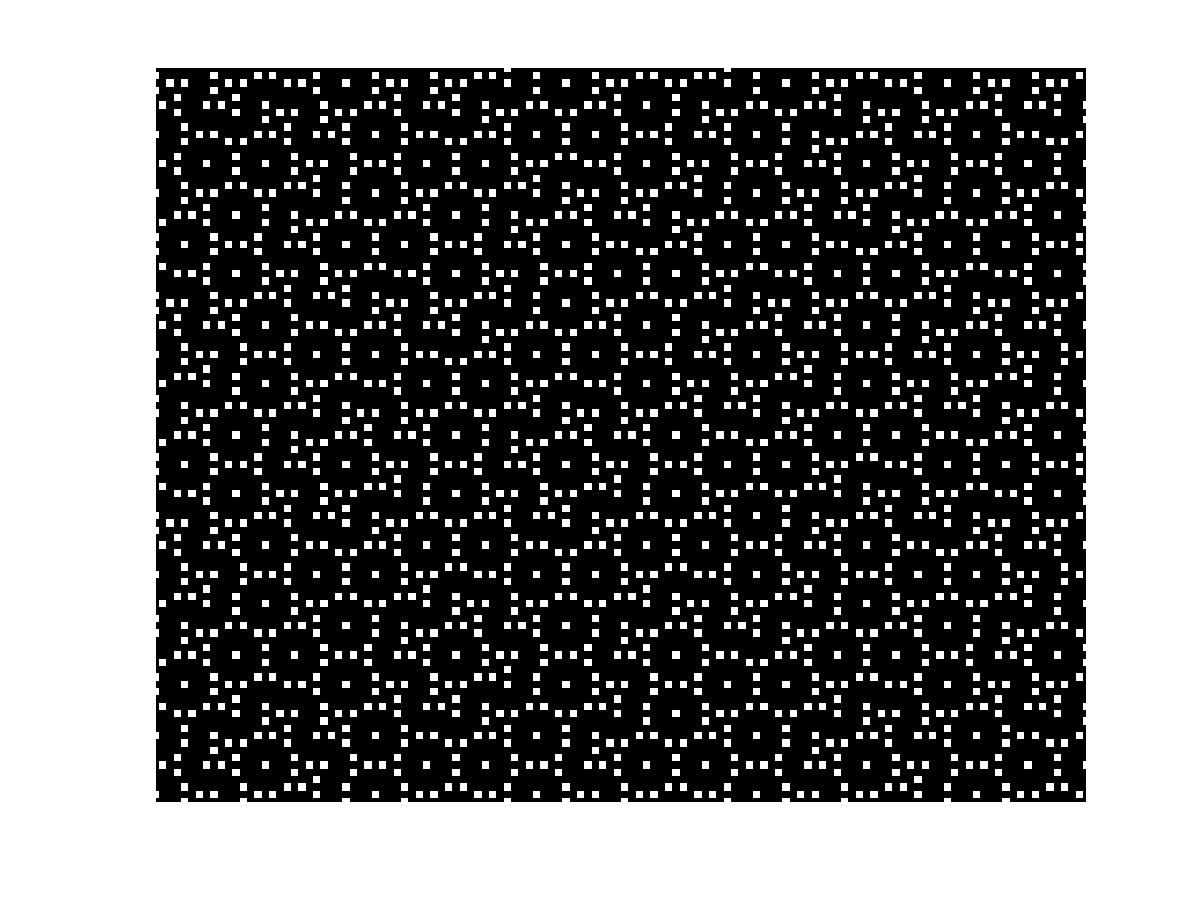}
\end{minipage}

\begin{minipage}[c]{.33\linewidth}
	\includegraphics[width=\linewidth, trim = 1.5cm .5cm 1.5cm .5cm]{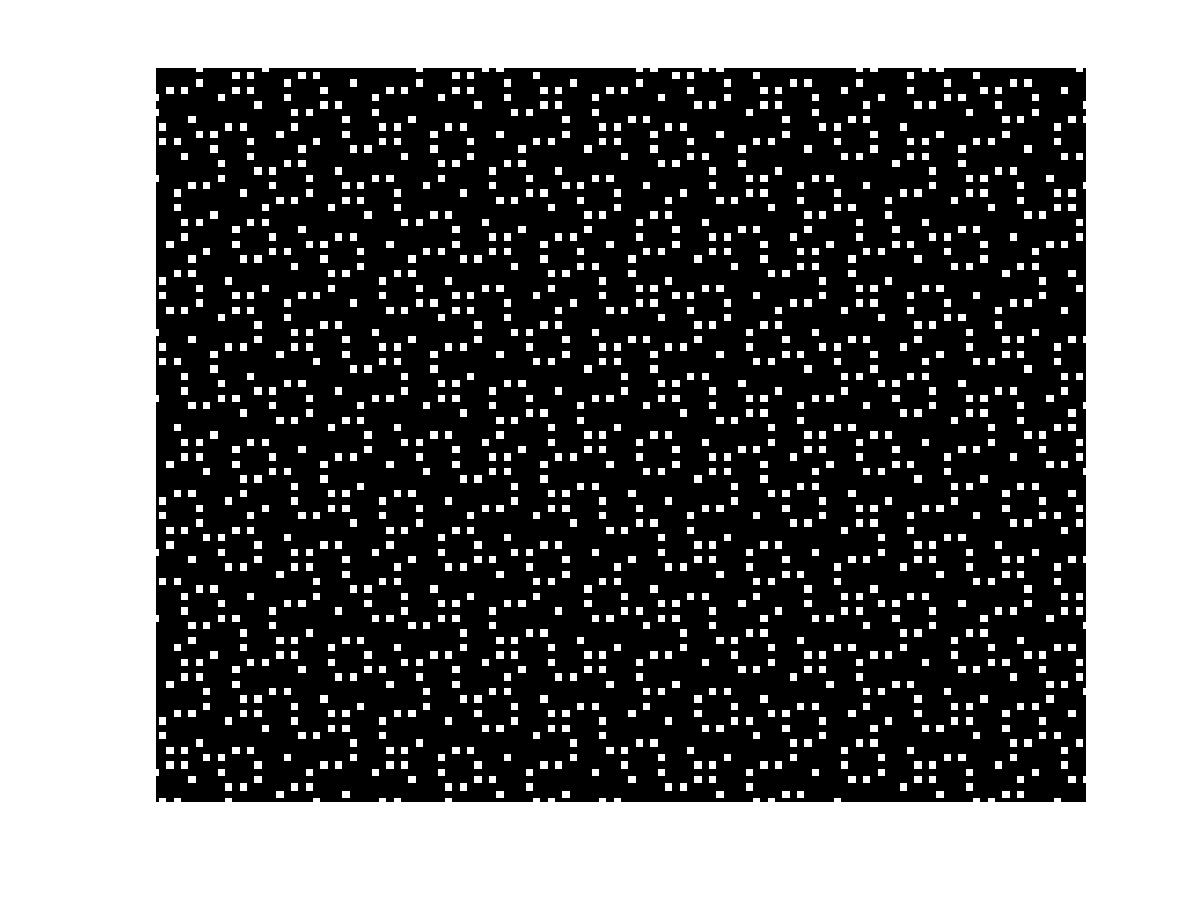}
\end{minipage}\hfill
\begin{minipage}[c]{.33\linewidth}
	\includegraphics[width=\linewidth, trim = 1.5cm .5cm 1.5cm .5cm]{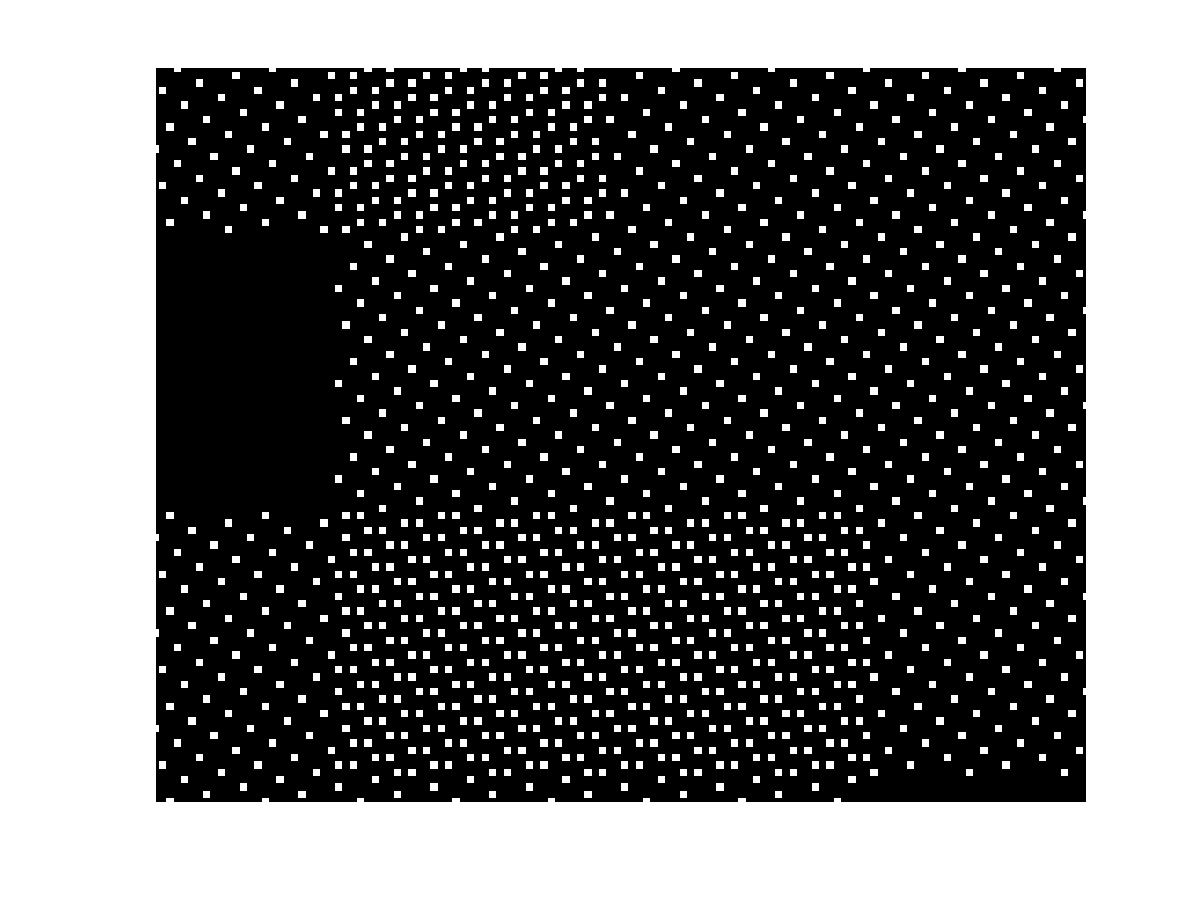}
\end{minipage}\hfill
\begin{minipage}[c]{.33\linewidth}
	\includegraphics[width=\linewidth, trim = 1.5cm .5cm 1.5cm .5cm]{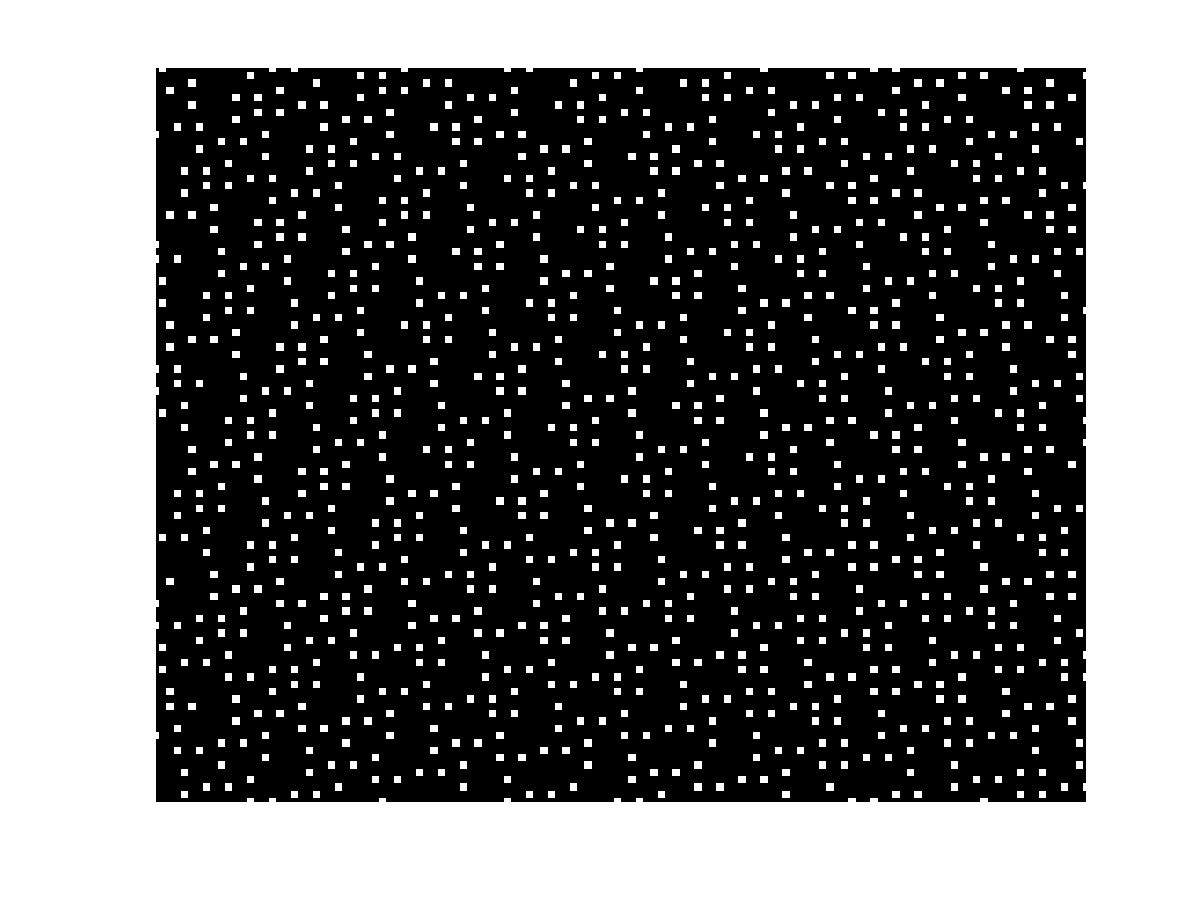}
\end{minipage}
\caption[Images of $\Z^2$ by discretizations of rotations]{Images of $\Z^2$ by discretizations of rotations, a point is black if it belongs to the image of $\Z^2$ by the discretization of the rotation. From left to right and top to bottom, angles $\pi/4$, $\pi/5$, $\pi/6$, $\pi/7$, $\pi/8$ and $\pi/9$.}\label{ImagesRotations}
\end{figure}

\begin{figure}[ht]
\begin{minipage}[c]{.33\linewidth}
	\includegraphics[width=\linewidth, trim = 1.5cm .5cm 1.5cm .5cm]{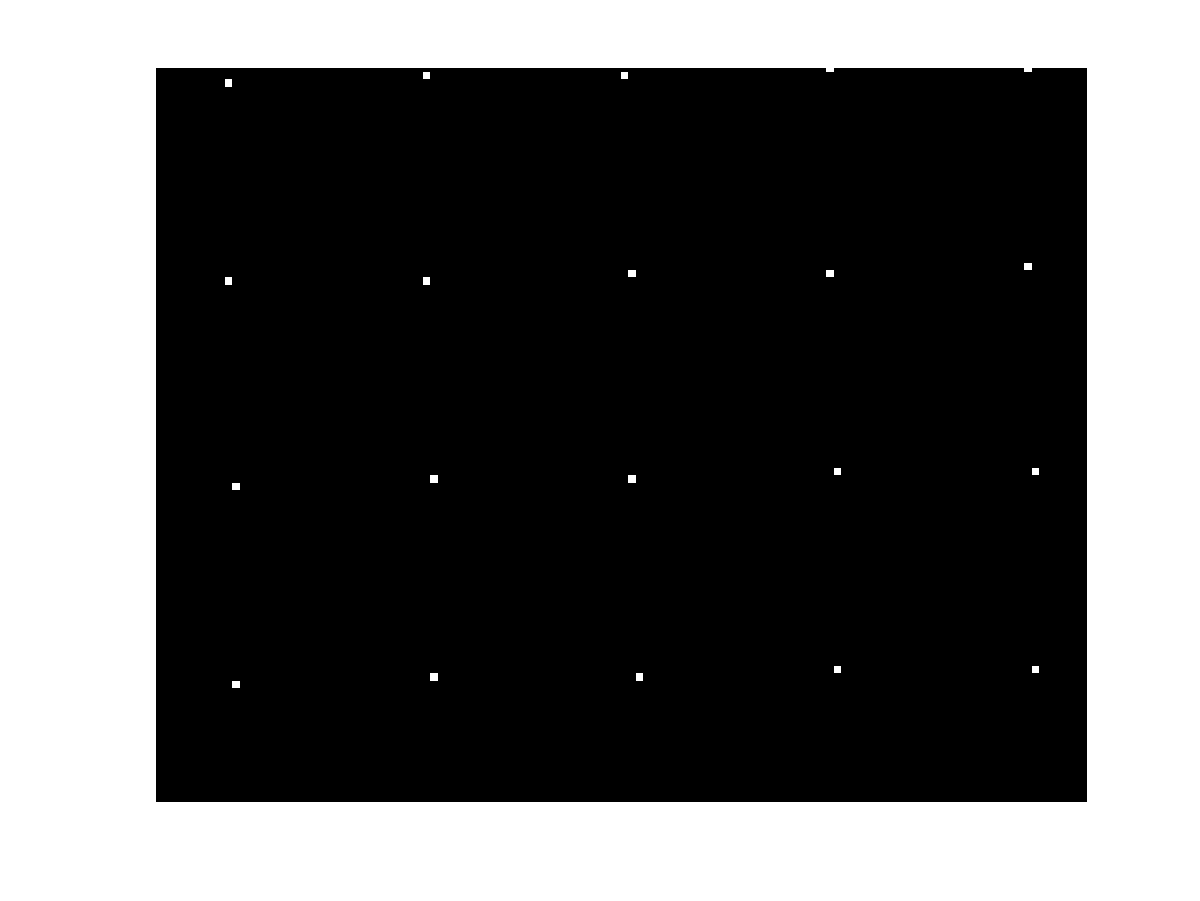}
\end{minipage}\hfill
\begin{minipage}[c]{.33\linewidth}
	\includegraphics[width=\linewidth, trim = 1.5cm .5cm 1.5cm .5cm]{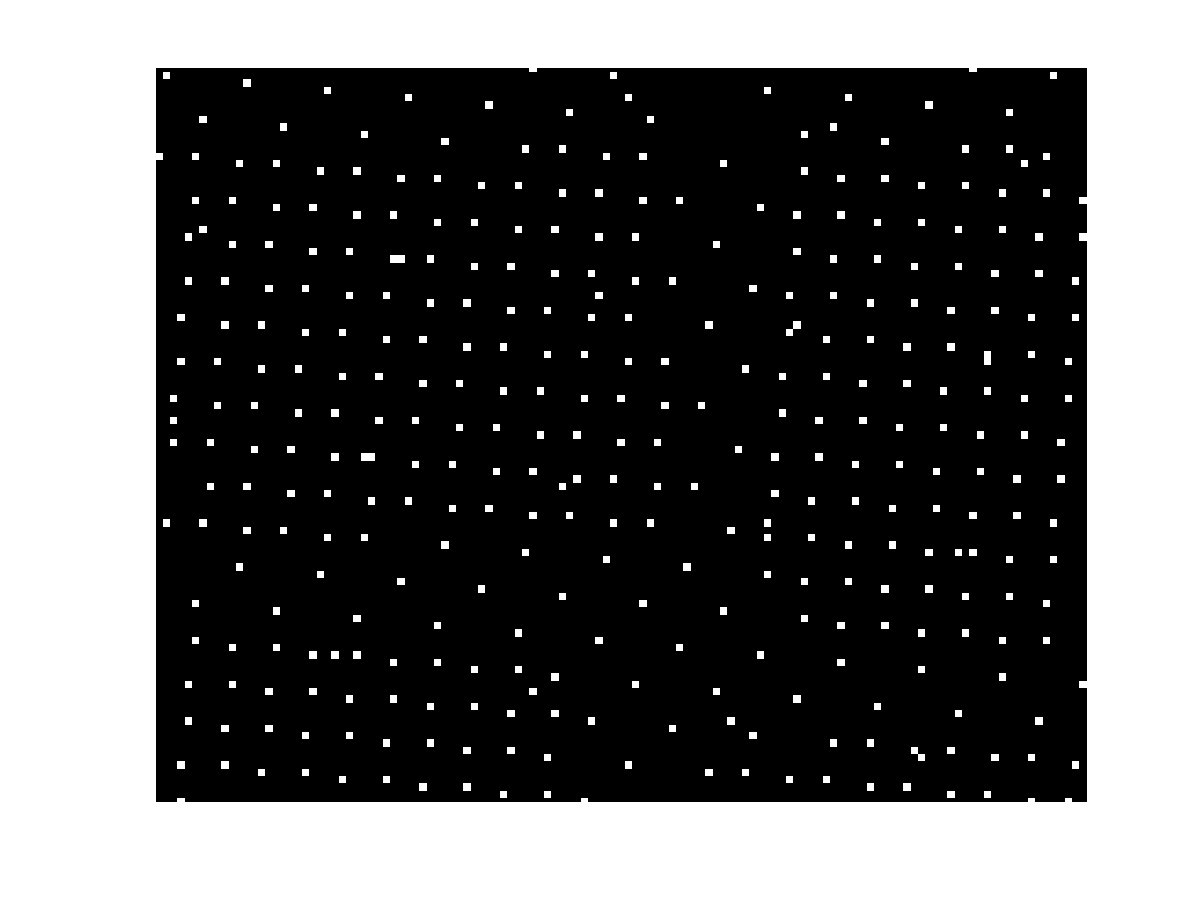}
\end{minipage}\hfill
\begin{minipage}[c]{.33\linewidth}
	\includegraphics[width=\linewidth, trim = 1.5cm .5cm 1.5cm .5cm]{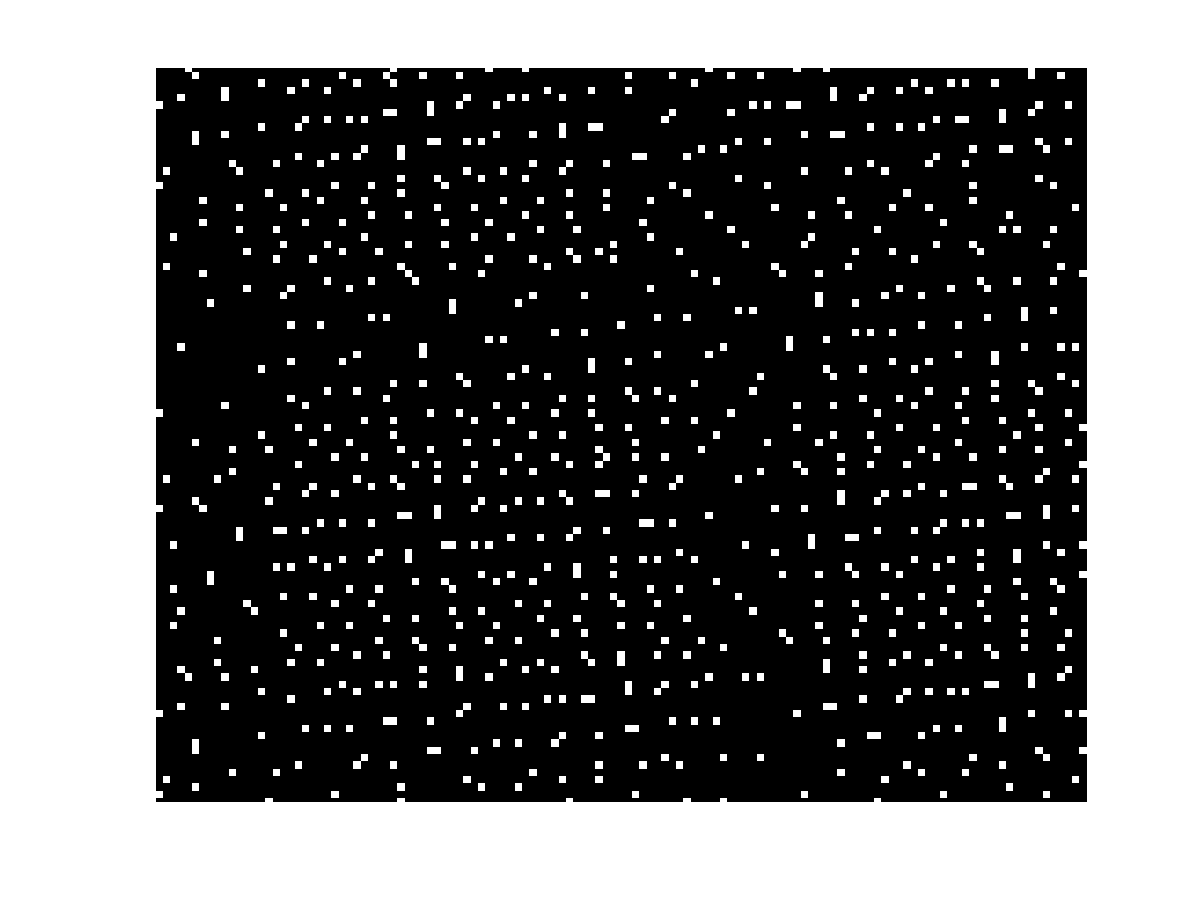}
\end{minipage}

\begin{minipage}[c]{.33\linewidth}
	\includegraphics[width=\linewidth, trim = 1.5cm .5cm 1.5cm .5cm]{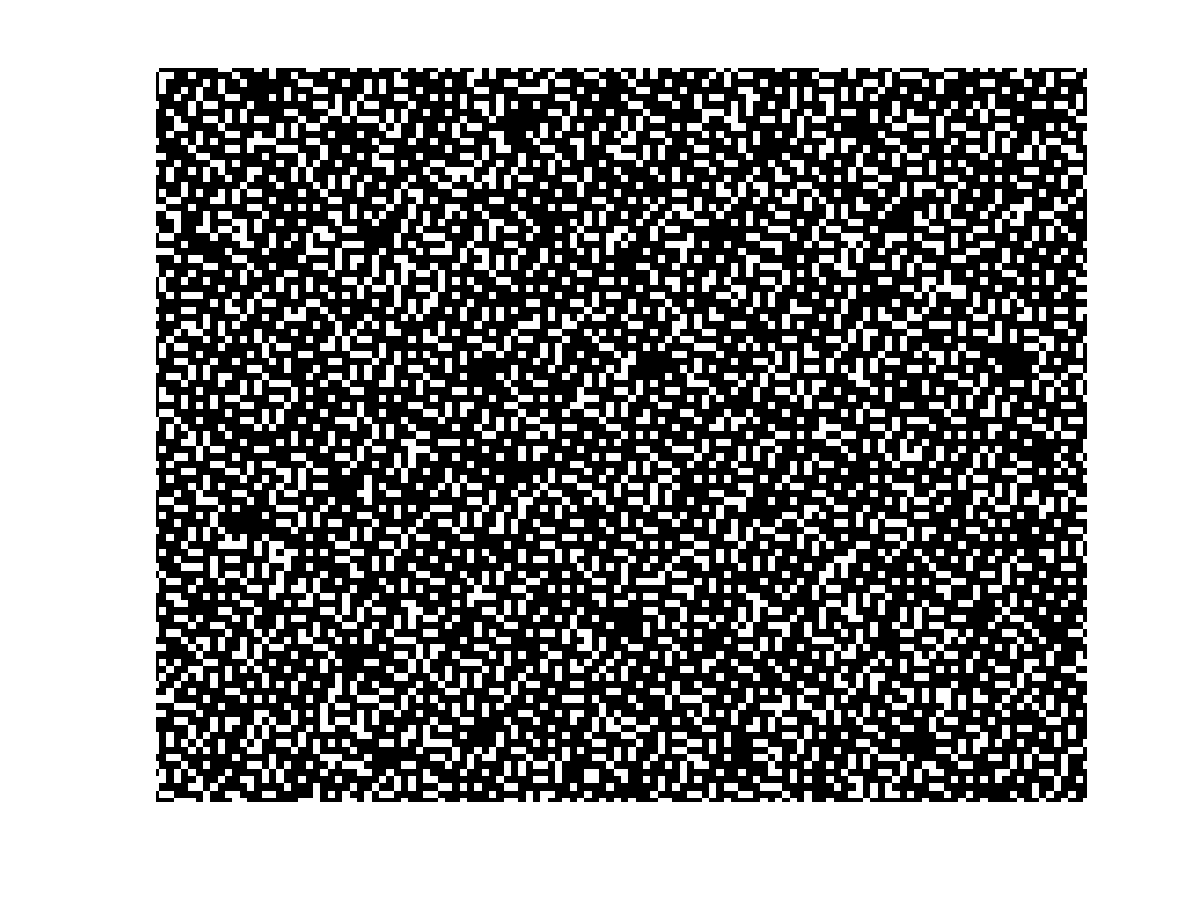}
\end{minipage}\hfill
\begin{minipage}[c]{.33\linewidth}
	\includegraphics[width=\linewidth, trim = 1.5cm .5cm 1.5cm .5cm]{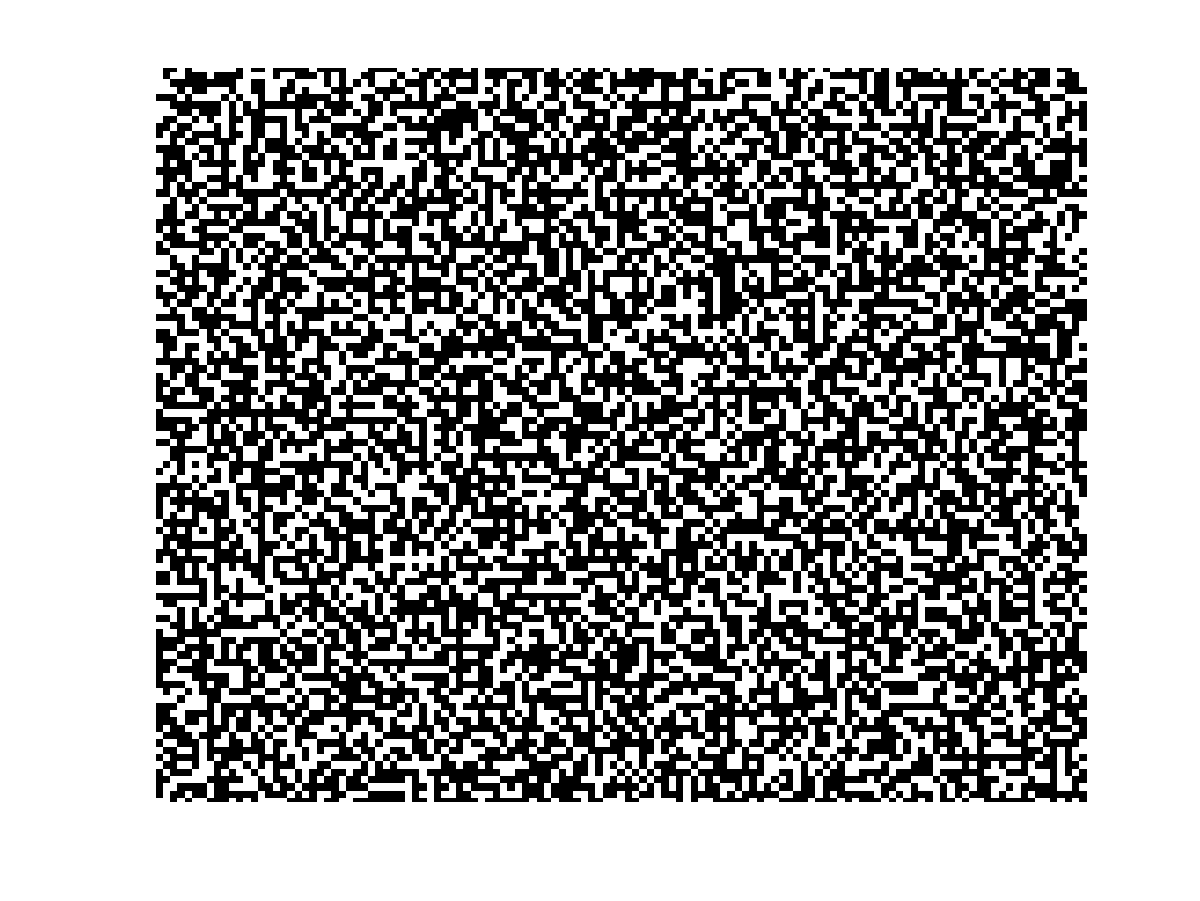}
\end{minipage}\hfill
\begin{minipage}[c]{.33\linewidth}
	\includegraphics[width=\linewidth, trim = 1.5cm .5cm 1.5cm .5cm]{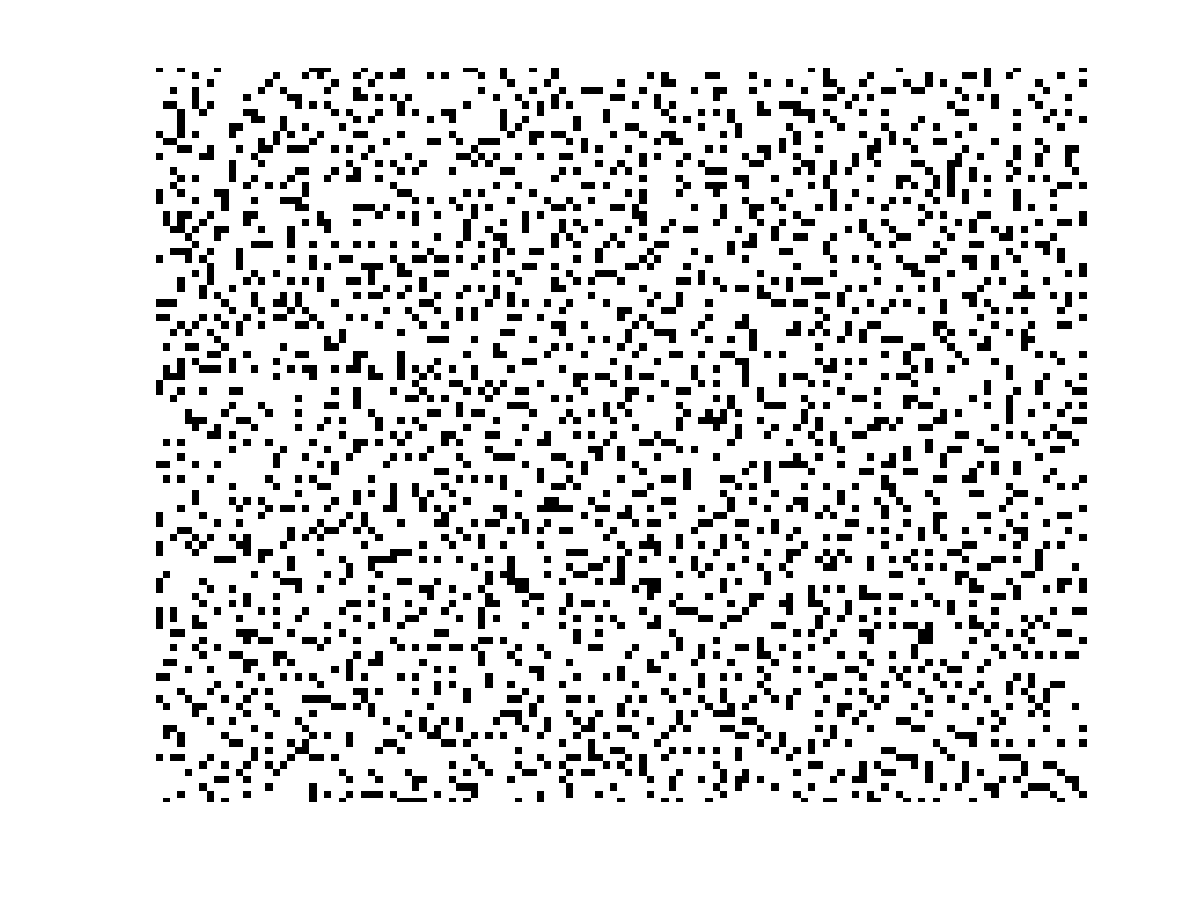}
\end{minipage}
\caption[Successive images of $\Z^2$ by discretizations of random rotations]{Successive images of $\Z^2$ by discretizations of random rotations, a point is black if it belongs to $(\widehat{R_{\theta_k}}\circ\cdots\circ\widehat{R_{\theta_1}})(\Z^2)$, where the $\theta_i$ are chosen uniformly randomly in $[0,2\pi]$. From left to right and top to bottom, $k=1,\, 2,\, 3,\, 5,\, 10,\, 50$.}\label{ImagesSuitesRotations}
\end{figure}

\begin{figure}[t]
\begin{center}
\includegraphics[width=.6\linewidth]{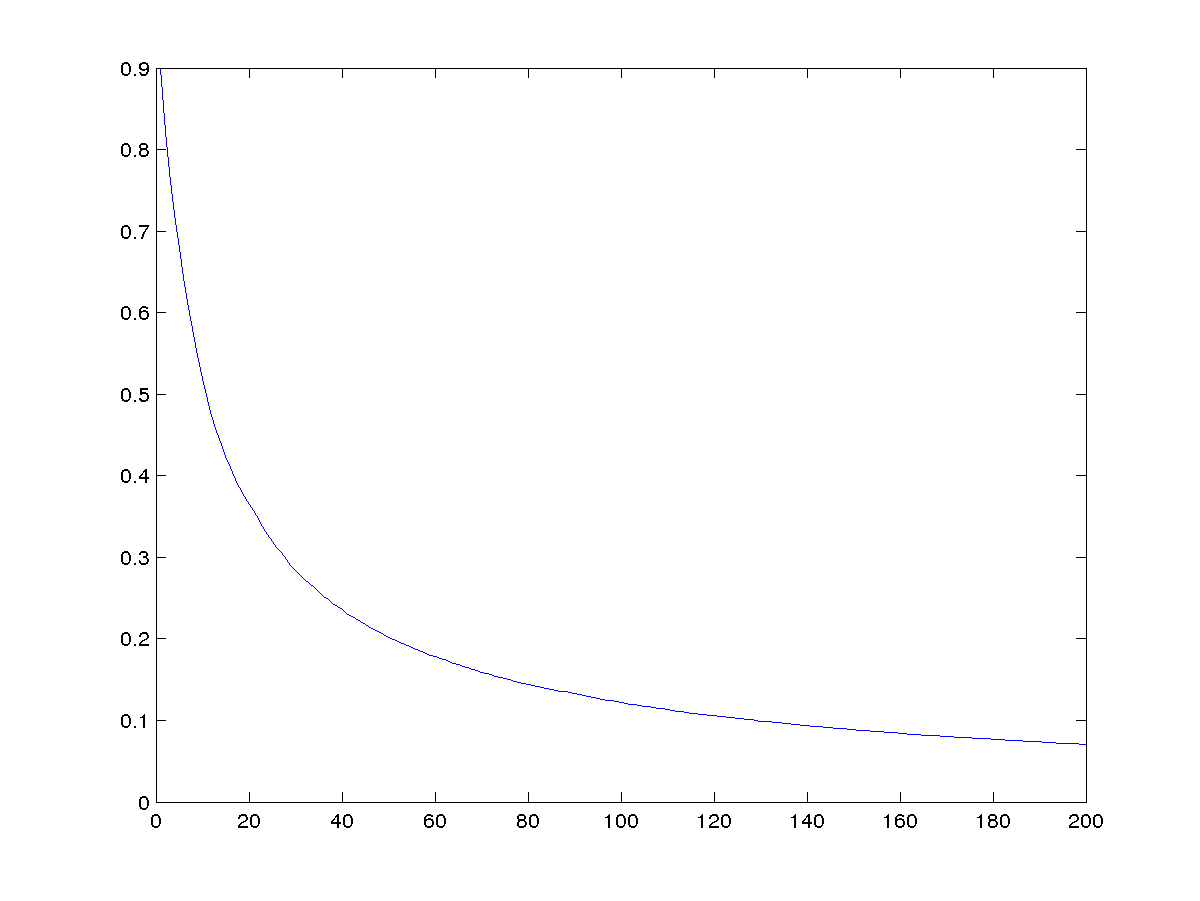}
\caption[Expectation of the rate of injectivity of a random sequences of rotations]{Expectation of the rate of injectivity of a random sequences of rotations: the graphic represents the mean of the rate of injectivity $\tau^k(R_{\theta_k},\cdots,R_{\theta_1})$ depending on $k$, $1\le k\le 200$, for 50 random draws of sequences of angles $(\theta_i)_i$, with each $\theta_i$ chosen independently and uniformly in $[0,2\pi]$. Note that the behaviour is not exponential.}\label{TauxSuiteRotations}
\end{center}
\end{figure}


In the previous section, the starting point of the proof was Lemma~\ref{majoration}, which ensures that when the rate of injectivity is bigger than $1/2$, the frequency of \emph{any} difference is bigger than a constant depending on the rate. Here, the starting point is the Minkowski theorem for almost periodic patterns (Theorem~\ref{MinkAlm}), which gives \emph{one} nonzero difference whose frequency is positive. The rest of the proof of Theorem~\ref{AnswerConjIsom} consists in using again an argument of equidistribution. More precisely, we apply successively the following lemma, which asserts that given an almost periodic pattern $\Gamma$ of density $D_0$, a sequence of isometries and $\delta>0$, then, perturbing each isometry of at most $\delta$ if necessary, we can make the density of the $k_0$-th image of $\Gamma$ smaller than $\lambda_0 D_0$, with $k_0$ and $\lambda_0$ depending only on $D_0$ and $\delta$. The proof of this lemma involves the study of the action of the discretizations on differences made in Proposition~\ref{ActionDiff}

\begin{lemme}\label{EstimPerteRot}
Let $(P_k)_{k\ge 1}$ be a sequence of matrices of $O_n(\R)$ and $\Gamma\subset\Z^n$ an almost periodic pattern. Given $\delta>0$ and $D>0$ such that $D(\Gamma)\ge D$, there exists $k_0 = k_0(D)$ (decreasing in $D$), $\lambda_0 = \lambda_0(D,\delta)<1$ (decreasing in $D$ and in $\delta$), and a sequence $(Q_k)_{k\ge 1}$ of totally irrational matrices of $O_n(\R)$, such that $\|P_k - Q_k\|\le \delta$ for every $k\ge 1$ and
\[D\big((\widehat{Q_{k_0}}\circ\dots\circ \widehat{Q_1})(\Gamma)\big)< \lambda_0 D.\]
\end{lemme}

We begin by proving that this lemma implies Theorem~\ref{AnswerConjIsom}.

\begin{proof}[Proof of Theorem \ref{AnswerConjIsom}]
Suppose that Lemma~\ref{EstimPerteRot} is true. Let $\tau_0\in ]0,1[$ and $\delta>0$. We want to prove that we can perturb the sequence $(P_k)_k$ into a sequence $(Q_k)_k$ of isometries, which is $\delta$-close to $(P_k)_k$ and is such that its asymptotic rate is smaller than $\tau_0$ (and that this remains true on a whole neighbourhood of these matrices).

Thus, we can suppose that $\tau^\infty((P_k))>\tau_0$. We apply Lemma~\ref{EstimPerteRot} to obtain the parameters $k_0 = k_0(\tau_0/2)$ (because $k_0(D)$ is decreasing in $D$) and $\lambda_0 = \lambda_0(\tau_0/2,\delta)$ (because $\lambda_0(D,\delta)$ is decreasing in $D$). Applying the lemma $\ell$ times, this gives a sequence $(Q_k)_k$ of isometries, which is $\delta$-close to $(P_k)_k$, such that, as long as $\tau^{\ell k_0}(Q_0,\cdots,Q_{\ell k_0})> \tau_0/2$, we have $\tau^{\ell k_0}(Q_1,\cdots,Q_{\ell k_0})<\lambda_0^\ell D(\Z^n)$. But for $\ell$ large enough, $\lambda_0^\ell<\tau_0$, which proves the theorem.
\end{proof}

\begin{proof}[Proof of Lemma \ref{EstimPerteRot}]
The idea of the proof is the following. Firstly, we apply the Minkowski-type theorem for almost periodic patterns (Theorem~\ref{MinkAlm}) to find a uniform constant $C>0$ and a point $u_0\in\Z^n\setminus\{0\}$ whose norm is ``not too big'', such that $\rho_\Gamma (u_0) > C D(\Gamma)$. Then, we apply Proposition \ref{ActionDiff} to prove that the difference $u_0$ in $\Gamma$ eventually goes to 0; that is, that there exists $k_0\in\N^*$ and an almost periodic pattern $\widetilde\Gamma$ of positive density (that we can compute) such that there exists a sequence $(Q_k)_k$ of isometries, with $\|Q_i-P_i\|\le\delta$, such that for every $x\in \widetilde\Gamma$,
\[(\widehat{Q_{k_0}}\circ\dots\circ \widehat{Q_1})(x) = (\widehat{Q_{k_0}}\circ\dots\circ \widehat{Q_1})(x+u_0).\]
\bigskip

We begin by applying the Minkowsky-like theorem for almost periodic patterns (Theorem \ref{MinkAlm}) to a \emph{Euclidean} ball $B'_R$\index{$B'_R$} with radius $R$ (recall that $[B]$ denotes the set of integer points inside $B$):
\begin{equation}\label{eqMinkChap2}
\frac{1}{\card[B'_R]}\sum_{u\in [B'_R]} \rho_\Gamma(u) \ge D(\Gamma) \frac{\card[B'_{R/2}]}{\card[B'_R]}.
\end{equation}
An easy estimation of the number of integer points inside of the balls $B'_R$ and $B'_{R/2}$ leads to\footnote{See for example \cite[page 5]{MR673938}, this book also performs a complete investigation on the subject of finding the number of integer points in a Euclidean ball.}:
\[\frac{\card[B'_{R/2}]}{\card[B'_R]} \ge \left( \frac{\sqrt{R/2} - \sqrt n/2}{\sqrt{R} + \sqrt n/2}\right)^{2n},\]
in particular, if $R\ge n$, we obtain
\[\frac{\card[B'_{R/2}]}{\card[B'_R]} \ge \left(\frac{\sqrt 2 -1}{3}\right)^{2n} \ge \frac{1}{53^n},\]
thus Equation~\eqref{eqMinkChap2} becomes
\begin{equation}\label{eqMinkChap2b}
\frac{1}{\card[B'_R]}\sum_{u\in [B'_R]} \rho_\Gamma(u) \ge D(\Gamma) \frac{1}{53^n}.
\end{equation}

We now want that $\rho_\Gamma(0)=1$ plays only a little role in Equation~\eqref{eqMinkChap2b}, that is
\[\frac{1}{\card[B'_R]} \le \frac{D(\Gamma)}{2\times 53^n},\]
using the same estimate as previously for $\card[B'_R]$, that is true if
\[\frac{1}{(\sqrt{R}-\sqrt{n}/2)^n} \le \frac{D(\Gamma)}{2\times 53^n};\]
a direct calculation shows that it is true for example if
\begin{equation}\label{superEq}
R = R_0 = n+\frac{212}{\sqrt[n]{D(\Gamma)}}.
\end{equation}
In this case, Equation~\eqref{eqMinkChap2b} implies
\[\frac{1}{\card[B'_R]}\sum_{u\in [B'_R]\setminus\{0\}} \rho_\Gamma(u) \ge D(\Gamma) \frac{1}{53^n} - \frac{\rho_\Gamma(0)}{\card[B'_R]} \ge \frac{D(\Gamma)}{2\times 53^n},\]
thus there exists $u_0\in [B'_R]\setminus \{0\}$ such that
\begin{equation}\label{EqDefU0}
\rho_\Gamma(u_0) \ge \frac{D(\Gamma)}{2\times 53^n}.
\end{equation}

We now take $R\ge R_0$, and perturb each matrix $P_k$ into a totally irrational matrix $Q_k$ such that for every point $x\in [B'_R]\setminus\{0\}$, the point $Q_k (x)$ is far away from the lattice $\Z^n$. More precisely, as the set of matrices $Q\in O_n(\R)$ such that $Q([B'_R]) \cap \Z^n \neq \{0\}$ is finite, there exists a constant $d_0(R,\delta)$ such that for every $P\in O_n(\R)$, there exists $Q\in O_n(\R)$ such that $\|P-Q\|\le\delta$ and for every $x\in [B'_R]\setminus\{0\}$, we have $d_\infty(Q(x),\Z^n)> d_0(R,\delta)$. Applying Lemma~\ref{passifacil} (which states that if the sequence $(Q_k)_k$ is generic, then the matrices $Q_k$ are ``non resonant''), we build a sequence $(Q_k)_{k\ge 1}$ of totally irrational matrices of $O_n(\R)$ such that for every $k\in\N^*$, we have:
\begin{itemize}
\item $\|P_k-Q_k\|\le\delta$;
\item for every $x\in [B'_R]\setminus\{0\}$, we have $d_\infty(Q_k(x),\Z^n)> d_0(R,\delta)$;
\item the set $Q_k\circ \widehat{Q_{k-1}} \circ \cdots\circ \widehat{Q_1}(\Gamma)$ is equidistributed modulo $\Z^n$.
\end{itemize}

We then consider the difference $u_0$ (given by Equation~\eqref{EqDefU0}). We denote by $\lfloor P \rfloor (u)$ the point of the smallest integer cube of dimension $n'\le n$ that contains $P(u)$ which has the smallest Euclidean norm (that is, the point of the support of $\varphi_{P(u)}$ with the smallest Euclidean norm). In particular, if $P(u)\notin \Z^n$, then $\|\lfloor P\rfloor (u)\|_2 < \|P(u)\|_2$ (where $\|\cdot\|_2$ is the Euclidean norm). Then, the point (ii) of Proposition \ref{ActionDiff} shows that 
\begin{align*}
\rho_{\widehat{Q_1}(\Gamma)}(\lfloor Q_1 \rfloor (u_0)) & \ge \frac{D(\Gamma)}{D(\widehat{Q_1}(\Gamma))} \varphi_{Q_1(\lfloor Q_1 \rfloor (u_0))} (u_0) \rho_\Gamma(u_0)\\
     & \ge \frac{\big(d_0(R,\delta)\big)^n}{2\times 53^n}D(\Gamma),
\end{align*}
and so on, for every $k\in\N^*$,
\[\rho_{(\widehat{Q_k}\circ\cdots\circ \widehat{Q_1})(\Gamma)}\big(\big(\lfloor Q_{k} \rfloor \circ \cdots \circ \lfloor Q_1 \rfloor \big)(u_0)\big) \ge \left(\frac{\big(d_0(R,\delta)\big)^n}{2\times 53^n}\right)^k D(\Gamma).\]

We then remark that the Euclidean norm of $\lfloor Q_1 \rfloor (u_0)$ can only take a finite number of values (it lies in $\sqrt\Z$). Then, there exists $k_0\le R^2$ such that
\[\big(\lfloor Q_{k_0} \rfloor \circ \cdots \circ \lfloor Q_1 \rfloor \big) (u_0) = 0;\]
in particular, by Equation \eqref{superEq}, we have
\[k_0\le 2n^2 + \frac{89\,888}{(D(\Gamma))^{n/2}}.\]
Then, point (ii) of Remark \ref{RemActionDiff} implies that the density of the image set satisfies
\[D\big((\widehat{Q_k}\circ\cdots\circ \widehat{Q_1})(\Gamma)\big) \le \left(1-\left(\frac{\big(d_0(R,\delta)\big)^n}{2\times 53^n}\right)^{k_0}\right) D(\Gamma).\]
We obtain the conclusions of the lemma by setting $\lambda_0 = 1-\left(\frac{(d_0(R,\delta))^n}{2\times 53^n}\right)^{k_0}$.
\end{proof}

\begin{rem}
We could be tempted to try to apply this strategy of proof to any sequence of matrices in $SL_n(\R)$. However, this does not work in the general case, because if we consider a non-Euclidean norm $N$, there could be some $x\in\R^n$ such that for every vertex $y$ of the cube containing $x$, we have $N(y)\ge N(x)$.
\end{rem}
\bigskip

We now set out a conjecture which states a generic dichotomy for the behaviour of the frequency of differences for discretizations of generic sequences of matrices, in the same vein as \cite{MR1944399} (for dimension 2) or \cite{ArturSylvain} (in the general case).

\begin{conj}\label{ConjSecond}
Let $(A_k)_{k\ge 1}$ be a generic sequence of matrices of $SL_n(\R)$. Then, if we note $\Gamma_k = (\widehat A_k\circ\cdots\circ\widehat A_1)(\Z^n)$, the standard deviation of $\rho_{\Gamma_k}$ is either very small, either as big as possible. More precisely, the standard deviation of
\[\frac{1}{\sqrt{D(\Gamma_k) \big(1-D(\Gamma_k)\big)}}\rho_{\Gamma_k}\]
either tends to 0 as $k$ goes to infinity (which corresponds to the case of zero Lyapunov exponent), or tends to 1 as $k$ goes to infinity (which corresponds to the hyperbolic case).
\end{conj}

%
The idea underlying the conjecture is that for for every $\varep>0$, for every generic sequence $(A_k)_{k\ge 1}$ of matrices and every $k$ big enough, we have two cases.
\begin{itemize}
\item Either
\[D\left\{ v\in\Z^n\ \Big\vert\ \left|\frac{1}{D(\Gamma_k)}\rho_{\Gamma_k}(v) - 1\right| \ge \varep \right\}<\varep,\]
which means that $\rho_{\Gamma_k}(v)$ is very close to its mean on a set of density almost 1. This corresponds to the zero Lyapunov exponent case, and the idea is that diffusion process wins on the hyperbolicity of the sequence of matrices.
\item Or
\[D\left\{ v\in\Z^n\ \Big\vert\ \rho_{\Gamma_k}(v) \in [\varep,1-\varep] \right\}<\varep,\]
which means that $\rho_{\Gamma_k}(v)$ is very close to 0 or 1 on a set of density almost 1. This corresponds to the hyperbolic case and the idea is that the hyperbolicity of the sequence of matrices wins on the diffusion process.
\end{itemize}

\section[Applications of model sets]{Applications of the notion of model set}

This section is devoted to the application of the notion of model set to the computation of the rate of injectivity of a finite number of matrices.

We take advantage of the rational independence between the matrices of a generic sequence to generalize to arbitrary times the geometric formulas of Section~\ref{ptgeom}. The end of the present section is devoted to the proof of the main theorem of this chapter (Theorem~\ref{ConjPrincip}).
\bigskip

Let us summarize the different notations we will use throughout this section. We will denote by $0^k$\index{$0^k$} the origin of the space $\R^k$, and $W^k = ]-1/2,1/2]^{nk}$ (unless otherwise stated). In this section, we will denote $D_c(E)$\index{$D_c$} the density of a ``continuous'' set $E\subset \R^n$, defined as (when the limit exists)
\[D_c(E) = \lim_{R\to+\infty} \frac{\Leb(B_R\cap E)}{\Leb(B_R)},\]
while for a discrete set $E\subset \R^n$, the notation $D_d(E)$\index{$D_d$} will indicate the discrete density of $E$, defined as (when the limit exists)
\[D_d(E) = \lim_{R\to+\infty} \frac{\card(B_R\cap E)}{\card(B_R\cap \Z^n)},\]

We will consider $(A_k)_{k\ge 1}$ a sequence of matrices of $SL_n(\R)$, and denote
\[\Gamma_k = (\widehat{A_k}\circ\dots\circ\widehat{A_1}) (\Z^n).\]
Also, $\Lambda_k$ will be the lattice $M_{A_1,\cdots,A_k} \Z^{n(k+1)}$, with
\begin{equation*}
M_{A_1,\cdots,A_k} = \left(\begin{array}{ccccc}
A_1 & -\Id &        &        & \\
    & A_2  & -\Id   &        & \\
    &      & \ddots & \ddots & \\
    &      &        & A_k    & -\Id\\
    &      &        &        & \Id
\end{array}\right)\in M_{n(k+1)}(\R),
\end{equation*}
and $\widetilde \Lambda_k$ will be the lattice $\widetilde M_{A_1,\cdots,A_k} \Z^{nk}$, with
\begin{equation}\label{DefMatTilde}
\widetilde M_{A_1,\cdots,A_k} = \left(\begin{array}{ccccc}
A_1 & -\Id &        &         & \\
    & A_2  & -\Id   &         & \\
    &      & \ddots & \ddots  & \\
    &      &        & A_{k-1} & -\Id\\
    &     &        &          & A_k
\end{array}\right)\in M_{nk}(\R).
\end{equation}
Finally, we will denote
\[\overline\tau^k(A_1,\cdots,A_k) = D_c\left( W^{k+1} + \Lambda_k \right).\]

These quantities will be related during this section.

\subsection[A geometric viewpoint to compute the rate of injectivity]{A geometric viewpoint to compute the rate of injectivity in arbitrary times}

We recall the (trivial) equality stated in Proposition~\ref{ImgModel}: let $A_1,\cdots,A_k\in GL_n(\R)$, then\index{$\Gamma_k$}
\begin{align}
\Gamma_k & = (\widehat{A_k}\circ\dots\circ\widehat{A_1}) (\Z^n)\nonumber\\
         & = \big\{p_2(\lambda_k)\mid \lambda_k\in \Lambda_k,\, p_1(\lambda)\in W^k\big\}\nonumber\\
				 & = p_2\Big(\Lambda_k \cap \big(p_1^{-1}(W^k)\big)\Big),\label{CalcGamma}
\end{align}
with\index{$M_{A_1,\cdots,A_k}$}
\begin{equation}\label{DefMat}
M_{A_1,\cdots,A_k} = \left(\begin{array}{ccccc}
A_1 & -\Id &        &        & \\
    & A_2  & -\Id   &        & \\
    &      & \ddots & \ddots & \\
    &      &        & A_k    & -\Id\\
    &      &        &        & \Id
\end{array}\right)\in M_{n(k+1)}(\R),
\end{equation}
$\Lambda_k = M_{A_1,\cdots,A_k} \Z^{n(k+1)}$,\index{$\Lambda_k$} $p_1$ the projection on the $nk$ first coordinates, $p_2$ the projection on the $n$ last coordinates and $W^k = ]-1/2,1/2]^{nk}$.

Here, we suppose that the set $p_1(\Lambda_k)$ is dense (thus, equidistributed) in the image set $\im p_1$ (note that this condition is generic among the sequences of invertible linear maps). In particular, the set $\{p_2(\gamma)\mid \gamma\in\Lambda_k\}$ is equidistributed in the window $W^k$.

The following property makes the link between the density of $\Gamma_k$ -- that is, the rate of injectivity of $A_1,\cdots,A_k$ -- and the density of the union of unit cubes centred on the points of the lattice $\Lambda_k$.

\begin{prop}\label{CalculTauxModel}
For a generic sequence of matrices $(A_k)_k$ of $SL_2(\R)$, we have
\[D_d(\Gamma_k) = D_c\left( W^k + \widetilde\Lambda_k \right).\]
Equivalently, we have
\[D_d(\Gamma_k) = D_c\left( W^{k+1} + \Lambda_k \right).\]
\end{prop}

Of course, this proposition generalizes Proposition~\ref{FormTau1} to an arbitrary number of matrices $A_1,\cdots,A_k$. Remark that this kind of ideas was already present in \cite[Theorem 5.3.1]{MR1440853}.

\begin{rem}
The density on the left of the equality is the density of a discrete set (that is, with respect to counting measure), whereas the density on the right of the equality is that of a continuous set (that is, with respect to Lebesgue measure). The two notions coincide when we consider discrete sets as sums of Dirac masses.
\end{rem}

This proposition leads to the definition of the mean rate of injectivity in time $k$.

\begin{definition}\label{DefTauBarrr}
Let $(A_k)_{k\ge 1}$ be a sequence of matrices of $SL_n(\R)$. The \emph{mean rate of injectivity in time $k$} of this sequence is defined by\index{$\overline\tau^k$}
\[\overline\tau^k(A_1,\cdots,A_k) = D_c\left( W^{k+1} + \Lambda_k \right).\]
\end{definition}

\begin{rem}\label{conttaukk}
As in Section~\ref{ptgeom} for the rate of injectivity in time 1, Proposition~\ref{CalculTauxModel} asserts that for a generic sequence of matrices, the rate of injectivity $\tau^k$ in time $k$ coincides with the mean rate of injectivity $\overline\tau^k$, which is continuous and piecewise polynomial of degree $\le nk$ in the coefficients of the matrix.
\end{rem}

\begin{rem}
The formula of Proposition~\ref{CalculTauxModel} could be used to compute numerically the mean rate of injectivity in time $k$ of a sequence of matrices: it is much faster to compute the area of a finite number of intersections of cubes (in fact, a small number) than to compute the cardinalities of the images of a big set $[-R,R]^n \times \Z^n$.
\end{rem}

\begin{proof}[Proof of Proposition \ref{CalculTauxModel}]
We want to determine the density of $\Gamma_k$. By Equation~\eqref{CalcGamma}, we have
\[x\in\Gamma_k \iff x\in\Z^n\ \text{and}\ \exists \lambda\in \Lambda_k : x=p_2(\lambda),\, p_1(\lambda) \in W^k.\]
But if $x=p_2(\lambda)$, then we can write $\lambda=(\widetilde\lambda,0^n) + (0^{(k-1)n},-x,x)$ with $\widetilde\lambda\in \widetilde\Lambda_k$. Thus,
\begin{align*}
x\in\Gamma_k & \iff x\in\Z^n\ \text{and}\ \exists \widetilde\lambda\in \widetilde\Lambda_k : (0^{(k-1)n},-x)-\widetilde\lambda	\in W^k\\
             & \iff x\in\Z^n\ \text{and}\ (0^{(k-1)n},x)\in \bigcup_{\widetilde\lambda\in\widetilde\Lambda_k} \widetilde\lambda - W^k.
\end{align*}
Thus, $x\in\Gamma_k$ if and only if the projection of $(0^{(k-1)n},x)$ on $\R^{nk}/\widetilde\Lambda_k$ belongs to $\bigcup_{\widetilde\lambda\in\widetilde\Lambda_k} \widetilde\lambda - W^k$. Then, the proposition follows directly from the fact that the points of the form $(0^{(k-1)n},x)$, with $x\in\Z^n$, are equidistributed in $\R^{nk} / \widetilde \Lambda_k$.
\bigskip

To prove this equidistribution, we compute the inverse matrix of $\widetilde M_{A_1,\cdots,A_k}$:
\[{\widetilde M_{A_1,\cdots,A_k}}^{-1} = \begin{pmatrix}
A_1^{-1} & A_1^{-1}A_2^{-1} &  A_1^{-1}A_2^{-1}A_3^{-1} & \cdots & A_1^{-1}\cdots A_k^{-1}\\
    & A_2^{-1}  & A_2^{-1}A_3^{-1} & \cdots & A_2^{-1}\cdots A_k^{-1}\\
    &      & \ddots & \ddots  & \vdots\\
    &      &        & A_{k-1}^{-1} & A_{k-1}^{-1}A_k^{-1}\\
    &     &        &          & A_k^{-1}
\end{pmatrix}.\]
Thus, the set of points of the form $(0^{(k-1)n},x)$ in $\R^{nk} / \widetilde \Lambda_k$ corresponds to the image of the action
\[\Z^n \ni x \longmapsto
\begin{pmatrix}
A_1^{-1}\cdots A_k^{-1}\\
A_2^{-1}\cdots A_k^{-1}\\
\vdots\\
A_{k-1}^{-1}A_k^{-1}\\
A_k^{-1}
\end{pmatrix}x\]
of $\Z^n$ on the canonical torus $\R^{nk}/\Z^{nk}$. But this action is trivially ergodic (even in restriction to the first coordinate) when the sequence of matrices is generic.
\end{proof}

As in the case of the time 1 (Proposition~\ref{ActionDiffGeom}), there is also a geometric method to compute the frequency of a difference in $\Gamma_k$.

\begin{figure}
\begin{minipage}[t]{.47\linewidth}
\begin{center}
\begin{tikzpicture}[scale=1.4]
\draw (-1,0) -- (3,0);
\draw[fill=white,opacity=.7] (-.5,-1.2) -- (.5,-.6) -- (.5,1.2) -- (-.5,.6) -- cycle;
\draw[->] (0,-.08) to (2,-.08);
\draw[fill=white,opacity=.7] (1.5,-1.2) -- (2.5,-.6) -- (2.5,1.2) -- (1.5,.6) -- cycle;
\draw[color=blue!40!black] (2-.167,-.4) -- (2.167,-.2) -- (2.167,.4) -- (2-.167,.2) -- cycle;
\draw[color=blue!40!black] (2.167,-.2) node[below]{$\scriptstyle W^k$};
\draw (0,0) -- (1,0);
\draw (2,0) -- (3,0);
\foreach\i in {-1,...,1}{
\foreach\j in {-1,...,1}{
\draw[color=green!40!black] (.3*\i-.05*\j,.3*\i+.6*\j) node{$\scriptscriptstyle\bullet$};
\draw[color=green!40!black] (.3*\i-.05*\j+2+.07,.3*\i+.6*\j+.05) node{$\scriptscriptstyle\bullet$};
}}
\draw[color=green!40!black] (0,1.4) node{$\scriptstyle p_2^{-1}(0)\cap \Lambda_k$};
\draw[color=green!40!black] (2,1.4) node{$\scriptstyle p_2^{-1}(x)\cap \Lambda_k$};
\draw (0,0) node{$\bullet$};
\draw (2,0) node{$\bullet$};
\draw (1,-.3) node {$x$};
\end{tikzpicture}
\caption{Construction of Proposition~\ref{CalculTauxModel}.}\label{FigCalculTauxModel}
\end{center}
\end{minipage}\hfill
\begin{minipage}[t]{.47\linewidth}
\begin{center}
\begin{tikzpicture}[scale=1.4]
\draw (-1,0) -- (3,0);
\draw[fill=white,opacity=.7] (-.5,-1.2) -- (.5,-.6) -- (.5,1.2) -- (-.5,.6) -- cycle;
\draw[->] (-.09,-.08) to (2-.09,-.08);
\draw[fill=white,opacity=.7] (1.5,-1.2) -- (2.5,-.6) -- (2.5,1.2) -- (1.5,.6) -- cycle;
\draw[color=blue!40!black] (-.167,-.4) -- (.167,-.2) -- (.167,.4) -- (-.167,.2) -- cycle;
\draw[color=blue!40!black] (2-.167,-.4) -- (2.167,-.2) -- (2.167,.4) -- (2-.167,.2) -- cycle;
\draw (0,0) -- (1,0);
\draw (2,0) -- (3,0);
\foreach\i in {-1,...,1}{
\foreach\j in {-1,...,1}{
\draw[color=green!40!black] (.3*\i-.05*\j-.09,.3*\i+.6*\j-.08) node{$\scriptscriptstyle\bullet$};
\draw[color=green!40!black] (.3*\i-.05*\j+2+.07,.3*\i+.6*\j+.05) node{$\scriptscriptstyle\bullet$};
}}
\draw[color=green!40!black] (0,1.4) node{$\scriptstyle p_2^{-1}(x)\cap \Lambda_k$};
\draw[color=green!40!black] (2,1.4) node{$\scriptstyle p_2^{-1}(x+v)\cap \Lambda_k$};
\draw (0,0) node{$\bullet$};
\draw (2,0) node{$\bullet$};
\draw (1,-.3) node {$v$};
\end{tikzpicture}
\caption{Construction of Proposition~\ref{CalculDiffModel}}\label{FigCalculDiffModel}
\end{center}
\end{minipage}
\end{figure}

\begin{prop}\label{CalculDiffModel}
For a generic sequence of matrices $(A_k)_k$, we have
\[\rho_{\Gamma_k}(v) = \Leb\Big( \big(W^{k+1} + \Lambda_k\big)\cap \big(W^{k+1} + (0^{kn},v)\big)\Big).\]
\end{prop}

\begin{proof}[Proof of Proposition~\ref{CalculDiffModel}]
To compute $\rho_{\Gamma_k}(v)$, we want to know for which proportion of $x\in\Z^n$ such that $x\in\Gamma_k$, we have $x+v\in\Gamma_k$, that is, knowing that $p_1\big(p_2^{-1}(x)\cap\Lambda_k\big)\cap W^k\neq \emptyset$, we have $p_1\big(p_2^{-1}(x+v)\cap\Lambda_k\big)\cap W^k\neq \emptyset$. But
\begin{align*}
p_1\big(p_2^{-1}(x+v)\cap\Lambda_k\big)\cap W^k\neq \emptyset \quad & \iff \quad \big(\widetilde \Lambda_k + (0^{(k-1)n},-x-v)\big) \cap W^k \neq\emptyset\\
     & \iff \quad (0^{(k-1)n},x+v)\in \widetilde \Lambda_k + W^k.\\
\end{align*}
As the sets $p_1(p_2^{-1}(x)\cap\Lambda_k)$ are equidistributed in $W^k$ (see the proof of Proposition~\ref{CalculTauxModel}), the proportion of such points $x$ is equal to the area of the intersection $(W^{k+1} + \Lambda_k)\cap (W^{k+1} + (0^{kn},v))$.
\end{proof}

As in the case of time 1 (see Section~\ref{ptgeom}), an argument of equidistribution combined with Proposition~\ref{IntRho} allows to deduce Proposition~\ref{CalculTauxModel} from Proposition~\ref{CalculDiffModel}.

There is also a dual method to compute the rate of injectivity of a sequence of matrices: we define $\psi : \R^{nk}\to \R$ by
\[\psi = \sum_{\lambda\in\Lambda_k} \1_{W^k + \lambda},\]
and obtain the following formula (which is a generalization of Proposition~\ref{FormTau2}).

\begin{prop}\label{CalculTauxModel2}
For a generic sequence of matrices $(A_k)_k$, we have
\[D_d(\Gamma_k) = \int_{B_{1/2}} \frac{1}{\psi(\lambda)}\, \ud \Leb(\lambda).\]
\end{prop}

The proof of this proposition is similar to that of Proposition~\ref{FormTau2} page~\pageref{FormTau2} (see in particular Lemma~\ref{DoubleComptage}).
\bigskip

Recall the problem raised by Theorem~\ref{ConjPrincip}: we want to make $\tau^k$ as small as possible. By an argument of equidistribution, generically, it is equivalent to make the mean rate of injectivity $\overline\tau^k$ as small as possible when $k$ goes to infinity, by perturbing every matrix in $SL_n(\R)$ of at most $\delta>0$ (fixed once for all). In the framework of model sets, Theorem~\ref{ConjPrincip} is motivated by the phenomenon of concentration of the measure on a neighbourhood of the boundary of the cubes in high dimension.

\begin{prop}\label{concentration}
Let $W^k$ be the infinite ball of radius $1/2$ in $\R^k$ and $v^k$ the vector $(1,\cdots,1)\in \R^k$. Then, for every $\varep,\delta>0$, there exists $k_0\in\N^*$ such that for every $k\ge k_0$, we have $\Leb\big(W^k \cap (W^k + \delta v^k)\big) < \varep$.
\end{prop}


Applying Haj\'os theorem (Theorem~\ref{hajos}), it is easy to see when the density of $\Gamma_k$ is equal to 1: combining this theorem with Proposition~\ref{CalculTauxModel}, we obtain that this occurs if and only if the lattice given by the matrix $M_{A_1,\cdots,A_k}$ satisfies the conclusions of Haj\'os theorem\footnote{Of course, this property can be obtained directly by saying that the density is equal to 1 if and only if the rate of injectivity of every matrix of the sequence is equal to 1.}. The heuristic suggested by the phenomenon of concentration of the measure is that if we perturb ``randomly'' any sequence of matrices, we will go ``far away'' from the lattices satisfying Haj\'os theorem and then the rate of injectivity will be close to 0.

\begin{rem}
The kind of questions addressed by Haj\'os theorem are in general quite delicate. For example, we can wonder what happens if we do not suppose that the centres of the cubes form a lattice of $\R^n$. O. H. Keller conjectured in \cite{zbMATH02567416} that the conclusion of Haj\'os theorem is still true under this weaker hypothesis. This conjecture was proven to be true for $n\le 6$ by O. Perron in \cite{MR0003041,MR0002185}, but remained open in higher dimension until 1992, when J. C. Lagarias and P. W. Shor proved in \cite{MR1155280} that Keller's conjecture is false for $n\ge 10$ (this result was later improved by \cite{MR1920144} which shows that it is false as soon as $n\ge 8$; the case $n=7$ is still open).
\end{rem}

\subsection[New proof that the rate is smaller than $1/2$]{New proof that the asymptotic rate of injectivity is generically smaller than $1/2$}

As a first application of the concept of model set, we give a new proof of the fact that rate of injectivity of a generic sequence of $SL_n(\R)$ is smaller than $1/2$. To begin with, we give a lemma estimating the sizes of intersections of cubes when the rate is bigger than $1/2$ (which is the geometric counterpart of Lemma~\ref{majoration}).

\begin{lemme}\label{EstimTaux}
Let $W^k = ]-1/2,1/2]^k$ and $\Lambda\subset \R^k$ be a lattice with covolume 1 such that $D_c(W^k + \Lambda)\ge 1/2$. Then, for every $v\in\R^k$, we have
\[D_c\big((W^k + \Lambda + v)\cap (W^k + \Lambda)\big) \ge 2D_c(W^k + \Lambda)-1.\]
\end{lemme}

\begin{proof}[Proof of Lemma~\ref{EstimTaux}]
We first remark that $D_c(\Lambda+W^k)$ is equal to the volume of the projection of $W^k$ on the quotient space $\R^k/\Lambda$. For every $v\in\R^k$, the projection of $W^k+v$ on $\R^k/\Lambda$ has the same volume; as this volume is greater than $1/2$, and as the covolume of $\Lambda$ is 1, the projections of $W^k$ and $W^k+v$ overlap, and the volume of the intersection is bigger than $2D_c(W^k + \Lambda)-1$. Returning to the whole space $\R^k$, we get the conclusion of the lemma.
\end{proof}

\begin{rem}	
In the case where $\Lambda$ is the lattice spanned by $M_{A_1,\cdots,A_k}$, where $A_1,\cdots,A_k$ is a generic family, Lemma~\ref{EstimTaux} can be deduced directly from Lemma~\ref{majoration} (more precisely, improves its conclusion of a factor 2).
\end{rem}

This lemma allows us to give another proof of Theorem~\ref{PerLin1} (page~\pageref{PerLin1}), which implies in particular that the asymptotic rate of injectivity of a generic sequence of matrices of $SL_n(\R)$ is smaller than $1/2$ (and even is smaller than a sequence converging exponentially fast to $1/2$).

\begin{proof}[Alternative proof of Theorem~\ref{PerLin1}]\label{Proof2PerLin1}
Let $(A_k)_{k\ge 1}$ be a bounded sequence of matrices of $SL_n(\R)$ and $\delta>0$. As in the previous proof of Theorem~\ref{PerLin1}, we proceed by induction on $k$ and suppose that the theorem is proved for a rank $k\in\N^*$. Let $\widetilde \Lambda_k$ be the lattice spanned by the matrix $\widetilde M_{B_1,\cdots,B_k}$ (defined by Equation~\eqref{DefMatTilde} page~\pageref{DefMatTilde}) and $W^k = ]-1/2, 1/2]^{nk}$ be the window corresponding to the model set $\Gamma_k$ modelled on $\Lambda_k$ (the lattice spanned by the matrix $M_{B_1,\cdots,B_k}$, see Equation~\eqref{DefMat}). By Proposition~\ref{CalculTauxModel}, we have
\[\tau^k(B_1,\cdots,B_k) = D_c\left(W^k+ \widetilde \Lambda_k \right).\]

We now choose a matrix $B_{k+1}$ satisfying $\|A_{k+1} - B_{k+1}\|\le\delta$, such that there exists $x_1\in \Z^n\setminus \{0\}$ such that $\| B_{k+1}x_1\|_\infty \le 1-\varep$, with $\varep>0$ depending only on $\delta$ and $\|(A_k)_k\|$ (and $n$): indeed, for every matrix $B\in SL_n(\R)$, Minkowski theorem implies that there exists $x_1\in \Z^n\setminus \{0\}$ such that $\| Bx_1\|_\infty \le 1$; it then suffices to modify slightly $B$ to decrease $\| Bx_1\|_\infty$. We can also suppose that the sequence of matrices $B_1,\cdots,B_{k+1}$ is generic. Again, Proposition~\ref{CalculTauxModel} reduces the calculation of the rate on injectivity $\tau^{k+1}(B_1,\cdots,B_{k+1}) $ to that of the density of $W^{k+1}+ \widetilde \Lambda_{k+1}$. By the form of the matrix $M_{B_1,\cdots,B_k}$, this set can be decomposed into
\[W^{k+1}+ \widetilde \Lambda_{k+1} = W^{k+1} + \begin{pmatrix} \widetilde \Lambda_k\\0^n \end{pmatrix} +
\begin{pmatrix} 0^{n(k-1)}\\ -\Id\\ B_{k+1} \end{pmatrix} \Z^n.\]
In particular, as $|\det(B_{k+1})|=1$, this easily implies that $D_c\left(W^{k+1}+ \widetilde \Lambda_{k+1}\right)\le D_c\left(W^k+ \widetilde \Lambda_{k} \right)$.

\begin{figure}[t]
\begin{minipage}[t]{\linewidth}
\emph{How to read these figures :} The top of the figure represents the set $W^k + \widetilde\Lambda_k$ by the 1-dimensional set $[-1/2,1/2] + \nu\Z$ (in dark blue), for a number $\nu>1$. The bottom of the figure represents the set $W^{k+1} + \widetilde\Lambda_{k+1}$ by the set $[-1/2,1/2]^2 + \Lambda$, where $\Lambda$ is the lattice spanned by the vectors $(0,\nu)$ and $(1,1-\varep)$ for a parameter $\varep>0$ close to 0. The dark blue cubes represent the ``old'' cubes, that is, the thickening in dimension 2 of the set $W^k + \widetilde\Lambda_k$, and the light blue cubes represent the ``added'' cubes, that is, the rest of the set $W^{k+1} + \widetilde\Lambda_{k+1}$.
\end{minipage}\vspace{25pt}

\begin{minipage}[t]{.48\linewidth}
\begin{center}
\begin{tikzpicture}[scale=.95]
\draw[color=gray] (-.8,3) -- (5,3);
\draw[very thick,blue, |-|] (-.5,3) -- (.5,3);
\draw[very thick,blue, |-|] (-.5+1.8,3) -- (.5+1.8,3);
\draw[very thick,blue, |-|] (-.5+3.6,3) -- (.5+3.6,3);

\clip (-.6,-1.3) rectangle (4.8,1.5);
\foreach\j in {0,...,2}{
\draw[fill=blue,opacity=.25] (-.5+1.8*\j,-.5) rectangle (.5+1.8*\j,.5);
}
\foreach\i in {-2,...,2}{
\foreach\j in {-1,...,3}{
\draw[fill=blue,opacity=.20] (-.5+\i+1.8*\j,-.5+.8*\i) rectangle (.5+\i+1.8*\j,.5+.8*\i);
\draw (-.5+\i+1.8*\j,-.5+.8*\i) rectangle (.5+\i+1.8*\j,.5+.8*\i);
}}
\end{tikzpicture}
\caption[Intersection of cubes, rate bigger than $1/2$]{In the case where the rate is bigger than $1/2$, some intersections of cubes appear automatically between times $k$ and $k+1$.}\label{FigTauDemi}
\end{center}
\end{minipage}\hfill
\begin{minipage}[t]{.48\linewidth}
\begin{center}
\begin{tikzpicture}[scale=.95]
\draw[color=gray] (-.8,3) -- (5,3);
\draw[very thick,blue, |-|] (-.5,3) -- (.5,3);
\draw[very thick,blue, |-|] (-.5+2.8,3) -- (.5+2.8,3);

\clip (-.6,-1.3) rectangle (4.8,1.5);
\foreach\j in {0,...,1}{
\draw[fill=blue,opacity=.25] (-.5+2.8*\j,-.5) rectangle (.5+2.8*\j,.5);
}
\foreach\i in {-3,...,2}{
\foreach\j in {-1,...,2}{
\draw[fill=blue,opacity=.20] (-.5+\i+2.8*\j,-.5+.8*\i) rectangle (.5+\i+2.8*\j,.5+.8*\i);
\draw (-.5+\i+2.8*\j,-.5+.8*\i) rectangle (.5+\i+2.8*\j,.5+.8*\i);
}}
\end{tikzpicture}
\caption[Intersection of cubes, rate smaller than $1/2$]{In the case where the rate is smaller than $1/2$, there is not necessarily new intersections between times $k$ and $k+1$.}\label{FigTauTiers}
\end{center}
\end{minipage}
\end{figure}

What we need is a more precise bound. We apply Corollary~\ref{CoroSansNom} to
\[\Lambda_1 = \big(\widetilde\Lambda_k,0^n\big), \qquad \Lambda_2 = \begin{pmatrix} 0^{n(k-1)}\\ -\Id\\ B_{k+1} \end{pmatrix} \Z^n \qquad \text{and} \qquad B = W^{k+1}.\]
Then, the decreasing of the rate of injectivity between times $k$ and $k+1$ is bigger than the $D_1$ defined in Corollary~\ref{CoroSansNom}: using Lemma~\ref{EstimTaux}, we have
\[D_c\left(\Big(W^k + \widetilde\Lambda_k\Big) \cap \Big(W^k + \widetilde\Lambda_k + \big(0^{n(k-1)},-x_1\big)\Big)\right) \ge 2D_d(\Gamma_k)-1;\]
thus, as $\|x_1\|_\infty <1-\varep$
\[D_1 \ge \varep^n \big(2D_d(\Gamma_k)-1\big).\]
From Corollary~\ref{CoroSansNom} we deduce that
\begin{align*}
D(\Gamma_{k+1}) & = D_c\left(W^{k+1}+ \widetilde \Lambda_{k+1}\right)\\
        & \le D_c\left(W^k + \widetilde\Lambda_k\right) - \frac12 D_1\\
				& \le D_d(\Gamma_k) - \frac12 \varep^n \big(2D_d(\Gamma_k)-1\big).
\end{align*}
This proves the theorem for the rank $k+1$.
\end{proof}

\subsection[Generically, the asymptotic rate is zero]{Proof of Theorem~\ref{ConjPrincip}: generically, the asymptotic rate is zero}

We now come to the proof of Theorem~\ref{ConjPrincip}. The strategy of proof is identical to that we used in the previous section to state that generically, the asymptotic rate is smaller than $1/2$: we will use an induction to decrease the rate step by step. Recall that $\overline\tau^k(A_1,\cdots,A_k)$ indicates the density of the set $W^{k+1} + \Lambda_k$ (see Definition~\ref{DefTauBarrr}).

Unfortunately, if the density of $W^k + \widetilde\Lambda_k$ -- which is generically equal to the density of the $k$-th image $\big(\widehat A_k\circ\cdots\circ \widehat A_1\big)(\Z^n)$ -- is smaller than $1/2$, then we can not apply exactly the strategy of proof of the previous section (see Figure~\ref{FigTauTiers}). For example, if we take
\[A_1 = \begin{pmatrix} 4 & \\ & 1/4 \end{pmatrix}\quad \text{and}\quad A_2 = \begin{pmatrix} 1/2 & \\ & 2 \end{pmatrix},\]
then $\big(\widehat{A_2}\circ \widehat{A_1}\big)(\Z^2) = (2\Z)^2$, and for every $B_3$ close to the identity, we have $\tau^3(A_1,A_2,B_3) = \tau^2(A_1,A_2) = 1/4$.

\begin{figure}[t]
\begin{center}
\includegraphics[width=\linewidth]{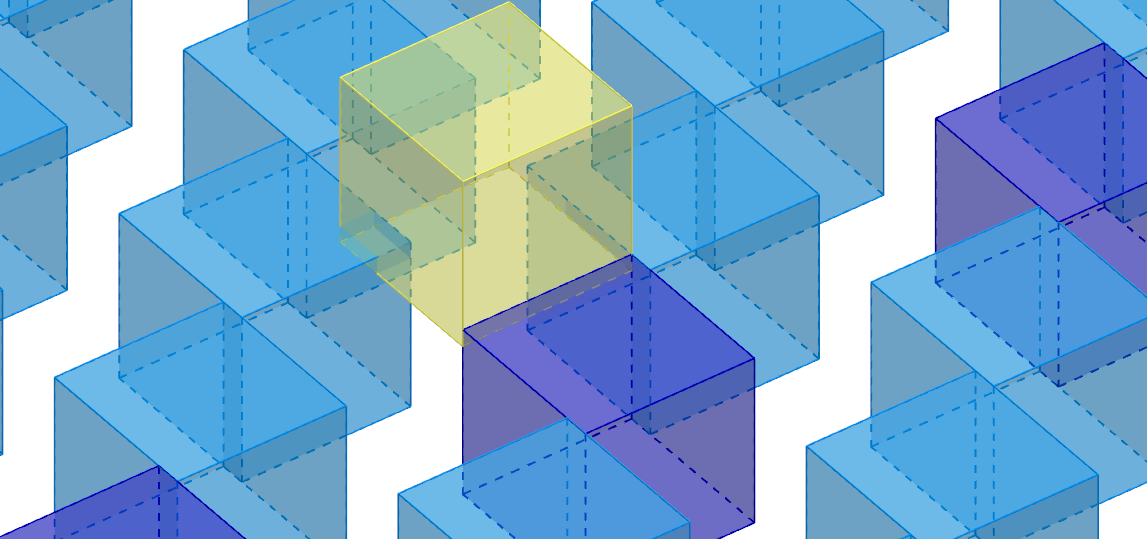}
\caption[Intersection of cubes, rate bigger than $1/3$]{Intersection of cubes in the case where the rate is bigger than $1/3$. The thickening of the cubes of $W^k + \widetilde\Lambda_k$ is represented in dark blue and the thickening of the rest of the cubes of $W^{k+1} + \widetilde\Lambda_{k+1}$ is represented in light blue; we have also represented another cube of $W^{k+2} + \widetilde\Lambda_{k+2}$ in yellow. We see that if the projection on the $z$-axis of the centre of the yellow cube is smaller than 1, then there is automatically an intersection between this cube and one of the blue cubes.}\label{FigInterCubes3DSupDemi}
\end{center}
\end{figure}

To overcome this difficulty, we prove that for a generic sequence $(A_k)_{k\ge 1}$, if $\overline\tau^k(A_1,\cdots,A_k)> 1/\ell$, then $\overline\tau^{k+\ell-1}(A_1,\cdots,A_{k+\ell-1})$ is strictly smaller than $\overline\tau^k(A_1,\cdots,A_k)$. An argument of equirepartition (in fact, Propostion~\ref{CalculTauxModel}) allows to see this problem in terms of area of intersections of cubes. More precisely, we consider the maximal number of disjoint translates of $W^k + \widetilde\Lambda_k$ in $\R^{nk}$: we easily see that if the density of $W^k + \widetilde\Lambda_k$ is bigger than $1/\ell$, then there can not be more than $\ell$ disjoint translates of $W^k + \widetilde\Lambda_k$ in $\R^{nk}$(Lemma~\ref{EstimTauxN}). Then, Lemma~\ref{LemDeFin} states that if the sequence of matrices is generic, either the density of $W^{k+1} + \widetilde\Lambda_{k+1}$ is smaller than that of $W^k + \widetilde\Lambda_k$ (see Figure~\ref{FigLemDeFin}), or there can not be more than $\ell-1$ disjoint translates of $W^{k+1} + \widetilde\Lambda_{k+1}$ in $\R^{n(k+1)}$(see Figure~\ref{FigLemDeFin2}). Applying this reasoning (at most) $\ell-1$ times, we obtain that the density of $W^{k+\ell-1} + \widetilde\Lambda_{k+\ell-1}$ is smaller than that of $W^k + \widetilde\Lambda_k$. For example if $D_c\big(W^k + \widetilde\Lambda_k\big) > 1/3$, then $D_c\big(W^{k+2} + \widetilde\Lambda_{k+2}\big) < D\big(W^k + \widetilde\Lambda_k\big)$ (see Figure~\ref{FigInterCubes3DSupDemi}). To apply this strategy in practice, we have to obtain quantitative estimates about the loss of density we get between times $k$ and $k+\ell-1$.

Remark that with this strategy, we do not need to make ``clever'' perturbations of the matrices: provided that the coefficients of the matrices are rationally independent, the perturbation of each matrix is made independently from that of the others. However, this reasoning does not tell when exactly the rate of injectivity decreases (likely, in most of cases, the speed of decreasing of the rate of injectivity is much faster than the one obtained by this method), and does not say either where exactly the loss of injectivity occurs in the image sets.

Firstly, we give a more precise statement of Theorem~\ref{ConjPrincip}.

{\renewcommand{\thelemme}{\ref{ConjPrincip}}
\begin{theoreme}
For a generic sequence of matrices $(A_k)_{k\ge 1}$ of $SL_n(\R)$, for every $\ell\in\N$, there exists $\lambda_\ell\in]0,1[$ such that for every $k\in\N$,
\begin{equation}\label{EstimTauxExp}
\tau^{\ell k}(A_1,\cdots,A_{\ell k}) \le \lambda_\ell^k + \frac{1}{\ell}.
\end{equation}
Also, for every $\nu<1$, we have
\begin{equation}\label{EqConvLog}
\tau^k(A_1,\cdots,A_k) = o\big(\ln(k) ^{-\nu}\big).
\end{equation}
In particular, the asymptotic rate of injectivity $\tau^\infty\big( (A_k)_{k\ge 1}\big)$ is equal to zero.
\end{theoreme}
\addtocounter{lemme}{-1}}

The following lemma is a generalization of Lemma~\ref{EstimTaux}. It expresses that if the density of $W^k + \widetilde\Lambda_k$ is bigger than $1/\ell$, then there can not be more than $\ell$ disjoint translates of $W^k + \widetilde\Lambda_k$, and gives an estimation on the size of these intersections.

\begin{lemme}\label{EstimTauxN}
Let $W^k = ]-1/2,1/2]^k$ and $\Lambda\subset \R^k$ be a lattice with covolume 1 such that $D_c(W^k + \Lambda)\ge 1/\ell$. Then, for every collection $v_1,\cdots,v_\ell\in\R^k$, there exists $i\neq i'\in \llbracket 1,\ell\rrbracket$ such that
\[D_c\big((W^k + \Lambda + v_i)\cap (W^k + \Lambda + v_{i'})\big) \ge 2\frac{\ell D_c(W^k + \Lambda)-1}{\ell(\ell-1)}.\]
\end{lemme}

\begin{proof}[Proof of Lemma~\ref{EstimTauxN}]
For every $v\in\R^k$, the density $D_c(W^k + \Lambda+v)$ is equal to the volume of the projection of $W^k$ on the quotient space $\R^k/\Lambda$. As this volume is greater than $1/\ell$, and as the covolume of $\Lambda$ is 1, the projections of the sets $(W^k + v_i)_{1\le i \le \ell}$ on $\R^k/\Lambda$ overlap, and the volume of the points belonging to at least two different sets is bigger than $\ell D_c(W^k + \Lambda)-1$. As there are $\ell(\ell-1)/2$ possibilities of intersection, there exists $i\neq i'$ such that the volume of the intersection between the projections of $W^k+v_i$ and $W^k+v_{i'}$ is bigger than $2(\ell D_c(W^k + \Lambda)-1)/(\ell(\ell-1))$. Returning to the whole space $\R^k$, we get the conclusion of the lemma.
\end{proof}

Recall that we denote $\widetilde\Lambda_k$ the lattice spanned by the matrix
\[\widetilde M_{A_1,\cdots,A_k} = \left(\begin{array}{ccccc}
A_1 & -\Id &        &         & \\
    & A_2  & -\Id   &         & \\
    &      & \ddots & \ddots  & \\
    &      &        & A_{k-1} & -\Id\\
    &     &        &          & A_k
\end{array}\right)\in M_{nk}(\R),
\]
and $W^k$ the cube $]-1/2,1/2]^{nk}$.


\begin{lemme}\label{LemDeFin}
For every $\delta>0$ and every $M>0$, there exists $\varep>0$ and an open set of matrices\footnote{Independent from the matrices $A_1,\cdots,A_k$.} $\mathcal O \subset SL_n(\R)$ which is $\delta$-dense in the set of matrices of norm $\le M$, such that if $\ell\ge 2$ and $D_0>0$ are such that for every collection of vectors $v_1,\cdots,v_{\ell} \in\R^n$, there exists $j,j'\in\llbracket 1,\ell\rrbracket$ such that
\[D_c\bigg(\Big(W^k + \widetilde\Lambda_k + (0^{(k-1)n} , v_j) \Big) \cap \Big(W^k + \widetilde\Lambda_k + (0^{(k-1)n} , v_{j'}) \Big)\bigg)\ge D_0,\]
then for every $B \in \mathcal O$, if we denote by $\widetilde\Lambda_{k+1}$ the lattice spanned by the matrix $\widetilde M_{A_1,\cdots,A_k,B}$,
\begin{enumerate}[(1)]
\item either $D_c(W^{k+1} + \widetilde\Lambda_{k+1}) \le D_c(W^{k} + \widetilde\Lambda_{k}) - \varep D_0/(4\ell)$;
\item or for every collection of vectors $w_1,\cdots,w_{\ell-1} \in\R^n$, there exists $i\neq i'\in\llbracket 1,\ell-1\rrbracket$ such that
\[D_c\bigg(\Big(W^{k+1} + \widetilde\Lambda_{k+1} + (0^{kn} , w_i) \Big) \cap \Big(W^{k+1} + \widetilde\Lambda_{k+1} + (0^{kn} , w_{i'}) \Big)\bigg)\ge \varep D_0 /\ell^2.\]
\end{enumerate}
\end{lemme}

\begin{rem}
If $\ell = 2$, then we have automatically the conclusion (1) of the lemma.
\end{rem}

In a certain sense, conclusion (1) corresponds to an hyperbolic case, and conclusion (2) expresses that there is a diffusion between times $k$ and $k+1$.

\begin{proof}[Proof of Lemma \ref{LemDeFin}]
Let $\mathcal O_\varep$ be the set of the matrices $B\in SL_n(\R)$ satisfying: for any collection of vectors $w_1,\cdots,w_{\ell-1} \in\R^n$, there exists a set $U\subset\R^n/B\Z^n$ of measure $>\varep$ such that every point of $U$ belongs to at least $\ell$ different cubes of the collection $(Bv + w_i + W^1)_{v\in\Z^n,\,1\le i\le\ell-1}$. In other words\footnote{Matrices that does not possess this property are such that the union of cubes form a $k$-fold tiling. This was the subject of Furtwängler conjecture, see \cite{MR1550530}, proved false by G.~Haj\'os. R.~Robinson gave a characterization of such $k$-fold tilings in some cases, see \cite{MR526466} or \cite[p. 29]{MR1311249}.}, every $x\in\R^n$ whose projection $\overline x$ on $\R^n/B\Z^n$ belongs to $U$ satisfies
\begin{equation}\label{EqLemDeFin}
\sum_{i=1}^{\ell-1} \sum_{v\in \Z^n}  \1_{x\in Bv + w_i + W^1} \ge \ell.
\end{equation}
We easily see that the sets $\mathcal O_\varep$ are open and that the union of these sets over $\varep>0$ is dense (it contains the set of matrices $B$ whose entries are all irrational). Thus, if we are given $\delta>0$ and $M>0$, there exists $\varep>0$ such that $\mathcal O = \mathcal O_\varep$ is $\delta$-dense in the set of matrices of $SL_n(\R)$ whose norm is smaller than $M$.

We then choose $B\in \mathcal O$ and a collection of vectors $w_1,\cdots,w_{\ell-1} \in\R^n$. Let $x\in\R^n$ be such that $\overline x\in U$. By hypothesis on the matrix $B$, $x$ satisfies Equation~\eqref{EqLemDeFin}, so there exists $\ell+1$ integer vectors $v_1,\cdots,v_{\ell}$ and $\ell$ indices $i_1,\cdots,i_{\ell}$ such that the couples $(v_j,i_j)$ are pairwise distinct and that
\begin{equation}\label{EqIxIn}
\forall j\in\llbracket 1,\ell+1\rrbracket,\quad x\in Bv_j + w_{i_j} + W^1.
\end{equation}

\begin{figure}
\begin{minipage}[t]{.48\linewidth}
\begin{center}
\begin{tikzpicture}[scale=.95]
\draw[color=gray] (-.8,3) -- (5,3);
\draw[very thick,blue, |-|] (-.5,3) -- (.5,3);
\draw[very thick,blue, |-|] (-.5+2.8,3) -- (.5+2.8,3);

\draw[dashed] (-.8,.4) -- (5,.4);
\draw(-.8,.4) node[left]{$x$};
\draw[dashed] (2.4,-1.5) -- (2.4,1.7);
\draw (2.4,-1.5) node[below]{$y$};

\clip (-.6,-1.3) rectangle (4.8,1.5);
\foreach\j in {0,...,1}{
\draw[fill=blue,opacity=.25] (-.5+2.8*\j,-.5) rectangle (.5+2.8*\j,.5);
}
\foreach\i in {-3,...,5}{
\foreach\j in {-1,...,3}{
\draw[fill=blue,opacity=.20] (-.5+\i-1+2.8*\j,-.5+.4*\i-.4) rectangle (.5+\i-1+2.8*\j,.5+.4*\i-.4);
\draw (-.5+\i-1+2.8*\j,-.5+.4*\i-.4) rectangle (.5+\i-1+2.8*\j,.5+.4*\i-.4);
}}
\end{tikzpicture}
\caption[First case of Lemma \ref{LemDeFin}]{First case of Lemma \ref{LemDeFin}, in the case $\ell=3$: the set $W^{k+1} + \widetilde\Lambda_{k+1}$ auto-intersects.}\label{FigLemDeFin}
\end{center}
\end{minipage}\hfill
\begin{minipage}[t]{.48\linewidth}
\begin{center}
\begin{tikzpicture}[scale=.95]
\draw[color=gray] (-.8,3) -- (5,3);
\draw[very thick,blue, |-|] (-.5,3) -- (.5,3);
\draw[very thick,blue, |-|] (-.5+2.8,3) -- (.5+2.8,3);

\draw[dashed] (-.8,.4) -- (5,.4);
\draw(-.8,.4) node[left]{$x$};
\draw[dashed] (1.4,-1.5) -- (1.4,1.7);
\draw (1.4,-1.5) node[below]{$y$};

\clip (-.6,-1.3) rectangle (4.8,1.5);
\foreach\j in {0,...,1}{
\draw[fill=blue,opacity=.25] (-.5+2.8*\j,-.5) rectangle (.5+2.8*\j,.5);
}
\foreach\i in {-3,...,2}{
\foreach\j in {-1,...,2}{
\draw[fill=blue,opacity=.20] (-.5+\i+2.8*\j,-.5+.8*\i) rectangle (.5+\i+2.8*\j,.5+.8*\i);
\draw (-.5+\i+2.8*\j,-.5+.8*\i) rectangle (.5+\i+2.8*\j,.5+.8*\i);
\draw[fill=gray,opacity=.15] (-.5+\i+2.8*\j,-.5+.8*\i+1.2) rectangle (.5+\i+2.8*\j,.5+.8*\i+1.2);
\draw (-.5+\i+2.8*\j,-.5+.8*\i+1.2) rectangle (.5+\i+2.8*\j,.5+.8*\i+1.2);
}}
\end{tikzpicture}
\caption[Second case of Lemma \ref{LemDeFin}]{Second case of Lemma \ref{LemDeFin}, in the case $\ell=3$: two distinct vertical translates of $W^{k+1} + \widetilde\Lambda_{k+1}$ intersect (the first translate contains the dark blue thickening of $W^k + \widetilde\Lambda_k$, the second is represented in grey).}\label{FigLemDeFin2}
\end{center}
\end{minipage}
\end{figure}

The following formula makes the link between what happens in the $n$ last and in the $n$ penultimates coordinates of $\R^{n(k+1)}$:
\begin{equation}\label{EqChangeDim}
W^{k+1} + \widetilde\Lambda_{k+1} + \big(0^{(k-1)n},0^n,w_{i_j}\big) = W^{k+1} + \widetilde\Lambda_{k+1} + \big(0^{(k-1)n},-v_j,w_{i_j} + Bv_j\big),
\end{equation}
(we add a vector belonging to $\widetilde\Lambda_{k+1}$). 

We now apply the hypothesis of the lemma to the vectors $-v_1,\cdots,-v_{\ell}$: there exists $j\neq j'\in\llbracket 1,\ell\rrbracket$ such that
\begin{equation}\label{EqInter}
D_c\bigg(\Big(W^k + \widetilde\Lambda_k + (0^{(k-1)n} , -v_j)\Big) \cap \Big(W^k + \widetilde\Lambda_k + (0^{(k-1)n} , -v_{j'})\Big)\bigg) \ge D_0.
\end{equation}
Let $y$ be a point belonging to this intersection. Applying Equations~\eqref{EqIxIn} and \eqref{EqInter}, we get	that
\begin{equation}\label{EqLemDeFinCentr}
(y,x) \in W^{k+1} + \big(\widetilde\Lambda_k,0^n\big) + \big(0^{(k-1)n}, -v_j, w_{i_j} + Bv_j\big)
\end{equation}
and the same for $j'$.

Two different cases can occur.
\begin{enumerate}[(i)]
\item Either $i_j = i_{j'}$ (that is, the translation vectors $w_{i_j}$ and $w_{i_{j'}}$ are equal). As a consequence, applying Equation~\eqref{EqLemDeFinCentr}, we have
\begin{align*}
(y,x) + \big(0^{(k-1)n}, v_j, -Bv_j-w_{i_j}\big) \in & \Big( W^{k+1} + \big(\widetilde\Lambda_k,0^n\big) \Big)\cap\\
                                                     & \Big( W^{k+1} + \big(\widetilde\Lambda_k,0^n\big) + v'\Big),
\end{align*}
with
\[v' = \big(0^{(k-1)n}, -(v_{j'}-v_j), B(v_{j'}-v_j)\big) \in \widetilde\Lambda_{k+1} \setminus\widetilde\Lambda_k.\]
This implies that the set $W^{k+1} + \widetilde\Lambda_{k+1}$ auto-intersects (see Figure~\ref{FigLemDeFin}).
\item Or $i_j \neq i_{j'}$ (that is, $w_{i_j}\neq w_{i_{j'}}$). Combining Equations~\eqref{EqLemDeFinCentr} and \eqref{EqChangeDim} (note that $\big(\widetilde\Lambda_k,0^n\big) \subset \widetilde\Lambda_{k+1}$), we get
\[(y,x)\in \Big(W^{k+1} + \widetilde\Lambda_{k+1} + \big(0^{kn}, w_{i_j}\big)\Big) \cap \Big(W^{k+1} + \widetilde\Lambda_{k+1} + \big(0^{kn},w_{i_{j'}}\big)\Big).\]
This implies that two distinct vertical translates of $W^{k+1} + \widetilde\Lambda_{k+1}$ intersect (see Figure~\ref{FigLemDeFin2}).
\end{enumerate}

We now look at the global behaviour of all the $x$ such that $\overline x\in U$. Again, we have two cases.
\begin{enumerate}[(1)]
\item Either for more than the half of such $x$ (for Lebesgue measure), we are in the case (i). To each of such $x$ corresponds a translation vector $w_i$. We choose $w_i$ such that the set of corresponding $x$ has the biggest measure; this measure is bigger than $\varep/\big(2(\ell-1)\big)\ge \varep/(2\ell)$. Reasoning as in the second proof of Theorem~\ref{PerLin1} (page~\pageref{Proof2PerLin1}), and in particular applying Corollary~\ref{CoroSansNom}, we get that the density $D_1$ of the auto-intersection of $W^{k+1} + \widetilde\Lambda_{k+1} + (0,w_i)$ is bigger than $D_0\varep/(2\ell)$. This leads to
\[D_c(W^{k+1} + \widetilde\Lambda_{k+1}) < D_c(W^{k} + \widetilde\Lambda_{k}) - \frac{D_0 \varep}{4\ell}.\]
In this case, we get the conclusion (1) of the lemma.
\item  Or for more than the half of such $x$, we are in the case (ii). Choosing the couple $(w_i,w_{i'})$ such that the measure of the set of corresponding $x$ is the greatest, we get 
\[D_c\bigg(\Big(W^{k+1} + \widetilde\Lambda_{k+1} + (0^{kn} , w_i) \Big) \cap \Big(W^{k+1} + \widetilde\Lambda_{k+1} + (0^{kn} , w_{i'}) \Big)\bigg)\ge \frac{D_0 \varep}{(\ell-1)(\ell-2)}.\]
In this case, we get the conclusion (2) of the lemma.
\end{enumerate}
\end{proof}

We can now prove Theorem~\ref{ConjPrincip}.

\begin{proof}[Proof of Theorem~\ref{ConjPrincip}]
As in the previous proofs of such results, we proceed by induction on $k$. Suppose that $\widetilde\Lambda_k$ is such that $D_c(W^k + \widetilde\Lambda_k)>1/\ell$. Then, Lemma~\ref{EstimTauxN} ensures that it is not possible to have $\ell$ disjoint translates of $W^k + \widetilde\Lambda_k$. Applying Lemma~\ref{LemDeFin}, we obtain that either $D_c(W^{k+1} + \widetilde\Lambda_{k+1})<D_c(W^k + \widetilde\Lambda_k)$, or 
it is not possible to have $\ell-1$ disjoint translates of $W^{k+1} + \widetilde\Lambda_{k+1}$. And so on, applying Lemma~\ref{LemDeFin} at most $\ell-1$ times, there exists $k'\in \llbracket k+1,k+\ell-1\rrbracket $ such that $W^{k'} + \widetilde\Lambda_{k'}$ has additional auto-intersections. Quantitatively, combining Lemmas~\ref{EstimTauxN} and \ref{LemDeFin}, we get
\[D_c\big(W^{k+\ell-1} + \widetilde\Lambda_{k+\ell-1}\big) \le D\big(W^k + \widetilde\Lambda_k\big) - \frac{\varep}{4\ell} \left(\frac{\varep}{\ell^2}\right)^{\ell-1} 2\frac{\ell D_c(W^k + \widetilde\Lambda_k)-1}{\ell(\ell-1)},\]
thus
\[D_c\big(W^{k+\ell-1} + \widetilde\Lambda_{k+\ell-1}\big) - 1/\ell \le  \left( 1 -\frac12 \left(\frac{\varep}{\ell^2}\right)^\ell\right) \Big(D_c\big(W^k + \widetilde\Lambda_k\big) - 1/\ell \Big),\]
in other words, if we denote $\overline\tau^k = \overline\tau^k(B_1,\cdots,B_k)$ and $\lambda_\ell = 1 - \frac12 \left(\frac{\varep}{\ell^2}\right)^\ell$,
\begin{equation}\label{EqFinChap}
\overline\tau^{k+\ell-1} - 1/\ell \le  \lambda_\ell\Big( \overline\tau^k - 1/\ell \Big).
\end{equation}
This implies that for every $\ell>0$, the sequence of rates $\overline\tau^k$ is smaller than a sequence converging exponentially fast to $1/\ell$: we get Equation~\eqref{EstimTauxExp}. In particular, the asymptotic rate of injectivity is generically equal to zero.
\bigskip

We now prove the estimation of Equation~\eqref{EqConvLog}. Suppose that $\overline\tau^k \in [1/(\ell-2),1/(\ell-1)]$, we compute how long it takes for the rate to be smaller than $1/(\ell-1)$. We apply Equation~\eqref{EqFinChap} $j$ times to $\ell$ and get
\[\overline\tau^{k+j\ell} -1/\ell \le \lambda_\ell^j(\overline\tau^{k} -1/\ell),\]
with $\lambda = 1 - \frac12 \left(\frac{\varep}{\ell^2}\right)^\ell$. In the worst case, we have $\overline\tau^k = 1/(\ell-2)$, thus if $j'$ satisfies
\begin{equation}\label{eqj'}
\frac{1}{\ell-1} - \frac{1}{\ell} = \lambda_\ell^{j'}\left(\frac{1}{\ell-2}-\frac{1}{\ell}\right),
\end{equation}
then $j=\lceil j' \rceil$ is such that $\overline\tau^{k+j\ell}\le 1/(\ell-1)$. Equation~\eqref{eqj'} is equivalent to
\[j' = \frac{-1}{\log \lambda_\ell}\left( \log 2 - \log\left(1-\frac{1}{\ell-1}\right) \right).\]
And for $\ell$ very big, we have the equivalent (recall that $\lambda_\ell = 1 - \frac12\left(\frac{\varep}{\ell^2}\right)^\ell$)
\[j' \sim \left(\frac{\ell^2}{\varep}\right)^\ell 2\log 2.\]
Thus, when $\ell$ is large enough, the time it takes for the rate to decrease from $1/(\ell-2)$ to $1/(\ell-1)$ is smaller than $\ell^{ 2(\ell+1)} = e^{ 2(\ell+1)\log \ell}$. 

On the other hand, if we set $f(k) = (\log k)^{-\nu}$, the time it takes for $f$ to go from $1/(\ell-2)$ to $1/(\ell-1)$ is equal to
\[e^{(\ell-1)^{1/\nu}} - e^{(\ell-2)^{1/\nu}} = e^{(\ell-1)^{1/\nu}}\left(1 - e^{(\ell-2)^{1/\nu} - (\ell-1)^{1/\nu}}\right) \sim e^{(\ell-1)^{1/\nu}}\]
when $l$ goes to infinity, thus smaller than $e^{(\ell-1)^{1/\nu}+1}$ when $\ell$ is large enough. But we have $2(\ell+1)\log \ell = o\big((\ell-1)^{1/\nu}+1\big)$. So, when $\ell$ is large enough, it takes arbitrarily much more time for $\overline\tau^k$ to decrease from $1/(\ell-2)$ to $1/(\ell-1)$ than for $f$ to decrease from $1/(\ell-2)$ to $1/(\ell-1)$. As a consequence, $\overline\tau(k) = o(f(k))$.
\end{proof}

\begin{rem}
Fix a probability measure $\mu$ on $SL_n(\R)$ whose support has nonempty interior. Then, for almost every sequence of independent $\mu$-identically distributed matrices, the asymptotic rate of injectivity is zero. This is an easy consequence of the fact that the (open and $\delta$-dense) set $\mathcal O$ in Lemma~\ref{LemDeFin} is independent from the matrices $A_1,\cdots,A_k$, from the 0-1 law and from the continuity of the rate of injectivity on a full Lebesgue measure set.
\end{rem}

In the next part, we will need another technical statement, whose proof is based on Lemma~\ref{ConjPrincip}.

\begin{lemme}\label{DerTheoPart2}
For every $R_0>0$ and $\delta>0$, there exists $k_0\in\N$ such that the set $\mathcal O_\varep^{k_0}$ of sequences $\{(A_k)_{k\ge 1}\in\ell^\infty(SL_n(\R))$ such that there exists a sequence $(w_k)_{k\ge 1}$ of translation vectors belonging to $[-1/2,1/2]^n$, and a vector $\widetilde y_0\in\Z^n$, with norm bigger than $R_0$, such that ($\pi(A+w)$\index{$\pi(A+w)$} denotes the discretization of the affine map $A+w$)
\[\big(\pi(A_{k_0} + w_{k_0})\circ \cdots \circ \pi(A_1 + w_1)\big)(\widetilde y_0) = \big(\widehat{A_{k_0} + w_{k_0}}\circ \cdots \circ \widehat{A_1 + w_1}\big)(0) =0.\]
Moreover, the point $\widetilde y_0$ being fixed, this property can be supposed to remain true on a whole neighbourhood of the sequence $(A_k)_{k\ge 1}\in\mathcal O_\varep^{k_0}$.
\end{lemme}

\begin{proof}[Proof of Lemma~\ref{DerTheoPart2}]
We set
\[ \mathcal O_\varep^k = \{(A_k)_{k\ge 1}\in\ell^\infty(SL_n(\R)) \mid \tau^k(A_1,\cdots,A_k) <\varep\}.\]
Lemma~\ref{ConjPrincip} states that for every $\varep>0$, the set $\bigcup_{k\ge 0} \mathcal O_\varep^k$ contains an open and dense subset of $\ell^\infty(SL_n(\R))$. Together with the continuity of $\tau^k$ at every generic sequence (Remark~\ref{conttaukk}), this implies that for every $\delta>0$, there exists $k_0>0$ such that $\mathcal O_\varep^{k_0}$ contains an open and $\delta$-dense subset of $\ell^\infty(SL_n(\R))$.

Then, if $\tau^{k_0}(A_1,\cdots,A_k) <\varep$, then there exists a point $x_0\in\Z^n$ such that
\[\card\big((A_{k_0}\circ\cdots\circ A_1)^{-1}(x_0)\big) \ge \frac{1}{\varep}\]
(and moreover if the sequence $(A_k)_{k\ge 1}$ is generic, then this property remains true on a whole neighbourhood of the sequence). The lemma follows from this statement by remarking that on the one hand, if we choose $w_k\in [-1/2,1/2]^n$ such that
\[w_k = A_k^{-1}\Big(\big(\widehat{A_{k_0}}\circ \cdots \circ \widehat{A_{k-1}}\big)^{-1}(x_0)\Big) \mod \Z^n,\]
then the properties of the cardinality of the inverse image of $x_0$ are transferred to the point 0, and that on the other hand, for every $R_0>0$, there exists $m\in\N$ such that every subset of $\Z^n$ with cardinality bigger than $m$ contains at least one point with norm bigger than $R_0$.
\end{proof}

\chapter{Statistics of roundoff errors}\label{Stat0}

In this (very) short chapter, we use the formalism of model sets to study the statistics of the errors made when we compute the images of a point of $\Z^n$ by the discretizations of a generic sequence of linear maps. Here, the main result is that when we consider all the points of $\Z^n$, the roundoff errors made at each iteration are equidistributed in $[-1/2,1/2]^n$ (Proposition~\ref{EquidistribErr}). This result had already been obtained by P.P.~Flockermann in \cite[Theorem~10 page 44]{Flocker} in dimension 1 (see Proposition~\ref{LinPourLinstant}), with quite different techniques. From this result, we deduce the statistics of the cumulative errors made after $k$ iterations. In particular, in dimension $n=1$, it allows us to compute the \emph{discrepancy} (see Definition~\ref{Discrepancy}) between the appropriate uniform measure on $\R$ and the image sets (Proposition~\ref{EstimDiscrepancy}). 

Ultimately, we hope that these notions can be used to tackle a conjecture of O.E. Lanford (Conjecture~\ref{Lalan}) concerning the physical measures of expanding maps of the circle.

\section{Roundoff errors}


Given $x\in \Z^n$ and a sequence $(A_m)_{m\ge 1}$ of invertible matrices of $\R^n$, we want to compute the sequence\index{$\varep_x^m$} $(\varep_x^m)_m$ of roundoff errors made at each iteration. They are defined by
\[\varep_x^m = \big(\widehat A_m - A_m\big)\big((\widehat A_{m-1} \circ \cdots \circ \widehat A_1)(x)\big)\in [-1/2,1/2]^n.\]
We also set\index{$v_x^m$}
\[v_x^m = (\widehat A_m\circ \widehat A_{m-1} \circ \cdots \circ \widehat A_1)(x)\in\Z^n.\]
We fix $k\ge 0$ and set $v_x = \big(v_x^1,\cdots,v_x^k\big)$\index{$v_x$}, $\varep_x = \big(\varep_x^1,\cdots,\varep_x^k\big)$\index{$\varep_x$}, $u_x = (A_1x,0^{(n-1)k})\in\R^{nk}$\index{$u_x$} and\index{$N_{A_1,\cdots,A_k}$}
\[N_{A_1,\cdots,A_k} = \begin{pmatrix}
-1  &     &        &     & \\
A_2 & -1  &        &     & \\
    & A_3 & \ddots &     & \\
    &     & \ddots & -1  & \\
    &     &        & A_k & -1
\end{pmatrix}\in M_{nk}(\R).\]

As 
\[N_{A_1,\cdots,A_k} v_x = \begin{pmatrix}
-v_x^1 \\
A_2 v_x^1 - v_x^2\\
A_3 v_x^2 - v_x^3\\
\vdots\\
A_k v_x^{k-1} - v_x^k
\end{pmatrix},\]
the vector $u_x$ can be decomposed into
\[u_x  = N_{A_1,\cdots,A_k}v_x - \varep_x,\]
with $v_x \in \Z^{nk}$ and $\varep_x \in W^k$ (recall that $W^k = ]-1/2,1/2]^{nk}$). The vector $u_x$ being fixed, this condition characterizes completely $\varep_x$ and $v_x$, as $W^k$ is a fundamental domain of $N_{A_1,\cdots,A_k}\Z^{nk}$ (remark that we are in the case of the conclusion of Haj\'os  theorem, see Theorem~\ref{hajos}). Thus, $\varep_x$ is equal to the projection of $u_x$ on $W^k$ modulo $N_{A_1,\cdots,A_k}\Z^{nk}$; equivalently, $N_{A_1,\cdots,A_k}^{-1}\varep_x$ is equal to the projection of $N_{A_1,\cdots,A_k}^{-1} u_x$ on $N_{A_1,\cdots,A_k}^{-1} W^k$ modulo $\Z^{nk}$.

This implies that the sequences of errors $\varep_x$ are equidistributed in $(\R^n/\Z^n)^k$ when $x$ ranges over $\Z^n$ if and only if the vectors $\big(N_{A_1,\cdots,A_k}^{-1} u_x\big)_{x\in\Z^n}$ are equidistributed modulo $\Z^{nk}$. For this purpose, the matrix $N_{A_1,\cdots,A_k}^{-1}$ can be easily computed:
\[N_{A_1,\cdots,A_k}^{-1} = \begin{pmatrix}
-1  &      &        &        & \\
-A_2 & -1 &    &        & \\
-A_3A_2  & -A_3 & \ddots &        & \\
\vdots   & \vdots & \ddots & -1 & \\
{}\ -A_k \cdots A_2\ {} & {}\ -A_k \cdots A_3\ {} & \cdots & -A_k & -1
\end{pmatrix},\]
thus
\begin{equation}\label{action!}
N_{A_1,\cdots,A_k}^{-1} u_x = -\begin{pmatrix}
A_1\\
A_2 A_1\\
A_3 A_2 A_1\\
\vdots\\
A_k \cdots A_1
\end{pmatrix}x.
\end{equation}
As a consequence, the sequences of errors $\varep_x$ are equidistributed in $(\R^n/\Z^n)^k$ if and only if the action of $\Z^n$ on $(\R^n/\Z^n)^k$ given by Equation~\eqref{action!} is ergodic. This leads to the following proposition.

\begin{prop}\label{EquidistribErr}
For a generic sequence $(A_k)_{k\ge 1}$ of matrices of $GL_n(\R)$, or $SL_n(\R)$, or $O_n(\R)$, for every fixed integer $k$, the finite sequence of errors $\varep_x= \big(\varep_x^1,\cdots,\varep_x^k\big)$ is equidistributed in $(\R^n/\Z^n)^k$ when $x$ ranges over $\Z^n$.
\end{prop}

In particular, this proposition implies that the errors $\varep^k$ are globally independent.

If we take $n=1$, we obtain an alternative proof of the following statement, which was already observed in the thesis \cite{Flocker} of P.P.~Flockermann (Theorem~10 page 44). 

\begin{prop}\label{LinPourLinstant}
If $n=1$, we denote $A_m = \lambda_m$ and for $m\le k$, set $\overline\lambda_m^k = \lambda_k\lambda_{k-1}\cdots\lambda_m$, with the convention $\overline\lambda_{k+1}^k = 1$. If the coefficients $(\overline\lambda_m^k)_{1\le m \le k+1}$ are linearly independent over $\Q$, then for every fixed integer $k$, the sequence $\big(\varep_x^1,\cdots,\varep_x^k\big)$ is equidistributed in $\R^k/\Z^k$ when $x$ ranges over $\Z$.
\end{prop}


\section{Cumulative errors and discrepancy}

From the previous study, it is possible to deduce the statistics of the global error
\[\mathcal E_x^k = \big(\widehat A_k \circ\widehat A_{k-1} \circ \cdots \circ \widehat A_1\big)(x) - \big(A_k \circ A_{k-1} \circ \cdots \circ A_1\big)(x)\]
made after $k$ iterations. Indeed, we have 
\[\mathcal E_x^{k+1} = A_{k+1} \mathcal E_x^k + \varep_x^{k+1}.\]
From this recurrence relation, we deduce that
\[\mathcal E_x^k = \sum_{m=1}^k B_m \varep_x^m,\]
where $B_m = A_k A_{k-1} \cdots A_{i+1}$. As the $\varep_x^m$ are independent and equidistributed, this gives the law of the global error $\mathcal E_k$. In particular, the covariance of $\mathcal E_k$ is equal to
\[\operatorname{Var} (\mathcal E_x^k) = \frac{1}{12}\sum_{m=1}^k B_m B_m^\top,\]
where $B^\top$ denotes the transpose of the matrix $B$.
\bigskip

We assume $n=1$, and that the linear maps $A_i$ are expanding. In this case, if we denote $A_i = \lambda_i\ge 1$, the variance of $\mathcal E_x^k$ is equal to
\begin{equation}\label{Variance}
\operatorname{Var} (\mathcal E_x^k) = \frac{1}{12}\sum_{m=1}^k (\lambda_k\lambda_{k-1}\cdots \lambda_{m})^2.
\end{equation}
In particular, if there exists $\alpha>1$ such that $\lambda_k\ge \alpha$ for every $k$, then 
\[\operatorname{Var} (\mathcal E_x^k) \ge \frac{1}{12}\sum_{m=1}^k \alpha^{2m} = \frac{\alpha^2 (\alpha^{2k}-1)}{12(\alpha-1)}.\]

These considerations allow us to compute the discrepancy of the image set $(\widehat A_k \circ \cdots \circ \widehat A_1)(x)$.
\begin{definition}\label{Discrepancy}
Let $E\subset\Z$ and $\mu$ a Borel measure on $\R$. We call \emph{discrepancy of $E$} with respect to the measure $\mu$ the quantity (when it is well defined)\index{$\Disc(E,\mu)$}
\[\Disc(E,\mu) = \lim_{R\to +\infty} \left(\frac{1}{2R} \int_{-R}^R \big( \card([-x,x]\cap E) - \mu([-x,x]) \big)^2 \ud x\right)^{1/2}.\]
\end{definition}

For a complete investigation of the subject of geometric discrepancy and an extensive bibliography, see the excellent book \cite{MR2683232} of J.~Matou{\v{s}}ek.

In dimension 1, the discrepancy of the set $(\widehat A_k \circ \cdots \circ \widehat A_1)(\Z)$ can be easily computed: the following proposition is a direct consequence of Equation~\eqref{Variance}.

\begin{prop}\label{EstimDiscrepancy}
The discrepancy of the set $(\widehat A_k \circ \cdots \circ \widehat A_1)(\Z)$ is equal to the standard deviation of $(\mathcal E_x^k)_x$ :
\[\Disc\big((\widehat A_k \circ \cdots \circ \widehat A_1)(\Z),\, \det(A_k^{-1}\cdots A_1^{-1})\Leb\big) = \operatorname{Var} (\mathcal E_x^k)^{1/2} = \left(\frac{1}{12}\sum_{m=1}^k (\lambda_k\lambda_{k-1}\cdots \lambda_{m})^2\right)^{1/2}.\]
\end{prop}

In particular, if there exists $\alpha>1$ such that $\lambda_k\ge \alpha$ for every $k$, then 
\[\Disc\big((\widehat A_k \circ \cdots \circ \widehat A_1)(\Z),\, \det(A_k^{-1}\cdots A_1^{-1})\Leb\big) \ge \left(\frac{\alpha^2 (\alpha^{2k}-1)}{12(\alpha-1)}\right)^{1/2}\mathop{\sim}\limits_{k\to+\infty} \frac{\alpha^{k+1}}{\sqrt{12(\alpha-1)}}.\]

As has already been said, we hope that these kind of considerations about the statistics of the deviation of the uniform measure on the image set with respect to Lebesgue measure can help to understand the behaviour of the image measures $(f_N^*)^k(\lambda_N)$ for $k\gg \log N$ (but not too big either). In particular, we would like to use discrepancy to improve the results of Section~\ref{SecTransOp} and possibly explain the behaviour observed in Figure~\ref{GrafDistMesTemps}.

\part{Discretization of maps in higher regularity}\label{PartTri}

\parttoc

\chapter*[Introduction]{Introduction}

In this third part of the manuscript, we will study the dynamics of the discretizations of generic $C^r$-maps, namely both conservative and dissipative diffeomorphisms, and expanding maps. We conducted this study for several reasons.

First of all, it is commonly accepted that the generic $C^1$ dynamics represents some physical systems better than the $C^0$ generic dynamics: a lot of these concrete systems are smooth, and the generic $C^0$ dynamics contains ``wild'' behaviours (for example, for both generic conservative and dissipative homeomorphisms, whether the set of periodic points of a given period is empty, whether it forms a Cantor set, thus it is uncountable). Moreover, in the $C^1$ generic case, a wider class of behaviours can occur, as there are ``generic dichotomy theorems'' (see for example \cite{ArturSylvain}).

Furthermore, some aspects of the behaviour of the discretizations of generic conservative homeomorphisms are a bit disappointing. The combinatorics of the discretizations varies a lot depending on the order of discretization, and does not tell anything about the actual dynamics of the continuous system. In addition, nothing in the results we have proved in the first part gives a method to to detect the ``physical'' behaviour of the homeomorphism on the discretizations. We would like to know if in the $C^1$ case, it is possible to obtain a method to recover the physical measures from some dynamical features of the discretizations.
\bigskip

In this part, we will consider a compact boundaryless manifold $X$ of dimension $n\ge 2$ (we will also study the specific case of the circle), equipped with a measure $\lambda$ derived froma volume form. We will use the following notations.

\begin{notation}
We denote by $\Diff^1(X)$\index{$\Diff^1(X)$} the set of $C^1$-diffeomorphisms of $X$, endowed by the metric $d_{C^1}$ defined by:
\[d_{C^1}(f,g) = \sup_{x\in X} d\big(f(x),g(x)\big) + \sup_{x\in X} \big\|Df_x - Dg_x\big\|.\]
We denote by $\Diff^1(X,\lambda)$\index{$\Diff^1(X,\lambda)$} the subset of $\Diff^1(X)$ made of the diffeomorphisms that preserve the measure $\lambda$, endowed with the same metric $d_{C^1}$.
\end{notation}

The metric $d_{C^1}$ on $\Diff^1(X)$ and $\Diff^1(X,\lambda)$ makes them Baire spaces. In this introduction, we will state all the results for the phase space $X=\T^n$, the measure $\lambda=\Leb$ and the uniform grids
\[E_N = \left\{\left(\frac{i_1}{N},\cdots,\frac{i_n}{N}\right)\in \T^n \big|\ \forall j,\, 0\le i_j\le {N}-1\right\},\]
However, these results are true in a more general setting (see the concerned chapters for precisions on these settings).
\bigskip

The first chapter of this part will be devoted to the direct application of classical perturbation lemmas for $C^1$-diffeomorphisms.

These perturbations lemmas are in general much more difficult to obtain than in the case of homeomorphisms (proposition of finite maps extension, see Proposition~\ref{extension}): when we perturb a diffeomorphism, we have to take care of the norm of the differentials of the perturbation. This is a crucial difference. Indeed, in the $C^0$ case, if we want to perturb the image of a point of $\varep$, we only have to make a perturbation on a ball of size $\varep$. If we want to do the same in topology $C^1$, the support of the perturbation has to contain a ball whose size is way bigger than $\varep$. Thus, for homeomorphisms, and for $N$ large enough, we are able to move independently each point of a grid $E_N$, while for $C^1$ diffeomorphisms we are only authorized to move a small proportion\footnote{This proportion is asymptotically independent from the order $N$ of the discretization, and depends on the size of the $C^1$ perturbation.} of the points of $E_N$.

Thus, if the $C^1$ generic dynamics is mainly considered as more interesting than the $C^0$ dynamics, the price to pay is that our results are in general weaker: they only concern a small proportion of the points of the grids. Note that the case of the regularity $C^1$ is somehow the ``limit case'': to our knowledge, there is no general statement of interesting perturbation result in regularity $C^r$, with $r>1$.

As a consequence, the results we obtain in this chapter -- by applying directly classical perturbation lemmas -- only concern a sub-dynamics of the discretizations. More precisely, in Part~1, we were able to determine the dynamics of \emph{all} the points of the discretizations of generic conservative homeomorphisms. Here, we will only control the dynamics of a few points of the grid. For example, even if the dynamical property we are interested in concern a $\varep$-dense subset of $X$, it may happen that the mesh of the grids which satisfies the conclusions of the theorems is very small compared to $\varep$; in this case we only know the dynamics (for the discretizations) of a small proportion of points of these grids.
\bigskip

The results of this first chapter of Part 3 will be based on the following statement (Corollary~\ref{PropShadowC1}).

\begin{propo}
Let $f \in \Diff^1(\T^n)$ (or $f \in \Diff^1(\T^n,\Leb)$) be a generic diffeomorphism. Then, for every $\tau\in\N$, every $N_0\in\N$ and every $\varep>0$, there exists $N\ge N_0$ such that every periodic orbit of $f$ of period smaller than $\tau$ is $\varep$-shadowed by a periodic orbit of $f_N$ with the same period.
\end{propo}

This proposition states that every periodic orbit of a generic diffeomorphism is detected by an infinite number of discretizations. It will also allow us to deduce a lot of results about the dynamics of the discretizations of generic diffeomorphisms by applying classical perturbation lemmas. For example, applying the connecting lemma of C.~Bonatti and S.~Crovisier, \cite{MR2090361}, we get the following result (Corollary~\ref{typlaxC1}).

\begin{propo}\label{BbB}
Let $f$ be a generic diffeomorphism of $\Diff^1(\T^n,\Leb)$. Then, for any $\varep>0$ and any $N_0>0$, there exists $N\ge N_0$ such that $f_N$ has a periodic orbit which is $\varep$-dense.

The statement is still true in the dissipative generic case if we restrict to a maximal invariant chain-transitive set.
\end{propo}

This statement can be seen as a weak version of Corollary~\ref{typlax}, which states that for a generic conservative homeomorphism, there exists an infinite number of discretizations which are cyclic permutations of the grids. We can also apply an ergodic closing lemma (due to R.~Mañé \cite{MR678479} and F.~Abdenur, C.~Bonatti and S.~Crovisier \cite{MR2811152}); we then get the following result for discretizations (Corollary~\ref{CoroMane}).

\begin{propo}
Let $f$ be a generic diffeomorphism of $\Diff^1(\T^n,\Leb)$. Then, for any $f$-invariant measure $\mu$, any $\varep>0$ and any $N_0>0$, there exists $N\ge N_0$ such that $f_N$ supports a periodic measure which is $\varep$-close to $\mu$.

The statement is still true in the dissipative generic case for every measure $\mu$ supported by a maximal invariant chain-transitive set.
\end{propo}

Again, this result can be seen as a weak version of Theorem~\ref{mesinv} for homeomorphisms.
\bigskip

Finally, we can say that all these results of Chapter~\ref{ChapPerturbLem} go in the direction of the following heuristic.

\emph{For a generic diffeomorphism $f\in \Diff^1(\T^n,\Leb)$, each ``sub-dynamics'' of $f$ (periodic orbits, chain-transitive invariant compact sets, invariant measure, rotation set\dots) can be detected by some discretizations $f_N$. However, we have no control over the orders of discretization $N$ for which it is true, and no control over the global dynamics of $f_N$. In the dissipative case, the same holds on each chain-recurrent class of $f$.}
\bigskip

Remark that for diffeomorphisms, the behaviours in the dissipative case are quite close to those in the conservative case, which was not true for generic homeomorphisms. Indeed, for a generic dissipative homeomorphism, the chain recurrent classes are totally disconnected, while for generic dissipative diffeomorphisms they can contain nonempty open sets (for example, there are open sets of transitive Anosov diffeomorphisms). On each of these chain recurrent classes, the dynamics generically resembles to that of a conservative diffeomorphism.
\bigskip
	
Recall that the class of $C^1$-diffeomorphisms bears open sets of Anosov maps. One of the important properties of such systems is that they satisfy the shadowing lemma: if $f$ is Anosov, then for every $\varep>0$, there exists $\delta>0$ such that every $\delta$-pseudo-orbit of $f$ is $\varep$-shadowed by a true orbit of $f$. As orbits of discretizations are in particular pseudo-orbits, orbits of discretizations of Anosov maps are shadowed by real orbits of the initial map. However, we do not control the dynamics of the shadowing orbit: for example, if the Anosov map is ergodic, we can not impose to the shadowing orbit to be typical with respect to Lebesgue measure. This phenomenon could seem anecdotic, but its consequences are in fact quite bad. Recall, for example, the behaviour of the discretizations of the linear automorphisms of the torus on the canonical grids (see Figure~\ref{Miaou}): in this case, the discretizations are permutations with a very small global period; thus the dynamics of the 
discretizations does not reflect the mixing properties of the initial dynamics. Thus, we can say that the shadowing lemma does not imply that the dynamics of $f_N$ looks like that of $f$.

In a certain sense, the statements we prove in this first chapter express that for a generic diffeomorphism $f$, any dynamical behaviour of $f$ is shadowed by a similar dynamical behaviour of the discretizations of $f$. More precisely, by a sub-dynamics of the discretizations; the dynamics of the rest of the grid is not totally random, as by the shadowing lemma it is close to some dynamics of the initial system. The bad news is that in the theorems we prove, we have no explicit control on the orders of discretization for which we detect this or that dynamical feature of $f$, whereas in the shadowing lemma the dependence of the parameter $\delta$ to $\varep$ is explicit.
\bigskip

The second chapter of this part is devoted to the study of the degree of recurrence of a generic $C^1$-diffeomorphism, in both conservative and dissipative cases, and of a generic expanding map.

Indeed, in Chapter~\ref{ChapPerturbLem}, we do not say anything about the global dynamics of the discretizations; we obtain results about sub-dynamics of the discretizations. In Chapter~\ref{ChapDeg}, we will study the simplest global combinatorial quantity for the discretizations: the degree of recurrence. We say that this invariant is the simplest to study because it is obtained as a decreasing limit in time. Thus, the study of the degree of recurrence reduces to that of finite time quantities. This study will also be the occasion to observe that the local behaviour of the discretizations in small time is governed by the behaviour of the differentials of the diffeomorphism; in particular we will use crucially the study of the discretizations of linear maps we have conducted in Part 2 of this manuscript.

Recall that the degree of recurrence $D(f_N)$ of a homeomorphism $f$ is defined as the ratio between the cardinality of the recurrent set of the discretization $f_N$ and that of the grid $E_N$ (Definition~\ref{DefDegree}). This degree of recurrence somehow represents the amount of information we lose when we iterate the discretization. The study of the degree of recurrence of a generic diffeomorphism easily reduces to that of the rates of injectivity
(see Definition~\ref{DefTauxDiffeo})
\[\tau^t(f) = \underset{N\to +\infty}{\overline\lim} \frac{\card\big((f_N)^t(E_N)\big)}{\card(E_N)}.\]

We first show a local-global formula for this rate of injectivity: the rate of injectivity of a generic diffeomorphism is linked with the rates of injectivity of its differentials. To do that, we define the discretization of a linear map $A : \R^n \to \R^n$ as the map $\widehat A = \pi\circ A : \Z^n \to \Z^n$, where $\pi : \R^n\to\Z^n$ is a projection on (one of) the nearest integer point for the euclidean distance (see Definition~\ref{DefDiscrLin}). Then, the rate of injectivity of a sequence $(A_k)_{k\ge 0}$ of matrices is defined as (see Definition~\ref{DefTaux})
\[\tau^k(A_1,\cdots,A_k) = \lim_{R\to +\infty} \frac{\card \big((\widehat{A_k}\circ\cdots\circ\widehat{A_1}) [B_R]\big)}{\card [B_R]}\in]0,1].\]
The local-global formula is the following (Theorem~\ref{convBis}, see also Theorem~\ref{convBisMieux}).

\begin{theorem}\label{Hue}
Let $r\ge 1$ and $f\in \Diff^r(\T^n)$ (or $f\in \Diff^r(\T^n,\Leb)$) be a generic diffeomorphism. Then $\tau^t(f)$ is well defined (that is, the limit exists) and satisfies:
\[\tau^t(f) = \int_{\T^n} \tau^t\left(Df_x, \cdots, D f_{f^{t-1}(x)}\right) \ud \Leb(x).\]
Moreover, the function $\tau^t$ is continuous in $f$.
\end{theorem}

Remark that the hypothesis of genericity is necessary to get this theorem (see Example~\ref{aoirfhaoijfaeij}). The proof of this result involves the local linearization of a diffeomorphism (Lemma~\ref{morceaux}), which is a difficult result in the conservative case (it uses the smoothing of a conservative $C^1$-diffeomorphism of A.~Avila, see \cite{MR2736152}). It also uses crucially the study of the continuity of the rate we have conducted in Part 2 of this manuscript.

In the last section of Chapter~\ref{ChapDeg}, we give a variation of Theorem~\ref{Hue} for expanding maps. This statement links the rate of injectivity of a generic $C^r$-expanding map with the probability of percolation $\overline D$ of a random graph associated to the derivatives of this map (Theorem~\ref{TauxExpand}).

\begin{theorem}\label{Perec}
For every $r\ge 1$, and every generic $C^r$-expanding map $f$ on $\Sp^1$, the rate of injectivity $\tau^k(f)$ satisfies
\[\tau^k(f) = \int_{\T^n} \overline D\big( (\det Df_x^{-1})_{\begin{subarray}{l} 1\le m\le k \\ x\in f^{-m}(y)\end{subarray}}\big) \ud \Leb(y)\]
(see Definition~\ref{Noel!} for the definition of $\overline D$).
\end{theorem}

The techniques developed to prove this statement also lead to the proof of Theorem~\ref{Hue} in the $C^r$ regularity, for every $r\ge 1$. Notice that even if theorem~\ref{Hue} is true for generic $C^r$-diffeomorphisms for every $r\ge 1$, it is not sufficient to get that the degree of recurrence is generically equal to 0 in the $C^r$ topology: to prove Theorem~\ref{Souvenir}, we need to perturb the derivative of a diffeomorphism on a big subset of $\T^n$. However, there is a hope that in the case of expanding maps, Theorem~\ref{Perec} leads to the fact that the degree of recurrence is generically equal to 0 in any topology of $C^r$ convergence.

The behaviour described by Theorem~\ref{Hue} is very different from the $C^0$ case. For a generic dissipative homeomorphism $f$, we have $\tau^1(f) = 0$ (Corollary~\ref{totsingdiscr}), while for a generic conservative homeomorphism $f$, for every $t>0$, the quantity $\card\big((f_N)^t(E_N)\big)/\card(E_N)$ accumulates on both $0$ and $1$ when $N$ goes to $+\infty$ (Corollaries~\ref{typlax} and \ref{crush}). The behaviour for generic diffeomorphisms does not depend on the assumption of preservation of a measure or not, and is in a certain sense smoother. 
\bigskip

In short, the rate of injectivity of a diffeomorphism can be computed by integrating the rates of injectivity of its differentials. This allows us to use the study of the rate of injectivity of linear maps we have conducted in the second part of this thesis, and in particular Theorem~\ref{ConjPrincip}. Applying classical techniques of $C^1$ dynamics (in particular Rokhlin tower lemma), the application of this theorem is quite straightforward and leads to the following theorem.

\begin{theorem}\label{Souvenir}
For a generic conservative diffeomorphism $f\in\Diff^1(\T^n,\Leb)$, we have
\[\lim_{t\to\infty}\tau^t(f)=0.\]
This implies in particular that $\lim_{N\to+\infty}D(f_N)=0$.
\end{theorem}

This theorem must be compared with the case of homeomorphisms (see Corollary~\ref{ConjEt}): for a generic conservative homeomorphism $f$, the sequence $D(f_N)$ accumulates on the whole segment $[0,1]$. Thus, in the case of diffeomorphisms, the behaviour of the rate of injectivity is less irregular than for homeomorphisms; despite this, the theorem shows that when we iterate the discretizations of a generic conservative diffeomorphism, we lose a great amount of information. Moreover, although $f$ is conservative, its discretizations tend to behave like dissipative maps. This can be compared with the work of P.~Lax \cite{MR0272983}: for any conservative homeomorphism $f$, there is a bijective finite map arbitrarily close to $f$ (see Theorem~\ref{Lax}). The previous theorem states that for a generic conservative $C^1$ diffeomorphism, the discretizations never possess this property.

To summarize, for a generic conservative $C^1$-diffeomorphism, the recurrent sets $\Omega(f_N)$ are such that their limit superior is $\T^n$ (see Proposition~\ref{BbB}), and their cardinality behaves as $o(\card(E_N))$.

The same behaviour of the rate of injectivity holds in the generic dissipative case (Corollary~\ref{CoroCoroArturJairo}): for a generic dissipative diffeomorphism $f\in\Diff^1(\T^n)$, we have $\lim_{N\to+\infty}D(f_N)=0$. This result is an easy consequence of a theorem of A.~Avila and J.~Bochi (Theorem~\ref{ArturJairo}, see also \cite{MR2267725}).
\bigskip

Finally, in Chapter~\ref{chapPhys}, we tackle the question of the physical behaviour of the discretizations. In section~\ref{SecPhys12}, we study the measures $\mu_{x}^{f_{N_k}}$ for ``a lot'' of points $x$. Recall that for every $x$, the orbit of $x_N$ under $f_N$ falls in a periodic orbit of $f_N$; we denote by $\mu_{x}^{f_{N_k}}$ the uniform measure on this periodic orbit. In the third chapter of this part, we will prove the following statement (Theorem~\ref{TheoMesPhysDiff}).

\begin{theorem}\label{TheomesPhysDiffIntroChap}
For a generic diffeomorphism $f\in\Diff^1(\T^n,\Leb)$, for a generic point $x\in X$, and for any $f$-invariant probability measure $\mu$, there exists a subsequence $(N_k)_k$ of discretizations such that
\[\mu_{x}^{f_{N_k}} \underset{k\to+\infty}{\longrightarrow} \mu.\]
\end{theorem}

Again, the study of the linear case (more precisely, Lemma~\ref{DerTheoPart2}) plays a pivotal role in the proof of this theorem. It also uses an ergodic closing lemma adapted from \cite{MR2811152} (Lemma~\ref{ErgoLemPlus}), and the connecting lemma for pseudo-orbits of \cite{MR2090361} (see also Theorem~\ref{PapySylvain}).

Recall that in the case of homeomorphisms, we have proved that for a generic conservative homeomorphism $f$, and for any $f$-invariant measure $\mu$, there exists a subsequence $(f_{N_k})_k$ of discretizations such that for every $x\in\T^n$, we have $\mu_{x}^{f_{N_k}} \underset{k\to+\infty}{\longrightarrow} \mu$. The result we have for diffeomorphisms -- even if it is much more difficult to prove -- is much weaker: it not only concerns the behaviour of a generic subset of the torus instead of all the points, but the orders of discretization for which the conclusions of the theorem are true strongly depend on the point $x$.

Despite this, this theorem says that in practice, we must be very careful when we want to find a physical measure. Suppose that we want to find the unique physical measure $\mu_0$ of a generic $C^1$-diffeomorphism $f$. We could think that it is sufficient to compute the measures $\mu_{x}^f$ for a certain number of points $x$; since $\mu_0$ is the unique physical measure, $\mu_{x}^f$ should be close to $\mu_0$ for most of the points $x$. In practice, of course, we will consider a finite number of such points $x$, and we will compute with a given precision, that is we will iterate a discretization $f_N$. The result says that for some $N$, we observe a measure $\mu_{x}^{f_{N_k}}$ that has nothing to do with the physical measure $\mu_0$.

Note that Theorem~\ref{TheomesPhysDiffIntroChap} does not say anything about the measures $\mu_{\T^n}^{f_N}$. Recall that these measures are defined as the average over $x\in E_N$ of the measures $\mu_x^{f_N}$; equivalently $\mu_{\T^n}^{f_N}$ is supported by the union of periodic orbits of $f_N$, and the total measure of each of these periodic orbits is proportional to the size of its basin of attraction under $f_N$. For now, the theoretical study of the measures $\mu_{\T^n}^{f_N}$ for generic conservative $C^1$-diffeomorphisms seems quite hard, as this kind of questions is closely related to the still open problem of genericity of ergodicity among these maps (see \cite{ArturSylvain} for the most recent advances on this topic).

On numerical simulations of these measures $\mu_{\T^n}^{f_N}$, it is not clear whether they converge towards Lebesgue measure or not (see Figures~\ref{MesC1IdCons2p}, \ref{MesC1AnoCons2p} and \ref{MesC1AnoConsSerie}). However, one can hope that their behaviour is not as erratic as for generic conservative homeomorphisms, where they accumulate on the whole set of $f$-invariant measures (Theorem~\ref{mesinv}). Indeed, Theorem~\ref{Souvenir} shows that the degree of recurrence $D(f_N)$ converges to 0 as $N$ tends to $+\infty$ for a generic conservative $C^1$-diffeomorphism, while it accumulates on the whole segment $[0,1]$ for generic conservative homeomorphisms (Corollary~\ref{ConjEt}). This shows that the \emph{global} dynamics of discretizations of generic conservative $C^1$-diffeomorphism and homeomorphisms might be very different.
\bigskip

Section~\ref{SecTransOp} is devoted to the study of the measures $(f_N^*)^m \lambda_N$ for ``small'' times $m$ (recall that $\lambda_N$ is the uniform measure on the grid $E_N$), in the case where $f$ is a $C^{1+\alpha}$ expanding map of the circle. In this setting, a classical result asserts that $f$ has a single physical measure $\mu_0$, which coincides with the SRB measure of $f$, and is also the unique invariant probability measure with a $C^\alpha$ density (see Theorem~\ref{Liverani}). Our theorem asserts that if $N$ goes to infinity much faster than $m$, the measures $(f_N^*)^m \lambda_N$ converge to $\mu_0$ (Theorem~\ref{MainMoche}).

\begin{theorem}\label{lalalalalalalalalala}
For every $0<\alpha\le 1$ and every $C^{1+\alpha}$ expanding map $f\in \mathcal E_d^{1+\alpha}(\Sp^1)$, there exists a constant $c_0 = c_0(f)>0$ such that if $(N_m)_m$ is a sequence of integers going to infinity and satisfying $\log N_m > c_0 m$, then the convergence $(f_{N_m}^*)^m(\lambda_{N_m}) \to \mu_0$ holds.
\end{theorem}

This result answers partially to a conjecture of O.E.~Lanford (see Conjecture~\ref{Lalan}). Its proof consists in a calculus of the difference of action on measures between the Ruelle-Perron-Frobenius operator and the pushforwards by the discretization $f_N^*$.
\bigskip

In Section~\ref{NumSimPhys}, we first present numerical experiments simulating the measures $\mu_x^{f_N}$ and $\mu_{\T^2}^{f_N}$ for some examples of conservative $C^1$-diffeomorphisms of the torus. The results of these simulations are quite striking: for all the diffeomorphisms we have tested (even a diffeomorphism $C^1$-close to a linear Anosov map, thus $C^0$ conjugated to it), for some large orders of discretization $N$, the measures $\mu_x^{f_N}$ do not look like Lebesgue measure (see for instance Figures~\ref{MesPhysIdC1}, \ref{MesPhysRotC1} and \ref{MesPhysAnoC1}). Thus, in practice, we have to be very careful when we compute numerically some measures $\mu_{x,T}^{f_N}$, it may happen that they do not reflect the physical behaviour of the initial diffeomorphism at all.

The end of this section is devoted to numerical simulations of invariant measures of expanding maps; in particular we observe on Figure~\ref{GrafDistMesTemps} that as predicted by Theorem~\ref{lalalalalalalalalala}, the distance $\dist(\mu_0,(f_{N_m}^*)^m(\lambda_{N_m}))$ is very close to 0 for small values of $m$, but then increases with $m$ to a significant value.
\bigskip

To summarize, we take stock about what we have proved to answer the question: is is possible to observe physical measures on discretizations? In practice, this question can be formalized in many ways\footnote{Here, we consider only the case of the regularity $C^r$, $r\ge 1$; the case of homeomorphisms has been treated in Part~\ref{PartOne}.}.
\begin{itemize}
\item First, we can wonder if the measures $\mu_{x}^{f_N}$ tend to a physical measure of $f$ for ``many'' points $x$. The answer is no, when by ``many points $x$'' we mean ``generic points $x\in X$'', for generic conservative and dissipative $C^1$-diffeomorphisms (Theorems~\ref{TheoMesPhysDiff} and \ref{TheoMesPhysDiffDissip}) and generic expanding maps (Proposition~\ref{TheoMesPhysDiffExp}). The behaviours described by these statements can even be observed in practice (see Section~\ref{NumSimPhys})
However, we do not know what happens when we consider a set of points $x$ which is typical for Lebesgue measure.
\item Second, we can ask whether the measures\footnote{Recall that the measures $\mu_{X}^{f_N}$ are supported by the union of periodic orbits of $f_N$, such that the total measure of each periodic orbit is proportional to the cardinality of its basin of attraction.} $\mu_{X}^{f_N}$ converge to a physical measure of $f$ or not. We have no theoretical statement about this question, even for expanding maps of the circle. Moreover, the results of numerical simulations are not clear, as they suggest that the answer to this question could be ``yes, if we consider averages among a wide range of orders $N$'' (see Figures~\ref{GrafMesInvGro} and \ref{MoyMeasPlInv}).
\item As we do not understand the behaviours of the measures\footnote{Note that the measures $\mu_{X}^{f_N}$ are obtained by averaging the measures $(f_N^*)^m\lambda_N$ over $m$.} $\mu_{X}^{f_N}$, as suggested by O.E.~Lanford (see Conjecture~\ref{Lalan}), we can ask if the measures $(f_N^*)^m\lambda_N$ converge to a physical measure for both $m$ and $N$ going to infinity, with $m$ not too big with respect to $N$. Theorem~\ref{lalalalalalalalalala} answers this question for some $m$ satisfying $m=O(\log N)$ and for a $C^\alpha$ expanding map on the circle: in this case, the measures $(f_N^*)^m\lambda_N$ converge exponentially fast to the unique physical measure of the map. For now, we do not know what happens for bigger times $m$. The numerical experiments suggest a quite surprising evolution of the measures $(f_N^*)^m\lambda_N$: for a fixed order of discretization $N$, these measures converge very fast towards the SRB measure, and then moves away from it slowly. This phenomenon is maybe due to the fact that the numerical errors due to the discretization process are equidistributed (see Proposition~\ref{EstimDiscrepancy}; this phenomenon could be specific to the dimension 1).
\item Finally, we can investigate whether the two first questions are connected or not. More precisely, we wonder if there exits a discrete Birkhoff's theorem, which would state that if the measures $\mu_{X}^{f_N}$ converge to a measure $\mu$, then the measures $\mu_{x}^{f_N}$ tend to $\mu$ for ``a lot'' of points $x$ (in a sense to define).
\end{itemize}

\chapter{Applications of perturbation lemmas}\label{ChapPerturbLem}

In this chapter, we apply classical perturbation lemmas of $C^1$-diffeomorphisms to deduce statements about the dynamics of the discretizations of $C^1$ generic conservative and dissipative diffeomorphisms.
\bigskip

In this chapter, $X$ will be a smooth compact boundaryless Riemannian manifold of dimension $n\ge 2$. It will be equipped with a measure $\lambda$, derived from a volume form (independent from the metric on $X$). We will be interested in the dynamics of the discretizations of generic $C^1$-diffeomorphisms, in both dissipative and conservative cases; \emph{i.e.} generic elements of $\Diff^1(X)$ and $\Diff^1(X,\lambda)$.

We will also fix a sequence $(E_N)_{N\ge 0}$ of discretizations grids on $X$. Recall that the very definition of discretization grid (Definition~\ref{grillmiam}) supposes that the mesh of these discrete sets tends to 0, that is: for every $\varep>0$, there exists $N_0$ such that for every $N\ge N_0$, the grid $E_N$ is $\varep$-dense in $X$. In this chapter, this is the only assumption we will make on the grids. Given a diffeomorphism $f\in \Diff(X,\lambda)$, we denote by $f_N$ the discretization of $f$ with respect to the grid $E_N$.
\bigskip

We will begin this chapter by stating that the periodic orbits of a generic diffeomorphism are shadowed by periodic orbits of its discretizations (Corollary~\ref{PropShadowC1}). This will allow us to apply classical closing lemmas, which state that some dynamical invariants of a generic diffeomorphism are shadowed by periodic orbits.

We first apply the connecting lemma for pseudo-orbits of C.~Bonatti and S.~Crovisier \cite{MR2090361}, to prove that for a generic diffeomorphism $f\in\Diff^1(X,\lambda)$ and for every $\varep>0$, an infinite number of discretizations $f_N$ have a periodic orbit which is $\varep$-dense in $X$ (Corollary~\ref{typlaxC1}). We also apply a closing lemma of F.~Abdenur and S.~Crovisier \cite{MR2975581} to obtain that an infinite number of discretizations $f_N$ are ``$\varep$-topologically mixing'' (Corollary~\ref{CoroFlavSylvain2}). The same results holds in the dissipative case for every chain-recurrent class.

We then apply an ergodic closing lemma of R.~Mañé \cite{MR678479} and F.~Abdenur, C.~Bonatti and S.~Crovisier \cite{MR2811152} to state that if $f\in\Diff^1(X,\lambda)$ is a generic conservative diffeomorphism, then every $f$-invariant measure is shadowed by periodic measures of the discretizations $f_N$ (Corollary~\ref{CoroMane}). The same kind of results holds for chain-transitive invariant compact sets (Corollary~\ref{CoroSylvIHES}), by applying a connecting lemma of S.~Crovisier \cite{MR2264835}. Again, the same results holds in the dissipative case when restricted to a chain-recurrent class.

Finally, we will use a theorem of J.~Franks \cite{MR958891} of realization of rotation vectors by periodic points to prove that the rotation set of a generic conservative diffeomorphism $f\in\Diff^1(X,\lambda)$ is well approximated by rotation sets of discretizations (Proposition~\ref{DiscrC1}). The same holds for generic dissipative $C^1$-diffeomorphisms, provided that we take convex hulls of the rotation sets of the discretizations (Proposition~\ref{DiscrC1Dissip}).

\section{Shadowing of periodic orbits}

First of all, we give the statement of an elementary perturbation lemma that we will use all along this chapter. This is a $C^1$ counterpart of the proposition of extension of finite maps for the $C^0$ topology (Proposition~\ref{extension}). 

\begin{lemme}[Elementary perturbation lemma in $C^1$ topology]\label{PerturbElem}
For every diffeomorphism $f\in\Diff^1(X)$ and every $\delta>0$, there exists $\eta>0$ and $r_0>0$ such that the following property holds: for every $x,y\in X$ such that $d(x,y)<r_0$, there exists a diffeomorphism $g\in\Diff^1(X)$ satisfying $d_{C^1}(f,g)<\delta$, such that $g(x) = f(y)$ and that $f$ and $g$ are equal out of the ball $B\big(\frac{x+y}{2},\frac{1+\eta}{2}d(x,y)\big)$.

Moreover, if $f\in\Diff^1(X,\lambda)$, then we can also suppose that $g\in\Diff^1(X,\lambda)$.
\end{lemme}

Before discussing the consequences of this lemma, let us underline the differences with the $C^0$ case. \emph{A priori}, we will apply the lemma to a small neighbourhood of $f$, thus a small parameter $\delta>0$. In facts, the number $\eta>0$ tends to $+\infty$ when $\delta$ goes to 0; thus the ball $B\big(\frac{x+y}{2},\frac{1+\eta}{2}d(x,y)\big)$ in which lies the support of the perturbation is large with respect to the distance between $x$ and $y$. This phenomenon does not happen in the $C^0$ topology, where the proposition of extension of finite maps (Proposition~\ref{extension}) states that the size of the support of the perturbation is almost equal to the distance between $x$ and $y$ (in fact, it can be supposed to be included in the neighbourhood of a given path joining $x$ to $y$). Remark that the situation would be even worse in the $C^r$ topology, for $r>1$: in this case, the $C^r$ size of the perturbation $\delta>0$ being fixed, the ratio between the size of the support of the perturbation and $d(x,y)$ would go to infinity when $d(x,y)$ goes to 0.

Thus, the difference between the statements we get for homeomorphisms and for $C^1$-diffeomorphisms lies in this difference between the elementary perturbation lemmas: in the $C^1$ case, perturbations are less local than in the $C^0$ case. In the view of obtaining generic properties of discretizations, the impossibility of making perturbations of all the points of the discretization grids has the effect that we can not obtain \emph{global} properties of discretizations of generic diffeomorphisms with this strategy of proof.

\begin{figure}[t]
\begin{center}
\begin{tikzpicture}[scale=1]

\draw[color=gray] (0,-2) circle (2);
\draw[color=gray] (0,-2) circle (2.3);
\draw[color=gray] (0,-2) circle (1.7);
\draw[color=gray] (0,-2) circle (1.4);
\draw[->,>=stealth,thick] (2,-2) arc (0:30:2);
\draw[->,>=stealth,thick] (2.3,-2) arc (0:20:2.2);
\draw[->,>=stealth,thick] (2.6,-2) arc (0:10:2.6);
\draw[->,>=stealth,thick] (1.7,-2) arc (0:20:1.7);
\draw[->,>=stealth,thick] (1.4,-2) arc (0:10:1.4);

\draw (1,-.268)  node[fill=white, above right] {$x$};
\draw (-1,-.268)  node[fill=white, above left] {$y$};
\draw (1,-.268) node{$\times$};
\draw (-1,-.268) node{$\times$};

\draw[color=gray] (0,-2) circle (2.6);

\end{tikzpicture}
\caption[Proof of Lemma~\ref{PerturbElem}]{Flow of the Hamiltonian used to prove Lemma~\ref{PerturbElem} (``local rotation'').}\label{FigLocRota}
\end{center}
\end{figure}
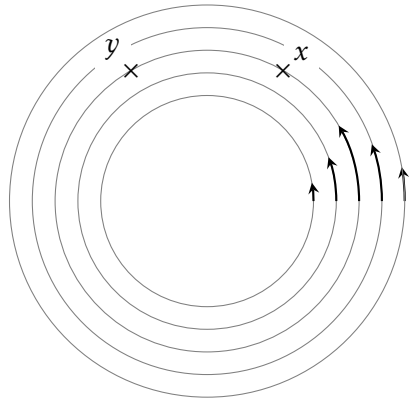

\begin{proof}[Sketch of proof of Lemma~\ref{PerturbElem}]
We sketch the proof in the conservative case; the dissipative case being obtained by the same arguments (and even, by simpler arguments).

The idea of the proof of this lemma is to apply ``local rotations'', as in the proof of the proposition of extension of finite maps in the $C^0$ topology.

More precisely, in dimension 2, $g$ is obtained by composing $f$ by the time 1 of a Hamiltonian flow whose orbits are circles with a common centre, such that one of these circles meets both $x$ and $y$ (see Figure~\ref{FigLocRota}). As the time 1 of this Hamiltonian maps $x$ into $y$, its $C^0$ norm is bigger than $d(x,y)$. But as we also want that the norm of its differential is smaller than $\delta$, its support should contain a ball or radius $(1+1/\delta)d(x,y)$. This explains why when $\delta$ is small, the size of the perturbation has to be big with respect to $d(x,y)$.

A precise  proof can be found for example in \cite[Proposition~5.1.1]{MR1662930}.
\end{proof}

From Lemma~\ref{PerturbElem}, we now deduce another perturbation result which will be at the basis of the rest of the chapter. It states that every periodic orbit of a generic conservative diffeomorphism is shadowed by periodic orbits of its discretizations.

\begin{coro}\label{PropShadowC1}
Let $f \in \Diff^1(X)$ (or $f \in \Diff^1(X,\lambda)$) be a generic diffeomorphism. Then, for every $\tau\in\N$, every $N_0\in\N$ and every $\varep>0$, there exists $N\ge N_0$ such that every periodic orbit of $f$ of period smaller than $\tau$ is $\varep$-shadowed by a periodic orbit of $f_N$ with the same period.
\end{coro}

\begin{proof}[Proof of Corollary~\ref{PropShadowC1}]
We prove the corollary in the dissipative case, the conservative case being identical.

We set $\mathcal S_{\tau,\varep,N}$ the set of $f\in \Diff^1(X)$ such that for every periodic orbit $\omega$ of $f$ of period smaller than $\tau$, there exists a periodic orbit of $f_N$ with the same period which is $\varep$-close to $\omega$. Then, the set of diffeomorphisms satisfying the conclusion of the lemma is the set
\[\bigcap_{\tau,\varep,N_0}\bigcup_{N\ge N_0} \mathcal S_{\tau,\varep,N}.\]
Thus, it suffices to prove that for every $\tau,\varep,N_0$, the set $\bigcup_{N\ge N_0} \mathcal S_{\tau,\varep,N}$ contains an open and dense subset of $\Diff^1(X)$.

Let $\tau>0$, $\varep>0$, $N_0\in\N$, $f\in \Diff^1(X,\lambda)$ and $\delta>0$; we want to find $g\in \Diff^1(X)$ such that $d(f,g)<\delta$ and $g\in \bigcup_{N\ge N_0} \mathcal S_{\tau,\varep,N}$.

We first use a classical result of C.~Robinson, which asserts that that for a generic diffeomorphism $f\in \Diff^1(X)$, the set of periodic points of period smaller than $\tau$ is finite, and moreover is continuous in the diffeomorphism $f$: if $g\in \Diff^1(X)$ is $C^1$-close to $f$, then its set of periodic points of period smaller than $\tau$ is close to that of $f$ for Hausdorff topology (see \cite{MR0273640,MR0279403}, \cite{MR0321141}, see also \cite{MR2288283}). Thus, perturbing a little $f$ if necessary, we can suppose that this property actually holds; taking a smaller $\delta$ if necessary, we can also suppose that there is no creation of periodic orbit when we perturb $f$ into the $\delta$-neighbourhood of $f$.

So we can enumerate $\omega_1,\cdots,\omega_\ell$ the periodic orbits of $f$ of period smaller than $\tau$. For each $i$ we denote $\omega_i = \{x_{i,1},\cdots,x_{i,\tau_i}\}$, with $f(x_{i,j}) = x_{i,j+1}$ ($j$ belonging to $\Z/\tau_i\Z$). Using Lemma~\ref{PerturbElem}, there exists a radius $r_1\in ]0,\varep[$ (which depends on $\delta$, on the $C^1$-norm of $f$ and on the minimal distance between two distinct points of the union of these orbits) such that for every collection of points $y_{i,j}$ satisfying $d(y_{i,j},x_{i,j})<r_1$ for every $i\in\llbracket 1,\ell\rrbracket$ and $j\in\llbracket 1,\tau_i\rrbracket$, there exists a diffeomorphism $g\in\Diff^1(X)$ such that $d_{C^1}(f,g)<\delta$ and that $g(y_{i,j}) = y_{i,j+1}$.

By the hypothesis on the grids $E_N$, there exists $N_1\ge N_0$ such that if $N\ge N_1$, then every ball of radius $r_1$ contains at least one point of $E_N$. We can apply the property stated in the last paragraph to the discretizations $(x_{i,j})_N$ of the points $x_{i,j}$ on the grid $E_N$, and get a diffeomorphism $g$ which belongs to the interior of $\bigcup_{n\ge N_0} \mathcal S_{\tau,\varep,N}$. This completes the proof of the lemma.
\end{proof}

Corollary~\ref{PropShadowC1} states that for a generic conservative diffeomorphism, any periodic orbit is shadowed by periodic orbits of the same period of an infinite number of discretizations. In particular, if we denote by $\operatorname{Ord}(f)$ the set of orders of periodic orbits of $f$, for a generic diffeomorphism $f\in \Diff^1(X,\lambda)$ and for any $M\in\N$, the supremum limit over $N$ of $[0,M] \cap \operatorname{Ord}(f_N)$ is equal $[0,M] \cap \operatorname{Ord}(f)$. Moreover, Corollary~\ref{PropShadowC1} indicates that it is theoretically possible to recover the set of periodic points of a generic conservative $C^1$-diffeomorphism by looking at the corresponding periodic points of its discretizations (see Figure~\ref{FigPtsPer} for numerical simulations corresponding to this property).

\section{Discrete counterparts of transitivity and topological mixing}

As a first application of Corollary~\ref{PropShadowC1}, we apply the connecting lemma of C.~Bonatti and S.~Crovisier \cite{MR2090361} to get a $C^1$ counterpart of Corollary~\ref{typlax}.

\begin{theoreme}[Bonatti-Crovisier]\label{PapySylvain}
A generic diffeomorphism $f\in\Diff^1(X,\lambda)$ is transitive. More precisely, for any generic diffeomorphism $f\in\Diff^1(X,\lambda)$ and any $\varep>0$, there exists a periodic orbit of $f$ which is $\varep$-dense.

In the dissipative case, for a generic diffeomorphism $f\in\Diff^1(X)$, every maximal invariant chain-transitive set is transitive. More precisely, for every maximal invariant chain-transitive set $K$ and every $\varep>0$, there exists a periodic orbit of $f$ which is $\varep$-dense in $K$. 
\end{theoreme}

In the conservative case, this difficult theorem is somehow a weak answer to the question asking whether a generic conservative $C^1$-diffeomorphism is ergodic or not.

Combining this theorem with Corollary~\ref{PropShadowC1}, we obtain directly the following corollary.

\begin{coro}\label{typlaxC1}
Let $f$ be a generic diffeomorphism of $\Diff^1(X,\lambda)$. Then, for any $\varep>0$ and any $N_0>0$, there exists $N\ge N_0$ such that $f_N$ has a periodic orbit which is $\varep$-dense. 

Let $f$ be a generic diffeomorphism of $\Diff^1(X)$. Then, for every maximal invariant chain-transitive set $K$, for every $\varep>0$ and every $N_0>0$, there exists $N\ge N_0$ such that $f_N$ has a periodic orbit which is $\varep$-dense in $K$. 
\end{coro}

In the conservative case, this corollary expresses that for every $\varep>0$, an infinite number of discretizations are ``$\varep$-topologically transitive'', thus contain a sub-dynamics which is similar to that of the initial diffeomorphism.

In the dissipative case, it implies that it is possible to recover the chain transitive set $\mathcal R(f)$ of a generic dissipative diffeomorphism $f$, by using the recurrent sets of the discretizations: $\mathcal R(f) = \overline\lim_{N\to +\infty} \Omega(f_N)$.

\bigskip

More recently, F.~Abdenur and S.~Crovisier have obtained in \cite{MR2975581} a more sophisticated closing lemma, which allows to prove that a generic diffeomorphism is topologically mixing\footnote{Recall that a continuous map $f : X\to X$ is \emph{topologically mixing} if for any non-empty open sets $U,V\subset X$, there exists $M\in\N$ such that $f^m(U)\cap V \neq\emptyset$ for every $m\ge M$.}.

\begin{theoreme}[Closing lemma with time control, Abdenur-Crovisier]\label{FlavSylv}
Let $f\in \Diff^1(X)$, $\ell\ge 2$ be an integer, $x$ be either a non-periodic point or a non-resonant periodic point\footnote{All the periodic points of a generic conservative diffeomorphism are non resonant, see Definition~3.1 of \cite{MR2975581}.}. Assume that each neighbourhood $V$ of $x$ intersects some iterate $f^m(V)$ such that $m$ is not a multiple of $\ell$. Then, there exists diffeomorphisms $g\in \Diff^1(X)$, arbitrarily $C^1$-close to $f$, such that $x$ is periodic under $g$ whose period is not a multiple of~$\ell$.

Moreover, if $f\in \Diff^1(X,\lambda)$, then $g$ can be supposed to be conservative too.
\end{theoreme}

This perturbation lemma can be seen as a weak version of Proposition~\ref{melfaibl}, which implies that under the same conditions, but in the $C^1$ topology, the period of $x$ can be supposed to be coprime with $\ell$. However, combining this lemma with arguments that are specific to diffeomorphisms, we get a stronger result of genericity, that is that a generic diffeomorphism is topologically strongly mixing.

From this theorem, it is possible to deduce the following statement.

\begin{coro}\label{CoroFlavSylvain}
For a generic $f\in \Diff^1(X,\lambda)$ and for any $\varep>0$, $f$ has two $\varep$-dense periodic points, whose periods are not multiples on to the other. The same property holds for a generic dissipative diffeomorphism on every maximal invariant chain-transitive set $K$.
\end{coro}

\begin{proof}[Proof of Corollary~\ref{CoroFlavSylvain}]
Let $\varep>0$; we show that the set of $f\in \Diff^1(X,\lambda)$ satisfying the conclusions of the corollary for $\varep$ is open and dense. By the connecting lemma of C.~Bonatti and S.~Crovisier (Theorem~\ref{PapySylvain}), for a generic $f\in \Diff^1(X,\lambda)$ and for any $\varep>0$, $f$ has a $\varep$-dense periodic point, that we denote by $p$, whose period is $\ell$; we can furthermore assume that the periodic point $p$ is persistent. Then, by Proposition~2.3 of \cite{MR2975581}, the map $f^\ell$ is transitive on the pointwise homoclinic class\footnote{The \emph{pointwise homoclinic class} is the closure of the set of transverse intersection points between the manifolds $W^s(p)$ and $W^u(p)$} of $p$. Thus, for any neighbourhood $V$ of $p$, there exists $k\in\N$ such that $f^{k\ell}\big(f(V)\big)\cap V \neq \emptyset$; this allows us to apply Theorem~\ref{FlavSylv} to perturb the diffeomorphism $f$ to a diffeomorphism $g$. Making both of the obtained periodic points of $g$ persistent (if necessary), a 
whole neighbourhood of $g$ satisfies the conclusions of the corollary.

The dissipative case is identical: it suffices to apply the dissipative version of Theorem~\ref{PapySylvain}.
\end{proof}

From Corollary~\ref{CoroFlavSylvain}, and using Corollary~\ref{PropShadowC1}, we deduce directly a statement about discretizations, which is somehow a $C^1$ counterpart of Corollary~\ref{méldiscr}.

\begin{coro}\label{CoroFlavSylvain2}
For a generic $f\in \Diff^1(X,\lambda)$, for any $\varep>0$ and any $N_0>0$, there exists $N\ge N_0$ such that $f_N$ has two $\varep$-dense periodic points, whose periods are not multiples on to the other. The same property holds for a generic dissipative diffeomorphism on every maximal invariant chain-transitive set $K$.
\end{coro}

Thus, in a certain sense, for every $\varep>0$, an infinite number of discretizations are $\varep$-topologically mixing on every maximal invariant chain-transitive set.

\section{Shadowing of invariant measures}

We now come to the ergodic properties of the discretizations of a generic diffeomorphism. We make use of an ergodic closing lemma, which allows to approximate every invariant measure by periodic invariant measures

\begin{theoreme}[Mañé, Abdenur-Bonatti-Crovisier]\label{Mane}
Let $f\in\Diff^1(X,\lambda)$ be a generic diffeomorphism. Then, every $f$-invariant measure $\mu$ is the weak limit of periodic measures.

Let $f\in\Diff^1(X)$ be a generic diffeomorphism. Then, every $f$-invariant measure $\mu$ which is supported by an invariant chain-transitive set is the weak limit of periodic measures.
\end{theoreme}

The case where the invariant measure is ergodic has been obtained by R.~Mañé in \cite{MR678479}, the general case has been treated by F.~Abdenur, C.~Bonatti and S.~Crovisier in \cite[Theorem 3.5]{MR2811152} (the theorem also holds for generic conservative $C^1$-diffeomorphisms).

Applying Corollary~\ref{PropShadowC1}, we get the following corollary, which is a $C^1$ counterpart of Theorems~\ref{mesinv} and \ref{EnsMesInv}.

\begin{coro}\label{CoroMane}
Let $f$ be a generic conservative diffeomorphism of $\Diff^1(X,\lambda)$. Then, for any $f$-invariant measure $\mu$, any $\varep>0$ and any $N_0>0$, there exists $N\ge N_0$ such that $f_N$ supports a periodic measure which is $\varep$-close to $\mu$.

More generally, for a generic diffeomorphism $f\in\Diff^1(X,\lambda)$, for any $\varep>0$ and $N_0>0$, there exists $N\ge N_0$ such that the set of $f_N$-invariant measures is $\varep$-close to the set of $f$-invariant measures (for the Hausdorff distance on the space of compact sets of probability measures on $X$).

The statement is still true for a generic dissipative diffeomorphism $f\in\Diff^1(X)$, if the measure $\mu$ is supported by an invariant chain-transitive set.
\end{coro}

Roughly speaking, every invariant measure supported by a maximal invariant chain-transitive set is ``seen'' by an infinite number of discretizations.

In Chapter~\ref{chapPhys}, we will obtain an improvement of this corollary (Theorem~\ref{TheoMesPhysDiff}), which will describe the basin of attraction of the periodic measure of the discretization.

\section{Shadowing of chain-transitive invariant sets}

We now treat the case of the shadowing of invariant sets of a generic diffeomorphism. In \cite[Theorem~4]{MR2264835}, S.~Crovisier stated the following result, which asserts that for a generic diffeomorphism, it is possible to recover the chain-transitive invariant sets by looking only at the periodic points.

\begin{theoreme}[Crovisier]\label{SylvIHES}
Any chain-transitive compact invariant set of a generic diffeomorphism $f\in\Diff^1(X)$ of $f\in\Diff^1(X,\lambda)$ is approximated in the Hausdorff topology by periodic orbits.
\end{theoreme}

Obviously, the converse is always true: any point of accumulation of periodic orbits is a chain-transitive compact invariant set.

%

Combined with Corollary~\ref{PropShadowC1}, Theorem~\ref{SylvIHES} leads to the following corollary, which is a $C^1$ counterpart of Theorem~\ref{CompactInvSimpl}, or Theorem~\ref{CompactInv}.

\begin{coro}\label{CoroSylvIHES}
Let $f$ be a generic diffeomorphism of $\Diff^1(X)$ or $\Diff^1(X,\lambda)$. Then, for any $f$-invariant chain-transitive compact set $K$, any $\varep>0$ and any $N_0>0$, there exists $N\ge N_0$ such that $f_N$ has an invariant compact set which is $\varep$-close to $K$ for Hausdorff distance on the set of compact subsets of $X$.

More generally, for a generic diffeomorphism $f\in\Diff^1(X)$ (or $f\in\Diff^1(X,\lambda)$), for any $\varep>0$ and $N_0>0$, there exists $N\ge N_0$ such that the set of $f_N$-invariant sets is $\varep$-close to the set of $f$-invariant chain-transitive compact sets (for Hausdorff distance on the set of  compact sets of compact subsets of $X$).
\end{coro}

\section{Rotation sets}

In this section, we are interested in the approximation of the rotation set of a generic diffeomorphism by the rotation sets of its discretizations (see Chapter~\ref{ChapRot}). Thus, we will take $X = \T^n$ and consider only diffeomorphisms in the connected component of the identity. We also take the same notations as in Chapter~\ref{ChapRot}.

We will see that in the $C^1$-case, applying the elementary perturbation lemma (Corollary~\ref{PropShadowC1}), it is possible to obtain a weaker result about discretized rotation sets (in this case, we can not control what happens on the whole grid $E_N$, but only on a subgrid of $E_N$): for a generic conservative diffeomorphism, the upper limit of the discretized rotation sets is equal to the rotation set of the diffeomorphism itself.

To begin with, we state an approximation lemma we will use in this section; this lemma is a quite direct consequence of Lemma~\ref{PerturbElem}.

\begin{lemme}\label{LemDiscrC1}
If $f$ is generic among $\Diff^1(\T^n)$ or $\Diff^1(\T^n,\Leb)$, then for every finite collection of rotation vectors $\{v_1,\cdots,v_n\}$, each one realized by a periodic orbit of $f$, there exists a subsequence $f_{N_i}$ of discretizations such that for every~$i$, $\{v_1,\cdots,v_n\} \subset \rho_{N_i}(f)$.
\end{lemme}

\begin{proof}[Proof of Lemma \ref{LemDiscrC1}]
The proof of this lemma is very similar to that of Lemma~\ref{EnsRotDiscrCons}; we take the same notations (in particular, $D_q \subset\Q$ is the set of fractions whose numerator is smaller than $q^2$ and whose denominator is smaller than $q$). We prove this lemma in the conservative setting, the dissipative case being identical.

Consider the set
\begin{equation}\label{equation}
\bigcap_{q,N_0}\bigcap_{D\in\mathcal D_q}\bigcup_{N\ge N_0} \left\{\begin{array}{r}
f\in\Diff^1(\T^n,\Leb) \mid(\forall v\in D,\,v\text{ is realised by a}\\
\text{persistent periodic point of $f$}) \implies D\subset \rho(F_N)
\end{array}\right\}.
\end{equation}
To prove the lemma, it suffices to prove that this set contains a $G_\delta$ dense subset of $\Diff^1(\T^n,\Leb)$.

Let $f\in\Hom(\T^n,\Leb)$, $\varep>0$, $q, N_0\in\N$ and $D\in \mathcal D_q$. We suppose that for all $v\in D$, $v$ is realizable by a persistent periodic orbit $\omega_i$ of $f$. Then, by the elementary perturbation lemma (Lemma~\ref{PerturbElem}), it is possible to perturb $f$ into a diffeomorphism $g$ such that $d_{C^1}(f,g)<\varep$ and that there exists $N\ge N_0$ such that for every $i$, there exists a periodic orbit $\omega'_i$ of $g$ which is close to $\omega_i$ (in particular, it has the same rotation vector) and such that $\omega'_i\subset E_N$. Moreover, perturbing a little $g$ if necessary (as in the proof of Proposition \ref{RotGeneCons}), we can suppose that the periodic orbits $\omega'_i$ are persistent. This proves that the set of \eqref{equation} contains a $G_\delta$ dense subset of $\Diff^1(\T^n,\Leb)$.
\end{proof}

To begin with, we treat the conservative case. We will combine Lemma~\ref{LemDiscrC1} with a realization theorem of rotation vectors by periodic orbits; but to use this realization theorem we need to prove that the rotation set of a generic conservative diffeomorphism is nonempty. It is a $C^1$ counterpart of Proposition~\ref{RotGeneCons}.

\begin{prop}\label{RotGeneConsC1}
On a open dense\footnote{In particular, if $f$ is generic.} subset of $\Diff^1(\T^n,\Leb)$, $\rho(F)$ has non-empty interior.
\end{prop}

\begin{rem}
We do not know the shape of the boundary of the rotation set of a generic conservative $C^1$-diffeomorphism. In particular we do not know if it is a polygon or not.
\end{rem}

\begin{proof}[Proof of Proposition \ref{RotGeneConsC1}]
It suffices to resume the proof of Proposition~\ref{RotGeneCons}, and to replace the $C^0$ ergodic closing lemma is replaced by the $C^1$-ergodic closing lemma which is stated in Theorem~\ref{Mane}. To make the obtained periodic point persistent, we just have to apply Franks lemma \cite{MR0283812} to perturb the differential of $g$ on the periodic orbit to avoid having the eigenvalue 1, so that the periodic point becomes persistent (see \cite[page 319]{MR1326374}). The rest of the proof is identical to the $C^0$ case.
\end{proof}

Combined with the realisation theorem of J. Franks (Theorem~\ref{Franks}, see also \cite[Theorem 3.2]{MR958891}) and the approximation lemma (Lemma~\ref{LemDiscrC1}), this proposition directly implies the following result on the discretizations.

\begin{prop}\label{DiscrC1}
If $f$ is generic among $\Diff^1(\T^n,\Leb)$, then there exists a subsequence $f_{N_i}$ of discretizations such that $\rho_{N_i}(F)$ tends to $\rho(F)$ for the Hausdorff topology; in particular, the asymptotic discretized rotation set tends to the rotation set of $f$.
\end{prop}

In particular, when $f$ is a generic conservative $C^1$-diffeomorphism, the asymptotic discretized rotation set
\[\rho^{discr}(F) = \bigcap_{M\in\N}\bigcup_{N\ge M} \rho(F_N)\]
coincides with the rotation set.

For a generic $C^1$-diffeomorphism, we do not know if the observable rotation set is reduced to a singleton or not. This is true when the diffeomorphism is ergodic, which is conjectured to be a generic property (see also page~\pageref{AvilaErgo}).

Remark that the behaviour described by Proposition~\ref{DiscrC1} can actually be observed on actual examples of conservative $C^\infty$-diffeomorphisms, see Section~\ref{SecNumRot}.
\bigskip

The dissipative case is less straightforward, as Proposition~\ref{RotGeneConsC1} is false for generic dissipative $C^1$-diffeomorphisms (there are open sets of diffeomorphisms whose rotation set has empty interior). However, the connected component of $\Id$ in $\Diff^1(\T^2)$, includes three disjoint open subsets, whose union is dense.
\begin{enumerate}
\item The first one is the open subset $\mathcal O_1$ of diffeomorphisms $f$ such that $\rho(F)$ has a nonempty interior (this set is open by the continuity\footnote{With the $C^0$ topology on the domain of $\rho$ and the Hausdorff distance on compact subsets of $\R^2$ on its image.} of $F\mapsto\rho(F)$ on every homeomorphism whose rotation set has nonempty interior, see \cite{MR1101087}). On $\mathcal O_1$, the same arguments than in the conservative case can be applied, thus if $f$ is generic among $\mathcal O_1$, then there exists a subsequence $f_{N_i}$ of discretizations such that $\rho_{N_i}(F)$ tends to $\rho(F)$.
\item The second one is the open subset $\mathcal O_2$ of diffeomorphisms $f$ such that $\rho(F)$ is stably a segment. By the ergodic closing lemma of R.~Mañé (Theorem~\ref{Mane}), for a generic $f\in \mathcal O_2$, every rotation vector of an ergodic measure (in particular, the vertices of the segment) is arbitrarily approximated by a sequence of rotation vectors of (persistent) periodic points. In particular, on an open and dense subset of $\mathcal O_2$, the rotation set is a segment with rational slope, which coincides with the closure of the convex hull of rotation vectors of (persistent) periodic points\footnote{Remark that in general, as the slope of the rotation set is rational, we can not hope to get a (for instance) dense subset of the rotation set on which each vector is realized by a periodic point, see for example \cite[Section~3.3]{Begu-ens}}. Thus, applying Lemma~\ref{LemDiscrC1}, we get that if $f$ is generic among $\mathcal O_2$, then there exists a subsequence $\big(f_{N_i}\big)_i$ of 
discretizations such that $\conv(\rho_{N_i}(F))$ tends to $\rho(F)$. 
\item The third one is the open subset $\mathcal O_3$ of diffeomorphisms $f$ such that $\rho(F)$ is stably a singleton. Then, as on an open dense subset of $\Diff^1(\T^2)$, every diffeomorphism possesses a periodic point, on an open and dense subset of $\mathcal O_3$, the rotation vector is realized by a periodic point. Then, trivially, we get that for $f$ belonging to an open and dense subset of $\mathcal O_3$, the rotation sets of $f_N$ converge globally to the rotation set of $f$.
\end{enumerate}

To summarize, we have proved the following property.

\begin{prop}\label{DiscrC1Dissip}
If $f\in\Diff^1(\T^n)$ is a generic dissipative diffeomorphism, then there exists a subsequence $\big(f_{N_i}\big)_i$ of discretizations such that $\conv\big(\rho_{N_i}(F)\big)$ tends to $\rho(F)$ for the Hausdorff topology; in particular, the convex hull of the asymptotic discretized rotation set tends to the rotation set of $f$. Moreover, if $\rho(F)$ has nonempty interior, or is reduced to a singleton, then there is no need to take convex hulls.
\end{prop}

Remark that even if the techniques of proof are quite different, this behaviour is identical to what happens for generic dissipative homeomorphisms (Proposition~\ref{RotDiscrDissip}).

\chapter{Degree of recurrence of a generic diffeomorphism}\label{ChapDeg}

In this chapter, we begin the study of the the global dynamics of the discretizations of a generic $C^1$-diffeomorphism (both conservative and dissipative) by focusing on the degree of recurrence of the discretizations.

We will consider that the space phase is the torus $\T^n$, the measure is Lebesgue measure and the grids are the uniform grids
\[E_N = \left\{ \left(\frac{i_1}{N},\cdots,\frac{i_n}{N}\right)\in \R^n/\Z^n \middle\vert\ 1\le i_1,\cdots,i_n \le N\right\}.\]
We will see in Section~\ref{AddendSett} that this quite restrictive framework can be generalized to arbitrary manifolds, provided that the discretizations grids behave locally (and almost everywhere) like the canonical grids on the torus.

Let us recall the definition of the degree of recurrence.

\begin{definition}\label{DefDegree}
Let $E$ be a finite set and $\sigma : E\to E$ be a finite map on $E$. The \emph{recurrent set} of $\sigma$ is the union $\Omega(\sigma)$ of the periodic orbits of $\sigma$; it is also equal to the set $\sigma^t(E)$ for every $t$ large enough.

The \emph{degree of recurrence} of the finite map $\sigma$ is the ratio $D(\sigma)$ between the cardinality of the recurrent set and the cardinality of $E$, that is\index{$D(f_N)$}
\[D(\sigma) = \frac{\card\big(\Omega(\sigma)\big)}{\card(E)}.\]
\end{definition}

The goal of this chapter is to study the behaviour of the degree of recurrence $D(f_N)$ as $N$ goes to infinity and for a generic conservative/dissipative diffeomorphism. This degree of recurrence somehow represents the amount of information we lose when we iterate the discretization.

As it can be obtained as the decreasing limit of finite time quantities, the degree of recurrence is maybe the easiest combinatorial invariant to study: we will deduce the behaviour of the degree of recurrence from that of the rate of injectivity.

\begin{definition}\label{DefTauxDiffeo}
Let $f\in \End(\T^n)$ be an endomorphism of the torus and $t\in\N$. The \emph{rate of injectivity} in time $t$ and for the order $N$ is the quantity\index{$\tau_N^t$}
\[\tau_N^t(f) = \frac{\card\big((f_N)^t(E_N)\big)}{\card(E_N)}.\]
Then, the \emph{upper rate of injectivity} of $f$ in time $t$ is defined as\index{$\tau^t$}
\begin{equation}\label{EqTauT}
\tau^t(f) = \underset{N\to +\infty}{\overline\lim} \tau_N^t(f).
\end{equation}
\end{definition}

The link between the degree of recurrence and the rates of injectivity is made by the following formula:
\[D(f_N) = \lim_{t\to+\infty}\tau_N^t(f).\]

The study of the rates of injectivity will be the opportunity to understand the local behaviour of the discretizations of diffeomorphisms: Theorem~\ref{convBis} asserts that that these rates of injectivity are obtained by averaging the corresponding quantities for the differentials of the diffeomorphism. The proof of this result involves the local linearization of a diffeomorphism (Lemma~\ref{LemExtension}), together with estimates of the lack of continuity of the rate of injectivity in the linear case (Proposition~\ref{oscill}). We generalize later this result to the $C^r$ topology, for generic diffeomorphisms (Theorem~\ref{convBisMieux}) and for generic expanding maps (Theorem~\ref{TauxExpand}).

In the conservative case, we will use the study of the rate of injectivity of matrices with determinant 1 in Part~2 of this manuscript (see Theorem~\ref{ConjPrincip}). It will lead to the proof of the fact that the sequence $D(f_N)$ of degrees of recurrence of a generic conservative $C^1$-diffeomorphism tends to 0 (Theorem~\ref{limiteEgalZero}).

For a generic dissipative $C^1$-diffeomorphism, the sequence $D(f_N)$ also converges to 0 (Corollary~\ref{CoroCoroArturJairo}); it is an easy consequence of a theorem of A.~Avila and J.~Bochi (Theorem~\ref{ArturJairo}, see also \cite{MR2267725}).

Note that the fact that the local-global formula is true for $C^r$-generic diffeomorphisms does not help to conclude about the degree of recurrence of such maps: a priori, we need to perturs the derivative of such maps on a large subset of the torus. However, we can hope that Theorem~\ref{limiteEgalZero} remains true for these higher regularities.
\bigskip

We now explain in more detail why the behaviour of $D(f_N)$ can be deduced from that of $\tau^t$. When $N$ is fixed, the sequence $(\tau_N^t(f))_t$ is decreasing in $t$, so $D(f_N) \le \tau_N^t(f)$ for every $t\in\N$. Taking the upper limit in $N$, we get
\[\underset{N\to +\infty}{\overline\lim} D(f_N) \le \tau^t(f)\]
for every $t\in\N$, so
\begin{equation}\label{intervLim}
\underset{N\to +\infty}{\overline\lim} D(f_N) \le \lim_{t\to +\infty} \tau^t(f)
\end{equation}
(as the sequence $(\tau^t(f))_t$ is decreasing, the limit is well defined). In particular, if we have an upper bound on $\lim_{t\to +\infty} \tau^t(f)$, this will give a bound on $\underset{N\to +\infty}{\overline\lim} D(f_N)$. Thus, the proof of Theorem~\ref{limiteEgalZero} is reduced to the study of the quantity $\lim_{t\to +\infty} \tau^t(f)$.

At the end of this chapter, we recall the results of the simulations we have conducted about the degree of recurrence of $C^1$-diffeomorphisms; it shows that in practice, the degree of recurrence tends to 0, at least for the examples of diffeomorphisms we have tested.

\section{Local-global formula}\label{ChapLocGlob}

In this section, we state a local-global formula which links the rate of injectivity of a generic diffeomorphism $f$ and the rates of injectivity of its differentials. The fact that the map $f$ is $C^1$ -- then possesses differentials -- introduces a mesoscopic scale for the study of the action of the discretizations in small time:
\begin{itemize}
\item at the macroscopic scale, the discretization of $f$ acts as $f$;
\item at the intermediate mesoscopic scale, the discretization of $f$ acts as a linear map;
\item at the microscopic scale, we are able to see that the discretization is a finite map and we see that the phase space is discrete.
\end{itemize}

Recall that the discretization of a linear map $A : \R^n \to \R^n$ is the map $\widehat A : \Z^n \to \Z^n$ defined by $\widehat A(x) = \pi(Ax)$, where $\pi : \R^n\to\Z^n$ is a projection on (one of) the nearest integer point for the euclidean distance (see Definition~\ref{DefDiscrLin}).

The rate of injectivity in time $k$ of the matrices $A_1,\cdots,A_k \in GL_n(\R)$ is then defined as (see Definition~\ref{DefTaux}; see also Corollary~\ref{corolimitexist} and Theorem~\ref{imgquasi} for the fact that this limit is well defined)
\[\tau^k(A_1,\cdots,A_k) = \lim_{R\to +\infty} \frac{\card \big((\widehat{A_k}\circ\cdots\circ\widehat{A_1}) [B_R]\big)}{\card [B_R]}\in]0,1];\]
and for an infinite sequence $(A_k)_{k\ge 1}$ of invertible matrices, as the previous quantity is decreasing, we define the asymptotic rate of injectivity
\[\tau^\infty\big((A_k)_{k\ge 1}\big) = \lim_{k\to +\infty}\tau^k(A_1,\cdots,A_k)\in[0,1].\]

So, the link between local and global behaviours of the rates of injectivity is given by the following theorem.

\begin{theoreme}\label{conv}
Let $f\in \Diff^1(\T^n)$ (or $f\in \Diff^1(\T^n,\Leb)$) be a generic diffeomorphism. Then $\tau^1(f)$ is well defined (that is, the superior limit is actually a limit) and satisfies:
\begin{equation}\label{eqInt}\tag{L-G}
\tau^1(f) = \int_{\T^n} \tau^1(D f_x) \ud x.
\end{equation}
Moreover, the function $\tau$ is continuous in $f$.
\end{theoreme}


\begin{ex}\label{aoirfhaoijfaeij}
Theorem~\ref{conv} becomes false if we do not suppose that the diffeomorphism is generic. For example, take a diffeomorphism $f\in\Diff^1(\T^n,\Leb)$, which is equal to
\[f_0 +v = \begin{pmatrix} \frac12 & -1 \\ \frac12 & 1 \end{pmatrix}+v.\]
in a ball of radius $r>0$, with a vector $v$ such that $v = (0.2/N, -0.1/N)$ modulo $\Z^2/N$, and ``generic'' where it does not coincide with this affine map. To construct this example more rigorously, we can for example apply Lemma~\ref{LemExtension} to have a $C^1$-diffeomorphism which is equal equal to
$f_0 +v$ in a ball of radius $r>0$, and apply the process of proof of Theorem~\ref{conv} outside of this ball, such that outside of this ball, a formula similar to Equation~\ref{eqInt} holds.

Then, by Lemma~\ref{CalcEx}, the rate of injectivity of $f_0$ is equal to $1/2$ and that of $f_0 +v$ is equal to $1$. Because $f$ is generic anywhere else, both sides of Equation~\eqref{eqInt} are different. Moreover, the rate of injectivity $\tau^1$ is not continuous in $f$.
\end{ex}

%

To show Theorem \ref{conv}, we have to show that Equation~\eqref{eqInt} holds on a dense $G_\delta$ set of diffeomorphisms. The ``density'' part is easily deduced from the following lemma.

\begin{figure}[t]
\begin{center}
\begin{tikzpicture}[scale=.85]
\draw[thick] (0,0) rectangle (4,4);
\draw (0,0) grid (4,4);
\foreach \k in {1,...,4}
 {\foreach \l in {1,...,4}
  {\fill[color=blue!8!white] (\k-.9,\l-.9) rectangle (\k-.1,\l-.1);
	\draw[color=blue!60!black] (\k-.9,\l-.9) rectangle (\k-.1,\l-.1);
}}
\end{tikzpicture}\caption[The set $\mathcal C$]{The set $\mathcal C$ (in blue) in the torus.}\label{FigAvila}
\end{center}
\end{figure}
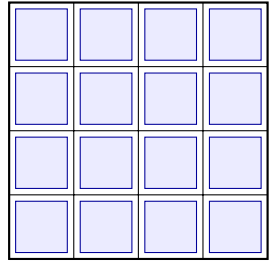

\begin{lemme}\label{morceaux}
Let $f_1\in\Diff^1(\T^n, \Leb)$ (respectively $f_1\in\Diff^1(\T^n)$). Then there exists $f_2\in\Diff^1(\T^n, \Leb)$ (respectively $f_2\in\Diff^1(\T^n)$) arbitrarily close to $f_1$, and a subset $\mathcal C\in\T^n$ (which is a finite union of cubes, see Figure~\ref{FigAvila}) of measure arbitrarily close to $1$, such that the differential of $f_2$ is piecewise constant and totally irrational\footnote{A matrix $A\in GL_n(\R)$ is \emph{totally irrational} if $A\Z^n$ is equidistributed modulo $\Z^n$, see Definition~\ref{DefMeanRate}.} on $\mathcal C$.
\end{lemme}

This lemma follows from the following result of A.~Avila, S.~Crovisier and A.~Wilkinson.

\begin{lemme}[Avila, Crovisier, Wilkinson]\label{LemExtension}
Let $C$ be the unit ball of $\R^n$ for $\|\cdot\|_\infty$ and $\varep>0$. Then, there exists $\delta>0$ such that for every $g_1\in\Diff^\infty(\R^n)$ such that $d_{C^1}({g_1}_{|C}, \Id_{|C})<\delta$, there exists $g_2\in\Diff^\infty(\R^n)$ such that:
\begin{enumerate}[(i)]
\item $d_{C^1}({g_2}_{|C}, {g_1}_{|C})<\varep$;
\item ${g_2}_{|(1-\varep)C} = \Id_{|(1-\varep)C}$;
\item ${g_2}_{|C^\complement} = {g_1}_{|C^\complement}$.
\end{enumerate}
Moreover, if $g_1$ preserves Lebesgue measure, then $g_2$ can be chosen to preserve it as well.
\end{lemme}

The dissipative case of this lemma is easily obtained by interpolating the diffeomorphism with identity, using a smooth bump function. The proof in the conservative case is more difficult and involves a result of J.~Moser \cite{MR0182927}. The reader may refer to \cite[Corollary 6.9]{ArturSylvain} for a complete proof of this lemma\footnote{Le 24/02/2015, cette version n'est pas encore en ligne\dots}.

\begin{proof}[Proof of Lemma \ref{morceaux}]
First of all, we regularize the diffeomorphism $f_1$ to get a $C^\infty$ diffeomorphism $f_3$ which is close to $f_1$ in the $C^1$ topology. In the dissipative case, this is easily obtained (for example) by convolving it by an approximation of the identity. The conservative case is much more difficult; this result has been obtained recently by A.~Avila \cite{MR2736152}.

We then obtain the lemma by applying Lemma~\ref{LemExtension} to the restriction of $f_3$ on each cube of a fine enough cubulation of $\T^n$ (see Figure~\ref{FigAvila}).
\end{proof}

The proof of Theorem \ref{conv} will primarily consist in applying Lemma~\ref{morceaux} together with the following lemma.

\begin{lemme}\label{LemPerLin}
Let $C$ be a cube of $\T^n$. Then, for every totally irrational matrix $A\in GL_n(\R)$, every $v\in \T^n$ and every $\varep>0$, there exists $\delta>0$ such that for every $C^1$ map $f : C \to\T^n$ satisfying $\sup_{x\in C} \|Df_x -A \|\le \delta$, there exists $N_0\in\N$ such that for every $N\ge N_0$,
\[\frac{\card\big( (A+v)_N(E_N\cap C)\, \Delta\, f_N(E_N\cap C) \big)}{\card(E_N\cap C)} \le \varep.\]
\end{lemme}

Before proving Lemma~\ref{LemPerLin}, we explain how it implies Theorem~\ref{conv}.

\begin{proof}[Proof of Theorem \ref{conv}]
We perform the proof in the dissipative setting, the conservative case being identical.
Let $f\in\Diff^1(\T^n)$. The idea is to cut the torus $\T^n$ into small pieces on which $f$ is very close to its Taylor expansion at order $1$.

Let $\mathcal U_\ell$ ($\ell\in\N^*$) be the set of diffeomorphisms $f$ such that the set of accumulation points of the sequence $\big(\tau^1_N(f)\big)_N$ is included in the ball of radius $1/\ell$ and centre $\int_{\T^n} \tau^1(D f_x) \ud x$ (that is, the right side of Equation~\eqref{eqInt}). We want to show that $\mathcal U_\ell$ contains a dense open subset of $\Diff^1(\T^n)$. In other words, we pick a diffeomorphism $f$, an integer $\ell$ and $\delta>0$, and we want to find another diffeomorphism $g$ which is $\delta$-close to $f$, and which belongs to the interior of $\mathcal U_\ell$.

To do that, we first use Lemma~\ref{morceaux}, which gives a diffeomorphism $g$ which is $\delta$-close to $f$ and whose differential is piecewise constant and irrational on a finite union of cubes $\mathcal C$ whose measure is bigger than $1-1/(10\ell)$. Then by construction, $g\in\mathcal U_{2\ell}$. Indeed,
as the differentials of $g$ are irrational on $C$, the rates of injectivity of the translates of the differentials are all equal to that of the differential itself (see Proposition~\ref{ThMeanRate}), thus
\[ \frac{\card\big(g_N(E_N\cap \mathcal C)\big)}{\card(E_N\cap \mathcal C)}\underset{N\to +\infty}{\longrightarrow}\int_{\mathcal C} \tau^1(D f_x) \ud x.\]
Moreover, as the measure of $\mathcal C$ is bigger than $1-1/(10\ell)$, what happens on the complement of $\mathcal C$ can be neglected, more precisely,
\begin{align*}
\left| \tau^1_N(g) - \frac{\card\big(g_N(E_N\cap \mathcal C)\big)}{\card(E_N\cap \mathcal C)}\right| \le & \left| \frac{\card\big(g_N(E_N)\big)}{\card(E_N)} - \frac{\card\big(g_N(E_N)\big)}{\card(E_N\cap \mathcal C)}\right| \\
      & + \left| \frac{\card\big(g_N(E_N)\big)}{\card(E_N\cap \mathcal C)} - \frac{\card\big(g_N(E_N\cap \mathcal C)\big)}{\card(E_N\cap \mathcal C)}\right|\\
	\le & \left| 1 - \frac{\card(E_N)}{\card(E_N\cap \mathcal C)}\right| + \frac{\card(E_N\cap \mathcal C)}{\card(E_N)}\\
	\underset{N\to +\infty}{\longrightarrow} & \, \frac{1}{1-1/(10\ell)}-1+\frac{1}{10\ell}\le \frac{1}{4\ell},
\end{align*}
and as for every $x$, $0\le \tau^1(D f_x) \le 1$, we also have
\[\left|\int_{\T^n} \tau^1(D f_x) \ud x - \int_{\mathcal C} \tau^1(D f_x) \ud x \right| \le \Leb\big(\mathcal C^\complement\big) \le \frac{1}{10\ell}.\]

To show that a whole neighbourhood of $g$ belongs to $\mathcal U_\ell$, it suffices to apply Lemma~\ref{LemPerLin} to each cube of $\mathcal C$; as the measure of $\mathcal C$ is bigger than $1-/(10\ell)$ (again, what happens on $\mathcal C^\complement$ does not count), we get the conclusion of the theorem.
\end{proof}

Thus, it remains to prove Lemma~\ref{LemPerLin}.

\begin{proof}[Proof of Lemma~\ref{LemPerLin}]
Let $A\in GL_n(\R)$ be a totally irrational matrix. Then, there exists a parameter $\delta>0$ and a radius $R_0>0$ such that every matrix $B\in Gl_n(\R)$ such that $\|A-B\|\le\delta$ is ``almost totally irrational'', so that:
\begin{enumerate}
\item[(H1)] the rates of injectivity of the translates of $B$ are close one to each other and can be seen on every ball of radius bigger than $R_0$: there exists $R_0>0$ such that for every $R>R_0$, every $B$ satisfying $\|A-B\|\le\delta$ and every $v\in\R^n$, we have
\[\left|\tau^1(B) - \det(B) D_R^+\big(\pi(B(\Z^n)+v)\big)\right| \le 1/(10\ell);\]
\item[(H2)] the image sets of the differentials are well distributed on every ball of radius bigger than $R_0$: if we note $(\Z^n)'$\index{$(\Z^n)'$} the set of points of $\R^n$ at least one coordinate of which belongs to $\Z+1/2$, then for every $R\ge R_0$ and every $B$ satisfying $\|A-B\|\le\delta$, we have
\[\det(B) D_R^+\big\{ y\in \Z^n \mid d(B(y), (\Z^n)') < 1/(40\ell n) \big\} \le 1/(10\ell).\]
\end{enumerate}
In other words, everything about the rate of injectivity is uniform. The fact that Property~(H1) is true on an open dense set of matrices is obtained by applying Proposition~\ref{ThMeanRate} and Proposition~\ref{oscill}, and that Property (H2) is true on an open dense set is a direct consequence of Lemma~\ref{passifacil}. Remark that the uniformity of the $R_0$ comes from Remark~\ref{Rempassifacil}.
\bigskip

Let $f : C \to\T^n$ be a $C^1$ map satisfying $\sup_{x\in C} \|Df_x -A \|\le \delta$. We write the Taylor expansion of order 1 of $f$ at the neighbourhood of $x\in C$; by compactness we obtain
\[\sup\left\{\frac{1}{\|h\|}\big\|f(x+h)-f(x)-Df_x(h)\big\|\ \middle|\ x\in C,h\in B(0,\varep)\right\}\underset{\varep\to 0}{\longrightarrow}0.\]
Thus, for every $\eta > 0$, there exists $\varep>0$ such that for all $x\in C$ and all $h\in B(0,\varep)$, we have
\begin{equation}\label{eqeta}
\left\|f(x+h)-\big(f(x)+Df_x(h)\big)\right\| < \eta\|h\| \le \eta\varep.
\end{equation}

We now take $R\ge R_0$ (given by hypotheses (H1) and (H2)). We want to find an order of discretization $N$ such that the error made by linearizing $f$ on $B(x,R/N)$ is small compared to $N$, that is, for every $h\in B(0,R/N)$, we have
\[\left\|f(x+h)-\big(f(x)+Df_x(h)\big)\right\| < \frac{1}{40\ell n N}.\]
To do that, we take $\eta = 1/(40 \ell n R)$, then get a $\varep>0$ given by Equation~\eqref{eqeta} (we can take $\varep$ as small as we want), and set $N = \lceil R/\varep \rceil$ (thus, we can take $N$ as big as we want). By Equation~\eqref{eqeta}, for every $y\in B(0,R)$, we obtain
\[\left\|f(x+y/N)-\big(f(x)+Df_x(y/N)\big)\right\| < \frac{1}{40\ell n N}.\]
Combining this result with Hypothesis~(H2), we obtain that the proportion of points of $E_N\cap B(x,R/N)$ whose images by the discretizations of $f$ and that of the linearization of $f$ in $x$ is smaller than $1/(10\ell)$. By Hypothesis~(H1), the density of the discretization of the linearization of $f$ in $x$ is close to $\tau^1(Df_x)$, that is  (recall that $P_N$ is the projection of $X$ on $E_N$)
\[\left|\tau^1(Df_x) - \det(Df_x) D_R^+\big(P_N(Df_x(\Z^n)+v)\big)\right| \le 1/(10\ell).\]

These two facts lead to
\begin{equation}\label{eqDiffTaux}
\left|\frac{\card\Big(P_N\big(f(x+B(0,R/N) \cap E_N)\big)\Big)}{\card\big(B(0,R/N) \cap E_N\big)} - D_{R}^+(Df_x)\right| < \frac{1}{2\ell}.
\end{equation}
\bigskip

We then tile $C$ by smaller cubes of size of the order of $\varep\simeq R/N$. On each of these squares, by Equation~\eqref{eqDiffTaux}, the rate of injectivity of $f$ and that of its differential are $1/2\ell$-close. More precisely, using $\|f\|_{C^1}$, we can find $\delta>0$ arbitrarily small so that the images by $f_N$ of the $\delta$-interiors of the cubes of the tiling are disjoint; so that the rates of injectivity on the cubes add. This proves the lemma.
\end{proof}

With the same kind of proof, we get the same result for arbitrary times $t$.

\begin{theoreme}\label{convBis}
Let $f\in \Diff^1(\T^n)$ (or $f\in \Diff^1(\T^n,\Leb)$) be a generic diffeomorphism. Then $\tau^k(f)$ is well defined (that is, the limit exists) and satisfies:
\[\tau^k(f) = \int_{\T^n} \tau^k\left(Df_x, \cdots, D f_{f^{k-1}(x)}\right) \ud \Leb(x).\]
Moreover, the function $\tau^k$ is continuous in $f$.
\end{theoreme}

\begin{coro}\label{ContTauBarDiff}
The mean rate of injectivity in time $t$ (see also Definition~\ref{DefMeanRate})
\[\overline\tau^k(f) = \int_{\T^n} \overline\tau^k\left(Df_x, \cdots, D f_{f^{k-1}(x)}\right) \ud \Leb(x)\]
is continuous (and even locally Lipschitz) on $\Diff^1(\T^n)$, and coincides with the rate of injectivity when the diffeomorphism $f$ is generic.
\end{coro}

\begin{coro}\label{derdesder?}
As the asymptotic rate of injectivity $\lim_{k\to +\infty} \tau^k(f)$ is upper semi continuous at every generic diffeomorphism $f$, that is, for every $\varep>0$ and every $g$ close enough to $f$, we have
\[\lim_{k\to +\infty}\tau^k(g) \le \lim_{k\to +\infty} \tau^k(f) + \varep.\]
\end{coro}

\section[A local-global formula in $C^r$ topology]{A local-global formula for $C^r$-generic expanding maps and $C^r$-generic diffeomorphisms}\label{SecRateExpand}

The goal of this section is to generalize the results of the previous section : we will obtain results for the $C^r$ regularity, and for both generic expanding maps and diffeomorphisms of the torus $\T^n$. Here, the term expanding map is taken from the point of view of discretizations: we say that a linear map $A$ is expanding if there does not exist two distinct integer points $x,y\in\Z^n$ such that $\widehat A(x) = \widehat A(y)$. The main result of this section is that the rate of injectivity of both generic $C^1$-diffeomorphisms and generic $C^r$-expanding maps of the torus $\T^n$ is obtained from a local-global formula (Theorems~\ref{convBisMieux} and \ref{TauxExpand}). Let us begin by explaining the case of diffeomorphisms.

\begin{theoreme}\label{convBisMieux}
Let $r\ge 1$, and $f\in \Diff^r(\T^n)$ (or $f\in \Diff^r(\T^n,\Leb)$) be a generic diffeomorphism. Then $\tau^k(f)$ is well defined (that is, the limit superior in \eqref{EqTauT} is a limit) and satisfies:
\[\tau^k(f) = \int_{\T^n} \tau^k\left(Df_x, \cdots, D f_{f^{k-1}(x)}\right) \ud \Leb(x).\]
Moreover, the function $\tau^k$ is continuous in $f$.
\end{theoreme}

The idea of the proof of this theorem is very simple: locally, the diffeomorphism is almost equal to a linear map. This introduces an intermediate \emph{mesoscopic} scale on the torus:
\begin{itemize}
\item at the macroscopic scale, the discretization of $f$ acts as $f$;
\item at the intermediate mesoscopic scale, the discretization of $f$ acts as a linear map;
\item at the microscopic scale, we are able to see that the discretization is a finite map and we see that the phase space is discrete.
\end{itemize}
This remark is formalized by Taylor's formula: for every $\varep>0$ and every $x\in\T^n$, there exists $\rho>0$ such that $f$ and its Taylor expansion at order 1 are $\varep$-close on $B(x,\rho)$. We then suppose that the derivative $Df_x$ is ``good'': the rate of injectivity of any of its $C^1$-small perturbations can be seen on a ball $B_R$ of $\R^n$ (with $R$ uniform in $x$). Then, the proof of the local-global formula is made in two steps.
\begin{itemize}
\item Prove that ``a lot'' of maps of $SL_n(\R)$ are ``good''. This is formalized by Lemma~\ref{LemTauxExpand}, which gives estimations of the size of the perturbations of the linear map allowed, and of the size of the ball $B_R$. Its proof is quite technical and uses crucially the formalism of model sets, and an improvement of Weyl's criterion.
\item Prove that for a generic diffeomorphism, the derivative satisfies the conditions of Lemma~\ref{LemTauxExpand} at almost every point. This follows easily from Thom's transversality theorem.
\end{itemize}

As the case of expanding maps is more complicated but similar, we will prove the local-global formula only for expanding maps; the adaptation of it for diffeomorphisms is straightforward.
\bigskip

Remark that the hypothesis of genericity is necessary to get Theorem~\ref{convBisMieux}. For example, it can be seen that if we set 
\[f_0 = \left(\begin{matrix} \frac12 & -1 \\ \frac12 & 1 \end{matrix}\right),\]
then $\tau(f_0)=1/2$ whereas $\tau(f_0+(1/4,3/4)) = 3/4$. Thus, if $g$ is a diffeomorphism of the torus which is equal to $f_0+v$ on an open subset of $\T^2$, with $v$ a suitable translation vector, then the conclusions of Theorem~\ref{convBisMieux} does not hold (see Example~\ref{aoirfhaoijfaeij} for more explanations).
\bigskip

The definition of the linear analogue of the rate of injectivity of an expanding map in time $k$ is more complicated than for diffeomorphisms: in this case, the set of preimages has a structure of $d$-ary tree. We define the rate of injectivity of a tree -- with edges decorated by linear expanding maps -- as the probability of percolation of a random graph associated to this decorated tree (see Definition~\ref{Noel!}). In particular, if all the expanding maps were equal, then the connected component of the root of this random graph is a Galton-Watson tree. We begin by the definition of the set of expanding maps.

\begin{definition}\label{DefExpan}
For $r\ge 1$ and $d\ge 2$, we denote by $\mathcal D^r(\T^n)$\index{$\mathcal D^{r}(\T^n)$} the set of $C^r$ ``$\Z^n$-expanding maps'' of $\T^n$ for the infinite norm. More precisely, $\mathcal D^r(\T^n)$ is the set of maps $f : \T^n\to \T^n$, which are local diffeomorphisms, such that the derivative $f^{(\lfloor r\rfloor)}$ is well defined and belongs to $C^{r-\lfloor r\rfloor}(\T^n)$ and such that for every $x\in \T^n$ and every $v\in\Z^n\setminus\{0\}$, we have $\|Df_x v\|_\infty \ge 1$.

In particular, for $f\in \mathcal D^r(\T^n)$, the number of preimages of any point of $\T^n$ is equal to a constant, that we denote by $d$.
\end{definition}

Remark that in dimension $n=1$, the set $\mathcal D^r(\Sp^1)$ coincides with the classical set of expanding maps: $f\in\mathcal D^r(\Sp^1)$ if and only if it belongs to $C^r(\Sp^1)$ and $f'(x)\ge 1$ for every $x\in\Sp^1$.
\bigskip

We now define the linear setting corresponding to a map $f\in\mathcal D(\T^n)$. 

\begin{definition}
We set (see also Figure~\ref{TreeBee})\index{$I_k$}
\[I_k = \bigsqcup_{m=1}^k \llbracket 1,d\rrbracket^m\]
the set of $m$-tuples of integers of $\llbracket 1,d\rrbracket$, for $1\le m\le k$.

For $\ind = (i_1,\cdots,i_m)\in \llbracket 1,d\rrbracket^m$, we set $\len(\ind) = m$\index{$\len(\ind)$} and $\fat(\ind) = (i_1,\cdots,i_{m-1})\in \llbracket 1,d\rrbracket^{m-1}$\index{$\fat(\ind)$} (with the convention $\fat(i_1) = \emptyset$).
\end{definition}

\begin{figure}[!b]
\begin{minipage}[c]{.35\linewidth}
\begin{center}
\begin{tikzpicture}[scale=.7]
\node (O) at (1,0){$\emptyset$};
\node (A) at (3,1){$(1)$};
\node (B) at (3,-1){$(2)$};
\node (C) at (6,1.5){$(1,1)$};
\node (D) at (6,.5){$(1,2)$};
\node (E) at (6,-.5){$(2,1)$};
\node (F) at (6,-1.5){$(2,2)$};
\draw (O) -- (A);
\draw (A) -- (C);
\draw (A) -- (D);
\draw (O) -- (B);
\draw (B) -- (E);
\draw (B) -- (F);
\end{tikzpicture}
\caption[The tree $T_2$ for $d=2$]{The tree $T_2$ for $d=2$.}\label{TreeBee}
\end{center}
\end{minipage}\hfill
\begin{minipage}[c]{.6\linewidth}
\begin{center}
\begin{tikzpicture}[scale=1.2]
\node (O) at (0,0){$y$};
\node (A) at (3,1){$x_{(1)}$};
\node (B) at (3,-1){$x_{(2)}$};
\node (C) at (6,1.5){$x_{(1,1)}$};
\node (D) at (6,.5){$x_{(1,2)}$};
\node (E) at (6,-.5){$x_{(2,1)}$};
\node (F) at (6,-1.5){$x_{(2,2)}$};
\draw (O) -- (A) node[sloped, midway, above]{$\det Df_{x_{(1)}}^{-1}$};
\draw (A) -- (C) node[sloped, midway, above]{$\det Df_{x_{(1,1)}}^{-1}$};
\draw (A) -- (D) node[sloped, midway, below]{$\det Df_{x_{(1,2)}}^{-1}$};
\draw (O) -- (B) node[sloped, midway, below]{$\det Df_{x_{(2)}}^{-1}$};
\draw (B) -- (E) node[sloped, midway, above]{$\det Df_{x_{(2,1)}}^{-1}$};
\draw (B) -- (F) node[sloped, midway, below]{$\det Df_{x_{(2,2)}}^{-1}$};
\end{tikzpicture}
\caption[Probability tree associated to the preimages of $y$]{The probability tree associated to the preimages of $y$, for $k=2$ and $d=2$. We have $f(x_{(1,1)}) = f(x_{(1,2)}) = x_{(1)}$, etc.}\label{ProbTree}
\end{center}
\end{minipage}
\end{figure}

The set $I_k$ is the linear counterpart of the set $\bigsqcup_{m=1}^k f^{-m}(y)$. Its cardinal is equal to $d(1-d^k)/(1-d)$.

\begin{definition}\label{Noel!}
Let $k\in\N$. The \emph{complete tree of order $k$} is the rooted $d$-ary tree $T_k$\index{$T_k$} whose vertices are the elements of $I_k$ together with the root, denoted by $\emptyset$, and whose edges are of the form $(\fat(\ind),\ind)_{\ind\in I_k}$ (see Figure~\ref{TreeBee}).

Let $(p_\ind)_{\ind\in I_k}$ be a family of numbers belonging to $[0,1]$. These probabilities will be seen as decorations of the edges of the tree $T_k$. We will call \emph{random graph associated to $(p_\ind)_{\ind\in I_k}$} the random subgraph $G_{(p_\ind)_\ind}$\index{$G_{(p_\ind)}$} of $T_k$, such that the laws of appearance the edges $(\fat(\ind),\ind)$ of $G_{(p_\ind)_\ind}$ are independent Bernoulli laws of parameter $p_\ind$. In other words, $G_{(p_\ind)_\ind}$ is obtained from $T_k$ by erasing independently each vertex of $T_k$ with probability $1-p_\ind$.

We define the \emph{mean density} $\overline D((p_\ind)_\ind)$\index{$\overline D((p_\ind)_\ind)$} of $(p_\ind)_{\ind\in I_k}$ as the probability that in $G_{(p_\ind)_\ind}$, there is at least one path linking the root to a leaf.
\end{definition}

Remark that if the probabilities $p_\ind$ are constant equal to $p$, the random graph $G_{(p_\ind)_\ind}$ is a Galton-Watson tree, where the probability for a vertex to have $i$ children is equal to $\binom{d}{i} p^i(1-p)^{d-i}$.

\begin{definition}\label{DefTreeMap}
By the notation $\overline D( (\det Df_x^{-1})_{x\in f^{-m}(y),\, 1\le m\le k})$, we will mean that the mean density is taken with respect to the random graph $G_{f,y}$ associated to the decorated tree whose vertices are the $f^{-m}(y)$ for $0\le m\le k$, and  whose edges are of the form $(f(x),x)$ for $x\in f^{-m}(y)$ with $1\le m \le k$, each one being decorated by the number $\det Df_x^{-1}$ (see Figure~\ref{ProbTree}).
\end{definition}

Recall that the rates of injectivity are defined by (see also Definition~\ref{DefTauxDiffeo})
\[\tau^k(f_N) = \frac{\card\big((f_N)^k(E_N)\big)}{\card(E_N)} \qquad \text{and} \qquad \tau^k(f) = \limsup_{N\to+\infty} \tau^k(f_N).\]

\begin{theoreme}\label{TauxExpand}
Let $r\ge 1$, $f$ a generic element of $\mathcal D^r(\T^n)$ and $k\in\N$. Then, $\tau^k(f)$ is a limit (that is, the sequence $(\tau^k(f_N))_N$ converges), and we have
\begin{equation}\label{EqIntMieux}
\tau^k(f) = \int_{\T^n} \overline D\big( (\det Df_x^{-1})_{\begin{subarray}{l} 1\le m\le k \\ x\in f^{-m}(y)\end{subarray}}\big) \ud \Leb(y).
\end{equation}
Moreover, the map $f\mapsto \tau^k(f)$ is continuous in $f$.
\end{theoreme}

The proof of Theorem~\ref{TauxExpand} is mainly based on the following lemma, which treats the linear corresponding case. Its statement is divided into two parts, the second one being a quantitative version of the first.

\begin{lemme}\label{LemTauxExpand}
Let $k\in\N$, and a family $(A_\ind)_{\ind\in I_k}$ of invertible matrices, such that for any $\ind\in I_k$ and any $v\in\Z^n\setminus\{0\}$, we have $\|A_\ind v\|_\infty \ge 1$.

If the image of the map
\[\Z^n\ni x \mapsto \bigoplus_{\ind\in I_k} A_\ind^{-1} A_{\fat(\ind)}^{-1} \cdots A_{\fat^{\len(\ind)}(\ind)}^{-1} x\]
projects on a dense subset of the torus $\R^{n\card I_k}/\Z^{n\card I_k}$, then we have
\[D\left(\bigcup_{\ind\in \llbracket 1,d\rrbracket^k}\big( \widehat A_{\fat^{k-1}(\ind)} \circ \cdots \circ \widehat A_\ind \big) (\Z^n)\right) = \overline D\big((\det A_\ind^{-1})_\ind\big).\]

More precisely, for every $\ell', c\in\N$, there exists a locally finite union of positive codimension submanifolds $V_q$ of $(GL_n(\R))^{\card I_k}$ such that for every $\eta'>0$, there exists a radius $R_0>0$ such that if $(A_\ind)_{\ind\in I_k}$ satisfies $d((A_\ind)_\ind, V_q)>\eta'$ for every $q$, then for every $R\ge R_0$, and every family $(v_\ind)_{\ind\in I_k}$ of vectors of $\R^n$, we have\footnote{The map $\pi(A+v)$\index{$\pi(A+v)$} is the discretization of the affine map $A+v$.}
\begin{equation}\label{EqLem22}
\left|D_R^+\left(\bigcup_{\ind\in \llbracket 1,d\rrbracket^k}\big( \pi(A_{\fat^{k-1}(\ind)} + v_{\fat^{k-1}(\ind)}) \circ \cdots \circ \pi(A_\ind + v_\ind) \big) (\Z^n) \right) - \overline D\big((\det A_\ind^{-1})_\ind\big)\right| < \frac{1}{\ell'}
\end{equation}
(the density of the image set is ``almost invariant'' under perturbations by translations), and for every $m\le k$ and every $\ind\in \llbracket 1,d\rrbracket^k$, we have\footnote{Where $(\Z^n)'$ stands for the set of points of $\R^n$ at least one coordinate of which belongs to $\Z+1/2$.}
\begin{equation}\label{EqLem222}
D_R^+\left\{x\in \big(A_{\fat^m(\ind)} + v_{\fat^m(\ind)}\big)(\Z^n)\ \middle\vert\ d\big(x,(\Z^n)' \big) < \frac{1}{{c}\ell'(2n+1)} \right\} < \frac{1}{{c}\ell'}
\end{equation}
(there is only a small proportion of the points of the image sets which are obtained by discretizing points close to $(\Z^n)'$).
\end{lemme}

The local-global formula~\eqref{EqIntMieux} will later follow from this lemma, an appropriate application of Taylor's theorem and Thom's transversality theorem (Lemma~\ref{PerturbCr}).

The next lemma uses the strategy of proof of Weyl's criterion to get a uniform convergence in Birkhoff's theorem for rotations of the torus $\T^n$ whose rotation vectors are outside of a neighbourhood of a finite union of hyperplanes.

\begin{lemme}[Weyl]\label{Weyl}
Let $\dist$ be a distance generating the weak-* topology on $\Prb$ the space of Borel probability measures on $\T^n$. Then, for every $\varep>0$, there exists a locally finite family of affine hyperplanes $H_i \subset \R^n$, such that for every $\eta>0$, there exists $M_0\in\N$, such that for every $\boldsymbol\lambda \in \R^n$ satisfying $d(\boldsymbol\lambda,H_q)>\eta$ for every $q$, and for every $M\ge M_0$, we have
\[\dist \left(\frac{1}{M}\sum_{m=0}^{M-1} \bar\delta_{m\boldsymbol\lambda }\,,\ \Leb_{\R^n/\Z^n}\right) < \varep,\]
where $\bar\delta_x$ is the Dirac measure of the projection of $x$ on $\R^n/\Z^n$.
\end{lemme}

\begin{proof}[Proof of Lemma~\ref{Weyl}]
As $\dist$ generates the weak-* topology on $\Prb$, it can be replaced by any other distance also generating the weak-* topology on $\Prb$. So we consider the distance $\dist_W$\index{$\dist_W$} defined by:
\[\dist_W(\mu,\nu) = \sum_{\boldsymbol k\in\N^n} \frac{1}{2^{k_1+\cdots+k_n}}\left| \int_{\R^n/\Z^n} e^{i2\pi\boldsymbol k\cdot x} \ud(\mu-\nu)(x)\right|;\]
there exists $K>0$ and $\varep'>0$ such that if a measure $\mu\in\Prb$ satisfies
\begin{equation}\label{EqWeyl}
\forall \boldsymbol k\in\N^n :\, 0<k_1+\cdots+k_n \le K, \quad \left| \int_{\R^n/\Z^n} e^{i2\pi\boldsymbol k\cdot x} \ud\mu(x)\right|<\varep',
\end{equation}
then $\dist(\mu,\Leb)<\varep$.

For every $\boldsymbol k\in\N^n\setminus\{0\}$ and $j\in\Z$, we set
\[H_{\boldsymbol k}^j = \{\boldsymbol\lambda\in\R^n \mid \boldsymbol k\cdot \boldsymbol\lambda = j\}.\]
Remark that the family $\{H_{\boldsymbol k}^j\}$, with $j\in\Z$ and $\boldsymbol k$ such that $0<k_1+\cdots+k_n \le K$, is locally finite. We denote by $\{H_q\}_q$ this family, and choose $\boldsymbol\lambda \in \R^n$ such that $d(\boldsymbol \lambda,H_q) > \eta$ for every $q$. We also take
\begin{equation}\label{DefMoe}
M_0 \ge \frac{2}{\varep' \left|1-e^{i2\pi\eta}\right|}.
\end{equation}
Thus, for every $\boldsymbol k\in\N^n$ such that $k_1+\cdots+k_n \le K$, and every $M\ge M_0$, the measure
\[\mu = \frac{1}{M}\sum_{m=0}^{M-1} \bar\delta_{m\boldsymbol\lambda}.\]
satisfies
\[\left| \int_{\R^n/\Z^n} e^{i2\pi\boldsymbol k x} \ud\mu(x)\right| = \frac{1}{M} \left|\frac{1-e^{i2\pi M\boldsymbol k\cdot \boldsymbol\lambda}}{1-e^{i2\pi \boldsymbol k\cdot \boldsymbol\lambda}}\right| \le \frac{2}{M_0} \frac{1}{\left|1-e^{i2\pi \boldsymbol k\cdot \boldsymbol\lambda}\right|}.\]
By \eqref{DefMoe} and the fact that $d(\boldsymbol k\cdot \boldsymbol\lambda,\Z)\ge \eta$, we deduce that
\[\left| \int_{\R^n/\Z^n} e^{i2\pi\boldsymbol k x} \ud\mu(x)\right| \le\varep'.\]
Thus, the measure $\mu$ satisfies the criterion \eqref{EqWeyl}, which proves the lemma.
\end{proof}

\begin{proof}[Proof of Lemma~\ref{LemTauxExpand}]
To begin with, let us treat the case $d=1$. Let $A_1,\cdots,A_k$ be $k$ invertible matrices. We want to compute the rate of injectivity of $\widehat A_k \circ \cdots \circ \widehat A_1$. Recall that we set
\[\widetilde M_{A_1,\cdots,A_k} = \left(\begin{array}{ccccc}
A_1 & -\Id &        &         & \\
    & A_2  & -\Id   &         & \\
    &      & \ddots & \ddots  & \\
    &      &        & A_{k-1} & -\Id\\
    &     &        &          & A_k
\end{array}\right)\in M_{nk}(\R),\]
$\widetilde \Lambda_k = \widetilde M_{\lambda_1,\cdots,\lambda_k} \Z^{nk}$ and $W^k = ]-1/2,1/2]^{nk}$. Resuming the proof of Proposition~\ref{CalculTauxModel}, we see that $x\in \big( \widehat A_k \circ \cdots \circ \widehat A_1 \big) (\Z^n)$ if and only if $(0^{n(k-1)},x) \in W^k + \widetilde \Lambda_k$. This implies the following statement.

\begin{lemme}\label{ppe}
We have
\begin{equation}\label{densitete}
\det(A_k\cdots A_1)D\big( \widehat A_k \circ \cdots \circ \widehat A_1 \big) (\Z^n) = \nu(\operatorname{pr}_{\R^{nk}/\widetilde \Lambda_k}(W^k)),
\end{equation}
where $\nu$ is the uniform measure on the submodule $\operatorname{pr}_{\R^{nk}/\widetilde \Lambda_k}(0^{n(k-1)},\Z^n)$ of $\R^{nk}/\widetilde \Lambda_k$.
\end{lemme}

In particular, if the image of the map
\[\Z^n\ni x \mapsto \bigoplus_{m=1}^k (A_m)^{-1}\cdots(A_k)^{-1} x\]
projects on a dense subset of the torus $\R^{nk}/\Z^{nk}$, then the quantity~\eqref{densitete} is equal to the volume of the intersection between the projection of $W^k$ on $\R^{nk}/\widetilde\Lambda_k$ and a fundamental domain of $\widetilde\Lambda_k$ (see the end of the proof of Proposition~\ref{CalculTauxModel} and in particular the form of the matrix $\widetilde M_{A_1,\cdots,A_k}^{-1}$). By the hypothesis made on the matrices $A_m$ -- that is, for any $v\in\Z^n\setminus\{0\}$, $\|A_m v\|_\infty \ge 1$ -- this volume is equal to $1$ (simply because the restriction to $W^k$ of the projection $\R^{nk}\mapsto \R^{nk}/\widetilde\Lambda_k$ is injective). Thus, the density of the set $\big( \widehat A_k \circ \cdots \circ \widehat A_1 \big) (\Z^n)$ is equal to $1/(\det(A_k \cdots A_1))$  .

\begin{figure}[t]
\begin{center}
\begin{tikzpicture}[scale=1.26]
\node (O) at (0,0){$\displaystyle D\left(\bigcup_{\ind\in \llbracket 1,d\rrbracket^k}\!\!\!\! \big( \widehat A_{\fat^{k-1}(\ind)} \circ \cdots \circ \widehat A_\ind \big) (\Z^n)\right)$};
\node (A) at (6.5,1.5){$\displaystyle D\left(\bigcup_{\substack{\ind\in \llbracket 1,d\rrbracket^k \\ \fat^{k-1}(\ind) = 1}}\!\!\!\!\!\!\!\!\!\! \big( \widehat A_{\fat^{k-2}(\ind)} \circ \cdots \circ \widehat A_\ind \big) (\Z^n)\right)$};
\node (B) at (6.5,-1.5){$\displaystyle D\left(\bigcup_{\substack{\ind\in \llbracket 1,d\rrbracket^k \\ \fat^{k-1}(\ind) = 2}}\!\!\!\!\!\!\!\!\!\! \big( \widehat A_{\fat^{k-2}(\ind)} \circ \cdots \circ \widehat A_\ind \big) (\Z^n)\right)$};
\draw (O.5) -- (A.185) node[sloped, midway, above]{$\det A_{(1)}^{-1}$};
\draw (O.-5) -- (B.175) node[sloped, midway, below]{$\det A_{(2)}^{-1}$};
\end{tikzpicture}
\end{center}
\caption[Calculus of the density at the level $k$]{Calculus of the density of the image set at the level $k$ according to the density of its sons.}\label{ProbTree2}
\end{figure}

We now consider the general case where $d$ is arbitrary. We take a family $(A_\ind)_{\ind\in I_k}$ of invertible matrices, such that for any $\ind\in I_k$ and any $v\in\Z^n\setminus\{0\}$, we have $\|A_\ind v\|_\infty \ge 1$. A point $x\in\Z^n$ belongs to
\[\bigcup_{\ind\in \llbracket 1,d\rrbracket^k}\big( \widehat A_{\fat^{k-1}(\ind)} \circ \cdots \circ \widehat A_\ind \big) (\Z^n)\]
if and only if there exists $\ind\in \llbracket 1,d\rrbracket^k$ such that $(0^{m-1},x) \in W^k + \widetilde \Lambda_\ind$. Equivalently, a point $x\in\Z^n$ does not belong to the set
\[\big(\widehat A_{\fat^{k-1}(\ind)} \circ \cdots \circ \widehat A_\ind\big)(\Z^n) \]
if and only if for every $\ind\in \llbracket 1,d\rrbracket^k$, we have $(0^{n(k-1)},x) \notin W^k + \widetilde \Lambda_\ind$. Thus, if the image of the map
\[\Z^n\ni x \mapsto \bigoplus_{\ind\in I_k} A_\ind^{-1} A_{\fat(\ind)}^{-1} \cdots A_{\fat^{\len(\ind)}(\ind)}^{-1} x\]
projects on a dense subset of the torus $\R^{n\card I_k}/\Z^{n\card I_k}$, then the events $x\in S_i$, with
\[S_i = \bigcup_{\substack{\ind\in \llbracket 1,d\rrbracket^k \\ \fat^{k-1}(\ind) = i}}\big( \widehat A_{\fat^{k-1}(\ind)} \circ \cdots \circ \widehat A_\ind \big) (\Z^n)\]
are independent (see Figure~\ref{ProbTree2}), meaning that for every $F\subset \llbracket 1,d\rrbracket$, we have
\begin{equation}\label{independant}
D\left(\bigcap_{i\in F} S_i\right) = \prod_{i\in F} D\big(S_i \big).
\end{equation}
Thus, by the inclusion-exclusion principle, we get
\[D\left(\bigcup_{i\in \llbracket 1,d\rrbracket} S_i\right) = \sum_{\emptyset \neq F \subset \llbracket 1,d\rrbracket} (-1)^{\card(F)+1} \prod_{i\in F} D\big(S_i \big).\]
Moreover, the fact that for any $\ind\in I_k$ and any $v\in\Z^n\setminus\{0\}$, we have $\|A_\ind v\|_\infty \ge 1$ leads to
\[D(S_i) = \det A_{\fat^{k-1}(\ind)}^{-1}\ D\left(\bigcup_{\substack{\ind\in \llbracket 1,d\rrbracket^k \\ \fat^{k-1}(\ind) = i}}\big( \widehat A_{\fat^{k-2}(\ind)} \circ \cdots \circ \widehat A_\ind \big) (\Z^n)\right).\]

These facts imply that the density we look for follows the same recurrence relation as $\overline D\big((\det A_\ind^{-1})_\ind\big)$, thus
\[D\left(\bigcup_{\ind\in \llbracket 1,d\rrbracket^k}\big( \widehat A_{\fat^{k-1}(\ind)} \circ \cdots \circ \widehat A_\ind \big) (\Z^n)\right) = \overline D\big((\det A_\ind^{-1})_\ind\big).\]
\bigskip

The second part of the lemma is an effective improvement of the first one. To obtain the bound~\eqref{EqLem22}, we combine Lemma~\ref{Weyl} with Lemma~\ref{ppe} to get that for every $\varep>0$, there exists a locally finite collection of submanifolds $V_q$ of $(GL_n(\R))^{\card I_k}$ with positive codimension, such that for every $\eta'>0$, there exists $R_0>0$ such that if $d((A_\ind)_\ind,V_q)>\eta'$ for every $q$, then Equation~\eqref{independant} is true up to $\varep$.

The other bound~\eqref{EqLem222} is obtained independently from the rest of the proof by a direct application of Lemma~\ref{Weyl} and of Lemma~\ref{ppe} applied to $k=1$.
\end{proof}

\begin{lemme}[Perturbations in $C^r$ topology]\label{PerturbCr}
Let $1\le r \le +\infty$ and $f$ a generic element of $\mathcal D^r(\T^n)$. Then, for every $k\in\N$, every $\ell'\in\N$ and every finite collection $(V_q)$ of submanifolds of positive codimension of $(GL_n(\R))^{dm}$, there exists $\eta>0$ such that the set
\[T_\eta = \left\{ y\in\T^n\ \middle\vert\ \forall q,\, d\Big( \big( Df_{x}\big)_{\begin{subarray}{l} 1\le m \le k \\ x\in f^{-m}(y)\end{subarray}}\, ,\, V_q\Big) > \eta \right\}\]
contains a finite disjoint union of cubes\footnote{here, a cube is just any ball for the infinite norm.}, whose union has measure bigger than $1-1/\ell'$.
\end{lemme}

\begin{proof}[Proof of Lemma~\ref{PerturbCr}]
By Thom's transversality theorem, for a generic map $f\in\mathcal D^r(\T^n)$, the set
\[\left\{ y\in\T^n\ \middle\vert\ \forall q,\, \big(Df_{x}\big)_{\begin{subarray}{l} 1\le m \le k \\ x\in f^{-m}(y)\end{subarray}}\in V_q \right\}\]
if finite. Thus, the sets $T_\eta^{\scriptscriptstyle\complement}$ are compact sets and their (decreasing) intersection over $\eta$ is a finite set. So, there exists $\eta>0$ such that $T_\eta^{\scriptscriptstyle\complement}$ is close enough to this finite set for Hausdorff topology to have the conclusions of the lemma.
\end{proof}

We can now begin the proof of Theorem~\ref{TauxExpand}.

\begin{proof}[Proof of Theorem~\ref{TauxExpand}]
Let $f\in\mathcal D^r(\T^n)$. The idea is to cut the torus $\T^n$ into small pieces on which $f$ is very close to its Taylor expansion at order $1$.

Let $m\in\N$, and $\mathcal U_\ell$ ($\ell\in\N^*$) be the set of maps $f\in\mathcal D^r(\T^n)$ such that the set of accumulation points of the sequence $(\tau^m_N(f))_N$ is included in the ball of radius $1/\ell$ and centre
\[\int_{\T^n} \overline D\big( (Df_x)_{\begin{subarray}{l} 1\le m\le k \\ x\in f^{-m}(y)\end{subarray}}\big) \ud \Leb(y).\]
(that is, the right side of Equation~\eqref{EqIntMieux}). We want to show that $\mathcal U_\ell$ contains an open and dense subset of $\mathcal D^r(\T^n)$. In other words, we pick a map $f$, an integer $\ell$ and $\delta>0$, and we want to find another map $g\in \mathcal D^r(\T^n)$ which is $\delta$-close to $f$ for the $C^r$ distance, and which belongs to the interior of $\mathcal U_\ell$.

To do that, we first set $\ell' = 3\ell$ and $c=d(1-d^k)/(1-d) = \card(I_k)$, and use Lemma~\ref{LemTauxExpand} to get a locally finite union of positive codimension submanifolds $V_q$ of $(GL_n(\R))^{\card(I_k)}$. We then apply Lemma~\ref{PerturbCr} to these submanifolds, to the $\delta$ we have fixed at the beginning of the proof and to $\ell' = 4\ell$; this gives us a parameter $\eta>0$ and a map $g\in \mathcal D^r(\T^n)$ such that $d_{C^r}(f,g)<\delta$, and such that the set
\[\left\{ y\in\T^n\ \middle\vert\ \forall q,\, d\Big( \big( Dg_{x}\big)_{\begin{subarray}{l} 1\le m\le k \\ x\in g^{-m}(y)\end{subarray}}\, ,\, V_q\Big)<\eta \right\}\]
is contained in a disjoint finite union $\mathcal C$ of cubes, whose union has measure bigger than $1-1/(4\ell)$. Finally, we apply Lemma~\ref{LemTauxExpand} to $\eta' = \eta/2$; this gives us a radius $R_0>0$ such that if $(A_\ind)_{\ind \in I_k}$ is a family of matrices of $GL_n(\R)$ satisfying $d((A_\ind)_\ind, V_q)>\eta/2$ for every $q$, then for every $R\ge R_0$, and every family $(v_\ind)_{\ind\in I_k}$ of vectors of $\R^n$, we have
\begin{equation}\label{DiffDensGene}
\left|D_R^+\left(\bigcup_{\ind\in \llbracket 1,d\rrbracket^k}\big( \pi(A_{\fat^{k-1}(\ind)} + v_{\fat^{k-1}(\ind)}) \circ \cdots \circ \pi(A_\ind + v_\ind) \big) (\Z^n) \right) - \overline D\big((\det A_\ind^{-1})_\ind\big)\right| < \frac{1}{3\ell},
\end{equation}
and for every $i,j$,
\begin{equation}\label{ModuloGene}
D_R^+\left\{x\in \big(A_{j^m(\ind)} + v_{j^m(\ind)}\big)(\Z^n)\ \middle\vert\ d\big(x,(\Z^n)' \big) < \frac{1}{3\ell(2n+1)\card I_k} \right\} < \frac{1}{3\ell\card I_k}.
\end{equation}
\bigskip

We now take a map $h\in \mathcal D^r(\T^n)$ such that $d_{C^1}(g,h)<\delta'$, and prove that if $\delta'$ is small enough, then $h$ belongs to the interior of $\mathcal U_\ell$. First of all, we remark that if $\delta'$ is small enough, then the set
\[\left\{ y\in\T^n\ \middle\vert\ \forall q,\, d\Big( \big( Dh_{x}\big)_{\begin{subarray}{l} 1\le m\le k \\ x\in h^{-m}(y)\end{subarray}}\, ,\, V_q\Big)>\eta/2 \right\}\]
contains a set $\mathcal C'$, which is a finite union of cubes whose union has measure bigger than $1-1/(3\ell)$.

Let $C$ be a cube of $\mathcal C'$, $y\in C$ and $x\in f^{-m}(y)$, with $1\le m \le k$. We write the Taylor expansion of order 1 of $h$ at the neighbourhood of $x$; by compactness we obtain
\[\sup\left\{\frac{1}{\|z\|}\big\|h(x+z)-h(x)-Dh_x(z)\big\|\ \middle|\ x\in C,\,z\in B(0,\rho)\right\}\underset{\rho\to 0}{\longrightarrow}0.\]
Thus, for every $\varep > 0$, there exists $\rho>0$ such that for all $x\in C$ and all $z\in B(0,\rho)$, we have
\begin{equation}\label{eqetata}
\left\|h(x+z)-\big(h(x)+Dh_x(z)\big)\right\| < \varep\|z\| \le \varep \rho.
\end{equation}

We now take $R\ge R_0$. We want to find an order of discretization $N$ such that the error made by linearizing $h$ on $B(x,R/N)$ is small compared to $N$, that is, for every $z\in B(0,R/N)$, we have
\[\left\|h(x+z)-\big(h(x)+Dh_x(z)\big)\right\| < \frac{1}{3\ell(2n+1)\card I_k}\cdot\frac{1}{N}.\]
To do that, we apply Equation~\eqref{eqetata} to
\[\varep = \frac{1}{3R\ell(2n+1)\card I_k},\]
to get a radius $\rho>0$ (we can take $\rho$ as small as we want), and we set $N = \lceil R/\rho \rceil$ (thus, we can take $N$ as big as we want). By~\eqref{eqetata}, for every $z\in B(0,R)$, we obtain the desired bound:
\[\left\|h(x+z/N)-\big(h(x)+Dh_x(z/N)\big)\right\| < \frac{1}{3\ell(2n+1)\card I_k}\cdot\frac{1}{N}.\]
Combined with \eqref{ModuloGene}, this leads to
\begin{equation}\label{eqetatata}
\frac{\card\Big(h_N\big(B(x,R/N)\big)\, \Delta\, P_N\big(h(x) + Dh_x (B(0,R/N))\big)\Big)}{\card\big(B(x,R/N) \cap E_N\big)} \le \frac{1}{3\ell \card I_k};
\end{equation}
in other words, on every ball of radius $R/N$, the image of $E_N$ by $h_N$ and the discretization of the linearization of $h$ are almost the same (that is, up to a proportion $1/(3\ell \card I_k)$ of points).

We now set $R_1 = R_0 \|f'\|_\infty^m$, and choose $R\ge R_1$, to which is associated a number $\rho>0$ and an order $N = \lceil R/\rho \rceil$, that we can choose large enough so that $2R/N \le \|f'\|_\infty$. We also choose $y\in C$. As
\[\card\big(h_N^m(E_N) \cap B(y,R/N) \big) = \card\left( \bigcup_{x\in h^{-m}(y)} h_N^m\big(B(x,R/N) \cap E_N\big) \cap B(y,R/N)\right),\]
and using the estimations~\eqref{DiffDensGene} and~\eqref{eqetatata}, we get
\[\left|\frac{\card\big(h_N^m(E_N) \cap B(y,R/N) \big)}{\card\big(B(y,R/N) \cap E_N\big)} -\overline D\big( (\det Df_x^{-1})_{\begin{subarray}{l} 1\le m\le k \\ x\in f^{-m}(y)\end{subarray}}\big)\right| < \frac{2}{3\ell}.\]
As such an estimation holds on a subset of $\T^n$ of measure bigger than $1-1/(3\ell)$, we get the conclusion of the theorem.
\end{proof}
\bigskip

We can easily adapt the proof of Lemma~\ref{LemTauxExpand} to the case of sequences of matrices, without the hypothesis of expansivity. This leads to the following improvement of Theorem~\ref{convBis}, for higher regularity.

\begin{lemme}\label{LemTauxDiff2}
For every $k\in\N$ and every $\ell', c\in\N$, there exists a locally finite union of positive codimension submanifolds $V_q$ of $(GL_n(\R))^{k}$ (respectively $(SL_n(\R))^{k}$) such that for every $\eta'>0$, there exists a radius $R_0>0$ such that if $(A_m)_{1\le m\le k}$ is a finite sequence of matrices of $(GL_n(\R))^{k}$ (respectively $(SL_n(\R))^{k}$) satisfying $d((A_m)_m, V_q)>\eta'$ for every $q$, then for every $R\ge R_0$, and every family $(v_m)_{1\le m \le k}$ of vectors of $\R^n$, we have
\[\left|D_R^+\left(\big( \pi(A_{k}+v_k) \circ \cdots \circ \pi( A_1+v_1)\big) (\Z^n) \right) - \det(A_k^{-1}\cdots A_1^{-1})\overline\tau^k(A_1,\cdots,A_k)\right| < \frac{1}{\ell'}\]
(the density of the image set is ``almost invariant'' under perturbations by translations), and for every $m\le k$, we have\footnote{Recall that $(\Z^n)'$ stands for the set of points of $\R^n$ at least one coordinate of which belongs to $\Z+1/2$.}
\[D_R^+\left\{x\in \big(A_{m} + v_{m}\big)(\Z^n)\ \middle\vert\ d\big(x,(\Z^n)' \big) < \frac{1}{{c}\ell'(2n+1)} \right\} < \frac{1}{{c}\ell'}\]
(there is only a small proportion of the points of the image sets which are obtained by discretizing points close to $(\Z^n)'$).
\end{lemme}

With the same proof as Theorem~\ref{TauxExpand}, Lemma~\ref{LemTauxDiff2} leads to the local-global formula for $C^r$-diffeomorphisms (Theorem~\ref{convBisMieux}).

\section{Asymptotic rate of injectivity for a generic dissipative diffeomorphism}\label{SecJairo}

First of all, we tackle the issue of the asymptotic rate of injectivity of generic dissipative diffeomorphisms. We will deduce it from a theorem of A.~Avila and J.~Bochi. Again, we will consider the torus $\T^n$ equipped with Lebesgue measure $\Leb$ and the canonical measures $E_N$, see Section~\ref{AddendSett} for a more general setting where the result is still true. The study of the rate of injectivity for generic dissipative diffeomorphisms is based on the following theorem of A.~Avila and J.~Bochi.


\begin{theoreme}[Avila, Bochi]\label{ArturJairo}
Let $f$ be a generic $C^1$ maps of $\T^n$. Then for every $\varep>0$, there exists a compact set $K\subset \T^n$ and an integer $m\in\N$ such that
\[\Leb(K) > 1-\varep \qquad \text{and} \qquad \Leb(f^m(K))<\varep.\]
\end{theoreme}

This statement is obtained by combining Lemma 1 and Theorem 1 of \cite{MR2267725}.

\begin{rem}
As $C^1$ expanding maps of $\T^n$ and $C^1$ diffeomorphisms of $\T^n$ are open subsets of the set of $C^1$ maps of $\T^n$, the same theorem holds for generic $C^1$ expanding maps and $C^1$ diffeomorphisms of $\T^n$ (this had already been proved in the case of $C^1$-expanding maps by A.~Quas in \cite{MR1688216}).
\end{rem}

This theorem can be used to compute the asymptotic rate of injectivity of a generic diffeomorphism.

\begin{coro}\label{CoroCoroArturJairo}
The asymptotic rate of injectivity of a generic dissipative diffeomorphism $f\in\Diff^1(\T^n,\Leb)$ is equal to 0. In particular, the degree of recurrence $D(f_N)$ of a generic dissipative diffeomorphism tends to 0 when $N$ goes to infinity.
\end{coro}

\begin{proof}[Proof of Corollary~\ref{CoroCoroArturJairo}]
The proof of this corollary mainly consists in stating which good properties can be supposed to possess the compact set $K$ of Theorem~\ref{ArturJairo}. Thus, for $f$ a generic diffeomorphism and $\varep>0$, there exists $m>0$ and a compact set $K$ such that $\Leb(K) > 1-\varep$ and $\Leb(f^m(K))<\varep$.

First of all, it can be easily seen that Theorem~\ref{ArturJairo} is still true when the compact set $K$ is replaced by an open set $O$: simply consider an open set $O'\supset f^m(K)$ such that $\Leb(O')<\varep$ (by regularity of the measure $\Leb$) and set $O = f^{-m}(O') \supset K$. We then approach the set $O$ by unions of dyadic cubes of $\T^n$: we define the cubes of order $2^M$ on $\T^n$
\[C_{M,i} = \prod_{j=1}^n\left[\frac{i_j}{2^M},\frac{i_j+1}{2^M} \right],\]
and set\label{pageCubes}
\[U_M = \operatorname{Int}\left(\overline{\bigcup_{C_{M,i} \subset O} C_{M,i} }\right),\]
where $\operatorname{Int}$ denotes the interior. Then, the union $\bigcup_{M\in\N} U_M$ is increasing in $M$ and we have $\bigcup_{M\in\N} U_M = O$. In particular, there exists $M_0\in\N$ such that $\Leb(U_{M_0})>1-\varep$, and as $U_{M_0}\subset O$, we also have $\Leb(f^m(U_{M_0}))<\varep$. We denote $U = U_{M_0}$. Finally, as $U$ is a finite union of cubes, and as $f$ is a diffeomorphism, there exists $\delta>0$ such that the measure of the $\delta$-neighbourhood of $f^m(U)$ is smaller than $\varep$. We call $V$ this $\delta$-neighbourhood.

As $U$ is a finite union of cubes, there exists $N_0\in\N$ such that if $N\ge N_0$, then the proportion of points of $E_N$ which belong to $U$ is bigger than $1-2\varep$, and the proportion of points of $E_N$ which belong to $V$ is smaller than $2\varep$. Moreover, if $N_0$ is large enough, then for every $N\ge N_0$, and for every $x_N\in E_N\cap U$, we have $f_N^m(x_N)\in V$. This implies that
\[\frac{\card(f_N^m(E_N))}{\card(E_N)}\le 4\varep,\]
which proves the corollary.
\end{proof}

\section{Asymptotic rate of injectivity for a generic conservative diffeomorphism}\label{ChapAsympto}

The goal of this section is to prove that the degree of recurrence of a generic conservative $C^1$-diffeomorphism is equal to 0. It will be obtained by using the local-global formula (Theorem~\ref{convBis}) and the result about the asymptotic rate of injectivity of a generic sequence of matrices (Theorem~\ref{ConjPrincip}). As a warm-up, we begin by a weak and very easy version of this statement, which already shows that the discretizations of generic conservative $C^1$-diffeomorphism do not behave in the same fashion as discretization of a generic conservative homeomorphism.

\begin{prop}\label{DegRecurDiffEasy}
For a generic conservative diffeomorphism $f\in\Diff^1(\T^n,\Leb)$, we have $\underset{N\to +\infty}{\overline\lim} \tau^1_N(f) < 1$, in particular $\underset{N\to +\infty}{\overline\lim} D(f_N) < 1$.
\end{prop}

\begin{proof}[Proof of Proposition \ref{DegRecurDiffEasy}]
To prove this proposition, we prove that on an open and dense set of diffeomorphisms $f$, we have $\underset{N\to +\infty}{\overline\lim} \tau^1_N(f) < 1$. This is easily obtained by perturbing locally any diffeomorphism. Indeed, we take a diffeomorphism $f\in\Diff^1(\T^n,\Leb)$ and a point $x\in\T^n$. Then, we perturb the differential of $f$ at $x$ such that $\overline \tau^1 (Df_x)<1$ (the set of matrices satisfying this property is locally included in a finite union of submanifolds of codimension at least 1, see Corollary~\ref{CoroHajos}). We then apply the local linearization theorem to get a diffeomorphism which is linear in a neighbourhood of $x$. Thus, $\overline \tau^1(f)<1$ by the local-global formula Theorem~\ref{convBis}, and this remains true on a whole neighbourhood of $f$ by the continuity of $\overline\tau^1$ (Corollary~\ref{ContTauBarDiff}).
\end{proof}

%

\begin{theoreme}\label{limiteEgalZero}
For a generic conservative diffeomorphism $f\in\Diff^1(\T^n,\Leb)$, we have
\[\lim_{t\to\infty}\tau^k(f)=0;\]
more precisely, for every $\varep>0$, the set of diffeomorphisms $f\in\Diff^1(\T^n,\Leb)$ such that $\lim_{t\to +\infty} \tau^k(f)<\varep$ is open and dense.

In particular\footnote{Using Equation~\eqref{intervLim} page~\pageref{intervLim}.}, we have $\lim_{N\to +\infty} D(f_N) = 0$.
\end{theoreme}

\begin{proof}[Proof of theorem \ref{limiteEgalZero}]
We show that for every $\ell\in \N$ and every $\varep>0$, the set of conservative diffeomorphisms such that $\overline\lim_{t\to\infty}\tau_\infty^t < 1/\ell+\varep$ contains an open dense subset of $\Diff^1(\T^n,\Leb)$. To begin with, we fix $f\in \Diff^1(\T^n,\Leb)$ and $\delta>0$ (which will be a size of perturbation of $f$). By Theorem~\ref{ConjPrincip}, and in particular Equation~\eqref{EstimTauxExp}, there exists a parameter $\lambda\in]0,1[$ (depending only on $\delta$, $\ell$ and $\|f\|_{C^1}$), such that for every sequence $(A_k)_{k\ge 1}$ of linear maps in $SL_n(\R)$, there exists a sequence $(B_k)_{k\ge 1}$ of (generic) linear maps in $SL_n(\R)$ such that for each $k$, we have $\|A_k-B_k\|\le\delta$ and $\tau^{\ell k}(B_1,\cdots,B_{\ell k})\le \lambda^k + 1/\ell$ (as the sequence is generic, this property remains true on a whole neighbourhood of $(B_k)_{k\ge 1}$, see Remark~\ref{conttaukk}). From that parameter $\lambda$, we deduce a time $k_0>0$ such that
\[\frac{1}{k_0} \sum_{k=1}^{k_0} \lambda^k = \frac{\lambda}{k_0} \frac{1-\lambda^{k_0}}{1-\lambda} < \varep/100.\]

Applying a classical technique in this context (see for example \cite{MR1944399}), we use a Rokhlin tower of height $k_0$ with an open basis $U$:
\begin{itemize}
\item The sets $U, f(U), \cdots, f^{k_0-1}(U)$ are pairwise disjoint;
\item the measure of the union of the ``floors'' $U \cup f(U) \cup \cdots \cup f^{k_0-1}(U)$ is bigger than $1-\varep/100$
\end{itemize}
For the existence of such towers, see for example \cite[Lemme~6.8]{MR2931648} or \cite[Chapter ``Uniform topology'']{MR0097489}.

We then approach the basis $U$ by a union of cubes of a dyadic subdivision of $\T^n$ (as in page~\pageref{pageCubes}, from now we suppose that $U$ is a union of such cubes). If that dyadic subdivision is fine enough, on each cube $C$, it is possible to perturb $f$ into a diffeomorphism $g$ such that on each set $(1-\varep/100)C$, $g$ is affine and irrational, using Lemma~\ref{LemExtension} and Franks lemma (see \cite{MR0283812} or \cite{MR2288283}), which is valid only in the $C^1$ topology. We do the same thing on the $k_0-1$ first images of each cube and perturb $g$ such that on each set $g^k((1-\varep/100)C)$, the perturbed diffeomorphism $h$ is linear and equal to $B_k$. By what we have said at the beginning of the proof, we can moreover suppose that the sequence $B_1, \cdots, B_{k_0}$ of linear maps is generic, and satisfies $\tau^k(B_1,\cdots,B_k)\le \lambda^k + 1/\ell$ for every $k\le k_0-1$. By the choice of $k_0$ we have made, this implies that
\begin{align*}
\tau^{k_0}(h) \le & \sum_{k=0}^{k_0-1} \Leb(h^{k}(U))\tau^{k_0} \big(h_{|h^{k}(U)}\big)\\
                  & + \Leb\big((U \cup \cdots \cup h^{k_0-1}(U))^\complement\big) \tau^{k_0} \big(h_{(U \cup \cdots \cup h^{k_0-1}(U))^\complement}\big)\\
              \le & \sum_{k=0}^{k_0-1} \Leb(h^{k}(U))\tau^{k_0-k} \big(h_{|h^{k}(U)}\big) + \Leb\big((U \cup \cdots \cup h^{k_0-1}(U))^\complement\big)\\
              \le & \Leb(U)\sum_{k=0}^{k_0-1} \left(\lambda^{k_0-k} + \frac{1}{\ell}\right) + \varep/100\\
              \le & \frac{1}{k_0}\sum_{k=1}^{k_0} (\lambda^{k} + 1)/2 + \varep/100\\
						  \le & 1/\ell + \varep/2.
\end{align*}
Moreover, the differentials of $h$ form generic sequences on a set of measure at least $1-\varep/10$. This implies that the rate of injectivity os continuous in $h$ when restricted to this subset of $\T^n$. Thus, the inequality $\tau^{k_0}(h) \le 1/\ell + \varep$ still holds on a whole neighbourhood of $h$. This proves the theorem.
\end{proof}

\subsection*{Numerical simulation}

\begin{figure}[t]
\begin{center}
\includegraphics[width=.5\linewidth,trim = .5cm .3cm .6cm .1cm,clip]{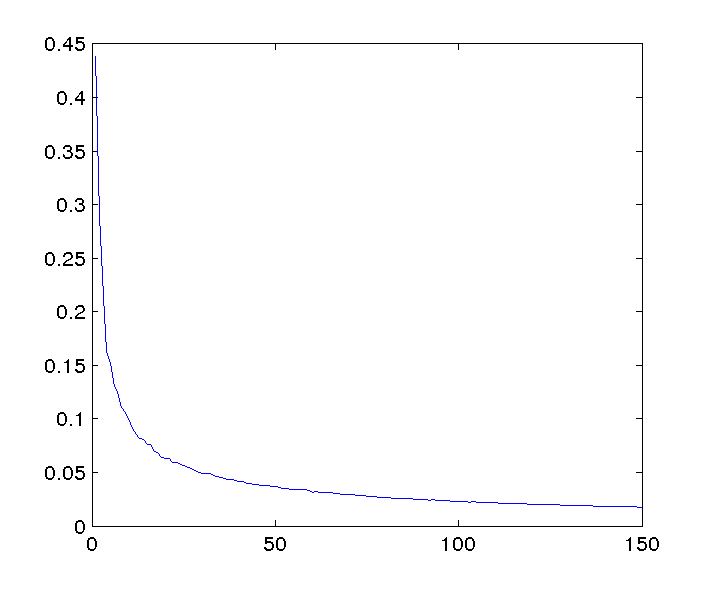}
\caption[Simulation of the degree of recurrence of $(f_5)_N$]{Simulation of the degree of recurrence $D((f_5)_N)$ depending on $N$, on the grids $E_N$ with $N=128k$, $k=1,\cdots,150$.}\label{GrafDnDiffeoCons}
\end{center}
\end{figure}

We recall the results of numerical simulations we have presented in Part~\ref{PartOne} (Figure~\ref{GrafCons}). We have computed numerically the degree of recurrence of the diffeomorphism $f_5$, which is $C^1$-close to $\Id$. Recall that $f_5 = Q\circ P$, with
\[P(x,y) = \big(x,y+p(x)\big)\quad\text{and}\quad Q(x,y) = \big(x+q(y),y\big),\]
\[p(x) = \frac{1}{209}\cos(2\pi\times 17x)+\frac{1}{271}\sin(2\pi\times 27x)-\frac{1}{703}\cos(2\pi\times 35x),\]
\[q(y) = \frac{1}{287}\cos(2\pi\times 15y)+\frac{1}{203}\sin(2\pi\times 27y)-\frac{1}{841}\sin(2\pi\times 38y).\]

On Figure~\ref{GrafDnDiffeoCons}, we have represented graphically the quantity $D((f_5)_{128k})$ for $k$ from 1 to $150$. It appears that, as predicted by Theorem~\ref{limiteEgalZero}, this degree of recurrence goes to 0. In fact, it is even decreasing, and converges quite fast to 0: as soon as $N=128$, the degree of recurrence is smaller than $1/2$, and if $N\gtrsim 1000$, then $D((f_5)_N) \le 1/10$.

\section{Asymptotic rate of injectivity of a generic $C^r$ expanding map}

In this section, we prove that the asymptotic rate of injectivity of a generic expanding map is equal to 0. Note that a local version of this result was already obtained by P.P.~Flockermann in his thesis (Corollary~2 page 69 and Corollary~3 page 71 of \cite{Flocker}), stating that for a generic $C^{1+\alpha}$ expanding map $f$ of the circle, the ``local asymptotic rate of injectivity'' is equal to 0 almost everywhere. Some of his arguments will be used in this section. Note also that in $C^1$ regularity, the equality $\tau^\infty(f) = 0$ for a generic $f$ is a consequence of Theorem~\ref{ArturJairo} of A.~Avila and J.~Bochi (see also Corollary~\ref{CoroCoroArturJairo}); the same theorem even proves that the asymptotic rate of injectivity of a generic $C^1$ endomorphism of the circle is equal to 0.

\begin{definition}
We define $Z_m$ as the number of children at the $m$-th generation in $G_{f,y}$ (see Definition~\ref{DefTreeMap}).
\end{definition}

\begin{prop}\label{PropTruc}
For every $r\in]1,+\infty]$, for every $f\in \mathcal D^r(\Sp^1)$ and every $y\in\Sp^1$, we have
\[\mathbf{P}(Z_m>0) \underset{m\to+\infty}{\longrightarrow}0.\]
Equivalently,
\[\overline D\big( (\det Df_x^{-1})_{\begin{subarray}{l} 1\le m\le k \\ x\in f^{-m}(y)\end{subarray}}\big) \underset{k\to+\infty}{\longrightarrow}0.\]
\end{prop}

\begin{lemme}\label{Espoir}
The expectation of $Z_m$ satisfies
\[\E(Z_m) = (\Ll^m 1)(y),\]
where $\Ll$ is the Ruelle-Perron-Frobenius associated to $f$ and $1$ denotes the constant function equal to 1 on $\Sp^1$.
In particular, there exists a constant $\Sigma_0>0$ such that $\E(Z_m) \le \Sigma_0$ for every $m\in\N$.
\end{lemme}

The second part of the lemma is deduced from the first one by applying the theorem stating that for every $C^r$ expanding map $f$ of $\Sp^1$ ($r>1$), the maps $\Ll^m 1$ converge uniformly towards a Hölder map, which is the density of the unique SRB measure of $f$ (see for example \cite{MR2504311}). The first assertion of the lemma follows from the convergence of the operators $f_N^*$ acting on $\Prb$ (the space of Borel probability measures) towards the Ruelle-Perron-Frobenius operator.

\begin{definition}
The \emph{transfer operator} associated to the map $f$ (usually called Ruelle-Perron-Fro\-be\-nius operator), which acts on densities of probability measures, will be denoted by $\Ll_f$\index{$\Ll_f$}. It is defined by
\[\Ll_f \phi(y) = \sum_{x\in f^{-1}(y)} \frac{\phi(x)}{f'(x)}.\]
\end{definition}

Lemma~\ref{Espoir} follows directly from the following lemma.

\begin{lemme}\label{LemRPF}
Denoting $\Leb_N$ the uniform measure on $E_N$, for every $C^1$ expanding map of $\Sp^1$ and every $m\ge 0$, we have convergence of the measures $(f_N^*)^m(\lambda_N)$ towards the measure of density $\Ll_f^m 1$ (where $1$ denotes the constant function equal to 1).
\end{lemme}

The proof of this lemma is straightforward but quite long. We sketch here this proof, the reader will find a complete proof using generating functions in Section~3.4 of \cite{Flocker} and a quantitative version of it in Section~12.2 of \cite{Guih-These}.

\begin{proof}[Sketch of proof of Lemma~\ref{LemRPF}]
As $f$ is $C^1$, by the mean value theorem, for every segment $I$ small enough, we have
\[\left| \Leb(I) - \frac{\Leb(f(I))}{f'(x_0)}\right| \le \varep.\]
Moreover, for every interval $J$, 
\[\left|\Leb(J) - \frac{\card(J\cap E_N)}{\card(E_N)}\right| \le \frac{1}{N}.\]
These two inequalities allows to prove the local convergence of the measures $f_N^*(\lambda_N)$ towards the measure with density $\Ll_f 1$. The same kind of arguments holds in arbitrary times, and allow to prove the lemma.
\end{proof}

\begin{proof}[Proof of Proposition~\ref{PropTruc}]
We fix $\varep>0$, and set $K=\lceil\Sigma_0/\varep\rceil$ (the constant $\Sigma_0$ being given by Lemma~\ref{Espoir}) and
\[a_m = \mathbf{P}(Z_m = 0)\,,\quad  b_m = \mathbf{P}(0 < Z_m\le K)\,,\quad  c_m = \mathbf{P}(Z_m > K).\]

\begin{figure}[b]
\begin{center}
\begin{tikzpicture}[scale=1]
\node[draw,ellipse] (C) at (6,0) {\small$c_m\le \varep$};
\node[draw,circle,minimum height=1.2cm] (B) at (3,0) {\large$b_m$};
\node[draw,circle,minimum height=1.5cm] (A) at (0,0) {\large$a_m$};
\draw[->,>=latex,shorten >=3pt, shorten <=3pt] (C) to[bend left] (B);
\draw[->,>=latex,shorten >=3pt, shorten <=3pt] (B) to[bend left] (C);
\draw[->,>=latex,shorten >=3pt, shorten <=3pt] (C) to[bend left] (A);
\draw[->,>=latex,shorten >=3pt, shorten <=3pt] (B) to[bend right] node[midway, above]{$\ge \alpha$} (A) ;
\draw[->,>=latex,shorten >=3pt, shorten <=3pt] (A) to[loop left,looseness=16,min distance=10mm] (A);
\end{tikzpicture}
\end{center}
\caption{Transition graph for $Z_m$: $Z_m=0$, $0<Z_m\le K$ and $Z_m>K$.}\label{GraphTransZ}
\end{figure}
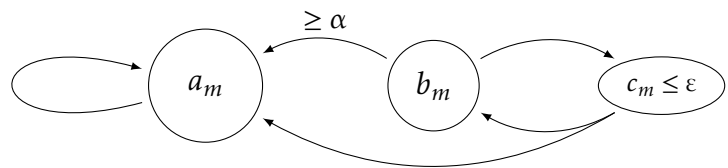

\noindent We want to prove that the sequence $(a_m)_{m\in\N}$ tends to 1.

Of course, $a_m+b_m+c_m=1$ for any $m$. We then remark that a generation with less than $K$ children will give birth to zero child with positive probability: we have the inequality
\[\mathbf{P}(Z_{m+1}=0 \mid 0 < Z_m\le K) \ge \left( 1 - \frac{1}{\|f'\|_\infty}\right)^{dK}.\]
In other words, setting $\alpha = ( 1 - \|f'\|_\infty^{-1})^{dK}$, we get $a_{m+1} \ge a_m + \alpha b_m$ (see also Figure~\ref{GraphTransZ}).

Furthermore, by Markov inequality and Lemma~\ref{Espoir}, we have
\[\mathbf{P}(Z_m\ge \Sigma_0/\varep) \le \varep,\]
so $c_m\le\varep$.

In summary, we have
\[\left\{\begin{array}{l}
a_m + b_m + c_m = 1\\
c_m\le \varep\\
a_{m+1} \ge a_m + \alpha b_m.
\end{array}\right.\]
An easy computation leads to $a_{m+1}\ge (1-\alpha) a_m + \alpha(1-\varep)$, which implies that $\liminf a_m \ge 1-\varep$. As this holds for any $\varep>0$, we get that $\lim a_m = 1$.
\end{proof}

\begin{coro}\label{Tau0Dilat}
For every $r\in]1,+\infty]$ and a generic map $f\in \mathcal D^r(\Sp^1)$, we have $\tau^\infty(f) = 0$. In particular, $\lim_{N\to +\infty} D(f_N) = 0$.
\end{coro}

\begin{proof}[Proof of Corollary~\ref{Tau0Dilat}]
It is an easy consequence of the local-global formula (Theorem~\ref{TauxExpand}), Proposition~\ref{PropTruc} and the dominated convergence theorem.
\end{proof}

\section{A more general setting where the theorems are still true}\label{AddendSett}

Here, we give weaker assumptions under which the theorems of Sections~\ref{ChapLocGlob}, \ref{SecJairo}, \ref{ChapAsympto} and \ref{Section!} (that is, Theorems~\ref{conv}, \ref{convBis}, \ref{limiteEgalZero} and \ref{TheoMesPhysDiff}, and Corollaries~\ref{ContTauBarDiff}, \ref{CoroCoroArturJairo} and \ref{derdesder?}) are still true: the framework ``torus $\T^n$ with grids $E_N$ and Lebesgue measure'' could be seen as a little too restrictive.
\bigskip

So, we take a compact smooth manifold $M$ (possibly with boundary) and choose a partition $M_1,\cdots,M_k$ of $M$ into closed sets\footnote{That is, $\bigcup_i M_i = M$, and for $i\neq j$, the intersection between the interiors of $M_i$ and $M_j$ are empty.} with smooth boundaries, such that for every $i$, there exists a chart $\varphi_i : M_i\to \R^n$. We endow $\R^n$ with the euclidean distance, which defines a distance on $M$ \emph{via} the charts $\phi_i$ (this distance is not necessarily continuous). From now, we study what happens on a single chart, as what happens on the neighbourhoods of the boundaries of these charts ``counts for nothing'' from the Lebesgue measure viewpoint.

Finally, we suppose that the uniform measures on the grids $E_N = \bigcup_i E_{N,i}$ converge to a smooth measure $\lambda$ on $M$ when $N$ goes to infinity.

This can be easily seen that these conditions are sufficient for Corollary~\ref{CoroCoroArturJairo} to be still true.
\bigskip

For the rest of the statements of this chapter, we need that the grids behave locally as the canonical grids on the torus.

For every $i$, we choose a sequence $(\kappa_{N,i})_N$ of positive real numbers such that $\kappa_{N,i}\underset{N\to +\infty}{\longrightarrow} 0$. This defines a sequence $E_{N,i}$ of grids on the set $M_i$ by $E_{N,i} = \varphi_i^{-1} (\kappa_{N,i}\Z^n)$. Also, the canonical projection $\pi : \R^n\to \Z^n$ (see Definition~\ref{DefDiscrLin}) allows to define the projection $\pi_{N,i}$, defined as the projection on $\kappa_{N,i}\Z^n$ in the coordinates given by $\varphi_i$:
\[\begin{array}{rcl}
\pi_{N,i} : & M_i & \longrightarrow E_{N,i}\\
            & x   & \longmapsto \varphi_i^{-1}\Big(\kappa_{N,i}\pi\big(\kappa_{N,i}^{-1}\varphi_i(x)\big)\Big).
\end{array}\]

We easily check that under these conditions, Theorems~\ref{conv}, \ref{convBis} and \ref{limiteEgalZero} are still true, that is if we replace the torus $\T^n$ by $M$, the uniform grids by the grids $E_N$, the canonical projections by the projections $\pi_{N,i}$, and Lebesgue measure by the measure $\lambda$.

\subsection*{Numerical simulation}

\begin{figure}[b]
\begin{center}
\includegraphics[width=.5\linewidth,trim = .5cm .3cm .6cm .1cm,clip]{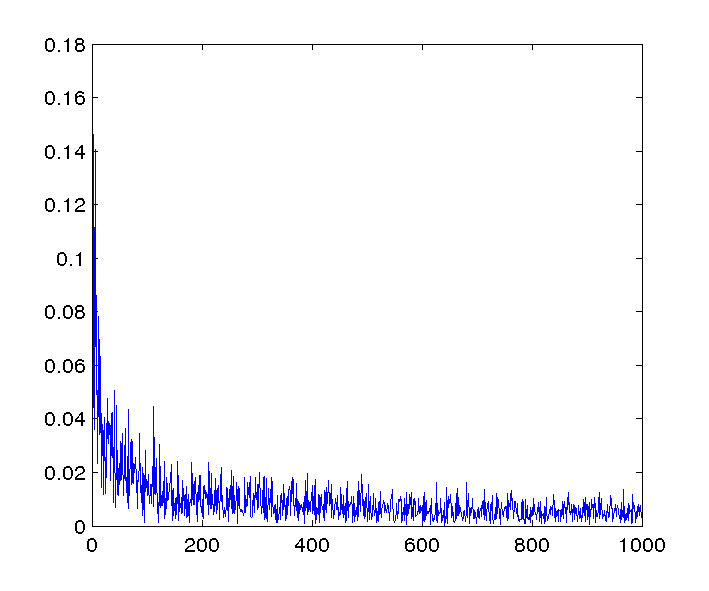}
\caption[Simulation of the degree of recurrence of $f_N$]{Simulation of the degree of recurrence $D((f_5)_N)$ of the expanding map $N$, depending on $N$, on the grids $E_N$ with $N=128k$, $k=1,\cdots,1\,000$.}\label{GrafDnExp}
\end{center}
\end{figure}

We present the results of the numerical simulation we have conducted for the degree of recurrence of the expanding map of the circle $f$, defines dy
\[f(x) = 2x + \varep_1 \cos(2\pi x) + \varep_2 \sin(6\pi x),\]
with $\varep_1 = 0.127\,943\,563\,72$ and $\varep_2 = 0.008\,247\,359\,61$.

On Figure~\ref{GrafDnExp}, we have represented the quantity $D((f)_{128k})$ for $k$ from 1 to $1\,000$. It appears that, as predicted by the theorem of P.P.~Flockermann (Corollary~2 page 69 and Corollary~3 page 71 of \cite{Flocker}), this degree of recurrence seems to tend to 0. In fact, it is even decreasing, and converges quite fast to 0: as soon as $N=128$, the degree of recurrence is smaller than $1/5$, and if $N\gtrsim 25\,000$, then $D((f)_N) \le 1/50$.

\chapter[Physical measures]{Physical measures of discretizations of generic diffeomorphisms and expanding maps}\label{chapPhys}

\section{Physical measures of discretizations of generic diffeomorphisms}\label{SecPhys12}

In this section, we will consider the torus $\T^n$, equipped with Lebesgue measure $\Leb$ and the uniform grids 
\[E_N = \left\{ \left(\frac{i_1}{N},\cdots,\frac{i_n}{N}\right)\in \R^n/\Z^n \middle\vert\ 1\le i_1,\cdots,i_n \le N\right\}.\]
Fore a more general setting where the theorems of this section are still true, see Section~\ref{AddendSett}.

This chapter is devoted to the study of the physical measures of the discretizations of a generic conservative diffeomorphism. Recall the classical definition of a physical measure for a map $f$: a Borel probability measure $\mu$ is called \emph{physical} for the map $f$ if its basin of attraction has positive Lebesgue measure, where the basin of attraction of $\mu$ is the set
\[\left\{x\in \T^n\ \Big\vert\ \frac{1}{M} \sum_{m=0}^{M-1} \delta_{f^m(x)} \underset{M\to+\infty}{\longrightarrow} \mu\right\}.\]
(see also Definition~\ref{sport}). Heuristically, the physical measures are the ones that can be observed in practical experiments, because they are ``seen'' by a set of points $x$ of positive Lebesgue measure. Here, our aim is to study similar concepts in the view of discretizations: which measures can be seen by the discretizations of generic conservative $C^1$-diffeomorphisms?

Remember that we denote by $\mu_{x}^{f_N}$ the limit of the Birkhoff sums
\[\frac{1}{M} \sum_{m=0}^{M-1} \delta_{f_N^m(x_N)}.\]
More concretely, $\mu_{x}^{f_N}$ is the $f_N$-invariant probability measure supported by the periodic orbit on which the positive orbit of $x_N=P_N(x)$ falls after a while (see also Definition~\ref{defmes} page~\pageref{defmes}). We would like to know the answer the following question: for a generic conservative $C^1$-diffeomorphism $f$, does the sequence of measure $\mu_{x}^{f_N}$ tend to a physical measure of $f$ for most of the points $x$ as $N$ goes to infinity?

Recall that in the $C^0$ case, we have proved Theorem~\ref{mesinv}, which implies in particular that for a generic homeomorphism $f\in\Hom(\T^n,\Leb)$ and every $x\in\T^n$, the measures $\mu^{f_N}_x$ accumulate on the whole set of $f$-invariant measures when $N$ goes to infinity (moreover, given an $f$-invariant measure $\mu$, the sequence $(N_k)_{k\ge 0}$ such that $\mu^{f_{N_k}}_x$ tends to $\mu$ can be chosen independently of $x$). In a certain sense, this theorem in the case of homeomorphisms expresses that from the point of view of the discretizations, all the $f$-invariant measures are physical.

In the $C^1$-case, we have already proved Corollary~\ref{CoroMane} (in Chapter~\ref{ChapPerturbLem}), which states that for a generic conservative $C^1$-diffeomorphism $f$, any $f$-invariant measure is the limit of a sequence of $f_N$-invariant measures. This corollary has been obtained as a simple consequence of an ergodic closing lemma of R.~Mañé and F.~Abdenur, C.~Bonatti and S.~Crovisier (Theorem~\ref{Mane}); however it does not say anything about the basin of attraction of these discrete measures.

In the theorem we prove in this chapter (Theorem~\ref{TheoMesPhysDiff}), we improve the previous statement for generic conservative $C^1$-diffeomorphisms, in order to describe the basin of attraction of the discrete measures. In particular, we prove that for points $x$ belonging to a generic subset of points of the torus, the measures $\mu^{f_N}_x$ accumulate on the whole set of $f$-invariant measures (Theorem~\ref{TheoMesPhysDiff}). Notice that given an $f$-invariant measure $\mu$, the sequence $(N_k)_k$ such that $\mu_{x}^{f_{N_k}}$ converges to $\mu$ depends on the point $x$, contrary to what happens in the $C^0$ case.

Moreover, if we fix a countable subset $D\subset \T^n$, then for a generic conservative $C^1$-diffeomorphism $f$ and for any $x\in D$, the measures $\mu^{f_N}_x$ accumulate on the whole set of $f$-invariant measures (Addendum~\ref{AddTheoMesPhysDiff}). This is a process that is usually applied in practice to detect the $f$-invariant measures: fix a finite set $D\subset \T^n$ and compute the measure $\mu^{f_N}_x$ for $x\in D$ and for a large order of discretization $N$. The theoretical result expresses that it is possible that the measure that we observe on numerical experiments is very far away from the physical measure.

\label{AvilaErgo}Note that in the space $\Diff^1(\T^n,\Leb)$, there are open sets where generic diffeomorphisms are ergodic: the set of Anosov diffeomorphisms is open in $\Diff^1(\T^n,\Leb)$, and a generic Anosov conservative $C^1$-diffeomorphism is ergodic (it is a consequence of the fact that any $C^2$ Anosov conservative diffeomorphism is ergodic, together with the theorem of regularization of conservative diffeomorphisms of A.~Avila \cite{MR2736152}). More generally, A.~Avila, S.~Crovisier and A.~Wilkinson have set recently in \cite{ArturSylvain} a generic dichotomy for a conservative diffeomorphism $f$: either $f$ is ergodic, either all the Lyapunov exponents of $f$ vanish. In short, there are open sets where generic conservative diffeomorphisms have only one physical measure; in this case, our result asserts that this physical measure is not detected on discretizations by computing the measures $\mu^{f_N}_x$.

Recall that results of stochastic stability are known to be true in various contexts (for example, expanding maps \cite{MR884892},\cite{MR874047}, \cite{MR685377}, uniformly hyperbolic attractors \cite{MR874047}, \cite{MR857204}, etc.). These theorems suggest that the physical measures can always be observed in practice, even if the system is noisy. Our Theorem~\ref{TheoMesPhysDiff} indicates that the effects of discretizations (i.e. numerical truncation) might be quite different from those of a random noise.

However, we shall remark that the same proof as for Theorem~\ref{TheoMesPhysDiff} implies that for a generic diffeomorphism $f\in\Diff^1(\T^n,\Leb)$ and a generic point $x\in\T^n$ (or equivalently, for any $x\in\T^n$ and for a generic $f\in\Diff^1(\T^n,\Leb)$), the measures
\[\mu^f_{x,m} = \frac 1m \sum_{i=0}^{m-1}{f}_*^i \delta_x\]
accumulate on the whole set of $f$-invariant measures.
\bigskip

The proof of Theorem~\ref{TheoMesPhysDiff} uses crucially the results of Part~\ref{PartII} on the fact that the asymptotic rate is null (in particular Lemma~\ref{DerTheoPart2}); it also uses two connecting lemmas (Theorem~\ref{PapySylvain} of \cite{MR2090361} and an improvement of Theorem~\ref{Mane} of \cite{MR2811152}).
\bigskip

At the end of this chapter, we present numerical experiments simulating the measures $\mu_x^{f_N}$ for some examples of conservative $C^1$-diffeomorphisms $f$ of the torus. The results of these simulations are quite striking for an example of $f$ $C^1$-close to $\Id$ (see Figure~\ref{MesPhysIdC1}): even for very large orders $N$, the measures $\mu_x^{f_N}$ do not converge to Lebesgue measure at all, and are very different ones from the others. This illustrates perfectly Theorem~\ref{TheoMesPhysDiff} (more precisely, Addendum~\ref{AddTheoMesPhysDiff}), which states that if $x$ is fixed, then for a generic $f\in\Diff^1(\T^2,\Leb)$, the measures $\mu_x^{f_N}$ accumulate on the whole set of $f$-invariant measures, but do not say anything about, for instance, the frequency of orders $N$ such that $\mu_x^{f_N}$ is not close to Lebesgue measure. Moreover, the same phenomenon (although less pronounced) occurs for diffeomorphisms close to a translation of $\T^2$ (Figure~\ref{MesPhysRotC1}) or a linear Anosov 
automorphism (Figure~\ref{MesPhysAnoC1}).

\subsection{Statement of the theorem}\label{Section!}

\begin{theoreme}\label{TheoMesPhysDiff}
For a generic diffeomorphism $f\in\Diff^1(\T^n,\Leb)$, for a generic point $x\in X$, for any $f$-invariant probability measure $\mu$, there exists a subsequence $(N_k)_k$ of discretizations such that
\[\mu_{x}^{f_{N_k}} \underset{k\to+\infty}{\longrightarrow} \mu.\]
\end{theoreme}

Remark that the theorem in the $C^0$ case is almost the same, except that here, the starting point $x\in\T^n$ is no longer arbitrary but has to be chosen in a generic subset of the torus, and that the sequence $(N_k)_k$ depends on the starting point $x$. The proof of this theorem will also lead to the two following statements.

\begin{add}\label{AddTheoMesPhysDiff}
For a generic diffeomorphism $f\in\Diff^1(\T^n,\Leb)$, for any $\varep>0$ there exists a $\varep$-dense subset $(x_1,\cdots,x_m)$ such that for any $f$-invariant probability measure $\mu$, there exists a subsequence $(N_k)_k$ of discretizations such that for every $j$,
\[\mu_{x_j}^{f_{N_k}} \underset{k\to+\infty}{\longrightarrow} \mu.\]

Also, for any countable subset $D\subset \T^n$, for a generic diffeomorphism $f\in\Diff^1(\T^n,\Leb)$, for any $f$-invariant probability measure $\mu$, and for any finite subset $E\subset D$, there exists a subsequence $(N_k)_k$ of discretizations such that for every $x\in E$, we have
\[\mu_{x}^{f_{N_k}} \underset{k\to+\infty}{\longrightarrow} \mu.\]
\end{add}

The first statement asserts that if $f$ is a generic conservative $C^1$-diffeomorphism, then for any $f$-invariant measure $\mu$, there exists an infinite number of discretizations $f_N$ which possess an invariant measure which is close tu $\mu$, and whose basin of attraction is $\varep$-dense. Basically, for an infinite number of $N$ any $f$-invariant will be seen from any region of the torus.

In the second statement, a countable set of starting points of the experiment is chosen ``by the user''. This is quite close to what happens in practice: we take a finite number of points $x_1,\cdots,x_m$ and compute the measures $\mu_{x_m,T}^{f_{N_k}}$ for all $m$, for a big $N\in\N$ and for ``large'' times $T$ (we can expect that $T$ is large enough to have $\mu_{x_m,T}^{f_{N_k}} \simeq \mu_{x_m}^{f_{N_k}}$). In this case, the result expresses that it may happen (in fact, for arbitrarily large $N$) that the measures $\mu_{x_m,T}^{f_{N_k}}$ are not close to the physical measure of $f$ but are rather chosen ``at random'' among the set of $f$-invariant measures.
\bigskip

We also have a dissipative counterpart of Theorem~\ref{TheoMesPhysDiff}, whose proof is easier.

\begin{theoreme}\label{TheoMesPhysDiffDissip}
For a generic dissipative diffeomorphism $f\in\Diff^1(\T^n)$, for any $f$-invariant probability measure $\mu$ such that the sum of the Lyapunov exponents of $\mu$ is negative (or equal to 0), for a generic point $x$ belonging to the same chain recurrent class as $\mu$, there exists a subsequence $(N_k)_k$ of discretizations such that
\[\mu_{x}^{f_{N_k}} \underset{k\to+\infty}{\longrightarrow} \mu.\]
\end{theoreme}

Remark that if we also consider the inverse $f^{-1}$ of a generic diffeomorphism $f\in\Diff^1(\T^n)$, we can recover any invariant measure $\mu$ of $f$ by looking at the measures $\mu_{x}^{f_{N_k}}$ for generic points $x$ in the chain recurrent class of $\mu$.

The proof of this result is obtained by applying Lemma~\ref{RemplaceDissip} during the proof of Theorem~\ref{TheoMesPhysDiff}.
\bigskip

We also have the same statement as Theorem~\ref{TheoMesPhysDiff} but for expanding maps of the circle. We denote $\mathcal E_d^1(\Sp^1)$ the set of $C^1$-expanding maps of the circle of degree $d$.

\begin{prop}\label{TheoMesPhysDiffExp}
For a generic expanding map $f\in\mathcal E_d^1(\Sp^1)$, for any $f$-invariant probability measure $\mu$, for a generic point $x\in\Sp^1$, there exists a subsequence $(N_k)_k$ of discretizations such that
\[\mu_{x}^{f_{N_k}} \underset{k\to+\infty}{\longrightarrow} \mu.\]
\end{prop}

The proof of this statement is far easier than that of Theorem~\ref{TheoMesPhysDiff} , as it can be obtained by coding any expanding map of class $C^1$ (that is, any $f\in\mathcal E_d^1(\Sp^1)$ is homeomorphic to a full shift on a set with $d$ elements).
\bigskip

We will use the connecting lemma for pseudo-orbits (see \cite{MR2090361}, see also Theorem~\ref{PapySylvain}), together with an ergodic closing lemma (adapted from \cite{MR2811152}) and the results of Part~\ref{PartII} on the fact that the asymptotic rate is null (in particular Lemma~\ref{DerTheoPart2}), to prove that any invariant measure of the diffeomorphism can be observed by starting at any point of a generic subset of $\T^n$.

By Baire theorem and the fact that for a generic conservative diffeomorphism, a generic invariant measure is ergodic, non periodic and has no zero Lyapunov exponent (see Theorem~3.5 of \cite{MR2811152}), the proof of Theorem~\ref{TheoMesPhysDiff} can be reduced easily to that of the following approximation lemma.

\begin{lemme}\label{LemMesPhysDiff}
For every $f\in\Diff^1(\T^n,\Leb)$, for every $f$-invariant measure $\mu$ which is ergodic, not periodic and has no zero Lyapunov exponent, for every open subset $U\subset \T^n$, for every $C^1$-neighbourhood $\mathcal V$ of $f$, for every $\varep>0$ and every $N_0\in\N$, there exists $g\in\Diff^1(\T^n,\Leb)$ such that $g\in \mathcal V$, there exists $y\in U$ and $N\ge N_0$ such that $\dist(\mu, \mu^{g_N}_y)<\varep$. Moreover, we can suppose that this property remains true on a whole neighbourhood of $g$.
\end{lemme}

First of all, we explain how to deduce Theorem~\ref{TheoMesPhysDiff} from Lemma~\ref{LemMesPhysDiff}.

\begin{proof}[Proof of Theorem~\ref{TheoMesPhysDiff}]
We consider a sequence $(\nu_\ell)_{\ell\ge 0}$ of Borel probability measures, which is dense in the whole set of probability measures. We also consider a sequence $(U_i)_{i\ge 0}$ of open subsets of $\T^n$ which spans the topology of $\T^n$. This allows us to set
\[\mathcal S_{\,N_0,k_0,\ell,i} = \left\{ f\in\Diff^1(\T^n,\Leb)\ \middle\vert\ \begin{array}{l}\exists \mu\ f\text{-inv.} : \dist(\mu,\nu_\ell)\le 1/k_0\implies\\
\exists N\ge N_0, y\in U_i : d(\mu^{f_N}_y ,\nu_\ell)<2/k_0 \end{array}\right\}.\]

We easily see that the set
\[\bigcap_{N_0,k_0,\ell,i \ge 0} \mathcal S_{\,N_0,k_0,\ell,i}\]
in contained in the set of diffeomorphisms satisfying the conclusions of the theorem.

It remains to prove that each set $\mathcal S_{\,N_0,k_0,\ell,i}$ contains an open and dense subset of $\Diff^1(\T^n,\Leb)$. Actually the interior of each set $\mathcal S_{\,N_0,k_0,\ell,i}$ is dense. This follows from\footnote{See the proof of Theorem \ref{EnsMesInv} for more details on the same kind of arguments.} the upper semi-continuity of the set of $f$-invariant measures with respect to $f$ and from the combination of Lemma~\ref{LemMesPhysDiff} with the fact that for a generic diffeomorphism, a generic invariant measure is ergodic, non periodic and has no zero Lyapunov exponent (see Theorem~3.5 of \cite{MR2811152}).
\end{proof}

The rest of this chapter is devoted to the proof of Lemma~\ref{LemMesPhysDiff}. We now outline the main arguments of this quite long and technical proof.

\paragraph{Sketch of proof of Lemma~\ref{LemMesPhysDiff}.}
First of all, we take a point $x\in\T^n$ which is typical for the measure $\mu$. In particular, by an ergodic closing lemma derived from that of F.~Abdenur, C.~Bonatti and S.~Crovisier \cite{MR2811152} (Lemma~\ref{ErgoLemPlus}), there is a perturbation of $f$ (still denoted by $f$) so that the orbit $\omega$ of $x$ is periodic of period $\tau_1$; moreover, $\omega$ can be supposed to bear an invariant measure close to $\mu$, to have an arbitrary large length, and to have Lyapunov exponents and Lyapunov subspaces close to that of $\mu$ under $f$. Applying the (difficult) connecting lemma for pseudo-orbits of C.~Bonatti and S.~Crovisier \cite{MR2090361}, we get another perturbation of the diffeomorphism (still denoted by $f$), such that the stable manifold of $x$ under $f$ meets the open set $U$ at a point that we denote by $y$.

So, we need to perturb the diffeomorphism $f$ so that:
\begin{itemize}
\item the periodic orbit $x$ is stabilized by $f_N$. This can be easily made by a small perturbation of $f$;
\item the positive orbit of $y$ under $f_N$ falls on the periodic orbit of $x$ under $f_N$. This is the difficult part of the proof: we can apply the previous strategy to put every point of the positive orbit of $y$ on the grid only during a finite time. It becomes impossible to perform perturbations to put the orbit of $y$ on the grid -- without perturbing the orbit of $x$ -- as soon as this orbit comes into a $C/N$- neighbourhood of the orbit of $x$ (where $C$ is a constant depending on $\mathcal V$).
\end{itemize}
To solve this problem, we need the results about the linear case we have proved in Chapter~\ref{Souris}. They allow us to find a point $z$ whose distance to the orbit of $x$ is bigger than $C/N$, and such that for $t$ large enough, $f_N^t(z)$ belongs to the orbit of $x$ under $f_N$.
\bigskip

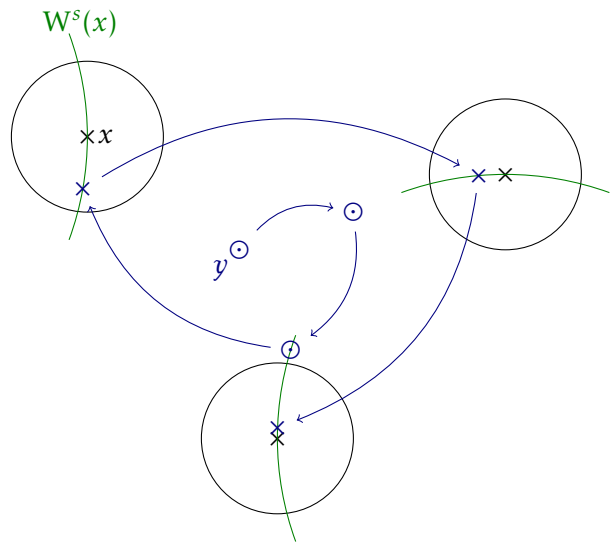
\begin{figure}[t]
\begin{center}
\begin{tikzpicture}[scale=1]

\draw (0,0) node{$\times$};
\draw (0,0) circle (1);
\draw[color=green!50!black] (0,0) arc (180:200:4);
\draw[color=green!50!black] (0,0) arc (180:160:4);

\draw (3,3.5) node{$\times$};
\draw (3,3.5) circle (1);
\draw[color=green!50!black] (3,3.5) arc (90:110:4);
\draw[color=green!50!black] (3,3.5) arc (90:70:4);

\draw (-2.5,4) node{$\times$};
\draw (-2.5,4) node[right]{$x$};
\draw (-2.5,4) circle (1);
\draw[color=green!50!black] (-2.5,4)  arc (0:20:4);
\draw[color=green!50!black] (-2.5,4)  arc (0:-20:4);
\draw[color=green!50!black] (-2.6,5.5)  node{$W^s(x)$};

\node[color=blue!50!black] (A) at (-.5,2.5) {$\odot$};
\node[color=blue!50!black] (B) at (1,3) {$\odot$};
\node[color=blue!50!black] (C) at (.175,1.17) {$\odot$};
\node[color=blue!50!black] (D) at (-2.561,3.31) {$\times$};
\node[color=blue!50!black] (E) at (2.651,3.485) {$\times$};
\node[color=blue!50!black] (F) at (0.002,0.140) {$\times$};

\draw[->,color=blue!50!black] (A) to[bend left] (B);
\draw[->,color=blue!50!black] (B) to[bend left] (C);
\draw[->,color=blue!50!black] (C) to[bend left] (D);
\draw[->,color=blue!50!black] (D) to[bend left] (E);
\draw[->,color=blue!50!black] (E) to[bend left] (F);

\draw[->,color=blue!50!black] (A) node[below left]{$y$};

\end{tikzpicture}
\caption[Perturbation of Lemma~\ref{LemMesPhysDiff}]{During the proof of Lemma~\ref{LemMesPhysDiff}, it is easy to perturb the first points of the orbit of $y$ (small disks) until the orbit meets the neighbourhoods of the orbits of $x$ where the diffeomorphism is linear (inside of the circles). The difficulty of the proof is to make appropriate perturbations in these small neighbourhoods.}\label{DessinGlobPhys}
\end{center}
\end{figure}

In more detail, we use Lemma~\ref{LemExtension} to linearize locally the diffeomorphism in the neighbourhood of the periodic orbit of $\omega$. In particular, the positive orbit of $y$ eventually belongs to this linearizing neighbourhood, from a time $T_1$. We denote $y' = f^{T_1}(y)$. To summarize, the periodic orbit $\omega$ bears a measure close to $\mu$, its Lyapunov exponents are close to that of $\mu$, its Lyapunov linear subspaces are close to that of $\mu$ (maybe not all along the periodic orbit, but at least for the first iterates of $x$). The diffeomorphism $f$ is linear around each point of $\omega$. Finally, the stable manifold of $\omega$ meets $U$ at $y$, and the positive orbit of $y$ is included in the neighbourhood of $\omega$ where $f$ is linear from the point $y' = f^{T_1}(y)$.

We then choose an integer $N$ large enough, and perturb the orbit of $x$ such that it is stabilized by the discretization $f_N$. We want to make another perturbation of $f$ such that the backward orbit of $x$ by $f_N$ also contains $y'$ (recall that $f_N$ is not necessarily one-to-one). This is done by a perturbation supported in the neighbourhood of
$\omega$ where $f$ is linear. First of all, during a time $t_4\ge 0$, we apply Lemma~\ref{DerTheoPart2} of Chapter~\ref{Souris} to find a point $z$ in the neighbourhood of $f^{-t_4}(x)$ where $f$ is linear, but far enough from $f^{-t_4}(x)$ compared to $1/N$, such that the $t_4$-th image of $z$ by the discretization $f_N$ is equal to $x$. Next, we perturb the orbit of $z$ under $f^{-1}$ during a time $t_3\ge 0$ such that $f^{-t_3}(z)$ belongs to the stable subspace of $f^{-t_4-t_3}(x)$. Note that the support of this perturbation must be disjoint from $\omega$; this is the reason why $z$ must be ``far enough from $x$''. Finally, we find another time $t_2$ such that the negative orbit $\{f^{-t}(z')\}_{t\ge 0}$ of $z' = f^{-t_3-t_2}(z)$ has an hyperbolic behaviour. We then perturb each point of the negative orbit of $z'$ (within the stable manifold of $\omega$), so that it contains an arbitrary point of the stable manifold of $\omega$, far enough from $\omega$. This allows us to meet the point $y'$, provided 
that the order of discretizations $N$ is large enough.

To complete the proof, we we consider the segment of $f$-orbit joining $y$ to $z$; we perturb each one of these points to put them on the grid $E_N$ (with a perturbation whose supports size is proportional to $1/N$).

Notice that we shall have chosen carefully the parameters of the first perturbations in order to make this final perturbation possible. Also, remark that the length of the periodic orbit $\omega$ must be very large compared to the times $t_2$, $t_3$ and $t_4$. This is why we will perform the proof in the opposite direction : we will begin by choosing the times $t_i$ and make the perturbation of the dynamics afterwards.

Note that the Addendum~\ref{AddTheoMesPhysDiff} can be proved by using a small variation on Lemma~\ref{LemMesPhysDiff}, that we will explain at the end of Section~\ref{AvDerSec}.

\subsection{An improved ergodic closing lemma}

The proof of Theorem~\ref{TheoMesPhysDiff} begins by the approximation of any invariant measure $\mu$ of any conservative $C^1$-diffeomorphism by a periodic measure of a diffeomorphism $g$ close to $f$. This is done by R.~Mañé's ergodic closing lemma, but we will need the fact that the obtained periodic measure inherits some of the properties of the measure $\mu$. More precisely, given a $C^1$-diffeomorphism $f$, we will have to approach any non periodic ergodic measure of $f$ with nonzero Lyapunov exponent by a periodic measure of a diffeomorphism $g$ close to $f$, such that the Lyapunov exponents and the Lyapunov subspaces of the measure are close to that of $f$ by $\mu$. We will obtain this result by modifying slightly the proof of a lemma obtained by F.~Abdenur, C.~Bonatti and S.~Crovisier in \cite{MR2811152} (Proposition~6.1).

\begin{lemme}[Ergodic closing lemma]\label{ErgoLemPlus}
Let $f\in\Diff^1(\T^n,\Leb)$. We consider
\begin{itemize}
\item a number $\varep>0$;
\item a $C^1$-neighbourhood $\mathcal V$ of $f$;
\item a time $\tau_0\in\N$;
\item an ergodic measure $\mu$ without zero Lyapunov exponent;
\item a point $x\in X$ which is typical for $\mu$ (see the beginning of the paragraph 6.1 of \cite{MR2811152});
\end{itemize}
moreover, we denote by $\lambda$ the smallest absolute value of the Lyapunov exponents of $\mu$, by $F^f_x$ the stable subspace at $x$ and by $G_x^f$ the unstable subspace\footnote{Stable and unstable in the sense of Oseledets splitting.} at $x$.
Then, there exists a diffeomorphism $g\in\Diff^1(\T^n,\Leb)$ and a time $\widetilde{t_0}>0$ (depending only in $f$, $\mu$ and $x$) such that:
\begin{enumerate}
\item $g\in\mathcal V$;
\item the point $x$ is periodic for $g$ of period $\tau\ge \tau_0$;
\item for any $t\le\tau$, we have $d\big(f^t(x),g^t(x)\big)<\varep$;
\item $x$ has no zero Lyapunov exponent for $g$ and the smallest absolute value of the Lyapunov exponents of $x$ is bigger than $\lambda/2$, we denote by $F^g_x$ the stable subspace and $G^g_x$ the unstable subspace;
\item the angles between $F^f_x$ and $F^g_x$, and between $G^f_x$ and $G^g_x$, are smaller than $\varep$;
\item for any $t\ge \widetilde{t_0}$, for any vectors of unit norm $v_F\in F^g_x$ and $v_G\in G^g_x$, we have
\[\frac{1}{t}\log \big(\|Dg^{-t}_x (v_F)\|\big)\ge \frac{\lambda}{4} \qquad \text{and} \qquad \frac{1}{t}\log \big(\|Dg^{t}_x (v_G)\|\big)\ge \frac{\lambda}{4}\]
\end{enumerate}
\end{lemme}

Remark that the proof of Proposition~6.1 of \cite{MR2811152} yields a similar lemma but with the weaker conclusion
\emph{\begin{itemize}
\item[5.] ``the angle between $G^f_x$ and $G^g_x$, is smaller than $\varep$''.
\end{itemize}}
\noindent Indeed, the authors obtain the linear space $G^g_x$ by a fixed point argument: Lemma~6.5 of \cite{MR2811152} states that the cone $C^s_{j,4C}$ is invariant by $Df_n^{-t_n}$, and thus contains both $G^f_x$ and $G^g_x$. Taking $C$ as big as desired, the cone $C^s_{j,4C}$ is as thin as desired and thus the angle between $G^f_x$ and $G^g_x$, is as small as desired. Unfortunately, in the original proof of Proposition~6.1 of \cite{MR2811152}, the linear space $F^g_x$ is not defined in the same way ; it is an invariant subspace which belongs to $C^u_{j,4C}$, which is an arbitrarily thick cone. Thus, the angle between $F^f_x$ and $F^g_x$, is not bounded by this method of proof. Our goal here is to modify the proof of Proposition~6.1 of \cite{MR2811152} to have simultaneously two thin cones $C'^u_{j,4C}$ and $C^s_{j,4C}$ which are invariant under respectively $Df_n^{t_n}$ and $Df_n^{-t_n}$

We begin by modifying the Lemma~6.2 of \cite{MR2811152}: we replace its forth point
\begin{itemize}
\item \emph{a sequence of linear isometries $P_n\in O_d(\R)$ such that $\|P_n - \Id\|< \varep$,}
\end{itemize}
by the point
\begin{itemize}
\item \emph{two sequences of linear isometries $P_n,Q_n\in O_d(\R)$ such that $\|P_n - \Id\|< \varep$ and $\|Q_n - \Id\|< \varep$,}
\end{itemize}
and its forth conclusion
\begin{itemize}
\item[\emph{d)}] \emph{For every $i\le j\in\{1,\cdots,k\}$ the inclination\footnote{The \emph{inclination} of a linear subspace $E\subset\R^n$ with respect to another subspace $E'\subset\R^n$ with the same dimension is the minimal norm of the linear maps $f : E\to E^\perp$ whose graph are equal to $E$.} of $Df_n^{t_n}. E_{i,j}$ with respect to $E_{i,j}$ is less than $C$.}
\end{itemize}
by the conclusion
\begin{itemize}
\item[\emph{d)}] \emph{For every $i\le j\in\{1,\cdots,k\}$ the inclination of $Df_n^{t_n}. E_{i,j}$ with respect to $E_{i,j}$ is less than $C$, and the inclination of $Df_n^{-t_n}. E_{i,j}$ with respect to $E_{i,j}$ is less than $C$.}
\end{itemize}

These replacements in the lemma are directly obtained by replacing Claim~6.4 of \cite{MR2811152} by the following lemma.

\begin{lemme}\label{Claim6.4}
For any $\eta>0$, there exists a constant $C>0$ such that for any matrix $A\in GL_n(\R)$ and any linear subspace $E\subset\R^n$, there exists two orthogonal matrices $P,Q\in O_n(\R)$ satisfying $\|P-\Id\|<\eta$ and $\|Q-\Id\|<\eta$, such that the inclinations of $(PAQ)(E)$ and $(PAQ)^{-1}(E)$ with respect to $E$ are smaller than $C$.
\end{lemme}


\begin{proof}[Proof of Lemma~\ref{Claim6.4}]
Given $\eta>0$, there exists a constant $C>0$ and a matrix $P_0\in O_n(\R)$ such that $\|P_0-\Id\|<\eta$, satisfying: for any linear subspace $E'\subset \R^n$, one of the two inclinations of $E'$ and of $P_0(E')$ with respect to $E$ is smaller then $C$.

We then choose an orthogonal matrix $Q\in O_n(\R)$ such that $\|Q-\Id\|<\eta$ and that (taking a bigger $C$ if necessary) both inclinations of $Q^{-1}\big(A^{-1}(E)\big)$ and $Q^{-1}\big((A^{-1} P_0^{-1})(E)\big)$ with respect to $E$ are smaller than $C$. There are two cases: either the inclination of $(AQ)(E)$ with respect to $E$ is smaller than $C$, and in this case we choose $P=\Id$, or the inclination of $(AQ)(E)$ with respect to $E$ is bigger than $C$, and in this case we can choose $P=P_0$. In both cases, the lemma is proved.
\end{proof}

The rest of the proof of Lemma~\ref{ErgoLemPlus} can be easily adapted from the proof of Proposition~6.1 of \cite{MR2811152}.

\subsection[Proof of the perturbation lemma]{Proof of the perturbation lemma (Lemma~\ref{LemMesPhysDiff})}\label{AvDerSec}

We now come to the proof of Lemma~\ref{LemMesPhysDiff}. We first do this proof in dimension 2, to simplify some arguments and to be able to make pictures.

\begin{proof}[Proof of Lemma~\ref{LemMesPhysDiff}]
Let $f$ be a conservative $C^1$-diffeomorphism, $\mathcal V$ a $C^1$-neighbourhood of $f$, $\varep>0$ and $N_0\in\N$. We denote $M = \max\big(\|Df\|_\infty, \|Df^{-1}\|_\infty\big)$. We also choose an $f$-invariant measure $\mu$ which is ergodic, not periodic and has no zero Lyapunov exponent, and an open set $U\subset \T^2$. We will make several successive approximations of $f$ in $\mathcal V$; during the proof we will need to decompose this neighbourhood: we choose $\delta>0$ such that the open $\delta$-interior $\mathcal V'$ of $\mathcal V$ is non-empty.

\paragraph{Step 0: elementary perturbation lemmas.}
During the proof of Lemma~\ref{LemMesPhysDiff}, we will use three different elementary perturbation lemmas.

The first one is the elementary perturbation lemma in $C^1$ topology we stated in Chapter~\ref{ChapPerturbLem} (Lemma~\ref{PerturbElem}). It allows to perturb locally the orbit of a diffeomorphism.

The second one is an easy corollary of the first one. We will use it to perturb a segment of orbit such that for any $N$ large enough, each point of this segment of orbit belongs to the grid $E_N$.

\begin{lemme}[Perturbation of a point such that it belongs to the grid]\label{Lem*}
For every open set $\mathcal V'$ of $\Diff^1(\T^n,\Leb)$, there exists $\eta'>0$ such that for $N$ large enough an every $x\in\T^n$, there exists $g\in\Diff^1(\T^n,\Leb)$ such that
\begin{itemize}
\item $g\in\mathcal V'$;
\item $g(x_N) = \big( f(x)\big)_N$;
\item $f=g$ outside of $B\big(x,(1+\eta')/N\big)$.
\end{itemize} 
\end{lemme}

Applying this lemma to several points $x_i\in \T^n$ which are far enough one from the others (for $i\neq j$, $d(x_i,x_j)\ge 2(1+\eta')/N$), it is possible to perturb $f$ into a diffeomorphism $g$ such that for every $i$, $g\big((x_i)_N\big) = \big( f(x_i)\big)_N$.

These two perturbations will be applied locally.
\bigskip

The third perturbation lemma is an improvement of Lemma~\ref{PerturbElem}; it states that the perturbation can be supposed to be a translation in a small neighbourhood of the perturbed point.

\begin{figure}[t]
\begin{center}
\begin{tikzpicture}[scale=1]

\fill[color=blue!10!white] (0,1.6) circle (.3);
\draw[color=blue!50!black] (0,1.6) circle (.3);
\draw[color=blue!50!black] (0,2.5) node{$B(0,r)$};

\draw[color=gray] (1,-1) arc (-90:90:1) -- (-1,1) arc (90:270:1) -- cycle;
\draw[color=gray] (1,-1.3) arc (-90:90:1.3) -- (-1,1.3) arc (90:270:1.3) -- cycle;
\draw[color=gray][color=gray] (1,-1.6) arc (-90:90:1.6) -- (-1,1.6) arc (90:270:1.6) -- cycle;
\draw[color=gray] (1,-1.9) arc (-90:90:1.9) -- (-1,1.9) arc (90:270:1.9) -- cycle;
\draw[color=gray] (1,-2.2) arc (-90:90:2.2) -- (-1,2.2) arc (90:270:2.2) -- cycle;
\draw[->,>=stealth,thick] (-.5,-1) -- (0,-1);
\draw[->,>=stealth,thick] (-.5,-1.3) -- (.5,-1.3);
\draw[->,>=stealth,thick] (-.5,-1.6) -- (.5,-1.6);
\draw[->,>=stealth,thick] (-.5,-1.9) -- (.5,-1.9);
\draw[->,>=stealth,thick] (-.5,-2.2) -- (0,-2.2);

\end{tikzpicture}
\caption[Proof of Lemma~\ref{LocTrans}]{Flow of the Hamiltonian used to prove Lemma~\ref{LocTrans} (``staduim'').}\label{FigStade}
\end{center}
\end{figure}
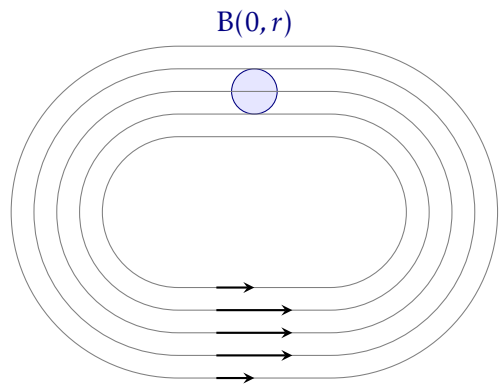

\begin{lemme}[Elementary perturbation with local translation]\label{LocTrans}
For every open set $\mathcal V'$ of $\Diff^1(\T^n,\Leb)$, and every $r>0$, there exists $N_1>0$ such that for every $N\ge N_1$ and every $\|v\|_\infty \le 1/(2N)$, there exists $g\in\Diff^1(\R^n,\Leb)$ such that:
\begin{itemize}
\item $g\in\mathcal V'$;
\item $\operatorname{Supp}(g)\subset B(0,10\,r)$;
\item for every $x\in B(0,r)$, $g(x) = x+v$.
\end{itemize}
\end{lemme}

\begin{proof}[Proof of Lemma~\ref{LocTrans}]
Take an appropriate Hamiltonian, see Figure~\ref{FigStade}.
\end{proof}

\paragraph{Step 1: choice of the starting point $x$ of the orbit.} Let $\lambda$ be the smallest absolute value of the Lyapunov exponents of $\mu$ (in particular, $\lambda > 0$).

We choose a point $x$ which is regular for the measure $\mu$: we suppose that it satisfies the conclusions of Oseledets and Birkhoff theorems, and Mañé's ergodic closing lemma (see Paragraph~6.1 of \cite{MR2811152}). We denote by $F_x^f$ the stable subspace and $G_x^f$ the unstable subspace for the Oseledets splitting at the point $x$. By Oseledets theorem, the growth of the angles $\angle\big(F_{f^i(x)}^f,G_{f^i(x)}^f\big)$ between the stable and unstable subspaces is subexponential (in both positive and negative times).

\paragraph{Step 2: choice of the parameters we use to apply the ergodic closing lemma.} In this second step, we determine the time during which we need an estimation of the angle between the stable and unstable subspaces of $f$ and its perturbations, and the minimal length of the approximating periodic orbit.

We first use the ``hyperbolic-like'' behaviour of $f$ near the orbit of $x$: for well chosen times $t_1$ and $t_2$, each vector which is not too close to $G_{f^{t_1}(x)}^f$ is mapped by $Df^{-t_2}$ into a vector which is close to $F_{f^{t_1-t_2}(x)}^f$. Given a vector $v\in T\T^n_{f^{t_1}(x)}$, it will allow us to perturb $f$ into $g$ such that an iterate of $v$ under $Dg^{-1}$ belongs to $F_{f^{t_1-t_2}(x)}^f$.

\begin{lemme}\label{LemAnglesOsel}
For every $\alpha>0$, there exists two times $t_1$ and $t_2 \ge 0$ such that if $v\in T \T^n_{f^{t_1}(x)}$ is such that the angle between $v$ and $G_{f^{t_1}(x)}$ is bigger than $\alpha$, then the angle between $Df^{-t_2}_{f^{t_1}(x)} v$ and $F_{f^{t_1 - t_2}(x)}$ is smaller than $\alpha$ (see Figure~\ref{FigRotaGlobC1}).
\end{lemme}

\begin{proof}[Proof of Lemma~\ref{LemAnglesOsel}]
It easily follows from Oseledets theorem, and more precisely from the fact that the function $\exp(t\lambda) / \angle(F_{f^t(x)},G_{f^t(x)})$ goes to $+\infty$ when $t$ goes to $+\infty$.
\end{proof}

So, we fix two times $t_1$ and $t_2 \ge 0$, obtained by applying Lemma~\ref{LemAnglesOsel} to $\alpha = \arcsin\big(1/(1+\eta)\big)$, where $\eta$ is the parameter obtained by applying the elementary perturbation lemma (Lemma~\ref{PerturbElem}) to $\delta/2$ (see Figure~\ref{FigRotaC1}).

We also choose a time $t_3\ge \widetilde t_0$ ($\widetilde t_0$ being given by Lemma~\ref{ErgoLemPlus}) such that
\begin{equation*}
e^{\lambda (t_3+t_2)/4} \ge M^{t_2}.
\end{equation*}
This estimation will be applied to point 6. of Lemma~\ref{ErgoLemPlus}. It will imply that for every $v\in F_{f^{t_1}(x)}^f$ and for every $t\ge t_2+t_3$, we have
\begin{equation}\label{Pastropgros}
\|Df_{f^{t_1}(x)}^{-t}(v)\| \ge \|Df_{f^{t_1}(x)}^{-t_2}(v)\| \ge \frac{1}{M^{t_2}} \|v\|.
\end{equation}
We then apply Lemma~\ref{DerTheoPart2} of Chapter~\ref{Souris} to
\begin{equation}\label{DefRr0}
R_0 = M^{t_2+t_3}(1+\eta'),
\end{equation}
where $\eta'$ is given by Lemma~\ref{Lem*} applied to the parameter $\delta/2$. This gives us a parameter $k_0 = t_4$. Note that $R_0$ is chosen so that if $v\in T_{f^{t_1}(x)}\T^n$ is such that $\|v\|\ge R_0/N$, then for any $t\in \llbracket 0, t_2+t_3\rrbracket$, we have
\begin{equation}\label{EqPasIdee}
\big\|Df^{-t}_{f^{t_1}(x)}(v)\big\| \ge (1+\eta')/N.
\end{equation}
Thus, we will be able to apply Lemma~\ref{Lem*} to the points $f^{-t}\big(f^{t_1}(x) + v\big)$, with $t\in \llbracket 0, t_2+t_3\rrbracket$, without perturbing the points of the orbit of $x$.

\paragraph{Step 3: global perturbation of the dynamics.} We can now apply the ergodic closing lemma we have stated in the previous subsection (Lemma~\ref{ErgoLemPlus}) to the neighbourhood $\mathcal V'$, the measure $\mu$, the point $x_1 = f^{t_1-t_2-t_3}(x)$ and $\tau_0 \ge t_2+t_3+t_4$ large enough so that $\tau_0 \lambda/4\ge 3$. We also need that the expansion of vectors $F^{g_1}$ along the segment of orbit $\big(x_2,g_2(x_2),\cdots,g_2^{\tau_0-t_2-t_3-t_4}(x_2)\big)$ is bigger than 3, but it can be supposed true by taking a bigger $\tau_0$ if necessary. This gives us a first perturbation $g_1$ of the diffeomorphism $f$, such that the point $x_1$ is periodic under $g_1$ with period $\tau_1\ge\tau_0$, and such that the Lyapunov exponents of $x_1$ for $g_1$ are close to that of $x_1$ under $f$, and the stable and unstable subspaces of $g_1$ at the point $g_1^t(x_1)$ are close to that of $f$ at the point $g_1^t(x_1)$ for every $t\in \llbracket 0,t_3+t_2\rrbracket$.

Remark that by the hypothesis on $\tau_0$, the Lyapunov exponent of $g_1^{\tau_1}$ at $x_1$ is bigger than 3, thus we will be able to apply Lemma~\ref{Lem*} to every point of the orbit belonging to $F^{g_1}_{x_1}$, even when the orbit returns several times near $x_1$. Also note that these properties are stable under $C^1$ perturbation.

We then use the connecting lemma for pseudo-orbits of C.~Bonatti and S.~Crovisier (Theorem~\ref{PapySylvain}, see \cite{MR2090361}), which implies that the stable manifolds of the periodic orbits of a generic conservative $C^1$-diffeomorphism are dense. This allows us to perturb the diffeomorphism $g_1$ into a diffeomorphism $g_2\in\mathcal V'$ such that there exists a point $x_2$ close to $x_1$ such that:
\begin{enumerate}[(1)]
\item $x_2$ is periodic for $g_2$ with the same period than that of $x_1$ under $g_1$, and moreover the periodic orbit of $x_2$ under $g_2$ shadows that of $x_1$ under $g_1$;
\item the Lyapunov exponents and the Lyapunov subspaces of $x_2$ for $g_2$ are very close to that of $x_1$ for $g_1$ (see the conclusions of Lemma~\ref{ErgoLemPlus}, in particular the Lyapunov subspaces are close during a time $t_3+t_2$);
\item the stable manifold of $x_2$ under $g_2$ meets the set $U$, at a point denoted by $y_2$.
\end{enumerate}

\paragraph{Step 4: linearization near the periodic orbit.} We then use Franks lemma (see \cite{MR0283812}) to perturb slightly the differentials of $g_2$ at the points $g_2^{t_2+ t_3}(x_2),$ $\cdots,g_2^{t_2+ t_3+t_4}(x_2)$, such that these differentials belong to the open set of matrices defined by Lemma~\ref{DerTheoPart2} of Chapter~\ref{Souris}. This gives us another diffeomorphism $g_3\in\mathcal V'$ close to $g_2$, such that the point $x_2$ still satisfies the nice properties (1), (2) and (3).

By Lemma~\ref{LemExtension}, there exists a parameter $r>0$ such that it is possible to linearize $g_3$ in the $r$-neighbourhood of the periodic orbit of $x_2$, without changing the nice properties (1), (2) and (3) of the periodic orbit of $x_2$. We can choose $r$ small enough so that the $10\,r$-neighbourhoods of the points of the periodic orbit of $x_2$ are pairwise disjoint. This gives us a diffeomorphism $g_4$, to which are associated two points $x_4$ and $y_4$, such that $x_4$ satisfies the properties (1), (2) and (3), and such that:
\begin{enumerate}[(1)]\setcounter{enumi}{3}
\item the differentials of $f$ at the points $g_4^{t_2+ t_3}(x_4),\cdots,g_4^{t_2+ t_3+t_4}(x_4)$ lie in the open dense set of matrices of Lemma~\ref{DerTheoPart2} of Chapter~\ref{Souris};
\item $g_4$ is linear in the $r$-neighbourhood of each point of the periodic orbit of $x_4$.
\end{enumerate}

\paragraph{Step 5: choice of the order of discretization.}

We choose a neighbourhood $\mathcal V''\subset\mathcal V'$ of $g_4$ such that properties (1) to (3) are still true for every diffeomorphism $g\in\mathcal V''$. We denote by $\omega_{x_4}$ the periodic orbit of $x_4$ under $g_4$, and by $B(\omega_{x_4},r)$ the $r$-neighbourhood of this periodic orbit. We also denote $T_1$ the smallest integer such that $g_4^t(y_4) \in B(\omega_{x_4},r/2)$ for every $t\ge T_1$, and set $y'_4 = g_4^{T_1}(y_4)$. Thus, the positive orbit of $y'_4$ will stay forever in the linearizing neighbourhood of $\omega_{x_4}$. Taking $T_1$ bigger if necessary, we can suppose that $y'_4$ belongs to the linearizing neighbourhood of the point $x_4$. We can also suppose that for every $t\in \llbracket 0, \tau_1\rrbracket$,
\begin{equation}\label{LoinOrbPer}
3d\big(g_4^{T_1-t}(y_4),g_4^{-t}(x_4)\big) \le \min_{\tau_1\le t'\le T_1}d\big(g_4^{T_1-t'}(y_4),g_4^{-t}(x_4)\big).
\end{equation}

We can now choose the order $N$ of the discretization, such that
\begin{enumerate}[(i)]
\item $N\ge N_0$ ($N_0$ has been chosen at the very beginning of the proof);
\item \label{itemlinea} applying Lemma~\ref{LocTrans} to the parameter $r$ and the neighbourhood $\mathcal V''$ to get an integer $N_1$, we have $N\ge N_1$, so that it is possible to choose the value of the points of $\omega_{x_4}$ modulo $E_N$ without changing the properties (1) to (5);
\item the distance between two distinct points of the segment of orbit $y_4,$ $g_4(y_4), \cdots, g_4^{T_1}(y_4) = y'_4$ is bigger than $2(1+\eta')/N + 2/N$, so that it will be possible to apply Lemma~\ref{Lem*} simultaneously to each of these points, even after the perturbation made during the point \eqref{itemlinea}, such that these points belong to $E_N$;
\item every $\sqrt 2/N$-pseudo-orbit\footnote{The constant $\sqrt n/N$ comes from the fact that an orbit of the discretization is a $\sqrt 2/N$-pseudo-orbit.} starting at a point of the periodic orbit $\omega_{x_4}$ stays during a time $T'\tau_1$ in the $d(y'_4,\omega_{x_4})$-neighbourhood of the periodic orbit, where $T'$ the smallest integer such that
\begin{equation}\label{hyphyphyp}
 \left(1 + \frac{1}{3(1+\eta)}\right)^{T'} \ge \nu,
\end{equation}
and $\nu$ is the maximal modulus of the eigenvalues of $(Dg_4)_{x_4}^{\tau_1}$. A simple calculus shows that this condition is true if for example
\begin{equation*}
N \ge \frac{2\sqrt n (M^{T'\tau_1}-1)}{r(M-1)}.
\end{equation*}
This condition will be used to apply the process described by Lemma~\ref{LemLyapPerturb}.
\end{enumerate}

\paragraph{Step 6: application of the linear theorem.}

By the hypothesis (ii) on $N$, we are able to use Lemma~\ref{LocTrans} (elementary perturbation with local translation) to perturb each point of the periodic orbit $\omega_{x_4}$ such that we obtain a diffeomorphism $g_5\in\mathcal V''$ and points $x_5$, $y_5$ and $y'_5$ satisfying properties (1) to (5) and moreover:
\begin{enumerate}[(1)]\setcounter{enumi}{5}
\item for every $t\in \llbracket t_2 + t_3, t_2 + t_3 + t_4\rrbracket$, the value of $g_5^t(x_5)$ modulo $E_N$ is equal to $w_k/N$, where $w_k$ is given by Lemma~\ref{DerTheoPart2} of Chapter~\ref{Souris};
\item for any other $t$, $g_5^t(x_5)$ belongs to $E_N$.
\end{enumerate}
In particular, the periodic orbit of $x_5$ under $g_5$ is stabilized by the discretization $(g_5)_N$ (indeed, recall that $w_k \in [-1/2,1/2]^k$).

By construction of the diffeomorphism $g_5$ (more precisely, the hypotheses (4), (5), (6) and (7)), it satisfies the conclusions of Lemma~\ref{DerTheoPart2} of Chapter~\ref{Souris}; thus there exists a point $z\in B\big(g_5^{t_2+t_3}(x_5),r\big)$ such that $(g_5)_N^{t_4} (z) = (g_5)_N^{t_2+t_3+t_4} (x_5)$ and that $\|z - g_5^{t_2+t_3}(x_5)\|\ge R_0/N$ (where $R_0$ is defined by Equation~\eqref{DefRr0}). Remark that hypothesis (iv) implies that $\|z-g_5^{t_2+t_3}(x_5)\|\ll r$.

\paragraph{Step 7: perturbations in the linear world.}

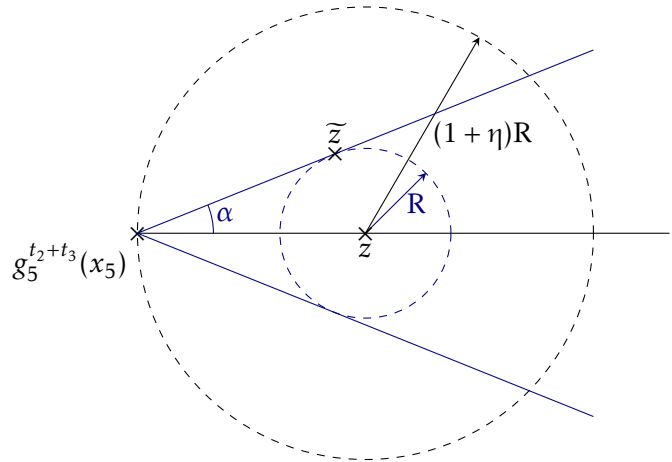
\begin{figure}[t]
\begin{center}
\begin{tikzpicture}[scale=1]
\draw (0,0) node {$\times$};
\draw (0,0) node[below left] {$g_5^{t_2+t_3}(x_5)$};

\draw (0,0) -- (7,0);
\draw[dashed] (3,0) circle (3);
\draw[->,>=stealth] (3,0) -- (4.5,2.6) node[pos=.5,right]{$(1+\eta)R$};
\draw[dashed, color=blue!50!black] (3,0) circle (9/8);
\draw[->,>=stealth, color=blue!50!black] (3,0) -- (3.8,0.8) node[pos=.5,right]{$R$};

\draw[color=blue!50!black] (0,0) -- (6,2.43);
\draw[color=blue!50!black] (0,0) -- (6,-2.43);
\draw[color=blue!50!black] (1,0) arc (0:22.05:1);
\draw[color=blue!50!black] (1.15,.25) node{$\alpha$};

\draw (2.6,1.053) node {$\times$};
\draw (2.6,1.053) node[above] {$\widetilde z$};

\draw (3,0) node {$\times$};
\draw (3,0) node[below] {$z$};
\end{tikzpicture}
\caption[Perturbation made to apply Lemma~\ref{LemAnglesOsel}]{Perturbation we make to apply Lemma~\ref{LemAnglesOsel} (see also Figure~\ref{FigRotaGlobC1}): we make an elementary perturbation in a neighbourhood of $z$ mapping $z$ into $\widetilde z$, such that the angle between the lines $\big(g_5^{t_2+t_3}(x_5)\ z\big)$ and $\big(g_5^{t_2+t_3}(x_5)\ \widetilde z\big)$ is bigger than $\alpha = \arcsin\big(1/(1+\eta)\big)$, and such that the support of the perturbation does not contain $g_5^{t_2+t_3}(x_5)$.}\label{FigRotaC1}
\end{center}
\end{figure}

In this step, our aim is to perturb the negative orbit of $z$ under $g_5$ such that it meets the point $y'_5$. Remark that by hypothesis (iv), every point of $z,g_5^{-1}(z),\cdots,g_5^{-t_2}(z)$ is in the linearizing neighbourhood of $\omega_{x_5}$.

From now, all the perturbations we will make will be local, and we will only care of the positions of a finite number of points. Thus, it will not be a problem if these perturbations make hypotheses (3) and (5) become false, provided that they have a suitable behaviour on this finite set of points.

First, if necessary, we make a perturbation in the way of Figure~\ref{FigRotaC1}, so that the angle between the lines $\big( g_5^{t_2+t_3}(x_5)\ z\big)$ and $G_{g_5^{t_2+t_3}(x_5)}^{g_5}$ is bigger than $\alpha$; this gives us a diffeomorphism $g_6$. More precisely, the support of the perturbation we apply is contained in a ball centred at $z$ and with radius $d(z,x_6)$, so that this perturbation does not change the orbit of $x_6$. Under these conditions, we satisfy the hypotheses of Lemma~\ref{LemAnglesOsel}, thus the angle between $\big( g_6^{t_3}(x_6)\ g_6^{-t_2}(z)\big)$ and $F_{{g_6}^{t_3}(x_6)}^{g_6}$ is smaller than $\alpha$. Another perturbation, described by Figure~\ref{FigRotaC1}, allows us to suppose that $g_6^{-t_2}(z)$ belongs to $F_{{g_6}^{t_3}(x_6)}^{g_6}$. This gives us a diffeomorphism that we still denote by $g_6$. Remark that it was possible to make these perturbations independently because the segment of negative orbit of the point $z$ we considered does not enter twice in the 
neighbourhood of a point of $\omega_{x_6}$ where the diffeomorphism is linear.

Thus, the points $z' = g_6^{-t_2}(z)$ and $y'_6 = y'_5$ both belong to the local stable manifold of the point $x_6 = x_5$ for $g_6$ (which coincides with the Oseledets linear subspace $F_{x_6}^{g_6}$ since $g_6$ is linear near $x_6$).
\bigskip

The next perturbation takes place in the neighbourhood of the point $x_6$ (and not in all the linearizing neighbourhoods of the points of $\omega_{x_6}$).

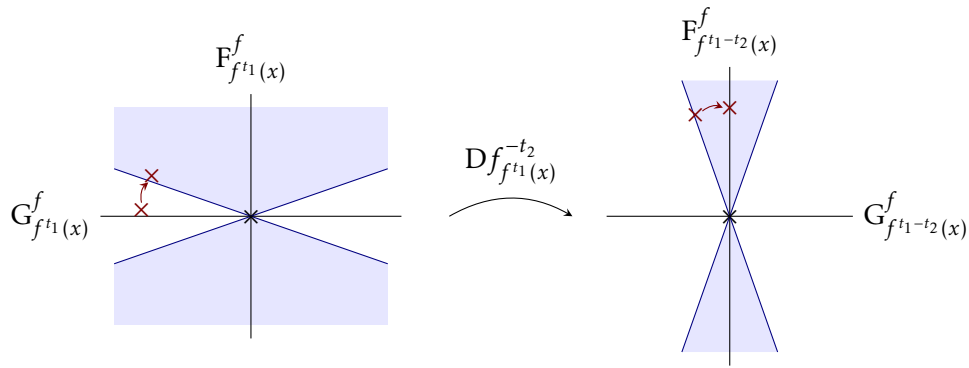
\begin{figure}[t]
\begin{center}
\begin{tikzpicture}[scale=.9]

\fill[color=blue!10!white] (-6,.7) -- (-2,-.7)  -- (-6,-.7) -- (-2,.7) -- cycle;
\fill[color=blue!10!white] (-6,.7) -- (-2,.7)  -- (-2,1.6) -- (-6,1.6) -- cycle;
\fill[color=blue!10!white] (-6,-.7) -- (-2,-.7)  -- (-2,-1.6) -- (-6,-1.6) -- cycle;
\draw (-4,0) node {$\times$};
\draw[color=blue!50!black] (-6,-.7) -- (-2,.7);
\draw[color=blue!50!black] (-6,.7) -- (-2,-.7);
\draw (-1.8,0) -- (-6.2,0);
\draw (-4,-1.8) -- (-4,1.8);
\draw (-4,1.8) node[above]{$F^f_{f^{t_1}(x)}$};
\draw (-6.2,0) node[left]{$G^f_{f^{t_1}(x)}$};
\draw[color=red!50!black] (-5.6,.1) node{$\times$};
\draw[color=red!50!black] (-5.45,.6) node{$\times$};
\draw[->,>=stealth, color=red!50!black] (-5.6,.2) to[bend left] (-5.5,.5);

\draw[->,>=stealth] (-1.1,0) to[bend left] (.7,0);
\draw (-.2,.8) node{$Df^{-t_2}_{f^{t_1}(x)}$};

\fill[color=blue!10!white] (2.3,-2) -- (3.7,2)  -- (2.3,2) -- (3.7,-2) -- cycle;
\draw (3,0) node {$\times$};
\draw[color=blue!50!black] (2.3,-2) -- (3.7,2);
\draw[color=blue!50!black] (2.3,2) -- (3.7,-2);
\draw (1.2,0) -- (4.8,0);
\draw (3,-2.2) -- (3,2.2);
\draw (3,2.2) node[above]{$F^f_{f^{t_1-t_2}(x)}$};
\draw (4.8,0) node[right]{$G^f_{f^{t_1-t_2}(x)}$};
\draw[color=red!50!black] (2.5,1.5) node{$\times$};
\draw[color=red!50!black] (3,1.6) node{$\times$};
\draw[->,>=stealth, color=red!50!black] (2.6,1.55) to[bend left] (2.9,1.6);

\end{tikzpicture}
\caption[Proof of Lemma~\ref{LemAnglesOsel}]{Proof of Lemma~\ref{LemAnglesOsel}: make a small perturbation at times $t_1$ and $t_1-t_2$ (in red), the hyperbolic-like behaviour of $f$ does the rest of the work for you. In red: the perturbation that we will make during step 7.}\label{FigRotaGlobC1}
\end{center}
\end{figure}

\begin{lemme}\label{LemLyapPerturb}
For every $y'\in F_{x_6}^{g_6}$ such that $d(y',x_6) > d(z',x_6)\nu^{T'\tau_1}$ ($T'$ being defined by Equation~\eqref{hyphyphyp}), there exists a diffeomorphism $g_7$ close to $g_6$ and $T''\in\N$ such that $g_7^{-\tau_1 T''}(z') = y'$. Moreover, the perturbations made to obtain $g_7$ are contained in the linearizing neighbourhood of $\omega_{x_6}$, do not modify the images of $\omega_{x_6}$, nor these of the negative orbit of $z'$ by the discretization or these of the positive orbit of $y'$ in the linearizing neighbourhood of $\omega_{x_6}$
\end{lemme}

\begin{figure}[t]
\begin{center}
\begin{tikzpicture}[scale=1]

\draw[->,>=stealth] (0,0) -- (10.5,0);
\draw (10.5,0) node[right]{$F_{x_6}^{g_6}$};
\draw (.5,0) node {$|$};
\draw (.5,-.5) node {$x_6$};

\draw (1.5,0) node {$\times$};
\draw (1.5,-.5) node {$z'$};
\draw[->,>=stealth,color=blue!50!black] (1.7,-.1) to[bend right] (2.8,-.1);
\draw[color=blue!50!black] (2.25,-.7) node {$g_6^{-\tau_1}$};

\draw[color=blue!50!black] (3,0) node {$\times$};
\draw[dotted] (3,0) circle (1.2);
\draw[dashed] (3,0) circle (.6);
\draw[->,>=stealth,color=blue!50!black] (3.2,-.1) to[bend right] (5.7,-.1);
\draw[color=blue!50!black] (4.45,-.9) node {$g_6^{-\tau_1}$};

\draw[color=blue!50!black] (5.9,0) node {$\times$};

\draw[->,>=stealth,color=red!50!black] (1.7,.1) to[bend left] (3.4,.1);
\draw[->,>=stealth,color=red!50!black,dotted] (1.7,.1) to[bend left] (2.8,.1);
\draw[->,>=stealth,color=red!50!black] (3.8,.1) to[bend left] (8.2,.1);
\draw[->,>=stealth,color=red!50!black,dotted] (3.8,.1) to[bend left] (7,.1);
\draw[color=red!50!black] (3.6,0) node {$\times$};
\draw[color=red!50!black] (8.4,0) node {$\times$};
\draw[dotted] (7.2,0) circle (2.4);
\draw[dashed] (7.2,0) circle (1.2);

\draw[color=red!50!black] (8.4,0) node {$\times$};

\end{tikzpicture}
\caption[Perturbation made to merge orbits]{Perturbation such that the point $y'_6$ belongs to the negative orbit of $z'$: the initial orbit is drawn in blue (below) and the perturbed orbit in red (above). From a certain time, the red orbit overtakes the blue orbit.}\label{FigLyapPerturb}
\end{center}
\end{figure}
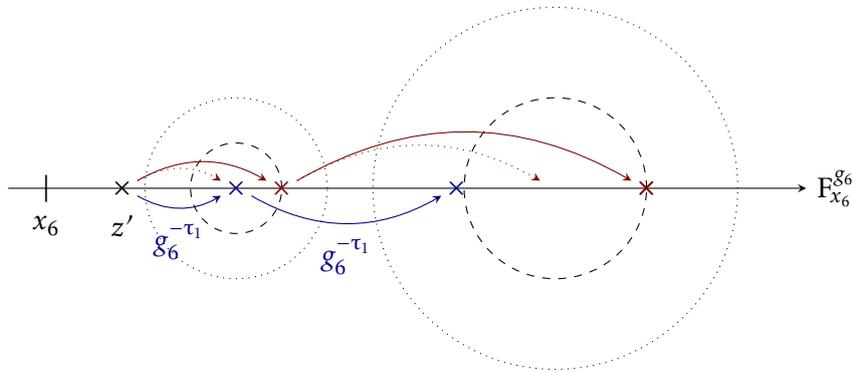

\begin{proof}[Proof of Lemma~\ref{LemLyapPerturb}]
During this proof, if $r$ and $s$ are two points of $W^s(x_6)$, we will denote by $[r,s]$ the segment of $W^s(x_6)$ between $r$ and $s$. Remark that if $r$ and $s$ lie in the neighbourhood of $x_6$ where $g_6$ is linear, then $[r,s]$ is a real segment, included in $F_{x_6}^{g_6}$ moreover, we will denore $[r,+\infty[$ the connected component of $W^s(x_6)\setminus\{r\}$ which does not contain $x_6$.

Consider the point $z'\in F_{x_6}^{g_6}$, and choose a point 
\[p\in \left[g_6^{-\tau_1}(z'), \left(1+\frac{1}{2(1+\eta)}\right)g_6^{-\tau_1}(z')\right].\]
By applying an elementary perturbation (Lemma~\ref{PerturbElem})  whose support is contained into $B\big(g_6^{-\tau_1}(z'),\, d(x_6,g_6^{-\tau_1}(z'))/2\big)$, it is possible to perturb $g_6$ into a diffeomorphism $g_7$ such that $g_6^{\tau_1}(z') = p$ (see Figure~\ref{FigLyapPerturb}). Applying this process $t$ times, for every
\[p\in \left[g_6^{-\tau_1t}(z'), \left(1+\frac{1}{2(1+\eta)}\right)^{t}g_6^{-\tau_1t}(z')\right],\]
it is possible to perturb $g_6$ into a diffeomorphism $g_7$ such that $g_7^{-\tau_1 t}(z') = p$ (the supports of the perturbations are disjoint because the expansion of $g^{-\tau_1}_{|F_{x_6}^{g_6}}$ is bigger than 3). But as $T'$ satisfies Equation~\eqref{hyphyphyp}, the union
\[\bigcup_{t\ge 0}\left[g_6^{-\tau_1t}(z'), \left(1+\frac{1}{2(1+\eta)}\right)^{t}g_6^{-\tau_1t}(z')\right]\]
covers all the interval $[g_6^{-\tau_1T'}(z'), +\infty[$. By the hypothesis made on $y'$, we also have $y'\in [g_6^{-\tau_1T'}(z'), +\infty[$; this proves the lemma.
\end{proof}
Thus, by hypothesis (iv), it is possible to apply Lemma~\ref{LemLyapPerturb} to our setting. This gives us a diffeomorphism $g_7$.

\paragraph{Step 8: final perturbation to put the segment of orbit on the grid.}

To summarize, we have a diffeomorphism $g_7\in\mathcal V'$, and periodic orbit $\omega_{x_7}$ of $g_7$, stabilized by $(g_7)_N$, which bears a measure close to $\mu$. We also have a segment of real orbit of $g_7$ which links the points $y_7\in U$ and $z$, where $z$ is such that $(g_7)_N^{t_4}(z) \in (\omega_{x_7})_N$. To finish the proof of the lemma, it remains to perturb $g_7$ so that the segment of orbit which links the points $y_7$ and $z$ is stabilized by the discretization $(g_7)_N$.
\bigskip

We now observe that by the construction we have made, the distance between two different points of the segment of orbit under $g_7$ between $y_7$ and $z$ is bigger than $2(1+\eta')/N$, and the distance between one point of this segment of orbit and a points of $\omega_{x_7}$ is bigger than $(1+\eta')/N$.

Indeed, if we take one point of the segment of forward orbit $z, g_7^{-1} (z),\cdots, g_7^{-t_2-t_3} (z)$, and one point in the periodic orbit $\omega_{x_7}$, this is due to the hypothesis $\|z - x_7\|\ge R_0/N$ ($R_0$ being defined by Equation~\eqref{DefRr0}) combined with Equation~\eqref{EqPasIdee}. If we take one point in this segment $z, g_7^{-1} (z),\cdots, g_7^{-t_2-t_3} (z)$, and one among the rest of the points (that is, the segment of orbit between $y_7$ and $z$), this is due to the fact that the Lyapunov exponent of $g_7^{\tau_1}$ in $x_7$ is bigger than 3, and to Equation~\eqref{LoinOrbPer}.

If we take one point of the form $g_7^{-t}(z)$, with $t> t_2+t_3$, but belonging to the neighbourhood of $\omega_{x_7}$ where $g_7$ is linear, and one point of $\omega_{x_7}$, this follows from the estimation given by Equation~\eqref{Pastropgros} applied to $\|v\|\ge R_0/N$. If for the second point, instead of considering a point of $\omega_{x_7}$, we take an element of the segment of orbit between $y_7$ and $z$, this follows from the fact that the Lyapunov exponent of $g_7^{\tau_1}$ in $x_7$ is bigger than 3.

Finally, for the points of the orbit that are not in the neighbourhood of $\omega_{x_7}$ where $g_7$ is linear, the property arises from hypothesis (iii) made on $N$.

Thus, by Lemma~\ref{Lem*}, we are able to perturb each of the points of the segment of orbit under $g_7$ between $y_7$ and $z$, such that each of these points belongs to the grid. This gives us a diffeomorphism $g_8\in \mathcal V$. 

To conclude, we have a point $y_8\in U$ whose orbit under $(g_8)_N$ falls on the periodic orbit $(\omega_{x_7})_N$, which bears a measure $\varep$-close to $\mu$. The lemma is proved.
\bigskip

The proof in higher dimensions is almost identical. The perturbation lemmas are still true\footnote{In particular, Lemma~\ref{LocTrans} can be obtained by considering a plane $(P)$ containing both $x$ and $y$ and taking a foliation of $\R^n$ by planes parallel to $(P)$. The desired diffeomorphism is then defined on each leave by the time-$\psi(t)$ of the Hamiltonian given in the proof of the lemma, with $\psi$ is a smooth compactly supported map on the space $\R^n/(P)$, equal to 1 in $0$ and with small $C^1$ norm.}, and the arguments easily adapts by considering the ``super-stable'' manifold of the orbit $\omega_x$, that is the set of points $y\in\T^n$ whose positive orbit is tangent to the Oseledets subspace corresponding to the maximal Lyapunov exponent. In particular, Lemma~\ref{ErgoLemPlus} is still true in this setting, and the connecting lemma (Theorem~\ref{PapySylvain}) implies that generically, this ``super-stable'' manifold is dense in $\T^n$.
\end{proof}

The proofs of the two statements of the addendum are almost identical.

For the first statement (the fact that for every $\varep>0$, the basin of attraction of the discrete measure can be supposed to contain a $\varep$-dense subset of the torus), we apply exactly the same proof than that of Lemma~\ref{LemMesPhysDiff}: making smaller perturbations of the diffeomorphism if necessary, we can suppose that the stable manifold of $y_8$ is $\varep$-dense. Thus, there exists a segment of backward orbit of $y_8$ which is $\varep$-dense, and we apply the same strategy of proof consisting in putting this segment of orbit on the grid. 

For the second statement, it suffices to apply the strategy of the first statement, and to conjugate the obtained diffeomorphism $g_9$ by an appropriate conservative diffeomorphism with small norm (this norm can be supposed to be as small as desired by taking $\varep$ small), so that the image of the $\varep$-dense subset of $\T^2$ by the conjugation contains the set $E$.
\bigskip

To obtain Theorem~\ref{TheoMesPhysDiffDissip} (dealing with the dissipative case) it suffices to replace the use of Lemmas~\ref{LemAnglesOsel} and~\ref{DerTheoPart2} by the following easier statement.

\begin{lemme}\label{RemplaceDissip}
For every $\alpha>0$ and every $R_0>0$, there exists three times $t_1$, $t_2 \ge 0$ and $t_4\ge 0$ such that:
\begin{itemize}
\item there exists $v\in T \T^n_{f^{t_1}(x)}\cap \Z^n$ such that $\|v\|\ge R_0$ and
\[\big(\widehat{Df}_{f^{t_1 + t_4}(x)}\circ \cdots \circ \widehat{Df}_{f^{t_1}(x)}\big)(v) = 0;\]
\item if $v\in T \T^n_{f^{t_1}(x)}$ is such that the angle between $v$ and $G_{f^{t_1}(x)}$ is bigger than $\alpha$, then the angle between $Df^{-t_2}_{f^{t_1}(x)} v$ and $F_{f^{t_1 - t_2}(x)}$ is smaller than $\alpha$ (see Figure~\ref{FigRotaGlobC1}).
\end{itemize}
\end{lemme}

\begin{proof}[Proof of Lemma~\ref{RemplaceDissip}]
This comes from Oseledets theorem and the hypotheses made on the Lyapunov exponents of $x$, and in particular that their sum is strictly negative.
\end{proof}

\section{Transfer operators and physical measures}\label{SecTransOp}

The aim of this section is to study the behaviour of the measures $(f_N^*)^m \lambda_N$ for ``small'' times $m$ (recall that $\lambda_N$ is the uniform measure on the grid $E_N$). We will focus on the case where $f$ is a $C^{1+\alpha}$ expanding map of the circle: a classical theorem asserts that in this case, the map $f$ has a single physical measure (which is also the SRB measure), that we will denote by $\mu_0$ (see Theorem~\ref{Liverani}). In this setting, O.E.~Lanford has stated the following conjecture.

\begin{conj}[Lanford]\label{Lalan}
Let $f : \Sp^1\to\Sp^1$ be a generic $C^2$ expanding map. Then the convergence $(f_N^*)^m(\lambda_N) \to \mu_0$ holds for $N,m$ both going to $+\infty$, with $\log N \ll m \ll \sqrt N$.
\end{conj}

In this section, we will prove a kind of weak version of this conjecture (Theorem~\ref{MainMoche}). We begin by some notations.

\begin{definition}
For $0<\alpha\le 1$ and $d\ge 2$, the set $\mathcal E_d^{1+\alpha}(\Sp^1)$\index{$\mathcal E_d^{1+\alpha}(\Sp^1)$} is the set of maps $f : \Sp^1\to \Sp^1$ of degree $d$, whose derivative $f'$ belongs to $C^\alpha(\Sp^1)$ and satisfies $f'(x)>1$ for every $x\in\Sp^1$.
\end{definition}

Here, we will only consider parameters $\alpha\le 1$. In particular, the set $\mathcal E_d^2$ will denote the set of maps $f : \Sp^1\to \Sp^1$ of degree $d$, whose derivative $f'$ is Lipschitz and satisfies $f'(x)>1$ for every $x\in\Sp^1$.

The main result of this section is the following.

\begin{theoreme}\label{MainMoche}
For every $0<\alpha\le 1$ and every $C^{1+\alpha}$ expanding map $f\in \mathcal E_d^{1+\alpha}(\Sp^1)$, there exists a constant $c_0 = c_0(f)>0$ such that if $(N_m)_m$ is a sequence of integers going to infinity and satisfying $\log N_m > c_0 m$, then the convergence $(f_{N_m}^*)^m(\lambda_{N_m}) \to \mu_0$ holds.
\end{theoreme}

This theorem will be obtained by combining the classical result of convergence of the pushforwards of any smooth measure towards the SRB measure $\mu_0$ (Theorem~\ref{Liverani}), together with the convergence of the discrete operators $f_N^*$ to the Ruelle-Perron-Frobenius operator (Theorem~\ref{Prop25}). Both of these theorems are effective: the constant $c_0$ of Theorem~\ref{MainMoche} can be computed in practical.
\bigskip

\emph{We fix once for all an expanding map $f\in \mathcal E_d^{1+\alpha}(\Sp^1)$, with $0<\alpha\le 1$. In the sequel, we will freely identify a measure with its density.}
\bigskip

First of all, we recall an effective result about the convergence of the pushforwards of any smooth measure towards the SRB measure.

The \emph{transfer operator} associated to the map $f$ (usually called Ruelle-Perron-Fro\-be\-nius operator), which acts on densities of probability measures, will be denoted by $\Ll_f$\index{$\Ll_f$}. It is defined by
\[\Ll_f \phi(y) = \sum_{x\in f^{-1}(y)} \frac{\phi(x)}{f'(x)}.\]

Then, for every map $\phi\in C^\alpha(\Sp^1)$ which is the density of a probability measure on $\Sp^1$, the sequence $\Ll^m_f \phi$ converges exponentially fast towards a map $\phi_0\in C^\alpha(\Sp^1)$; this map is the density of the unique $f$-invariant probability measure $\mu_0$ having a $C^\alpha$ density. This measure is also the unique SRB measure of $f$. More precisely, we have Theorem~\ref{Liverani}. To state it, we will need the following notation.

\begin{definition}
For $0<\alpha\le 1$ and $g\in C^\alpha(\Sp^1)$, we denote by $[g]_\alpha$\index{ $[g]_\alpha$} the Hölder norm of $g$, i.e.
\[[g]_\alpha = \sup_{x\neq y} \frac{|g(x)-g(y)|}{|x-y|^\alpha}\]
\end{definition}

\begin{theoreme}[\cite{MR2504311}, see also Theorem~2.5 of \cite{MR1343323} for the $C^2$ case]\label{Liverani}
Let $\phi\in C^\alpha(\Sp^1)$ be a positive function whose integral is equal to 1. Then, there exists an integer $M = M(\phi) > 0$ such that for every $m\ge M$,
\[\| \Ll_f^m(\phi) - \phi_0\|_\infty \le C_0 \Lambda^{m-M} \|\phi\|_\infty,\]
where $\phi_0$ denotes the density of the measure $\mu_0$,
\[C_0 = e^\Delta \|\phi_0\|_\infty \left(\|\phi_0\|_\infty + \frac{2}{\min_{x\in\Sp^1}\phi_0(x)} \right),\]
and $\Lambda = \tanh(\Delta/4)$, with
\[\Delta \le 2\log\left(\frac{1+\tau}{1-\tau}\right) + 2\big(3-\tau+2\|f'\|_\infty\big) \frac{2[\log f']_\alpha}{\lambda^{\alpha}-1} \|f'\|_\infty^{-\alpha},\]
where $\lambda = \min_{x\in \Sp^1}f'(x)$, and $\tau = (1+\lambda^{-\alpha})/2$.

Moreover, if $\phi$ is the constant function equal to 1, then $M(\phi) = 0$.
\end{theoreme}
This theorem is obtained by taking $K = 2 K_0 /(\lambda^\alpha-1)$ in the Theorem 10.1 of \cite{MR2504311} (where $K$ is defined at the beginning of section 10 of \cite{MR2504311}).
\bigskip

We now come to the original part of the proof of Theorem~\ref{MainMoche}, and compare the actions of the transfer operator and of its discrete counterpart on smooth densities of measures. The discrete operator will be denoted by $f_N^*$\index{$f_N^*$}; it acts on the probability measures supported by $E_N$:
\[(f_N^*\mu)(y) = \sum_{x\in f_N^{-1}(y)} \mu(x).\]
We will show that in a certain sense, we have
\[f_N^* \underset{N\to+\infty}{\longrightarrow} \Ll_f.\]
Thus, for a fixed time $m\in\N$, the measures $(f_N^*)^m(\lambda_N)$ tend to the measure $\mu_0$ when $N$ goes to $+\infty$. The problem is to have precise estimates about the speed of this convergence. This speed will be measured by the following distance on measures.

\begin{definition}
For $\mu$ and $\nu$ two probability Borel measures, we define $d_{\Lip}(\mu,\nu)$\index{$d_{\Lip}$} the distance 
\[d_{\Lip}(\mu,\nu) = \sup\left\{\left| \int_{\Sp^1} \psi\, \ud (\mu-\nu)\right| \mid \psi\in\Lip_1(\Sp^1) \right\},\]
where $\Lip_1(\Sp^1)$\index{$\Lip_1(\Sp^1)$} is the space of Lipschitz functions whose Lipschitz constant is smaller than 1.
\end{definition}

This distance $d_{\Lip}$ spans the weak-* topology on $\Prb$. Remark that by a theorem of L.~Kantorovich and G.~Rubinstein, this distance coincides with Wasserstein distance $W_1$ (see \cite{MR0102006}). Note that the distance between Lebesgue measure and the uniform measure $\lambda_N$ on $E_N$ is $d_{\Lip}(\Leb,\lambda_N) = 1/(4N)$.

\begin{theoreme}\label{Prop25}
Let $m\in\N$. To every map $\phi\in C^\alpha(\Sp^1)$ which is the density of a probability measure, we associate the measure $\varphi_N$\index{$\varphi_N$} supported by $E_N$: $\varphi_N(A) = \int_{P_N^{-1}(A)} \phi$ (where $P_N$ is the projection on $E_N$). Then, for every map $f\in \mathcal E_d^{1+\alpha}(\Sp^1)$,
\[{f_N^*}^m (\varphi_N) \underset{N\to+\infty}{\longrightarrow} \Ll^m_f(\phi).\]
More precisely, we have
\[d_{\Lip}\Big	((f_N^*)^m(\varphi_N) , \Ll_f^m\phi\Big)\le \eta\]
as soon as
\[N > (4/\eta)^{1+1/\alpha}\big(2 + \|f'\|_\infty \big)^{m(1+1/\alpha)}\left(d \big(1 + [f']_\alpha\big) \left([\phi]_\alpha + \frac{1}{d}\right)\right)^{m/\alpha}.\]
\end{theoreme}

\begin{coro}
For example, we have (recall that $\lambda_N$ is the uniform measure on $E_N$)
\[d_{\Lip}\Big	((f_N^*)^m(\lambda_N) , \Ll_f^m\Leb\Big)\le \eta\]
as soon as
\[N> \left(\frac{4}{\eta}\right)^{1+1/\alpha}\left(\big(2 + \|f'\|_\infty \big)^{1+1/\alpha}\big(1 + [f']_\alpha\big) ^{1/\alpha}\right)^m.\]
In particular, if $f$ is $C^2$, this gives
\[d_{\Lip}\Big	((f_N^*)^m(\lambda_N) , \Ll_f^m\Leb\Big)\le \eta\]
as soon as
\[N> \frac{16}{\eta^2}\left(\big(2 + \|f'\|_\infty \big)^{2}\big(1 + [f']_\alpha\big)\right)^m.\]
\end{coro}

Notice that this theorem is only a partial answer to the Conjecture~\ref{Lalan} of O.E.~Lanford. In particular, it does not answer to the following question.

\begin{ques}
Let $f : \Sp^1\to\Sp^1$ be a generic $C^{1+\alpha}$ expanding map. What is the best function $m(N)$ such that the convergence $(f_N^*)^m(\lambda_N) \to \mu_0$ holds for $N,m$ both going to $+\infty$, with $\log N \ll m \ll m(N)$?
\end{ques}

It is reasonable to think (see Conjecture~\ref{Lalan}) that $m(N)\le \sqrt N$. Indeed, if $\sigma : E\to E$ is a typical random map of a set $E$ with $q$ elements, and $x$ is a typical point of $E$, then the smallest integer $m$ such that $\sigma^m(x)\in \Omega(\sigma)$ (recall that $\Omega(\sigma)$ is the union of the periodic orbits of $\sigma$) is of order $\sqrt q$ (see \cite[XIV.5]{Boll-rand} or the Theorem~2.3.1 of \cite{Mier-dyna}).

Notice that in Theorem~\ref{Prop25}, there is no hypothesis of genericity. We can hope that the uniform distribution of the roundoff errors of a generic map (see Proposition~\ref{LinPourLinstant}) allows to have finer results about the best function $m(N)$.

We do not know what happens in the $C^1$ generic case. It is possible that the behaviour is very different from the $C^{1+\alpha}$ case. Indeed, a theorem of J.~Campbell and A.~Quas (see \cite{MR1845327}) asserts that a generic $C^1$ expanding map of the circle has a single physical measure, but that this measure is singular with respect to Lebesgue measure.

We could also imagine similar studies in more general settings. For example, we can expect that the results on the circle generalize to $C^{1+\alpha}$ expanding maps of the torus $\T^n$, with $n\ge 1$. It could be less straightforward to have generalizations to Anosov diffeomorphisms of $\T^n$: in this case, the transfer operator is quite different from that used in the case of expanding maps.
\bigskip

We now begin the proof of Theorem~\ref{Prop25} by the easier case where the time $m$ is equal to $1$.

\begin{lemme}\label{LemDistTps1}
For $\varep<1/5$ and
\[N>\varep^{-1-1/\alpha}\max\left\{1,\ [\phi]_\alpha^{1/\alpha},\ {[f']_\alpha^{1/\alpha}}\right\},\]
we have
\[d_{\Lip}\big(f_N^*(\varphi_N),\,\Ll_f\phi\big) \le \varep \big(5+2\|f'\|_\infty\big).\]
\end{lemme}

\begin{proof}[Proof of Lemma~\ref{LemDistTps1}]
Let $\phi$ be the $C^\alpha(\Sp^1)$ density of a probability measure on $\Sp^1$, i.e. $\phi\ge 0$ and $\int_{\Sp^1} \phi = 1$. Let $\varep>0$ and $f\in \mathcal E_d^{1+\alpha}$. As $f\in C^{1+\alpha}(\Sp^1)$, we can estimate its lack of linearity on an interval depending on its length. This will allow us to estimate the difference between $\Ll_f$ and $f_N$ on every interval based on the length of this interval. The global difference between $\Ll_f$ and $f_N$ will be obtained by a summation on a partition of $\Sp^1$ of intervals of appropriate length.

\begin{lemme}\label{Titigro}
Let $I$ be an interval of $\Sp^1$ with length smaller than
\[\left(\frac{\varep}{[f']_\alpha}\right)^{1/\alpha},\]
and $x_0\in I$. Then
\begin{equation*}
\left|\frac{\card(I\cap E_N)}{\card(E_N)} - \frac{\card(f(I)\cap E_N)}{f'(x_0)\card(E_N)}\right| \le \Leb(I)\left(\varep + \frac{2}{N}\frac{f'(x_0)}{\Leb(f(I))} \frac{1-\varep}{1-2\varep} \right).
\end{equation*}
\end{lemme}

\begin{proof}[Proof of Lemma~\ref{Titigro}]
The hypothesis on the length of $I$ implies that if $x,y\in I$, then $|f'(x)-f'(y)|<\varep$. Using the mean value inequality, we get
\[\left| \Leb(I) - \frac{\Leb(f(I))}{f'(x_0)}\right| \le \varep;\]
we deduce the bound
\begin{equation}\label{ESTD}
\left| \Leb(I) - \frac{\Leb(f(I))}{f'(x_0)}\right| \le \frac{\varep\Leb(I)}{f'(x_0)}.
\end{equation}
This implies that
\[\left| \Leb(I) - \frac{\Leb(f(I))}{f'(x_0)}\right| \le \frac{\varep/f'(x_0)}{1-\varep/f'(x_0)}\Leb(f(I)),\]
thus (because $f'(x_0)\ge 1$)
\[\left| \Leb(I) - \frac{\Leb(f(I))}{f'(x_0)}\right| \le \frac{\varep}{f'(x_0)(1-\varep)}\Leb(f(I)) ;\]
in particular,
\begin{equation}\label{TAF}
\Leb(I) \ge \frac{\Leb(f(I))}{f'(x_0)}\,\frac{1-2\varep}{1-\varep}.
\end{equation}

Moreover, for every interval $J$, 
\[\left|\Leb(J) - \frac{\card(J\cap E_N)}{\card(E_N)}\right| \le \frac{1}{N},\]
as a result, using Equation~\eqref{ESTD},
\[\left|\frac{\card(I\cap E_N)}{\card(E_N)} - \frac{\card(f(I)\cap E_N)}{f'(x_0)\card(E_N)}\right| \le \frac{1}{N} + \frac{\varep\Leb(I)}{f'(x_0)} + \frac{1}{Nf'(x_0)}.\]
We deduce that (still because $f'(x_0)\ge 1$)
\[\left|\frac{\card(I\cap E_N)}{\card(E_N)} - \frac{\card(f(I)\cap E_N)}{f'(x_0)\card(E_N)}\right| \le \Leb(I)\left(\varep + \frac{2}{N \Leb(I)} \right),\]
using Equation~\eqref{TAF}, this leads to
\[\left|\frac{\card(I\cap E_N)}{\card(E_N)} - \frac{\card(f(I)\cap E_N)}{f'(x_0)\card(E_N)}\right| \le \Leb(I)\left(\varep + \frac{2}{N}\frac{f'(x_0)}{\Leb(f(I))} \frac{1-\varep}{1-2\varep} \right).\]
\end{proof}

Let $\psi\in \Lip_1(\Sp^1)$ be a test function. We want to compute the difference
\[\big|\langle f_N^*(\varphi_N),\psi\rangle - \langle\Ll\phi^*,\psi\rangle \big|.\]
So, we compute:
\begin{align*}
\big|\langle f_N^*(\varphi_N),\psi\rangle - & \langle\Ll_f\phi,\psi\rangle\big|\\
                                          = & \left|\sum_{y_N\in E_N} \int_{P_N^{-1}(y_N)}\psi(y)\big(f_N^*(\varphi_N) - \Ll_f\phi\big)(y) \ud y\right|\\
                                          = & \left|\sum_{y_N\in E_N} \left(\psi(y_N)\!\!\sum_{x_N\in f_N^{-1}(y_N)}\!\!\varphi_N(x_N) - \int_{P_N^{-1}(y_N)}\!\psi(y)\! \sum_{x\in f^{-1}(y)}\frac{\phi(x)}{f'(x)} \ud y\right)\right|.
\end{align*}
We remark that by preservation of the total mass (the operators $f_N^*$ and $\Ll$ map probability measures to probability measures), this expression is independent from the mean of $\psi$. We can therefore freely add or remove a constant to the function $\psi$.

We introduce an intermediate scale between that of the grid $E_N$ and that of the circle $\Sp^1$: the grid $E_M$. When $M$ is large enough, we have:
\begin{enumerate}
\item[(C1)] $|\psi(y)-\psi(y')|<\varep$ for every points $y,y'$ such that $|y-y'|<1/M$,
\item[(C2)] $|f'(x)-f'(x')|<\varep$ et $|\phi(x)-\phi(x')|<\varep$ for every points $x,x'$ in the same connected component of $f^{-1}\big(P_M^{-1}(y_M)\big)$,
\item[(C3)] $1/M \le (\varep/[f']_\alpha)^{1/\alpha}$.
\item[(C4)] $M/N \le \varep$.
\end{enumerate}

In particular, for every points $x,x'$ in the same connected component of the preimage $f^{-1}(P_M^{-1}(y_M))$, we have
\begin{equation}\label{estphi}
\left|\varphi_N(x_N) - \frac{1}{N}\phi(x')\right| \le \frac{\varep}{N},
\end{equation}
and
\begin{align}
\left|\frac{\phi(x)}{f'(x)} - \frac{\phi(x')}{f'(x')}\right| & \le \phi(x)\left|\frac{1}{f'(x)}-\frac{1}{f'(x')}\right| + \frac{1}{f'(x')}\left|\phi(x)-\phi(x')\right|\nonumber\\
                       & \le \frac{\phi(x)}{f'(x)}\left|\frac{f'(x')-f'(x)}{f'(x')}\right| + \frac{\varep}{f'(x)}\frac{f'(x)}{f'(x')}\nonumber\\
											 & \le \frac{\phi(x)}{f'(x)}\frac{\varep}{\min f'} + \frac{\varep}{f'(x)}\left(1+\frac{\varep}{\min f'}\right)\nonumber\\
											 & \le \varep\left(\frac{\phi(x)}{f'(x)} + \frac{1+\varep}{f'(x)}\right)\label{estquot}
\end{align}
(we want to have a $f'(x)$ at the denominator to be able to integrate properly and keep working with $L^1$ norms instead of $L^\infty$ norms).

We begin by cutting out the distance we want to compute by using the intermediate scale $E_M$:
\begin{align*}
\big|\langle f_N^*(\varphi_N),\psi\rangle - & \langle\Ll_f\phi,\psi\rangle\big|\\
																				\le & \sum_{y'_M\in E_M} \left|\sum_{y_N\in E_N\cap P_M^{-1}(y'_M)}\! \left(\psi(y_N)\!\!\!\!\!\!\!\sum_{x_N\in f_N^{-1}(y_N)}\!\!\!\!\!\!\!\varphi_N(x_N) - \!\!\!\!\!\!\!\int\limits_{P_N^{-1}(y_N)}\!\!\!\!\!\psi(y)\!\!\! \sum_{x\in f^{-1}(y)}\frac{\phi(x)}{f'(x)} \ud y\right)\right|.
\end{align*}
Using condition (C1) and the fact that $\langle f_N^*(\varphi_N),1\rangle = \langle\Ll\phi^*,1\rangle = 1$, we get:
\begin{align*}
\big|\langle f_N^*(\varphi_N),\psi\rangle - & \langle\Ll_f\phi,\psi\rangle\big|\\
																				\le & 2\varep \\
																				    & + \|\psi\|_\infty\sum_{y'_M\in E_M} \left|\sum_{y_N\in E_N\cap P_M^{-1}(y'_M)} \left(\sum_{x_N\in f_N^{-1}(y_N)}\!\!\!\!\!\!\!\!\varphi_N(x_N) - \!\!\!\!\!\!\!\int\limits_{P_N^{-1}(y_N)}\! \sum_{x\in f^{-1}(y)}\frac{\phi(x)}{f'(x)} \ud y\right)\right|.\\
\end{align*}
For $y'_M\in E_M$, we denote by $(I_{y'_M,k})_{1\le k\le d}$ the connected components of $f^{-1}\big(P_M^{-1}(y'_M)\big)$. For every $k$, we also denote by  $x_{y'_M,k}$ the unique point of $I_{y'_M,k}\cap f^{-1}(y'_M)$. Using the bound~\eqref{estphi}, we get:
\begin{align*}
\big|\langle f_N^*(\varphi_N) & ,\psi\rangle - \langle\Ll_f\phi,\psi\rangle\big|\\
                                        \le & \varep \left(2 + \|\psi\|_\infty\right)\\
																					  & + \|\psi\|_\infty\!\!\!\!\sum_{y'_M\in E_M} \left| \sum_{k=1}^d \card(I_{y'_M,k}\cap E_N) \frac{\phi(x_{y'_M,k})}{N} - \!\!\!\!\!\!\!\!\!\!\!\!\sum_{y_N\in E_N\cap P_M^{-1}(y'_M)}\ \int\limits_{P_N^{-1}(y_N)} \sum_{x\in f^{-1}(y)}\frac{\phi(x)}{f'(x)} \ud y\right|.
\end{align*}
Combined with a change of variables, the bound~\eqref{estquot} leads to:
\begin{align*}
\big|\langle f_N^*(\varphi_N) & ,\psi\rangle - \langle\Ll_f\phi,\psi\rangle\big|\\
                                        \le & \varep \Big(2 + \|\psi\|_\infty\big(1+1+(1+\varep)\big)\Big)\\
																					  & + \|\psi\|_\infty\!\!\!\sum_{y'_M\in E_M} \left|\sum_{k=1}^d \card(I_{y'_M,k}\cap E_N)\frac{\phi(x_{y'_M,k})}{N} - \!\!\!\!\!\!\!\!\!\!\!\!\sum_{y_N\in E_N\cap P_M^{-1}(y'_M)}\! \frac{1}{N}\sum_{k=1}^d\frac{\phi(x_{y'_M,k})}{f'(x_{y'_M,k})} \right|\\
																				\le & \varep \big(2 + \|\psi\|_\infty (3+\varep)\big)\\
																					  & + \|\psi\|_\infty\sum_{y'_M\in E_M} \left| \sum_{k=1}^d \frac{\phi(x_{y'_M,k})}{N} \left(\card(I_{y'_M,k}\cap E_N) - \frac{\card(P_M^{-1}(y'_M))}{f'(x_{y'_M,k})}\right)\right|\\
																				\le &  \varep \big(2 + \|\psi\|_\infty (3+\varep)\big)\\
																			      & + \|\psi\|_\infty\sum_{y'_M\in E_M} \sum_{k=1}^d \phi(x_{y'_M,k}) \left|\frac{\card(I_{y'_M,k}\cap E_N)}{N} - \frac{\card(f(I_{y'_M,k})\cap E_N)}{f'(x_{y'_M,k}) N} \right|.
\end{align*}
Using Lemma~\ref{Titigro} (which is valid by Condition~(C3)), we deduce that:
\begin{align*}
\big|\langle f_N^*(\varphi_N),\psi\rangle - & \langle\Ll_f\phi,\psi\rangle\big|\\
																				\le &  \varep \big(2 + \|\psi\|_\infty (3+\varep)\big)\\
																			      & + \|\psi\|_\infty\sum_{y'_M\in E_M} \sum_{k=1}^d \phi(x_{y'_M,k}) \Leb(I_{y'_M,k})\left(\varep + \frac{2}{N}\frac{f'(x_{y'_M,k})}{\Leb(f(I_{y'_M,k}))} \frac{1-\varep}{1-2\varep} \right).
\end{align*}
As a result, using Condition~(C2) and the fact that $\|\phi\|_{L^1} = 1$, we get
\begin{align*}
\big|\langle f_N^*(\varphi_N),\psi\rangle - & \langle\Ll_f\phi,\psi\rangle\big|\\
																				\le &  \varep \big(2 + \|\psi\|_\infty (3+\varep)\big)\\
																				    & + \|\psi\|_\infty(1+\varep)\left( \varep + 2\frac{M}{N} \| f'\|_\infty \frac{1-\varep}{1-2\varep}\right)\\
\end{align*}
finally, condition~(C4) gives
\begin{align*}
\big|\langle f_N^*(\varphi_N),\psi\rangle - & \langle\Ll_f\phi,\psi\rangle\big|\\
																				\le &  \varep \left(2 + \|\psi\|_\infty \left(4+2\varep + 2 \| f'\|_\infty \frac{1+\varep}{1-2\varep}\right)\right).
\end{align*}

We now use the fact that $\psi\in\Lip_1(\Sp^1)$. In this case (recall that we are allowed to add any constant to $\psi$), we can suppose that $\|\psi\|_\infty\le 1/2$. We have a bound on the minimal order $M$ to have Conditions~(C1) to (C3) :
\[\frac{1}{M}<\min\left\{\varep,\ \left(\frac{\varep}{[\phi]_\alpha}\right)^{1/\alpha},\ \left(\frac{\varep}{[f']_\alpha}\right)^{1/\alpha}\right\}.\]
So, for $\varep \le 1$, Condition~(C4) gives
\[\frac{1}{N}<\varep^{1+1/\alpha}\min\left\{1,\ \frac{1}{[\phi]_\alpha^{1/\alpha}},\ \frac{1}{[f']_\alpha^{1/\alpha}}\right\},\]
and in this case,
\[d_{\Lip}\big(f_N^*(\varphi_N),\,\Ll_f\phi\big) \le \varep \left(4+\varep + \| f'\|_\infty \frac{1+\varep}{1-2\varep}\right).\]
If moreover, we suppose that $\varep<1/5$, we get the conclusion of the lemma:
\[d_{\Lip}\big(f_N^*(\varphi_N),\,\Ll_f\phi\big) \le \varep \big(5 + 2\| f'\|_\infty \big).\]
\end{proof}

We now focus on the general case $m>1$. It follows from two easy lemmas.

\begin{lemme}\label{LemFacilIter}
Let $\varphi_1,\varphi_2 : E_N \to \R_+$. Then, for every $j\in\N$, we have
\[d_{\Lip}\big( (f_N^*)^j\varphi^1 , (f_N^*)^j \varphi^2\big) \le \big( 2+\|f'\|_\infty \big)^j d_{\Lip}( \varphi^1 , \varphi^2).\]
\end{lemme}

\begin{proof}[Proof of Lemma~\ref{LemFacilIter}]
For $\psi\in\Lip_1(\Sp^1)$, we compute
\begin{align*}
\big \langle f_N^*\varphi^1 - f_N^* \varphi^2 , \psi \big\rangle = & \sum_{y_N\in E_N}\big( f_N^*\varphi^1(y_N) - f_N^*\varphi^2(y_N) \big) \psi(y_N)\\
    = & \sum_{y_N\in E_N} \sum_{x_N\in f_N^{-1}(y)}\big( \varphi^1(x_N) - \varphi^2(x_N) \big) \psi(y_N).
\end{align*}
But, for $x_N, x'_N \in E_N$, we have (remark that either $x_N=x'_N$, or $|x_N-x'_N|\ge 1/N$) 
\[\big|f_N(x_N) - f_N(x'_N)\big| \le \big(2 + \|f'\|_\infty\big) |x_N-x'_N|.\]
As a result, the map $x_N \mapsto \psi(f_N(x_N))$ is $\big(2 + \|f'\|_\infty\big)$-Lipschitz, so
\[d_{\Lip}\big( f_N^*\varphi^1 , f_N^* \varphi^2\big) \le \big( 2+\|f'\|_\infty \big) d_{\Lip}( \varphi^1 , \varphi^2).\]
The lemma easily follows from an induction.
\end{proof}

\begin{lemme}\label{LemFacilOpe}
\[\big[\Ll^j\phi\big]_\alpha \le d^j \big(1 + [f']_\alpha\big)^j \left([\phi]_\alpha + \frac{1}{d}\right).\]
\end{lemme}

\begin{proof}[Proof of lemma~\ref{LemFacilOpe}]
We have
\[\big|(\Ll\phi)(y) - (\Ll\phi)(y')\big| \le \left|\sum_{x\in f^{-1}(y)} \frac{\phi(x)}{f'(x)} - \sum_{x'\in f^{-1}(y')} \frac{\phi(x')}{f'(x')} \right|.\]
We write $f^{-1}(y) = \{x_1,\cdots,x_d\}$ and $f^{-1}(y') = \{x'_1,\cdots,x'_d\}$, so that for every $k$, we have $|x_k - x'_k|\le |y-y'|$. Thus
\begin{align*}
\big|(\Ll\phi)(y) - (\Ll\phi)(y')\big| & \le \sum_{k=1}^d \left|\frac{\phi(x_k)}{f'(x_k)} - \frac{\phi(x'_k)}{f'(x'_k)} \right|\\
         & \le \sum_{k=1}^d |\phi(x_k)| \left|\frac{1}{f'(x_k)} - \frac{1}{f'(x'_k)} \right| +  \frac{1}{f'(x'_k)}\left|\phi(x_k) - \phi(x'_k) \right|.
\end{align*}
Using the fact that $\min_{x\in\Sp^1} f'(x) \ge 1$, this leads to
\begin{align*}
\big|(\Ll\phi)(y) - (\Ll\phi)(y')\big| & \le \sum_{k=1}^d |\phi(x_k)| \left|f'(x_k) - f'(x'_k) \right| +  \left|\phi(x_k) - \phi(x'_k) \right|\\
          & \le \sum_{k=1}^d \|\phi\|_\infty [f']_\alpha |x_k - x'_k|^\alpha +  [\phi]_\alpha |x_k - x'_k|^\alpha\\
					& \le d \big(\|\phi\|_\infty [f']_\alpha + [\phi]_\alpha \big) |y-y'|^\alpha.
\end{align*}
Moreover, as $\phi\ge 0$ and as there exists $x\in\Sp^1$ such that $\phi(x)\le 1$, we have $\|\phi\|_\infty \le 1+[\phi]_\alpha$. So
\[\big[\Ll\phi\big]_\alpha \le d \Big([\phi]_\alpha \big(1 + [f']_\alpha\big) + [f']_\alpha\Big).\]
By iterating, we get
\begin{align*}
\big[\Ll^j\phi\big]_\alpha & \le d^j \big(1 + [f']_\alpha\big)^j \left([\phi]_\alpha + \frac{[f']_\alpha}{d \big(1 + [f']_\alpha\big) - 1}\right)\\
      & \le d^j \big(1 + [f']_\alpha\big)^j \left([\phi]_\alpha + \frac{1}{d}\right).
\end{align*}
\end{proof}

\begin{proof}[Proof of Theorem~\ref{Prop25}]
To begin with, we decompose the distance between $(f_N^*)^m(\varphi_N)$ and $\Ll_f^m\phi$:
\[d_{\Lip}\Big	((f_N^*)^m(\varphi_N) , \Ll_f^m\phi\Big) \le \sum_{j=1}^m d_{\Lip}\left((f_N^*)^j\big(\Ll_f^{m-j}(\phi)\big)_N, (f_N^*)^{j-1}\big(\Ll_f^{m+1-j}(\phi)\big)_N\right).\]
We then use Lemma~\ref{LemFacilIter}:
\begin{align*}
d_{\Lip}\Big	((f_N^*)^m(\varphi_N) & , \Ll_f^m\phi\Big)\\
              \le & \sum_{j=1}^m \big(2+\|f'\|_\infty \big)^{j-1}d_{\Lip}\left(f_N^*\big(\Ll_f^{m-j}(\phi)\big)_N, \Ll_f\big(\Ll_f^{m-j}(\phi)_N\big)\right).
\end{align*}
Applying Lemma~\ref{LemDistTps1} to $\varep<1/5$ and $N$ such that for every $j$,
\begin{equation}\label{Suspens}
N>\varep^{-1-1/\alpha}\max\left\{1,\ \big[\Ll_f^{m-j}(\phi)\big]_\alpha^{1/\alpha},\ {[f']_\alpha^{1/\alpha}}\right\},
\end{equation}
we deduce that
\[d_{\Lip}\big(f_N^*(\varphi_N),\,\Ll_f\phi\big) \le \varep \big(5+2\|f'\|_\infty\big).\]
\begin{align*}
d_{\Lip}\Big	((f_N^*)^m(\varphi_N) , \Ll_f^m\phi\Big) \le & \sum_{j=1}^m \big(2+\|f'\|_\infty \big)^{j-1} \big(5+2\|f'\|_\infty\big)\\
							\le & \frac{\big(2+\|f'\|_\infty \big)^{m}-1}{\big(2+\|f'\|_\infty \big)-1} \big(5+2\|f'\|_\infty\big).
\end{align*}
As $\|f'\|_\infty \ge 2$, we get
\[d_{\Lip}\Big	((f_N^*)^m(\varphi_N) , \Ll_f^m\phi^*\Big)\le 3\varep \big(2 + \|f'\|_\infty \big)^{m}.\]

By Lemma~\ref{LemFacilOpe}, \eqref{Suspens} is satisfied if
\[N>\varep^{-1-1/\alpha}\max\left\{1,\ \left(d \big(1 + [f']_\alpha\big) \left([\phi]_\alpha + \frac{1}{d}\right)\right)^{m/\alpha},\ {[f']_\alpha^{1/\alpha}}\right\},\]
i.e.
\[N> \varep^{-1-1/\alpha}\left(d \big(1 + [f']_\alpha\big) \left([\phi]_\alpha + \frac{1}{d}\right)\right)^{m/\alpha}.\]
Thus, for every $\eta\le 2$, we have
\[d_{\Lip}\Big	((f_N^*)^m(\varphi_N) , \Ll_f^m\phi^*\Big)\le \eta\]
as soon as
\[N> \eta^{-1-1/\alpha} 4^{1+1/\alpha}\big(2 + \|f'\|_\infty \big)^{m(1+1/\alpha)}\left(d \big(1 + [f']_\alpha\big) \left([\phi]_\alpha + \frac{1}{d}\right)\right)^{m/\alpha}.\]
\end{proof}

\section{Numerical simulations}\label{NumSimPhys}

In this section, we present the results of the numerical simulations we have conducted in connection with Theorems~\ref{TheoMesPhysDiff} and \ref{MainMoche}.

\subsection[Simulations of the measures $\mu^{f_N}_{x}$]{Simulations of the measures $\mu^{f_N}_{x}$ for conservative torus diffeomorphisms}

We have computed numerically the measures $\mu^{f_N}_x$ for conservative diffeomorphisms $f\in\Diff^1(\T^2,\Leb)$, for the uniform grids
\[E_N = \left\{\left(\frac{i}{N},\frac{j}{N}\right)\in \T^2 \big|\  0\le i,j\le {N}-1\right\},\]
and for starting points $x$ either equal to $(1/2,1/2)$, or chosen at random. We present images of sizes $128\times 128$ pixels representing in logarithmic scale the density of the measures $\mu^{f_N}_x$: each pixel is coloured according to the measure carried by the set of points of $E_N$ it covers. Blue corresponds to a pixel with very small measure and red to a pixel with very high measure. Scales on the right of each image corresponds to the measure of one pixel on the $\log 10$ scale: if green corresponds to $-3$, then a green pixel will have measure $10^{-3}$ for $\mu^{f_N}_x$. For information, when Lebesgue measure is represented, all the pixels have a value about $-4.2$.

We have carried out the simulations on three different diffeomorphisms.\label{DefDiffeoPhys}
\begin{itemize}
\item The first conservative diffeomorphism $h_1$ is of the form $h_1= Q\circ P$, where both $P$ and $Q$ are homeomorphisms of the torus that modify only one coordinate:
\[P(x,y) = \big(x,y+p(x)\big)\quad\text{and}\quad Q(x,y) = \big(x+q(y),y\big),\]
with
\[p(x) = \frac{1}{209}\cos(2\pi\times 17x)+\frac{1}{471}\sin(2\pi\times 29x)-\frac{1}{703}\cos(2\pi\times 39x),\]
\[q(y) = \frac{1}{287}\cos(2\pi\times 15y)+\frac{1}{403}\sin(2\pi\times 31y)-\frac{1}{841}\sin(2\pi\times 41y).\]
This $C^\infty$-diffeomorphism is in fact $C^1$-close to the identity. This allows $h_1$ to admit periodic orbits with not too large periods. Note that $h_1$ is also chosen so that it is not $C^2$-close to the identity. 

\item The second conservative diffeomorphism $h_2$ is the composition $h_2 = h_1\circ R$, with the translation of the torus
\[R(x,y) = \big(x+1/10 , y+1/15\big).\]
Again, for $R$, we have chosen a translation with a relatively small order (here 30) to ensure that the discretizations can have periodic orbits with small periods.
\item The third conservative diffeomorphism $h_3$ is the composition $h_2 = h_1\circ A$, with $A$ the linear Anosov map
\[A = \begin{pmatrix} 2 & 1 \\ 1 & 1 \end{pmatrix}.\]
As $h_1$ is $C^1$-close to $\Id$, the diffeomorphism $h_3$ is $C^0$-conjugated to the linear automorphism $A$, which is in particular ergodic.
\end{itemize}
\bigskip

\begin{figure}[h]
\begin{minipage}[c]{.33\linewidth}
	\includegraphics[width=\linewidth,trim = .5cm .3cm .6cm .1cm,clip]{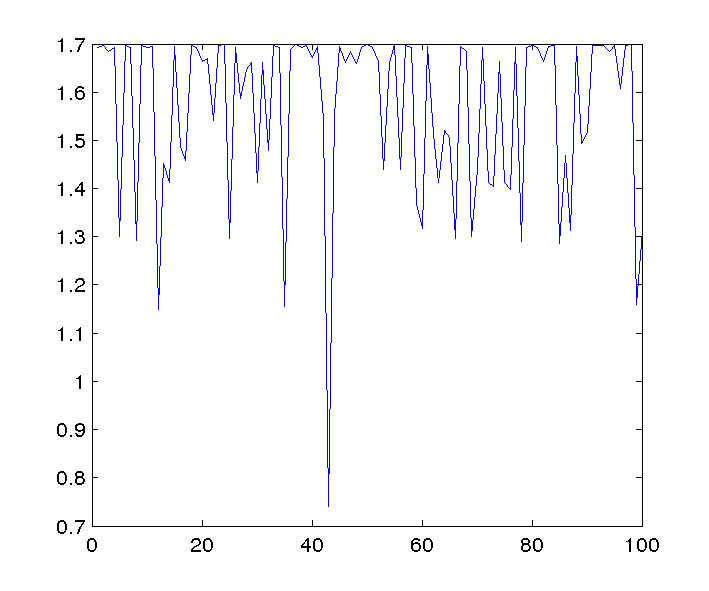}
\end{minipage}\hfill
\begin{minipage}[c]{.33\linewidth}
	\includegraphics[width=\linewidth,trim = .5cm .3cm .6cm .1cm,clip]{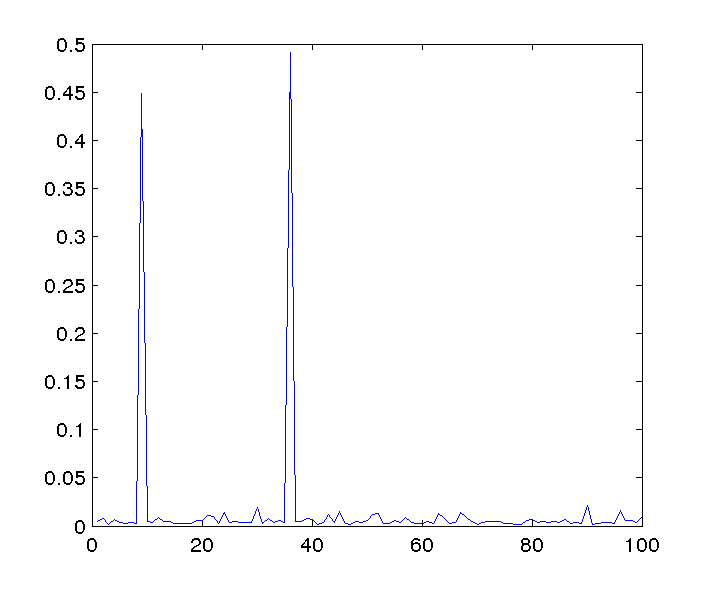}
\end{minipage}\hfill
\begin{minipage}[c]{.33\linewidth}
	\includegraphics[width=\linewidth,trim = .5cm .3cm .6cm .1cm,clip]{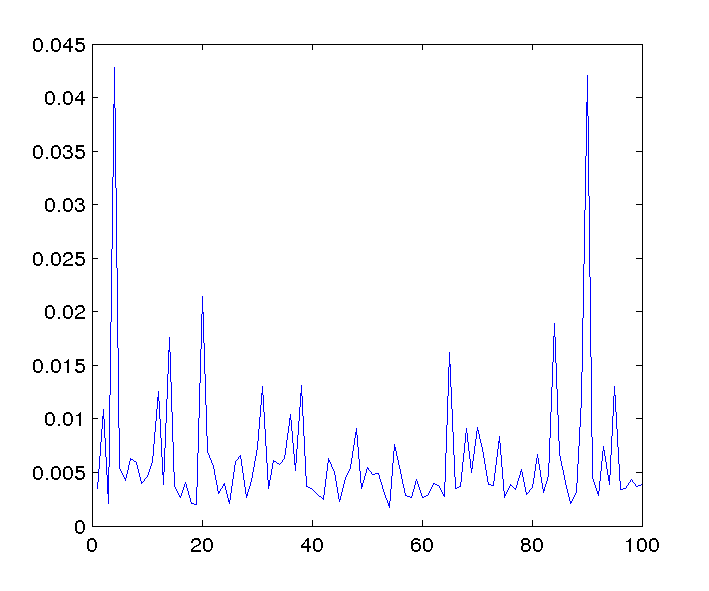}
\end{minipage}
\caption[Simulation of the distance between $\Leb$ and $\mu^{(h_i)_N}_{(1/2,1/2)}$ for 3 examples of conservative diffeomorphisms]{Distance between Lebesgue measure and the measure $\mu^{(h_i)_N}_{(1/2,1/2)}$ depending on $N$ for $h_1$ (left), $h_2$ (middle) and $h_3$ (right), on the grids $E_N$ with $N=2^{20}+k$, $k=1,\cdots,100$.}\label{GrafDistLebPhys}
\end{figure}

To compute these measures, we used Floyd's Algorithm (or the ``tortoise and the hare algorithm''). It has appeared that on the examples of diffeomorphisms we have tested, we were able to test orders of discretization $N\simeq 2^{20}$. Thus, the first figures represent the measures $\mu^{f_N}_x$ for $N\in\llbracket 2^{20}+1,2^{20}+9\rrbracket$. We have also computed the distance between the measure $\mu^{f_N}_x$ and Lebesgue measure (see Figure~\ref{GrafDistLebPhys}). The distance we have chosen is given by the formula
\[d(\mu,\nu) = \sum_{k=0}^\infty \frac{1}{2^k} \sum_{i,j=0}^{2^k-1} \big| \mu(C_{i,j,k}) - \nu(C_{i,j,k})\big|\in[0,2],\]
where
\[C_{i,j,k} = \left[\frac{i}{2^k},\frac{i+1}{2^k}\right] \times \left[\frac{j}{2^k},\frac{j+1}{2^k}\right].\]
In practice, we have computed an approximation of this quantity by summing only on the $k\in\llbracket 0,7 \rrbracket$.

\bigskip

\begin{figure}[ht]
\begin{minipage}[c]{.31\linewidth}
	\includegraphics[height=4.8cm,trim = 1.5cm .95cm 2.8cm .5cm,clip]{Fichiers/MesPhys/IdC1Bis/MesureLog20-1.png}
\end{minipage}\hfill
\begin{minipage}[c]{.31\linewidth}
	\includegraphics[height=4.8cm,trim = 1.5cm .95cm 2.8cm .5cm,clip]{Fichiers/MesPhys/IdC1Bis/MesureLog20-2.png}
\end{minipage}\hfill
\begin{minipage}[c]{.37\linewidth}
	\includegraphics[height=4.8cm,trim = 1.5cm .95cm 1cm .5cm,clip]{Fichiers/MesPhys/IdC1Bis/MesureLog20-3.png}
\end{minipage}

\begin{minipage}[c]{.31\linewidth}
	\includegraphics[height=4.8cm,trim = 1.5cm .95cm 2.8cm .5cm,clip]{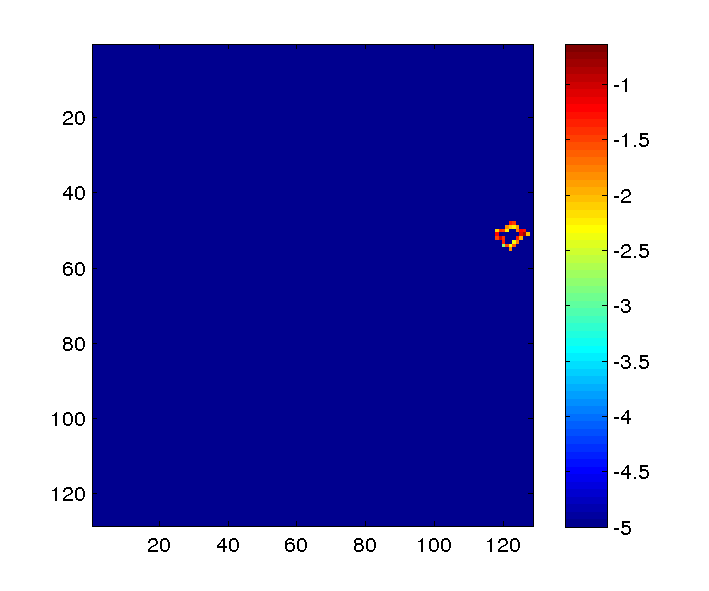}
\end{minipage}\hfill
\begin{minipage}[c]{.31\linewidth}
	\includegraphics[height=4.8cm,trim = 1.5cm .95cm 2.8cm .5cm,clip]{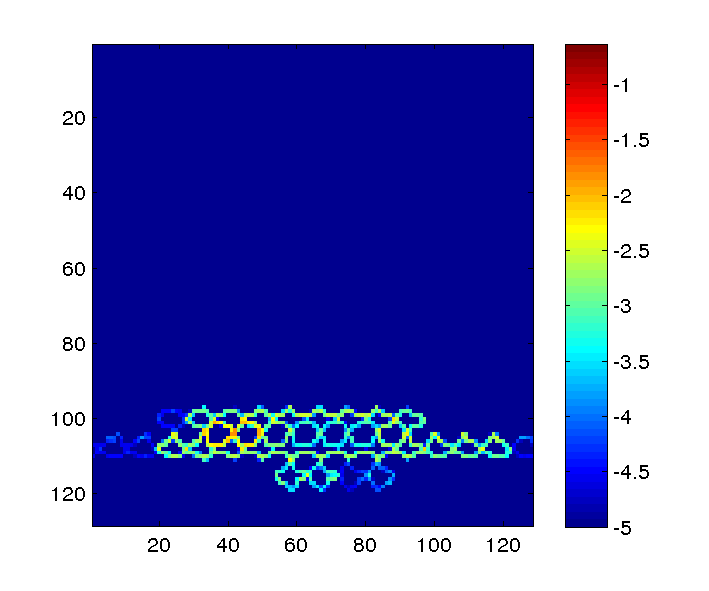}
\end{minipage}\hfill
\begin{minipage}[c]{.37\linewidth}
	\includegraphics[height=4.8cm,trim = 1.5cm .95cm 1cm .5cm,clip]{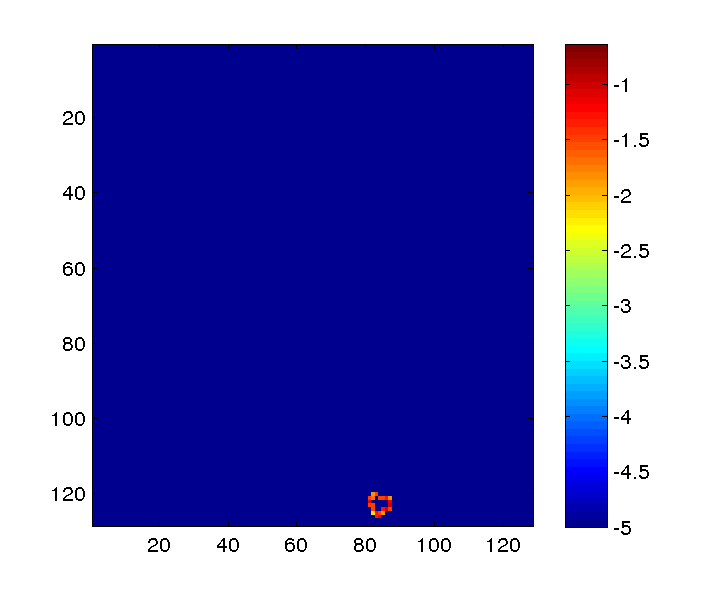}
\end{minipage}

\begin{minipage}[c]{.31\linewidth}
	\includegraphics[height=4.8cm,trim = 1.5cm .95cm 2.8cm .5cm,clip]{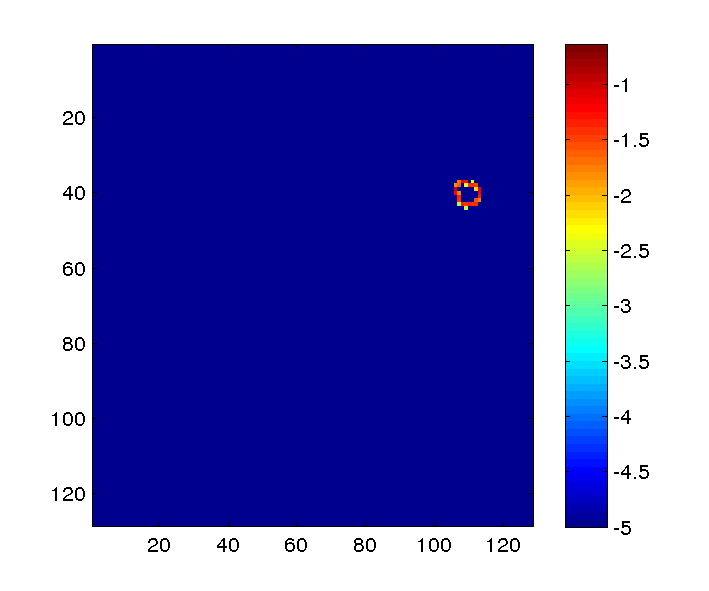}
\end{minipage}\hfill
\begin{minipage}[c]{.31\linewidth}
	\includegraphics[height=4.8cm,trim = 1.5cm .95cm 2.8cm .5cm,clip]{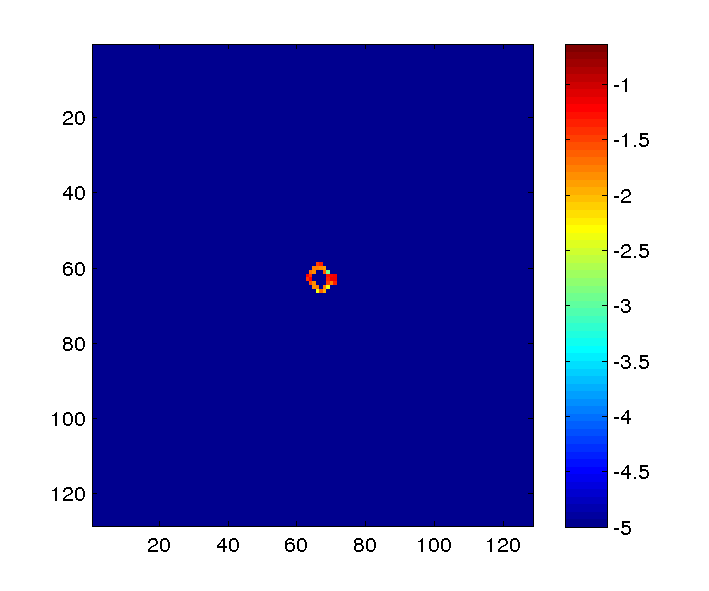}
\end{minipage}\hfill
\begin{minipage}[c]{.37\linewidth}
	\includegraphics[height=4.8cm,trim = 1.5cm .95cm 1cm .5cm,clip]{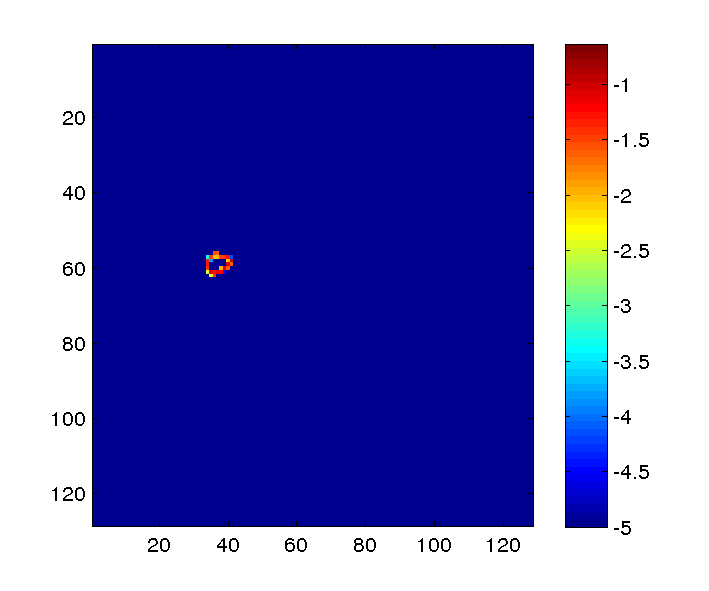}
\end{minipage}
\caption[Simulations of $\mu^{(h_1)_N}_x$ on the grids $E_N$, with $N=2^{20}+i$, $i=1,\cdots,9$]{Simulations of invariant measures $\mu^{(h_1)_N}_x$ on the grids $E_N$, with $N=2^{20}+i$, $i=1,\cdots,9$ and $x=(1/2,1/2)$ (from left to right and top to bottom).}\label{MesPhysIdC1}
\end{figure}

\begin{figure}[ht]
\begin{minipage}[c]{.31\linewidth}
	\includegraphics[height=4.8cm,trim = 1.5cm .95cm 2.8cm .5cm,clip]{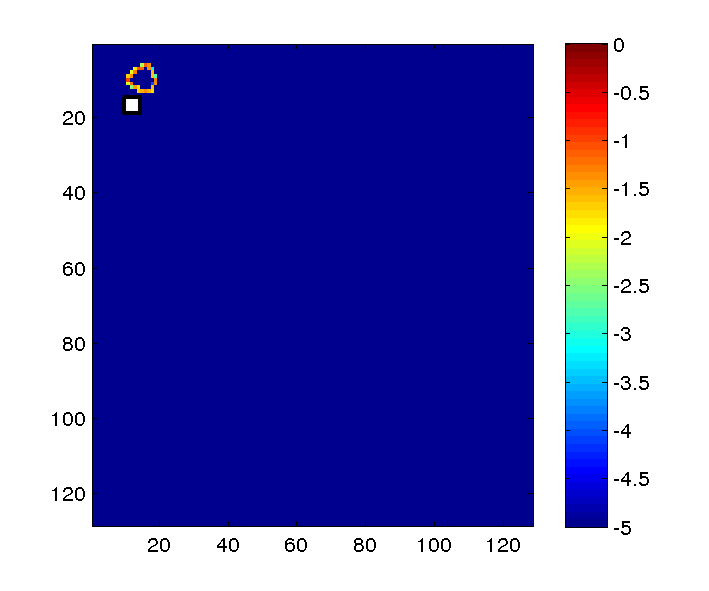}
\end{minipage}\hfill
\begin{minipage}[c]{.31\linewidth}
	\includegraphics[height=4.8cm,trim = 1.5cm .95cm 2.8cm .5cm,clip]{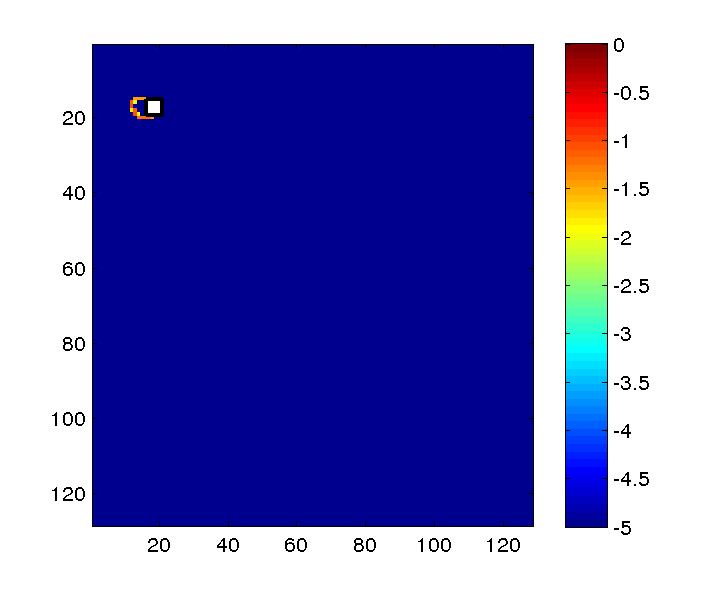}
\end{minipage}\hfill
\begin{minipage}[c]{.37\linewidth}
	\includegraphics[height=4.8cm,trim = 1.5cm .95cm 1cm .5cm,clip]{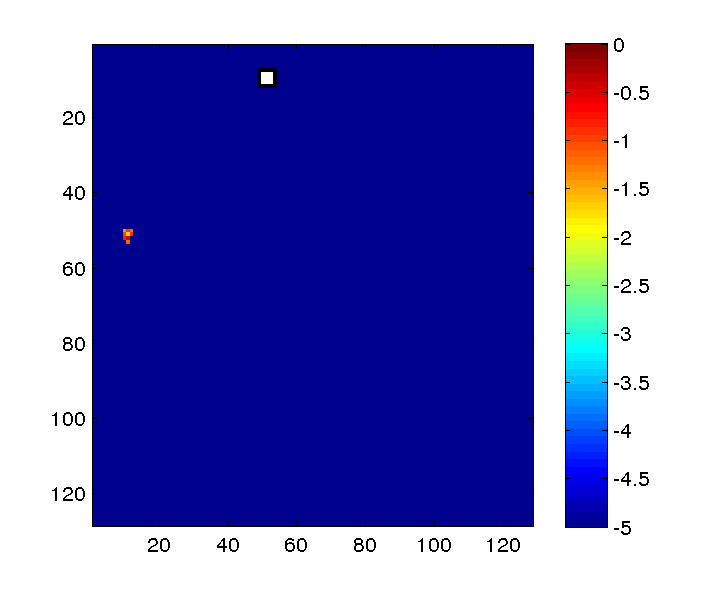}
\end{minipage}

\begin{minipage}[c]{.31\linewidth}
	\includegraphics[height=4.8cm,trim = 1.5cm .95cm 2.8cm .5cm,clip]{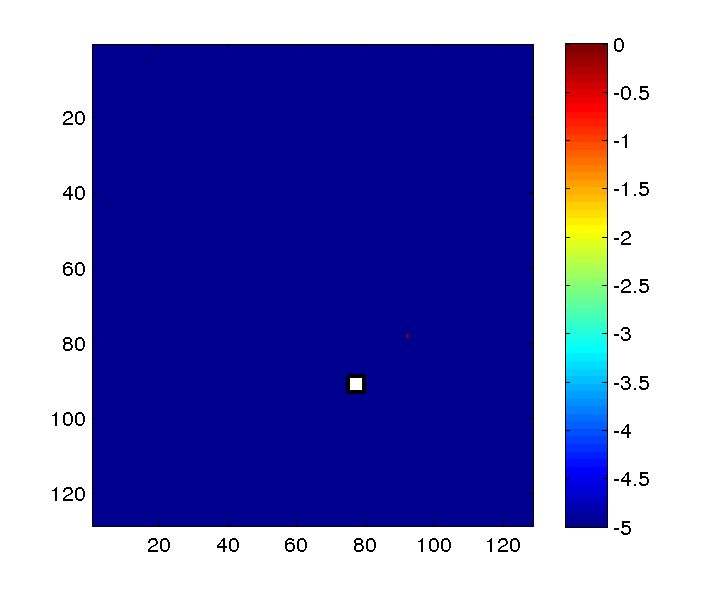}
\end{minipage}\hfill
\begin{minipage}[c]{.31\linewidth}
	\includegraphics[height=4.8cm,trim = 1.5cm .95cm 2.8cm .5cm,clip]{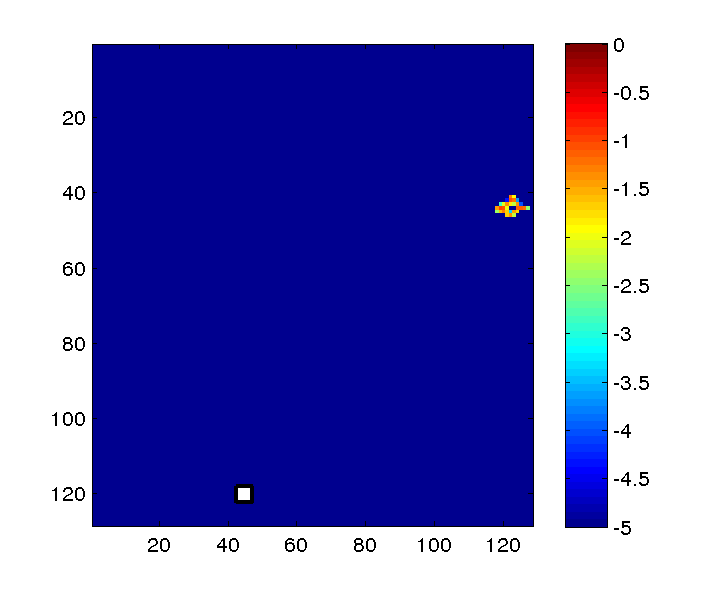}
\end{minipage}\hfill
\begin{minipage}[c]{.37\linewidth}
	\includegraphics[height=4.8cm,trim = 1.5cm .95cm 1cm .5cm,clip]{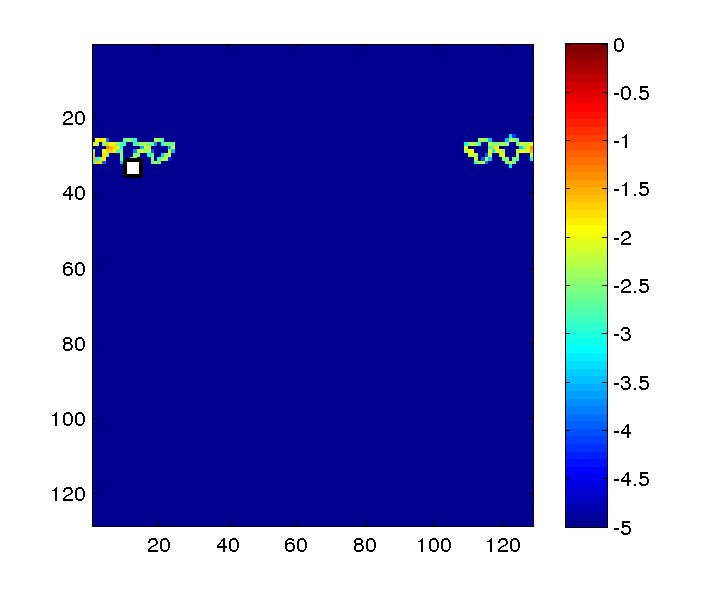}
\end{minipage}
\caption[Simulations of $\mu^{(h_1)_N}_x$ on the grid $E_N$, with $N=2^{23}$ and a random point $x$]{Simulations of invariant measures $\mu^{(h_1)_N}_x$ on the grid $E_N$, with $N=2^{23}$, and $x$ a random point of $\T^2$, represented by the black and white box. The behaviour observed on the top left picture is the most frequent, but we also observe other kind of measures: for example, the measures has a very small support like on the bottom left picture on about 10 of the 100 random draws we have made; we even see appearing the strange behaviour of the last picture once.}\label{MesPhysIdC1Pt}
\end{figure}

In the case of the diffeomorphism $h_1$, which is close to the identity, we observe a strong variation of the measure $\mu^{(h_1)_N}_x$ depending on $N$ (left of Figure~\ref{GrafDistLebPhys} and Figure~\ref{MesPhysIdC1}). More precisely, for 7 on the 9 orders of discretization represented, these measures seems to be supported by a small curve; for $N = 2^{20}+3$, this measure seems to be supported by a figure-8 curve, and for $N = 2^{20}+5$, the support of the measure is quite complicated and looks like an interlaced curve. The fact that the measures $\mu^{(h_1)_N}_x$ strongly depend on $N$ reflects the behaviour predicted by Theorem~\ref{TheoMesPhysDiff}: in theory, for a generic $C^1$ diffeomorphism, the measures $\mu^{f_N}_x$ should accumulate on the whole set of $f$-invariant measures; here we see that these measures strongly depend on $N$ (moreover, we can see on Figure~\ref{GrafDistLebPhys} that on the orders of discretization we have tested, these measures are never close to Lebesgue measure). We have no satisfying explanation to the specific shape of the supports of the measures. When we fix the order of discretization and make vary the starting point $x$, the behaviour is very similar: the measures $\mu^{(h_1)_N}_x$ widely depend on the point $x$ (see Figure~\ref{MesPhysIdC1Pt}). We also remark that increasing the order of discretizations does not make the measures $\mu^{(h_1)_N}_x$ evolve more smoothly.
\bigskip

\begin{figure}[ht]
\begin{minipage}[c]{.31\linewidth}
	\includegraphics[height=4.8cm,trim = 1.5cm .95cm 2.8cm .5cm,clip]{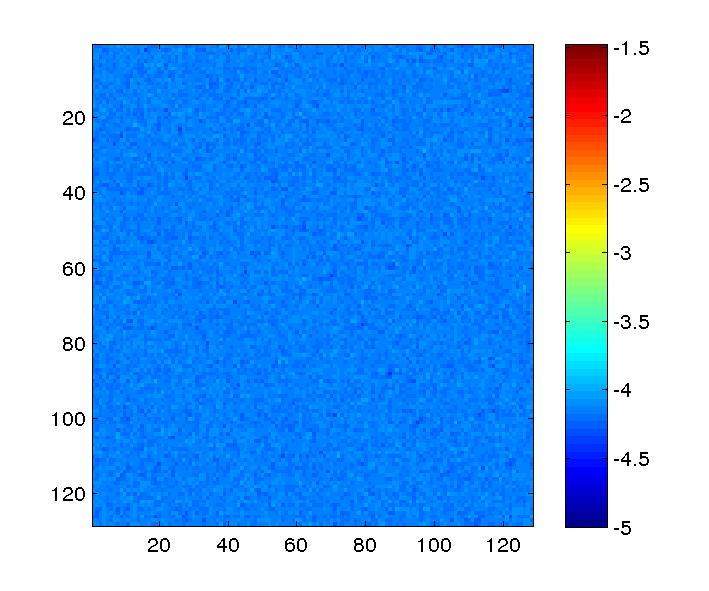}
\end{minipage}\hfill
\begin{minipage}[c]{.31\linewidth}
	\includegraphics[height=4.8cm,trim = 1.5cm .95cm 2.8cm .5cm,clip]{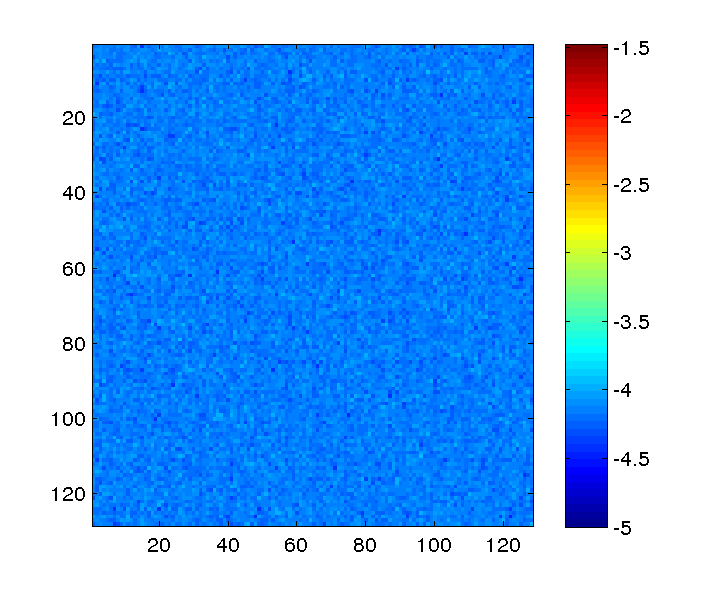}
\end{minipage}\hfill
\begin{minipage}[c]{.37\linewidth}
	\includegraphics[height=4.8cm,trim = 1.5cm .95cm 1cm .5cm,clip]{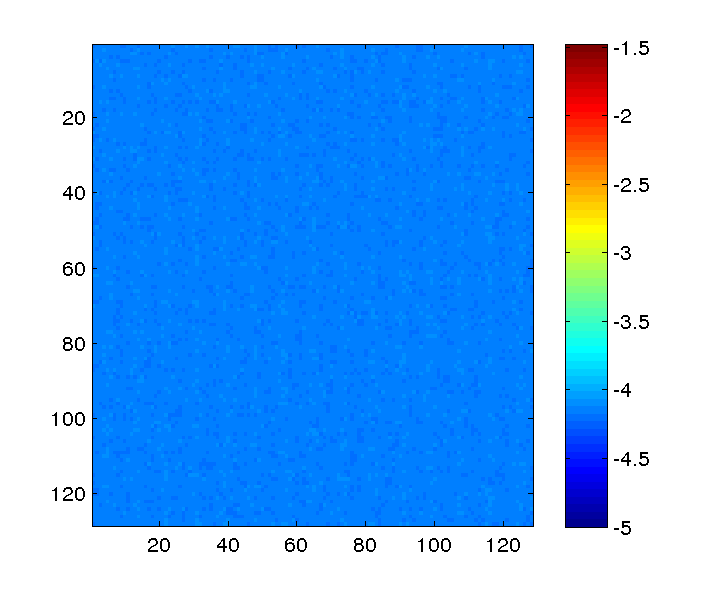}
\end{minipage}

\begin{minipage}[c]{.31\linewidth}
	\includegraphics[height=4.8cm,trim = 1.5cm .95cm 2.8cm .5cm,clip]{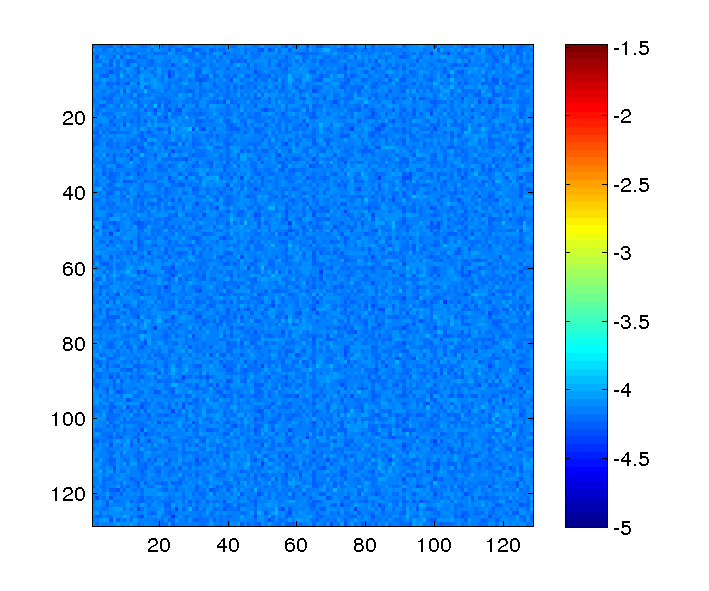}
\end{minipage}\hfill
\begin{minipage}[c]{.31\linewidth}
	\includegraphics[height=4.8cm,trim = 1.5cm .95cm 2.8cm .5cm,clip]{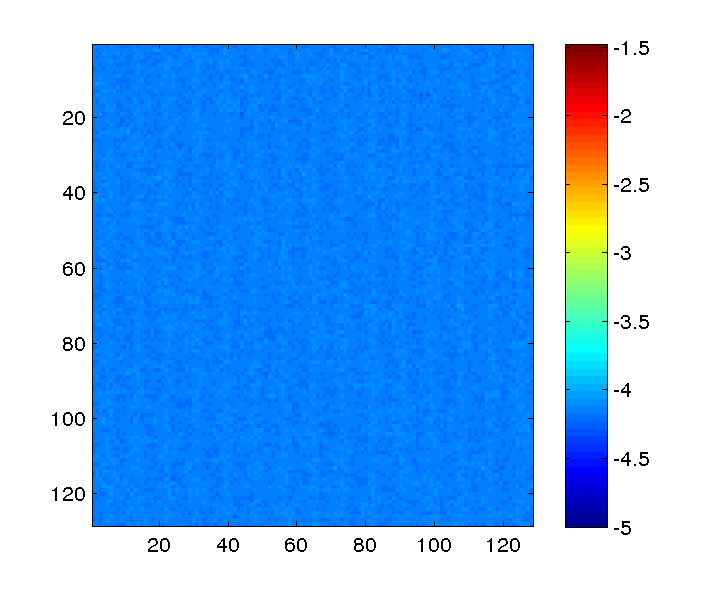}
\end{minipage}\hfill
\begin{minipage}[c]{.37\linewidth}
	\includegraphics[height=4.8cm,trim = 1.5cm .95cm 1cm .5cm,clip]{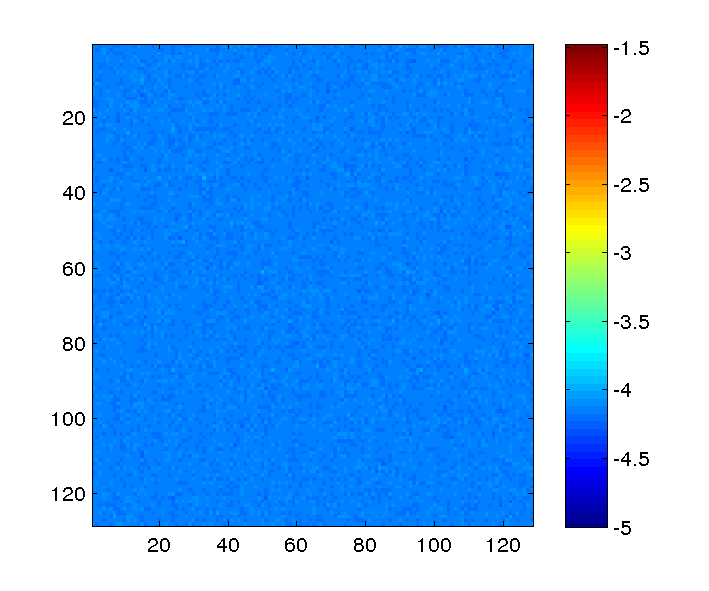}
\end{minipage}

\begin{minipage}[c]{.31\linewidth}
	\includegraphics[height=4.8cm,trim = 1.5cm .95cm 2.8cm .5cm,clip]{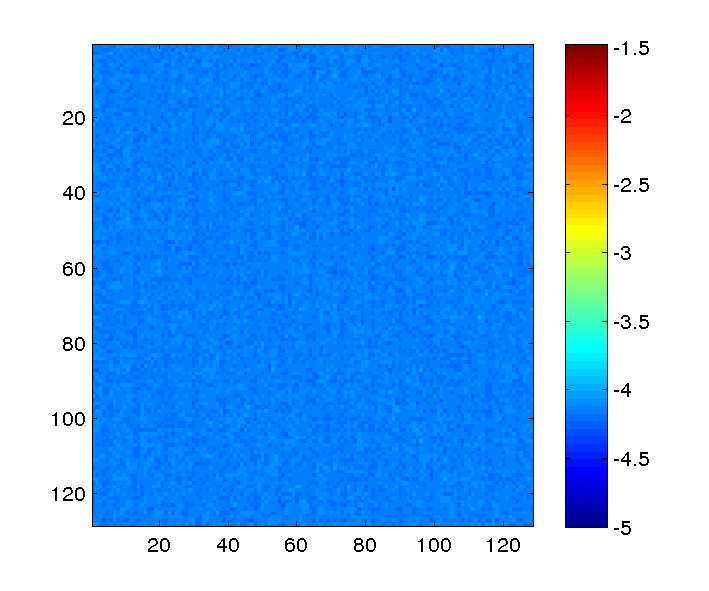}
\end{minipage}\hfill
\begin{minipage}[c]{.31\linewidth}
	\includegraphics[height=4.8cm,trim = 1.5cm .95cm 2.8cm .5cm,clip]{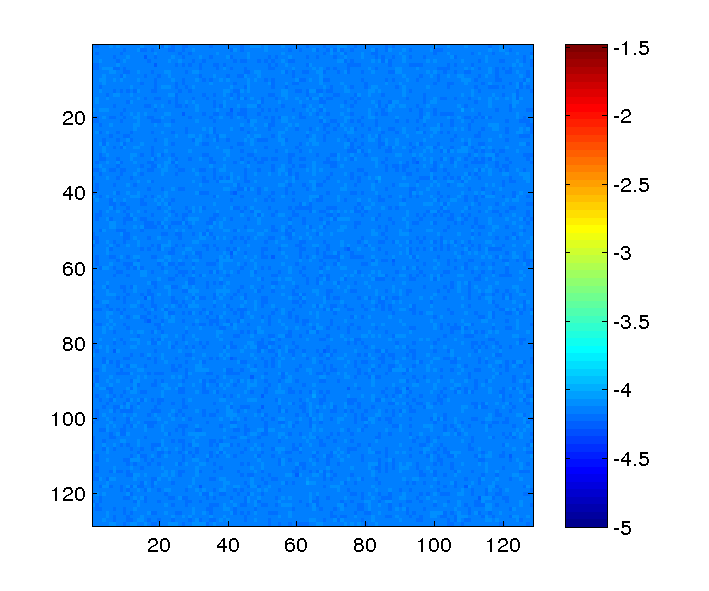}
\end{minipage}\hfill
\begin{minipage}[c]{.37\linewidth}
	\includegraphics[height=4.8cm,trim = 1.5cm .95cm 1cm .5cm,clip]{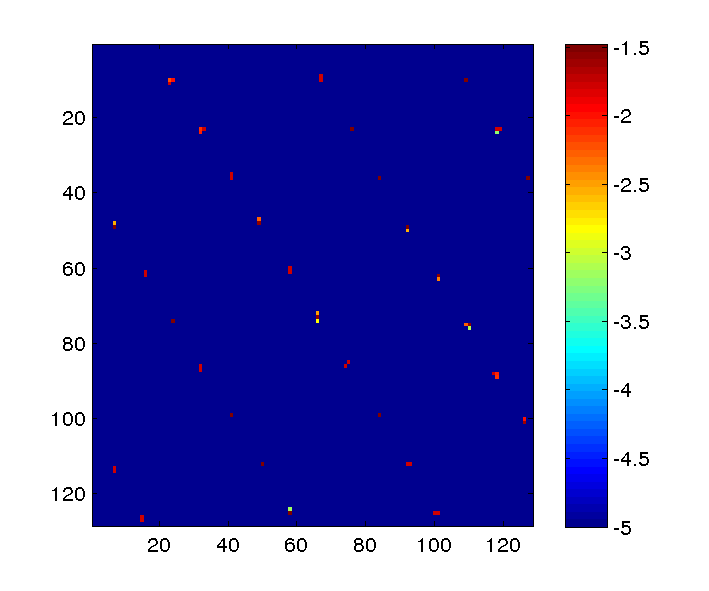}
\end{minipage}
\caption[Simulations of $\mu^{(h_2)_N}_x$ on the grids $E_N$, with $N=2^{20}+i$, $i=1,\cdots,9$]{Simulations of invariant measures $\mu^{(h_2)_N}_x$ on the grids $E_N$,  with $N=2^{20}+i$, $i=1,\cdots,9$ and $x=(1/2,1/2)$ (from left to right and top to bottom).}\label{MesPhysRotC1}
\end{figure}

\begin{figure}[ht]
\begin{minipage}[c]{.31\linewidth}
	\includegraphics[height=4.8cm,trim = 1.5cm .95cm 2.8cm .5cm,clip]{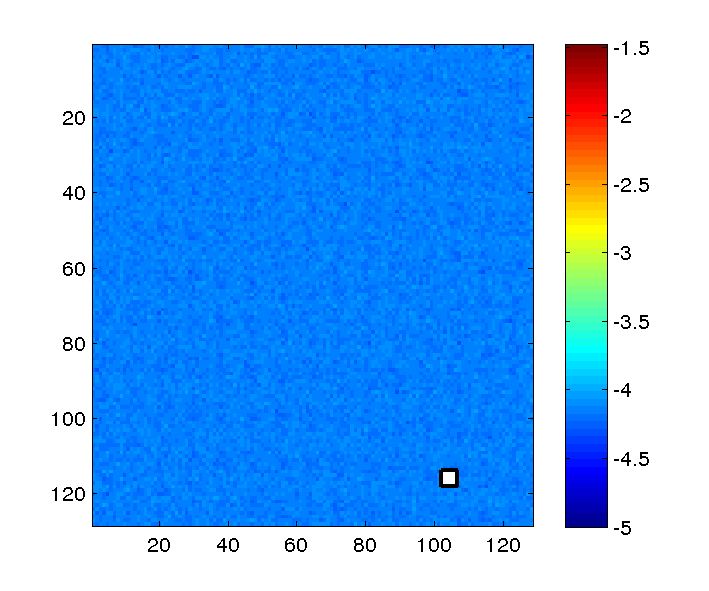}
\end{minipage}\hfill
\begin{minipage}[c]{.31\linewidth}
	\includegraphics[height=4.8cm,trim = 1.5cm .95cm 2.8cm .5cm,clip]{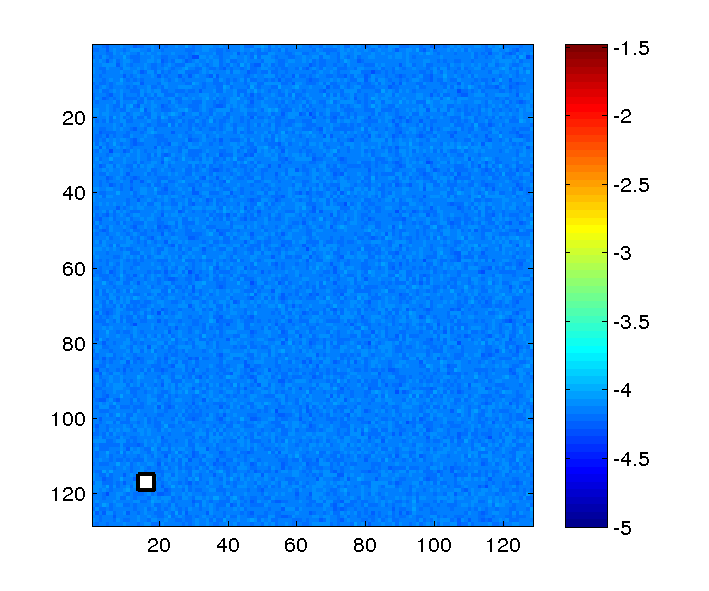}
\end{minipage}\hfill
\begin{minipage}[c]{.37\linewidth}
	\includegraphics[height=4.8cm,trim = 1.5cm .95cm 1cm .5cm,clip]{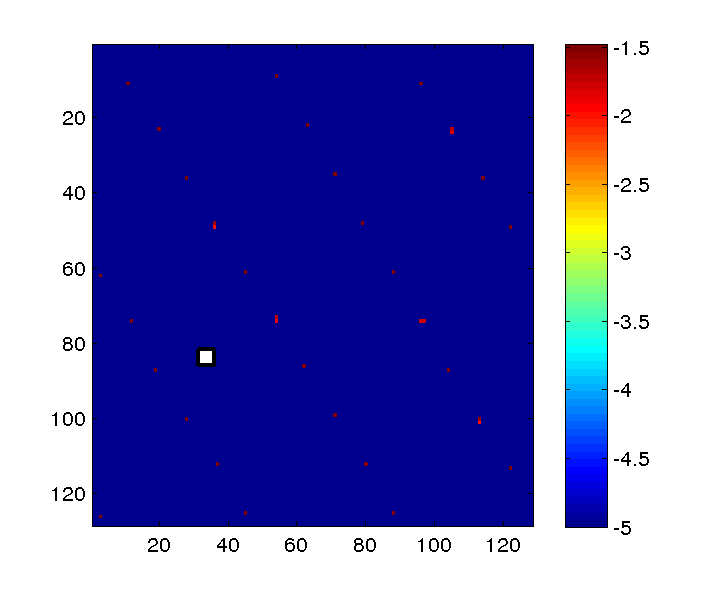}
\end{minipage}
\caption[Simulations of $\mu^{(h_2)_N}_x$ on the grid $E_N$, with $N=2^{23}+5$ and a random point $x$]{Simulations of invariant measures $\mu^{(h_2)_N}_x$ on the grid $E_N$, with $N=2^{23}+5$, and $x$ a random point of $\T^2$, represented by the black and white box. The behaviour observed on the two firsts pictures is the most frequent, but sometimes (in fact, twice on $1\,000$ random draws), we also observe a measure whose support is very close to a periodic point with small period, like on the right picture.}\label{MesPhysRotC1Pt}
\end{figure}

The measures $\mu^{(h_2)_N}_x$ vary less than the measures $\mu^{(h_1)_N}_x$ (recall that $h_2$ is a small perturbation of a rotation of order $30$). For the 8 first measures of the Figure~\ref{MesPhysRotC1}, we obtain a measure which is very close to Lebesgue measure. The fact that we obtain measures closer to Lebesgue measures was predictable: the dynamics $h_2$ is close to the rotation $R$, whose orbits are better distributed in the torus than that of the identity. But for the order $N=2^{20}+9$, the measure we obtain is very different from the previous one: its support seems quite close to a real orbit of $R$ (of order 30), in particular this support covers a very small proportion of $\T^2$. This is what is predicted by Theorem~\ref{TheoMesPhysDiff}: at least sometimes, the measures $\mu^{(h_2)_N}_x$ should not be close to Lebesgue measure. We observe exactly this behaviour when we make simulations for more different orders of discretization $N$ (middle of Figure~\ref{GrafDistLebPhys}): for two orders $N$ between $2^{20}+1$ and $2^{20}+100$, the measure is far away from Lebesgue measure. Remark that the same behaviour holds when we fix the order of discretization and make the point $x$ vary (Figure~\ref{MesPhysRotC1Pt}); however, we observe that the frequency of occurrence of the event ''the measure $\mu^{(h_2)_N}_x$ is close to a periodic measure with small period'' is smaller in the case of Figure~\ref{MesPhysRotC1Pt} than in the case of Figure~\ref{MesPhysRotC1}. We think that it is more due to the fact that the order of discretization is bigger in the first case, than to a fundamental difference of the processes used to produce these simulations.
\bigskip

\begin{figure}[ht]
\begin{minipage}[c]{.31\linewidth}
	\includegraphics[height=4.8cm,trim = 1.5cm .95cm 2.8cm .5cm,clip]{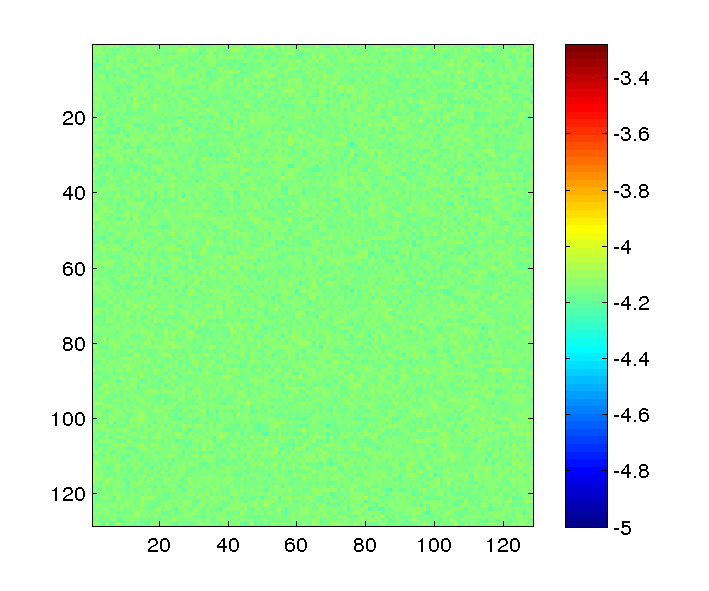}
\end{minipage}\hfill
\begin{minipage}[c]{.31\linewidth}
	\includegraphics[height=4.8cm,trim = 1.5cm .95cm 2.8cm .5cm,clip]{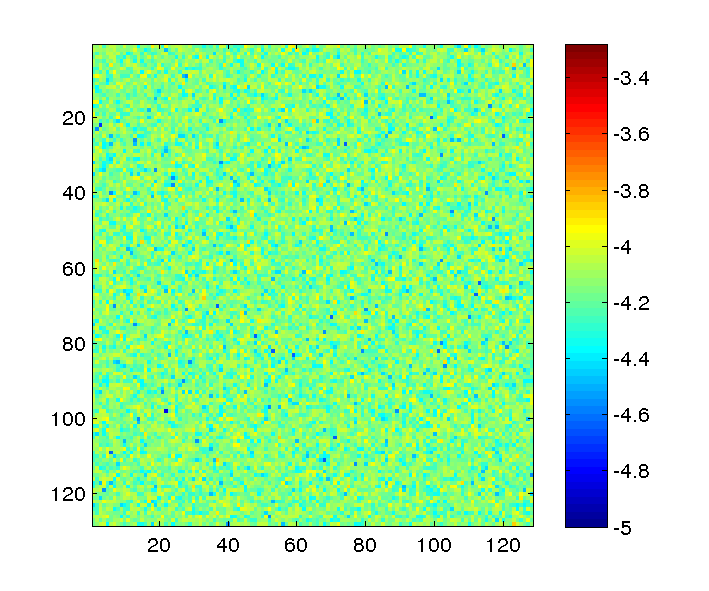}
\end{minipage}\hfill
\begin{minipage}[c]{.37\linewidth}
	\includegraphics[height=4.8cm,trim = 1.5cm .95cm 1cm .5cm,clip]{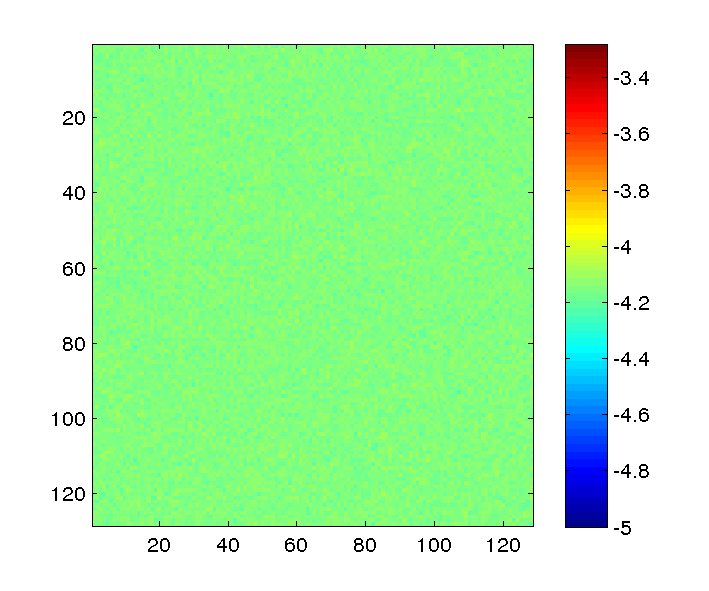}
\end{minipage}

\begin{minipage}[c]{.31\linewidth}
	\includegraphics[height=4.8cm,trim = 1.5cm .95cm 2.8cm .5cm,clip]{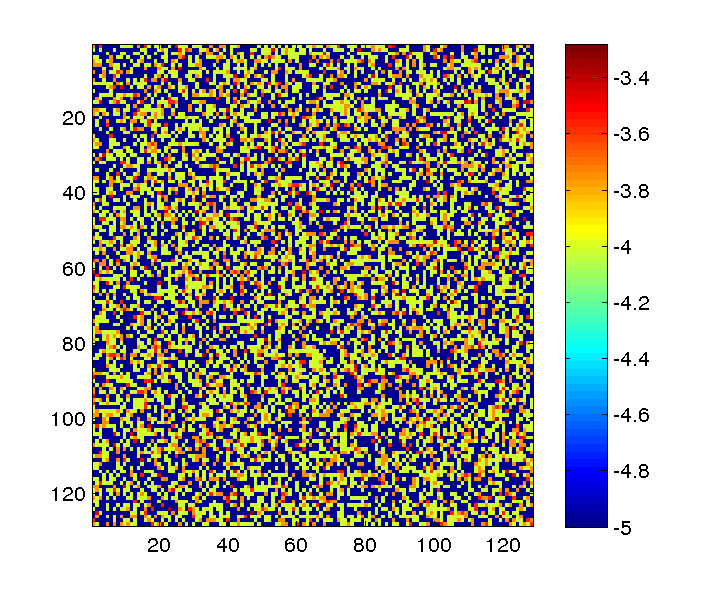}
\end{minipage}\hfill
\begin{minipage}[c]{.31\linewidth}
	\includegraphics[height=4.8cm,trim = 1.5cm .95cm 2.8cm .5cm,clip]{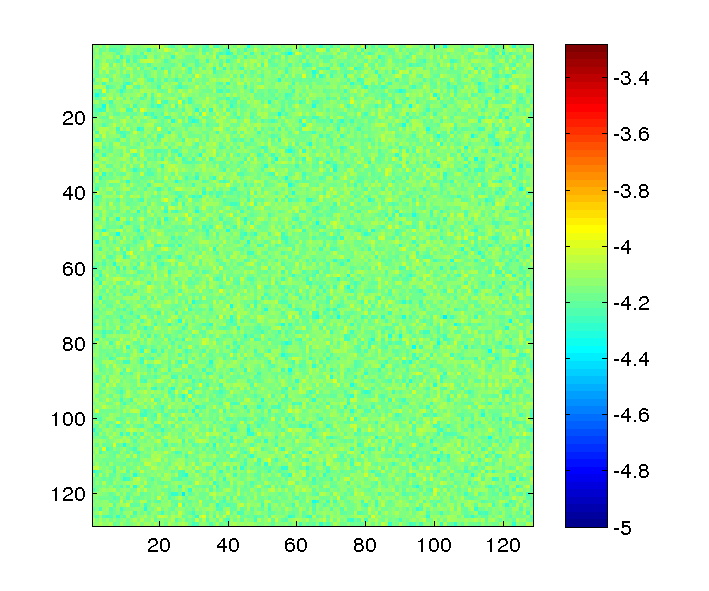}
\end{minipage}\hfill
\begin{minipage}[c]{.37\linewidth}
	\includegraphics[height=4.8cm,trim = 1.5cm .95cm 1cm .5cm,clip]{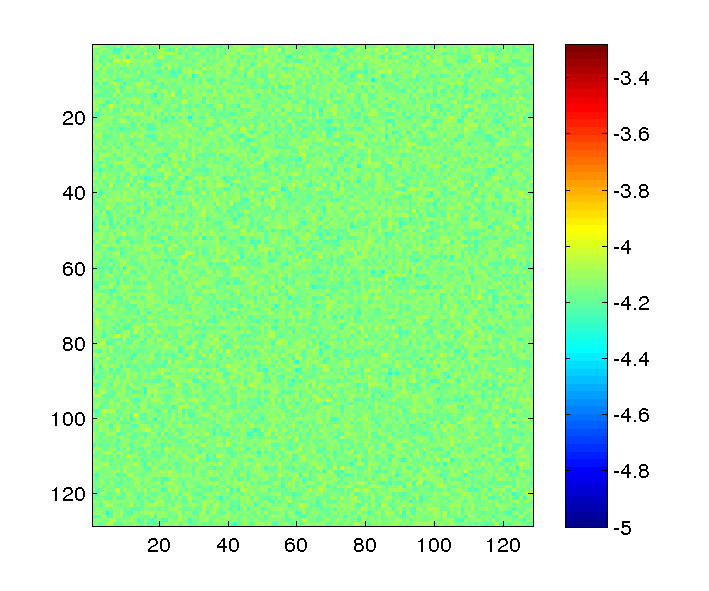}
\end{minipage}

\begin{minipage}[c]{.31\linewidth}
	\includegraphics[height=4.8cm,trim = 1.5cm .95cm 2.8cm .5cm,clip]{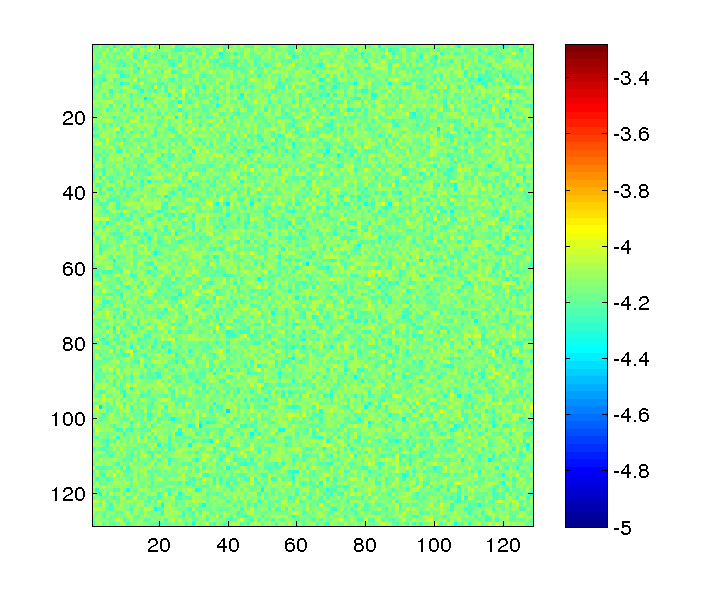}
\end{minipage}\hfill
\begin{minipage}[c]{.31\linewidth}
	\includegraphics[height=4.8cm,trim = 1.5cm .95cm 2.8cm .5cm,clip]{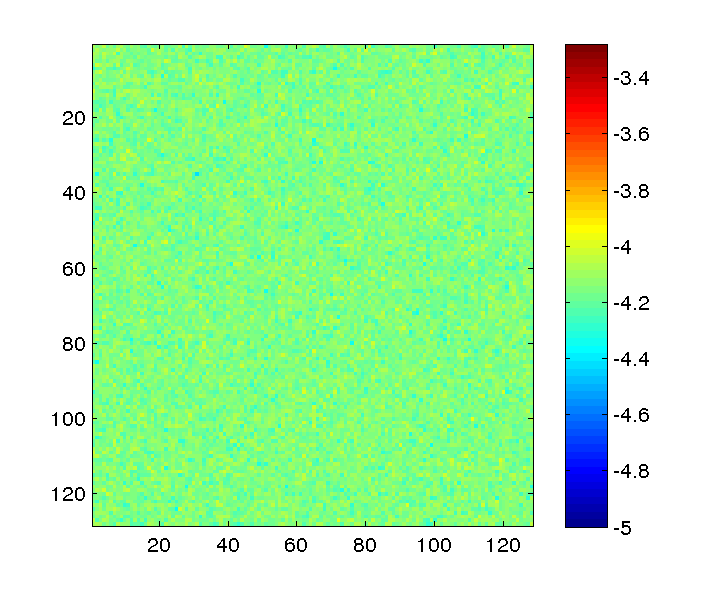}
\end{minipage}\hfill
\begin{minipage}[c]{.37\linewidth}
	\includegraphics[height=4.8cm,trim = 1.5cm .95cm 1cm .5cm,clip]{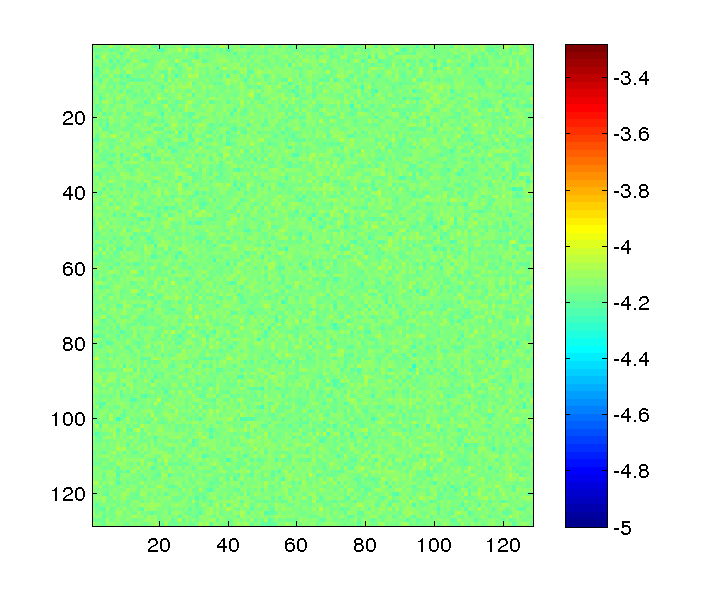}
\end{minipage}
\caption[Simulations of $\mu^{(h_3)_N}_x$ on the grids $E_N$, with $N=2^{20}+i$, $i=1,\cdots,9$]{Simulations of invariant measures $\mu^{(h_3)_N}_x$ on the grids $E_N$,  with $N=2^{20}+i$, $i=1,\cdots,9$ and $x=(1/2,1/2)$ (from left to right and top to bottom).}\label{MesPhysAnoC1}
\end{figure}

\begin{figure}[ht]
\begin{minipage}[c]{.31\linewidth}
	\includegraphics[height=4.8cm,trim = 1.5cm .95cm 2.8cm .5cm,clip]{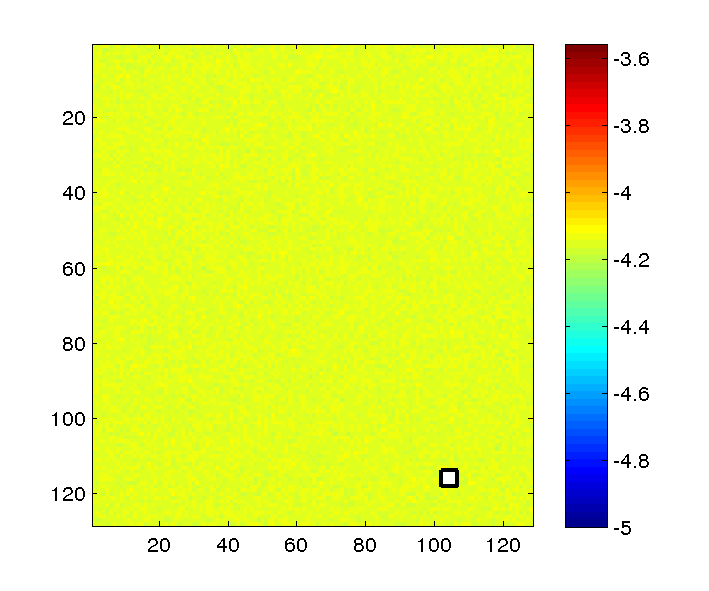}
\end{minipage}\hfill
\begin{minipage}[c]{.31\linewidth}
	\includegraphics[height=4.8cm,trim = 1.5cm .95cm 2.8cm .5cm,clip]{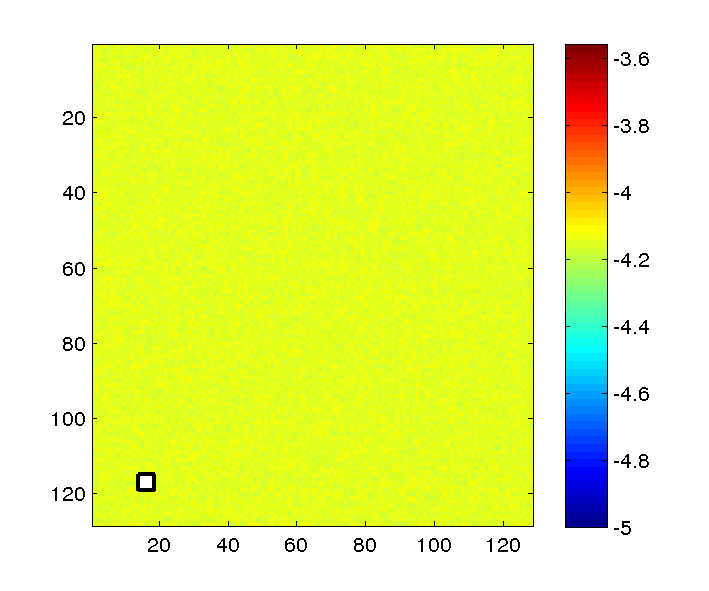}
\end{minipage}\hfill
\begin{minipage}[c]{.37\linewidth}
	\includegraphics[height=4.8cm,trim = 1.5cm .95cm 1cm .5cm,clip]{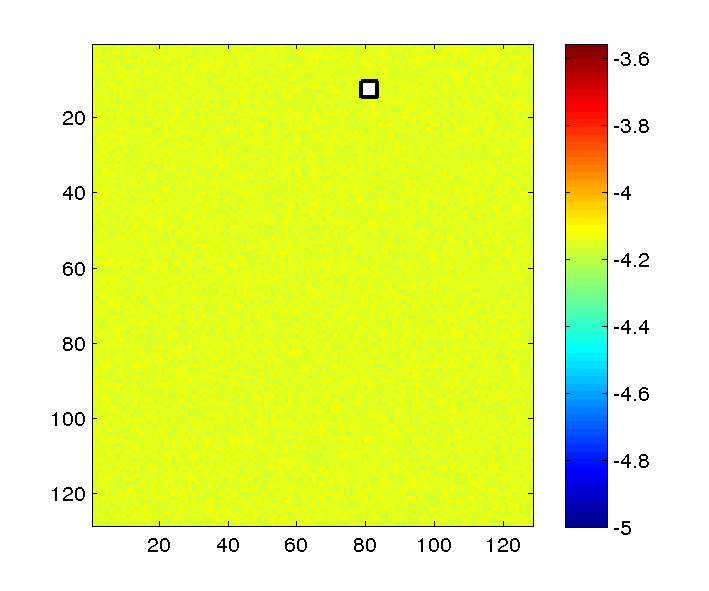}
\end{minipage}

\begin{minipage}[c]{.31\linewidth}
	\includegraphics[height=4.8cm,trim = 1.5cm .95cm 2.8cm .5cm,clip]{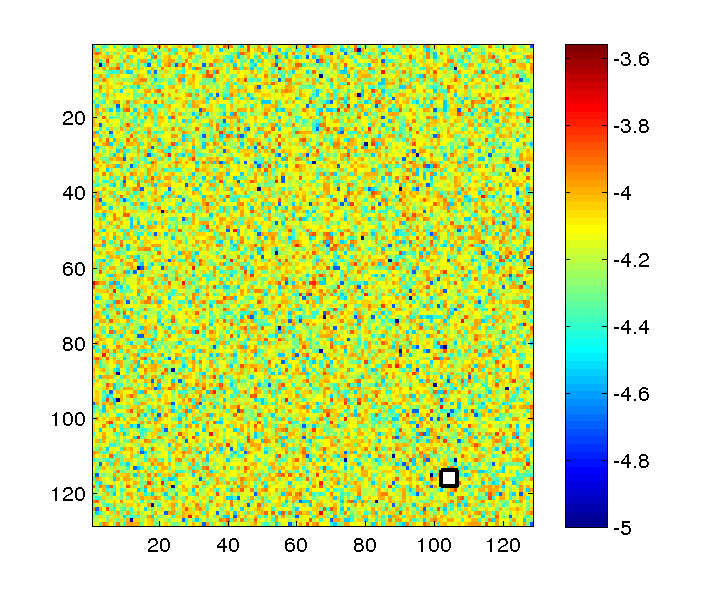}
\end{minipage}\hfill
\begin{minipage}[c]{.31\linewidth}
	\includegraphics[height=4.8cm,trim = 1.5cm .95cm 2.8cm .5cm,clip]{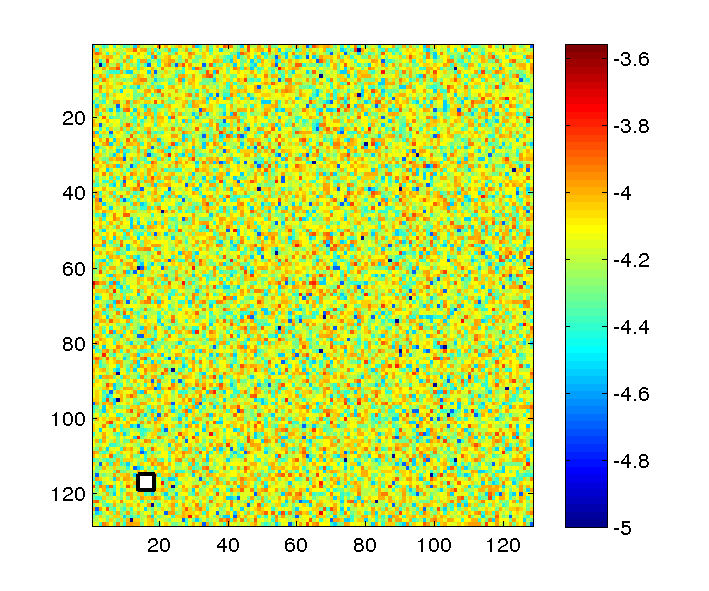}
\end{minipage}\hfill
\begin{minipage}[c]{.37\linewidth}
	\includegraphics[height=4.8cm,trim = 1.5cm .95cm 1cm .5cm,clip]{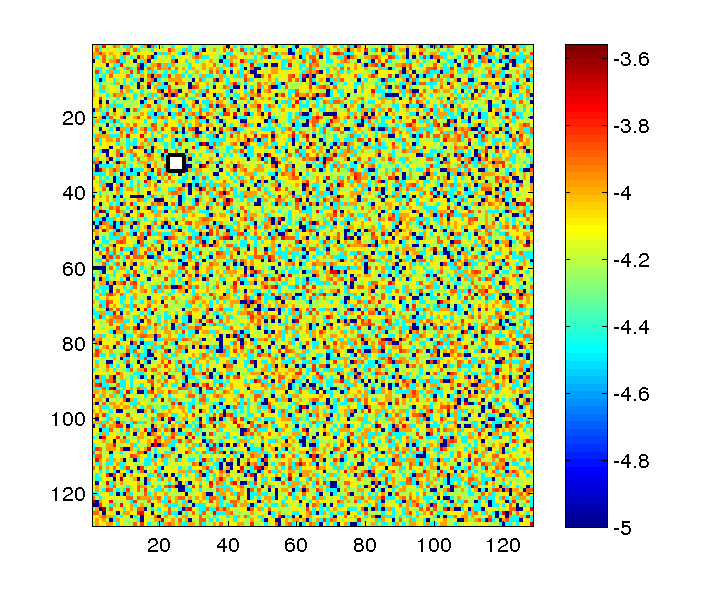}
\end{minipage}

\caption[Simulations of $\mu^{(h_3)_N}_x$ on the grid $E_N$, with $N\simeq 2^{23}$ and a random point $x$]{Simulations of invariant measures $\mu^{(h_3	)_N}_x$ on the grid $E_N$, with $N=2^{23}+1$ (top) and $N=2^{23}+17$ (bottom), and $x$ a random point of $\T^2$, represented by the black and white box. The behaviour observed on the three top pictures is the most frequent (for 17 over the 20 orders $N=2^{23}+i$, $i=0,\cdots,19$, all the 100 random draws we have made gave a measure very close to $\Leb$), but seldom we also observe measures further from Lebesgue measure, like what happens for $N=2^{23}+17$ (bottom), where $99$ over the $100$ random draws of $x$ produce a measure identical to the two first pictures, and the other random draw gives a measure a bit more singular with respect to Lebesgue measure (bottom right).}\label{MesPhysAnoC1Pt}
\end{figure}

The behaviour of the measures $\mu^{(h_3)_N}_x$, where $h_3$ is a small $C^1$-perturbation of the linear Anosov map $A$, is quite close to that of the the measures $\mu^{(h_2)_N}_x$ (see Figure~\ref{MesPhysAnoC1}): most of the time, these measures are close to Lebesgue measure, but for one order of discretization $N$ (here, $N = 2^{20}+4$), the measure becomes very different from Lebesgue measure (we can see on the right of Figure~\ref{GrafDistLebPhys} that this phenomenon appears twice when $N\in\llbracket 2^{20}+1,2^{20}+100\rrbracket$). The difference with the case of $h_2$ is that here, the ``exceptional'' measure is much better distributed than for $h_2$ (the maximal measure of a region of size $1/128\times 1/128$ is close to $10^{-3.5}$ for $h_3$ and close to $10^{-1.5}$ for $h_2$). The same phenomenon holds when we fix the order of discretization but change the starting point $x$ (see Figure~\ref{MesPhysAnoC1Pt}), except that (as for $f_2$) the number of apparition of measures that are singular with respect to Lebesgue measure is smaller than in Figure~\ref{MesPhysAnoC1}. Again, we think that this follows from the fact that the orders of discretizations tested are bigger. In this case, the simulations suggest the following behaviour: when the order of discretization $N$ increases, the frequency of apparition of measures $\mu^{(h_3)_N}_x$ far away from Lebesgue measures tends to $0$.

Recall that Addendum~\ref{AddTheoMesPhysDiff} states that if $x$ is fixed, then for a generic $f\in\Diff^1(\T^2,\Leb)$, the measures $\mu_x^{f_N}$ accumulate on the whole set of $f$-invariant measures, but do not say anything about, for instance, the frequency of orders $N$ such that $\mu_x^{f_N}$ is not close to Lebesgue measure. It is natural to think that this frequency depends a lot on $f$; for example that such $N$ are very rare close to an Anosov diffeomorphism and more frequent close to an ``elliptic'' dynamics like the identity. The results of numerical simulations seem to confirm this heuristic.

\clearpage

\subsection[Simulations of the measures $\mu^{f_N}_{\T^2}$]{Simulations of the measures $\mu^{f_N}_{\T^2}$ for conservative torus diffeomorphisms}\label{RedW}

We now present the results of numerical simulations of the measures $\mu_{\T^2}^{f_N}$. Recall that the measure $\mu_{\T^2}^{f_N}$ is supported by the union of periodic orbits of $f_N$, and is such that the total measure of each periodic orbit is equal to the cardinality of its basin of attraction.

First, we simulate a conservative diffeomorphism $f_5$ which is close to the identity in the $C^1$ topology. We have chosen $f_5 = Q\circ P$, with
\[p(x) = \frac{1}{209}\cos(2\pi\times 17x)+\frac{1}{271}\sin(2\pi\times 27x)-\frac{1}{703}\cos(2\pi\times 35x),\]
\[q(y) = \frac{1}{287}\cos(2\pi\times 15y)+\frac{1}{203}\sin(2\pi\times 27y)-\frac{1}{841}\sin(2\pi\times 38y).\]

We have also simulated the conservative diffeomorphism $f_6 = f_5\circ A$, with $A$ the standard Anosov automorphism
\[A = \begin{pmatrix} 2 & 1 \\ 1 & 1 \end{pmatrix},\]
thus $f_6$ is a small $C^1$ perturbation of $A$; in particular the theory asserts that it is topologically conjugated to $A$. We can test whether this property can be observed on simulations or not.
\bigskip

\begin{figure}[h!]
\begin{center}
\makebox[0.8\textwidth]{\parbox{0.8\textwidth}{%
\begin{minipage}[c]{.49\linewidth}
	\includegraphics[width=\linewidth,trim = .5cm .3cm .6cm .1cm,clip]{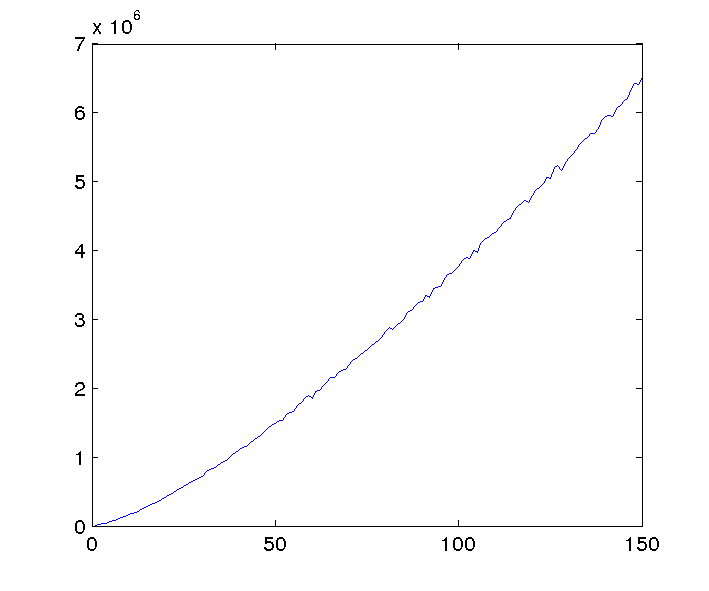}
\end{minipage}\hfill
\begin{minipage}[c]{.49\linewidth}
	\includegraphics[width=\linewidth,trim = .5cm .3cm .6cm .1cm,clip]{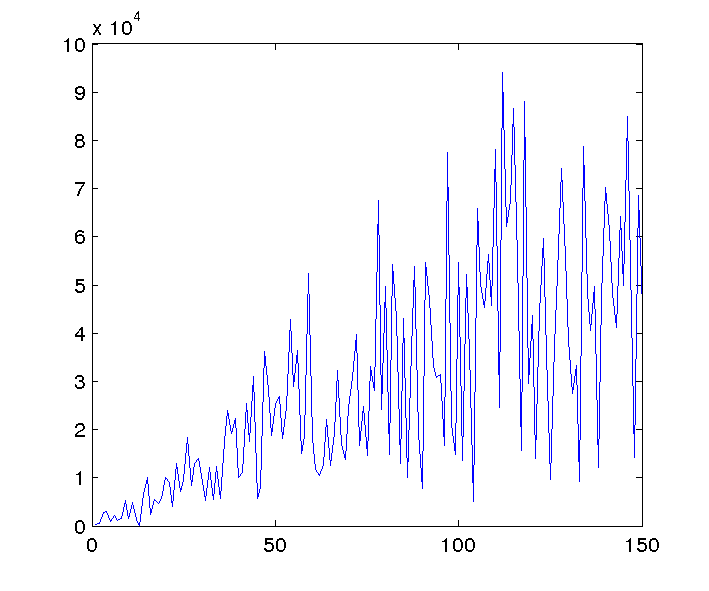}
\end{minipage}

\begin{minipage}[c]{.49\linewidth}
	\includegraphics[width=\linewidth,trim = .5cm .3cm .6cm .1cm,clip]{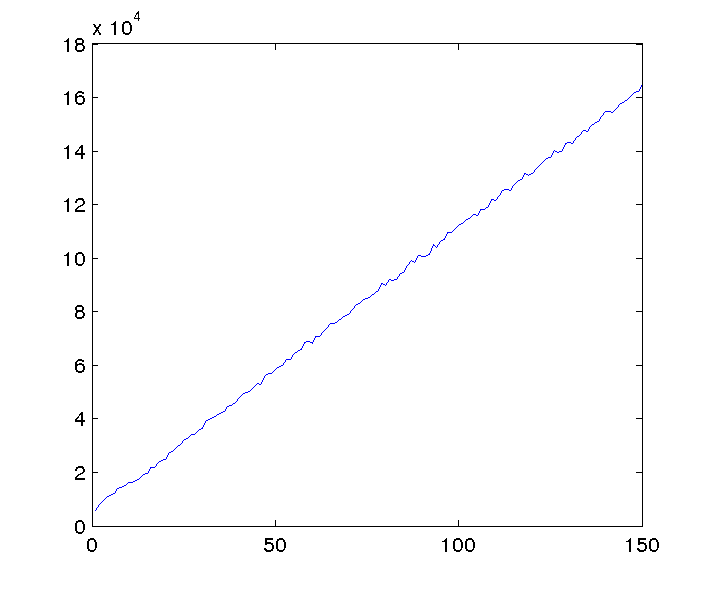}
\end{minipage}\hfill
\begin{minipage}[c]{.49\linewidth}
	\includegraphics[width=\linewidth,trim = .5cm .3cm .6cm .1cm,clip]{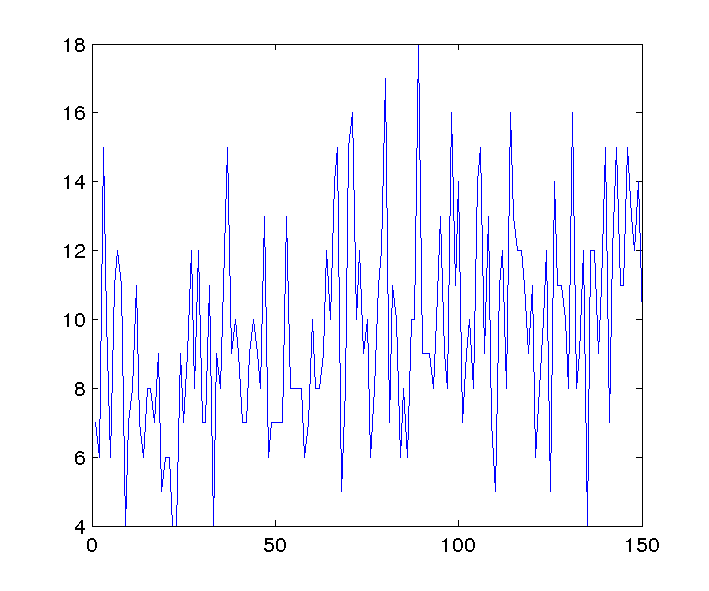}
\end{minipage}

\begin{minipage}[c]{.49\linewidth}
	\includegraphics[width=\linewidth,trim = .5cm .3cm .6cm .1cm,clip]{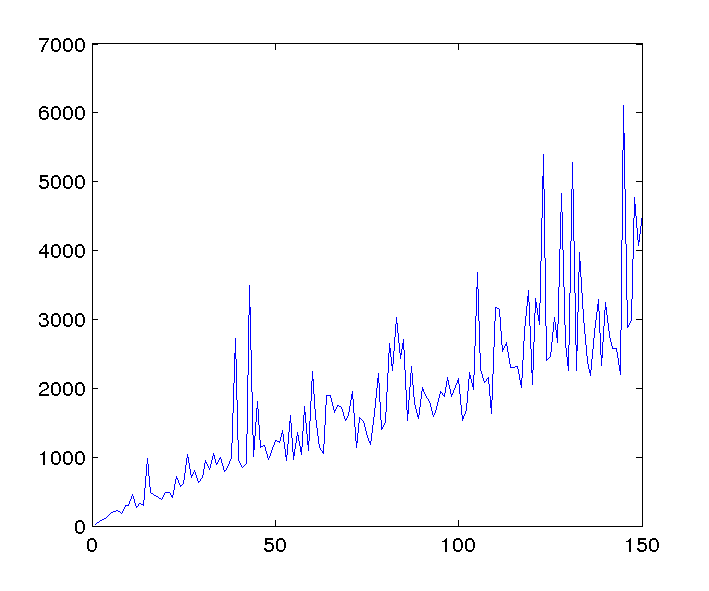}
\end{minipage}\hfill
\begin{minipage}[c]{.49\linewidth}
	\includegraphics[width=\linewidth,trim = .5cm .3cm .6cm .1cm,clip]{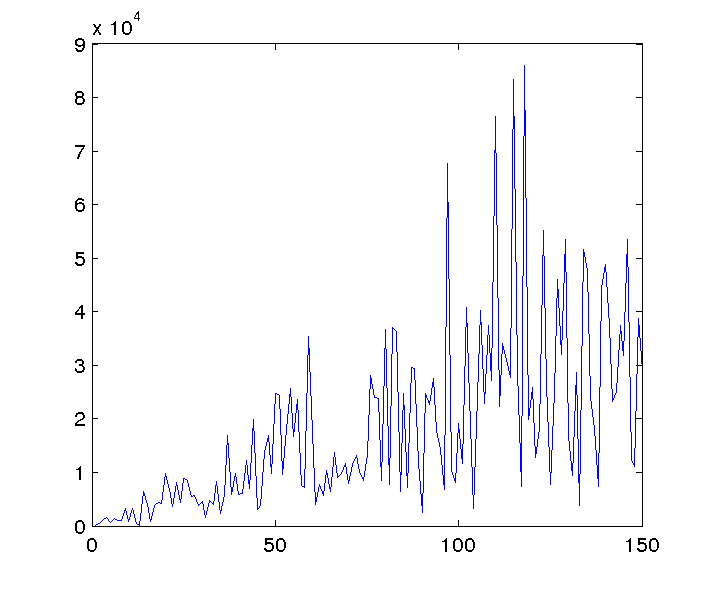}
\end{minipage}
}}\end{center}
\caption[Simulations of the combinatorial behaviour of 2 examples of conservative diffeomorphisms]{Size of the recurrent set $\Omega((f_i)_N)$ (top), number of periodic orbits of $(f_i)_N$ (middle) and length of the largest periodic orbit of $(f_i)_N$ (bottom) depending on $N$, for $f_5$ (left) and $f_6$ (right), on the grids $E_N$ with $N=128k$, $k=1,\cdots,150$.}\label{GrafConsDiffeo}
\end{figure}

First of all, we present the results of the simulations of the size of the recurrent set $\Omega((f_i)_N)$, the number of periodic orbits of $(f_i)_N$ (middle) and the length of the largest periodic orbit of $(f_i)_N$, for $N=128k$ and $k$ going from $1$ to $150$ (Figure~\ref{GrafConsDiffeo}).

For $f_5$, the cardinality of the recurrent set is clearly increasing; it looks as if it behaves like $N\ln N$. Fact remains that this behaviour is very different from the one we observe for the simulations of the homeomorphisms $f_3$ and $f_4$ (which are conservative homeomorphisms, which are respectively small $C^0$ perturbations of $\Id$ and $A$, see Section~\ref{partietroisb}, and more precisely Figure~\ref{GrafCons}). For its part, the evolution of the cardinality of the recurrent set of $f_6$ is much more erratic, and is quite similar to that observed in the corresponding $C^0$ case (top right of Figure~\ref{GrafCons}). Anyway, the behaviours are completely different at the neighbourhood of the identity and of the linear automorphism $A$.

As for the cardinality of the recurrent set, the number of periodic orbits of $(f_5)_N$ is very smooth and seems to evolve linearly in $N$. Obviously, this behaviour is very different from the case of homeomorphisms; however, we have no explanation to this ``smooth'' behaviour. In the neighbourhood of $A$, the number of periodic orbits (middle right of Figure~\ref{GrafConsDiffeo}) seems to be uniformly bounded in $N$ by $18$. Again, this behaviour is very similar to what happens to $f_4$, which is a small $C^0$ perturbation of $A$ (Figure~\ref{GrafCons}). For now, we do not have explanation to this behaviour.

The maximal period of a periodic orbit of $(f_5)_N$ evolves less smoothly then the number of periodic orbits. It is even quite similar to the corresponding $C^0$ case (Figure~\ref{GrafCons}). The same holds in the neighbourhood of $A$ for the diffeomorphism $f_6$. Seemingly, the maximal period of a periodic orbit of $f_N$ does not give a good criterion to test if a map behaves like a $C^1$ system or not.
\bigskip

\begin{figure}[h!]
\begin{center}
\makebox[0.8\textwidth]{\parbox{0.8\textwidth}{%
\begin{minipage}[c]{.49\linewidth}
	\includegraphics[width=\linewidth,trim = .5cm .3cm .6cm .1cm,clip]{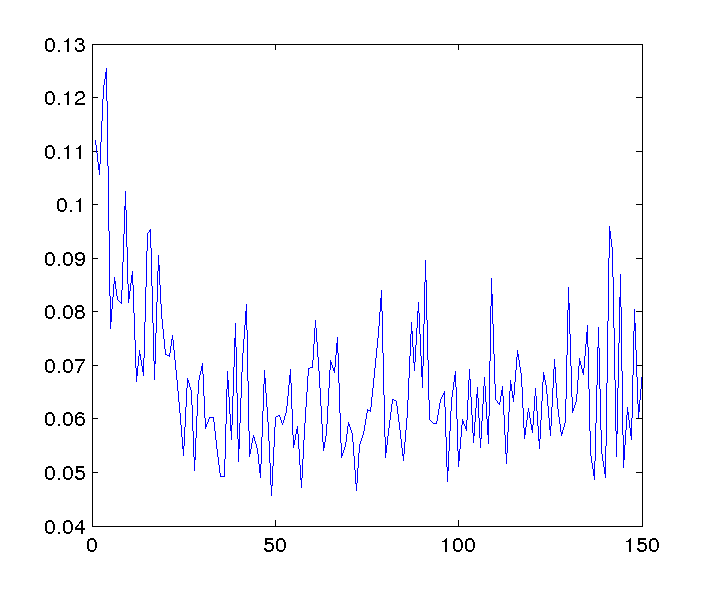}
\end{minipage}\hfill
\begin{minipage}[c]{.49\linewidth}
	\includegraphics[width=\linewidth,trim = .5cm .3cm .6cm .1cm,clip]{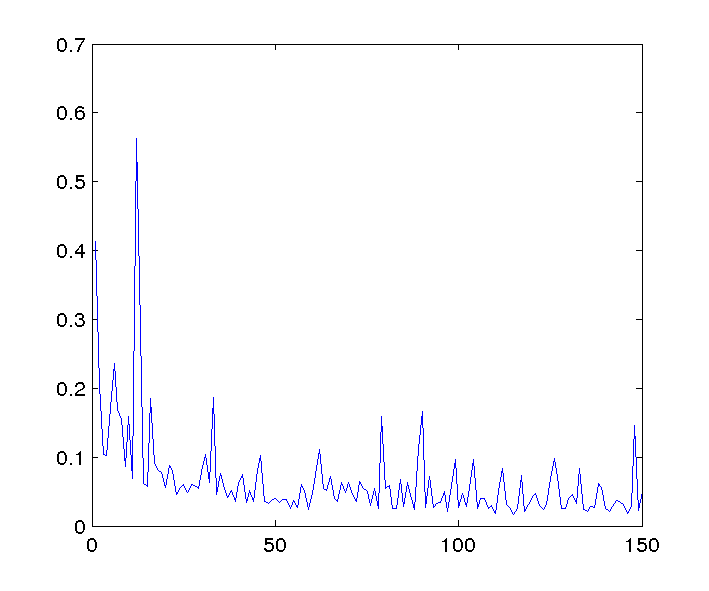}
\end{minipage}
}}\end{center}
\caption[Simulation of the distance between $\Leb$ and $\mu^{(f_i)_N}_{\T^2}$ for 3 examples of conservative diffeomorphisms]{Distance between Lebesgue measure and the measure $\mu^{(f_i)_N}_{\T^2}$ depending on $N$ for $f_5$ (left) and $f_6$ (right), on the grids $E_N$ with $N=128k$, $k=1,\cdots,150$.}\label{GrafDistLebC1Cons}
\end{figure}

For $f_5$, the distance $d(\mu^{f_N}_{\T^2},\Leb)$ is quite quickly smaller than $0.1$, and oscillates between $0.05$ and $0.1$ from $N=128\times 30$. It is not clear if in this case, the sequence of measures $(\mu^{f_N}_{\T^2})_N$ converge towards Lebesgue measure or not (while for the $C^0$ perturbation of the identity we have tested, it is clear that these measures do not converge to anything, see Figure~\ref{GrafDistLebCons}). The distance $d(\mu^{f_N}_{\T^2},\Leb)$ even seem to increase slowly (in average -- there are a lot of oscillations) from $N=50\times 128$. We have the same kind of conclusion for $f_6$: by looking at Figure~\ref{GrafDistLebC1Cons}, we can not decide if the sequence of measures $(\mu^{f_N}_{\T^2})$ seem to tend to Lebesgue measure or not.
\bigskip

\begin{figure}[ht]
\begin{minipage}[c]{.31\linewidth}
	\includegraphics[height=4.8cm,trim = 1.5cm .95cm 2.8cm .5cm,clip]{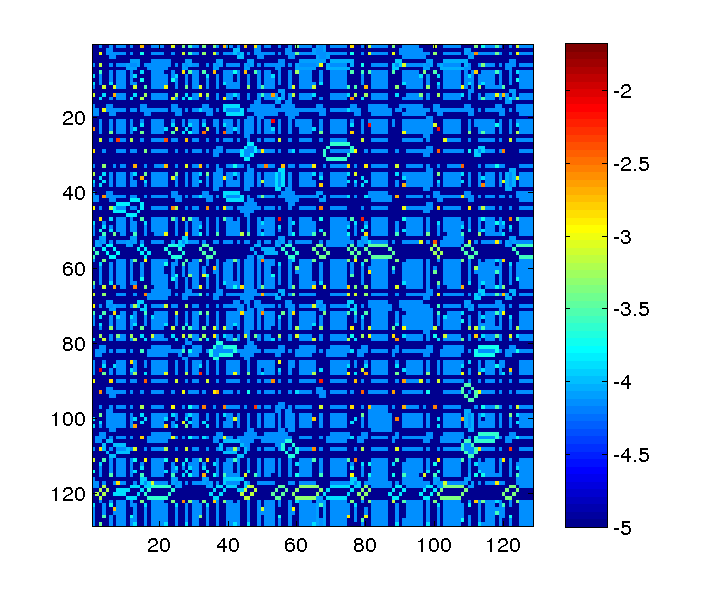}
\end{minipage}\hfill
\begin{minipage}[c]{.31\linewidth}
	\includegraphics[height=4.8cm,trim = 1.5cm .95cm 2.8cm .5cm,clip]{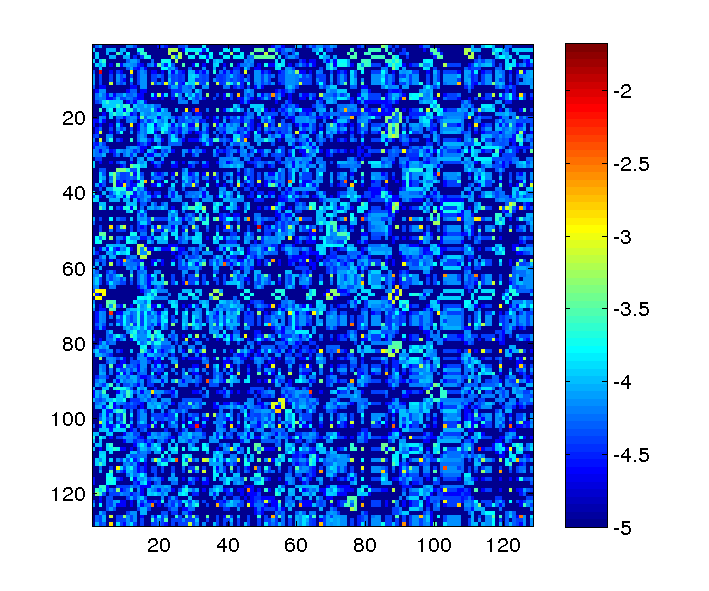}
\end{minipage}\hfill
\begin{minipage}[c]{.37\linewidth}
	\includegraphics[height=4.8cm,trim = 1.5cm .95cm 1cm .5cm,clip]{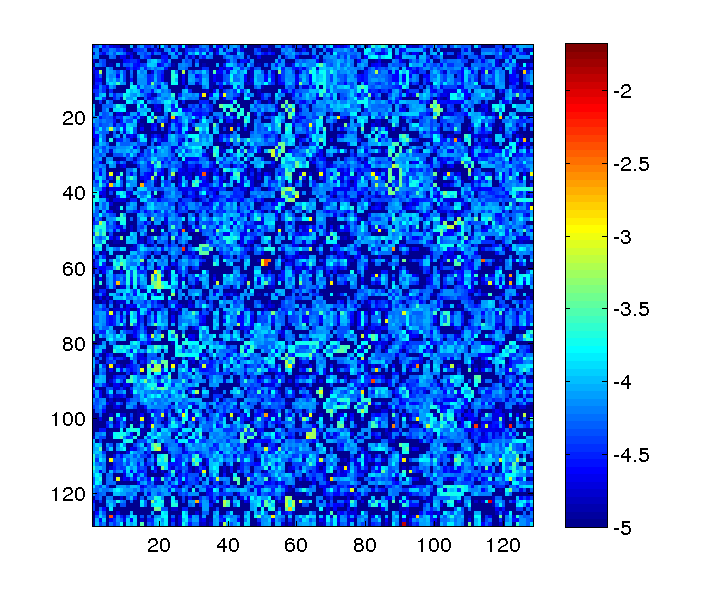}
\end{minipage}

\begin{minipage}[c]{.31\linewidth}
	\includegraphics[height=4.8cm,trim = 1.5cm .95cm 2.8cm .5cm,clip]{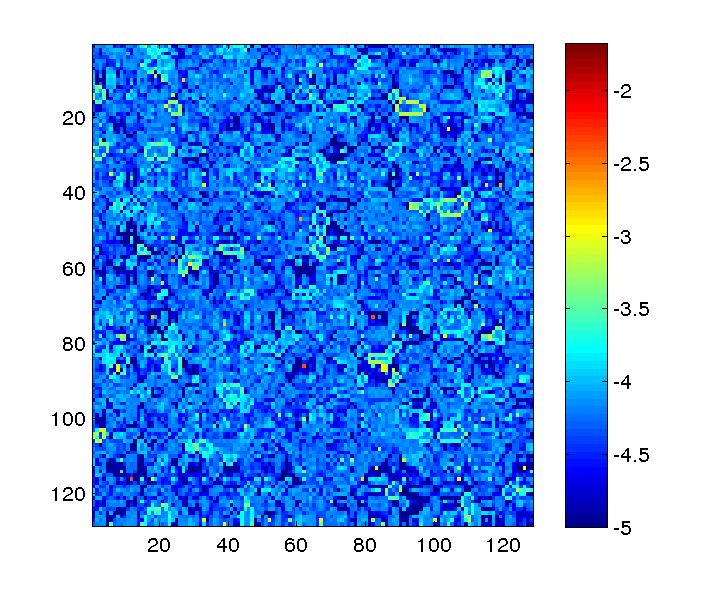}
\end{minipage}\hfill
\begin{minipage}[c]{.31\linewidth}
	\includegraphics[height=4.8cm,trim = 1.5cm .95cm 2.8cm .5cm,clip]{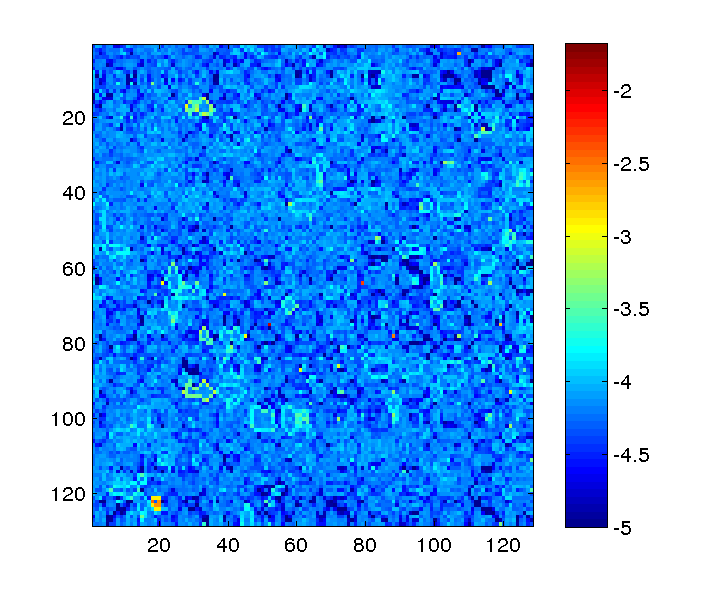}
\end{minipage}\hfill
\begin{minipage}[c]{.37\linewidth}
	\includegraphics[height=4.8cm,trim = 1.5cm .95cm 1cm .5cm,clip]{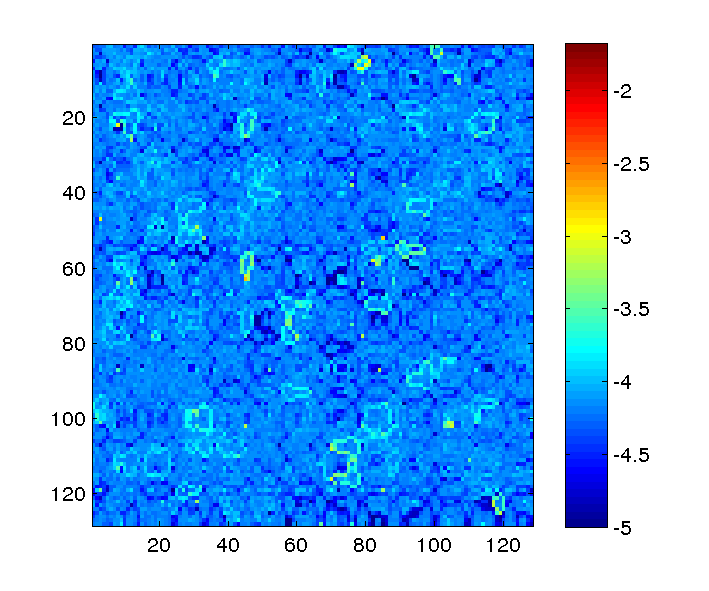}
\end{minipage}

\begin{minipage}[c]{.31\linewidth}
	\includegraphics[height=4.8cm,trim = 1.5cm .95cm 2.8cm .5cm,clip]{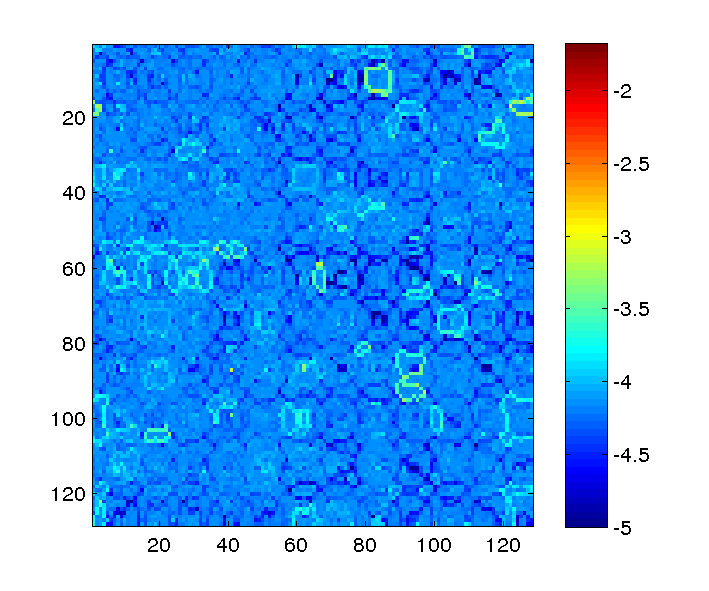}
\end{minipage}\hfill
\begin{minipage}[c]{.31\linewidth}
	\includegraphics[height=4.8cm,trim = 1.5cm .95cm 2.8cm .5cm,clip]{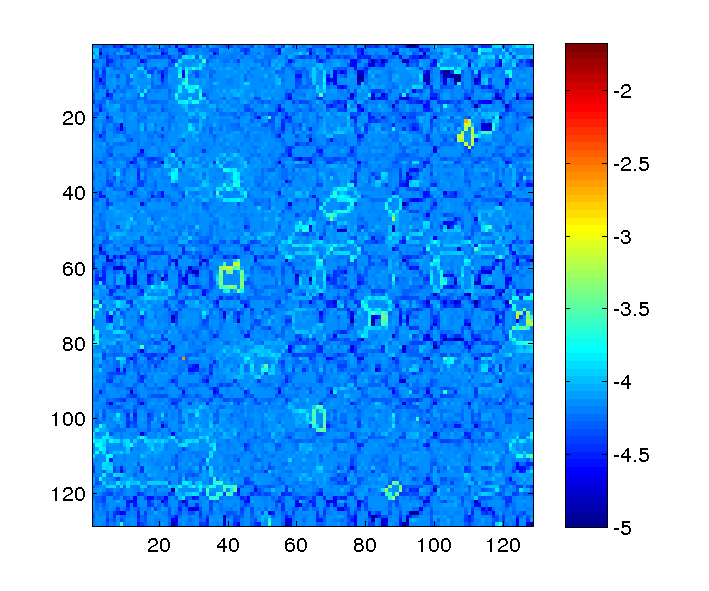}
\end{minipage}\hfill
\begin{minipage}[c]{.37\linewidth}
	\includegraphics[height=4.8cm,trim = 1.5cm .95cm 1cm .5cm,clip]{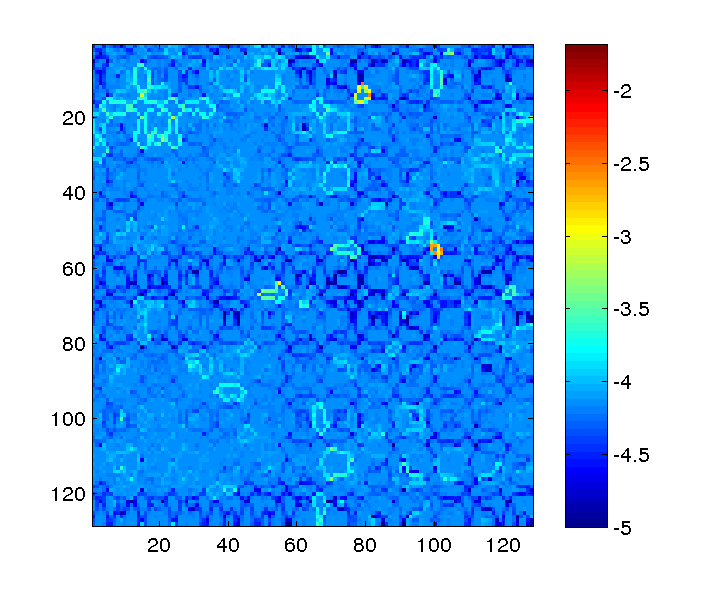}
\end{minipage}
\caption[Simulations of $\mu^{(f_5)_N}_{\T^2}$ on the grids $E_N$, with $N=2^k$, $k= 7,\cdots,15$]{Simulations of invariant measures $\mu^{(f_5)_N}_{\T^2}$ on the grids $E_N$, with $N=2^k$, $k= 7,\cdots,15$ (from left to right and top to bottom).}\label{MesC1IdCons2p}
\end{figure}

The behaviour of the computed invariant measures $\mu^{(f_5)_N}_{\T^2}$, where $f_5$ is a small $C^1$ perturbation of the identity, is way smoother than in the $C^0$ case (compare Figure~\ref{MesC1IdCons2p} with Figure~\ref{MesC0IdCons2p}). Indeed, the measure $\mu^{(f_5)_N}_{\T^2}$ has quickly a big component which is close to Lebesgue measure: the images contain a lot of light blue. Thus, we could be tempted to conclude that these measures converge to Lebesgue measure. However, there are still little regions that have a big measure: in the example of Figure~\ref{MesC1IdCons2p}, it seems that there are invariant curves that attract a lot of the points of the grid (as can also be observed on Figure~\ref{MesPhysIdC1}). We have no explanation to this phenomenon, and we do not know if it still occurs when the order of discretization is very large.
\bigskip

\begin{figure}[ht]
\begin{minipage}[c]{.31\linewidth}
	\includegraphics[height=4.8cm,trim = 1.5cm .95cm 2.8cm .5cm,clip]{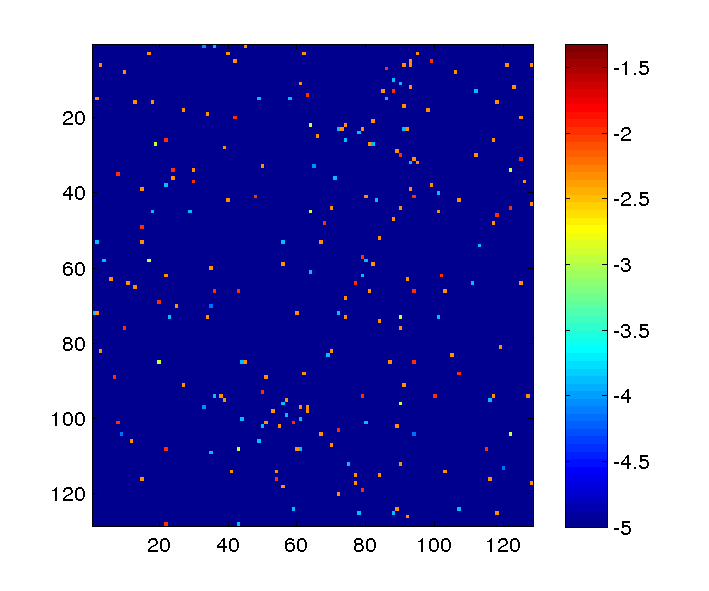}
\end{minipage}\hfill
\begin{minipage}[c]{.31\linewidth}
	\includegraphics[height=4.8cm,trim = 1.5cm .95cm 2.8cm .5cm,clip]{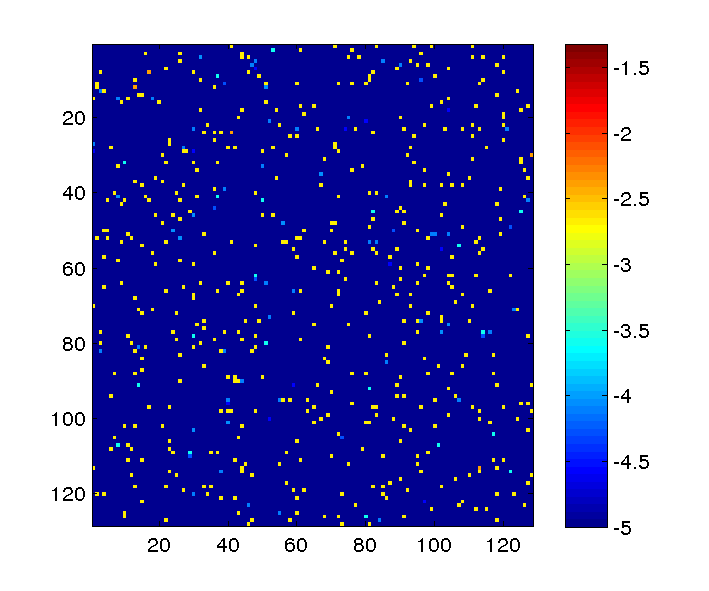}
\end{minipage}\hfill
\begin{minipage}[c]{.37\linewidth}
	\includegraphics[height=4.8cm,trim = 1.5cm .95cm 1cm .5cm,clip]{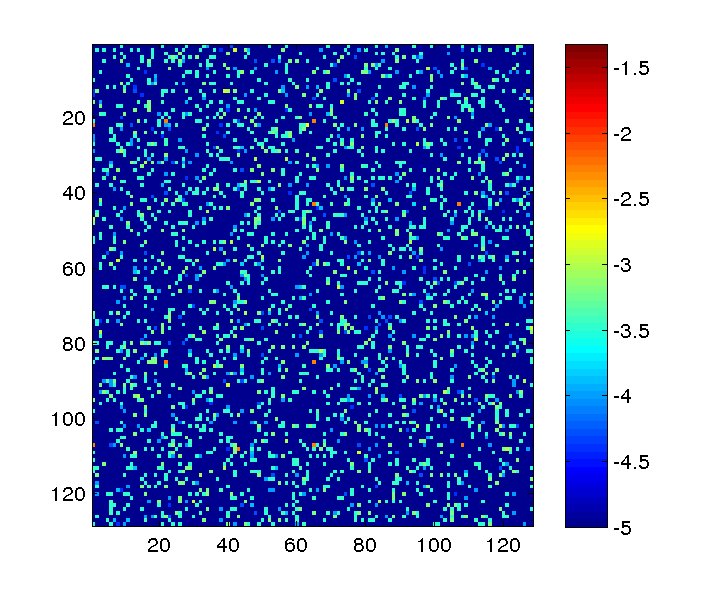}
\end{minipage}

\begin{minipage}[c]{.31\linewidth}
	\includegraphics[height=4.8cm,trim = 1.5cm .95cm 2.8cm .5cm,clip]{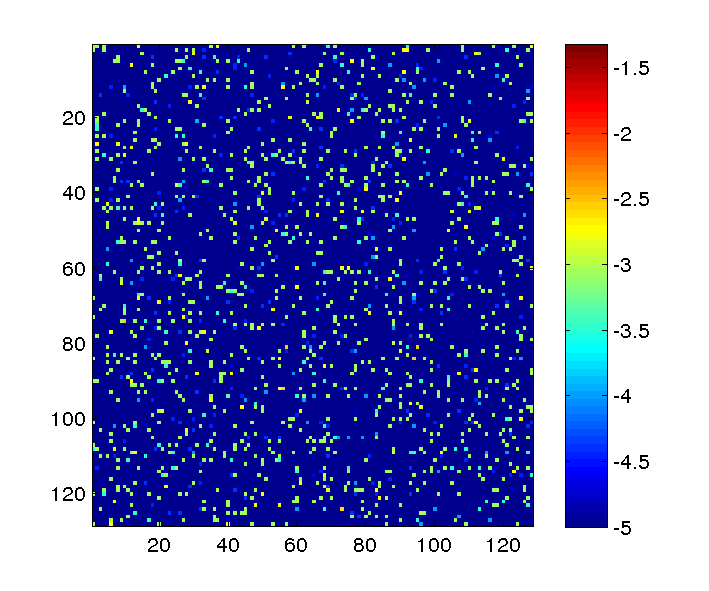}
\end{minipage}\hfill
\begin{minipage}[c]{.31\linewidth}
	\includegraphics[height=4.8cm,trim = 1.5cm .95cm 2.8cm .5cm,clip]{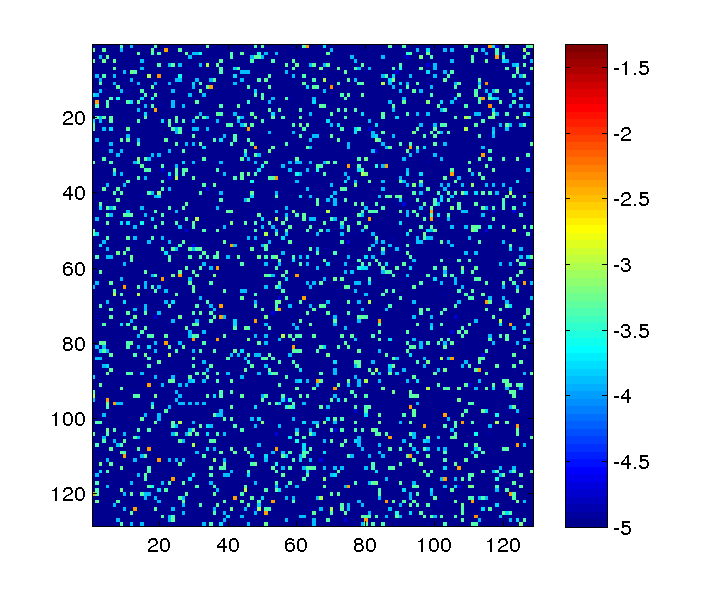}
\end{minipage}\hfill
\begin{minipage}[c]{.37\linewidth}
	\includegraphics[height=4.8cm,trim = 1.5cm .95cm 1cm .5cm,clip]{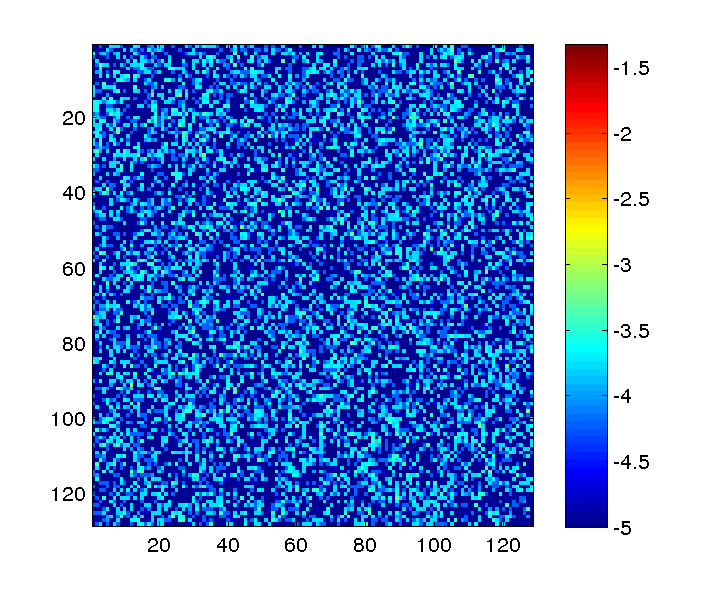}
\end{minipage}

\begin{minipage}[c]{.31\linewidth}
	\includegraphics[height=4.8cm,trim = 1.5cm .95cm 2.8cm .5cm,clip]{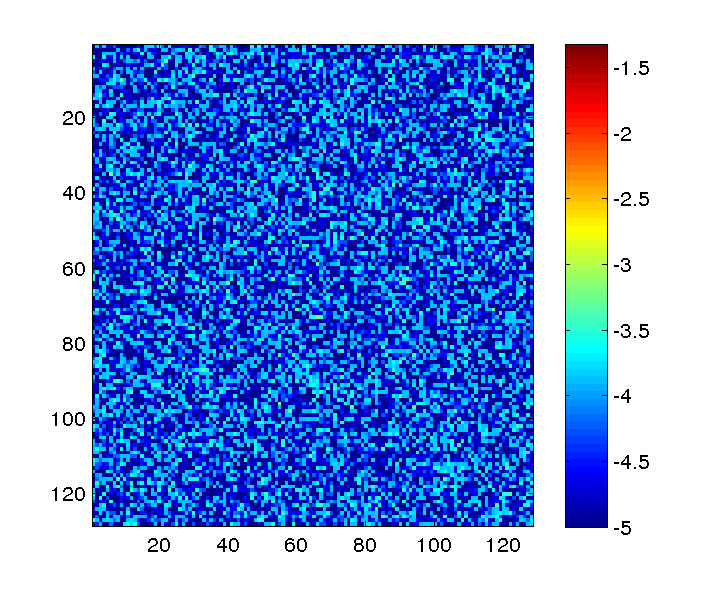}
\end{minipage}\hfill
\begin{minipage}[c]{.31\linewidth}
	\includegraphics[height=4.8cm,trim = 1.5cm .95cm 2.8cm .5cm,clip]{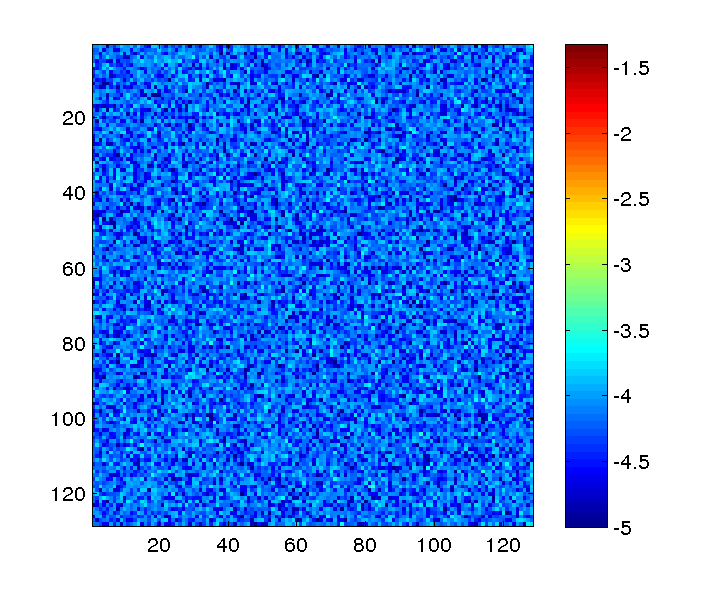}
\end{minipage}\hfill
\begin{minipage}[c]{.37\linewidth}
	\includegraphics[height=4.8cm,trim = 1.5cm .95cm 1cm .5cm,clip]{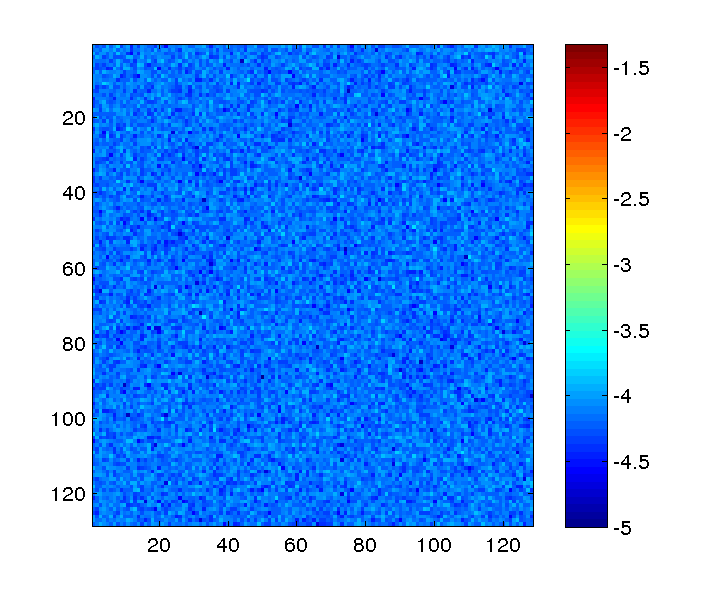}
\end{minipage}
\caption[Simulations of $\mu^{(f_6)_N}_{\T^2}$ on the grids $E_N$, with $N=2^k$, $k= 7,\cdots,15$]{Simulations of invariant measures $\mu^{(f_6)_N}_{\T^2}$ on the grids $E_N$, with $N=2^k$, $k= 7,\cdots,15$ (from left to right and top to bottom).}\label{MesC1AnoCons2p}
\end{figure}

\begin{figure}[ht]
\begin{minipage}[c]{.31\linewidth}
	\includegraphics[height=4.8cm,trim = 1.5cm .95cm 2.8cm .5cm,clip]{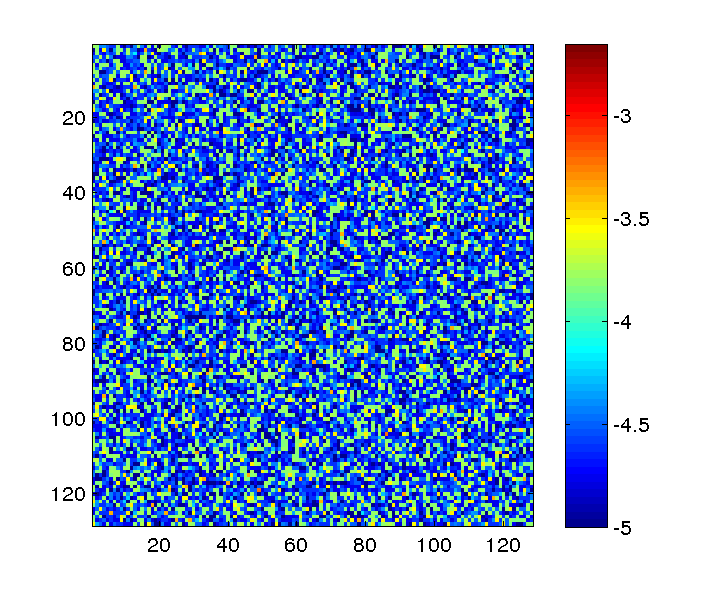}
\end{minipage}\hfill
\begin{minipage}[c]{.31\linewidth}
	\includegraphics[height=4.8cm,trim = 1.5cm .95cm 2.8cm .5cm,clip]{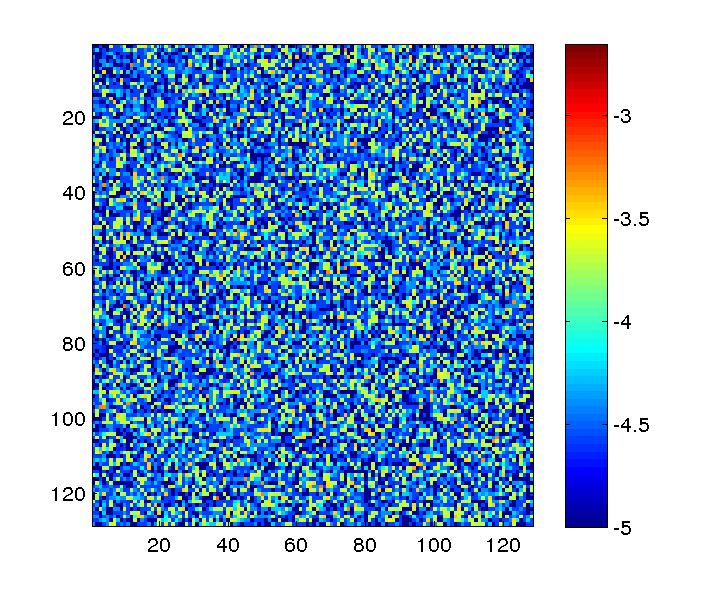}
\end{minipage}\hfill
\begin{minipage}[c]{.37\linewidth}
	\includegraphics[height=4.8cm,trim = 1.5cm .95cm 1cm .5cm,clip]{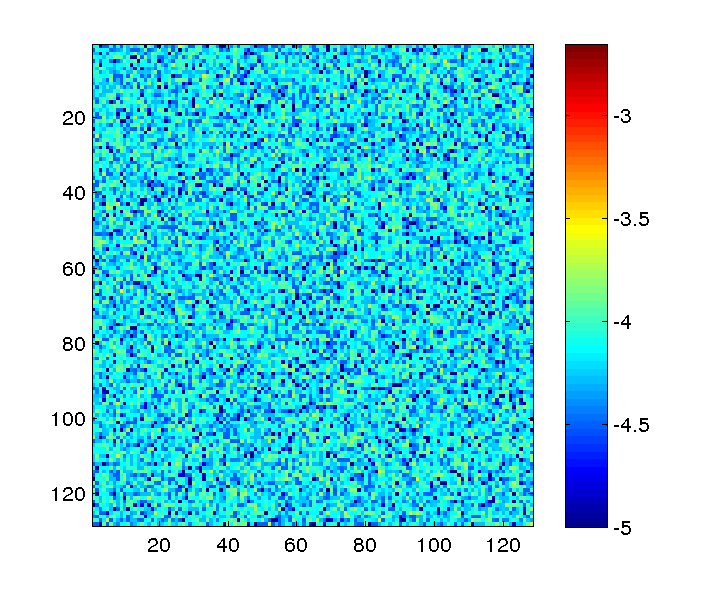}
\end{minipage}

\begin{minipage}[c]{.31\linewidth}
	\includegraphics[height=4.8cm,trim = 1.5cm .95cm 2.8cm .5cm,clip]{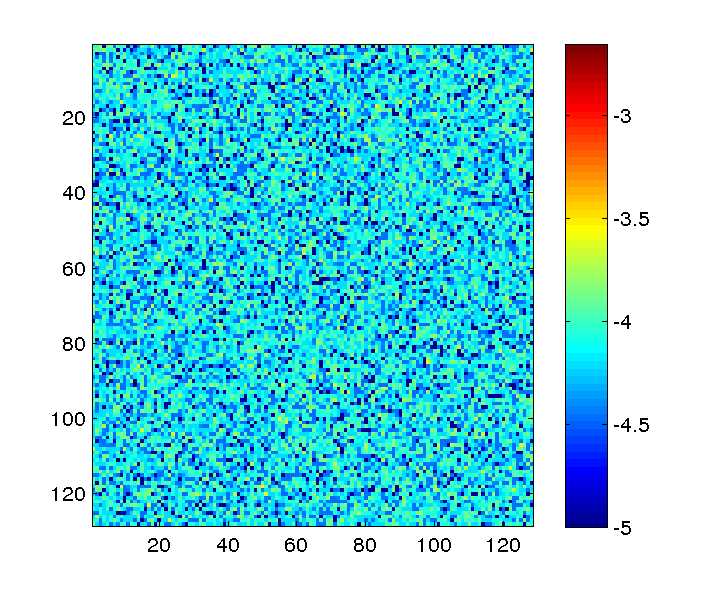}
\end{minipage}\hfill
\begin{minipage}[c]{.31\linewidth}
	\includegraphics[height=4.8cm,trim = 1.5cm .95cm 2.8cm .5cm,clip]{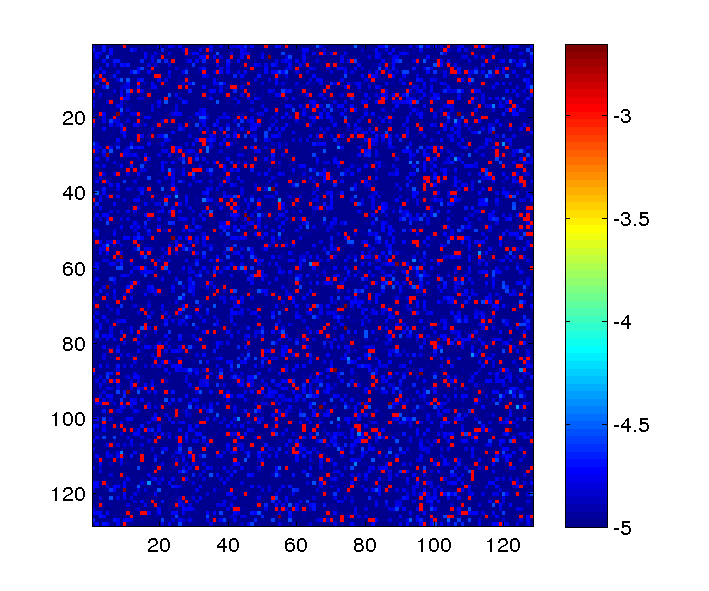}
\end{minipage}\hfill
\begin{minipage}[c]{.37\linewidth}
	\includegraphics[height=4.8cm,trim = 1.5cm .95cm 1cm .5cm,clip]{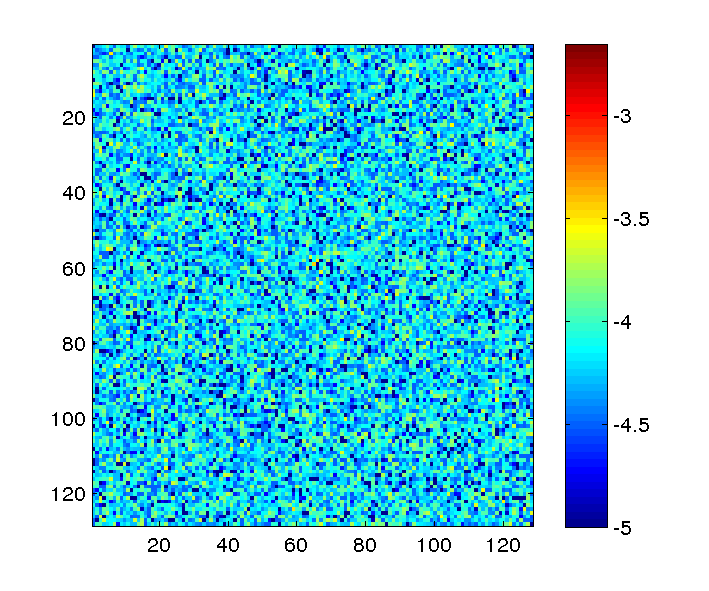}
\end{minipage}

\begin{minipage}[c]{.31\linewidth}
	\includegraphics[height=4.8cm,trim = 1.5cm .95cm 2.8cm .5cm,clip]{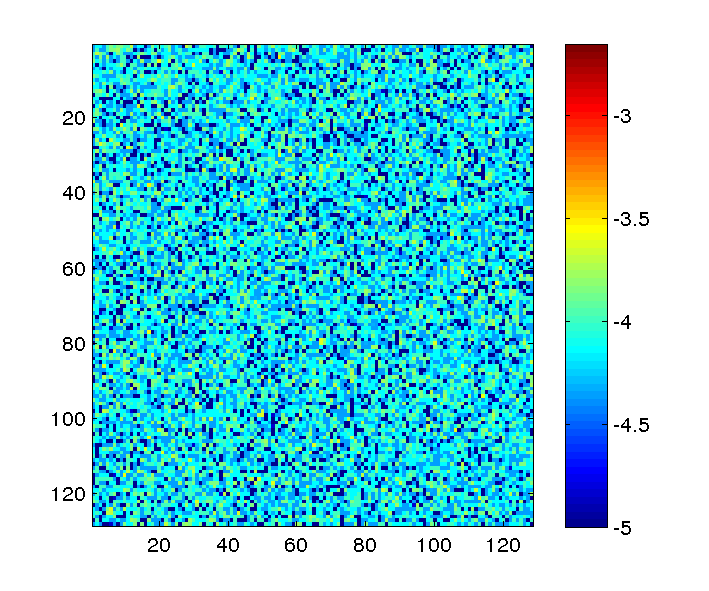}
\end{minipage}\hfill
\begin{minipage}[c]{.31\linewidth}
	\includegraphics[height=4.8cm,trim = 1.5cm .95cm 2.8cm .5cm,clip]{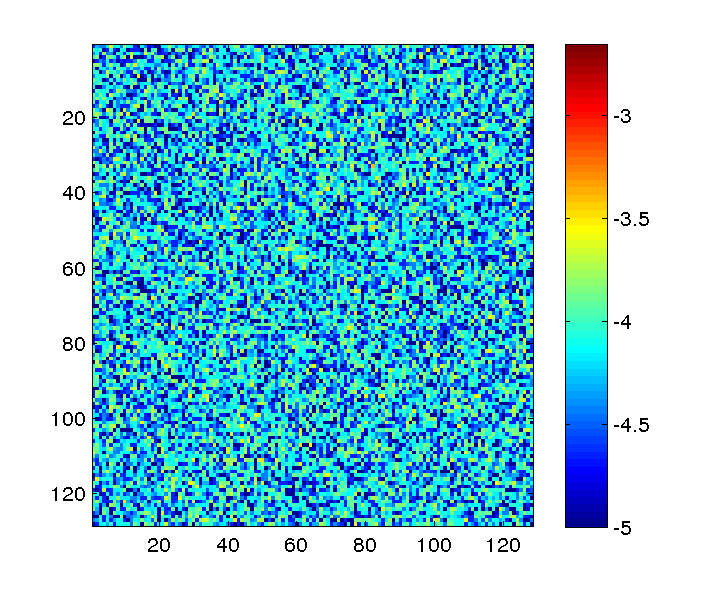}
\end{minipage}\hfill
\begin{minipage}[c]{.37\linewidth}
	\includegraphics[height=4.8cm,trim = 1.5cm .95cm 1cm .5cm,clip]{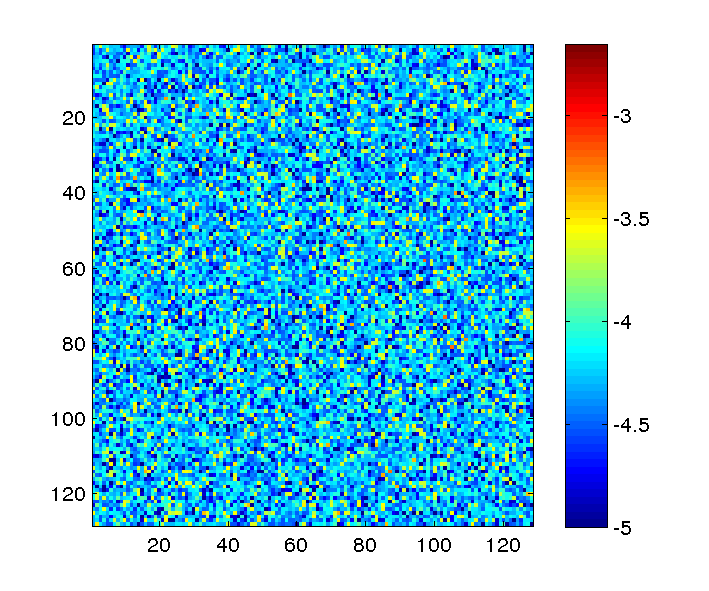}
\end{minipage}
\caption[Simulations of $\mu^{(f_6)_N}_{\T^2}$ on the grids $E_N$, with $N=11\,516, \cdots, 11\,524$]{Simulations of invariant measures $\mu^{(f_6)_N}_{\T^2}$ on the grids $E_N$, with $N=11\,516, \cdots, 11\,524$ (from left to right and top to bottom).}\label{MesC1AnoConsSerie}
\end{figure}

For the discretizations of $f_6$, the simulations on grids of size $2^k\times 2^k$ might suggest that the measures $\mu^{(f_6)_N}_{\T^2}$ tend to Lebesgue measure (Figure~\ref{MesC1AnoCons2p}). Actually, when we perform a lot of simulations, we realize that there are also big variations of the behaviour of the measures (Figure \ref{MesC1AnoConsSerie}): the measure is often well distributed in the torus, and sometimes quite singular with respect to Lebesgue measure (as it can be seen in Figure \ref{GrafDistLebC1Cons}). This behaviour is almost identical to that observed in the $C^0$ case in the neighbourhood of $A$ (see Figures~\ref{MesC0AnoCons2p} and \ref{MesC0AnoConsSer}, and also the corresponding discussion).

\clearpage

\subsection{The case of expanding maps}

We now present the results of the numerical simulations we have conducted for expanding maps of the circle. We have tested the following expanding map of degree 2:
\[f(x) = 2x + \varep_1 \cos(2\pi x) + \varep_2 \sin(6\pi x),\]
with $\varep_1 = 0.127\,943\,563\,72$ and $\varep_2 = 0.008\,247\,359\,61$.
\bigskip

We focus on the simulations of both measures $(f_N^*)^k\lambda_N$ and $\mu^{f_N}_{\Sp^1}$. Recall that the former is simply the push-forward of the uniform measure on $E_N$ by the iterate $f_N^k$ of the discretization. The latter is the measure supported by the union of periodic orbits of $f_N$, such that the total measure of a periodic orbit is equal to the size of its basin of attraction.
\bigskip

\begin{figure}[ht]
\begin{minipage}[c]{.33\linewidth}
	\includegraphics[width=\linewidth,trim = .5cm .3cm .6cm .1cm,clip]{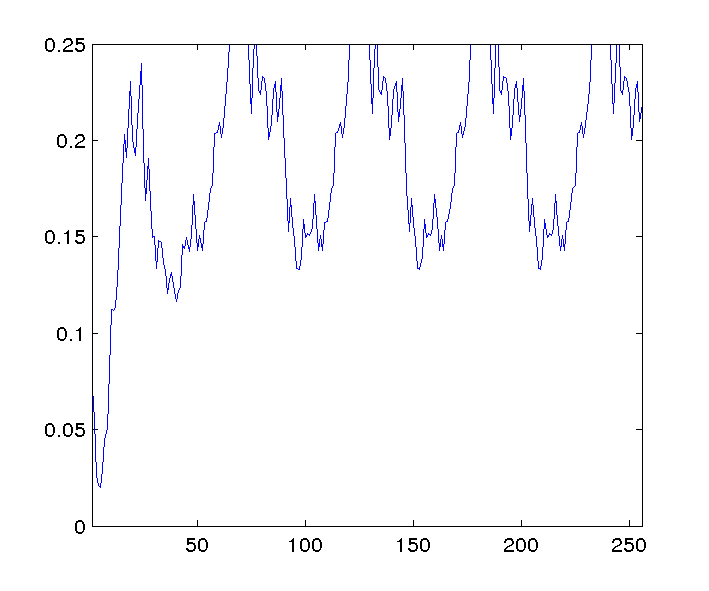}
\end{minipage}\hfill
\begin{minipage}[c]{.33\linewidth}
	\includegraphics[width=\linewidth,trim = .5cm .3cm .6cm .1cm,clip]{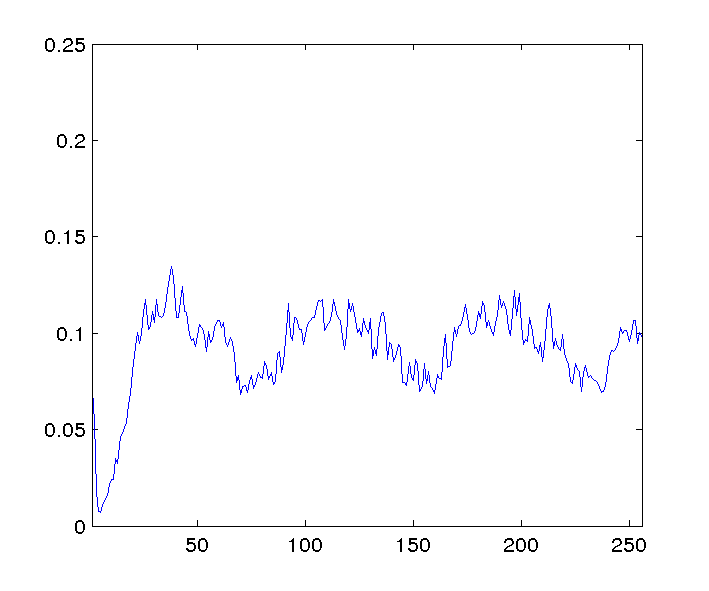}
\end{minipage}\hfill
\begin{minipage}[c]{.33\linewidth}
	\includegraphics[width=\linewidth,trim = .5cm .3cm .6cm .1cm,clip]{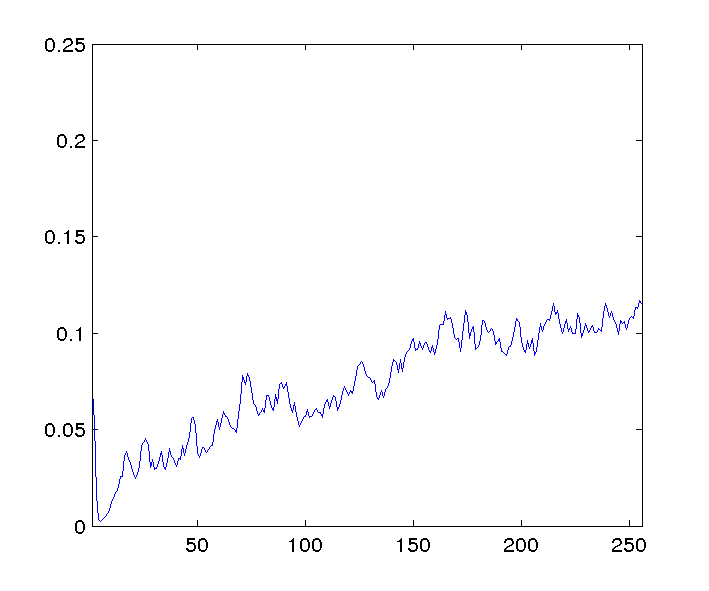}
\end{minipage}

\begin{minipage}[c]{.33\linewidth}
	\includegraphics[width=\linewidth,trim = .5cm .3cm .6cm .1cm,clip]{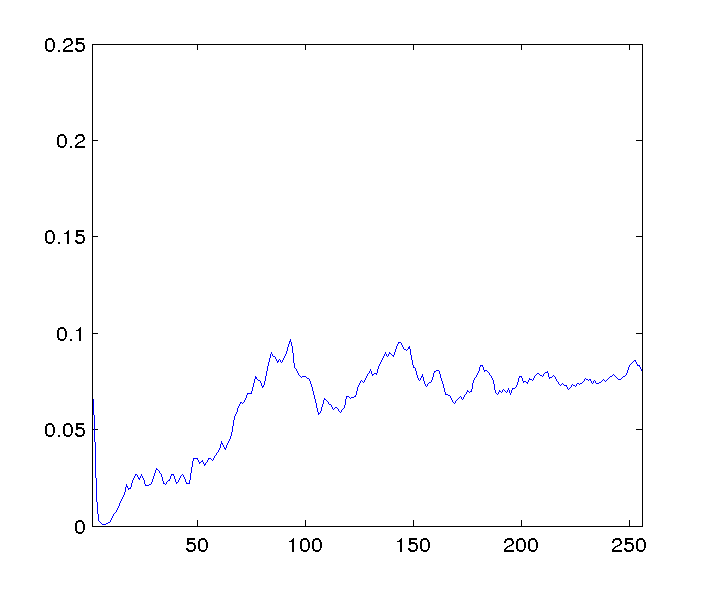}
\end{minipage}\hfill
\begin{minipage}[c]{.33\linewidth}
	\includegraphics[width=\linewidth,trim = .5cm .3cm .6cm .1cm,clip]{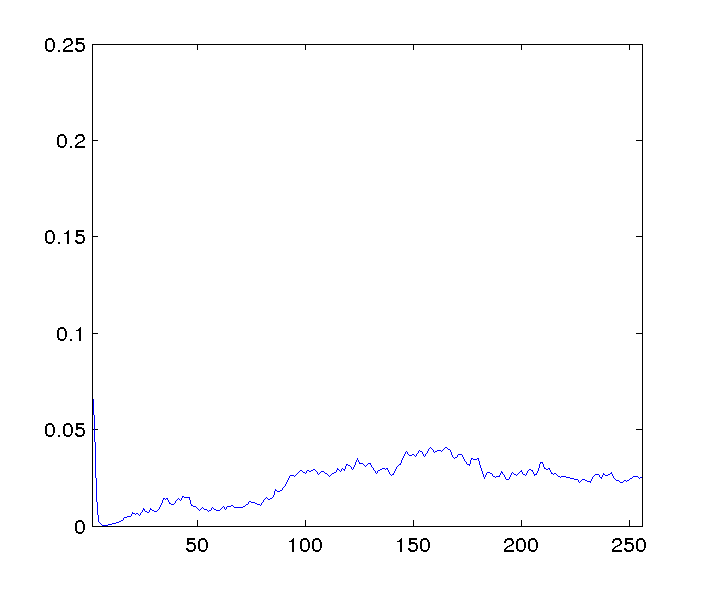}
\end{minipage}\hfill
\begin{minipage}[c]{.33\linewidth}
	\includegraphics[width=\linewidth,trim = .5cm .3cm .6cm .1cm,clip]{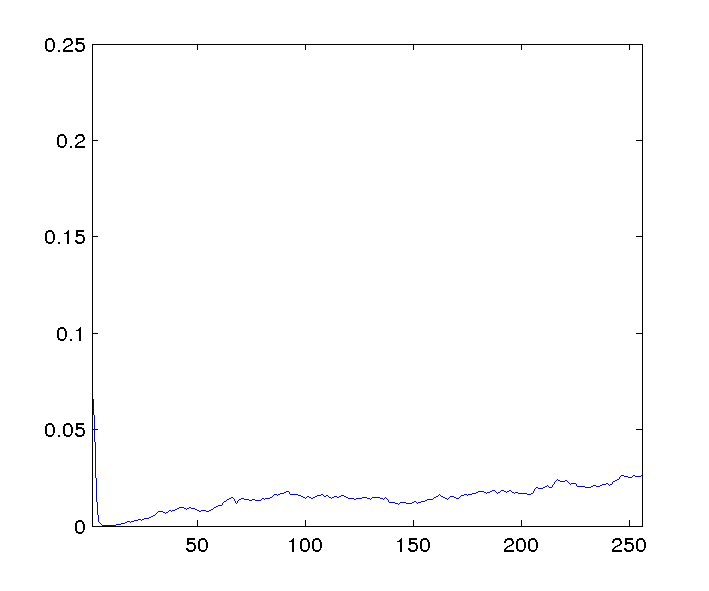}
\end{minipage}
\caption[Distance between the SRB measure of $f$ and the measure $(f_N^*)^k\lambda_N$]{Distance between the SRB measure of $f$ and the measure $(f_N^*)^k\lambda_N$ depending on the time $k$ ($1\le k \le 1\,024$), for $N=2^j$, with (from left to right and top to bottom) $j=9$, $11$, $13$, $15$, $17$ and $19$.}\label{GrafDistMesTemps}
\end{figure}

Recall that as $f$ is expanding and belongs to $C^2(\Sp^1)$, it has a single SRB measure, that we denote by $\mu_0$. This measure can be computed quite easily using the Ruelle-Perron-Frobenius transfer operator. In Figure~\ref{GrafDistMesTemps}, we have represented the distance $d(\mu_0,\,(f_N^*)^k\lambda_N)$ depending on the time $k$, for various orders of discretization $N$. The distance $d$ is defined by
\[d(\mu,\nu) = \sum_{k=0}^\infty \frac{1}{2^k} \sum_{i=0}^{2^k-1} \left| \mu\left(\left[\frac{i}{2^k},\frac{i+1}{2^k}\right]) - \nu(\left[\frac{i}{2^k},\frac{i+1}{2^k}\right]\right)\right|\in[0,2].\]
This distance spans the weak-* topology, which makes compact the set of probability measures on $\T^2$. In practice, we have computed an approximation of this quantity by summing only on the $k\in\llbracket 0,7 \rrbracket$.
\bigskip

\begin{figure}[ht]
\begin{center}
\includegraphics[width=.4\linewidth,trim = .5cm .3cm .6cm .1cm,clip]{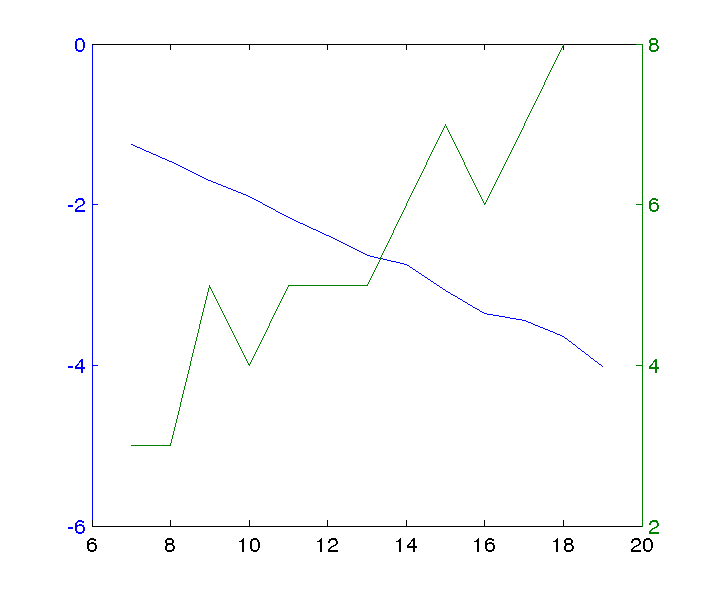}
\caption[Minimum of $d(\mu_0,\,(f_N^*)^k\lambda_N)$ depending on $N$, and time $k$ realizing this minimal distance]{Logarithm in base 10 of the minimum over $k$ of $d(\mu_0,\,(f_N^*)^k\lambda_N)$ depending on $N$ (in blue, left axis), and time $k$ realizing this minimal distance (in green, right axis), for $N=2^j$, with $j\in\llbracket 7,19\rrbracket$.}\label{minimumLog}
\end{center}
\end{figure}

On Figure~\ref{GrafDistMesTemps}, we observe that for $N$ being fixed, the distance between the SRB measure $\mu_0$ and the image measure $(f_N^*)^k\lambda_N$ reaches quite quickly a value quite close to 0, to increase thereafter. We shall notice that for the smaller values of $N$ ($N=2^7$ and $N=2^9$), the distance $d(\mu_0,\,(f_N^*)^k\lambda_N)$ looks eventually periodic: it seems that the stabilization time of the discretization is attained, or at least that for almost all the points $x\in E_N$, we have $f_N^{1\, 024}(x)\in \Omega(f_N)$.

The fact that the image measure $(f_N^*)^k\lambda_N$ is very close to the SRB measure for short times $k$ is not surprising, as it was predicted by Theorem~\ref{MainMoche}. For example, in Figure~\ref{minimumLog}, we can see that for $N=2^{19}$, the distance between $\mu_0$ and $(f_N^*)^8\lambda_N$ is approximatively $10^{-4}$

We think that growth of the distance $d(\mu_0,\,(f_N^*)^k\lambda_N)$ is increasing from a certain point is due to the fact that the roundoff errors arising from the discretization process imply that locally, the discrepancy between the image measure and the appropriate uniform measure grows exponentially with $k$ (see Proposition~\ref{EstimDiscrepancy}). However, for now, we have neither made the precise study of the discrepancy of a union of ``independent'' image sets (what happens for maps of degree $d\ge 2$), nor the application of the discrepancy to the non-linear case.

\begin{figure}[ht]
\begin{minipage}[c]{.33\linewidth}
	\includegraphics[width=\linewidth,trim = .5cm .3cm .6cm .1cm,clip]{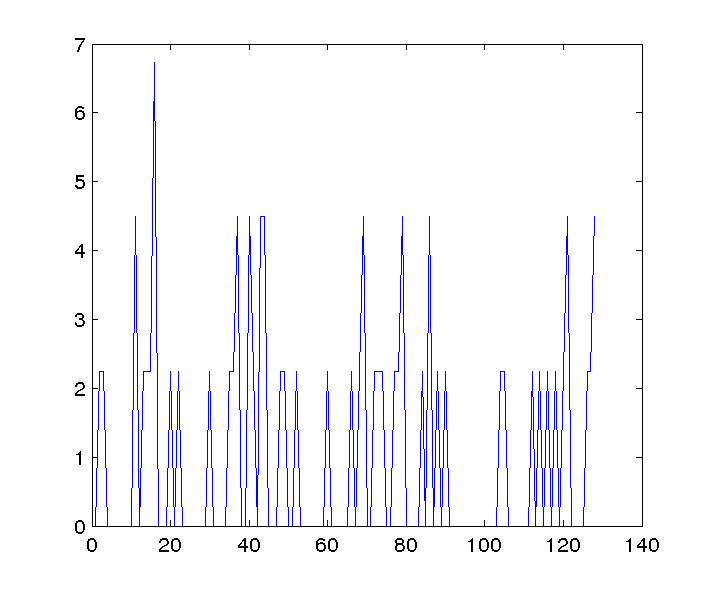}
\end{minipage}\hfill
\begin{minipage}[c]{.33\linewidth}
	\includegraphics[width=\linewidth,trim = .5cm .3cm .6cm .1cm,clip]{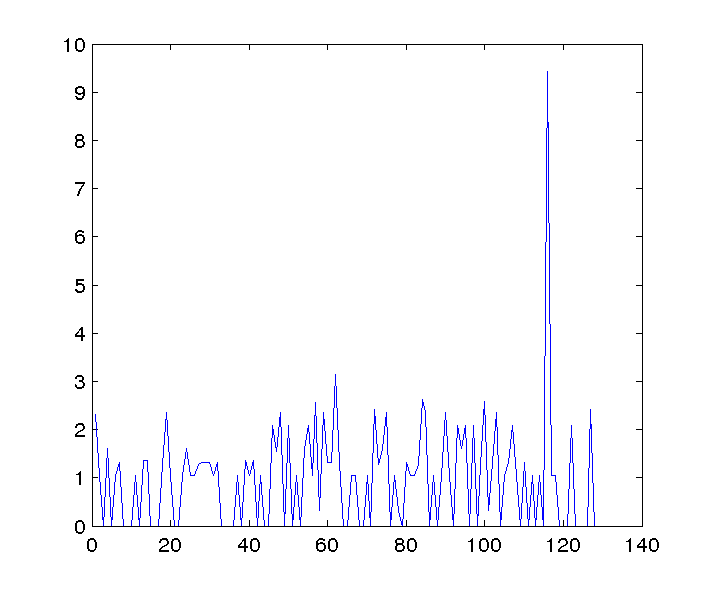}
\end{minipage}\hfill
\begin{minipage}[c]{.33\linewidth}
	\includegraphics[width=\linewidth,trim = .5cm .3cm .6cm .1cm,clip]{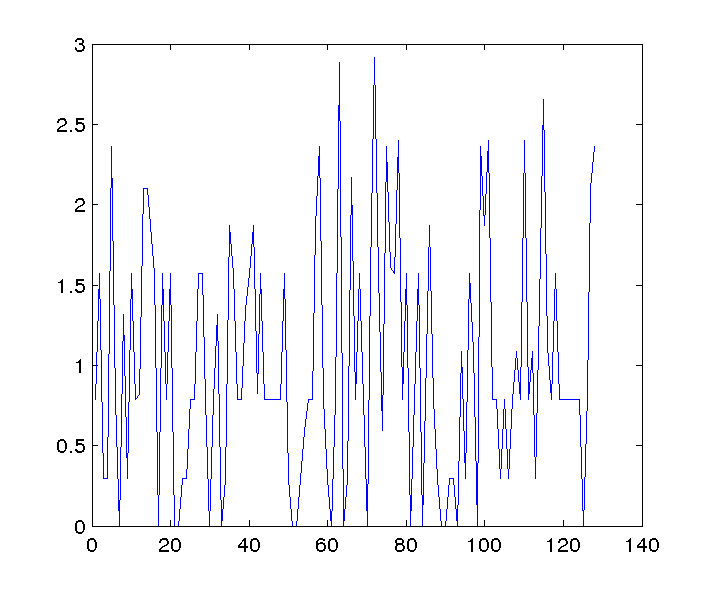}
\end{minipage}

\begin{minipage}[c]{.33\linewidth}
	\includegraphics[width=\linewidth,trim = .5cm .3cm .6cm .1cm,clip]{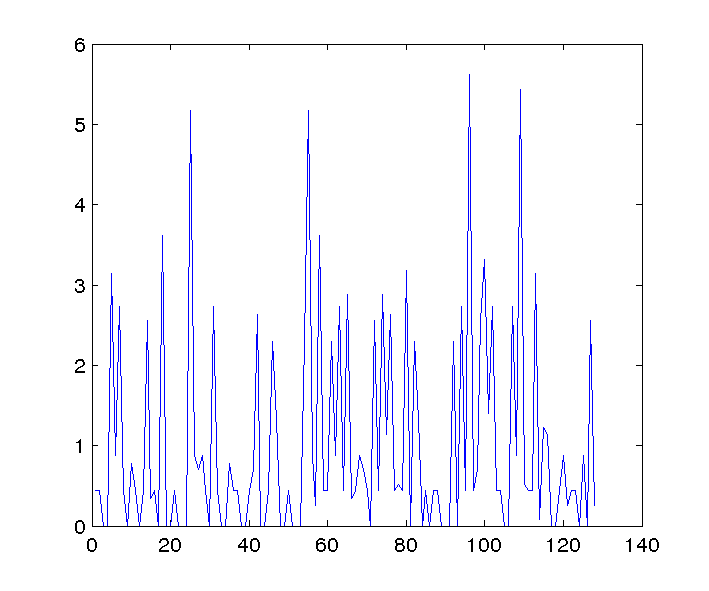}
\end{minipage}\hfill
\begin{minipage}[c]{.33\linewidth}
	\includegraphics[width=\linewidth,trim = .5cm .3cm .6cm .1cm,clip]{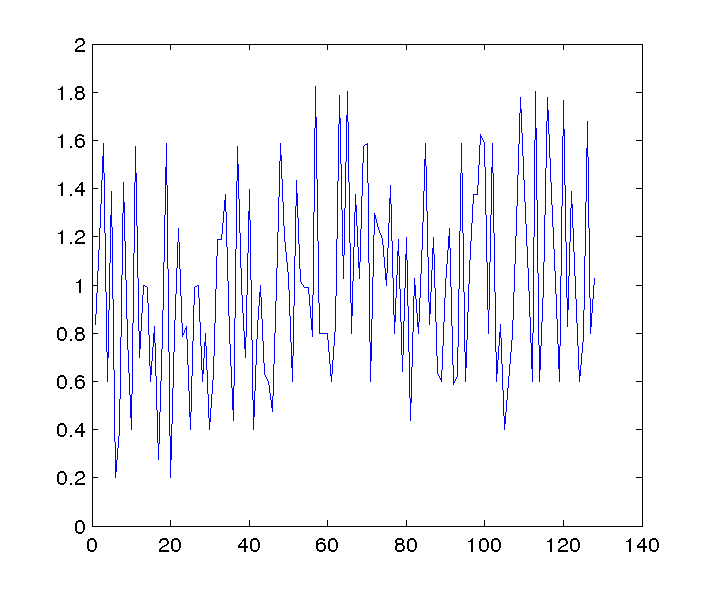}
\end{minipage}\hfill
\begin{minipage}[c]{.33\linewidth}
	\includegraphics[width=\linewidth,trim = .5cm .3cm .6cm .1cm,clip]{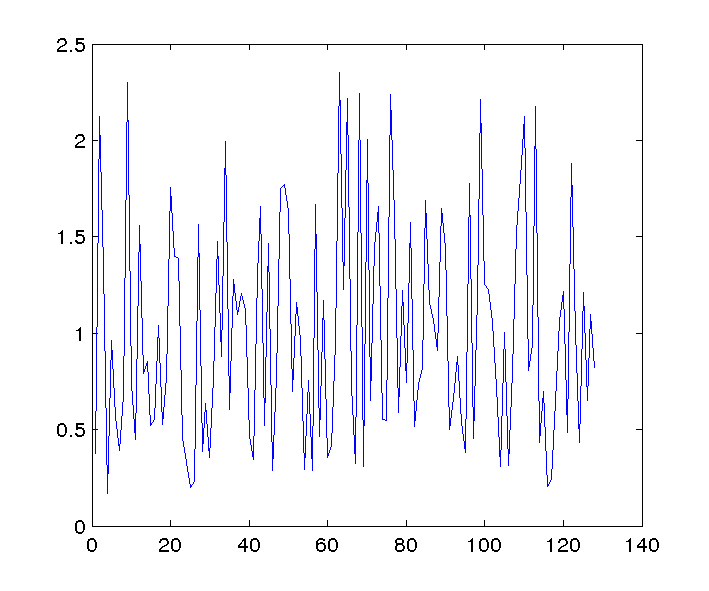}
\end{minipage}
\caption[Density of the measure $\mu^{f_N}_{\Sp^1}$ for $N=128\cdot 4^j$, $j\le5$]{Density of the measure $\mu^{f_N}_{\Sp^1}$, for $N=2^j$, with (from left to right and top to bottom) $j=7$, $9$, $11$, $13$, $15$, $17$ and $19$.}\label{GrafMestiti19}
\end{figure}

\begin{figure}[ht]
\begin{minipage}[c]{.33\linewidth}
	\includegraphics[width=\linewidth,trim = .5cm .3cm .6cm .1cm,clip]{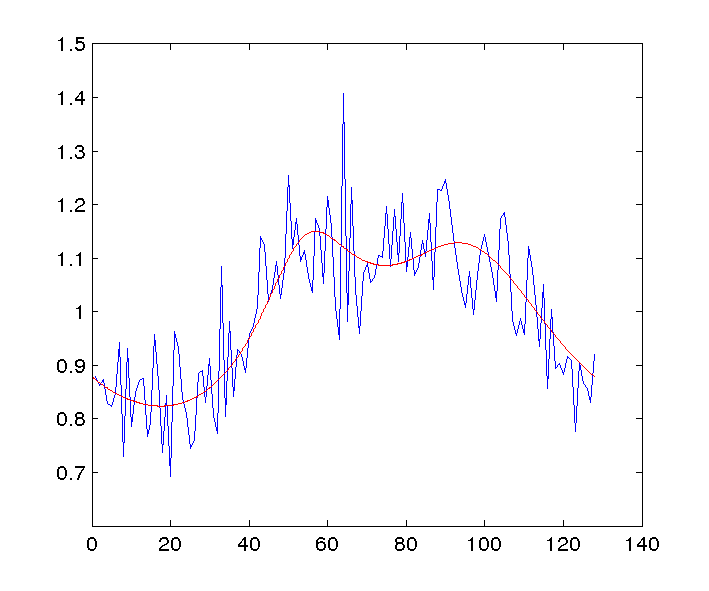}
\end{minipage}\hfill
\begin{minipage}[c]{.33\linewidth}
	\includegraphics[width=\linewidth,trim = .5cm .3cm .6cm .1cm,clip]{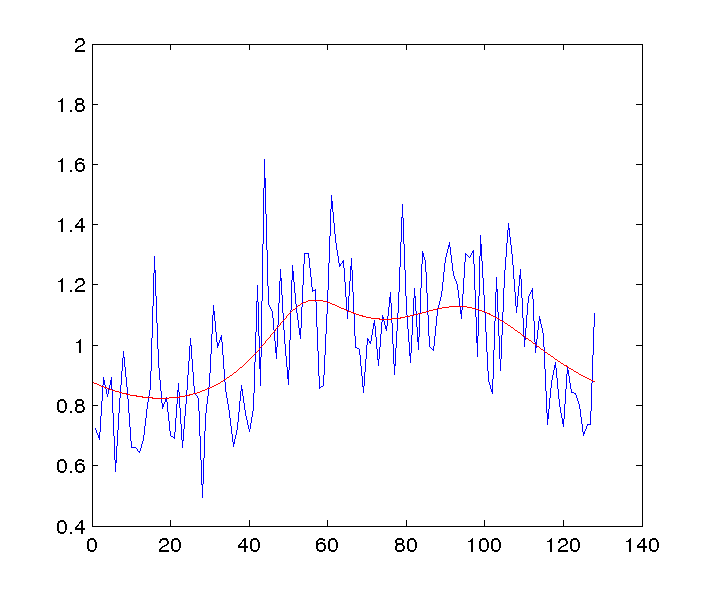}
\end{minipage}\hfill
\begin{minipage}[c]{.33\linewidth}
	\includegraphics[width=\linewidth,trim = .5cm .3cm .6cm .1cm,clip]{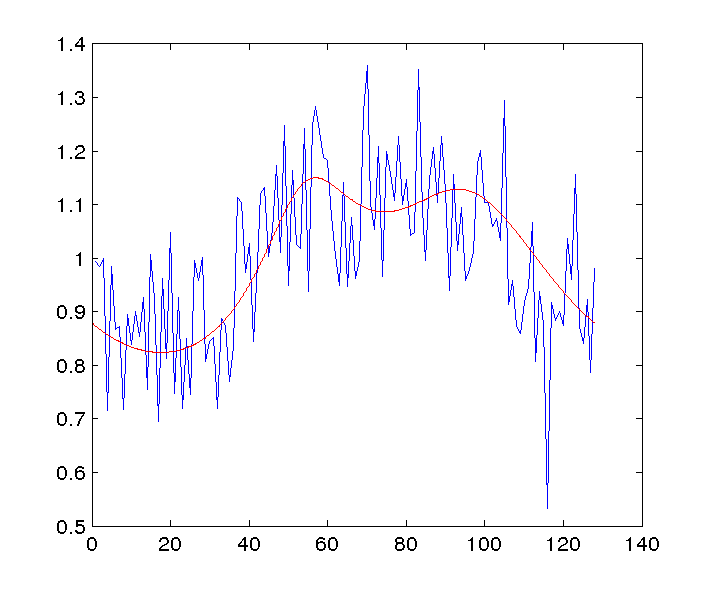}
\end{minipage}
\caption[Density of the measure $\mu^{f_N}_{\Sp^1}$, for $N=5.10^8 + j$, $j=1,2,3$.]{Density of the measure $\mu^{f_N}_{\Sp^1}$, for $N=5.10^8 + j$, $j=1,2,3$. The red curves represent the density of the SRB measure $\mu_0$.}\label{GrafMesInvGro}
\end{figure}

Figures~\ref{GrafMestiti19} and~\ref{GrafMesInvGro}display the density of the measure $\mu^{f_N}_{\Sp^1}$ for various values of $N$. Observe that even when $N$ is quite large (for $N=2^{19}$ on Figure~\ref{GrafMestiti19} and for $N\ge 2^{28}$ on Figure~\ref{GrafMesInvGro}), the measure $\mu^{f_N}_{\Sp^1}$ is quite far away from the SRB measure $\mu_0$. Thus, there is no evidence (at least on this example) that the measures $\mu^{f_N}_{\Sp^1}$ (or the measures $(f_N^*)^k\lambda_N$ for ``large'' times $k$) converge towards the SRB measure when $N$ goes to infinity.

\begin{figure}[ht]
\begin{center}
\includegraphics[width=.5\linewidth,trim = .5cm .3cm .6cm .1cm,clip]{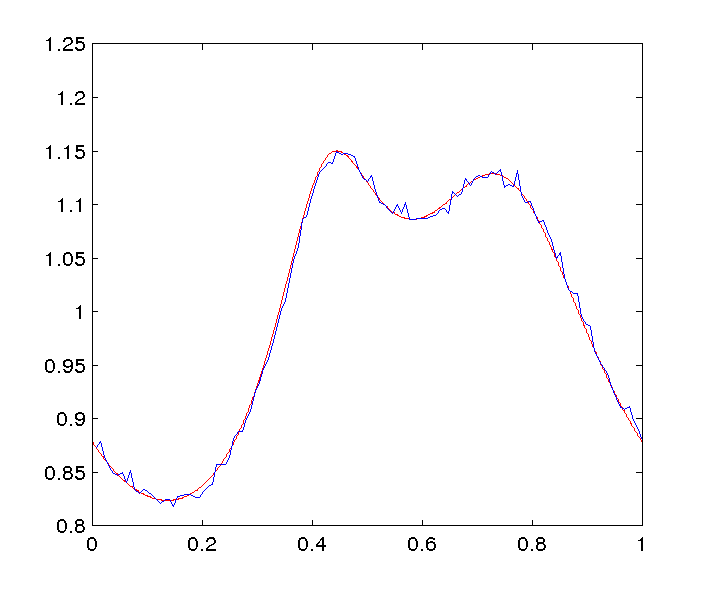}
\caption[Density of $\mu_0$ and mean of the densities of $\mu^{f_N}_{\Sp^1}$]{Density of the SRB measure $\mu_0$ of $f$, and mean of the densities of the measures $\mu^{f_N}_{\Sp^1}$, for $N=5.10^8 + j$, $j\in\llbracket1,1\,000\rrbracket$.}\label{MoyMeasPlInv}
\end{center}
\end{figure}

However, we can see on Figure~\ref{MoyMeasPlInv} that when we take the average of these measures on a large number of different grids, the measure we obtain is very close to the SRB measure. This suggests that on this example, the averages
\[\frac{1}{M}\sum_{N=0}^M \mu^{f_N}_{\Sp^1}\]
might converge towards the measure $\mu_0$.

\clearpage

\part{Discussion of results and research directions}\label{partconcl}

\parttoc

\chapter{Various factors influencing results}\label{Blablablabla}

This thesis is devoted to the study of discretizations in several contexts (conservative or dissipative, continuous or differentiable dynamics, etc.) and according to several viewpoints (combinatorics, topology, ergodic theory; study of a single discretization or all discretizations, etc.). In this chapter, we discuss the influence of these factors on the dynamics of discretizations. This summary seemed essential, as some statements concerning similar issues can differ from almost 150 pages in manuscript, simply because they concern different dynamics (see for example corollaries~\ref{typlax} et~\ref{typlaxC1}).

As in the introduction, we will consider the easy case where the phase space $X$ is the torus $\T^2 = \R^2/\Z^2$, endowed with Lebesgue measure $\Leb$ and uniform discretization grids
\[E_N = \left\{ \left(\frac{i}{N},\frac{j}{N}\right)\in \R^2/\Z^2 \middle\vert\ 1\le i,j \le N\right\}.\]
The results we will present are in fact true in various more general contexts; we refer to the relevant parts of the manuscript for the precise definitions of these contexts.

\section{Specific system versus generic system}

At the beginning of the introduction (page~\pageref{Miaou}), we have seen that the dynamics of the discretizations of a specific system -- for example a linear Anosov torus automorphism -- can be very particular. Specifically, in the case of this example (see Figure~\ref{Miaou} and \cite{Ghys-vari}), on one hand the discretizations are all bijective, and on the other hand there is a very quick recurrence. We interpreted this as a consequence of a resonance between the linear automorphism and the grids\footnote{There is a similar phenomenon for the doubling map $x\mapsto 2x$ on the circle $\R/\Z$. When $\R/\Z$ is equipped with the grids $E_N = \Z /N$, if $N$ is a power of 2, then the discretization $f_N$ has 0 as unique fixed point, which attracts all the points of the grid. If, instead, $N$ is coprime with 2, then $f_N$ is a cyclic permutation of the grid. In none of these cases the dynamics of the discretization does reflect the true dynamics of $x\mapsto 2x$.}.

This contrasts with the results concerning discretizations of generic conservative homeomorphisms obtained in this thesis. For example, let us recall Theorem~\ref{IntroAccuDN} of the introduction (Corollary~\ref{ConjEt} page~\pageref{ConjEt}).

\begin{theo}\label{ConclAccuDN}
For a generic conservative homeomorphism $f\in\Hom(\T^2,\Leb)$, the sequence $(D(f_N))_{N\ge 0}$ accumulates on the whole segment $[0,1]$.
\end{theo}

In particular, there is a subsequence of discretizations having a recurrent set\footnote{Recall that the recurrent set is the union of periodic orbits.} with cardinality arbitrarily small compared with the cardinality of the grid: an infinite number of discretizations are highly non injective. Thus, the generic dynamical behaviour of discretizations is very different from that of some particular system.

Nevertheless, let us stress the fact that considering generic homeomorphisms once grids are fixed induces some resonance (although less marked than in the example above) between the dynamics and the grid. Indeed, our results express that for a generic conservative homeomorphism $f$, each part of the dynamics of $f$ ``resonates'' with an infinite number of grids. For example, this is illustrated by Theorem~\ref{TheoMesIntro} of the introduction (Theorem~\ref{EnsMesInvSimpl} page~\pageref{EnsMesInvSimpl}).

\begin{theo}\label{TheoMesConcl}
For a generic conservative homeomorphism $f\in \Hom(\T^2,\Leb)$, for any $f$-invariant probability measure $\mu$, there exists a subsequence $(f_{N_k})_k$ of discretizations such that $f_{N_k}$ has a unique invariant measure $\mu_k$, which tends to $\mu$. In other words, for any $f$-invariant probability measure $\mu$, there is a sequence $(N_k)_{k\ge 0}$ of integers such that for any $x\in \T^2$,
\[\mu_{x}^{f_{N_k}}\underset{k\to+\infty}{\longrightarrow}  \mu.\]
\end{theo}

The same phenomenon occurs for generic conservative $C^1$-diffeomorphisms (but in a weaker sense): see for example Theorem~\ref{TheoMesIntroC1} of introduction (Theorem~\ref{TheoMesPhysDiff} page~\pageref{TheoMesPhysDiff}).

\begin{theo}\label{TheoMesConclC1}
For a generic conservative $C^1$-diffeomorphism $f\in\Diff^1(\T^2,\Leb)$, for a generic point $x\in \T^2$ (depending on $f$) and any $f$-invariant measure $\mu$, there exists a subsequence $(f_{N_k})_k$ such that
\[\mu_{x}^{f_{N_k}} \underset{k\to+\infty}{\longrightarrow} \mu.\]
\end{theo}

Note also that many of the results about genericity, concerning abstract dynamics, can be observed on concrete examples (see Section~\ref{TheoPrat} for a more detailed discussion).

\section{Conservative versus dissipative}

Of course, it is possible to detect attractive open sets of a dissipative dynamics $f$ on its discretizations. This is a well-known property, also verified by pseudo-orbits.

In the case of a generic dissipative homomorphism, this convergence of the dynamics of discretizations to that of the homeomorphism is generally observed on many properties. For example, a fairly regular behaviour is observed when looking at the combinatorial dynamics of discretizations of a generic dissipative homeomorphism; among others, we have shown the following statement (Corollary~\ref{totsingdiscr} page~\pageref{totsingdiscr}).

\begin{proposition}\label{Maybe}
For a generic dissipative homeomorphism $f\in \Hom(\T^2)$, the rate of injectivity $\card\big(f_N(E_N)\big) / \card (E_N)$ of $f_N$ tends to 0 when $N$ goes to infinity.
\end{proposition}

This trivially implies that thee degree of recurrence $D(f_N)=\card\big(\Omega(f_N)\big)/\card(E_N)$ tends to 0. This result follows directly from a similar property of the original dynamics: a generic dissipative homeomorphism $f$ is \emph{totally singular}, i.e. for every $\varep> 0$, there exists an open set with Lebesgue measure bigger than $1 - \varep$ whose image by $f$ has Lebesgue measure less than $\varep$ (see \cite{MR3027586} and Definition~\ref{strange} page~\pageref{strange} of this manuscript).

From the ergodic point of view, the dynamics of discretizations also converges to that of $f$. More precisely, we have shown the following result (Theorem~\ref{convmesdissip} page~\pageref{convmesdissip}).

\begin{theo}\label{TheoDissipIntro}
For a generic dissipative homeomorphism $f\in \Hom(\T^2)$, the measures $\mu_{\T^2}^{f_N}$ (defined as the Cesàro limits of the pushforwards of uniform measure on $E_N$ by iterates of $f_N$) tend to $\mu^f_{\T^2}$ (which is the Cesàro limit of pushforwards of Lebesgue measure by iterates of $f$).
\end{theo}

Thus, for a generic dissipative homeomorphism, the ``physical'' dynamics of discretizations (i.e. the dynamics of an arbitrarily large proportion of grid points) converges to the ``physical'' dynamics of the homeomorphism (i.e. the dynamics of almost every point with respect to the Lebesgue measure). This was predictable, given that at a certain scale, the dynamics of a generic dissipative homeomorphism is essentially a dynamics of sinks.
\bigskip

In the conservative case, the behaviour of discretizations is much less regular. For example, Theorem~\ref{ConclAccuDN} asserts that unlike the generic dissipative case, the degree of recurrence of discretizations of a generic conservative homeomorphism accumulates on all $ [0,1] $. From the ergodic point of view, Theorem~\ref{TheoMesConcl} says that any invariant measure (and not just the physical ones) of a generic conservative homeomorphism is seen by an infinity of discretizations.
\bigskip

The differences between conservative and dissipative behaviours for generic homeomorphisms can be summarized in the following moral. If we do not assume that the homeomorphism preserves the Lebesgue measure, then the dynamics of discretizations converges to the dynamics of the homeomorphism relative to Lebesgue measure. If, however, we assume that the homeomorphism preserves Lebesgue measure, then the whole dynamics of discretizations reflect all the possible dynamics of the homeomorphism (periodic orbit, compact invariant, invariant measure, rotation vector, etc.) and not just that of almost every point for Lebesgue measure.
\bigskip

This big difference between the behaviour of discretizations in the conservative and the dissipative cases seems rather specific to homeomorphisms: the case of the generic $C^1$-diffeomorphisms seems less contrasted. This is explained by the fact that every chain-recurrence class of a generic dissipative homeomorphism is totally disconnected, while for a generic $C^1$-diffeomorphism they can be stably with nonempty interior. Overall, the results on the dynamics of discretizations that are valid for generic conservative $C^1$-diffeomorphisms are also true for generic dissipative $C^1$-diffeomorphisms, as long as one restricts to a chain-recurrent class. For example, we have shown the following result (Corollary~\ref{typlaxC1} page~\pageref{typlaxC1}).

\begin{proposition}
If $f\in\Diff^1(\T^2,\Leb)$ is a generic conservative diffeomorphism, then for every $\varep>0$ and every $N_0\in\N$, there exists $N\ge N_0$ such that $f_N$ has a $\varep$-dense periodic orbit.

If $f\in\Diff^1(\T^2)$ is a generic dissipative diffeomorphism, then for every chain-recurrence class $K$, every $\varep>0$ and every $N_0\in\N$, there exists $N\ge N_0$ such that $f_N$ has a periodic orbit which is $\varep$-dense in $K$.
\end{proposition}

Thus, the behaviour of discretizations on a chain-recurrence class of a generic dissipative diffeomorphism seems very similar to that observed on discretizations of a generic conservative diffeomorphism. This can be imputed to the fact that real dynamics are themselves relatively close in this case.

\section{Dimension $1$ versus dimension $\ge 2$}

As one can imagine, the dynamics of discretizations of homeomorphisms in dimension 1 on the one hand, and in higher dimension in the other hand, are quite different. This reflects the fact that the dynamics themselves of such systems are very different.
\bigskip

In dimension 1, the length of periodic cycles of discretizations has a uniform behaviour at $+\infty$. For a generic circle homeomorphism  -- thus with rational rotation number $\rho(f) = p/q$ -- the length of these periodic cycles is eventually constant equal to $q$ (Proposition~\ref{mier1} page~\pageref{mier1}, due to T.~Miernowski). However, in dimension bigger than 2, an application of Baire theorem gives the following result.

\begin{proposition}
For a generic dissipative homeomorphism $f\in\Hom(\T^2)$, for any $m\in\N$, there exists a subsequence $(f_{N_k})_{k\ge 0}$ of discretizations such that for any $k$, $f_{N_k}$ has at least $m$ periodic orbits, with pairwise distinct lengths.
\end{proposition}
\bigskip

For $C^r$ circle diffeomorphisms, generic among those with irrational rotation number -- which corresponds more or less to the conservative case\footnote{In the sense that in regularity $C^{1+\alpha}$, Denjoy theorem asserts that the dynamics is transitive.} in higher dimensions -- Theorem~\ref{mier3} (page~\pageref{mier3}) of T.~Miernowski asserts that the common length $q_N(f)$ of periodic orbits of $f_N$ tends to $+\infty$ as $N$ goes to $+\infty$, but arbitrarily slowly. In higher dimension, the following result holds (obtained by combining proofs of Theorems~\ref{corovar2} and \ref{CompactInvSimpl} pages~\pageref{corovar2} and \pageref{CompactInvSimpl}).

\begin{theo}\label{TheoXIntro}
Let $f\in\Hom(\T^2,\Leb)$ be a generic conservative homeomorphism. Then for every $m\in\N$, there exists pairwise distinct integers $p_1,\cdots,p_m$ and a subsequence $(f_{N_k})_{k\ge 0}$ of discretizations such that the set of periods of periodic orbits of $f_{N_k}$ is exactly $\{p_1,\cdots,p_m\}$.
\end{theo}

In particular, for a generic conservative homeomorphism, there is a subsequence of discretizations having a number of periodic orbits tending to $+\infty$, and a subsequence of discretizations whose periodic orbits have uniformly bounded lengths. These phenomena are strongly opposed to what happens on the circle.
\bigskip

The dynamical behaviour of discretizations of a generic \emph{expanding} circle map is fairly close to that of a homeomorphism/diffeomorphism of a higher dimensional manifold. Indeed, as for generic diffeomorphisms (see Theorem~\ref{LocGlobIntro} page~\pageref{LocGlobIntro} of the introduction), we have a result linking local and global behaviours of discretizations.

\begin{theo}
For every $r\ge 1$ and every generic $C^r$ expanding circle map $f$, the rate of injectivity $\tau^k(f)$satisfies
\[\tau^k(f) = \int_{\T^n} \overline D\big( (\det Df_x^{-1})_{\begin{subarray}{l} 1\le m\le k \\ x\in f^{-m}(y)\end{subarray}}\big) \ud \Leb(y)\]
(see Definition~\ref{Noel!} page~\pageref{Noel!} for a definition of $\overline D$).
\end{theo}

Similarly, the ergodic behaviour of generic $C^1$ expanding circle maps strongly resembles that of a generic conservative $C^1$-diffeomorphism of the torus $\T^2$, as can be seen in comparing Theorem~\ref{TheoMesPhysDiff} page~\pageref{TheoMesPhysDiff} with Proposition~\ref{TheoMesPhysDiffExp} page~\pageref{TheoMesPhysDiffExp}.

\section{Regularity $C^0$ versus regularity $C^r$, $r\ge 1$}

The behaviour of generic homeomorphisms is often regarded as not very relevant from a physical point of view: it is locally too complicated to represent a large class of concrete systems. Despite this, the study of such dynamics is not without interest.
\begin{itemize}
\item Firstly, one gets results expressing what can happen when the considered system is not very regular. This is what happens in our case: one sees phenomena highlighted by the theory on quite simple examples of $C^\infty$-diffeomorphisms with large derivatives (see Figure~\ref{FigMesIntro} page~\pageref{FigMesIntro}).
\item Secondly, the generic $C^0$ dynamics are somehow toys models. Their study is a first step in understanding phenomena occurring in higher regularity.
\item Finally, the evolution of the the properties of discretizations of generic systems between regularities $C^0$ and $C^1$ can indicate what should happen in higher regularities.
\end{itemize}

This thesis is not an exception to the rule: the results we obtain in $C^1$ regularity are much weaker than those obtained in $C^0$ regularity. Indeed, to prove genericity results, the difficult part is often to get density statements. The perturbation lemma we use in $C^0$ regularity (Proposition~\ref{extension} page~\pageref{extension}) has a relatively simple proof and allows to perturb independently an arbitrary finite number of points of our space. For now, such a result does not exist in regularity $C^r$ for $r\ge 1$; we only have partial results in the case $r = 1$, called \emph{closing} or \emph{connecting} lemmas. These classical lemmas, whose proofs are much harder than in the $C^0$ case, allow us to transfer directly the results one has about the global behaviour of discretizations of generic homeomorphisms, to results concerning the behaviour of a sub-dynamics of discretizations of a generic $C^1$-diffeomorphism.

For example, one get the following statement (Corollary~\ref{typlax} page~\pageref{typlax}) for generic conservative homeomorphisms, which is a slight improvement of Theorem~\ref{TheoMierno} page~\pageref{TheoMierno} of T.~Miernowski.

\begin{proposition}\label{TypeLaxeIntro}
For a generic conservative torus homeomorphisms, there exists a subsequence of discretizations which are cyclic permutations of the grids.
\end{proposition}

This theorem states that for a generic conservative homeomorphism and for every $\varep> 0$, an infinite number of discretizations is topologically transitive at the scale $\varep$.

For diffeomorphisms, applying the (hard) connecting lemma of C.~Bonatti and S.~Crovisier \cite{MR2090361}, this theorem becomes the following statement (Corollary~\ref{typlaxC1} page~\pageref{typlaxC1}), which only concerns a small part of the point of the grid.

\begin{proposition}\label{TypeLaxeIntroC1}
For a generic conservative diffeomorphism $f\in\Diff^1(\T^2,\Leb)$, for any $\varep>0$, there exists a subsequence of discretizations having at least one $\varep$-dense periodic orbit.
\end{proposition}

Note that the difference between Theorems~\ref{TypeLaxeIntro} and \ref{TypeLaxeIntroC1} is similar to the difference between Theorems~\ref{TheoMesConcl} and \ref{TheoMesConclC1}.
\bigskip

The global behaviour of generic conservative diffeomorphisms is much more difficult to capture. We still get a result (with difficulty!) by considering a ``semi-dynamic'' quantity (since decreasing with time): the degree of recurrence. If it accumulates on all $ [0,1] $ for a generic conservative homeomorphism (Theorem ~\ref{ConclAccuDN}), it tends to 0 for a generic conservative diffeomorphism: let us recall Theorem~\ref{DnZeroIntro} of Introduction (Theorem~\ref{limiteEgalZero} page~\pageref{limiteEgalZero}).

\begin{theo}\label{DnZeroConcl}
For a generic diffeomorphism $f\in\Diff^1(\T^2,\Leb)$, 
\[\lim_{N\to +\infty} D(f_N) = 0.\]
\end{theo}

The proof of this result concerning diffeomorphisms uses crucially the fact that they possess differentials. It is based on an extensive study of the linear case (Part~\ref{PartII}) and on the formula linking the local and global behaviours of discretizations (Theorem~\ref{conv} page~\pageref{conv}, already stated in the introduction as Theorem~\ref{LocGlobIntro}).

\begin{theo}\label{LocGlobConcl}
For every $r\ge 1$ and for a generic conservative $C^r$-diffeomorphism $f\in\Diff^r(\T^2,\Leb)$, for every $t\in{\N^*}$,
\[\lim_{N\to +\infty} \tau_N^t(f) = \int_{\T^2} \tau(D f_{f^{t-1}(x)},\cdots,D f_x) \ud x,\]
where the rate of injectivity of a sequence of matrices is defined similarly to that of a diffeomorphism (see Definition~\ref{DefTaux} page~\pageref{DefTaux}).
\end{theo}
\bigskip

Note that there are open sets of $C^1$-diffeomorphisms made of Anosov diffeomorphisms, which satisfy the shadowing lemma: if $f$ is such a diffeomorphism, then for every $\varep> 0$, there exists $\delta> 0$ such that every $\delta$-pseudo-orbit is $\varep$-shadowed by a true orbit of $f$. In particular, this lemma is verified by the orbits of all discretization on a fine enough grid. One would hope that this very strong property implies that the dynamics of discretizations reflects that of the original diffeomorphism. The reality is somewhat more complex. For example, we saw that the discretizations of Arnold's cat map -- which is a well-known example of Anosov diffeomorphism -- do not reflect at all of the mixing properties of the dynamics. In this case, the dynamics of discretizations only catches the periodic points of small periods (see Figure~\ref{Miaou} page~\pageref{Miaou}). This is due to the fact that a priori, one knows nothing about the orbit that the pseudo-orbit tracks. It can represent rather badly the global dynamics of the application.

Somehow, statements that are shown on the dynamics of discretizations of a generic $C^1$-diffeomorphism $f$ express that every dynamical feature of $f$ is shadowed by similar behaviours of a sub-dynamics of some discretizations $f_{N_k}$ (see for example Proposition~\ref{TypeLaxeIntroC1}). The problem is that these results say nothing about the rest of the dynamics of $f_{N_k}$. The good news is that if $f$ satisfies the shadowing lemma, then it ensures that the rest of the dynamics is not arbitrary: this dynamics is actually tracked by real orbits of $f$.
\bigskip

Theorem~\ref{DnZeroConcl} also becomes interesting when compared to the corresponding case in $C^0$ regularity. It indicates that the bad behaviours observed for generic conservative homeomorphisms should disappear in higher regularity: the behaviour of the rate of injectivity is less irregular for generic conservative $C^1$-diffeomorphisms than for generic conservative homeomorphisms (compare with Theorem~\ref{ConclAccuDN}). This big difference suggests that the wild behaviours of the global dynamics of discretizations of generic conservative homeomorphisms observed in Chapter~\ref{ChapCons} may not appear in higher regularities. For example, contrary to what happens for generic conservative homeomorphisms, one can hope recovering the physical measures of a generic $C^1$ conservative diffeomorphism by looking at the measures $\mu_{\T^2}^{f_N}$, although the ``local measures'' $\mu_{x}^{f_N}$ accumulate on the whole set of $f$-invariant measures (Theorem~\ref{TheoMesConclC1}). This behaviour is also suggested by the numerical experiments we have conducted (see Subsection~\ref{RedW}).
\bigskip

Finally, we have shown that in short time, the pushforwards of the uniform measure on the grid by the discretizations of a smooth enough circle expanding map converge towards a single measure, which is the SRB measure $\mu_0$ of the expanding map ($\mu_0$ is also the unique absolutely continuous invariant measure). Specifically, the following theorem holds (Theorem~\ref{MainMoche}).

\begin{theo}\label{RepLanfordConlc}
For every $\alpha>0$ and every circle $C^{1+\alpha}$ expanding map $f$, there exists a constant $c_0 = c_0(f)>0$ such that if $(N_m)_m$ is a sequence of integers tending to infinity but such that $\ln N_m > c_0 m$, then the following convergence holds
\[(f_{N_m}^*)^m(\lambda_{N_m}) \to \mu_0.\]
\end{theo}

Observe that this result does not uses genericity. Note also that this theorem says nothing about what happens in the case of generic $C^1$ expanding maps: in this case, the techniques used in the proof (i.e. the Ruelle-Perron-Frobenius transfer operator) no longer work at all (and the behaviour of physical measures is very different, see for example \cite{MR1845327}).

\section{One discretization versus most of the discretizations versus all the discretizations}

For generic conservative homeomorphisms, the philosophy of the results is not really the same whether one considers a single discretization or all discretizations\footnote{In this paragraph, we will not treat the dissipative case, since the dynamics of discretizations converges uniformly to the dynamics of the starting homeomorphism.}. If one looks at the dynamics of a single discretization, its behaviour wildly depends on the order of the discretization, as expressed by the results of Section~\ref{partie 1.3} of Chapter~\ref{ChapCons}. For example, for a generic conservative homeomorphism, there exists a constant $P>$ 0 such that the following properties are all satisfied by an infinite number of discretization orders $ N $:

\begin{itemize}
\item $f_N$ is a cyclic permutation (Corollary \ref{typlax} page~\pageref{typlax});
\item $f_N$ has a single cycle, which has period smaller than $P$; this cycle attracts all the points of the grid $E_N$ (Corollary \ref{corovar2} page~\pageref{corovar2});
\item $f_N$ has at least $\sqrt{\card(E_N)}$ different cycles (Corollary \ref{corovar3} page~\pageref{corovar3}).
\end{itemize}
\bigskip

Worse, this variability of behaviours of discretizations persists if one tries to make statistics on the frequency of appearance of the properties among the discretizations; (see Section~\ref{bofbof} of Chapter~\ref{ChapCons}). For example, one shows the following result (Theorem \ref{propdemin} page~\pageref{propdemin} and Corollary~\ref{CoroAver1} page~\pageref{CoroAver1}).

\begin{theo}\label{MoyCesaroIntro}
Let $C$ be the set of $N\in\N$ such that $f_N$ is a cyclic permutation. Then, for a generic conservative homeomorphism $f\in\Hom(\T^2,\Leb)$, the sequence $\card(C\cap \llbracket 1,M\rrbracket)/M$ accumulates on both 0 and 1.
\end{theo}
\bigskip

Nevertheless, things are much better when considering the dynamics of \emph{all} the discretizations of the homeomorphism: in this case, what emerges from the results proved in Section~\ref{Sec8} is that it is possible to detect many dynamical properties of the homeomorphism on discretizations. For example, one can detect the periodic orbits of $f$ and their periods (Theorem~\ref{corovar2} page~\pageref{corovar2}).

\begin{theo}\label{TheoOrbPerC0Intro}
Let $f\in\Hom(\T^2,\Leb)$ be a generic conservative homeomorphism. Then, for every $f$-periodic orbit $\omega$ of period $p$, and every $\delta>0$, there exists a subsequence of discretization $f_{N_k}$ such that for every $k$, $f_{N_k}$ has a unique periodic orbit\footnote{Thus this periodic orbit attracts all the points of the grid.}, of length $p$, and which $\delta$-shadows the orbit $\omega$.
\end{theo}

Similarly, Theorem~\ref{TheoMesConcl} on invariant measures says that it is possible to retrieve all the invariant measures of the original homeomorphism by looking at the invariant measures of all discretizations of this homeomorphism. In a way, these statements express that the dynamics of the homeomorphism is tracked by that of its discretizations.
\bigskip

Concerning generic conservative $C^1$-diffeomorphisms, the results are yet too partial to decide whether there is a convergence of the dynamics of discretizations or not.

On the one hand, one knows that some dynamical invariants of the diffeomorpihsm can be detected on a sub-dynamics of an infinite number of discretizations. For example, periodic orbits can be recovered from the discretizations (see Lemma~\ref{PropShadowC1} page~\pageref{PropShadowC1}).

\begin{proposition}\label{TheoOrbPerC1Intro}
Let $f\in\Diff^1(\T^2,\Leb)$ be a generic conservative $C^1$-diffeomorphism. Then, for every $f$-periodic orbit $\omega$ of period $p$, and every $\delta>0$, there exists a subsequence of discretizations $f_{N_k}$ such that for every $k$, $f_{N_k}$ has at least one periodic orbit\footnote{Remark that contrary to the $C^0$ case, this orbit is in general not unique.} of length $p$, which $\delta$-shadows the orbit $\omega$.
\end{proposition}

Similarly, Theorem~\ref{TheoMesConclC1} indicates that it is possible to retrieve invariant measures of a generic conservative diffeomorphism (see also Corollary~\ref{CoroMane} for a proof of a simpler result).

On the other hand the global behaviour of the discretizations is more regular than in the $C^0$ case, as emphasized by the fact that the degree of recurrence tends to 0 for a generic conservative diffeomorphism (Theorem \ref{DnZeroConcl}). Unfortunately, for now, no one knows how to describe better the nature of the global dynamics of discretizations. For example, we have no idea of how is the asymptotic behaviour of the sequence of measures\footnote{Recall that given a discretization $f_N$, $\mu_{\T^2}^{f_N}$  is defined as the Cesàro limit of pushforwards by iterates of $f_N$ of the uniform measure on the grid $E_N$.} $\big(\mu_{\T^2}^{f_N}\big)_{N\in\N}$.

\section{Combinatorial properties versus topological and ergodic properties}

As already seen, some combinatorial quantities associated to discretizations of a generic conservative homeomorphism evolve very erratically depending on the order of $N$ discretization. For example, the degree of recurrence accumulates on all $[0,1]$ (Theorem~\ref{ConclAccuDN}), the lower limit of the number of periodic orbits is $1$ (Theorem~\ref{TypeLaxeIntro}) and its upper limit is $+\infty$ (Theorem~\ref{TheoXIntro}, see also Corollary~\ref{corovar3} page~\pageref{corovar3}), etc. We have even shown a general theorem in this sense: a property about finite maps of the grids $E_N$ is \emph{dense} if any conservative homeomorphism is arbitrarily well approximated by finite maps with this property. Theorem~\ref{génécycl} page~\pageref{génécycl} says that if a property is dense, then it appears on an infinite number of discretizations of a generic homeomorphism.

So we can say that these combinatorial quantities are not very useful to detect the dynamics of the initial application. Moreover, these results imply that the discretizations do not behave at all like a typical random map of a finite set: for example, such a map $\sigma : E \to E$, where $E$ is a finite set with $q$ elements, is such that the cardinality of its recurrent set is of the order of $\sqrt N$ (see \cite[XIV.5]{Boll-rand} or Theorem~2.3.1 of \cite{Mier-dyna}). Worse, it is still not true when one considers for example the average of degree of recurrence on long segments of discretization orders (see Theorem~\ref{MoyCesaroIntro}).

Regarding diffeomorphisms, we know very little about the global combinatorial behaviour of discretizations. For this problem, the main theorem we get is Theorem~\ref{DnZeroConcl}, which illustrates the fact that the discretizations of a generic conservative $C^1$-diffeomorphism behave differently from those of a generic conservative homeomorphism.
\bigskip

To recover the original dynamics, we have to use the geometry of the grids $E_N$, and focus on topological or ergodic properties of all discretizations. As already stated, both in $C^0$ and $C^1$ topologies, this allows to recover for example the set of invariant measures of the initial dynamics (Theorems\ref{TheoMesConcl} and \ref{TheoMesConclC1}) or all periodic points (Theorem~\ref{TheoOrbPerC0Intro} and Proposition\ref{TheoOrbPerC1Intro}), or the rotation set, as shown by Theorem~\ref{EnsRotIntro} of Introduction (Theorem~\ref{CoroRotDiscrCons} page~\pageref{CoroRotDiscrCons}).

\begin{theo}\label{EnsRotConcl}
For a generic conservative homeomorphism $f\in\Hom(\T^2,\Leb)$,
\begin{itemize}
\item the observable rotation set is reduced to a single point;
\item the supremum limit of the rotation sets of discretizations coincides with the rotation set of $f$.
\end{itemize}
\end{theo}

\section{Microscopic versus mesoscopic versus macroscopic}

A generic homomorphism (conservative or dissipative) has a very chaotic local behaviour. For instance, when it has a periodic point of period $p$, then its set of periodic points of period $p$ is a Cantor set (so it is uncountable). Moreover, for any multiple of $q$ of $p$, the set of periodic points of period $q$ is nonempty; and any periodic point of period $p$ is the limit of a sequence of periodic points of period $q$. This lack of local regularity of the map makes the local behaviour of its discretizations also very chaotic. Moreover, this local behaviour of discretizations strongly depends on the order of the discretization. For example, a discretization is sometimes locally bijective (Proposition~\ref{TypeLaxeIntro}), and sometimes locally  highly non injective, as shown by the following result (Corollary~\ref{crush} page~\pageref{crush}).

\begin{theo}\label{crushIntro}
Let $\vartheta : \N\to\R_+^*$ be a map tending to $+\infty$ in $+\infty$. Then, for a generic conservative homeomorphsim $f\in\Hom(\T^2,\Leb)$, we have\footnote{Remark that this result is an improvement of the fact that the degree of recurrence of a generic conservative homeomorphism accumulates on 0.}
\[\underset{N\to+\infty}{\underline\lim}\ \frac{\card\big(f_{N}(E_{N})\big)}{\vartheta({N})}=0.\]
\end{theo}

Therefore, the local behaviour of a discretization of a generic conservative homomorphism is rather wild, and therefore somehow irrelevant. This is reflected at a temporal level: for example, comparing Proposition~\ref{TypeLaxeIntro} and Theorem~\ref{crushIntro}, we see that discretizations of two (arbitrarily large) different orders may have opposite behaviours from the first iteration .
\bigskip

Things are quite different for diffeomorphisms. The existence of (continuous) differentials to a diffeomorphism $f$ dictates the local behaviour of its discretizations. This is illustrated by Theorem~\ref{LocGlobConcl}, which expresses that generically, a time $t$ being fixed, for large enough discretization orders, the cardinalities of the images of the grid by the discretizations in time smaller than $t$  is determined by the differentials of $f$. There are therefore three different space scales:
\begin{itemize}
\item Firstly, at the scale of the torus -- the macroscopic scale -- in finite time, and for large enough orders of discretization, the  dynamics of the discretizations resembles that of the original diffeomorphism.
\item Secondly, at the scale of the differentials of $f$ -- the mesoscopic scale -- one perceive the action of derivatives of $f$. In finite time, and for large enough orders of discretization, the dynamics of the discretizations resembles that of a sequence of matrices.
\item Finally, at the scale of the grid -- the microscopic scale -- we are able to see the points of the grids, the ``atoms''.
\end{itemize}

As for homeomorphisms, these spatial scales are transposed to temporal level. Here, having an additional mesoscopic scale induces a transitional regime for the local behaviour of discretizations, as expressed by Theorem~\ref{LocGlobConcl} (which is not really a purely  dynamical statement because it concerns only a finite number of iterations). To derive an asymptotic result we then use the decreasing of the rate of injectivity over time.
\bigskip

\begin{figure}[t]
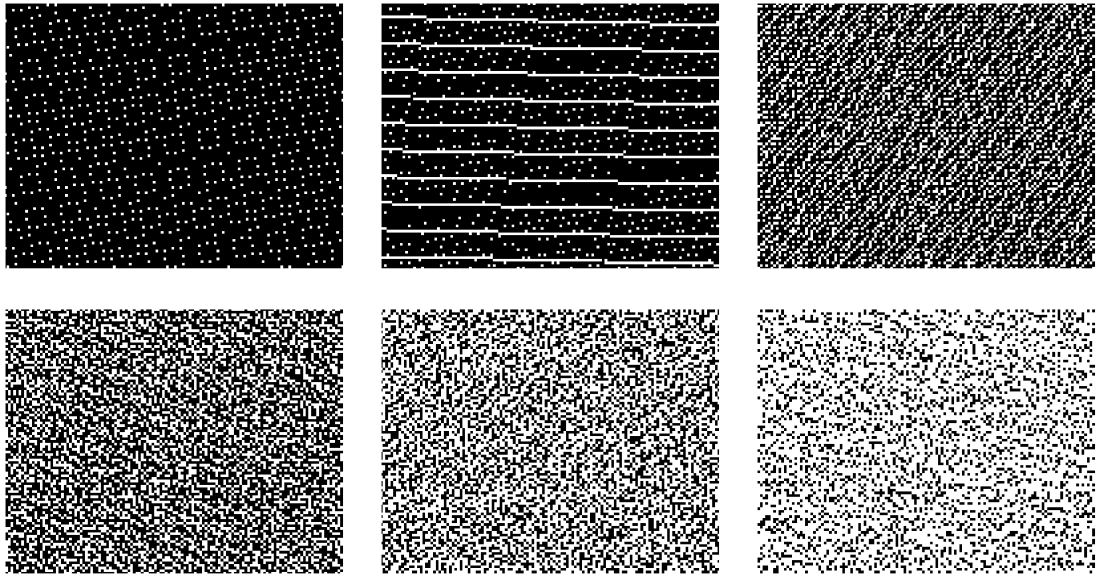

\begin{minipage}[c]{.33\linewidth}
	\includegraphics[width=\linewidth, trim = 1.5cm .5cm 1.5cm .5cm,clip]{Fichiers/ImagesEnsQuasiPer/ImageSlTemps1.png}
\end{minipage}\hfill
\begin{minipage}[c]{.33\linewidth}
	\includegraphics[width=\linewidth, trim = 1.5cm .5cm 1.5cm .5cm,clip]{Fichiers/ImagesEnsQuasiPer/ImageSlTemps2.png}
\end{minipage}\hfill
\begin{minipage}[c]{.33\linewidth}
	\includegraphics[width=\linewidth, trim = 1.5cm .5cm 1.5cm .5cm,clip]{Fichiers/ImagesEnsQuasiPer/ImageSlTemps3.png}
\end{minipage}

\begin{minipage}[c]{.33\linewidth}
	\includegraphics[width=\linewidth, trim = 1.5cm .5cm 1.5cm .5cm,clip]{Fichiers/ImagesEnsQuasiPer/ImageSlTemps5.png}
\end{minipage}\hfill
\begin{minipage}[c]{.33\linewidth}
	\includegraphics[width=\linewidth, trim = 1.5cm .5cm 1.5cm .5cm,clip]{Fichiers/ImagesEnsQuasiPer/ImageSlTemps10.png}
\end{minipage}\hfill
\begin{minipage}[c]{.33\linewidth}
	\includegraphics[width=\linewidth, trim = 1.5cm .5cm 1.5cm .5cm,clip]{Fichiers/ImagesEnsQuasiPer/ImageSlTemps20.png}
\end{minipage}
\caption[Successive images of $\Z^2$ by discretizations of matrices of $SL_2(\R)$]{Successive images of $\Z^2$ by discretizations of $k$ matrices of $SL_2(\R)$, chosen randomly (and independently) in a compact of $SL_2(\R)$ (see Figure~\ref{ImagesSuitesMat} for a precise explanation). A point of $\Z^2$ is coloured in black if it belongs to the image set. From left to right and top to bottom, $k=1,\, 2,\, 3,\, 5,\, 10,\, 20$.}\label{ImagesSuitesMatIntro}
\end{figure}

Thus, to study the dynamics of diffeomorphisms, we are led to consider the dynamics of discretizations of sequences of linear maps. Specifically, the discretization of a matrix $A\in GL_2(\R)$ is the map $\widehat A = P \circ A_{|\Z^2} $, where $P$ is a projection of $\R^2$ on the closest point of $\Z^2$. We want to study the sets $(\widehat{A_K} \circ \cdots \circ \widehat{A_1}) (\Z^n) $ for a generic \footnote{Generic in the sense of the topology spanned by the supremum norm on sequences, relative to a standard norm on $M_2(\R)$ fixed once for all.} sequence of matrices $(A_k)_{k \ge 1}$ (possibly with determinant 1). One can see in Figure~\ref{ImagesSuitesMatIntro} typical images of the image sets obtained. Among others, we can observe a phenomenon of almost-periodicity of the firsts sets. Indeed, we prove that this is the case in any time \footnote{Although it is far from being obvious on the latest images of the figure.} (Theorem~\ref{imgquasi}). Roughly speaking, a set $\Gamma$ is \emph{almost periodic} if for $R$ large enough, the set $B(0,R) \cap \Gamma$ determines the set $\Gamma$ up to a set density of less than $\varep$ (see Definition~\ref{DefAlmPer} page~\pageref{DefAlmPer}).

\begin{theo}
For any sequence $(A_k)_{k\ge 1}$ of invertible matrices, the sets $(\widehat{A_k}\circ\cdots\circ \widehat{A_1})(\Z^n)$ are almost periodic.
\end{theo}

Observing Figure~\ref{ImagesSuitesMatIntro}, we also remark that the density of the image sets seems to decrease with time. This behaviour is explained by the main theorem of Part~\ref{PartII} (Theorem~\ref{ConjPrincip} page~\pageref{ConjPrincip}).

\begin{theo}\label{CallMe}
For a generic sequence $(A_k)_{k\ge 1}$ of matrices with determinant 1, the density of the sets $(\widehat{A_k}\circ\cdots\circ \widehat{A_1})(\Z^n)$ tends to $0$ as $k$ tends to $+\infty$.
\end{theo}

The proof of this theorem, which is quite long and difficult, is mainly based on arguments of equidistribution, which allow to reduce the computation of the density of image sets (which is not very practical to handle, in particular it is nonlinear) to that of the area of an intersection of cubes in high dimensions.

\section{In theory versus in practice}\label{TheoPrat}

For each context studied in this manuscript (homeomorphisms and $C^1$-diffeomorphisms, both dissipative as conservative), we have compared the statements obtained with the results of numerical simulations carried out on examples supposed to represent the generic case.
\bigskip

\begin{figure}[t]
\begin{center}
\includegraphics[width=.35\linewidth,trim = .5cm .3cm .6cm .1cm,clip]{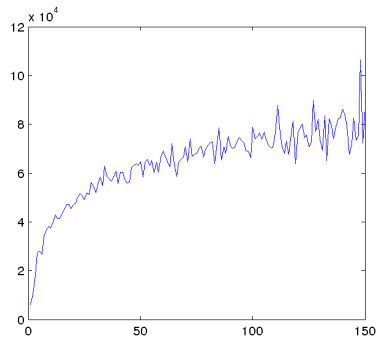}
\caption[Degree of recurrence of a generic conservative homeomorphism]{Cardinality of the recurrent set $\Omega(f_N)$ depending on $N$, for $N = 128k$, $k=1,\cdots,150$, and $f$ a conservative homeomorphism close to identity (see page~\pageref{PageDefSimulCons} for the definition of this homeomorphism).}\label{GrafIntro}
\end{center}
\end{figure}

Let us start with the case of conservative homeomorphisms. At first glance, simulations of quantities such as the degree of recurrence can be seen as disappointing: in practice, the degree of recurrence does not accumulate on all $[0,1]$, but tends to 0 relatively quickly (see Figure~\ref{GrafIntro}, see also Section~\ref{partietroisb}). This contradicts the conclusions of Theorem~\ref{ConclAccuDN}: in practice, the behaviour described by the statements regarding the generic case cannot be observed in practice, even when the definition of the map is made to mimic the generic case.

\begin{figure}[t]
\begin{center}
\begin{minipage}[t]{.48\linewidth}
\begin{center}
\includegraphics[width=.8\linewidth,trim = 1.5cm .95cm 1cm .5cm,clip]{Fichiers/Serie128/IdC0Dissip2/Mesure32768.png}
\caption[Measure $\mu_{\T^2}^{f_N}$ for a dissipative homeomorphism]{Density of the measure $\mu_{\T^2}^{f_N}$, for $f$ an example of dissipative homeomorphism close to identity with small attractive sets, and $N = 32\,768$ (see page~\pageref{defbis} for the definition of this homeomorphism).}\label{MesDissipIntro1}
\end{center}
\end{minipage}\hfill
\begin{minipage}[t]{.48\linewidth}
\begin{center}
	\includegraphics[width=.8\linewidth,trim = 1.5cm .95cm 1cm .5cm,clip]{Fichiers/Serie128/IdC0Dissip/Mesure32768.png}
		\caption[Measure $\mu_{\T^2}^{f_N}$ for a dissipative homeomorphism]{Density of the measure $\mu_{\T^2}^{f_N}$, for $f$ an example of dissipative homeomorphism close to identity with large attractive sets, and $N = 32\,768$ (see page~\pageref{defbis} for the definition of this homeomorphism).}\label{MesDissipIntro2}
\end{center}
\end{minipage}
\end{center}
\end{figure}

To observe in practice the phenomena predicted by theoretical results, one has to stop looking at the combinatorial dynamics of discretizations. For example, we have simulated the measures $\mu_{\T^2}^{f_N}$ (recall that they are the Cesàro limits of pushing forwards of the uniform measure on $E_N$ by the iterates of $f_N$). In practice, it can not really be expected that the measures $\mu_{\T^2}^{f_N}$ converge towards all invariant measures of $f$, because this set is an infinite dimensional convex set. Nevertheless, one can determine whether these measures seem to converge or not. In Figure~\ref{FigMesIntro} page~\pageref{FigMesIntro}, we see that the measures $\mu_{\T^2}^{f_N}$ do not converge at all towards Lebesgue measure, but in addition they have nothing to do with each other, even when the discretization orders are very close. As predicted by Theorem~\ref{TheoMesConcl}, we do not recover the physical measure of the initial homeomorphism on discretizations.
\bigskip

For dissipative homomorphisms, we also have simulated the measures $\mu_{\T^2}^{f_N}$ on two different examples. The first one is a dissipative homeomorphism close to identity, but with small sinks (see Figure~\ref{MesDissipIntro1}, see also Section~\ref{partietrois}). It turns out that the behaviour of these measures is identical to that observed in the conservative case: it does not detect the dissipative nature of the dynamics at all. This had already been pointed out by J.-M. Gambaudo and C. Tresser in \cite{MR700317}: if the sinks are too small -- this can happen quite easily, even if the definition of the application is rather simple -- then they are not detected by simulations.

To avoid this phenomenon, we have also tested an example of dissipative homeomorphism close to identity, but having much larger sinks (see Figure~\ref{MesDissipIntro2}, see also Section~\ref{partietrois}). In this case, as predicted by Theorem~\ref{TheoDissipIntro}, the measures $\mu_{\T^2}^{f_N}$ seem to converge rather quickly to a measure supported  by the attractive sets of the homeomorphism.
\bigskip

\begin{figure}[h]
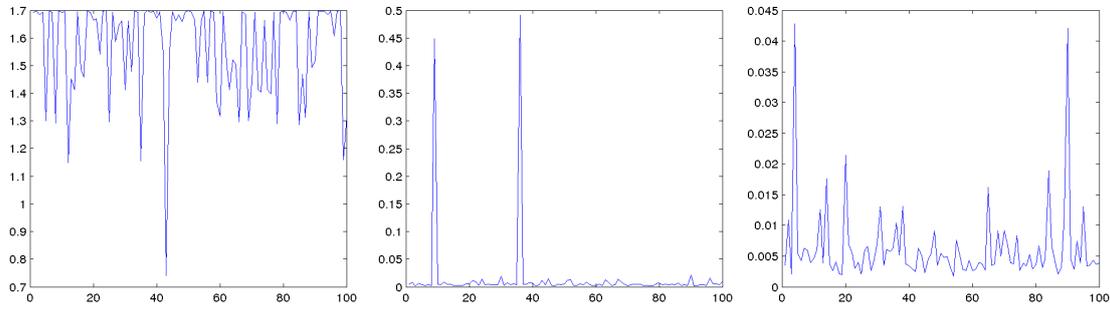

\begin{minipage}[c]{.33\linewidth}
	\includegraphics[width=\linewidth,trim = .5cm .3cm .6cm .1cm,clip]{Fichiers/MesPhys/IdC1Bis/DistLeb100.png}
\end{minipage}\hfill
\begin{minipage}[c]{.33\linewidth}
	\includegraphics[width=\linewidth,trim = .5cm .3cm .6cm .1cm,clip]{Fichiers/MesPhys/RotC1Bis/DistLeb100.png}
\end{minipage}\hfill
\begin{minipage}[c]{.33\linewidth}
	\includegraphics[width=\linewidth,trim = .5cm .3cm .6cm .1cm,clip]{Fichiers/MesPhys/AnoC1Bis/DistLeb100.png}
\end{minipage}
\caption[Distance between $\Leb$ and $\mu^{(f_i)_N}_x$ for 3 different examples of conservative diffeomorphisms]{Distance between Lebesgue measure ans measures $\mu^{(f_i)_N}_x$ depending on $N$, for $f_1$ a small $C^1$-perturbation of $\Id$ (left), $f_2$ a small $C^1$-perturbation of a translation of order 30 (middle) and $f_3$ a small $C^1$-perturbation of a linear Anosov automorphism (right), on grids $E_N$, with $N=2^{20}+k$, $k=1,\cdots,100$ (see page~\pageref{DefDiffeoPhys} for a precise definition of these diffeomorphisms).}\label{GrafDistLebPhysConcl}
\end{figure}

Note also that the frequency of occurrence of the phenomena described by the theorems can change dramatically depending on the considered type of dynamics. For example, Figure~\ref{GrafDistLebPhysConcl} illustrates Theorem~\ref{TheoMesConclC1}, which expresses that in some cases (in fact, most of the cases in the topological sense), if $x\in\T^2$ and $f\in\Diff^1(\T^2,\Leb)$ are fixed, then the measures\footnote{Recall that $\mu_x^{f_N}$ is the uniform measure on the periodic orbit where falls the positive orbit of $x_N$ under $f_N$.} $\mu_x^{f_N}$ accumulate on all the invariant measures of $f$. Figure~\ref{GrafDistLebPhysConcl} represents the distance $\dist(\mu_x^{f_N},\Leb)$ depending on $N$ for three different conservative diffeomorphisms (see Section~\ref{NumSimPhys} for precise definitions). For the first diffeomorphism, which is close to an ``elliptical'' dynamics -- namely identity --, this distance is still quite high: the obtained measure is always far away from Lebesgue measure (see Figure~\ref{FigMesIntro2} page~\pageref{FigMesIntro2}). The second tested diffeomorphism is a small $C^1$-perturbation of a translation of order 30 of the torus. In this case, the measure is almost always very close to Lebesgue measure, except for some few orders $N$, where it is quite far away. The third diffeomorphism is a small $C^1$-perturbation of a linear Anosov automorphism. Again, the calculated measure is almost always very close to Lebesgue measure, except for some few orders $N$ (although this is less marked than for the second diffeomorphism). As could be expected, it seems that in practice, the more the initial dynamics is chaotic, the less there are orders of discretization for which the conclusions of the obtained statements are true.

\chapter{Research directions}\label{ouaf}

\section{Combinatorial dynamics of discretizations of generic diffeomorphisms}

So far, we know very little about the combinatorial dynamics of discretizations of generic conservative diffeomorphisms. As already said, we know that the degree of recurrence tends to 0 (Theorem~\ref{DnZeroConcl}), but it is actually the only result we have about the global behaviour of discretizations. We still do not know answer simple questions such as the following.

\begin{ques}
Does the cardinality of the recurrent set of discretizations $f_N$ of a generic conservative diffeomorphism $f$ goes to infinity when $N$ goes to infinity? If so, do the length of the longest periodic orbit of $f_N$, or the number of periodic orbits of $f_N$, tend to infinity?
\end{ques}

The next step would be to determine whether discretizations of generic diffeomorphisms behave like random maps (see Section~\ref{AppAlea}). This kind of questions seems to me out of reach for the moment.

\section{Convergence of canonical invariant measures}\label{SecQuesPhys}

We still do not know the behaviour of the ``canonical'' invariant measures $\mu_{\T^2}^{f_N}$ of $f_N$  when $f$ is a generic conservative $C^1$-diffeomorphism. Recall that each of these measures $\mu_{\T^2}^{f_N}$ is supported by the recurrent set of $f_N$; the total measure of each periodic orbit being proportional to the size of its basin of attraction. In the case of generic conservative homeomorphisms, these measures accumulate on all $f$-invariant measures (Theorem~\ref{TheoMesConcl}). On the contrary, in the case of generic dissipative homeomorphisms, they converge towards a single measure (Theorem~\ref{TheoDissipIntro}). In the case of generic conservative diffeomorphisms, their behaviour is so far not known.

\begin{ques}
For a generic conservative diffeomorphism $f\in \Diff^1(\T^2,\Leb)$, do the measures $\mu_{\T^2}^{f_N}$ tend to Lebesgue measure? Do they accumulate on the whole set of $f$-invariant measures?
\end{ques}

The numerical simulations we present do not allow us to guess an answer to this question (see Figure~\ref{MesC1IdCons2p}). If such measures are much closer to Lebesgue measure that in the $C^0$ case, one still observe small regions of the torus with large measure.

Still in the case of generic conservative diffeomorphisms, one can also wonder what happens to the measures  $\mu_{x}^{f_N}$. Recall that the measure  $\mu_{x}^{f_N}$ is the uniform measure on the periodic orbit in which the orbit of $x$ under $f_N$ falls. Theorem~\ref{TheoMesConclC1} says that for a generic set points $x$, these measures accumulate on all $f$-invariant measures, but the following questions remain open.

\begin{ques}
For a generic conservative diffeomorphism $f\in \Diff^1(\T^2,\Leb)$, what is the behaviour of the measures $\mu_{x}^{f_N}$ for almost every point $x$ (for Lebesgue measure) ? From the point of view of the grids, what is the behaviour of the measures $\mu_{x}^{f_N}$ for ``most of '' the points $x_N\in E_N$ ?
\end{ques}

\section{Detection of some dynamical properties}

In this manuscript, we do not address the detection of a number of dynamic invariants. For example, our techniques provide no way to determine whether a dynamics is transitive or not. In fact, the systems studied are very specific: a generic dissipative homeomorphism is never transitive, while a generic conservative homeomorphism -- or $C^1$-diffeomorphism -- is always. We could also try to detect the conservativeness of a dynamics. Unfortunately, we have no algorithm to answer this type of question. We do not know either how to answer questions like the following.

\begin{ques}
If for every $\varep>0$, an infinite number of discretizations has a $\varep$-dense periodic orbit, does it imply that the initial dynamics is transitive?
\end{ques}

A classical problem in dynamics is the effective calculation of entropy or Lyapunov exponents of a system. The results obtained in this thesis do not provide answers to such questions, which are only relevant in the case of generic diffeomorphisms (for a generic conservative homeomorphism, the metric entropy is zero and topological entropy is infinite). We just can have a lower bound (possibly very bad) of the topological entropy of a conservative diffeomorphism of the torus with the calculation of its rotation set (with the bounds obtained by J.~Kwapisz in \cite{MR1213082}).

\begin{ques}
Only knowing the dynamics of discretizations, is it possible to calculate the entropy of the initial system?
\end{ques}

To handle this problem, the best method is probably not looking at the dynamics of discretizations, that is to say, looking at the behaviour of discretizations ``in infinite time''. It may be more efficient to stop computations at a well chosen time, presumably logarithmic in $N$ (see next section).

\section{Explicit estimations about convergence of dynamics}

One of the shortcomings of the results presented in this manuscript is their lack of effectiveness (using Baire theorem, we use countable axiom of choice). For example, when obtaining results of simulations of rotation sets using coarse discretizations (e.g. Figure~\ref{FigRotIntro} page~\pageref{FigRotIntro}), there is no theoretical result estimating whether the obtained sets are close to true rotation set or not. It would probably be quite easy to have estimates such as ``for $\delta> 0$, there exists an explicit $N_0>0$ such that if the grid order is bigger than $N_0$, then the rotation set of the discretization is contained in the $\delta$-neighbourhood of the real rotation set''. By cons, there is currently no way to know when the rotation set of the discretization ``fills'' the rotation set of the homeomorphism. One can only get an indication of the convergence when the obtained sets seem to stabilize; however this does not ensure that the simulated rotation set is actually close to the rotation set of the homeomorphism. It would be very interesting to be able to provide answers to the following question.

\begin{ques}
Is it possible to find a reasonably fast algorithm of approximation of the rotation set, such that there are rigorous estimates of the distance between the calculated and the actual rotation set, depending on the parameters of the problem?
\end{ques}

\section{Different notions notions of genericity}

In this thesis, we discuss results of genericity in only two specific cases: the $C^0$ and $C^1$ topologies. One could imagine study what happens in different contexts. One could of course be interested in what happens in the intermediate case of Hölder generic applications. It would also be natural to look at what happens in $C^r$ topology, with $r>1$, $r = \infty $ or $r = \omega$ (real analytic applications). However, the results we obtained are rather partial in $C^1$ topology, it may be very difficult to have satisfactory statements in these more rigid contexts.

In this manuscript, we distinguish the behaviours of discretizations of conservative and dissipative generic systems. One could imagine studying what happens to other types of dynamics like generic chain-transitive homeomorphisms\footnote{For a study of the dynamics of generic chain-transitive continuous maps (necessarily homeomorphisms), see the book of E.~Akin  \cite{MR1627928}. It is possible that the techniques presented there can be used for the study of discretizations of such applications.}.

There are other concepts of genericity than that of Baire. For example, \emph{prevalence}, which is supposed to be closer to the concept of full measure set. In this direction, part of the article \cite{MR2279269} of T.~ Miernowski about discretizations of circle diffeomorphisms focuses on prevalent diffeomorphisms.

The most most natural idea for the study of discretizations is perhaps that of Kolmogorov genericity: for a function $f$ given by an explicit formula, the Kolmogorov complexity describes the number of elementary operations needed to calculate the images $f(x)$. Indeed, to make numerical simulations, one has necessarily an explicit expression of the evolution law; so we want a notion of genericity for maps definable by a formula. One can then ask the following question.

\begin{ques}
What are the properties of discretizations of most of the functions with Kolmogorov complexity less than a given number?
\end{ques}

This problem seems far too complex to get satisfactory answers in the short term.

Note that some properties we got are true on open and dense sets of applications. For example, for every $\tau_0> $ 0, the property ``$\lim_{t\to +\infty}\tau^t(f) \le \tau_0$'' is satisfied on an open and dense subset of conservative $C^1$-diffeomorphisms. In this case, the concepts of genericity and prevalence more or less coincide.

\section{Generic behaviours among conjugates to a given system}

Finally, one may wonder what happens when one fixes a dynamics $f$, and considers its discretizations over a generic sequence of grids. A partial answer to this question is given in Section~\ref{AddendParti1} of Chapter~\ref{ChapCons}: if we consider a conservative homeomorphism $f$, and if ``generic sequence of grids'' means ``image of a fixed good sequence of grids by a generic conservative homeomorphism'', then it happens exactly the same as for generic conservative homeomorphisms with respect to a fixed sequence of grids, provided that the fixed points set of $f$ has empty interior. It would be interesting to know what happens in other cases, in particular to have answers to the following question.

\begin{ques}
A homomorphism $f$ being fixed, if we consider discretizations of $f$ with respect to the image of a sequence of grids by a generic conservative homeomorphism, is it possible to retrieve the classes of recurrence by chains of $f$, or even better, the transition graph between recurrence classes? What about a practical point of view? What happens on each recurrence class?
\end{ques}

\section{Characteristic time}

In general, when a computer is asked to calculate the orbit of a point $x$ by a map $f$, it works in double precision, with $52$ binary digits. If one digitally iterates $f$ a reasonable number of times $t$, it is very unlikely that the orbit of $x$ falls into a periodic orbit of the discretization in time smaller than $t$, that is that there is a time $t'<t$ such that $f_N^{t'}(x_N) = f_N^t(x_N)$. Indeed, for a random map of a set with $2^{52}$ elements, this time $t'$ is typically of the order of $84$ million (see e.g. Theorem 2.3.1 of \cite{Mier-dyna}). However, the results obtained in this manuscript concern precisely this recurring event; so a priori they cannot explain what happens in practice.

\begin{center}
\begin{tikzpicture}[every text node part/.style={align=center}]
\draw[->,>=stealth] (7,0) -- (12,0);
\draw (0,0) -- (5,0);
\draw[dashed] (5,0) -- (7,0);
\draw (0,0) node{$|$};
\draw (0,-.5) node{$t=0$};
\draw (3.5,0) node{$|$};
\draw (3.5,-.5) node{$t_1$};
\draw (3.5,-.9) node{small};
\draw[decorate,decoration={brace}]
(0,.3) -- (3.5,.3) node[above,pos=0.5] {orbits of $f$ and\\ $f_N$ are close};
\draw (8.5,0) node{$|$};
\draw (8.5,-.5) node{$t_2$};
\draw (8.5,-.9) node{stabilisation time};
\draw (6,.3) node[above] {???};
\draw (12,-.5) node{$\infty$};
\draw[decorate,decoration={brace}]
(8.5,.3) -- (12,.3) node[above,pos=0.5] {dynamics of\\ discretization};
\end{tikzpicture}
\end{center}

The timeline above concerns the temporal behaviour of a given discretization $f_N$. It shows a time $t_1$ corresponding to the time up to which the orbit of any starting point $x$ under $f_N$ shadows that of $x$ under $f$. The theoretical estimates we have for the time up to which there is a ``strong shadowing'' are rather bad. For example, if the map $f$ is Lipschitz, this time $t_1$ is logarithmic in the size of the grid. In practice, when performing simulations, one iterate much longer than this usually rather short time. The time $t_2$ is the stabilization time of $f_N$, that is to say the time from which the orbit of any point on the grid fall into the recurrent set of $f_N$. From this time, the dynamics of $f_N$ is strictly periodic; one can consider that it is only from $ t_2 $ that the dynamics of the discretization emerges.

This leaves a long time interval $[t_1, t_2]$ where we do not really know how evolve the discretizations.

In this manuscript, we are interested primarily in the \emph{dynamics} of discretizations. The finite time behavior of discretizations has already been studied by P.P.~Flockermann (see Section~\ref{Floque}). It would be extremely interesting to determine an ``intermediate''  characteristic time $t_1\ll t_1'\ll t_2$, until where we see for sure the dynamics of the initial map, without suffering recurrence phenomena induced by the discretizations. In this subject, O.E.~Lanford has stated the following conjecture.

\begin{conjec}[Lanford]
Let $f : \Sp^1\to\Sp^1$ be a generic $C^2$ expanding map of the circle. If $\Leb_N$ denotes the uniform measure on $E_N$ and $\mu$ the unique physical measure of $f$, then $(f_N^*)^k(\Leb_N) \rightharpoonup \mu$ when $N,k$ tend to $+\infty$ with $\ln N \ll k \ll \sqrt N$.
\end{conjec}

In this manuscript, we prove a theorem in this sense (see Theorem~\ref{RepLanfordConlc}) but unfortunately this only concerns times $\ln N = O(k) $. I hope that the techniques presented in Chapter~\ref{Stat0} will allow to get more results in the direction of this conjecture.

\section{Comparison with random maps}\label{AppAlea}

In this thesis, we (almost) do not compare the dynamics of discretizations with that of a typical random map of a set with $q$ elements. Here, ``typical'' means the following. An integer $q$ being fixed, one considers all maps from a set with $q$ elements into itself; this finite set is endowed with the uniform probability. One then wonders what properties are satisfied by most of these applications. This type of questions was studied by (among others) P.~Erd\H{o}s and P.~Tur\'an (see for example \cite{Erdo-prob1} and \cite{Erdo-prob2}), an overview of these results can be found in \cite{Boll-rand}. For example, the degree of recurrence of such applications (recall that the degree of recurrence is the ratio between the cardinality of the recurrent set and $q$) behaves like $q^{-1/2}$. Of course, the study we have conducted in the case of discretizations of generic conservative homeomorphisms shows that such discretizations do not behave at all like random maps, even from an experimental point of view. By cons, what happens to discretizations of $C^r$-diffeomorphisms, with $r\ge 1$, is not clear at all. Unfortunately, the precise study of the combinatorics of such discretizations seems for the moment out of reach.

\section{Spatial discretizations versus stochastic perturbations versus multivalued maps}

In this manuscript, we do not treat the problems of stochastic perturbations or multivalued discretizations.

The first problem is the following: consider a dynamics $f$, and look at the sequences $(x_k)_{k\ge 0}$, where $x_{k+1}$ is obtained by taking a small random perturbation of $f(x_k)$. This idea of ``small random perturbation'' can be formalized in many ways, for example one can choose to take the image of $x_k$ by a map randomly chosen from those that are close to $f$, or take a point randomly in a neighbourhood of $f(x_k) $, etc. The reader could consult the lecture notes \cite{VianaStoch} of M.~Viana to get an idea of the variety of results in this area (e.g. in the case of expanding maps \cite{MR884892},\cite{MR874047}, \cite{MR685377}, or uniformly hyperbolic maps \cite{MR874047}, \cite{MR857204}, etc.).

In general, these results need the dynamics to be quite regular (typically at least $C^{1+\alpha}$): a conventional method to get them is to use some transfer operators related to the dynamics. They all suggest that the physical dynamics of stochastic perturbations converges towards the physical dynamics of the initial map. As already said, what happens for discretizations of generic systems is really different: after a while, the orbit of each point of the grid falls into a periodic orbit. This recurring phenomenon destroys any hope that the discretizations behave like stochastic perturbations during long time scales. This may affect the ergodic dynamics of discretizations, as there is no uniform convergence in Birkhoff theorem. This does not prevent the discretizations of some dynamics being  close stochastic perturbations in the short term, as explained by P.P.~Flockermann in his thesis \cite{Flocker} in the case of circle expanding maps. However, P.P.~Flockermann does not get any practical estimate on how long this behaviour is true.

\begin{ques}
How long the discretizations of generic circle expanding maps on the circle, or generic diffeomorphisms, behave like stochastic perturbations?
\end{ques}

\bigskip

In this manuscript, we do not treat either approximation of the dynamics by multivalued maps. This viewpoint can be quite interesting to determine which dynamical properties of a system can be detected by considering discrete approximations of the dynamics, see for instance \cite{MR2776399} (see also \cite{MR1403460} et \cite{MR1387977}). A priori, approaching a dynamics by a multivalued map, one lose much less information than approaching it by a finite map. Despite this, for generic conservative homeomorphisms, we show that in fact we can retrieve every dynamical feature of the system by considering all its discretizations (see for example Theorem~\ref{TheoMesConcl}). Therefore, in theory, one does not lose more information in discretizing than in considering multivalued maps. In fact, in Chapter~\ref{ChapCons}, one implicitly uses multivalued maps: the proof of Lax theorem consists in extracting an actual map from a multivalued map. The rest of Chapter~\ref{ChapCons} uses the same type of arguments: the candidates for discretizations are contained in a multivalued map; then Baire theorem allows us to choose among all possible extractions those with properties that one wants to see appearing an infinite number of times on discretizations.

However, in the case of multivalued maps, it may be possible to get actual estimates of the order of discretization up to which we have to go to be able to retrieve dynamical invariants with a given accuracy. This is clearly not the case for discretizations of generic systems, at least with the methods used in this thesis. Note that the fact of considering multivalued maps poses additional practical problems: for example, it is much longer to digitally detect the cycles for a multivalued map than for a classical map.

\newpage
\appendix

\addtolength{\cftfigurenumwidth}{8pt}
\listoffigures

\printindex

\small
\nocite{Guih-discr}
\nocite{Guih-rot}

\bibliographystyle{smfalpha}
\bibliography{../../Biblio}

\end{document}